%% file: main.tex
\documentclass[11pt]{article}
\usepackage[T1]{fontenc}
\usepackage{lmodern}
\usepackage[margin=1in]{geometry}
\usepackage{amsmath,amssymb,amsthm,mathtools}
\usepackage{booktabs,array}
\usepackage{microtype}
\usepackage{tikz}
\usepackage{xurl}
\usepackage[colorlinks=true,linkcolor=blue!55!black,citecolor=blue!55!black,urlcolor=blue!55!black]{hyperref}
\newcommand{\papertopic}{Gaussian Belief Propagation and Periodic Jacobi Matrices on Trees}
\newcommand{\papersubtitle}{Spectral Curve, Antibound Poles, and Fixed Points}
\hypersetup{pdftitle={\papertopic: \papersubtitle},pdfauthor={Yusuke Watanabe}}

\newtheorem{maintheorem}{Theorem}
\newtheorem{theorem}{Theorem}[section]
\newtheorem{lemma}[theorem]{Lemma}
\newtheorem{proposition}[theorem]{Proposition}
\newtheorem{corollary}[theorem]{Corollary}
\newtheorem{conjecture}[theorem]{Conjecture}
\theoremstyle{definition}
\newtheorem{definition}[theorem]{Definition}
\newtheorem{example}[theorem]{Example}
\newtheorem{problem}[theorem]{Problem}
\theoremstyle{remark}
\newtheorem{remark}[theorem]{Remark}
\numberwithin{equation}{section}

\DeclareMathOperator{\diag}{diag}
\DeclareMathOperator{\spec}{spec}
\DeclareMathOperator{\Res}{Res}
\DeclareMathOperator{\ord}{ord}
\DeclareMathOperator{\rank}{rank}
\DeclareMathOperator{\sgn}{sgn}
\DeclareMathOperator{\Hess}{Hess}
\DeclareMathOperator{\ind}{ind}
\DeclareMathOperator{\tr}{tr}
\DeclareMathOperator{\lcm}{lcm}
\DeclareMathOperator{\im}{Im}
\DeclareMathOperator{\re}{Re}
\newcommand{\R}{\mathbb R}
\newcommand{\C}{\mathbb C}
\newcommand{\Z}{\mathbb Z}
\newcommand{\N}{\mathbb N}
\newcommand{\PP}{\mathbb P}
\newcommand{\ones}{\mathbf 1}
\newcommand{\abs}[1]{\lvert#1\rvert}
\newcommand{\norm}[1]{\lVert#1\rVert}
\newcommand{\T}{\mathsf T}                 
\newcommand{\tH}{\widetilde H}             
\newcommand{\tree}{\mathcal T}             
\newcommand{\sigm}{\sigma_-}               
\newcommand{\sigp}{\sigma_+}               
\newcommand{\Sigm}{\Sigma_-}               
\newcommand{\Sigp}{\Sigma_+}               
\newcommand{\cL}{\mathcal L}               
\newcommand{\LG}{L(G)}                     
\newcommand{\cZ}{\mathcal Z}               
\newcommand{\Xc}{X}                        
\newcommand{\Cphys}{C_{\mathrm{phys}}}     
\newcommand{\Xp}{X_+}                      
\newcommand{\Vphys}{V^{\mathrm{phys}}}
\newcommand{\taups}{\tau_{\mathrm{ps}}}    
\title{\papertopic:\\\papersubtitle}
\author{Yusuke Watanabe\thanks{\href{mailto:ywatanabex@gmail.com}{\texttt{ywatanabex@gmail.com}}}}
\date{}

\begin{document}
\maketitle
\begin{abstract}
Diagonal Green functions of a periodic Jacobi operator on a universal
covering tree satisfy the scalar Gaussian belief propagation covariance
equations on the finite connected loopless base multigraph.
At cycle rank at least two, the Green functions continue to simple
antibound poles at both Perron values for arbitrary real potentials
and nonzero real weights. This proves the conjecture of Christiansen,
Simon and Zinchenko. For positive weights the pole result includes
self-loops. On loopless graphs the lower pole is the walk-summability
threshold. Each energy strictly between it and the spectral bottom has
at least two fixed points with positive definite local beliefs. On a
connected simple graph, some coefficients yield a compact component of
positive solutions exactly when the graph contains two vertex-disjoint
connected subgraphs of cycle rank at least two. The covariance curve is
generically irreducible, with one sheet per edge subset with no vertex of
degree one. For cycle rank at least two and generic real coefficients, the
diagonal Green functions, infinite half-tree $m$-functions and the
density-of-states Stieltjes transform all have this degree, answering a
question of Avni, Breuer and Simon. Then all eigenvalues of the signed
Jacobi matrices are antibound poles. Signed principal submatrices
classify finite points with nonzero residue of $(\sum_i V_i) dz$,
where $V_i$ are the vertex variances. We give a genus
formula. Near zero potentials and unit weights on $K_4$, independently
varying the coefficients leaves the cover spectrum an interval. This
refutes Avni, Breuer and Simon's conjecture on the density of parameters
with all gaps open. In block models, generic irreducibility and a
uniqueness criterion extend, and at uniform dimension $p$ the degree can
exceed $L^p$ for scalar degree $L$. Part of the scalar theory is formally
verified in Lean 4 with Mathlib.
\end{abstract}

\tableofcontents

\part{Scalar variables}
\input{sections/introduction}
\input{sections/setup}
\input{sections/degree}
\input{sections/local-poles}
\input{sections/perron}
\input{sections/irreducibility}
\input{sections/classification}
\input{sections/real}
\input{sections/genus}
\input{sections/free-products}
\input{sections/degenerations}
\input{sections/abs}

\part{Vector-valued variables}
\input{sections/block-real}
\input{sections/block-complex}

\part{Formalization and open problems}
\input{sections/formalization}
\input{sections/open-problems}

\part{Appendices}
\appendix
\input{appendices/many-solutions}
\input{appendices/weak-bridge}
\input{appendices/bowtie}
\input{appendices/certificates}

\input{sections/declarations}

\input{references}
\end{document}

%% file: sections/introduction.tex
\section{Introduction}
\label{sec:intro}

\subsection{One curve, three correspondences}

Let $G=(V,E)$ be a finite connected graph with a real potential $D_i$ on
each vertex and a nonzero real weight $w_e$ on each edge, and let $\tH$ be
the periodic Jacobi operator that these data define on the universal
covering tree of $G$. Its diagonal Green function
$G_i(z)=((\tH-z)^{-1})_{ii}$ is an algebraic function of~$z$
\cite{Aomoto84,ABS}. This paper studies the algebraic curve on which it
lives, and it does so through a system of equations that comes from
statistical inference.

Gaussian belief propagation is the message-passing algorithm for
marginals of a Gaussian graphical model
\cite{WeissFreeman,MJW}. Its covariance fixed-point equations for the
quadratic form $\tfrac12x^{\T}(D+A-z)x$, where $A$ is the weighted
adjacency matrix, can be written as
\begin{equation}\label{eq:intro-system}
(D_i-z)\,V_i=1+\sum_{e\ni i}f_e\quad(i\in V),\qquad
f_e(1+f_e)=w_e^2\,V_iV_j\quad(e=ij\in E),
\end{equation}
with one unknown $V_i$ per vertex, the variance of the belief at $i$, and
one unknown $f_e$ per edge (Section~\ref{sec:setup}). A solution is
\emph{admissible} if all $V_i$ are nonzero and $z\notin\{D_i\}$, and a real
solution is \emph{positive} if all $V_i$ and all $f_e$ are positive; the
positive solutions are exactly the fixed points whose vertex and edge
beliefs are positive definite. The Green function is one solution: for
real $z$ below the spectrum of $\tH$, $V_i=G_i(z)$ is a positive solution,
the \emph{physical branch} (Proposition~\ref{prop:physical-branch}). The
admissible solutions of \eqref{eq:intro-system}, as $z$ varies, determine
the \emph{spectral curve} $\Xc=\Xc(G;D,w)$ by taking the Zariski closure
and retaining the components on which $z$ varies. For generic real
coefficients, away from a finite exceptional set of energies, every
admissible covariance fixed point is reached by analytic continuation of
the physical branch in the full $(V,f)$ coordinates. Write
$\omega=(\sum_i V_i)\,dz$. We describe the curve by three correspondences.

\begin{center}
\begin{tabular}{@{}lll@{}}
\toprule
On the curve & Corresponds to & Result \\
\midrule
sheets (points over $z=\infty$) & edge subsets with no vertex of degree one & Theorem~\ref{thm:main-sheets} \\
finite points with $\Res\omega\ne0$ & eigenvalues of signed principal submatrices & Theorems~\ref{thm:main-perron},~\ref{thm:main-poles} \\
real points with $V>0$, $f>0$ & fixed points with positive definite beliefs & Theorem~\ref{thm:main-real} \\
\bottomrule
\end{tabular}
\end{center}

\noindent The signed matrices are the matrices $D+A^s$ obtained from the
finite Jacobi matrix $D+A$ by changing the signs of the weights on a subset
of the edges; equivalently, $D+A^s$ is the part of the Jacobi matrix of a
two-fold cover of $G$ that is antisymmetric under the deck transformation.
The three correspondences meet at one point. When the cycle rank of $G$ is
at least two, the smallest eigenvalue over all signings,
$\sigm=\lambda_{\min}(D-\abs A)$, is at the same time a pole of the Green
function on a nonphysical sheet, the threshold of walk-summability of the
Gaussian model, and the end of a real arc of positive fixed points that
starts on the physical branch (Theorem~\ref{thm:main-perron}).

\subsection{Main results}

We write $n=\abs V$, $\beta=\abs E-n+1$ for the cycle rank of a connected
graph, $\Sigm$ and $\Sigp$ for the bottom and the top of the spectrum of
$\tH$, and
\[
\sigm=\lambda_{\min}(D-\abs A),\qquad\sigp=\lambda_{\max}(D+\abs A)
\]
for the \emph{Perron values} of the finite graph, with positive unit
eigenvectors $\bar u^-$ and $\bar u^+$. Graphs are finite multigraphs
without self-loops unless self-loops are explicitly allowed; parallel
edges are distinct edges.

Christiansen, Simon and Zinchenko \cite{CSZ} consider a finite connected
leafless graph with positive edge weights whose universal covering tree
is different from the line. They conjectured that each Green function
continues to a higher sheet on which it has a pole at the largest
eigenvalue of the finite Jacobi matrix, and, when the graph is bipartite,
also at its smallest eigenvalue. The first theorem proves their
conjecture and identifies a lower pole for nonbipartite graphs as well.

\begin{maintheorem}[Perron antibound poles]\label{thm:main-perron}
Let $G$ be a finite connected multigraph without self-loops and of cycle
rank $\beta\ge2$, with arbitrary real potentials and nonzero real edge
weights.
\begin{enumerate}
\item[(a)] \emph{Gap.} $\sigm<\Sigm$ and $\Sigp<\sigp$.
\item[(b)] \emph{Poles.} The continuation of each Green function $G_i$
along a suitable path in the $z$-plane avoiding finitely many points
reaches a branch that is single valued near $z=\sigp$ and has there a
simple pole with residue $(\abs E-n)(\bar u_i^+)^2$ in the variable $z$;
the same holds at $z=\sigm$ with $\bar u^-$. These poles do not lie on the
physical sheet. For positive weights the statement also holds for graphs
with self-loops, with the convention that a loop of weight $a$ contributes
$2a$ to the diagonal of the finite Jacobi matrix
(Proposition~\ref{prop:self-loops}).
\item[(c)] \emph{Threshold and real arc.} The comparison matrix
$D-z-\abs A$ of the model at parameter $z$ is positive definite exactly
when $z<\sigm$, and this is equivalent to positive definiteness of all
two-fold covers of the model. Positive solutions of
\eqref{eq:intro-system} exist only for $z\le\Sigm$. For every
$z\in(\sigm,\Sigm)$ there are at least two positive solutions: the
physical one, and one on a real arc of positive solutions that leaves the
physical branch at its turning point over $z=\Sigm$ and runs into the pole
at $\sigm$.
\end{enumerate}
For $\beta=1$ and a loopless graph one has $\sigm=\Sigm$ and $\sigp=\Sigp$.
Each Green function has inverse-square-root behavior there for arbitrary
real potentials and nonzero weights, including when trees are attached
to the cycle (Proposition~\ref{prop:unicyclic}).
\end{maintheorem}

For positive weights and a bipartite graph, $\sigm=\lambda_{\min}(D+A)$, so
(b) contains Conjecture~3.4 of \cite{CSZ}; the lower pole for nonbipartite
graphs, the residues, and graphs with leaves are additions. The gap (a) is
the theorem of Sy and Sunada \cite{SySunada} for simple graphs. We know of
no earlier proof of the conjecture. Part (c) concerns the local beliefs:
the finite model at parameter $z$ is itself positive definite only for
$z<\lambda_{\min}(D+A)$, and for a bipartite graph with positive weights
this excludes the whole gap. The equivalence between walk-summability and
positive definiteness of all two-fold covers is due to Ruozzi and Tatikonda
\cite{RuozziTatikonda}. The condition $\beta\ge2$ cannot be replaced by the
condition that the covering tree is not the line: a cycle with pendant
trees has $\beta=1$.

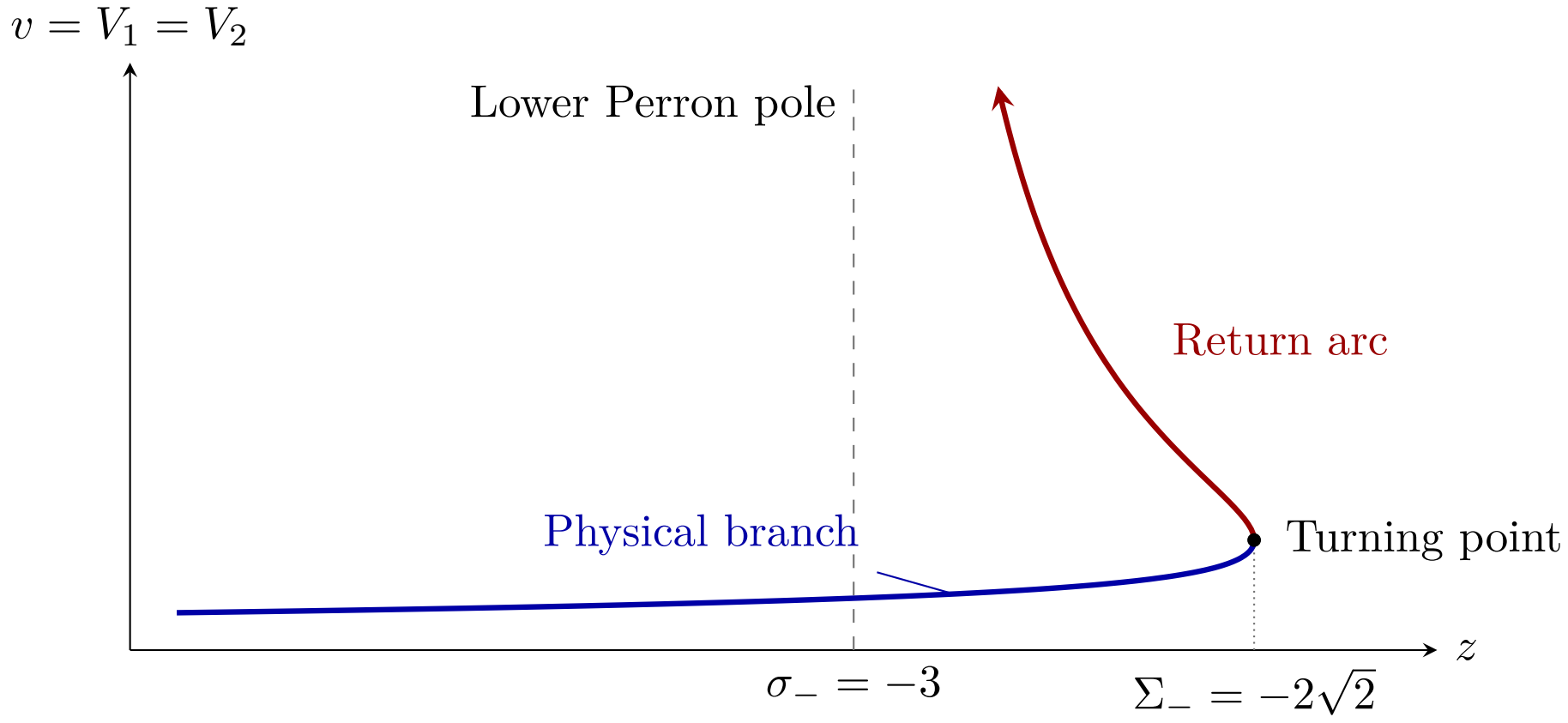
\begin{figure}[tbp]
\centering
\begin{tikzpicture}[x=18cm,y=0.60cm,>=stealth,font=\small]
  \pgfmathsetmacro{\edgeposition}{3-2*sqrt(2)}
  \pgfmathsetmacro{\turnvariance}{sqrt(2)}
  \draw[->] (-0.31,0) -- (0.25,0) node[right] {$z$};
  \draw[->] (-0.31,0) -- (-0.31,7.55) node[above] {$v=V_1=V_2$};
  \draw[black!55,dashed] (0,0) -- (0,7.25);
  \draw[black!55,densely dotted] (\edgeposition,0)
    -- (\edgeposition,\turnvariance);
  \node[below] at (0,0) {$\sigm=-3$};
  \node[below] at (\edgeposition,0) {$\Sigm=-2\sqrt{2}$};
  \draw[very thick,blue!65!black]
    plot[domain=0.48:\turnvariance,samples=85,variable=\vplot]
    ({3+(1-3*sqrt(1+4*\vplot*\vplot))/(2*\vplot)},\vplot);
  \draw[very thick,red!60!black,->]
    plot[domain=\turnvariance:7.25,samples=85,variable=\vplot]
    ({3+(1-3*sqrt(1+4*\vplot*\vplot))/(2*\vplot)},\vplot);
  \fill (\edgeposition,\turnvariance) circle[radius=1.5pt]
    node[right=3pt] {Turning point};
  \node[above left,blue!65!black] at (0.01,1.0) {Physical branch};
  \draw[blue!65!black,thin] (0.01,1.0) -- (0.043,0.72);
  \node[right,red!60!black] at (0.129,4.0) {Return arc};
  \node[left] at (0,7.0) {Lower Perron pole};
\end{tikzpicture}
\caption{The arc of positive solutions for the two-vertex graph with
three unit-weight edges and zero potentials. For the common variance
$v=V_1=V_2$, the equations give $z=(1-3\sqrt{1+4v^2})/(2v)$.
Increasing $v$ takes the physical branch to the turning point at $\Sigm$
and then follows the return arc as $z$ approaches $\sigm$ from above.}
\label{fig:perron-arc}
\end{figure}

Figure~\ref{fig:perron-arc} illustrates the arc in part~(c) on the
three-parallel-edge graph.

An edge subset $S\subseteq E$ is \emph{leafless} if no vertex is incident
to exactly one edge of~$S$. Let $\LG$ be the number of leafless subsets.
This is a known graph invariant, the number of sub-coregraphs or
generalized loops \cite{WatanabeFukumizuZeta,ChertkovChernyak}; it
satisfies $2^\beta\le\LG$ for connected $G$, with further bounds recalled
in Proposition~\ref{prop:L-bounds}.

\begin{maintheorem}[Sheets]\label{thm:main-sheets}
Let $G$ be a finite multigraph without self-loops.
\begin{enumerate}
\item[(a)] For all coefficients, the solutions of \eqref{eq:intro-system}
near $z=\infty$ form $2^{\abs E}$ branches, analytic in $1/z$ and indexed by
the edge subsets $S$, with $f_e\to-1$ for $e\in S$ and $f_e\to0$ otherwise;
the branch of $S$ is admissible if and only if $S$ is leafless. Hence the
spectral curve has degree $\LG$ over the $z$-line. For connected $G$ with
real coefficients the physical branch corresponds to $S=\emptyset$.
\item[(b)] For generic coefficients the spectral curve is irreducible, and
for a general value of $z$ the system \eqref{eq:intro-system} has exactly
$\LG$ admissible solutions, all simple. For connected $G$ and generic real
coefficients, continuation of the physical branch reaches every such
solution over energies outside a finite set.
\item[(c)] Let $G$ be connected with $\beta\ge2$. For generic real
coefficients each Green function $G_i$, each $m$-function of an infinite
half-tree, and the Stieltjes transform of the density of states are
algebraic functions of degree exactly $\LG$ over $\C(z)$.
\item[(d)] For all coefficients the number of nondegenerate admissible
solutions at a fixed $z$ is at most $\LG$.
\end{enumerate}
\end{maintheorem}

Avni, Breuer and Simon proved algebraicity of the Green function for
periodic Jacobi matrices on trees, bounded the number of sheets by
$2^{2\abs E}$, observed that all computed examples have two sheets, and
asked whether this is always so \cite[Problem~10.2]{ABS}. Part (c) answers
this: generically the number of sheets is $\LG\ge2^\beta\ge4$ as soon as
$\beta\ge2$, for example $5$ for two vertices joined by three edges and
$15$ for $K_4$. The two-sheeted examples are explained by symmetry: for
uniform coefficients on an edge-transitive graph the physical component is
a two-sheeted curve (Proposition~\ref{prop:uniform}). The count in (b) is the Gaussian
analogue of the generic number of fixed points of binary belief
propagation \cite{WatanabeCounting}, which is a sum over the same leafless
subsets with a different vertex weight.

Conjecture~\ref{open:monodromy} asks whether the generic monodromy
of this scalar curve is the full symmetric group on its $\LG$ sheets.

\begin{maintheorem}[Poles]\label{thm:main-poles}
Let $G$ be a finite connected multigraph without self-loops and with
$\beta\ge2$.
\begin{enumerate}
\item[(a)] \emph{Local existence, all coefficients.} Let $s$ be a signing,
$\lambda$ a simple eigenvalue of $D+A^s$ with $\lambda\notin\{D_i\}$ and
with an eigenvector $u$, $u^{\T}u=1$, without zero entries. Then the
spectral curve has a unique unramified branch near $z=\lambda$ with the
following vertex and edge expansions:
\[
V_i=\frac{(\abs E-n)u_i^2}{z-\lambda}+O(1),\qquad
f_e=-\frac{(\abs E-n)s_e w_e u_i u_j}{z-\lambda}+O(1)\quad(e=ij).
\]
\item[(b)] \emph{Generic coefficients.} For generic real coefficients, for
every signing $s$ and every eigenvalue $\lambda$ of $D+A^s$, the
continuation of the Green functions reaches the branch of (a). These
$n2^\beta$ poles, one for each switching class and eigenvalue, have
positive residues and are antibound. Conversely, a continuation of some
$G_i$ can be unbounded near a finite point $z_0$ only if $z_0$ is an
eigenvalue of a signed principal submatrix $(D+A^s)[U]$ with $G[U]$
connected.
\item[(c)] \emph{Classification.} For generic coefficients, the points of
the normalization of the spectral curve over finite $z$ at which the
differential $\omega=(\sum_iV_i)\,dz$ has a nonzero residue are exactly
the poles supported on vertex subsets $U$ described in
Theorem~\ref{thm:partial-pole}; their number and their residues are given
by Theorem~\ref{thm:classification}.
\item[(d)] \emph{Floquet identity.} On every sheet and for all
coefficients, $\omega=-d\log\Phi$ with
$\Phi=\prod_e(1+f_e)/\prod_iV_i$. All residues of $\omega$ are integers,
and for generic coefficients their sum over the finite points is
$\sum_{S\in\cL(G)}(2\abs S-n)$.
\item[(e)] \emph{Cycle rank two.} If $\beta=2$, then for generic real
coefficients the only points with nonzero residue of $\omega$ that lie
outside the convex hull $[\Sigm,\Sigp]$ of the spectrum of $\tH$ are the
two Perron poles.
\end{enumerate}
\end{maintheorem}

On the physical branch $\Phi$ is the Floquet function of
\cite{BBGSS}. Part (b) answers a question raised after the conjecture in
\cite{CSZ}, which asks for families of finite Jacobi matrices whose only
higher-sheet poles are at $\pm\sigma$: no such family exists among generic
coefficients. Reading the sign of a residue to separate physical from
nonphysical poles goes back, for two special cases of free products, to
Gutkin \cite{Gutkin}; poles on the physical sheet are the point spectrum
studied in \cite{Aomoto91,BGVM}.

For the last theorem we fix $z$ and absorb it into $D$.

\begin{maintheorem}[Positive fixed points]\label{thm:main-real}
Consider the positive solutions of \eqref{eq:intro-system} at a fixed real
$z$, that is, the fixed points of Gaussian belief propagation with
positive definite local beliefs, for the comparison matrix
$M=D-z-\abs A$.
\begin{enumerate}
\item[(a)] If $M$ is positive definite, there is exactly one positive
solution, the physical one.
\item[(b)] Suppose that no principal submatrix of $M$ on a nonempty
vertex set inducing a connected subgraph has smallest eigenvalue zero.
Then the set of positive solutions is compact. The Brouwer degree of the
gradient of the reduced Bethe function at zero, taken on a bounded open
set containing all positive solutions and with closure in the positive
orthant, is $1$ if $M$ is positive definite and $0$ otherwise. When the
solutions are finite, the degree is the sum of their local indices. If
all positive solutions are nondegenerate, their number is odd in the
first case and even in the second. Under these assumptions there is
exactly one positive solution if and only if the model is walk-summable.
\item[(c)] Let $G$ be simple and connected. If a positive solution exists,
there is a maximal one, with the largest cavity precisions. For $z<\Sigm$
it is the physical branch. At every other positive solution the Hessian of
the reduced Bethe function has a negative eigenvalue. At any positive
solution, the sum of the Morse index and nullity is at most the pseudoforest
deletion number of~$G$.
\end{enumerate}
\end{maintheorem}

Convergence of the algorithm for walk-summable models \cite{MJW}, the
criterion for the existence of a positive fixed point and convergence to
the maximal one \cite{SuWu14}, and uniqueness inside smaller admissible
sets \cite{MoallemiVanRoy,DuEtAl} are known. Theorem~\ref{thm:main-real}
counts all positive solutions and gives the converse direction.
Theorem~\ref{real:thm:packing} characterizes when a circle can occur on a
connected simple graph. If there are no two vertex-disjoint connected
subgraphs of cycle rank at least two, each energy has at most two positive
solutions and there are no circles for any coefficients. If there are
two such subgraphs, some coefficients give a regular circle and at
least four positive solutions at one energy. Section~\ref{sec:real}
also gives families with pseudoforest deletion number two in which the
number of solutions is unbounded, and describes the positive solution
set for generic coefficients and $\beta\ge2$ as one arc and finitely
many circles.

\subsection{Open problems addressed}

Table~\ref{tab:problems} lists the open problems from
\cite{CSZ,ABS,AomotoKato} on which the paper has something to say.

\begin{table}[htbp]
\centering
\small
\begin{tabular}{@{}>{\raggedright\arraybackslash}p{4.3cm}>{\raggedright\arraybackslash}p{8.1cm}>{\raggedright\arraybackslash}p{3.0cm}@{}}
\toprule
Problem & Result here & Where \\
\midrule
\cite[Conjecture~3.4]{CSZ}: poles at the Perron values & Proved under the original leafless hypothesis, with simple poles and explicit residues & Theorem~\ref{thm:main-perron}(b) \\
\cite{CSZ}, question after Conj.~3.4: families with poles only at $\pm\sigma$ & For generic coefficients every eigenvalue of every signing is a pole & Theorem~\ref{thm:main-poles}(b) \\
\cite[Problem~10.2]{ABS}: are the Green function and the $m$-functions two-sheeted? & Generically $\LG$ sheets & Theorem~\ref{thm:main-sheets}(c) \\
\cite[Problem~10.1]{ABS}: are all branch points of square-root type? & No (a known cube-root example, placed on the curve) & Section~\ref{sec:degenerations} \\
\cite[Problem~10.3]{ABS}: correspondence of poles of half-tree functions & The correspondence valid on the line fails at branching vertices; a local criterion and examples. The global question stays open & Section~\ref{sec:abs} \\
\cite[Conjectures~10.4, 10.5]{ABS}: gaps open generically & An open set of coefficients on $K_4$ with no gap; the exceptional set has codimension zero & Theorem~\ref{thm:gapless} \\
\cite[Conjectures~B.4, B.5]{ABS} (arXiv version): Borg-type statements & Counterexamples in the same open set & Corollary~\ref{cor:abs-conjectures} \\
\cite{AomotoKato}: irreducibility of the curve of all branches for free products of cyclic groups & Irreducible for generic weights; reducible for three factors $\Z/2$ with equal weights & Theorem~\ref{thm:aomoto-kato} \\
\bottomrule
\end{tabular}
\caption{Open problems addressed in this paper. The numbering of the
conjectures of \cite{ABS} follows the arXiv version.}
\label{tab:problems}
\end{table}

\subsection{An example}

Let $G=D_3$ be the graph with two vertices joined by three parallel edges
of positive weights $w_1,w_2,w_3$, with potentials $D_1,D_2$. It has
$\beta=2$ and five leafless edge subsets, namely the empty set, the three
pairs of edges, and the full set. For generic coefficients the spectral
curve is irreducible with five sheets and genus six. There are four
switching classes of signings, and Theorem~\ref{thm:main-poles} gives the
eight antibound poles
\[
z=\frac{D_1+D_2}2\pm\sqrt{\Bigl(\frac{D_1-D_2}2\Bigr)^2+(w_1\pm w_2\pm w_3)^2},
\]
each with residue $u_i^2$ at the vertex $i$, where $u$ is the unit
eigenvector. The choice $w_1+w_2+w_3$ for the signed sum gives the two
Perron poles of Theorem~\ref{thm:main-perron}. When the three weights are
equal, the three nontrivial switching classes give the same pair of
eigenvalues and the physical component has two sheets. The only finite
points on its normalization where $(V_1+V_2)\,dz$ has nonzero residue
are the two Perron poles. If $D_1\ne D_2$, there is also a branch point
over each $D_i$, where $V_i$ diverges as a nonzero multiple of
$(D_i-z)^{-1/2}$ and $V_i\,dz$ has zero residue
(Proposition~\ref{prop:uniform}). The model in \cite[Example~3.3]{CSZ}
has equal weights and potentials $b,-b$. For $D_1=D_2=0$ the curve is
the curve of Aomoto and Kato for the free product of three copies of
$\Z/2$ \cite{AomotoKato}. The example is worked out in
Section~\ref{free:subsec:D3}.

\subsection{Further results}

\emph{Genus.} Section~\ref{sec:genus} computes the genus of the spectral
curve for generic coefficients as an explicit combinatorial sum over
the pole configurations. The cycle $C_n$ gives the classical hyperelliptic curve of genus
$n-1$; $D_3$ has genus $6$ and $K_4$ has genus $94$.

\emph{Free products.} Section~\ref{sec:free} treats the free convolution
of finitely supported measures and the nearest-neighbor operators on free
products of cyclic groups, where the curve of all branches was written
down by Aomoto and Kato \cite{AomotoKato} and the sign of the residues was
used by Gutkin \cite{Gutkin}. We prove generic irreducibility, give the
exact degree and the residues at sums of atoms, and read poles supported
on subgraphs through operator-valued free convolution.

\emph{Degenerations.} Section~\ref{sec:degenerations} shows that for
unweighted graphs the degree can drop below the count for the degree
refinement, that attaching a rooted tree to every vertex pulls the curve
back along a rational map without losing sheets, and that the complete
pole data of all vertices do not determine the gaps of the spectrum.

\emph{Problems of Avni, Breuer and Simon.} Section~\ref{sec:abs} contains
the results on the half-tree functions and on spectral gaps listed in
Table~\ref{tab:problems}.

\emph{Vector-valued variables.} Part~II extends the theory to block
models. Section~\ref{blr:subsec:scalar-comparison} compares the two
theories through the uniqueness threshold, the complex sheets and the
poles. The block curve is generically irreducible. Its degree can grow
strictly beyond the count supplied by independent scalar channels, as
the two-dimensional bowtie with $36>\LG^2=16$ shows
(Theorems~\ref{thm:block-irreducible} and~\ref{thm:block-degree}).

\subsection{Related work}

\emph{Periodic Jacobi matrices on trees.} The framework is that of
\cite{ABS,ABKS,CSZ,BBGSS}; see also \cite{GVK,BGVM}. Aomoto
\cite{Aomoto84} established algebraicity in the free-group setting and
counted $2^g$ real branches at large $z$ for a bouquet of $g$ loops.
The equations for more general trees are due to Aomoto~\cite{Aomoto88};
see his later account \cite[\S1]{Aomoto91} for the periodic case.
Aomoto and Kato \cite{AomotoKato} extended the calculation of the curve
and the branch count to free products of cyclic groups and left
irreducibility open. For a general finite graph we are not aware of
earlier results on the exact degree, on irreducibility, on the poles on
nonphysical sheets, or on the genus. For
periodic operators on $\Z^d$ the corresponding irreducibility statements
for Bloch varieties are classical and recent \cite{Kuchment,FaustLiu};
there irreducibility excludes embedded eigenvalues, whereas on trees an
irreducible curve can carry point spectrum on the physical sheet.

\emph{Gaussian belief propagation.} Besides the works cited after
Theorem~\ref{thm:main-real}, see \cite{JBD,CsekeHeskes,SuWu15,KSdP}. The
loop series of Chertkov and Chernyak expresses the determinant as a sum
over generalized loops around any fixed point \cite{ChernyakChertkov};
Theorem~\ref{thm:main-sheets} says that the same index set counts the
fixed points themselves.

\emph{Binary models.} The count $\LG$ and the role of the pseudoforest
deletion number parallel the results of
\cite{WatanabeCounting,WatanabeResponse} for binary belief propagation;
the comparison is made in Sections~\ref{sec:degree} and~\ref{sec:real}.

The attributions in this paper rest on the sources we could read. In a
few cases, listed where they occur, we could consult only an abstract, a
review, or a restatement by the same authors.

\subsection{Formal verification and organization}

The Lean~4 formalization includes the fixed-coefficient sheet count away
from exceptional energies, continuation from one physical germ to both
Perron poles, analytic normalization of the physical component, existence
of two positive solutions throughout the lower gap, and the signed degree
and parity results. Section~\ref{sec:formal} specifies the hypotheses and
the scope of the formal statements.

Section~\ref{sec:setup} fixes the setting. Sections~\ref{sec:degree}
and~\ref{sec:local-poles} contain the two local analyses, at $z=\infty$ and
at signed eigenvalues, on which everything else rests.
Section~\ref{sec:perron} proves Theorem~\ref{thm:main-perron},
Section~\ref{sec:irreducibility} proves Theorem~\ref{thm:main-sheets},
Section~\ref{sec:classification} proves Theorem~\ref{thm:main-poles}, and
Section~\ref{sec:real} proves Theorem~\ref{thm:main-real}.
Sections~\ref{sec:genus}--\ref{sec:abs} contain the further results on
scalar models, and Part~II treats vector-valued variables.
Section~\ref{sec:open} collects open problems. The appendices contain the
longer constructions and the computer-assisted parts:
families with many positive solutions (Appendix~\ref{app:many}), weak
bridges and compact components (Appendix~\ref{app:bridge}), the count of
sheets of the two-dimensional model on the bowtie
(Appendix~\ref{app:bowtie}), and a catalog of all computer-assisted
statements with their certificates (Appendix~\ref{app:cert}).

%% file: sections/setup.tex
\section{Setting: the covariance system and its physical branch}
\label{sec:setup}

\subsection{Graphs, signings and Perron values}
\label{subsec:graphs}

Throughout, $G=(V,E)$ is a finite multigraph. Parallel edges are distinct
elements of $E$. Unless stated otherwise $G$ has no self-loops; self-loops
are treated in Proposition~\ref{prop:self-loops}. We write $n=\abs V$, and
for connected $G$ the cycle rank is $\beta=\abs E-n+1$. For $S\subseteq E$
and $i\in V$, $d_S(i)$ is the number of edges of $S$ at $i$, and
$d_i=d_E(i)$. We write $e=ij$ for an edge with ends $i$ and $j$, and
$\sum_{e\ni i}$ for the sum over the edges at $i$, each parallel edge
counted separately.

The coefficients are a potential $D\in\R^V$ and edge weights
$w\in(\R\setminus\{0\})^E$; complex coefficients $D\in\C^V$,
$w\in(\C^*)^E$ are used in the algebraic statements. A \emph{signing} is
$s\in\{\pm1\}^E$. We set
\[
A_{ij}=\sum_{e=ij}w_e,\qquad A^s_{ij}=\sum_{e=ij}s_ew_e,\qquad
\abs A_{ij}=\sum_{e=ij}\abs{w_e},
\]
and call $D+A^s$ the \emph{signed Jacobi matrices} of $(G,D,w)$. Note that
$\abs A$ is the sum of the absolute values edge by edge; it differs from
the entrywise absolute value of $A$ when parallel edges carry weights of
both signs. For $\varepsilon\in\{\pm1\}^V$ the signing
$s'_e=\varepsilon_i\varepsilon_js_e$ satisfies
$D+A^{s'}=\diag(\varepsilon)(D+A^s)\diag(\varepsilon)$. Two signings
related in this way are \emph{switching equivalent}; a connected graph has
$2^\beta$ switching classes, and the spectrum of $D+A^s$ and the squares
$u_i^2$ of the entries of its eigenvectors depend only on the class of~$s$.
The matrix $(D+A)\oplus(D+A^s)$ is, after an orthogonal change of basis,
the Jacobi matrix of the two-fold cover of $G$ in which the edge $e$ is
lifted crosswise exactly when $s_e=-1$.

For connected $G$ and real coefficients the \emph{Perron values} are
\begin{equation}\label{eq:perron-values}
\sigm=\lambda_{\min}(D-\abs A),\qquad \sigp=\lambda_{\max}(D+\abs A).
\end{equation}
Both are simple eigenvalues, with eigenvectors $\bar u^{-},\bar u^{+}$ that
we take positive and of unit length. They are the extreme eigenvalues among
all signings: $\sigm\le\lambda\le\sigp$ for every eigenvalue $\lambda$ of
every $D+A^s$.

\subsection{The periodic Jacobi operator on the universal cover}
\label{subsec:cover}

Let $\pi\colon\tree\to G$ be the universal cover of $G$, a tree, infinite
when $\beta\ge1$. The coefficients lift to $\tree$, and
\[
(\tH\psi)(x)=D_{\pi(x)}\psi(x)+\sum_{\tilde e=xy}w_{\pi(\tilde e)}\psi(y)
\qquad(\psi\in\ell^2(\tree))
\]
is a bounded self-adjoint operator, the \emph{periodic Jacobi operator} of
$(G,D,w)$ in the sense of
\cite{ABS,CSZ}. On a tree the signs of the
weights can be removed by a diagonal unitary, so $\tH$ is unitarily
equivalent to the operator with weights $\abs{w_e}$. We write
\[
\Sigm=\inf\sigma(\tH),\qquad\Sigp=\sup\sigma(\tH).
\]
For $z\in\C\setminus\sigma(\tH)$ and $i\in V$ the \emph{Green function} is
$G_i(z)=\langle\delta_x,(\tH-z)^{-1}\delta_x\rangle$ for any $x\in\pi^{-1}(i)$.
For an edge $e=ij$, removing a lift of $e$ from $\tree$ leaves two
half-trees; let $m_{i|e}(z)$ be the Green function, at its root over $i$, of
the restriction of $\tH$ to the half-tree containing that root. (When $e$
is a loopless edge both orientations are determined by the end; for a
parallel class the label $e$ distinguishes the edges.) The Schur complement
formula on the tree gives, for $z\notin\R$ and for real $z$ outside
$[\Sigm,\Sigp]$,
\begin{equation}\label{eq:m-recursion}
\frac1{G_i}=D_i-z-\sum_{e\ni i,\ e=ij}w_e^2\,m_{j|e},\qquad
\frac1{m_{i|e}}=D_i-z-\sum_{\substack{e'\ni i,\ e'\ne e\\ e'=ij'}}w_{e'}^2\,m_{j'|e'} ,
\end{equation}
and therefore
\begin{equation}\label{eq:G-from-m}
G_i=\frac{m_{i|e}}{1-w_e^2\,m_{i|e}\,m_{j|e}}\qquad(e=ij).
\end{equation}

\subsection{Gaussian belief propagation and the covariance system}
\label{subsec:bp}

Gaussian belief propagation on $G$ for the quadratic form
$\tfrac12x^{\T}(D+A)x$ iterates one precision $q_{i\to j}$ for each
oriented edge. A fixed point satisfies the cavity equations
\begin{equation}\label{eq:cavity}
q_{i\to j}=D_i-\sum_{\substack{e'\ni i,\ e'\ne e\\ e'=ik}}\frac{w_{e'}^2}{q_{k\to i}}
\qquad(e=ij),
\end{equation}
where for parallel edges the messages are indexed by the oriented edge.
Its vertex variance is $V_i=(D_i-\sum_{e\ni i}w_e^2/q_{k\to i})^{-1}$, and
the belief of the edge $e=ij$ is the Gaussian with precision matrix
$\bigl(\begin{smallmatrix}q_{i\to j}&w_e\\ w_e&q_{j\to i}\end{smallmatrix}\bigr)$.
We only consider the precision (covariance) part of the algorithm; the
means do not enter the equations below.

We replace $D$ by $D-z$ for a spectral parameter $z$ and pass to the
variables $V_i$ and one variable $f_e$ per edge.

\begin{definition}\label{def:system}
The \emph{covariance system} of $(G,D,w)$ is the system of $n+\abs E$
equations
\begin{equation}\label{eq:system}
(D_i-z)\,V_i=1+\sum_{e\ni i}f_e\quad(i\in V),\qquad
f_e(1+f_e)=w_e^2\,V_iV_j\quad(e=ij\in E)
\end{equation}
in the unknowns $(z,V,f)\in\C\times\C^V\times\C^E$. A solution is
\emph{admissible} if $z\notin\{D_i\}$ and $V_i\ne0$ for all $i$. A real
solution is \emph{positive} if $V_i>0$ and $f_e>0$ for all $i$ and $e$.
We write $\cZ$ for the set of admissible solutions and $\Xp$ for the set of
positive solutions.
\end{definition}

A positive solution has $z<\min_iD_i$, because the right-hand side of the
vertex equation is positive, and satisfies $f_e=f_+(w_e^2V_iV_j)$ with
\begin{equation}\label{eq:fplus}
f_+(x)=\tfrac12\bigl(\sqrt{1+4x}-1\bigr).
\end{equation}
Positivity of the $V_i$ alone does not imply positivity of the $f_e$.

\begin{lemma}\label{lem:bp-correspondence}
Fix $z\notin\{D_i\}$. The map
\[
V_i=\Bigl(D_i-z-\sum_{e\ni i}\frac{w_e^2}{q_{k\to i}}\Bigr)^{-1},\qquad
f_e=\frac{w_e^2}{q_{j\to i}}\,V_i\quad(e=ij)
\]
is a bijection from the fixed points of \eqref{eq:cavity} for the diagonal
$D-z$ whose cavity precisions, vertex precisions and edge determinants
$q_{i\to j}q_{j\to i}-w_e^2$ are all nonzero, onto the admissible solutions
of \eqref{eq:system} with this~$z$. Its inverse is
$q_{i\to j}=(1+f_e)/V_i$. For real coefficients and real $z$ it restricts
to a bijection between the fixed points all of whose vertex and edge
beliefs are positive definite and the positive solutions.
\end{lemma}

\begin{proof}
Write $t_i=D_i-z$, $m_{k\to i}=w_{e'}^2/q_{k\to i}$ and $R_i=1/V_i=t_i-\sum_km_{k\to i}$.
The cavity equation reads $q_{i\to j}=R_i+m_{j\to i}$. Multiplying by
$m_{i\to j}$ and using $m_{i\to j}q_{i\to j}=w_e^2=m_{j\to i}q_{j\to i}$ gives
$m_{i\to j}R_i=m_{j\to i}R_j$, so $f_e=m_{j\to i}V_i=m_{i\to j}V_j$ does not
depend on the orientation. Then $t_iV_i=V_i(R_i+\sum_km_{k\to i})=1+\sum_{e\ni i}f_e$,
$1+f_e=q_{i\to j}V_i$, and
$f_e(1+f_e)=m_{i\to j}q_{i\to j}V_iV_j=w_e^2V_iV_j$. Conversely, let
$(z,V,f)$ be admissible and put $q_{i\to j}=(1+f_e)/V_i$. If $f_e=-1$ the
edge equation gives $V_iV_j=0$, so $q_{i\to j}\ne0$; the edge equation
gives $w_e^2/q_{k\to i}=f_{e'}/V_i$, so that
$t_i-\sum_{e'\ne e}w_{e'}^2/q_{k\to i}=(1+f_e)/V_i=q_{i\to j}$ by the vertex
equation, the vertex precision is $1/V_i\ne0$, and the edge determinant is
$(1+f_e)/(V_iV_j)\ne0$. The two maps are mutually inverse. In the real case,
$V>0$ and $f>0$ give positive cavity precisions, vertex precisions and edge
determinants; conversely, if the edge belief is positive definite with
covariance $\bigl(\begin{smallmatrix}V_i&C_e\\ C_e&V_j\end{smallmatrix}\bigr)$,
then $f_e=C_e^2/(V_iV_j-C_e^2)>0$.
\end{proof}

By \eqref{eq:G-from-m} and Lemma~\ref{lem:bp-correspondence}, for any
admissible solution the quantity
\begin{equation}\label{eq:m-from-Vf}
m_{i|e}:=\frac{V_i}{1+f_e}=\frac1{q_{i\to j}}=\frac{f_e}{w_e^2V_j}\qquad(e=ij)
\end{equation}
plays the role of the half-tree $m$-function, and
$1+f_e=(1-w_e^2m_{i|e}m_{j|e})^{-1}$.

\begin{proposition}[Physical branch]\label{prop:physical-branch}
Let $G$ be connected with real coefficients.
For $z\in\C\setminus\sigma(\tH)$ set
\[
\Vphys_i(z)=G_i(z),\qquad
f^{\mathrm{phys}}_e(z)=\frac{w_e^2\,m_{i|e}(z)\,m_{j|e}(z)}{1-w_e^2\,m_{i|e}(z)\,m_{j|e}(z)}\quad(e=ij).
\]
Then $(z,\Vphys(z),f^{\mathrm{phys}}(z))$ solves \eqref{eq:system} wherever
the half-tree functions are finite, in particular for $z\notin\R$ and for
real $z\notin[\Sigm,\Sigp]$. It is admissible for $z\notin\R$ and for real
$z\notin[\Sigm,\Sigp]$, and it is a positive solution for real $z<\Sigm$.
\end{proposition}

\begin{proof}
By \eqref{eq:G-from-m}, $w_e^2m_{j|e}G_i=f^{\mathrm{phys}}_e$ and
$f^{\mathrm{phys}}_e(1+f^{\mathrm{phys}}_e)=w_e^2m_{i|e}m_{j|e}(1-w_e^2m_{i|e}m_{j|e})^{-2}=w_e^2G_iG_j$.
Multiplying the first identity of \eqref{eq:m-recursion} by $G_i$ gives the
vertex equation. For $z\notin\R$ the Green function is a nonconstant
Herglotz function, so $G_i(z)\ne0$ and $z\ne D_i$; for real
$z<\Sigm$ the operators $\tH-z$ and its half-tree restrictions are positive
definite, so $G_i>0$, $m_{i|e}>0$ and
$1-w_e^2m_{i|e}m_{j|e}=m_{i|e}/G_i>0$, whence $f^{\mathrm{phys}}_e>0$; the case
$z>\Sigp$ is symmetric, with $G_i<0$.
\end{proof}

We call this solution the \emph{physical branch}. For real $z<\Sigm$ it is
the fixed point to which the precision iteration of Gaussian belief
propagation converges from the standard initialization; its cavity
precisions are the largest, and its vertex variances the smallest, among
all positive solutions (Section~\ref{sec:real}).

\subsection{The spectral curve}
\label{subsec:curve}

\begin{definition}\label{def:curve}
The \emph{spectral curve} $\Xc=\Xc(G;D,w)$ is the union of the irreducible
components of the Zariski closure of $\cZ$ in $\C^{1+n+\abs E}$ on which $z$
is not constant. Its \emph{sheets} are the local branches over a
general~$z$; their number is the degree of $z\colon\Xc\to\C$. For connected
$G$ and real coefficients, the \emph{physical component} $\Cphys$ is the
irreducible component of $\Xc$ that contains the physical branch.
\end{definition}

That $\Xc$ is a curve, that its degree over the $z$-line is the same for
all coefficients, and that $\Cphys$ is well defined are proved in
Section~\ref{sec:degree}.

\begin{definition}\label{def:leafless}
An edge subset $S\subseteq E$ is \emph{leafless} if $d_S(i)\ne1$ for all
$i\in V$. We write $\cL(G)$ for the set of leafless edge subsets and
$\LG=\abs{\cL(G)}$.
\end{definition}

The empty set and, when no vertex has degree one, the full set $E$ are
leafless. The leafless subsets are the sub-coregraphs of
\cite{WatanabeFukumizuZeta} and the generalized loops of
\cite{ChertkovChernyak}.

\subsection{Conventions}
\label{subsec:conventions}

A statement holds for \emph{generic coefficients} if it holds for all
$(D,w)$ in a nonempty Zariski open subset of $\C^V\times(\C^*)^E$; such a
set contains a dense open subset of the real coefficients. We do not
describe the exceptional sets. A point of the normalization of $\Xc$ over a
finite value $z_0$ is \emph{unramified} if $z-z_0$ is a local coordinate
there. A pole of the continued Green function is \emph{antibound} if it
does not lie on the physical sheet $\C\setminus\sigma(\tH)$. We write
$\kappa_0=\abs E-n=\beta-1$.

%% file: sections/degree.tex
\section{Sheets at infinity and the generic degree}
\label{sec:degree}

In this section $G=(V,E)$ is a finite loopless multigraph, not necessarily
connected, and the coefficients are complex, $D\in\C^V$ and
$w\in(\C^*)^E$. We determine all solutions of the covariance system
\eqref{eq:system} for large $\abs z$, for every choice of the coefficients:
they form $2^{\abs E}$ branches indexed by the edge subsets, and a branch is
admissible exactly when its edge subset is leafless
(Lemma~\ref{lem:infinity}). An argument with finite covers of the space of
coefficients then shows that the number of admissible solutions over a fixed
$z$ is $\LG$, with all solutions simple, for generic coefficients
(Theorem~\ref{thm:generic-degree}), and that for arbitrary coefficients the
same count holds over all but finitely many $z$
(Corollary~\ref{cor:degree-all}). This gives the degree statement of
Theorem~\ref{thm:main-sheets} and the facts about the spectral curve that
were postponed in Section~\ref{subsec:curve}. The number $\LG$ itself is a
known graph invariant; its bounds are recalled in
Proposition~\ref{prop:L-bounds}.

\subsection{The reduced system}
\label{subsec:deg-reduced}

For $f\in\C^E$ and $i\in V$ put
\[
g_i(f)=1+\sum_{e\ni i}f_e .
\]
If $z\notin\{D_i\}$, the vertex equations of \eqref{eq:system} are solved by
$V_i=g_i(f)/(D_i-z)$, and the edge equations become a system in $f$ alone.

\begin{definition}\label{deg:def-reduced}
For $c\in\C^E$ the \emph{reduced system} is
\begin{equation}\label{deg:reduced}
f_e(1+f_e)=c_e\,g_i(f)\,g_j(f)\qquad(e=ij\in E)
\end{equation}
in the unknown $f\in\C^E$. A solution is \emph{simple} if the Jacobian
matrix of the $\abs E$ equations with respect to $f$ is invertible at it.
\end{definition}

\begin{lemma}\label{deg:lem-reduction}
Let $z\notin\{D_i\}$ and
\begin{equation}\label{deg:c-of-z}
c_e(z)=\frac{w_e^2}{(D_i-z)(D_j-z)}\qquad(e=ij).
\end{equation}
The map $(V,f)\mapsto f$ is a bijection from the solutions of
\eqref{eq:system} with this value of $z$ onto the solutions of
\eqref{deg:reduced} with $c=c(z)$; its inverse is $V_i=g_i(f)/(D_i-z)$. A
solution of \eqref{eq:system} is admissible if and only if $g_i(f)\ne0$ for
all $i$. The Jacobian determinant of the $n+\abs E$ equations
\eqref{eq:system} with respect to $(V,f)$ equals
$\prod_i(D_i-z)$ times the Jacobian determinant of \eqref{deg:reduced} with
respect to $f$; in particular simple solutions correspond to simple
solutions.
\end{lemma}

\begin{proof}
The first two statements are immediate from $(D_i-z)V_i=g_i(f)$. For the
last one, order the equations and the unknowns of \eqref{eq:system} with the
vertices first. The Jacobian matrix has the block form
$\bigl(\begin{smallmatrix}\diag(D_i-z)&B\\ C&\diag(1+2f_e)\end{smallmatrix}\bigr)$,
and its Schur complement with respect to the invertible first block is the
Jacobian matrix of the system obtained by substituting
$V_i=g_i(f)/(D_i-z)$ into the edge equations, which is \eqref{deg:reduced}.
\end{proof}

By Lemma~\ref{lem:bp-correspondence} the admissible solutions with $z=0$
are the fixed points of the cavity equations \eqref{eq:cavity} with nonzero
cavity precisions, vertex precisions and edge determinants. The two maps in
that lemma are regular near admissible points and are inverse to each other
modulo the respective equations, so that simple admissible solutions
correspond to nondegenerate fixed points of \eqref{eq:cavity}.

\subsection{Weak coupling}
\label{subsec:deg-weak}

For $S\subseteq E$ let $-\ones_S\in\C^E$ be the vector with entries $-1$ on
$S$ and $0$ on $E\setminus S$. For $c=0$ the reduced system is
$f_e(1+f_e)=0$, with the $2^{\abs E}$ solutions $-\ones_S$, and
$g_i(-\ones_S)=1-d_S(i)$.

\begin{lemma}[Weak coupling]\label{deg:lem-weak}
There is $\varepsilon>0$, depending only on $G$, with the following
properties for the polydisc
$P_\varepsilon=\{c\in\C^E:\abs{c_e}<\varepsilon\text{ for all }e\}$.
\begin{enumerate}
\item[(a)] There are holomorphic maps $f^S\colon P_\varepsilon\to\C^E$
$(S\subseteq E)$ with $f^S(0)=-\ones_S$ such that, for every
$c\in P_\varepsilon$, the solutions of \eqref{deg:reduced} are the
$2^{\abs E}$ distinct points $f^S(c)$, and all of them are simple.
\item[(b)] If $S$ is leafless, then $g_i(f^S(c))\ne0$ for all $i\in V$ and
all $c\in P_\varepsilon$.
\item[(c)] If $d_S(i)=1$, then $g_i(f^S(c))=0$ for all $c\in P_\varepsilon$.
\end{enumerate}
\end{lemma}

\begin{proof}
If $E=\emptyset$, the reduced system has one solution in zero variables,
with empty Jacobian determinant $1$, and every $g_i=1$. Assume $E\ne\emptyset$.
(a) Write $F_i=\sum_{e\ni i}f_e$. The part of highest degree of the $e$-th
equation of \eqref{deg:reduced} is $f_e^2-c_eF_iF_j$. Suppose that
$\abs{c_e}d_id_j<1$ for all $e=ij$, and let $f$ be a common zero of these
forms. If $e$ is an edge for which $\abs{f_e}$ is maximal, then
$\abs{F_i}\le d_i\abs{f_e}$ and $\abs{F_j}\le d_j\abs{f_e}$, so
$\abs{f_e}^2\le\abs{c_e}d_id_j\abs{f_e}^2$ and $f=0$. Hence for such $c$ the
closure in $\PP^{\abs E}$ of the solution set of \eqref{deg:reduced} does
not meet the hyperplane at infinity; the solution set is therefore finite,
and by B\'ezout's theorem it consists of $2^{\abs E}$ points counted with
multiplicity. At $c=0$ the Jacobian matrix at $-\ones_S$ is
$\diag(1+2f_e)=\diag(\pm1)$. By the implicit function theorem there are
holomorphic solutions $f^S(c)$ near $c=0$ with $f^S(0)=-\ones_S$, and for
small $\varepsilon$ they are defined on $P_\varepsilon$, pairwise distinct
and simple. Since the total number is $2^{\abs E}$, there are no other
solutions.

(b) We have $g_i(f^S(c))\to1-d_S(i)\ne0$ as $c\to0$; shrink $\varepsilon$.

(c) Let $I=\{i:d_S(i)=1\}$ and let $E_I$ be the set of edges with at least
one end in $I$. Consider the system in the unknowns $f_e$,
$e\in E\setminus E_I$, obtained from the equations of \eqref{deg:reduced}
for $e\in E\setminus E_I$ by substituting the constants
$f_e=-\ones_S(e)$ for $e\in E_I$. It has the same form as
\eqref{deg:reduced}, with the constant $1$ in $g_i$ replaced by an integer,
and at $c=0$ the restriction of $-\ones_S$ is a solution with Jacobian
matrix $\diag(\pm1)$. Hence it has a holomorphic solution $f'(c)$ near
$c=0$ with this value at $0$. Let $\hat f(c)\in\C^E$ be equal to
$-\ones_S$ on $E_I$ and to $f'(c)$ on $E\setminus E_I$. For $i\in I$ all
edges at $i$ belong to $E_I$, so $g_i(\hat f(c))=1-d_S(i)=0$. Hence for
$e\in E_I$ both sides of the $e$-th equation of \eqref{deg:reduced} vanish
at $\hat f(c)$, and the equations for $e\notin E_I$ hold by construction.
Thus $\hat f(c)$ is a solution of \eqref{deg:reduced} that tends to
$-\ones_S$, so $\hat f(c)=f^S(c)$ for small $c$ by (a), and
$g_i(f^S(c))=0$ for $i\in I$ near $c=0$, hence on $P_\varepsilon$ by the
identity theorem.
\end{proof}

\subsection{The branches at infinity}
\label{subsec:deg-infinity}

\begin{lemma}[Branches at infinity]\label{lem:infinity}
Let $D\in\C^V$ and $w\in(\C^*)^E$ be arbitrary. There is $R>0$ such that
the solutions of \eqref{eq:system} with $\abs z>R$ are exactly the points
of $2^{\abs E}$ pairwise disjoint branches
\[
\beta_S=\bigl\{(z,V^S(\tau),f^S(\tau)):\ \tau=1/z,\ 0<\abs\tau<1/R\bigr\}
\qquad(S\subseteq E),
\]
where $V^S$ and $f^S$ are holomorphic for $\abs\tau<1/R$ and
\begin{equation}\label{deg:asymptotics}
f^S_e(\tau)=-\ones_S(e)+O(\tau^2),\qquad
V^S_i(\tau)=-\bigl(1-d_S(i)\bigr)\tau+O(\tau^2).
\end{equation}
Every point of every branch is a simple solution of \eqref{eq:system} for
its value of~$z$. If $S$ is leafless, all points of $\beta_S$ are
admissible. If $S$ is not leafless, no point of $\beta_S$ is admissible:
$V^S_i$ vanishes identically for every $i$ with $d_S(i)=1$. In particular,
for $\abs z>R$ there are exactly $\LG$ admissible solutions over $z$, one on
each branch $\beta_S$ with $S\in\cL(G)$.
\end{lemma}

\begin{proof}
In terms of $\tau=1/z$ the functions \eqref{deg:c-of-z} are
$c_e=w_e^2\tau^2/((1-D_i\tau)(1-D_j\tau))$, holomorphic near $\tau=0$ with a
double zero there. Choose $R>\max_i\abs{D_i}$ such that
$c(z)\in P_\varepsilon$ for $\abs z>R$, with $\varepsilon$ as in
Lemma~\ref{deg:lem-weak}, and put $f^S(\tau)=f^S(c(1/\tau))$ and
$V^S_i(\tau)=-\tau\,g_i(f^S(\tau))/(1-D_i\tau)$. The statement follows from
Lemmas~\ref{deg:lem-reduction} and~\ref{deg:lem-weak}.
\end{proof}

Since each branch is parametrized by $\tau$, the corresponding point over
$z=\infty$ of the normalization of the projective closure of the solution
set is unramified: $\tau$ is a local coordinate there.

\begin{example}\label{deg:ex-D3}
Let $G=D_3$ be the graph with two vertices joined by three parallel edges
$e_1,e_2,e_3$. With $y_e=1+2f_e$ the system \eqref{eq:system} reads
\[
y_e^2=1+4w_e^2V_1V_2,\qquad (D_1-z)V_1=(D_2-z)V_2=\tfrac12(y_1+y_2+y_3-1).
\]
The eight branches at infinity are labeled by the signs of
$y_e\to\pm1$, and $S$ is the set of edges with $y_e\to-1$. The leafless
subsets are $\emptyset$, the three pairs and $E$, so $\LG=5$. On the three
branches with $\abs S=1$ we have $V_1=V_2=0$ identically and
$y$ is constant. On the branch $S=\emptyset$ all $y_e$ tend to $1$.
\end{example}

\subsection{Generic coefficients}
\label{subsec:deg-generic}

\begin{lemma}\label{deg:lem-open}
There is a Zariski open set $U'\subseteq\C^E$ containing the polydisc
$P_\varepsilon$ of Lemma~\ref{deg:lem-weak} such that for every $c\in U'$
the reduced system \eqref{deg:reduced} has exactly $2^{\abs E}$ solutions,
all of them simple, and exactly $\LG$ of them satisfy $g_i(f)\ne0$ for
all~$i$.
\end{lemma}

\begin{proof}
If $E=\emptyset$, take $U'=\C^E$. Otherwise proceed as follows.
Let $U_0\subseteq\C^E$ be the set of $c$ for which the forms
$f_e^2-c_eF_iF_j$ have no common zero in $\PP^{\abs E-1}$. Its complement
is the image of a closed subset of $\C^E\times\PP^{\abs E-1}$ under the
proper projection to $\C^E$, so $U_0$ is Zariski open. Over $U_0$ the
solution set $Z_0=\{(c,f):\eqref{deg:reduced}\}$ is closed in
$U_0\times\PP^{\abs E}$ and does not meet the hyperplane at infinity; the
projection $Z_0\to U_0$ is therefore finite, and every fiber has
$2^{\abs E}$ points counted with multiplicity. The set of $c\in U_0$ over
which some solution is not simple is the image of a closed subset of $Z_0$,
hence closed. Let $U$ be its complement in $U_0$ and $Z_U$ the part of
$Z_0$ over $U$. Then $Z_U\to U$ is a finite unramified covering of degree
$2^{\abs E}$ of the irreducible smooth variety $U$, and $U\supseteq
P_\varepsilon$ by Lemma~\ref{deg:lem-weak}(a).

Let $B=\{(c,f)\in Z_U:\prod_ig_i(f)=0\}$. The connected components
$Z_1,\dots,Z_r$ of $Z_U$ are irreducible, and each of them is a finite
unramified covering of $U$. For each $a$, either $Z_a\subseteq B$, or
$B\cap Z_a$ is a proper closed subset of $Z_a$, and then its image $Y_a$ in
$U$ is a proper closed subset. Let $Y$ be the union of the sets $Y_a$ of the
second kind and $b$ the sum of the degrees of the components $Z_a$ of the
first kind. The fiber of $B$ over $c\in U$ has exactly $b$ points if
$c\notin Y$, and more than $b$ points if $c\in Y$. By
Lemma~\ref{deg:lem-weak}(b),(c) the fiber of $B$ has $2^{\abs E}-\LG$
points for every $c\in P_\varepsilon$. Since $P_\varepsilon$ is not
contained in the proper closed set $Y$, we get $b=2^{\abs E}-\LG$ and then
$P_\varepsilon\cap Y=\emptyset$. The set $U'=U\setminus Y$ has the required
properties.
\end{proof}

\begin{theorem}[Generic degree]\label{thm:generic-degree}
Let $G$ be a finite loopless multigraph and fix $z_0\in\C$. For generic
coefficients $(D,w)\in\C^V\times(\C^*)^E$ the system \eqref{eq:system} with
$z=z_0$ has exactly $2^{\abs E}$ solutions $(V,f)$, all of them simple.
Exactly $\LG$ of them are admissible, and each of the remaining
$2^{\abs E}-\LG$ solutions has $V_i=0$ for some $i$. In particular, for
generic coefficients the cavity equations \eqref{eq:cavity} of Gaussian
belief propagation have exactly $\LG$ complex fixed points with nonzero
cavity precisions, vertex precisions and edge determinants, and all of them
are nondegenerate.
\end{theorem}

\begin{proof}
The map $(D,w)\mapsto c(z_0)$ given by \eqref{deg:c-of-z} is a morphism from
the open set $\{D_i\ne z_0\text{ for all }i\}\times(\C^*)^E$ to $\C^E$, and
its image contains $(\C^*)^E$ (take $D_i=z_0+1$). Hence the preimage of the
set $U'$ of Lemma~\ref{deg:lem-open} is a nonempty Zariski open set of
coefficients. For coefficients in this set the statement follows from
Lemmas~\ref{deg:lem-open} and~\ref{deg:lem-reduction}, and the last
sentence is the case $z_0=0$ together with
Lemma~\ref{lem:bp-correspondence} and the remark after
Lemma~\ref{deg:lem-reduction}.
\end{proof}

The exceptional set depends on $z_0$ only through the translation
$D\mapsto D-z_0$. Since it is a proper Zariski closed set, the theorem
holds in particular for all real coefficients outside a proper real
algebraic subset.

\subsection{The degree for all coefficients}
\label{subsec:deg-all}

Lemma~\ref{deg:lem-open} gives more than the count for generic
coefficients: for arbitrary coefficients it controls all but finitely many
values of $z$.

\begin{corollary}\label{cor:degree-all}
Let $D\in\C^V$ and $w\in(\C^*)^E$ be arbitrary.
\begin{enumerate}
\item[(a)] There is a finite set $\Xi\subset\C$ containing $\{D_i\}$ such
that for every $z\in\C\setminus\Xi$ the system \eqref{eq:system} has exactly
$\LG$ admissible solutions over $z$, and all of them are simple.
\item[(b)] The spectral curve $\Xc$ is a curve, that is, all its
irreducible components have dimension one. Its points over
$z\in\C\setminus\Xi$ are exactly the $\LG$ admissible solutions of
\textup{(a)}; they are smooth points of $\Xc$ at which $z$ is a local
coordinate. In particular $z\colon\Xc\to\C$ has degree $\LG$.
\item[(c)] The normalization of the projective closure of $\Xc$ has exactly
$\LG$ points over $z=\infty$. They are unramified and correspond to the
branches $\beta_S$, $S\in\cL(G)$, of Lemma~\ref{lem:infinity}.
\item[(d)] The branch $\beta_\emptyset$ is contained in exactly one
irreducible component of $\Xc$. If $G$ is connected and the coefficients are
real, $\beta_\emptyset$ is the physical branch for $\abs z>R$, so that the
physical component $\Cphys$ is well defined.
\end{enumerate}
\end{corollary}

For complex coefficients, and for disconnected $G$, we define $\Cphys$ as
the irreducible component of $\Xc$ that contains $\beta_\emptyset$; by (d)
this agrees with Definition~\ref{def:curve}.

\begin{proof}
(a) The map $z\mapsto c(z)$ of \eqref{deg:c-of-z} is a morphism from
$\C\setminus\{D_i\}$ to $\C^E$, so the preimage of $U'$ is Zariski open in
$\C\setminus\{D_i\}$. It contains all $z$ with $\abs z>R$, because then
$c(z)\in P_\varepsilon\subseteq U'$. Hence its complement $\Xi$ in $\C$ is
finite, and (a) follows from Lemmas~\ref{deg:lem-open}
and~\ref{deg:lem-reduction}.

(b) Let $\cZ^\circ$ be the set of admissible solutions over
$\C\setminus\Xi$. It is the complement in $\cZ$ of the admissible solutions
over finitely many values of $z$. At a point $p\in\cZ^\circ$ the Jacobian
matrix of \eqref{eq:system} with respect to $(V,f)$ is invertible, so near
$p$ the solution set of \eqref{eq:system} is the graph of a holomorphic map
of $z$. Since $\cZ$ is open in its Zariski closure, $p$ is a smooth point of
the closure, it lies on exactly one irreducible component, and that
component has dimension one and nonconstant $z$. Conversely, an irreducible
component of the closure of $\cZ$ on which $z$ is not constant meets
$\cZ^\circ$ in a dense subset. Hence $\Xc$ is the closure of $\cZ^\circ$,
and it is a curve. It remains to show that $\Xc$ has no points over
$z_0\in\C\setminus\Xi$ other than those of $\cZ^\circ$. Such a point $p$
would be a solution of \eqref{eq:system} that is not admissible. By
Lemmas~\ref{deg:lem-open} and~\ref{deg:lem-reduction} it is a simple
solution, so the solutions near $p$ form the graph of a holomorphic map of
$z$, and since $p$ is a limit of points of $\cZ^\circ$, the points of this
graph over $z\ne z_0$ near $z_0$ are admissible. The $\LG$ admissible
solutions over $z_0$ are also simple and continue to admissible solutions
over nearby $z$. This gives $\LG+1$ admissible solutions over $z$ near
$z_0$, contradicting~(a).

(c) By Lemma~\ref{lem:infinity} and (b) the part of $\Xc$ over
$\{\abs z>R\}$ is the disjoint union of the punctured discs $\beta_S$,
$S\in\cL(G)$, each mapped isomorphically to $\{0<\abs\tau<1/R\}$ by
$\tau=1/z$.

(d) The branch $\beta_\emptyset$ is connected and consists of smooth points
of $\Xc$, so it lies in one irreducible component. Let $G$ be connected with
real coefficients. By Proposition~\ref{prop:physical-branch} the physical
branch is an admissible solution with $f_e>0$ for all $e$, for real
$z<\Sigm$ and for real $z>\Sigp$. For real $z$ with $\abs z>R$ it lies on
one of the branches $\beta_S$, and for $\abs z$ large the points of
$\beta_S$ with $S\ne\emptyset$ have $f_e$ close to $-1$ for $e\in S$. Hence
for real $z$ of large absolute value the physical branch is
$\beta_\emptyset$. Both are holomorphic on the connected set
$\{\abs z>R'\}$ for some $R'\ge R$, so they agree there.
\end{proof}

\begin{remark}\label{deg:rem-constant-components}
The Zariski closure of $\cZ$ may have further irreducible components over
the finitely many points of $\Xi\setminus\{D_i\}$, on which $z$ is constant;
they may have positive dimension. By definition they are not part of $\Xc$.
Corollary~\ref{cor:degree-all} does not say that $\Xc$ is irreducible.
For generic coefficients it is (Theorem~\ref{thm:irreducible}), but for
special coefficients $\Cphys$ can have fewer than $\LG$ sheets; for
example, for a regular graph with uniform coefficients the physical branch
lies on a curve with two sheets
(Example~\ref{pole:ex-regular}).
\end{remark}

\begin{remark}\label{deg:rem-real}
For real coefficients and real $z$ with $\abs z>R$, all $2^{\abs E}$
solutions are real, since complex conjugation permutes the branches
$\beta_S$ and preserves their limits. A branch $\beta_S$ with
$S\ne\emptyset$ has $f_e<0$ for $e\in S$. Hence for $z<-R$ the only positive
solution is the physical one. For bounded $z$ the positive solutions are
studied in Section~\ref{sec:real}.
\end{remark}

\begin{proposition}\label{prop:nondegenerate-bound}
Let $D\in\C^V$, $w\in(\C^*)^E$ and $z_0\in\C\setminus\{D_i\}$ be arbitrary.
The number of simple admissible solutions of \eqref{eq:system} with
$z=z_0$ is at most $\LG$. In particular, for real coefficients and every
real $z_0$ the number of positive solutions over $z_0$ at which the
Jacobian matrix of
\begin{equation}\label{deg:vertex-map}
V\longmapsto\Bigl((D_i-z_0)V_i-1-\sum_{e\ni i,\ e=ij}f_+(w_e^2V_iV_j)\Bigr)_{i\in V}
\end{equation}
is invertible is at most $\LG$.
\end{proposition}

\begin{proof}
A simple solution over $z_0$ continues holomorphically to solutions over
all $z$ near $z_0$ by the implicit function theorem. Distinct solutions
have distinct continuations, and admissibility is an open condition. If
there were more than $\LG$ simple admissible solutions over $z_0$, there
would be more than $\LG$ admissible solutions over all nearby $z$,
contradicting Corollary~\ref{cor:degree-all}(a). For the second statement,
recall that a positive solution has $z_0<\min_iD_i$ and
$f_e=f_+(w_e^2V_iV_j)$. In the Jacobian matrix of \eqref{eq:system} the
block $\diag(1+2f_e)$ of the edge equations is invertible because
$f_e>0$, and its Schur complement is the Jacobian matrix of
\eqref{deg:vertex-map}. Hence a positive solution is simple if and only if
the latter is invertible.
\end{proof}

By the conservation of the local intersection number under perturbation,
the same argument bounds by $\LG$ the number of isolated admissible
solutions over $z_0$ counted with multiplicity. The bound of
Proposition~\ref{prop:nondegenerate-bound} for positive solutions is far
from sharp in general; the positive solutions are counted in
Section~\ref{sec:real}.

\subsection{The number of leafless subsets}
\label{subsec:deg-L}

The leafless edge subsets are the edge sets of the sub-coregraphs of
\cite{WatanabeFukumizuZeta} and the generalized loops of the loop series
\cite{ChertkovChernyak}. For the determinant of a Gaussian model, a loop
series around an arbitrary fixed point of Gaussian belief propagation, with
the same index set, is given in \cite{ChernyakChertkov}; fixed points are
not counted there. Theorem~\ref{thm:generic-degree} says that the index set
of these expansions also counts the fixed points themselves. What is new
here is this identification and the comparison with the binary case in
Proposition~\ref{prop:L-bounds}(c); the bounds in (b) are known.

Let $\psi(1)=0$ and $\psi(k)=1$ for $k\ne1$, so that
\begin{equation}\label{deg:L-sum}
\LG=\sum_{S\subseteq E}\ \prod_{i\in V}\psi\bigl(d_S(i)\bigr).
\end{equation}
The generic number of fixed points of belief propagation for binary
variables on $G$ is given by the same sum with another vertex weight
\cite{WatanabeCounting}:
\begin{equation}\label{deg:Q-sum}
Q(G)=\sum_{S\subseteq E}\ \prod_{i\in V}\varphi\bigl(d_S(i)\bigr),\qquad
\varphi(0)=1,\quad\varphi(k)=k-1\ \ (k\ge1).
\end{equation}

\begin{proposition}\label{prop:L-bounds}
Let $G$ be a finite loopless multigraph.
\begin{enumerate}
\item[(a)] $\LG$ is multiplicative over the connected components of $G$,
and it does not change when a vertex of degree one is removed together with
its edge. Hence $\LG=1$ for forests, and $\LG$ depends only on the $2$-core
of~$G$.
\item[(b)] \textup{(\cite[Theorem~4]{WatanabeFukumizuZeta})} If $G$ is
connected with cycle rank $\beta$, then
\begin{equation}\label{deg:WF-bounds}
2^\beta\ \le\ \LG\ \le\
\Bigl(\frac{5+\sqrt5}2\Bigr)^{\beta-1}+\Bigl(\frac{5-\sqrt5}2\Bigr)^{\beta-1}.
\end{equation}
Both bounds are equalities for $\beta\le1$. For $\beta\ge2$, equality
holds on the left if and only if the $2$-core of $G$ is a subdivision of a
bouquet of cycles, and on the right if and only if the $2$-core of $G$ is a
subdivision of a $3$-regular graph, where multiple edges and self-loops are
allowed in the $3$-regular graph.
\item[(c)] $\LG\le Q(G)$, with equality if and only if every connected
component of $G$ has cycle rank at most one.
\end{enumerate}
\end{proposition}

\begin{proof}
(a) Leafless subsets of a disjoint union are the unions of leafless subsets
of the parts, and an edge at a vertex of degree one belongs to no leafless
subset.

(b) We refer to \cite{WatanabeFukumizuZeta} for the upper bound and the
cases of equality. The lower bound is elementary: the $2^\beta$ edge
subsets in which every vertex has even degree are leafless.

(c) We have $\psi\le\varphi$, with equality exactly at the arguments
$0,1,2$, so each term of \eqref{deg:L-sum} is at most the corresponding
term of \eqref{deg:Q-sum}, and the terms of non-leafless $S$ vanish in
both sums. Equality holds if and only if no leafless subset has a vertex of
degree at least three. If every component has cycle rank at most one, the
leafless subsets are unions of disjoint cycles. If a component has cycle
rank at least two, the edge set of its $2$-core is leafless and connected
and is not a cycle, so it has a vertex of degree at least three.
\end{proof}

The right-hand side of \eqref{deg:WF-bounds} is $1,2,5,15,50,175,625$ for
$\beta=0,\dots,6$. Together with Theorem~\ref{thm:generic-degree} and
Proposition~\ref{prop:nondegenerate-bound}, it bounds the number of
nondegenerate fixed points of Gaussian belief propagation, complex or
positive, in terms of the cycle rank alone.

Table~\ref{deg:table} lists $\LG$ and $Q(G)$ for some small graphs; the
values are obtained by enumeration of the edge subsets. The theta graph,
the dumbbell, $D_3$, $K_4$, the prism, $K_{3,3}$, the cube and the Petersen
graph are subdivisions of $3$-regular graphs and attain the upper bound in
\eqref{deg:WF-bounds}. For the complete graph, the leafless subsets of
$E(K_n)$ are by definition the labeled graphs on $n$ vertices without
vertices of degree one, so
\[
L(K_n)=1,\ 1,\ 2,\ 15,\ 314,\ 13757,\ 1142968\qquad(n=1,\dots,7)
\]
is the sequence A059167 of \cite{OEIS}; see also \cite{HararyPalmer}.

\begin{table}[ht]
\centering
\begin{tabular}{lrrrr}
\toprule
$G$ & $\beta$ & $2^\beta$ & $\LG$ & $Q(G)$\\
\midrule
tree & $0$ & $1$ & $1$ & $1$\\
cycle, or any connected graph with $\beta=1$ & $1$ & $2$ & $2$ & $2$\\
bowtie (two triangles with a common vertex) & $2$ & $4$ & $4$ & $6$\\
$D_3$, theta graph, dumbbell, diamond, $K_{2,3}$ & $2$ & $4$ & $5$ & $8$\\
$K_4$ & $3$ & $8$ & $15$ & $48$\\
$K_4$ with one edge doubled & $4$ & $16$ & $38$ & $204$\\
$3\times3$ grid & $4$ & $16$ & $43$ & $252$\\
$K_{3,3}$, prism & $4$ & $16$ & $50$ & $320$\\
wheel with five spokes & $5$ & $32$ & $118$ & $1180$\\
cube & $5$ & $32$ & $175$ & $2176$\\
$K_5$ & $6$ & $64$ & $314$ & $4976$\\
Petersen graph & $6$ & $64$ & $625$ & $14848$\\
octahedron & $7$ & $128$ & $982$ & $27304$\\
\bottomrule
\end{tabular}
\caption{The number $\LG$ of leafless edge subsets and the generic number
$Q(G)$ of fixed points of binary belief propagation.}
\label{deg:table}
\end{table}

\begin{remark}[Comparison with binary variables]\label{deg:rem-binary}
The sums \eqref{deg:L-sum} and \eqref{deg:Q-sum} run over the same index
set, the leafless edge subsets, and differ only in the vertex weight. The
proofs have the same form as well: in \cite{WatanabeCounting} the fixed
points of binary belief propagation are counted by following them to a
limit of the couplings in which they are indexed by leafless subsets $S$,
with $\prod_i(d_S(i)-1)$ fixed points for each $S$; here each leafless $S$
carries exactly one solution, the branch $\beta_S$.
\end{remark}

%% file: sections/local-poles.tex
\section{Poles at signed eigenvalues and the Floquet identity}
\label{sec:local-poles}

This section contains the local analysis of the covariance system
\eqref{eq:system} at the points where some of the $V_i$ become infinite. The
statements are local and algebraic: $G=(V,E)$ is a finite loopless
multigraph, not necessarily connected, and the coefficients $D\in\C^V$,
$w\in(\C^*)^E$ are complex unless real coefficients are specified. We show
that at a simple eigenvalue of a signed Jacobi matrix $D+A^s$ with a
nowhere vanishing eigenvector there is a unique branch of admissible
solutions on which all $V_i$ have simple poles, with residues
$\kappa_0u_i^2$ (Theorem~\ref{thm:local-pole}); that the same construction
on a vertex subset $U$, coupled to the rest of the graph through a boundary
system, gives poles supported on $U$ (Theorem~\ref{thm:partial-pole}); that
the total residue of $(\sum_iV_i)\,dz$ over the finite part of the spectral
curve is $\sum_{S\in\cL(G)}(2\abs S-n)$
(Proposition~\ref{prop:residue-sum}); and that this differential is, on all
sheets, the negative logarithmic derivative of a rational function which
on the physical branch is the Floquet function of \cite{BBGSS}
(Proposition~\ref{prop:floquet}). These are the local ingredients of
Theorems~\ref{thm:main-perron} and~\ref{thm:main-poles}. Whether a given
pole is reached from the physical branch is a global question, treated in
Sections~\ref{sec:perron} and~\ref{sec:irreducibility}, and the fact that
there are no other poles with nonzero residue is
Theorem~\ref{thm:classification}.

\subsection{Poles at eigenvalues of signed Jacobi matrices}
\label{subsec:pole-full}

Let $s\in\{\pm1\}^E$ be a signing and let $\lambda$ be an eigenvalue of
the complex symmetric matrix $M=D+A^s$ of algebraic multiplicity one, with
eigenvector $u$. Then $u^{\T}u\ne0$: the image of $M-\lambda$ is
$\{x:u^{\T}x=0\}$, so if $u^{\T}u=0$ there is $v$ with $(M-\lambda)v=u$ and
the algebraic multiplicity is at least two. We may therefore normalize
$u^{\T}u=1$; for real coefficients this is $\norm u=1$. We write $u^2$ for
the vector $(u_i^2)_{i\in V}$.

\begin{theorem}[Poles at signed eigenvalues]\label{thm:local-pole}
Let $G$ be a finite loopless multigraph with $\abs E\ne n$, and let
$D\in\C^V$, $w\in(\C^*)^E$. Let $s$ be a signing, $\lambda$ an eigenvalue of
$D+A^s$ of algebraic multiplicity one, and $u$ an eigenvector with
$u_i\ne0$ for all $i$ and $u^{\T}u=1$. Put $\kappa_0=\abs E-n$. Then there
are $r>0$ and holomorphic maps $W\colon\{\abs\zeta<r\}\to\C^V$ and
$y\colon\{\abs\zeta<r\}\to\C^E$ with
\[
W_i(0)=\kappa_0u_i^2,\qquad y_e(0)=-2\kappa_0\,s_ew_eu_iu_j\quad(e=ij),
\]
such that for $0<\abs{z-\lambda}<r$
\begin{equation}\label{pole:branch}
V_i(z)=\frac{W_i(z-\lambda)}{z-\lambda},\qquad
f_e(z)=\frac12\Bigl(\frac{y_e(z-\lambda)}{z-\lambda}-1\Bigr)
\end{equation}
is an admissible solution of \eqref{eq:system}. Thus
\begin{equation}\label{pole:expansion}
V_i(z)=\frac{\kappa_0u_i^2}{z-\lambda}+O(1),\qquad
f_e(z)=-\frac{\kappa_0\,s_ew_eu_iu_j}{z-\lambda}+O(1)\quad(e=ij).
\end{equation}
This branch is unique in the following sense: there is a neighborhood
$\mathcal N$ of $(W(0),y(0))$ in $\C^V\times\C^E$ such that every solution
$(z,V,f)$ of \eqref{eq:system} with $0<\abs{z-\lambda}<r$ and
$\bigl((z-\lambda)V,(z-\lambda)(1+2f)\bigr)\in\mathcal N$ is given by
\eqref{pole:branch}. In particular, every sequence of solutions with
$z\to\lambda$, $z\ne\lambda$, along which $(z-\lambda)(V,1+2f)$ converges
to $(W(0),y(0))$ eventually lies on the branch.
\end{theorem}

The branch is contained in the spectral curve $\Xc$ and defines one point
of its normalization over $z=\lambda$. This point is unramified, each $V_i$
has a simple pole there, and
$\sum_i\Res_{z=\lambda}V_i=\kappa_0u^{\T}u=\kappa_0$. The case of interest
is $\abs E>n$, that is, $\beta\ge2$ for connected $G$; the eigenvector of a
simple eigenvalue can be nowhere zero only if $G$ is connected. The data
$u_i^2$ and $s_ew_eu_iu_j$ depend only on the switching class of $s$, so
the theorem attaches a branch to each pair of a switching class and an
eigenvalue: at most $n2^\beta$ branches in all.

\begin{proof}
Put $\zeta=z-\lambda$, $W=\zeta V$ and $y_e=\zeta(1+2f_e)$. Since
$f(1+f)=x$ is equivalent to $(1+2f)^2=1+4x$, for $\zeta\ne0$ the system
\eqref{eq:system} is equivalent to the polynomial system
\begin{equation}\label{pole:rescaled}
\begin{aligned}
p_i&:=(D_i-\lambda-\zeta)W_i-\zeta\Bigl(1-\frac{d_i}2\Bigr)-\frac12\sum_{e\ni i}y_e=0
&&(i\in V),\\
q_e&:=y_e^2-\zeta^2-4w_e^2W_iW_j=0&&(e=ij\in E).
\end{aligned}
\end{equation}
Let $y^0_e=-2s_ew_eu_iu_j$. For $\zeta=0$ and every $\kappa\in\C$ the point
$(W,y)=\kappa(u^2,y^0)$ solves \eqref{pole:rescaled}: the edge equations
hold identically, and
$p_i=\kappa u_i\bigl((D_i-\lambda)u_i+\sum_{e\ni i,\,e=ij}s_ew_eu_j\bigr)=0$.

\emph{Elimination of $y$.} Fix $\kappa_*\ne0$. At
$(\kappa_*u^2,\kappa_*y^0,0)$ we have
$\partial q_e/\partial y_e=2\kappa_*y^0_e\ne0$, so near this point the edge
equations are equivalent to $y_e=Y_e(W,\zeta)$ with $Y_e$ holomorphic,
$Y_e(\kappa u^2,0)=\kappa y^0_e$, and at the points $(\kappa u^2,0)$
\[
\frac{\partial Y_e}{\partial W_i}=\frac{2w_e^2W_j}{y_e}=-s_ew_e\frac{u_j}{u_i},
\qquad\frac{\partial Y_e}{\partial\zeta}=\frac\zeta{y_e}=0\qquad(e=ij).
\]
Let $\tilde p(W,\zeta)=p(W,Y(W,\zeta),\zeta)$. Using
$\sum_{e\ni i,\,e=ij}s_ew_eu_j=-(D_i-\lambda)u_i$ we find that the Jacobian
matrix of $\tilde p$ with respect to $W$ at $(\kappa u^2,0)$ is
\begin{equation}\label{pole:J}
J=\tfrac12\diag(u)\,(D+A^s-\lambda)\,\diag(u)^{-1},
\end{equation}
independently of $\kappa$. Its kernel is spanned by $u^2=\diag(u)u$, its
image is the hyperplane $H=\{x\in\C^V:\ones^{\T}x=0\}$, and
$\ones^{\T}u^2=u^{\T}u=1$, so $\C^V=\C u^2\oplus H$ and $J$ maps $H$
isomorphically onto $H$.

\emph{Reduction to one equation.} Write $W=\kappa u^2+\eta$ with
$\eta\in H$, and let $\pi(x)=x-(\ones^{\T}x)u^2$ be the projection onto $H$
along $u^2$. The equation $\pi\tilde p(\kappa u^2+\eta,\zeta)=0$ has
derivative $\pi J|_H=J|_H$ with respect to $\eta$ at $(\kappa_*,0,0)$, so by
the implicit function theorem it is equivalent, near this point, to
$\eta=\eta(\kappa,\zeta)$ with $\eta$ holomorphic; and
$\eta(\kappa,0)=0$, because the points $\kappa(u^2,y^0)$ are solutions. The
remaining equation is $h(\kappa,\zeta)=0$ with
\[
h(\kappa,\zeta)=\ones^{\T}\tilde p\bigl(\kappa u^2+\eta(\kappa,\zeta),\zeta\bigr).
\]
Since $h(\kappa,0)=0$, we have $h=\zeta h_1$ with $h_1$ holomorphic, and
because $\ones^{\T}J=0$ and $\partial Y/\partial\zeta=0$ at $\zeta=0$,
\[
h_1(\kappa,0)=\ones^{\T}\frac{\partial p}{\partial\zeta}(\kappa u^2,\kappa y^0,0)
=\sum_{i\in V}\Bigl(-\kappa u_i^2-1+\frac{d_i}2\Bigr)=\abs E-n-\kappa .
\]

\emph{Conclusion.} Take $\kappa_*=\kappa_0=\abs E-n$, which is not zero by
hypothesis. Then $h_1(\kappa_0,0)=0$ and $\partial_\kappa h_1(\kappa_0,0)=-1$,
so $h_1=0$ is equivalent near $(\kappa_0,0)$ to $\kappa=\kappa(\zeta)$ with
$\kappa$ holomorphic and $\kappa(0)=\kappa_0$. Hence, near
$(\kappa_0u^2,\kappa_0y^0,0)$, the solutions of \eqref{pole:rescaled}
consist of the line $\{\kappa(u^2,y^0)\}$ in $\{\zeta=0\}$ and of the
holomorphic curve
\[
W(\zeta)=\kappa(\zeta)u^2+\eta(\kappa(\zeta),\zeta),\qquad
y(\zeta)=Y(W(\zeta),\zeta).
\]
For $\zeta\ne0$ this curve gives the solution \eqref{pole:branch} of
\eqref{eq:system}, and the uniqueness statement follows. For small
$\zeta\ne0$ we have $V_i=W_i(\zeta)/\zeta\ne0$, since
$W_i(0)=\kappa_0u_i^2\ne0$, and $z=\lambda+\zeta\notin\{D_i\}$; so the
solution is admissible.
\end{proof}

\begin{remark}\label{pole:rem-hypotheses}
(a) The hypothesis $\abs E\ne n$ is used only to have $\kappa_0\ne0$. If
$\abs E=n$, for example for a cycle, then $\kappa_0=0$, the candidate limit
point is $(W,y)=0$, where the elimination of $y$ breaks down, and the
theorem gives no branch. For the cycle $C_n$ with $D=0$ and $w=1$ the
solutions with all $V_i$ equal and all $f_e$ equal satisfy $-zV=1+2f$ and
$(1+2f)^2=1+4V^2$, so that $V=\pm(z^2-4)^{-1/2}$: at the Perron values
$\pm2$, which here are the edges of the spectrum of $\tH$, the function $V$
has a square root singularity and no pole.

(b) For a tree, $\kappa_0=-1$ and all signings are switching equivalent.
The theorem then gives $V_i=-u_i^2/(z-\lambda)+O(1)$. This agrees with the
residues of $((D+A-z)^{-1})_{ii}$, which for real coefficients is the
physical branch, the universal cover of a tree being the tree itself, and
is the only sheet because $\LG=1$.

(c) The value $\lambda$ may coincide with some $D_i$; this does not affect
admissibility, since the points of the branch have $z\ne\lambda$.
\end{remark}

For connected $G$ the hypotheses on $\lambda$ and $u$ hold for all signings
and all eigenvalues when the coefficients are generic. We state this for
induced subgraphs, as needed in Section~\ref{subsec:pole-partial}. For
$U\subseteq V$ let $E(U)$ be the set of edges with both ends in $U$,
$G[U]=(U,E(U))$, and for a signing $s$ of $E(U)$ let $(D+A^s)[U]$ be the
signed Jacobi matrix of $G[U]$ with the coefficients restricted to $U$ and
$E(U)$.

\begin{lemma}\label{pole:lem-generic}
For generic coefficients the following holds for every nonempty
$U\subseteq V$ with $G[U]$ connected and every signing $s$ of $E(U)$: all
eigenvalues of $(D+A^s)[U]$ are simple and all its eigenvectors have no zero
entry.
\end{lemma}

\begin{proof}
There are finitely many pairs $(U,s)$, and the substitution
$w_e\mapsto s_ew_e$ preserves genericity, so it suffices to treat $U=V$,
$G$ connected and $s\equiv1$; let $M=D+A$ and $n\ge2$. If $\lambda$ is a
simple eigenvalue of a symmetric matrix $M$ with eigenvector $u$, the
adjugate of $M-\lambda$ is a nonzero symmetric matrix of rank one whose
columns are multiples of $u$, hence equal to $a\,uu^{\T}$ with $a\ne0$; its
$i$-th diagonal entry is $\det(M[V\setminus\{i\}]-\lambda)$. So $u_i\ne0$ if
and only if $\lambda$ is not an eigenvalue of $M[V\setminus\{i\}]$. The
statement therefore holds where the discriminant of the characteristic
polynomial $\chi_M$ and the resultants of $\chi_M$ and
$\chi_{M[V\setminus\{i\}]}$, $i\in V$, do not vanish. These are polynomials
in $(D,w)\in\C^V\times\C^E$, and it is enough to find one point, possibly
with some $w_e=0$, where none of them vanishes. Choose a spanning tree $T$
of $G$, real distinct $D_i$, and $w_e=t$ for $e\in T$, $w_e=0$ otherwise.
For small $t$ the eigenvalues $\lambda_k(t)\to D_k$ are simple. Let $u(t)$
be the eigenvector of $\lambda_k(t)$ with $u_k=1$; it is holomorphic in $t$
with $u(0)$ the $k$-th unit vector, and
$(\lambda_k-D_j)u_j=t\sum_lu_l$ for $j\ne k$, the sum being over the
neighbors $l$ of $j$ in $T$. By induction on $m$, $u_j=O(t^m)$ for all $j$
at distance at least $m$ from $k$ in $T$. Since a vertex at distance $m\ge1$
has exactly one neighbor at distance $m-1$ and its other neighbors are at
distance $m+1$, a second induction gives
$u_j=t^m\prod_l(D_k-D_l)^{-1}+O(t^{m+1})$, the product over the vertices
$l\ne k$ of the path from $k$ to $j$. Hence $u(t)$ has no zero entry for
small $t\ne0$.
\end{proof}

Consequently, for connected $G$ with $\beta\ge2$ and generic coefficients,
Theorem~\ref{thm:local-pole} applies to each of the $n$ eigenvalues of each
of the $2^\beta$ switching classes.

\begin{example}[The graph $D_3$]\label{pole:ex-D3}
Let $G=D_3$ as in Example~\ref{deg:ex-D3}, with weights $w_1,w_2,w_3$. Here
$\kappa_0=1$, there are four switching classes, represented by the signings
$(+,+,+)$, $(-,+,+)$, $(+,-,+)$, $(+,+,-)$, and $D+A^s$ is the
$2\times2$ matrix with diagonal $(D_1,D_2)$ and off-diagonal entry
$a_s=s_1w_1+s_2w_2+s_3w_3$. For real coefficients with all $a_s\ne0$ the
hypotheses of Theorem~\ref{thm:local-pole} hold, and we obtain eight
branches with simple poles, at
\[
z=\frac{D_1+D_2}2\pm\sqrt{\Bigl(\frac{D_1-D_2}2\Bigr)^2+a_s^2}\,,
\]
with residues $\Res V_i=u_i^2$ and $\Res V_1+\Res V_2=1$. The Perron values
$\sigma_\pm$ are the two poles with $a_s=\abs{w_1}+\abs{w_2}+\abs{w_3}$. If
$w_1=w_2=w_3$, the three nontrivial classes have the same eigenvalues.
\end{example}

\begin{example}[Regular graphs with uniform coefficients]\label{pole:ex-regular}
Let $G$ be connected and $d$-regular with $d\ge3$, $D=0$ and $w\equiv1$.
The solutions of \eqref{eq:system} with all $V_i$ equal and all $f_e$ equal
satisfy $-zV=1+df$ and $f(1+f)=V^2$, which gives the two-sheeted curve
\[
V=\frac{(d-2)z\mp d\sqrt{z^2-4(d-1)}}{2(z^2-d^2)} .
\]
With the square root that behaves like $z$ at infinity, the upper sign is
the physical branch, the Green function of the $d$-regular tree, and it is
regular at $z=\pm d$. The lower sign has simple poles at
$z=\pm d=\sigma_\pm$ with residue $(d-2)/2$. This agrees with
Theorem~\ref{thm:local-pole} for the Perron eigenvalue of $D\pm\abs A$:
$\bar u_i^2=1/n$ and $\kappa_0/n=(nd/2-n)/n$. These are the second-sheet
poles of \cite[Example~3.1]{CSZ}.
\end{example}

\subsection{Poles supported on a vertex subset}
\label{subsec:pole-partial}

The construction of Theorem~\ref{thm:local-pole} can be carried out on a
subset of the vertices. Let $U\subseteq V$ be nonempty. We use the
following notation.
\begin{itemize}
\item $N$ is the set of vertices outside $U$ that are adjacent to $U$, and
$T=V\setminus(U\cup N)$.
\item $E(U)$ and $E(T)$ are the sets of edges with both ends in $U$,
respectively in $T$; $\partial U$ is the set of edges between $U$ and $N$;
and $E_N$ is the set of edges with at least one end in $N$ and no end in
$U$. These four sets partition $E$.
\item $d^U_i$ is the number of edges of $E(U)$ at $i\in U$.
\end{itemize}
The \emph{limit data} consist of: a signing $s$ of $E(U)$; an eigenvalue
$\lambda$ of $(D+A^s)[U]$ of algebraic multiplicity one with an eigenvector
$u\in\C^U$ without zero entries, normalized by $u^{\T}u=1$; a
\emph{freezing} $\phi\in\{0,-1\}^{E_N}$; numbers $\kappa\in\C$,
$f^*_e$ $(e\in\partial U)$ and $\nu_x$ $(x\in N)$; and numbers $V^*_y$
$(y\in T)$, $f^*_e$ $(e\in E(T))$. They are required to satisfy the
\emph{boundary system}
\begin{equation}\label{pole:boundary}
\begin{aligned}
\kappa&=\abs{E(U)}-\abs U-\sum_{e\in\partial U}f^*_e,\\
f^*_e(1+f^*_e)&=w_e^2\,\kappa\,u_i^2\,\nu_x&&(e=ix\in\partial U,\ i\in U,\ x\in N),\\
0&=1+\sum_{e\in\partial U,\ e\ni x}f^*_e+\sum_{e\in E_N,\ e\ni x}\phi_e&&(x\in N),
\end{aligned}
\end{equation}
in the unknowns $(\kappa,(f^*_e)_{e\in\partial U},(\nu_x)_{x\in N})$, and
the \emph{remaining system}
\begin{equation}\label{pole:remaining}
\begin{aligned}
(D_y-\lambda)V^*_y&=1+\sum_{e\in E(T),\ e\ni y}f^*_e+\sum_{e\in E_N,\ e\ni y}\phi_e&&(y\in T),\\
f^*_e(1+f^*_e)&=w_e^2\,V^*_yV^*_{y'}&&(e=yy'\in E(T)),
\end{aligned}
\end{equation}
in the unknowns $((V^*_y)_{y\in T},(f^*_e)_{e\in E(T)})$. The remaining
system is the covariance system of $G[T]$ at $z=\lambda$ with the constant
$1$ of each vertex equation changed by the frozen edges. Both systems are
square, and we call a solution \emph{nondegenerate} if the Jacobian matrix
with respect to the unknowns is invertible at it.

\begin{theorem}[Poles supported on a vertex subset]\label{thm:partial-pole}
Let $G$ be a finite loopless multigraph, $D\in\C^V$, $w\in(\C^*)^E$, and
let $U\subseteq V$ be nonempty. Suppose that limit data as above are given
such that
\begin{enumerate}
\item[(i)] $\lambda\notin\{D_x:x\notin U\}$;
\item[(ii)] $(\kappa,f^*|_{\partial U},\nu)$ is a nondegenerate solution of
the boundary system \eqref{pole:boundary} with $\kappa\ne0$ and
$\nu_x\ne0$ for all $x\in N$;
\item[(iii)] $(V^*,f^*|_{E(T)})$ is a nondegenerate solution of the
remaining system \eqref{pole:remaining} with $V^*_y\ne0$ for all $y\in T$.
\end{enumerate}
Then there are $r>0$ and a branch of admissible solutions of
\eqref{eq:system}, defined and holomorphic for $0<\abs{z-\lambda}<r$, with
\begin{equation}\label{pole:partial-expansion}
\begin{aligned}
V_i&=\frac{\kappa u_i^2}{z-\lambda}+O(1)\quad(i\in U),&
f_e&=-\frac{\kappa\,s_ew_eu_iu_j}{z-\lambda}+O(1)\quad(e=ij\in E(U)),\\
V_x&=\nu_x(z-\lambda)+O\bigl((z-\lambda)^2\bigr)\quad(x\in N),&
f_e&=f^*_e+O(z-\lambda)\quad(e\in\partial U\cup E(T)),\\
V_y&=V^*_y+O(z-\lambda)\quad(y\in T),&
f_e&=\phi_e+O(z-\lambda)\quad(e\in E_N).
\end{aligned}
\end{equation}
The functions $(z-\lambda)V_i$ $(i\in U)$, $(z-\lambda)(1+2f_e)$
$(e\in E(U))$, $V_x/(z-\lambda)$ $(x\in N)$ and the remaining $V$ and $f$
extend holomorphically to $z=\lambda$, and the branch is the only one with
this property and these limits: every solution of \eqref{eq:system} with
$0<\abs{z-\lambda}<r$ for which these rescaled quantities are sufficiently
close to their limits lies on it. The corresponding point of the
normalization of $\Xc$ over $z=\lambda$ is unramified, and
\begin{equation}\label{pole:kappa-integer}
\Res_{z=\lambda}\sum_{i\in V}V_i=\kappa
=\abs{E(U)}-\abs U+\abs N-\#\{e\in E_N:\phi_e=-1\}-\#\{e\in E_N\cap E(N):\phi_e=-1\},
\end{equation}
where $E(N)$ is the set of edges with both ends in $N$. In particular
$\kappa$ is an integer determined by $U$ and the freezing.
\end{theorem}

For $U=V$ the sets $N$, $T$, $\partial U$ and $E_N$ are empty, the boundary
system is $\kappa=\abs E-n$, and the statement is
Theorem~\ref{thm:local-pole}.

\begin{proof}
Let $\zeta=z-\lambda$. We use the unknowns $W_i=\zeta V_i$ $(i\in U)$,
$y_e=\zeta(1+2f_e)$ $(e\in E(U))$, $\nu_x=V_x/\zeta$ $(x\in N)$, and the
unscaled $f_e$ $(e\notin E(U))$ and $V_y$ $(y\in T)$. For $\zeta\ne0$ the
system \eqref{eq:system} is equivalent to
\begin{align}
&(D_i-\lambda-\zeta)W_i-\zeta\Bigl(1-\frac{d^U_i}2+\sum_{e\in\partial U,\,e\ni i}f_e\Bigr)-\frac12\sum_{e\in E(U),\,e\ni i}y_e=0&&(i\in U),\label{pole:pU}\\
&y_e^2-\zeta^2-4w_e^2W_iW_j=0&&(e=ij\in E(U)),\label{pole:qU}\\
&f_e(1+f_e)-w_e^2W_i\nu_x=0&&(e=ix\in\partial U),\label{pole:bd-edge}\\
&1+\sum_{e\ni x}f_e-(D_x-\lambda-\zeta)\,\zeta\,\nu_x=0&&(x\in N),\label{pole:bd-vertex}\\
&f_e(1+f_e)-w_e^2\,\zeta^{a_e}\,\nu_xP_e=0&&(e\in E_N,\ e\ni x),\label{pole:frozen}
\end{align}
together with the equations of \eqref{eq:system} for the vertices of $T$
and the edges of $E(T)$, with $D_y-\lambda-\zeta$ in place of
$D_y-\lambda$. In \eqref{pole:frozen}, $x\in N$ is an end of $e$, and
$a_e=2$, $P_e=\nu_{x'}$ if the other end $x'$ is in $N$, while $a_e=1$,
$P_e=V_y$ if the other end $y$ is in $T$. All these equations are
polynomial, and the limit data, with $W=\kappa u^2$ and
$y_e=-2\kappa s_ew_eu_iu_j$, solve them at $\zeta=0$; note that the first
equation of \eqref{pole:boundary} is not among them.

The equations \eqref{pole:pU}, \eqref{pole:qU} differ from
\eqref{pole:rescaled} for the graph $G[U]$ only by the term
$\zeta\sum_{e\in\partial U,\,e\ni i}f_e$. We eliminate $y$ and split
$W=\kappa u^2+\eta$ as in the proof of Theorem~\ref{thm:local-pole}, near
the given value $\kappa\ne0$. The boundary variables
$f_\partial=(f_e)_{e\in\partial U}$ enter \eqref{pole:pU} only through a
term divisible by $\zeta$, so the Jacobian matrix with respect to $W$ at
$\zeta=0$ is again \eqref{pole:J}, for the matrix $(D+A^s)[U]$. Hence the
projection of \eqref{pole:pU} to $H$ is equivalent to
$\eta=\eta(\kappa,f_\partial,\zeta)$ with $\eta$ holomorphic and
$\eta(\kappa,f_\partial,0)=0$, and the remaining scalar equation is
$\zeta h_1(\kappa,f_\partial,\zeta)=0$ with $h_1$ holomorphic and
\[
h_1(\kappa,f_\partial,0)
=\sum_{i\in U}\Bigl(-\kappa u_i^2-1+\frac{d^U_i}2-\sum_{e\in\partial U,\,e\ni i}f_e\Bigr)
=\abs{E(U)}-\abs U-\sum_{e\in\partial U}f_e-\kappa .
\]
For $\zeta\ne0$ the system is thus equivalent, near the limit point, to the
system consisting of $h_1=0$, of \eqref{pole:bd-edge} with
$W_i=\kappa u_i^2+\eta_i(\kappa,f_\partial,\zeta)$, of
\eqref{pole:bd-vertex}, \eqref{pole:frozen} and of the equations of $T$, in
the unknowns
\[
\bigl((f_e)_{e\in E_N};\ (V_y)_{y\in T},(f_e)_{e\in E(T)};\ \kappa,f_\partial,(\nu_x)_{x\in N}\bigr),
\]
and this system depends holomorphically on $\zeta$ near $\zeta=0$. At
$\zeta=0$ it reads: $f_e(1+f_e)=0$ for $e\in E_N$; the remaining system
\eqref{pole:remaining} with $\phi$ replaced by $(f_e)_{e\in E_N}$; and the
boundary system \eqref{pole:boundary}, again with $\phi$ replaced by
$(f_e)_{e\in E_N}$. With the three groups of unknowns and equations in this
order, its Jacobian matrix at the limit data is block lower triangular. The
diagonal blocks are $\diag(1+2\phi_e)=\diag(\pm1)$, the Jacobian matrix of
the remaining system, and the Jacobian matrix of the boundary system, and
they are invertible by (ii) and (iii). By the implicit function theorem
there is exactly one solution near the limit data for each small $\zeta$,
and it is holomorphic in $\zeta$. This gives the branch, the expansions
\eqref{pole:partial-expansion} and the uniqueness statement.

For small $\zeta\ne0$ the solution is admissible: $V_i=W_i/\zeta\ne0$ for
$i\in U$ because $\kappa u_i^2\ne0$; $V_x=\zeta\nu_x\ne0$ for $x\in N$;
$V_y\ne0$ for $y\in T$ by (iii); and $z=\lambda+\zeta\notin\{D_i\}$. Only
the $V_i$ with $i\in U$ have poles, so the residue of $\sum_iV_i$ is
$\kappa\,u^{\T}u=\kappa$. Summing the last equation of \eqref{pole:boundary}
over $x\in N$, where an edge of $E_N$ with both ends in $N$ occurs twice,
and inserting the result in the first equation gives
\eqref{pole:kappa-integer}.
\end{proof}

\begin{remark}\label{pole:rem-nu}
Condition (i) is not used in the proof; it holds in all our applications,
and without it the remaining system is no longer the covariance system of
$G[T]$ at an admissible value of $z$. The condition $\nu_x\ne0$ cannot be omitted: a nondegenerate solution of the
boundary system with $\nu_x=0$ still continues to a branch of solutions of
\eqref{eq:system}, but on it $V_x$ vanishes identically, and the branch is
not admissible.
\end{remark}

The simplest case is a single vertex $U=\{v\}$. Then $E(U)=\emptyset$,
$\lambda=D_v$, $u=1$, and $N$ is the set of neighbors of $v$.

\begin{proposition}[Poles at the potential of a vertex]\label{prop:degree-two-vertex}
Let $v$ be a vertex whose edges $vx$ go to pairwise distinct neighbors
$x\in N$, and suppose $D_v\notin\{D_x:x\ne v\}$. Let $E_N$ be the set of
edges that have an end in $N$ and are not incident to $v$, let
$\phi\in\{0,-1\}^{E_N}$, and let $m_x$ be the number of edges $e\in E_N$ at
$x$ with $\phi_e=-1$. Suppose that $m_x\ge2$ for all $x\in N$, and that the
remaining system \eqref{pole:remaining} on $T=V\setminus(N\cup\{v\})$ has a
nondegenerate solution with $V^*_y\ne0$ for all $y\in T$. Then there is a
branch of admissible solutions near $z=D_v$, unique in the sense of
Theorem~\ref{thm:partial-pole}, on which
\[
V_v=\frac{\kappa}{z-D_v}+O(1),\qquad\kappa=-1-\sum_{x\in N}(m_x-1),
\]
$V_x\to0$ for $x\in N$, and the other $V_y$ stay finite. In particular, if
$v$ has degree two with neighbors $a\ne b$, then
$\Res_{z=D_v}V_v=-(m_a+m_b-1)$. If moreover every edge not incident to $v$
has an end in $\{a,b\}$, the remaining system has no edge unknowns, and the
condition on it is that every $y\in T$ is incident to a number different
from one of edges $e$ with $\phi_e=-1$.
\end{proposition}

\begin{proof}
For $U=\{v\}$ the last equation of \eqref{pole:boundary} gives
$f^*_{vx}=m_x-1$, the first gives $\kappa=-1-\sum_x(m_x-1)\le-1-\abs N$,
and the second gives $\nu_x=m_x(m_x-1)/(w_{vx}^2\kappa)$, which is nonzero
exactly when $m_x\ge2$. The Jacobian matrix of the boundary system is
triangular in the order $(f^*,\kappa,\nu)$ with diagonal entries $1$, $1$
and $-w_{vx}^2\kappa$, so the solution is nondegenerate. If $E(T)=\emptyset$,
the remaining system is $(D_y-D_v)V^*_y=1-m_y$ with $m_y$ the number of
edges at $y$ frozen at $-1$. Now apply Theorem~\ref{thm:partial-pole}.
\end{proof}

\begin{example}[Diamond, $K_{2,3}$ and subdivided $D_3$]\label{pole:ex-diamond}
Let $G$ be the diamond, $K_4$ with one edge removed, with vertices $a,b$ of
degree three and $v,v'$ of degree two. For $U=\{v\}$ we have $N=\{a,b\}$,
$T=\{v'\}$ and $E_N=\{ab,av',bv'\}$. The condition $m_a,m_b\ge2$ forces
$\phi\equiv-1$, the vertex $v'$ is then incident to two edges frozen at
$-1$, and Proposition~\ref{prop:degree-two-vertex} gives one branch on which
$V_v$ has a simple pole at $z=D_v$ with residue $-3$, provided
$D_v\notin\{D_a,D_b,D_{v'}\}$. The same holds at $v'$. In the same way each
vertex of degree two of $K_{2,3}$, and the subdividing vertex of the graph
obtained from $D_3$ by subdividing one edge, carries one such branch with
residue $-3$.
\end{example}

\begin{example}[Bowtie]\label{pole:ex-bowtie}
Let $G$ be the bowtie, with triangles $\{c,a,b\}$ and $\{c,a',b'\}$ sharing
the vertex $c$, and let $U=\{a,b\}$. Then $E(U)=\{ab\}$, $N=\{c\}$,
$\partial U=\{ac,bc\}$, $T=\{a',b'\}$, $E_N=\{ca',cb'\}$ and
$E(T)=\{a'b'\}$. By \eqref{pole:kappa-integer},
$\kappa=-\#\{e\in E_N:\phi_e=-1\}$. For $\phi\equiv-1$ we get $\kappa=-2$
and the boundary system
\[
f^*_{ac}+f^*_{bc}=1,\qquad f^*_{ac}(1+f^*_{ac})=-2w_{ac}^2u_a^2\,\nu_c,\qquad
f^*_{bc}(1+f^*_{bc})=-2w_{bc}^2u_b^2\,\nu_c,
\]
which after elimination of $\nu_c$ is a quadratic equation; the remaining
system has the solutions $f^*_{a'b'}=0$, which has $V^*=0$, and
$f^*_{a'b'}=1/(c'-1)$ with
$c'=w_{a'b'}^2/((D_{a'}-\lambda)(D_{b'}-\lambda))$. For generic coefficients
this gives two branches with $\Res\sum_iV_i=-2$ for each of the two
eigenvalues of the $2\times2$ matrix $(D+A)[U]$. For $\phi\equiv0$ we have
$\kappa=0$, and if exactly one edge is frozen at $-1$, the remaining system
forces $f^*_{a'b'}\in\{0,-1\}$ for generic coefficients and has no solution
with $V^*_{a'},V^*_{b'}\ne0$. Together with the symmetric choice
$U=\{a',b'\}$ we obtain eight branches with total residue $-2$ each.
\end{example}

\subsection{The total residue}
\label{subsec:pole-residue-sum}

Let $\widetilde\Xc$ be the normalization of the projective closure of the
spectral curve $\Xc$, a compact Riemann surface, possibly disconnected. The
coordinates $z$, $V_i$ and $f_e$ are meromorphic functions on
$\widetilde\Xc$, and
\[
\omega=\Bigl(\sum_{i\in V}V_i\Bigr)dz
\]
is a meromorphic differential on $\widetilde\Xc$.

\begin{proposition}[Total residue]\label{prop:residue-sum}
For all coefficients $D\in\C^V$, $w\in(\C^*)^E$, the residue of $\omega$ at
the point of $\widetilde\Xc$ over $z=\infty$ on the branch $\beta_S$,
$S\in\cL(G)$, is $n-2\abs S$, and the sum of the residues of $\omega$ at
the points of $\widetilde\Xc$ over finite values of $z$ is
\begin{equation}\label{pole:residue-sum}
\sum_{S\in\cL(G)}\bigl(2\abs S-n\bigr).
\end{equation}
\end{proposition}

\begin{proof}
By Corollary~\ref{cor:degree-all}(c) the points of $\widetilde\Xc$ over
$z=\infty$ are the $\LG$ points defined by the branches $\beta_S$,
$S\in\cL(G)$, with local coordinate $\tau=1/z$. On $\beta_S$ we have
$V_i=-\tau g_i/(1-D_i\tau)$ with $g_i=1+\sum_{e\ni i}f^S_e$ holomorphic and
$g_i(0)=1-d_S(i)$, and $dz=-d\tau/\tau^2$, so
\[
\omega=\sum_{i\in V}\frac{g_i(\tau)}{\tau(1-D_i\tau)}\,d\tau,\qquad
\Res_{\tau=0}\omega=\sum_{i\in V}\bigl(1-d_S(i)\bigr)=n-2\abs S .
\]
The sum of all residues of a meromorphic differential on a compact Riemann
surface is zero.
\end{proof}

\begin{example}\label{pole:ex-bookkeeping}
We compare \eqref{pole:residue-sum} with the poles constructed above. The
comparison uses Theorem~\ref{thm:classification}, according to which, for
generic coefficients, the poles of Theorem~\ref{thm:partial-pole} are the
only points over finite $z$ where $\omega$ has a nonzero residue, and which
also gives the number of solutions of the boundary systems.

For a cycle, $\cL(G)=\{\emptyset,E\}$ and the sum is $-n+n=0$.

For $D_3$ the sum is $-2+3\cdot2+4=8$, the contribution of the eight poles
of Example~\ref{pole:ex-D3}, each with total residue $\kappa_0=1$.

For the diamond the leafless subsets are $\emptyset$, two triangles, the
$4$-cycle and $E$, and the sum is $-4+2\cdot2+4+6=10$. The sixteen poles of
Theorem~\ref{thm:local-pole} contribute $16$, and the two poles of
Example~\ref{pole:ex-diamond} contribute $-6$.

For the bowtie the sum is $-5+2\cdot1+7=4$. The twenty poles of
Theorem~\ref{thm:local-pole} contribute $20$ and the eight poles of
Example~\ref{pole:ex-bowtie} contribute $-16$.

For $K_4$ the leafless subsets are $\emptyset$, four triangles, three
$4$-cycles, six subsets with five edges and $E$, and the sum is
$-4+4\cdot2+3\cdot4+6\cdot6+8=60$. The poles are as follows. For $U=V$
there are $4\cdot2^3=32$ poles with $\kappa=\kappa_0=2$. For each of the
four triangles $U$ we have $\abs N=1$, $E_N=\emptyset$ and $\kappa=1$; the
boundary system has one solution with $\nu\ne0$ for each of the
$3\cdot2$ pairs of an eigenvalue and a switching class, which gives $24$
poles. For each of the six edges $U$, $N$ consists of the other two
vertices and $E_N$ of the edge joining them; this edge must be frozen at
$-1$, then $\kappa=-1$, and there is one solution of the boundary system
for each of the two eigenvalues, which gives $12$ poles. For each of the
four vertices $U=\{v\}$, Proposition~\ref{prop:degree-two-vertex} with the
opposite triangle frozen at $-1$ gives one pole with $\kappa=-4$. The total
is $32\cdot2+24\cdot1-12\cdot1-4\cdot4=60$.
\end{example}

In general the poles of Theorem~\ref{thm:local-pole} contribute at most
$\kappa_0n2^\beta$ to \eqref{pole:residue-sum}, and the difference is
carried by poles supported on proper subsets, such as the poles at
$z=D_v$ of Proposition~\ref{prop:degree-two-vertex}.

\subsection{The Floquet identity}
\label{subsec:pole-floquet}

\begin{proposition}[Floquet identity]\label{prop:floquet}
For all coefficients $D\in\C^V$, $w\in(\C^*)^E$, the function
\begin{equation}\label{pole:Phi}
\Phi=\frac{\prod_{e\in E}(1+f_e)}{\prod_{i\in V}V_i}
\end{equation}
is a meromorphic function on $\widetilde\Xc$, not identically $0$ or
$\infty$ on any component, and
\begin{equation}\label{pole:floquet}
\Bigl(\sum_{i\in V}V_i\Bigr)dz=-d\log\Phi
\end{equation}
on $\widetilde\Xc$, that is, on all sheets. If $G$ is connected and the
coefficients are real, then on the physical branch, for
$z\in\C\setminus[\Sigm,\infty)$,
\begin{equation}\label{pole:Phi-phys}
\Phi(z)=\frac{\prod_{e=ij\in E}\bigl(1-w_e^2m_{i|e}(z)m_{j|e}(z)\bigr)^{-1}}{\prod_{i\in V}G_i(z)}
=\exp\Bigl(n\int\log(t-z)\,dk(t)\Bigr),
\end{equation}
where $dk=\frac1n\sum_{i\in V}d\mu_i$ is the density of states of $\tH$,
$\mu_i$ is the spectral measure of $\tH$ at a vertex over $i$, and
$\log$ is the principal branch.
\end{proposition}

The right-hand side of \eqref{pole:Phi-phys} is the Floquet function of
\cite{BBGSS}, and the second equality in \eqref{pole:Phi-phys} is their
Floquet formula, proved there for connected leafless graphs by a comparison
of logarithmic derivatives. In terms of belief propagation, $1+f_e$ is
$(1-\rho_e^2)^{-1}$ with $\rho_e$ the correlation coefficient of the edge
belief, and $\log\Phi=-\sum_i\log V_i-\sum_e\log(1-\rho_e^2)$ is the
Bethe approximation of $\log\det(D+A-z)$ at the fixed point. The identity \eqref{pole:floquet}
extends the relation between $\Phi$ and $\sum_iG_i$ from the physical
branch to all fixed points.

\begin{proof}
On each irreducible component of $\Xc$ the admissible points are dense, and
at an admissible point $V_i\ne0$ and
$f_e(1+f_e)=w_e^2V_iV_j\ne0$. Hence $\Phi$ is a well-defined meromorphic
function, not identically $0$ or $\infty$. Since $z$ is not constant on any
component, we can differentiate with respect to $z$; write $'=d/dz$ and
$r_i=V_i'/V_i$. The edge equations give
\[
f_e'(1+2f_e)=w_e^2(V_iV_j)'=f_e(1+f_e)(r_i+r_j)\qquad(e=ij),
\]
and the vertex equations give $-V_i+(D_i-z)V_i'=\sum_{e\ni i}f_e'$, that is,
\[
V_i=r_i\Bigl(1+\sum_{e\ni i}f_e\Bigr)-\sum_{e\ni i}f_e' .
\]
Summing over $i$, where each edge occurs at its two ends,
\[
\sum_iV_i=\sum_ir_i+\sum_{e=ij}\bigl((r_i+r_j)f_e-2f_e'\bigr)
=\sum_ir_i-\sum_e\frac{f_e'}{1+f_e},
\]
because
$(r_i+r_j)f_e(1+f_e)-2f_e'(1+f_e)=f_e'(1+2f_e)-2f_e'(1+f_e)=-f_e'$. The
right-hand side is $-(\log\Phi)'$.

On the physical branch $V_i=G_i$ and
$1+f_e=(1-w_e^2m_{i|e}m_{j|e})^{-1}$ by
Proposition~\ref{prop:physical-branch}, which is the first equality in
\eqref{pole:Phi-phys}. For the second, $G_i(z)=\int d\mu_i(t)/(t-z)$ gives
$\frac{d}{dz}\,n\!\int\log(t-z)\,dk(t)=-\sum_iG_i(z)$, so by
\eqref{pole:floquet} the two sides of the second equality have the same
logarithmic derivative on $(-\infty,\Sigm)$. As $z\to-\infty$ the physical
branch is $\beta_\emptyset$ (Corollary~\ref{cor:degree-all}(d)), so
$1+f_e\to1$ and $-zG_i\to1$ by \eqref{deg:asymptotics}, and both sides are
$(-z)^n(1+o(1))$. Hence they agree for $z<\Sigm$, and then on the connected
set $\C\setminus[\Sigm,\infty)$, where both are holomorphic, by analytic
continuation.
\end{proof}

\begin{corollary}\label{pole:cor-orders}
At every point $P$ of $\widetilde\Xc$ the residue of $\omega$ is the
integer $-\ord_P\Phi$. At the point over $z=\infty$ on the branch
$\beta_S$, $\ord\Phi=2\abs S-n$, and at the point defined by a branch of
Theorem~\ref{thm:partial-pole}, $\ord\Phi=-\kappa$.
\end{corollary}

\begin{proof}
The first statement is \eqref{pole:floquet}. The other two follow from
Proposition~\ref{prop:residue-sum} and \eqref{pole:kappa-integer}; they can
also be read off from the orders of the factors of $\Phi$. On $\beta_S$
each $V_i$ has a simple zero in $\tau$, and for $e=ij\in S$ the factor
$1+f_e=f_e^{-1}c_e\,g_ig_j$ has a double zero, by \eqref{deg:c-of-z}. At a
pole of Theorem~\ref{thm:partial-pole}, $V_i$ $(i\in U)$ and $1+f_e$
$(e\in E(U))$ have simple poles, $V_x$ $(x\in N)$ has a simple zero,
$1+f_e=f_e^{-1}w_e^2V_xV_y$ has a zero of order one or two for an edge
$e\in E_N$ frozen at $-1$, according to whether its other end is in $T$ or
in $N$, and the other factors are finite and nonzero.
\end{proof}

In particular Proposition~\ref{prop:residue-sum} says that the principal
divisor of $\Phi$ has degree zero.

\begin{corollary}[Sign of the residues on the physical sheet]\label{cor:antibound-sign}
Let $G$ be connected with real coefficients.
\begin{enumerate}
\item[(a)] If $\lambda$ is an isolated point of $\sigma(\tH)$, the physical
branch $G_i$ has at most a simple pole at $\lambda$, with residue
$-\mu_i(\{\lambda\})\le0$, and
$\Res_{z=\lambda}\sum_iG_i=-n\,k(\{\lambda\})$.
\item[(b)] Let $(V,f)$ be a branch of Theorem~\ref{thm:local-pole} with
$\kappa_0>0$, or a branch of Theorem~\ref{thm:partial-pole} with $\kappa>0$,
at a real $\lambda$. Then for every $i$ in the support of the pole
$(i\in V$, respectively $i\in U)$ and all sufficiently small
$\epsilon>0$, the values $V_i(\lambda\pm i\epsilon)$ differ from
$G_i(\lambda\pm i\epsilon)$. Hence the branch does not coincide with the
physical branch on either side of the real axis near $\lambda$: a pole with
positive residue is never a pole of the Green function on the physical
sheet, and if it is reached by analytic continuation of the Green function,
it is an antibound pole.
\end{enumerate}
\end{corollary}

\begin{proof}
(a) follows from $G_i(z)=\int d\mu_i(t)/(t-z)$. (b) Since the matrix
$(D+A^s)[U]$ is real symmetric, $u$ is real, and by
\eqref{pole:expansion}, \eqref{pole:partial-expansion},
\[
\im V_i(\lambda\pm i\epsilon)=\mp\frac{\kappa u_i^2}\epsilon+O(1),
\]
which is negative in the upper and positive in the lower half-plane for
small $\epsilon$, whereas $\im G_i(z)$ has the sign of $\im z$.
\end{proof}

Corollary~\ref{cor:antibound-sign} does not say on which sheet a pole with
positive residue lies, and it does not exclude poles with $\kappa<0$ from
the physical sheet. By (a) and Corollary~\ref{pole:cor-orders}, a branch of
Theorem~\ref{thm:partial-pole} that coincides with the physical branch near
an isolated eigenvalue $\lambda$ of $\tH$ has
\[
n\,k(\{\lambda\})=-\kappa=\abs U-\abs{E(U)}-\abs N+\#\{e\in E_N:\phi_e=-1\}+\#\{e\in E_N\cap E(N):\phi_e=-1\}.
\]

\begin{remark}[Aomoto's index theorem]\label{pole:rem-aomoto}
Aomoto's index theorem, in the form stated in \cite{BBGSS}, expresses the
mass $n\,k(\{\lambda\})$ of an eigenvalue of $\tH$ as
$\#X_1-\#\partial X_1-\#E_1$, where $X_1$ is the set of vertices of $G$
over which some eigenfunction does not vanish, $\partial X_1$ the set of
its neighbors outside, and $E_1$ the set of edges inside $X_1$; see also
\cite{BGVM}. The first three terms of the formula above have the same form,
with $U$, $N$ and $E(U)$ in these roles. We do not claim more than this
formal correspondence: the formula above contains the additional terms of
the frozen edges, it is derived for the branches of
Theorem~\ref{thm:partial-pole} under the nondegeneracy assumptions made
there, and it concerns all sheets; it is neither a consequence nor a
generalization of Aomoto's theorem.
\end{remark}

%% file: sections/perron.tex
\providecommand{\perum}{\bar u^{-}}
\providecommand{\perup}{\bar u^{+}}
\providecommand{\perC}{\widetilde C}

\section{The Perron antibound poles}\label{sec:perron}

In this section $G$ is a connected multigraph without self-loops and the
coefficients $D\in\R^V$, $w\in(\R\setminus\{0\})^E$ are real and otherwise
\emph{arbitrary}: no genericity and no absence of symmetry is assumed,
except where we say so explicitly. We prove
Theorem~\ref{thm:main-perron}: when $\beta\ge2$ the physical component
$\Cphys$ contains the two pole branches of Theorem~\ref{thm:local-pole}
at the Perron values $\sigm$ and $\sigp$, the lower one is reached from
the physical branch along a real path of positive solutions, the value
$\sigm$ is the threshold of walk-summability, and for every $z$ strictly
between $\sigm$ and $\Sigm$ the system \eqref{eq:system} has at least two
positive solutions. We recall that $\abs A$ is formed edge by edge
(Section~\ref{subsec:graphs}); with the entrywise absolute value of $A$
the statements below are false for multigraphs whose parallel edges
carry weights of both signs.

\subsection{Statement of the results}\label{per:subsec:statements}

We write $\kappa_0=\abs E-n=\beta-1$ and use Theorem~\ref{thm:local-pole}
in the following form. Let $\beta\ge2$. The value $\sigm$ is a simple
eigenvalue of $D-\abs A=D+A^s$ for the signing $s_e=-\sgn w_e$, its
eigenvector $\perum$ has no zero entry, and $\sigm<\min_iD_i$ because
$\sigm\le\min_i(D-\abs A)_{ii}$ with equality only if a coordinate vector
is an eigenvector. Put $\zeta=z-\sigm$, $W=\zeta V$ and
$y_e=\zeta(1+2f_e)$. Theorem~\ref{thm:local-pole} gives $r>0$ and analytic
functions $W(\zeta)$, $y(\zeta)$ on $\abs\zeta<r$ with
\begin{equation}\label{per:eq:pole-data}
W_i(0)=\kappa_0(\perum_i)^2,\qquad y_e(0)=2\kappa_0\abs{w_e}\,\perum_i\perum_j
\quad(e=ij),
\end{equation}
such that $V=W(\zeta)/\zeta$, $f=(y(\zeta)/\zeta-1)/2$ is an admissible
solution of \eqref{eq:system} for $0<\abs\zeta<r$, and such that every
solution of \eqref{eq:system} with $\zeta\ne0$ for which $(\zeta,\zeta V,\zeta(1+2f))$
is sufficiently close to $(0,W(0),y(0))$ is of this form. We call this
family the \emph{lower Perron pole branch} $P_-$. On it
\[
V_i(z)=\frac{\kappa_0(\perum_i)^2}{z-\sigm}+O(1),
\]
so that each $V_i$ has a simple pole at $z=\sigm$, with positive residue,
and is single valued near $\sigm$. The \emph{upper Perron pole branch}
$P_+$ is defined in the same way from $\sigp$, the signing $s_e=\sgn w_e$
and $\perup$; here $\sigp>\max_iD_i$ and
$V_i(z)=\kappa_0(\perup_i)^2/(z-\sigp)+O(1)$.

By Lemma~\ref{lem:infinity} there is $R>0$ such that for $\abs z>R$ all
solutions of \eqref{eq:system} lie on $2^{\abs E}$ branches which are
analytic in $\tau=1/z$ and are indexed by the subsets $S\subseteq E$,
with $f_e\to-1$ for $e\in S$ and $f_e\to0$ for $e\notin S$. Along the
physical branch $G_i(z)\to0$ and $m_{i|e}(z)\to0$ as $\abs z\to\infty$, so
$f^{\mathrm{phys}}\to0$, and the physical branch is the branch $S=\emptyset$
both for $z<-R$ and for $z>R$. Distinct branches have distinct limits, so
near a point of the branch $S=\emptyset$ with $\abs z>R$ the solution set
of \eqref{eq:system} is this branch alone. Thus $\Cphys$ is the unique
irreducible component of the closure of $\cZ$ containing it, complex
conjugation maps $\Cphys$ to itself, and the branch $S=\emptyset$ defines
one point $P_\infty$ of the normalization $\perC$ of the projective closure
of $\Cphys$, at which $\tau$ is a real local coordinate. The surface
$\perC$ is a compact connected Riemann surface with an antiholomorphic
involution; its fixed point set $\perC(\R)$ is a finite disjoint union of
circles, and the coordinates $z,V_i,f_e$ are meromorphic functions on
$\perC$ which are real on $\perC(\R)$. Let $F\subset\perC$ be the finite
set of their poles and $\nu\colon\perC\setminus F\to\C^{1+n+\abs E}$ the
normalization map.

\begin{theorem}\label{thm:perron-component}
Let $G$ be a connected loopless multigraph with $\beta\ge2$ and let
$(D,w)$ be arbitrary real coefficients. Then the physical component
$\Cphys$ contains the Perron pole branches $P_-$ and $P_+$. More
precisely, let $Q_\pm\in\perC$ be the points defined by $P_\pm$ and let
$O$ be the circle of $\perC(\R)$ through $P_\infty$.
\begin{enumerate}
\item[\textup{(a)}] The points $Q_-$ and $Q_+$ lie on $O$. At $Q_\pm$ the
function $z-\sigma_\pm$ is a local coordinate and each $V_i$ has a simple
pole with residue $\kappa_0(\bar u^{\pm}_i)^2>0$.
\item[\textup{(b)}] Let $\Gamma_-$ be the open arc of $O$ that leaves
$P_\infty$ on the side $\tau<0$ and ends at the first point of $F$. Its
other end is $Q_-$, and $\nu$ maps $\Gamma_-$ into the set $\Xp$ of
positive solutions.
\item[\textup{(c)}] The image $\nu(\Gamma_-)$ lies in $\{z\le\Sigm\}$. In
the order of $\Gamma_-$ it consists of the physical branch over
$(-\infty,\Sigm)$, the point
$p^*=\lim_{z\uparrow\Sigm}(z,\Vphys(z),f^{\mathrm{phys}}(z))$, and positive
solutions that do not lie on the physical branch; along the last part
$z\to\sigm$ from above as the point tends to $Q_-$.
\item[\textup{(d)}] The arc $\Gamma_+$ of $O$ that leaves $P_\infty$ on
the side $\tau>0$ and ends at the first point of $F$ ends at $Q_+$; on it
$z\ge\Sigp$, $V<0$ and $f>0$.
\end{enumerate}
\end{theorem}

The proof occupies Sections~\ref{per:subsec:domain}--\ref{per:subsec:proof}.
It uses only two properties of the positive solutions, valid for all
coefficients: they exist only for $z\le\Sigm$
(Lemma~\ref{lem:domain}), and they can become unbounded only at the
Perron pole (Lemma~\ref{lem:blowup}). For generic coefficients $\Xp$ is a
smooth curve and the theorem can be proved by following a connected
component of $\Xp$ (Corollary~\ref{cor:arc-circles}). For symmetric
coefficients $\Xp$ has singular points at which a real path could pass to
another irreducible component; working on the real points of the
normalization of $\Cphys$ avoids this. The consequences are drawn in
Sections~\ref{per:subsec:csz}--\ref{per:subsec:negative}: the conjecture of
Christiansen, Simon and Zinchenko (Corollary~\ref{cor:csz},
Proposition~\ref{prop:self-loops}), the case $\beta=1$
(Proposition~\ref{prop:unicyclic}), the two positive solutions in the gap
(Corollary~\ref{cor:two-roots}) and the threshold
(Proposition~\ref{prop:two-lift-threshold}).

\subsection{Where positive solutions exist and where they blow up}
\label{per:subsec:domain}

Let $(z,V,f)$ be a positive solution and put $t_i=D_i-z>0$. By
\eqref{eq:m-from-Vf} the numbers $m_{i|e}=V_i/(1+f_e)=f_e/(w_e^2V_j)$ are
positive, and the vertex equation of \eqref{eq:system} together with
$f_e/V_i=w_e^2m_{j|e}$ gives
\begin{equation}\label{per:eq:cavity}
\frac1{V_i}=t_i-\sum_{e\ni i,\ e=ij}w_e^2m_{j|e},\qquad
\frac1{m_{i|e}}=\frac{1+f_e}{V_i}=t_i-\sum_{\substack{e'\ni i,\ e'\ne e\\ e'=ij'}}w_{e'}^2m_{j'|e'} .
\end{equation}

\begin{lemma}[Domain of positive solutions]\label{lem:domain}
Let $G$ be connected with $n\ge2$ and let the coefficients be real. Then
$\Sigm<\min_iD_i$, and every positive solution $(z,V,f)$ of
\eqref{eq:system} has $z\le\Sigm$.
\end{lemma}

\begin{proof}
For the first claim let $x\in\tree$ lie over $i$, let $y$ be a neighbor
of $x$ along an edge of weight $w$, lying over $j$, and
$\psi=\delta_x-\varepsilon\sgn(w)\delta_y$. Then
$\langle\psi,\tH\psi\rangle/\norm\psi^2=(D_i-2\varepsilon\abs w+\varepsilon^2D_j)/(1+\varepsilon^2)<D_i$
for small $\varepsilon>0$.

Let $(z,V,f)$ be a positive solution and lift the numbers $m_{i|e}$ to the
oriented edges of $\tree$. Let $T'$ be a finite subtree of $\tree$ with a
root $\varrho$, and eliminate the vertices of $(\tH-z)[T']$ from the
leaves towards the root. The pivot at a vertex $v$ over $i$ is
$\pi_v=t_i-\sum_cw_{vc}^2/\pi_c$, the sum over the children $c$ of $v$ in
$T'$. We claim that $\pi_v\ge1/m_{v}>0$ for every $v\ne\varrho$, where
$m_v$ is the lifted value $m_{i|e}$ for the edge $e$ joining $v$ to its
parent. Indeed, if this holds for the children of $v$, then
$w_{vc}^2/\pi_c\le w_{vc}^2m_c$, and adding the positive terms of the
neighbors of $v$ that are not in $T'$ and are not the parent gives
$\pi_v\ge1/m_v$ by the second identity of \eqref{per:eq:cavity}. In the
same way the first identity gives $\pi_\varrho\ge1/V_\varrho>0$. All
pivots are positive, hence $(\tH-z)[T']$ is positive definite. (The
elimination along a tree produces no fill-in, and the weights enter only
through $w^2$.) Every finitely supported $\psi$ is supported in a finite
subtree, so $\langle\psi,(\tH-z)\psi\rangle\ge0$ on a dense subspace and
$z\le\inf\sigma(\tH)$.
\end{proof}

Next we rewrite \eqref{eq:system} for positive solutions. With
$f_+$ as in \eqref{eq:fplus} we have $f_+(x)=\sqrt x-\tfrac12+r(x)$,
\[
r(x)=\frac1{2\bigl(\sqrt{1+4x}+2\sqrt x\bigr)},\qquad
0<r(x)\le\min\Bigl\{\frac12,\frac1{8\sqrt x}\Bigr\}.
\]
Put $v_i=\sqrt{V_i}$ and $r_e=r(w_e^2V_iV_j)$. Since
$f_e=f_+(w_e^2V_iV_j)$, the vertex equation reads
\begin{equation}\label{per:eq:sqrt-form}
v_i\bigl((D-z-\abs A)v\bigr)_i=1-\frac{d_i}2+\sum_{e\ni i}r_e\qquad(i\in V).
\end{equation}

\begin{lemma}[Localization of blow-up]\label{lem:blowup}
Let $G$ be connected with $n\ge2$ and real coefficients.
\begin{enumerate}
\item[\textup{(a)}] Every positive solution with $z<\sigm$ satisfies
$\sum_iV_i\le n/(\sigm-z)$.
\item[\textup{(b)}] Let $(z_k,V^{(k)},f^{(k)})$ be positive solutions with
$\norm{V^{(k)}}\to\infty$, and put $\zeta_k=z_k-\sigm$. Then
$\zeta_k\to0$ and $\zeta_kV^{(k)}_i\to\kappa_0(\perum_i)^2$ for all $i$.
If $\beta\ge2$, then $\zeta_k>0$ for large $k$, and
$\zeta_k(1+2f^{(k)}_e)\to2\kappa_0\abs{w_e}\perum_i\perum_j$ for every edge
$e=ij$.
\end{enumerate}
\end{lemma}

\begin{proof}
(a) Summing \eqref{per:eq:sqrt-form} over $i$ and using $r_e\le\tfrac12$
gives $v^{\T}(D-z-\abs A)v=n-\abs E+2\sum_er_e\le n$, while the left-hand
side is at least $(\sigm-z)\norm v^2=(\sigm-z)\sum_iV_i$.

(b) By Lemma~\ref{lem:domain}, $z_k\le\Sigm<\min_iD_i$, and by (a)
$\liminf z_k\ge\sigm$; in particular $t_i=D_i-z_k$ stays in a compact
subset of $(0,\infty)$. Pass to a subsequence along which $z_k\to z^*$
and $v^{(k)}/\norm{v^{(k)}}\to\eta\ge0$, $\norm\eta=1$. Dividing
\eqref{per:eq:sqrt-form} by $\norm{v^{(k)}}^2$ and letting $k\to\infty$ we
get $\eta_i((D-z^*-\abs A)\eta)_i=0$ for all $i$. If $j$ is joined to $i$
by an edge $e$, then $t_jV_j\ge f_e\ge\abs{w_e}v_iv_j-\tfrac12$ and
$V_j\ge1/t_j$, whence $v_j\ge(\abs{w_e}/t_j)v_i-1/(2\sqrt{t_j})$. Thus
$\eta_j\ge c\,\eta_i$ with a constant $c>0$, and $\eta>0$ since $G$ is
connected. Therefore $(D-\abs A)\eta=z^*\eta$ with $\eta>0$, and by the
Perron--Frobenius theorem $z^*=\sigm$ and $\eta=\perum$. The limit does not
depend on the subsequence, so $\zeta_k\to0$ and
$v^{(k)}/\norm{v^{(k)}}\to\perum$ for the whole sequence. In particular
every $V^{(k)}_i\to\infty$ and $r_e\to0$.

Let $M=D-\sigm-\abs A$. It is positive semidefinite with kernel spanned
by $\perum$; let $\mu_2>0$ be its second smallest eigenvalue. Write
$a_k=\zeta_k\norm{v^{(k)}}^2$ and drop the index $k$. Summing
\eqref{per:eq:sqrt-form} gives
\[
v^{\T}Mv=a+n-\abs E+2\sum_er_e=a-\kappa_0+o(1).
\]
Since $v^{\T}Mv\ge0$ we get $a\ge\kappa_0-o(1)$. For the reverse
inequality, \eqref{per:eq:sqrt-form} gives
$(Mv)_i=\zeta v_i+(1-d_i/2+\sum_{e\ni i}r_e)/v_i$, so
$\norm{Mv}^2\le2\zeta^2\norm v^2+o(1)=2\zeta a+o(1)$, and
$v^{\T}Mv\le\norm{Mv}^2/\mu_2$ because $M\ge\mu_2$ on the orthogonal
complement of its kernel. Hence
$a(1-2\zeta/\mu_2)\le\kappa_0+o(1)$, and $\zeta\to0$ gives
$a\to\kappa_0$. Consequently
$\zeta V_i=a\,(v_i/\norm v)^2\to\kappa_0(\perum_i)^2$. If $\beta\ge2$
then $\kappa_0>0$, so $\zeta_k>0$ eventually, and
$\zeta(1+2f_e)=\zeta\sqrt{1+4w_e^2V_iV_j}\to2\abs{w_e}\kappa_0\perum_i\perum_j$.
\end{proof}

\begin{lemma}\label{per:far}
Let $G$ be connected with real coefficients. There is $R'>0$ such that
for every $z<-R'$ the physical branch is the only positive solution
of \eqref{eq:system}.
\end{lemma}

\begin{proof}
For $z<\min_iD_i$ and $V\in[0,\infty)^V$ let
\begin{equation}\label{per:eq:Phi}
\Phi_z(V)_i=\frac1{D_i-z}\Bigl(1+\sum_{e\ni i,\ e=ij}f_+(w_e^2V_iV_j)\Bigr).
\end{equation}
The positive solutions with this $z$ are the fixed points of $\Phi_z$
with positive entries. Since $0<f_+'\le1$, the sum of the $i$th row of
the Jacobian matrix $D\Phi_z$ is at most
$\sum_{e\ni i}w_e^2(V_i+V_j)/(D_i-z)$. By Lemma~\ref{lem:blowup}(a) all
positive solutions with $z<\sigm$ lie in the box
$[0,n/(\sigm-z)]^V$, on which the row sums are $O(z^{-2})$. For
$z<-R'$ they are less than $1$, so $\Phi_z$ is Lipschitz with constant less
than $1$ on this convex box in the maximum norm and has at most one fixed
point there. The physical
branch is one by Proposition~\ref{prop:physical-branch}.
\end{proof}

\subsection{The gap}\label{per:subsec:gap}

\begin{corollary}[Gap]\label{cor:gap}
Let $G$ be a connected loopless multigraph with $\beta\ge2$. For all real
coefficients
\[
\sigm<\Sigm\qquad\text{and}\qquad\Sigp<\sigp .
\]
\end{corollary}

\begin{proof}
On the pole branch $P_-$ take $\zeta\in(0,r)$ real. The coefficients are
real and the branch is unique, so $W(\zeta)$ and $y(\zeta)$ are real, and
by \eqref{per:eq:pole-data} we have $V_i=W_i(\zeta)/\zeta>0$ and
$1+2f_e=y_e(\zeta)/\zeta>0$ for small $\zeta$; together with
$f_e(1+f_e)=w_e^2V_iV_j>0$ this gives $f_e>0$. So $P_-$ consists of
positive solutions for small $\zeta>0$, and Lemma~\ref{lem:domain} gives
$\sigm+\zeta\le\Sigm$. For the upper inequality we use the map
\begin{equation}\label{per:eq:iota}
\iota\colon(D,z,V,f)\longmapsto(-D,-z,-V,f),
\end{equation}
which preserves both equations of \eqref{eq:system}. On the tree
$-\tH(D,w)=\tH(-D,-w)$ is unitarily equivalent, by a diagonal unitary, to
$\tH(-D,w)$; hence $\Sigm(-D,w)=-\Sigp(D,w)$, and clearly
$\sigm(-D,w)=\lambda_{\min}(-D-\abs A)=-\sigp(D,w)$. The lower inequality
for $(-D,w)$ is the upper one for $(D,w)$.
\end{proof}

\begin{remark}\label{per:rem:sysunada}
For simple graphs with positive weights the strict inequality
$\Sigp<\sigp$ is a consequence of the theorem of Sy and
Sunada~\cite{SySunada}, since the fundamental group of $G$ is
non-amenable exactly when $\beta\ge2$; it is quoted in this form
in~\cite{CSZ}. Corollary~\ref{cor:gap} gives a different proof, which
covers multigraphs and both ends, and which exhibits the reason in the
present language: just above $\sigm$ there are fixed points with positive
definite local beliefs, and such fixed points cannot exist above $\Sigm$.
For $\beta=1$ the gap closes
(Proposition~\ref{prop:unicyclic}).
\end{remark}

\subsection{The turning point at the bottom of the spectrum}
\label{per:subsec:turning}

For $z<\min_iD_i$ and $V$ in a complex neighborhood of $(0,\infty)^V$ put
\begin{equation}\label{per:eq:g}
g(z,V)_i=(D_i-z)V_i-1-\sum_{e\ni i,\ e=ij}f_+(w_e^2V_iV_j),\qquad
J_V=\partial_Vg=\diag(D-z)\,(I-D\Phi_z),
\end{equation}
with $\Phi_z$ as in \eqref{per:eq:Phi}. At a solution of \eqref{eq:system}
with $1+2f_e\ne0$ for all $e$, the two roots $f_e$ and $-1-f_e$ of the
edge equation are distinct, so nearby solutions are obtained by
eliminating $f_e$ as an analytic function of $V$; near a positive
solution this function is $f_+(w_e^2V_iV_j)$. Hence near a positive
solution the set $\cZ$ is biholomorphic to $\{g=0\}\subset\C^{1+n}$. For
$V>0$ the matrix $D\Phi_z$ is nonnegative, with positive entries
$(i,j)$ for adjacent $i\ne j$, hence irreducible. Note that
$\partial_zg=-V$.

The physical branch is increasing in $z$ on $(-\infty,\Sigm)$, because
$\frac{d}{dz}G_i(z)=\langle\delta_x,(\tH-z)^{-2}\delta_x\rangle>0$.

\begin{lemma}[Regularity of the turning point]\label{lem:turning}
Let $G$ be a connected loopless multigraph with $\beta\ge2$ and real
coefficients.
\begin{enumerate}
\item[\textup{(i)}] For $z<\Sigm$ the spectral radius of
$D\Phi_z$ at $\Vphys(z)$ is less than $1$. Hence $J_V$ is invertible there
and the points of the physical branch over $(-\infty,\Sigm)$ are smooth
points of $\cZ$ at which $z$ is a local coordinate.
\item[\textup{(ii)}] The limit $p^*=\lim_{z\uparrow\Sigm}\Vphys(z)$
exists and $(\Sigm,p^*)$, with $f=f_+$, is a positive solution. At
$p^*$ the spectral radius of $D\Phi_{\Sigm}$ equals $1$, the kernel and
the left kernel of $J_V$ are spanned by positive vectors, and the
$n\times(n+1)$ matrix $(\partial_zg\mid J_V)$ has rank $n$. Hence $p^*$ is
a smooth point of $\cZ$ at which $z$ has a critical point, and near $p^*$
the set $\Xp$ is a smooth real arc on which $z\le\Sigm$; one side of
$p^*$ on this arc is the physical branch and the other side consists of
positive solutions that are not on the physical branch.
\end{enumerate}
\end{lemma}

\begin{proof}
(i) Let $\rho(z)$ be the spectral radius of $D\Phi_z$ at $\Vphys(z)$; it
is continuous, and $\rho(z)\to0$ as $z\to-\infty$ because
$\Vphys=O(1/\abs z)$. Suppose $\rho(z_1)=1$ for some $z_1<\Sigm$, and let
$\ell>0$ be a left Perron vector of $D\Phi_{z_1}$. Then
$\tilde\ell=\diag(D-z_1)^{-1}\ell>0$ satisfies $\tilde\ell^{\T}J_V=0$.
Differentiating $g(z,\Vphys(z))=0$ gives $J_V\,\frac{d}{dz}\Vphys=\Vphys$,
and therefore $0=\tilde\ell^{\T}\Vphys>0$, a contradiction. So
$\rho<1$ on $(-\infty,\Sigm)$, and $I-D\Phi_z$ is invertible.

(ii) The physical branch is increasing in $z$. If it were unbounded as
$z\uparrow\Sigm$, Lemma~\ref{lem:blowup}(b) would give $z\to\sigm$,
contradicting $\sigm<\Sigm$ (Corollary~\ref{cor:gap}). So the limit $p^*$
exists. It satisfies $g(\Sigm,p^*)=0$, and $p^*_i\ge1/(D_i-\Sigm)>0$, so
with $f=f_+>0$ it is a positive solution. By continuity the spectral
radius at $p^*$ is at most $1$. If it were less than $1$, the implicit
function theorem would give real solutions of $g=0$ for $z$ on both sides
of $\Sigm$, positive near $p^*$, contradicting Lemma~\ref{lem:domain}.
So the spectral radius is $1$; it is a simple eigenvalue of the
irreducible nonnegative matrix $D\Phi_{\Sigm}$, with right and left
eigenvectors $\varrho,\ell>0$. Thus $\ker J_V=\R\varrho$ and the left
kernel is spanned by $\tilde\ell=\diag(D-\Sigm)^{-1}\ell>0$. Since
$\tilde\ell^{\T}\partial_zg=-\tilde\ell^{\T}p^*<0$, the vector
$\partial_zg$ is not in the image of $J_V$ and the rank is $n$. Therefore
$\{g=0\}$ is, near $p^*$, a smooth complex curve with tangent
$(0,\varrho)$, so $dz$ vanishes on it at $p^*$, and its real points form
a smooth real arc through $p^*$. The points of this arc near $p^*$ have
$V>0$, so they are positive solutions and $z\le\Sigm$ on the arc by
Lemma~\ref{lem:domain}. The points $(z,\Vphys(z))$ with $z$ close to
$\Sigm$ lie on the arc and converge to $p^*$; being a connected set that
does not contain $p^*$, they fill one side of $p^*$. The points of the
other side are not on the physical branch, since the physical branch is
a graph over $z$ and its points with $z\le\Sigm-\varepsilon$ stay at a
positive distance from $p^*$ by monotonicity.
\end{proof}

We do not claim that the turning point is nondegenerate, that is, that
$z$ has local degree exactly $2$ at $p^*$; the degree is even because
$z\le\Sigm$ on both sides.

\subsection{Proof of Theorem~\ref{thm:perron-component}}
\label{per:subsec:proof}

We first treat the lower end. Let $\Gamma$ be the connected component of
$O\setminus F$ that contains the points of $O$ near $P_\infty$ with
$\tau<0$. Since $P_\infty\in F$ and $F$ is finite, $\Gamma$ is an open
arc; one of its ends is $P_\infty$, approached with $\tau<0$, and we call
the other end $Q\in F$. The points of the physical branch over
$(-\infty,\Sigm)$ are smooth points of $\cZ$ by Lemma~\ref{lem:turning}(i),
so each has exactly one preimage under $\nu$; these preimages form an
open sub-arc $\Gamma_1\subset\Gamma$ starting at $P_\infty$, on which $z$ is
a coordinate with values in $(-\infty,\Sigm)$.

\emph{Step 1: $\nu(\Gamma)\subset\Xp$.} Let
$A=\{q\in\Gamma:\nu(q)\in\Xp\}$. A real point of the closure of $\Cphys$
satisfies \eqref{eq:system}, so
$A=\{q\in\Gamma: V_i(q)>0,\ f_e(q)>0\text{ for all }i,e\}$, which is open
in $\Gamma$; it is nonempty because it contains $\Gamma_1$. Let
$q_k\in A$, $q_k\to q\in\Gamma$. By Lemma~\ref{lem:domain},
$z(q_k)\le\Sigm$, so $z(q)\le\Sigm<\min_iD_i$. The vertex equation gives
$V_i(q_k)\ge1/(D_i-z(q_k))$, so $V_i(q)>0$; and $f_e(q)\ge0$ with
$f_e(1+f_e)=w_e^2V_iV_j>0$ gives $f_e(q)>0$. Thus $A$ is closed in
$\Gamma$, and $A=\Gamma$ since $\Gamma$ is connected.

\emph{Step 2: $z$ is bounded near the end $Q$.} Otherwise, as
$z\le\Sigm$ on $\Gamma$, there are $q_k\in\Gamma$ tending to the end $Q$
with $z(q_k)\to-\infty$. By Lemma~\ref{per:far} the positive solution
$\nu(q_k)$ is a point of the physical branch for large $k$, so
$q_k\in\Gamma_1$ and $q_k$ tends to the end $P_\infty$ of $\Gamma$. This is
impossible, because a sequence in an open arc cannot converge to both of
its ends.

\emph{Step 3: $Q=Q_-$.} Since $Q\in F$ and $z$ is bounded near $Q$,
some $V_i$ or $f_e$ is unbounded as $q\to Q$ along $\Gamma$; because
$0<f_e=f_+(w_e^2V_iV_j)\le\abs{w_e}\sqrt{V_iV_j}$, in fact
$\max_iV_i(q)\to\infty$. Lemma~\ref{lem:blowup}(b) shows that
$\zeta=z-\sigm\to0^+$ and that $(\zeta V,\zeta(1+2f))$ converges to the
point \eqref{per:eq:pole-data}. By the uniqueness part of
Theorem~\ref{thm:local-pole}, for $q\in\Gamma$ close to the end $Q$ the
point $\nu(q)$ lies on the pole branch $P_-$. The points of $P_-$ with
$\zeta\ne0$ are smooth points of the closure of $\cZ$, since all
solutions near them lie on $P_-$; they form a connected set, hence lie in
one irreducible component, and this component is $\Cphys$ because
infinitely many of them are in $\nu(\Gamma)$. The punctured disc
$0<\abs\zeta<r$ therefore lifts to $\perC$, the lift extends over the
puncture to a point $Q_-\in\perC$ at which $\zeta$ is a local coordinate,
and $q\to Q_-$; so $Q=Q_-$. At $Q_-$ the function $V_i=W_i(\zeta)/\zeta$
has a simple pole with residue $\kappa_0(\perum_i)^2$.

This proves (b) with $\Gamma_-=\Gamma$, and the statements of (a) about
$Q_-$. For (c): $z\le\Sigm$ on $\Gamma_-$ by Step~1 and
Lemma~\ref{lem:domain}. The point $(\Sigm,p^*)$ is a smooth point of
$\cZ$ by Lemma~\ref{lem:turning}(ii), in the closure of the physical
branch, so it has a single preimage $q^*$, which is the end of $\Gamma_1$
in $\Gamma_-$. The remaining part
$\Gamma_2=\Gamma_-\setminus(\Gamma_1\cup\{q^*\})$ is an open arc from
$q^*$ to $Q_-$. Its points are not mapped to the physical branch or to
$p^*$, because those points have their only preimages in
$\Gamma_1\cup\{q^*\}$. By Step~3, $z\to\sigm$ from above at $Q_-$.

\emph{Step 4: the upper end.} The map $\iota$ of \eqref{per:eq:iota}
sends the solutions of \eqref{eq:system} for $(D,w)$ to those for
$(-D,w)$. By the unitary equivalence used in the proof of
Corollary~\ref{cor:gap}, $G_i^{(-D,w)}(-z)=-G_i^{(D,w)}(z)$ and likewise
for the half-tree functions, so $\iota$ maps the physical branch of
$(D,w)$ over $z>\Sigp$ to the physical branch of $(-D,w)$ over
$z<\Sigm(-D,w)$. Hence $\iota$ maps $\Cphys(D,w)$ onto $\Cphys(-D,w)$ and
induces an isomorphism of the normalizations, compatible with the real
structures, which fixes $P_\infty$ and changes the sign of $\tau$. Since
$\sigm(-D,w)=-\sigp(D,w)$ with the same Perron vector, applying
Steps~1--3 to $(-D,w)$ and transporting back by $\iota$ gives (d) and the
statements of (a) about $Q_+$; the residue of $V_i$ at $Q_+$ is
$\kappa_0(\perup_i)^2$ because $\iota$ changes the signs of both $V$ and
$z-\sigma$. Finally $Q_\pm\in O$ because they are ends of arcs of $O$.
\qed

\subsection{The conjecture of Christiansen, Simon and Zinchenko}
\label{per:subsec:csz}

Christiansen, Simon and Zinchenko~\cite{CSZ} consider a finite leafless
graph, possibly with self-loops and multiple edges, with real potential
and positive edge weights, its Jacobi matrix $J$, and the periodic Jacobi
operator on the universal covering tree. Their Conjecture~3.4 states that
if the tree is not the line $\Z$, then the Green function, continued
analytically from $\C\setminus\sigma(\tH)$, has a pole on a higher sheet
at $z=\sigma$, the largest eigenvalue of $J$, and, if the graph is
bipartite, also at the smallest eigenvalue $\sigma_-$ of $J$. For a
leafless connected graph with at least one edge the tree is not $\Z$
exactly when $\beta\ge2$. For positive weights $\abs A=A$, so
$\sigma=\sigp$; and for a bipartite graph $D-A$ is conjugate to $D+A=J$ by
a diagonal matrix of signs, so $\sigm=\lambda_{\min}(J)$.

\begin{corollary}[The conjecture of Christiansen, Simon and Zinchenko]
\label{cor:csz}
Let $G$ be a connected loopless multigraph with $\beta\ge2$, with
arbitrary real coefficients. For every vertex $i$ there are paths in $\C$,
starting on the physical sheet and avoiding a finite set, along which the
Green function $G_i$ continues analytically to germs that extend to
functions meromorphic near $\sigp=\lambda_{\max}(D+\abs A)$, respectively
near $\sigm=\lambda_{\min}(D-\abs A)$, with
\[
G_i(z)=\frac{\kappa_0(\bar u^{\pm}_i)^2}{z-\sigma_\pm}+O(1),\qquad\kappa_0=\abs E-n .
\]
These simple poles are not on the physical sheet: $G_i$ itself is analytic
at $\sigma_\pm$. In particular Conjecture~3.4 of~\cite{CSZ} holds for all
loopless graphs with $\beta\ge2$ and all coefficients, and the pole at
$\lambda_{\min}(D-A)$ is present whether or not $G$ is bipartite.
\end{corollary}

\begin{proof}
Let $B\subset\C$ be the finite set consisting of the critical values of
$z\colon\perC\to\PP^1$ and of the values $z(q)\ne\infty$ for $q\in F$,
and let $\perC^\circ=\perC\setminus z^{-1}(B\cup\{\infty\})$. This is a
connected Riemann surface, $z\colon\perC^\circ\to\C\setminus B$ is a finite
unramified covering, and $V_i$ is holomorphic on $\perC^\circ$. Choose
$q_0\in\perC^\circ$ over a point $z_0<\Sigm$ of the physical branch and
$q_1\in\perC^\circ$ in the coordinate disc of $Q_-$. A path from $q_0$ to
$q_1$ in $\perC^\circ$ projects to a path $\gamma$ in $\C\setminus B$, and
the analytic continuation of the germ of $G_i=\Vphys_i$ at $z_0$ along
$\gamma$ is the germ at $z(q_1)$ of $W_i(\zeta)/\zeta$, which is
meromorphic in $\abs{z-\sigm}<r$ with the stated simple pole
(Theorem~\ref{thm:perron-component}(a)). The same argument applies to
$Q_+$. By Corollary~\ref{cor:gap}, $\sigma_\pm\notin[\Sigm,\Sigp]$, so
the function $G_i$ on $\C\setminus\sigma(\tH)$ is analytic at
$\sigma_\pm$. (Alternatively, a pole with positive residue is never on
the physical sheet, by Corollary~\ref{cor:antibound-sign}.)
\end{proof}

To our knowledge the conjecture had not been resolved, in whole or in
part, before. Three points lie outside its original statement: the lower
pole for non-bipartite graphs, graphs with leaves (which are allowed in
Corollary~\ref{cor:csz}, as long as $\beta\ge2$), and the residue. The
framework of~\cite{CSZ} allows self-loops, which we treat in
Proposition~\ref{prop:self-loops}; the case $\beta=1$ with leaves, in
which the tree is not $\Z$ either, is treated in
Proposition~\ref{prop:unicyclic}. In that case the Perron values coincide
with the outer band edges and give inverse-square-root singularities
with zero trace residue. The cycle-rank hypothesis for the extended
statement is therefore $\beta\ge2$, the same as for the gap.

Corollary~\ref{cor:csz} asserts the existence of a path of continuation.
Theorem~\ref{thm:perron-component} gives an explicit one for the lower
pole.

\begin{corollary}[A real path]\label{cor:real-path}
Under the hypotheses of Corollary~\ref{cor:csz}, let $\gamma$ be a path in
$\perC$ obtained from the arc $\Gamma_-$ by replacing a small sub-arc
around each critical point of $z$ on $\Gamma_-$ by a half circle in
$\perC$ around that point. The continuation of $G_i$ along
$z\circ\gamma$ leads from the physical sheet at $z<\Sigm$ to the pole
branch $P_-$. The path $z\circ\gamma$ runs along the real axis within
$(-\infty,\Sigm]$, except that near the value of $z$ at a
critical point of local degree $k$ it winds $k/2$ times around this
value. At the turning point $p^*$ over $\Sigm$ the degree $k$ is even;
if the turning point is nondegenerate, $k=2$ and the path makes one full
turn around $\Sigm$.
\end{corollary}

\begin{proof}
There are no poles of the $V_i$ on $\Gamma_-$, and at the points of
$\gamma$ the function $z$ is a local coordinate of $\perC$, so the germs
$V_i\circ(z|_{\text{loc}})^{-1}$ along $\gamma$ form an analytic
continuation along $z\circ\gamma$. Points of $\nu(\Gamma_-)$ that are
singular points of $\cZ$ need no detour as long as $z$ is a local
coordinate at the corresponding point of $\perC$. If $z-z_1=c\,t^k+\dots$,
$c\neq0$, in a real local coordinate $t$ at a critical point, a half
circle in the $t$-plane from $t<0$ to $t>0$ is mapped to a path winding
$k/2$ times around $z_1$.
\end{proof}

We do not claim that $z$ is monotone on $\Gamma_2$, the part of
$\Gamma_-$ after the turning point; a half circle in the $z$-plane
around a critical value would lead to a non-real point of $\perC$ and
not to the continuation of $\Gamma_-$. For generic coefficients
$\nu(\Gamma_-)$ is the non-compact component of $\Xp$
(Corollary~\ref{cor:arc-circles}). For symmetric coefficients
$\nu(\Gamma_-)$ may pass through singular points of $\Xp$: numerical
continuation for two identical copies of a three-edge dipole joined by a
weak bridge indicates that the second part of the arc crosses a closed
curve of non-symmetric positive solutions at two pitchfork points, while
$\Gamma_-$ stays on the symmetric component and reaches the pole.

\begin{example}\label{per:ex:regular}
Let $G$ be $d$-regular, $d\ge3$, with $D=0$ and $w\equiv1$, so that
$\sigma_\pm=\pm d$ and $(\bar u^{\pm}_i)^2=1/n$. The symmetric solutions
$V_i\equiv V$, $f_e\equiv f$ of \eqref{eq:system} satisfy $-zV=1+df$ and
\[
(1-cd^2)f^2+(1-2cd)f-c=0,\qquad c=z^{-2}.
\]
The physical branch is the root with $f\to0$ as $z\to\infty$. The
discriminant $1-4c(d-1)$ vanishes at $z=\pm2\sqrt{d-1}=\Sigma_\pm$, the
turning points, and the other root has poles at $z=\pm d$ with
$V=(d-2)/(2(z\mp d))+O(1)$. The residue $(d-2)/2$ equals
$\kappa_0/n$, in agreement with Theorem~\ref{thm:perron-component}(a) and
with the second-sheet poles of the regular tree in
\cite[Example~3.1]{CSZ}.
\end{example}

\subsubsection*{Self-loops}

We now allow self-loops, with the conventions of \cite{ABS,CSZ}: a loop
$e$ at $v$ with weight $a_e>0$ contributes two edges at each lift of $v$
in the covering tree, counts as one edge in $\abs E$ and in
$\beta=\abs E-n+1$, and adds $2a_e$ to the diagonal entry of the finite
Jacobi matrix $J$ at $v$ (apply the tree operator to the lift of a
function on $V$: both orientations of the loop contribute).

\begin{proposition}[Self-loops]\label{prop:self-loops}
Let $G$ be a finite connected graph, possibly with self-loops and
multiple edges, with $\beta\ge2$, real potential and positive edge
weights $a_e$. Let $J$ be its Jacobi matrix with the convention above and
$J^-$ the matrix obtained from $J$ by replacing every $a_e$ by $-a_e$
(so a loop adds $-2a_e$ to the diagonal). Then the Green function $G_i$
of the periodic Jacobi operator on the universal covering tree continues,
along paths avoiding a finite set, to germs with simple poles on higher
sheets at $z=\lambda_{\max}(J)$ and at $z=\lambda_{\min}(J^-)$, with
residues $(\abs E-n)\bar u_i^2$, where $\bar u$ is the corresponding
positive unit eigenvector.
\end{proposition}

\begin{proof}
If $G$ has no loop this is Corollary~\ref{cor:csz}. Otherwise let
$G_2\to G$ be the two-fold cover in which \emph{every} loop is lifted
crosswise and every other edge is lifted parallel. A loop at $v$ lifts to
two parallel edges between the two lifts $v_1,v_2$ of $v$, so $G_2$ is a
loopless multigraph; it is connected because $G$ has a loop, and
$n_2=2n$, $\abs{E_2}=2\abs E$, $\beta(G_2)=2\beta-1\ge3$. With the
coefficients pulled back, $G_2$ and $G$ have the same weighted universal
covering tree, hence the same operator $\tH$ and the same Green function
at $i$ and at a lift $i_1$. Lifting the positive eigenvector $\bar u$ of
$J$ to a function constant on the fibers gives a positive eigenvector of
the Jacobi matrix $D_2+A_2$ of $G_2$ with the same eigenvalue (the two
parallel lifts of a loop contribute $2a_e\bar u_v$), and a positive
eigenvector is the Perron vector; so
$\lambda_{\max}(D_2+A_2)=\lambda_{\max}(J)$ with unit Perron vector
$\bar u/\sqrt2$ on each sheet. Likewise
$\lambda_{\min}(D_2-A_2)=\lambda_{\min}(J^-)$. Corollary~\ref{cor:csz}
for $G_2$ gives poles with residues
$(\abs{E_2}-n_2)(\bar u_i/\sqrt2)^2=(\abs E-n)\bar u_i^2$.
\end{proof}

Lifting only one loop crosswise would leave the other loops as loops of
$G_2$. With Proposition~\ref{prop:self-loops}, the whole framework of
Conjecture~3.4 of~\cite{CSZ} is covered: leafless graphs whose covering
tree is not $\Z$, with self-loops and multiple edges, and positive
weights. For non-bipartite graphs the pole at $\lambda_{\min}(J^-)$ is an
additional pole not predicted there.

\subsection{Cycle rank one}\label{per:subsec:unicyclic}

For $\beta=1$ we have $\kappa_0=0$, so the unramified pole construction of
Theorem~\ref{thm:local-pole} does not apply. The graph is then a cycle of
length $\ell\ge3$, or a pair of parallel edges ($\ell=2$), with finitely
many trees attached, and $\LG=2$.

\begin{proposition}[Cycle rank one]\label{prop:unicyclic}
Let $G$ be a connected loopless multigraph with $\beta=1$, real potentials
and nonzero real edge weights.
\begin{enumerate}
\item[\textup{(a)}] For all real coefficients $\sigm=\Sigm$ and
$\sigp=\Sigp$: the Perron values are the outer edges of the spectrum of
$\tH$.
\item[\textup{(b)}] The spectral curve has two sheets, exchanged at the
simple branch points over $\sigma_\pm$.
Each $G_i$ has a simple pole in a local coordinate of the normalization
at these points and grows as a nonzero multiple of
$(z-\sigma_\pm)^{-1/2}$ in the energy variable. The differential
$(\sum_i V_i)\,dz$ is holomorphic there.
\end{enumerate}
\end{proposition}

\begin{proof}
(a) Both sides depend only on the $\abs{w_e}$
(Section~\ref{subsec:cover}), so we may assume $w_e>0$. Let $e_0=pq$ be an
edge of the cycle. The universal cover is a $\Z$-periodic graph with
fundamental domain the tree $G-e_0$, and by the Floquet--Bloch
decomposition~\cite{Kuchment}
$\sigma(\tH)=\bigcup_{\theta\in[0,2\pi)}\spec J(\theta)$, where
$J(\theta)$ is the Hermitian matrix obtained from $D+A$ by replacing the
contribution $w_{e_0}$ of $e_0$ to the entry $(p,q)$ by
$w_{e_0}e^{i\theta}$ and to the entry $(q,p)$ by $w_{e_0}e^{-i\theta}$
(for $\ell=2$ the contribution of the other parallel edge is unchanged).
For every $\psi\in\C^V$,
\[
\abs\psi^{\T}(D-A)\abs\psi\le\psi^*J(\theta)\psi\le\abs\psi^{\T}(D+A)\abs\psi ,
\]
so $\sigm\le\lambda_{\min}(J(\theta))$ and
$\lambda_{\max}(J(\theta))\le\sigp$. Since $J(0)=D+A$ we get
$\Sigp=\sigp$. For the lower edge, the tree $G-e_0$ admits
$\varepsilon\in\{\pm1\}^V$ with $\varepsilon_i\varepsilon_j=-1$ on each of
its edges; then $\varepsilon_p\varepsilon_q=(-1)^{\ell-1}$, and
$\diag(\varepsilon)J(\theta)\diag(\varepsilon)=D-A$ when
$e^{i\theta}=(-1)^\ell$. Hence $\Sigm=\sigm$.

(b) Let $F$ be the vertices outside the cycle. The determinant has the
form
\[
 \det(J(\theta)-zI)=P(z)+Q(z)\cos\theta,\qquad
 Q(z)=\pm 2\Bigl(\prod_{e\text{ on the cycle}}w_e\Bigr)
           \det\bigl(J(0)[F]-zI\bigr).
\]
In the determinant expansion a term with uncanceled phase must traverse
the unique cycle, and the two orientations give the cosine. For a
two-edge cycle the same formula follows from the cross term of
$\abs{w_1+w_2e^{i\theta}}^2$. An empty $F$ gives determinant one.

At either outer edge $z_*=\sigma_\pm$, the matrix realizing the edge
at $\theta_*\in\{0,\pi\}$ has a simple extremal eigenvalue and an
eigenvector with all entries nonzero by Perron--Frobenius and the sign
switch in (a). Its principal submatrix on $F$ cannot have the eigenvalue
$z_*$: extending such an eigenvector by zero would attain the same
extremal Rayleigh quotient with a zero coordinate. Thus $Q(z_*)\ne0$,
and exactly one of $P+Q$ and $P-Q$ has a simple zero there.

The diagonal cofactor at $i$ also has the form
$P_i(z)+Q_i(z)\cos\theta$. By Cramer's rule and the Bloch decomposition
the diagonal Green function is the average of this cofactor divided by
$P+Q\cos\theta$. In a neighborhood of $z_*$ this gives
\[
 G_i(z)=\frac{Q_i(z)}{Q(z)}
  +\left(P_i(z)-\frac{Q_i(z)P(z)}{Q(z)}\right)
    \frac{1}{\sqrt{P(z)^2-Q(z)^2}},
\]
where the square root is chosen using the average of
$(P+Q\cos\theta)^{-1}$ on the physical branch. At $z_*$ the coefficient
of the inverse square root is the diagonal cofactor of
$J(\theta_*)-z_* I$, which is nonzero because the eigenvalue is simple
and the corresponding eigenvector has nonzero $i$-th entry. Hence each
$G_i$ has inverse-square-root behavior. By
Theorem~\ref{thm:main-sheets}(a) the total degree is two for every choice
of coefficients, so the ramification of the physical component
accounts for both sheets. In the local coordinate
$\eta=\sqrt{z-z_*}$, each $G_i$ has a simple pole and $dz=2\eta\,d\eta$,
so $(\sum_i V_i)\,dz$ is holomorphic.
\end{proof}

For a cycle without trees this is the classical picture of the
hyperelliptic curve $y^2=\Delta(z)^2-4$ of a periodic Jacobi matrix on
$\Z$, in which $\sigma_\pm$ are the outermost roots of
$\Delta=\pm2$~\cite{vMM}. A
unicyclic graph with at least one leaf has a covering tree different
from $\Z$ and lies outside the framework of~\cite{CSZ}, which requires
leafless graphs. Proposition~\ref{prop:unicyclic} gives the band-edge
behavior for these graphs.

\subsection{Two fixed points in the gap, and the threshold}
\label{per:subsec:gap-roots}

\begin{corollary}[Two positive solutions in the gap]\label{cor:two-roots}
Let $G$ be a connected loopless multigraph with $\beta\ge2$ and arbitrary
real coefficients. For every $z\in(\sigm,\Sigm)$ the system
\eqref{eq:system} has at least two positive solutions: the physical
branch and a point of $\nu(\Gamma_2)$. Equivalently, by
Lemma~\ref{lem:bp-correspondence}, Gaussian belief propagation for the
diagonal $D-z$ has at least two fixed points all of whose vertex and edge
beliefs are positive definite.
\end{corollary}

\begin{proof}
With the notation of Section~\ref{per:subsec:proof}, $z$ is a continuous
real function on the open arc $\Gamma_2$, which tends to $\Sigm$ at the
end $q^*$ and to $\sigm$ at the end $Q_-$. By the intermediate value
theorem it takes every value in $(\sigm,\Sigm)$. The points of
$\nu(\Gamma_2)$ are positive solutions
(Theorem~\ref{thm:perron-component}(b)) and do not lie on the physical
branch (Theorem~\ref{thm:perron-component}(c)). Critical points of $z$ on
$\Gamma_2$ and intersections of $\nu(\Gamma_2)$ with other families of
solutions do not affect the argument.
\end{proof}

The lower Perron value is also a threshold in the finite graph.

\begin{proposition}[Threshold]\label{prop:two-lift-threshold}
Let $G$ be a connected loopless multigraph with real coefficients, and
let $z\in\R$. The following are equivalent.
\begin{enumerate}
\item[\textup{(i)}] $z<\sigm$, that is, the comparison matrix
$D-z-\abs A$ is positive definite.
\item[\textup{(ii)}] $D-z+A^s$ is positive definite for every signing $s$.
\item[\textup{(iii)}] For every two-fold cover of $G$, the Jacobi matrix
of the cover with the pulled back coefficients $D-z$ and $w$ is positive
definite.
\end{enumerate}
\end{proposition}

\begin{proof}
For real $\psi$ and any $s$ we have
$\psi^{\T}(D+A^s)\psi\ge\abs\psi^{\T}(D-\abs A)\abs\psi\ge\sigm\norm\psi^2$,
with equality for $s_e=-\sgn w_e$ and $\psi=\perum$. Hence
$\min_s\lambda_{\min}(D+A^s)=\sigm$, which is the equivalence of (i) and
(ii). The Jacobi matrix of the two-fold cover with crossing pattern $s$
is orthogonally equivalent to $(D+A)\oplus(D+A^s)$
(Section~\ref{subsec:graphs}), and $D+A$ is the case $s\equiv1$, so (ii)
and (iii) are equivalent.
\end{proof}

For a simple graph, condition~(i) says that the matrix $D-z+A$ is
walk-summable in the sense of Malioutov, Johnson and
Willsky~\cite{MJW}, and the equivalence of (i) and (iii) is due to
Ruozzi and Tatikonda~\cite[Theorem~11]{RuozziTatikonda}; we include the
short argument because we use the edgewise convention for multigraphs.
Together with Theorem~\ref{thm:perron-component} and
Corollary~\ref{cor:two-roots} this gives the statement that connects the
three aspects of the spectral curve: \emph{for $\beta\ge2$ the lower
Perron antibound pole of the Green function sits exactly at the boundary
of walk-summability, and for every $z$ between it and the bottom of the
spectrum of the covering tree there are at least two fixed points with
positive definite local beliefs.} For $z<\sigm$ the positive solution is
unique under the hypotheses of Theorem~\ref{thm:unique-ws}.

\begin{remark}\label{per:rem:caveats}
Three qualifications are necessary. First,
Corollary~\ref{cor:two-roots} concerns the local beliefs. The finite
quadratic form $\tfrac12x^{\T}(D-z+A)x$ defines a normalizable Gaussian
model only for $z<\lambda_{\min}(D+A)$, where
$\lambda_{\min}(D+A)\ge\sigm$. For a bipartite graph with positive
weights $\lambda_{\min}(D+A)=\sigm$, so that for no $z$ in the open gap
$(\sigm,\Sigm)$ is there a normalizable finite model; the two fixed
points are then fixed points of the covariance equations without an
underlying Gaussian distribution on $G$. Second, we make no statement here on the
exact number or the parity of the positive solutions for a given $z$ in
the gap (see Section~\ref{sec:real}). Third, the
convergence of the belief propagation iteration is a different question:
walk-summability is sufficient for convergence~\cite{MJW}, but the
iteration is known to converge for some models that are not
walk-summable~\cite{MJW}, and the second fixed point of
Corollary~\ref{cor:two-roots} is not the limit of the iteration from the
standard initialization.
\end{remark}

\subsection{The all-negative real branch and the shape of \texorpdfstring{$\Xp$}{X+}}
\label{per:subsec:negative}

The pole branch $P_+$ at $\sigp$ is real for real $z$. On the side
$z<\sigp$ it is the end of the arc $\Gamma_+$. We now identify the side
$z>\sigp$, when $G$ has no leaves, with the branch $S=E$ of
Lemma~\ref{lem:infinity}.

\begin{proposition}[All-negative branch]\label{prop:all-negative-branch}
Let $G$ be a connected loopless multigraph in which every vertex has
degree at least $2$, with real coefficients. Call a real solution of
\eqref{eq:system} \emph{all-negative} if $V_i>0$ for all $i$ and
$f_e=-\tfrac12\bigl(1+\sqrt{1+4w_e^2V_iV_j}\bigr)$ for all $e=ij$.
\begin{enumerate}
\item[\textup{(a)}] For every $z>\sigp$ there is exactly one all-negative
solution. It is real analytic in $z$, and for large $z$ it is the branch
$S=E$ of Lemma~\ref{lem:infinity}.
\item[\textup{(b)}] For $\max_iD_i<z\le\sigp$ there is no all-negative
solution.
\item[\textup{(c)}] As $z\downarrow\sigp$ the all-negative solution is
unbounded. If $\beta\ge2$ it coincides, for $z$ close to $\sigp$, with
the pole branch $P_+$; hence the branch $S=E$ continues along the real
axis to a simple pole at $\sigp$ with
$V_i=\kappa_0(\perup_i)^2/(z-\sigp)+O(1)$.
\end{enumerate}
By the symmetry \eqref{per:eq:iota}, for $z<\sigm$ there is exactly one
real solution with $V<0$ and $f_e=-\tfrac12(1+\sqrt{1+4w_e^2V_iV_j})$; it
is the branch $S=E$ for $z\to-\infty$ and, if $\beta\ge2$, coincides near
$\sigm$ with the pole branch $P_-$ on the side $z<\sigm$.
\end{proposition}

\begin{proof}
Let $z>\max_iD_i$, $t_i=z-D_i>0$ and $F_e(V)=\sqrt{1+4w_e^2V_iV_j}$. The
all-negative solutions are the fixed points with positive entries of
\[
\Psi_z(V)_i=\frac1{t_i}\Bigl(\frac12\sum_{e\ni i}F_e(V)+\frac{d_i}2-1\Bigr),
\]
which is monotone increasing on $[0,\infty)^V$ with
$\Psi_z(0)_i=(d_i-1)/t_i>0$.

\emph{Uniqueness and analyticity.} In the coordinates $x=\log V$ let
$T(x)=\log\Psi_z(e^x)$. Its Jacobian matrix has nonnegative entries, and
since $V_i\partial_{V_i}F_e=V_j\partial_{V_j}F_e=(F_e^2-1)/(2F_e)$ its
$i$th row sum is
\[
\frac{\sum_{e\ni i}\bigl(\tfrac12F_e-\tfrac1{2F_e}\bigr)}{\sum_{e\ni i}\tfrac12F_e+\tfrac{d_i}2-1}<1,
\]
because $d_i\ge2$. On the segment joining two fixed points the row sums
are bounded by a constant less than $1$, so
$\norm{T(x)-T(x')}_\infty<\norm{x-x'}_\infty$ unless $x=x'$: there is at
most one fixed point. At a fixed point $D\Psi_z$ is conjugate to $DT$ by
$\diag(V)$, so its spectral radius is less than $1$, and the implicit
function theorem gives real analytic dependence on $z$.

\emph{Existence for $z>\sigp$.} Let $V=c\,(\perup)^2$ with $c>0$. From
$F_e\le1+2\abs{w_e}\sqrt{V_iV_j}$ and
$\sum_{e\ni i}\abs{w_e}\perup_j=(\sigp-D_i)\perup_i$ we get
$\Psi_z(V)_i\le\bigl(d_i-1+c(\perup_i)^2(\sigp-D_i)\bigr)/t_i$, which is at
most $V_i$ as soon as
$c\ge c(z):=\max_i(d_i-1)/((\perup_i)^2(z-\sigp))$. The iterates of $\Psi_z$
starting from $0$ increase and stay below $c(z)(\perup)^2$, so they
converge to a fixed point $V(z)$ with $\Psi_z(0)\le V(z)\le c(z)(\perup)^2$.

\emph{Nonexistence for $z\le\sigp$.} At an all-negative solution
$F_e>2\abs{w_e}\sqrt{V_iV_j}$ and $d_i\ge2$ give
$t_iV_i>\sum_{e\ni i}\abs{w_e}\sqrt{V_iV_j}$, that is,
$((z-D-\abs A)v)_i>0$ for $v=\sqrt V$. Pairing with $\perup$ gives
$(z-\sigp)(\perup)^{\T}v>0$.

\emph{Identification at infinity.} Since $c(z)=O(1/z)$ we have $V(z)\to0$,
$F_e\to1$ and $f_e\to-1$ for all $e$ as $z\to\infty$. By
Lemma~\ref{lem:infinity} the solution lies on the branch $S=E$ for large
$z$.

\emph{Behavior at $\sigp$.} The map $\Psi_z$ decreases in $z$, hence so do
its iterates from $0$ and $V(z)$. If $V(z)$ stayed bounded as
$z\downarrow\sigp$, its limit would be an all-negative solution at
$z=\sigp$. Let $\beta\ge2$. On the pole branch $P_+$ the analogue of
\eqref{per:eq:pole-data} is $W_i(0)=\kappa_0(\perup_i)^2>0$ and
$y_e(0)=-2\kappa_0\abs{w_e}\perup_i\perup_j<0$, so for small real
$\zeta=z-\sigp>0$ it has $V=W/\zeta>0$ and $1+2f_e=y_e/\zeta<0$: it is an
all-negative solution, and by uniqueness it is $V(z)$.
\end{proof}

\begin{remark}\label{per:rem:full-circle}
Combining Proposition~\ref{prop:all-negative-branch} with
Theorem~\ref{thm:perron-component} we see that, when $\beta\ge2$ and $G$
has no leaves, the physical component $\Cphys$ contains the branch $S=E$
at infinity for all real coefficients: the real-analytic all-negative
family is contained in $\Cphys$ near $\sigp$, hence everywhere. In
Example~\ref{per:ex:regular} the two roots of the quadratic equation at
$z=\infty$ are indeed $f=0$ and $f=-1$.
\end{remark}

Finally we describe $\Xp$ for generic coefficients.

\begin{proposition}[Generic smoothness]\label{prop:generic-smooth}
Let $G$ be a loopless multigraph. For generic coefficients the set
$\cZ$ of admissible solutions is a smooth complex curve. For real
coefficients in this generic set, the real points of $\cZ$ form a smooth
one-dimensional real manifold, and $\Xp$ is an open and closed subset of
it.
\end{proposition}

\begin{proof}
Let $Y$ be the set of $(D,w,z,V,f)$ with $w_e\ne0$ such that $(z,V,f)$ is
an admissible solution of \eqref{eq:system} for $(D,w)$. The derivative
of the $i$th vertex equation with respect to $D_i$ is $V_i\ne0$ and that
of the equation of the edge $e=ij$ with respect to $w_e$ is
$-2w_eV_iV_j\ne0$; each of these coefficients occurs in one equation
only, so the Jacobian matrix of the defining equations has full rank and
$Y$ is a smooth variety of dimension $n+\abs E+1$. By generic smoothness
in characteristic zero, there is a dense Zariski open set of
coefficients over which the projection $Y\to\C^V\times(\C^*)^E$ is a
smooth morphism; its fibers are the sets $\cZ$, which are therefore
smooth of dimension $1$ (or empty). The real points of a smooth complex
curve defined over $\R$ form a smooth real curve. In it $\Xp$ is open,
being defined by strict inequalities, and closed: a limit of positive
solutions has $z\le\Sigm<\min_iD_i$ by Lemma~\ref{lem:domain},
$V_i\ge1/(D_i-z)>0$, and $f_e\ge0$ with $f_e(1+f_e)>0$.
\end{proof}

\begin{corollary}[One arc and finitely many circles]\label{cor:arc-circles}
Let $G$ be a connected loopless multigraph with $\beta\ge2$ and let the
real coefficients be generic in the sense of
Proposition~\ref{prop:generic-smooth}. Then $\Xp$ is a smooth real curve
contained in $\{z\le\Sigm\}$ with finitely many connected components.
Exactly one component is non-compact; it is the arc
$\nu(\Gamma_-)$, which runs from the physical branch at $z\to-\infty$
through the turning point over $\Sigm$ to the pole at $\sigm$, approached
from $z>\sigm$. All other components are diffeomorphic to circles.
\end{corollary}

\begin{proof}
The set $\Xp$ is semialgebraic, so it has finitely many components, each
diffeomorphic to $\R$ or to a circle and closed in $\R^{1+n+\abs E}$. Along
an end of a non-compact component the point leaves every compact set;
since $z\le\Sigm$ and $0<f_e\le\abs{w_e}\sqrt{V_iV_j}$, either $z$ is
unbounded below along the end, or $z$ stays bounded and
$\norm V\to\infty$. In the first case the end meets $\{z<-R'\}$, where
$\Xp$ is the physical branch (Lemma~\ref{per:far}), so it is the end
$z\to-\infty$ of the physical branch. In the second case
Lemma~\ref{lem:blowup}(b) and the uniqueness part of
Theorem~\ref{thm:local-pole} show that the end coincides with the part
$\zeta>0$ of the pole branch $P_-$. Hence the non-compact components
have at most two ends in total. Each has two ends, and the component
containing the physical branch is non-compact, so it is the only one.
By Theorem~\ref{thm:perron-component} it contains $\nu(\Gamma_-)$, which
is connected and approaches both of its ends; so it equals
$\nu(\Gamma_-)$.
\end{proof}

For such coefficients the positive solutions for a fixed $z$ are the
intersection points of the arc and of the circles with the fiber over
$z$. The right-hand side $[z<\sigm]$ of the degree identity of
Theorem~\ref{thm:degree-identity} records the position of the end of the
arc: for $z<\sigm$ a generic fiber meets the arc an odd number of times,
and for $z>\sigm$ an even number of times. Compact components do occur;
see Section~\ref{sec:real}.

%% file: sections/irreducibility.tex
\section{Irreducibility and the degree of the Green function}
\label{sec:irreducibility}

By Theorem~\ref{thm:generic-degree} the spectral curve has $\LG$ sheets
over the $z$-line. In this section we show that for generic coefficients
these sheets form a single irreducible curve, so that the analytic
continuation of the Green function visits all of them, and that no single
diagonal Green function, half-tree $m$-function or the Stieltjes transform
of the density of states lives on a smaller cover. This proves
Theorem~\ref{thm:main-sheets}.

Avni, Breuer and Simon proved that the Green functions of a periodic
Jacobi operator on a tree are algebraic \cite[Theorem~6.6]{ABS}; in the
remark following that theorem they bound the number of sheets by B\'ezout's
theorem by $2^{2q}$, where $q$ is the number of edges of the base graph,
and note that in the examples they could compute there are only two
sheets. Problem~10.2 of \cite{ABS} asks whether the Green function and the
$m$-functions are two-sheeted. For generic coefficients the answer is given by
Theorems~\ref{thm:green-degree} and~\ref{thm:m-function-degree}: the number
of sheets is exactly $\LG$, which is at least $2^\beta\ge4$ as soon as
$\beta\ge2$ (Proposition~\ref{prop:L-bounds}). Two sheets occur for
symmetric coefficients (Proposition~\ref{prop:uniform}).

\subsection{Statements}
\label{irr:subsec:statements}

\begin{theorem}\label{thm:irreducible}
Let $G$ be a finite loopless multigraph. For generic coefficients the
spectral curve $\Xc(G;D,w)$ is irreducible, of degree $\LG$ over the
$z$-line. Equivalently, the monodromy of $z\colon\Xc\to\C$ acts
transitively on the $\LG$ admissible solutions over a general point~$z$.
\end{theorem}

\begin{corollary}\label{cor:all-sheets}
Let $G$ be connected. For generic real coefficients $\Cphys=\Xc$: there is
a finite set $B\subset\C$ such that the analytic continuation of the
physical branch $(\Vphys,f^{\mathrm{phys}})$ along paths in $\C\setminus B$
reaches every one of the $\LG$ admissible solutions over every point of
$\C\setminus B$.
\end{corollary}

The word ``generic'' cannot be removed: for edge-transitive graphs with
uniform coefficients $\Cphys$ has two sheets
(Proposition~\ref{prop:uniform}). The theorem concerns complex admissible
solutions; it says nothing about the number of real or positive solutions
over a real~$z$, and we have no criterion that decides whether a given real
model belongs to the generic set. Transitivity of the monodromy in the
single variable $z$ is stronger than transitivity of the monodromy in all
the coefficients, which is what the degree count of
Theorem~\ref{thm:generic-degree} would give by a standard argument.

The proof of Theorem~\ref{thm:irreducible} is by induction on the number of
edges and occupies
Sections~\ref{irr:subsec:constants}--\ref{irr:subsec:proof}. When the weight
of one edge $e_0$ tends to zero, the curve of $G$ degenerates to the union
of two curves of $G-e_0$, on which $f_{e_0}\to0$ and $f_{e_0}\to-1$
respectively. The second one is the curve of a system in which the constant
$1$ of the vertex equations at the two ends of $e_0$ is replaced by~$0$, so
that the induction has to be carried out for systems with integer vertex
constants. Each of the two curves is irreducible by induction, and the
point is to show that for a small nonzero weight they are joined. We do
this near a point where the variances of a whole connected component of
$G-e_0$ diverge in the direction of an eigenvector of a signed Jacobi
matrix: in suitable coordinates the two families become the two sheets of
an explicit conic bundle over a small annulus, with two simple branch
points.

\subsection{Vertex constants and the passage to generic coefficients}
\label{irr:subsec:constants}

Throughout
Sections~\ref{irr:subsec:constants}--\ref{irr:subsec:proof}, $G=(V,E)$ is a
finite loopless multigraph, not necessarily connected, and
$p=(D,w)\in P_G:=\C^V\times(\C^*)^E$. For $\mu\in\C^V$ we consider the
system
\begin{equation}\label{irr:eq:system-mu}
(D_i-z)V_i=\mu_i+\sum_{e\ni i}f_e\quad(i\in V),\qquad
f_e(1+f_e)=w_e^2V_iV_j\quad(e=ij\in E),
\end{equation}
which is \eqref{eq:system} for $\mu=\ones$. As before, a solution is
admissible if $z\notin\{D_i\}$ and all $V_i\ne0$. For $z\notin\{D_i\}$ we
eliminate $V_i=g_i/(D_i-z)$, where $g_i=\mu_i+\sum_{e\ni i}f_e$, and obtain
the equivalent system
\begin{equation}\label{irr:eq:f-system}
f_e(1+f_e)=c_e\,g_ig_j\quad(e=ij\in E),\qquad
c_e=c_e(p,z)=\frac{w_e^2}{(D_i-z)(D_j-z)},
\end{equation}
in the unknowns $f\in\C^E$; admissibility means $g_i\ne0$ for all $i$. A
solution with $z\notin\{D_i\}$ is \emph{nondegenerate} if the Jacobian
matrix of \eqref{irr:eq:f-system} with respect to $f$ is invertible there.
A nondegenerate solution lies on a unique local analytic branch
$z\mapsto f(z)$ of solutions, hence on a unique irreducible component of
the solution variety of \eqref{irr:eq:system-mu}.

For $\mu\in\Z^V$ put
\[
L_\mu(G)=\#\{S\subseteq E:\ d_S(i)\ne\mu_i\ \text{for all }i\in V\},
\]
so that $L_{\ones}(G)=\LG$, and let $X_\mu(G;p)$ be the union of the
irreducible components of the Zariski closure of the set of admissible
solutions of \eqref{irr:eq:system-mu} on which $z$ is not constant. Thus
$X_{\ones}(G;p)=\Xc(G;D,w)$. We write $\mu\le\ones$ if $\mu_i\le1$ for
all~$i$, and $\ones_e=\ones_i+\ones_j$ for an edge $e=ij$.

\begin{theorem}\label{irr:thm:general}
Let $G$ be a finite loopless multigraph and $\mu\in\Z^V$ with
$\mu\le\ones$. For generic $p\in P_G$ the curve $X_\mu(G;p)$ is
irreducible, or empty if $L_\mu(G)=0$.
\end{theorem}

Theorem~\ref{thm:irreducible} is the case $\mu=\ones$. The restriction
$\mu\le\ones$ is used only in the combinatorial
Lemma~\ref{irr:lem:effective}; the constants that occur in the induction
starting from $\ones$ only decrease.

The next lemma extends Lemma~\ref{lem:infinity} to integer constants.

\begin{lemma}\label{irr:lem:infinity}
Let $\mu\in\Z^V$ and $p\in P_G$. There is $R_p>0$, which can be chosen
locally uniformly in $p$, such that the following holds.
\begin{enumerate}
\item[\textup{(i)}] For $\abs z>R_p$ the solutions of
\eqref{irr:eq:f-system} are $2^{\abs E}$ pairwise distinct points $f^S(z)$,
$S\subseteq E$, which are analytic in $1/z$, satisfy $f^S(\infty)=-\ones_S$,
and are nondegenerate.
\item[\textup{(ii)}] If $d_S(i)\ne\mu_i$ for all $i$, the branch $f^S$ is
admissible for $\abs z>R_p$. Otherwise $g_i$ vanishes identically on $f^S$
for some~$i$.
\item[\textup{(iii)}] The degree of $z\colon X_\mu(G;p)\to\C$ is
$L_\mu(G)$; the branches of $X_\mu(G;p)$ over $\abs z>R_p$ are the
admissible branches of~\textup{(ii)}, and the general point of every
component of $X_\mu(G;p)$ is a nondegenerate admissible solution.
\end{enumerate}
\end{lemma}

\begin{proof}
(i) At $c=0$ the system \eqref{irr:eq:f-system} reads $f_e(1+f_e)=0$, with
the $2^{\abs E}$ solutions $-\ones_S$ and Jacobian matrix
$\diag(1+2f_e)=\diag(\pm1)$. Its leading forms $f_e^2-c_e(\sum_{e'\ni i}f_{e'})(\sum_{e'\ni j}f_{e'})$
have only the trivial common zero at $c=0$, so for $c$ in a neighborhood
of $0$ there are no solutions at infinity and, by B\'ezout's theorem and
the implicit function theorem, the solutions are exactly the $2^{\abs E}$
analytic continuations $f^S(c)$ of the points $-\ones_S$. Since
$c(p,z)\to0$ as $z\to\infty$, locally uniformly in $p$, this gives~(i).

(ii) On $f^S$ we have $g_i\to\mu_i-d_S(i)$ as $z\to\infty$. Suppose
$d_S(i)=\mu_i$. Fix $f_e=-\ones_S(e)$ for the edges at~$i$ and consider
the system \eqref{irr:eq:f-system} of $G-i$ in which the constant of each
neighbor $j$ of $i$ is $\mu_j-\#\{e\in S:e=ij\}$. By~(i) it has a solution
near $-\ones_{S\setminus\{e\ni i\}}$ for small~$c$. Together with the fixed
values this is a solution of the system of $G$, because $g_i=\mu_i-d_S(i)=0$
identically and both sides of the equations of the edges at~$i$ vanish. By
the uniqueness in~(i) it is $f^S$.

(iii) Let $C$ be an irreducible component of the closure of the set of
admissible solutions on which $z$ is not constant. The image $z(C)$ is
cofinite, so $C$ contains points over $\abs z>R_p$, where it is a union of
branches $f^S$. If one of them is not admissible, some $g_i$ vanishes on an
infinite subset of $C$, hence on $C$, and $C$ contains no admissible point,
which is absurd. Conversely an admissible branch is an arc of admissible
solutions along which $z$ is not constant, so it lies in $X_\mu(G;p)$, and
since its points are nondegenerate it lies in only one component. Hence the
$L_\mu(G)$ admissible branches are distributed among the components of
$X_\mu(G;p)$. Nondegeneracy and admissibility are Zariski open conditions
on each component and hold over $\abs z>R_p$.
\end{proof}

We now show that irreducibility at one coefficient implies irreducibility
at generic coefficients. Let $\mathcal U\subseteq\C^E$ be the set of
$c$ for which \eqref{irr:eq:f-system} has $2^{\abs E}$ distinct solutions
in $\C^E$; it is the set of $c$ off the zero locus of the resultant of the
leading forms at which no solution is multiple, hence Zariski open as in
the proof of Lemma~\ref{deg:lem-open}, and for
$c\in\mathcal U$ all solutions are nondegenerate. By the proof of
Lemma~\ref{irr:lem:infinity} it contains a neighborhood of~$0$.

\begin{lemma}\label{irr:lem:generic}
Let $\mu\in\Z^V$.
\begin{enumerate}
\item[\textup{(i)}] There is a Zariski open set
$\mathcal U_{\mathrm{good}}\subseteq\mathcal U$ containing a neighborhood
of~$0$ such that for $c\in\mathcal U_{\mathrm{good}}$ exactly $L_\mu(G)$ of
the $2^{\abs E}$ solutions of \eqref{irr:eq:f-system} are admissible.
Consequently the set
$\Omega=\{(p,z):z\notin\{D_i\},\ c(p,z)\in\mathcal U_{\mathrm{good}}\}$ is
Zariski open in $P_G\times\C$ and $\Omega_p=\{z:(p,z)\in\Omega\}$ is
cofinite for every~$p$.
\item[\textup{(ii)}] For $z\in\Omega_p$ the points of $X_\mu(G;p)$ over $z$
are exactly the $L_\mu(G)$ admissible solutions over~$z$. Hence
$X_\mu(G;p)$ is irreducible if and only if the unramified covering
$X_\mu(G;p)|_{\Omega_p}\to\Omega_p$ of degree $L_\mu(G)$ is connected, that
is, if and only if its monodromy is transitive.
\item[\textup{(iii)}] The set
$I=\{p\in P_G:X_\mu(G;p)\text{ is irreducible}\}$ is open in the Euclidean
topology.
\item[\textup{(iv)}] If $I\ne\emptyset$, then $I$ contains a nonempty
Zariski open subset of $P_G$.
\end{enumerate}
\end{lemma}

\begin{proof}
(i) The solutions of \eqref{irr:eq:f-system} over $\mathcal U$ form a
finite \'etale covering $\mathcal S\to\mathcal U$ of degree $2^{\abs E}$.
The property that $\prod_ig_i$ does not vanish identically near a point of
$\mathcal S$ is preserved by analytic continuation, so the points with this
property form a union $\mathcal S'$ of connected components of $\mathcal S$.
Near $c=0$ the sheets of $\mathcal S'$ are the branches $f^S(c)$ with
$d_S(i)\ne\mu_i$ for all $i$, by the proof of
Lemma~\ref{irr:lem:infinity}\,(ii), so $\mathcal S'\to\mathcal U$ has degree
$L_\mu(G)$. The norm of $\prod_ig_i$ along $\mathcal S'\to\mathcal U$ is a
regular function on $\mathcal U$ that does not vanish near~$0$; let
$\mathcal U_{\mathrm{good}}$ be the complement of its zero locus. For
$c\in\mathcal U_{\mathrm{good}}$ the admissible solutions are the
$L_\mu(G)$ points of $\mathcal S'$ over~$c$. The remaining assertions follow
because $c(p,z)$ is regular on $\{z\notin\{D_i\}\}$ and tends to $0$ as
$z\to\infty$.

(ii) An admissible solution over $z\in\Omega_p$ is nondegenerate, hence
lies on an analytic branch of admissible solutions over a neighborhood
of~$z$, and therefore on $X_\mu(G;p)$. Conversely, let $P\in X_\mu(G;p)$ lie
over $z\in\Omega_p$. Near $z$ all solutions lie on $2^{\abs E}$ disjoint
analytic branches, and $P$ is a limit of admissible points of
$X_\mu(G;p)$, which lie on the branch through~$P$; this branch therefore
belongs to $\mathcal S'$, and $P$ is admissible because
$c(p,z)\in\mathcal U_{\mathrm{good}}$. An irreducible curve remains
connected after removal of finitely many fibers, and two distinct
components of $X_\mu(G;p)$ are disjoint over $\Omega_p$; this proves the
equivalence.

(iii) Let $p_1\in I$ and fix $z_0\in\Omega_{p_1}$ and finitely many loops
that generate $\pi_1(\Omega_{p_1},z_0)$. Since $\Omega$ is open and the
loops are compact, they lie in $\Omega_p$ for $p$ near $p_1$, and the
permutations of the fiber over $z_0$ that they induce do not depend on $p$
by continuity. These permutations generate a transitive group, so
$X_\mu(G;p)|_{\Omega_p}$ is connected.

(iv) The complement $\Delta$ of $\Omega$ in $P_G\times\C$ is Zariski closed
with finite fibers over $P_G$. Its irreducible components that dominate
$P_G$ form a hypersurface $\{h(p,z)=0\}$ with $h$ squarefree in~$z$, and the
other components project into a proper closed subset $Z'\subset P_G$. Let
$P_0$ be the set of $p\notin Z'$ at which the leading coefficient and the
discriminant of $h$ with respect to $z$ do not vanish. On $P_0$ the finite
set $\C\setminus\Omega_p$ consists of a constant number of points that
depend continuously on $p$. Hence, for $p$ near $p_1\in P_0$, loops around
small discs centered at these points generate $\pi_1(\Omega_p)$ and induce
the same permutations as for $p_1$, as in~(iii); transitivity of the
monodromy is therefore locally constant on $P_0$. As $P_0$ is connected in
the Euclidean topology, transitivity holds on all of $P_0$ or nowhere
on~$P_0$. By~(iii) the set $I$ contains a nonempty open set, which meets
the dense set $P_0$; thus $P_0\subseteq I$.
\end{proof}

\subsection{Separation of signed spectra}
\label{irr:subsec:separation}

The genericity conditions on eigenvalues that we use here and in
Sections~\ref{sec:classification} and~\ref{sec:genus} all follow from one
statement. For a nonempty $U\subseteq V$ with $G[U]$ connected and a
signing $s$ of the edges $E(U)$ of $G[U]$, let
$p_{U,s}(z)=\det\bigl(z-(D+A^s)[U]\bigr)$. It depends only on the
switching class $[s]$ of $s$ on $G[U]$.

\begin{lemma}\label{irr:lem:separation}
Regard $D_i$ and $w_e$ as indeterminates. The polynomials $p_{U,s}$ are
irreducible in $\C(D,w)[z]$, and $p_{U,s}\ne p_{U',s'}$ unless $U=U'$ and
$[s]=[s']$. Consequently, for generic coefficients, all eigenvalues of all
matrices $(D+A^s)[U]$, $G[U]$ connected, are simple, the spectra of
$(D+A^s)[U]$ and $(D+A^{s'})[U']$ are disjoint unless $U=U'$ and
$[s]=[s']$, every eigenvector $u$ of $(D+A^s)[U]$ has no zero entry and
satisfies $u^{\T}u\ne0$, and no eigenvalue of $(D+A^s)[U]$ equals $D_x$
unless $U=\{x\}$.
\end{lemma}

\begin{proof}
Put $x_i=z-D_i$ and $\mathbb K=\C(w)$. The polynomial
$P=\det(\diag(x_U)-A^s[U])\in\mathbb K[x_i:i\in U]$ has degree one in each
variable, with leading monomial $\prod_{i\in U}x_i$. Suppose $P=QR$ with
nonconstant $Q,R$. Comparing degrees in each variable, $Q$ and $R$ involve
disjoint sets of variables $I$ and $J$ with $I\sqcup J=U$, and we may take
both with leading coefficient one. Since the diagonal of the matrix is
$(x_i)$ with no constant term, the coefficient of $\prod_{k\ne i}x_k$ in $P$
is zero for every $i$; hence the coefficient of $\prod_{k\in I\setminus i}x_k$
in $Q$ vanishes for $i\in I$, and similarly for $R$. For $i\in I$, $j\in J$
the coefficient of $\prod_{k\ne i,j}x_k$ in $P$ is the product of these two
coefficients, so it is zero. On the other hand it equals $-(A^s_{ij})^2$,
which is nonzero in $\mathbb K$ when $i$ and $j$ are adjacent. As $G[U]$ is
connected, some $i\in I$ is adjacent to some $j\in J$, a contradiction.
Thus $P$ is irreducible in $\mathbb K[x_U]$, hence in
$\mathbb K[x_U][z]=\mathbb K[D_U][z]$, and, being monic in $z$, in
$\C(D,w)[z]$ by Gauss's lemma.

If $U\ne U'$ the two polynomials involve different variables $D_i$. Let
$U=U'$ and $[s]\ne[s']$. Then some cycle $C$ of $G[U]$ has different sign
products for $s$ and $s'$. If $C$ has length $\ell\ge3$, the monomial
$\prod_{e\in C}w_e\prod_{i\in U\setminus V(C)}x_i$ arises in $P$ only from
the two cyclic permutations along $C$ and has coefficient
$-2\prod_{e\in C}s_e$; if $C$ consists of two parallel edges $e,e'$ joining
$i$ and $j$, the monomial $w_ew_{e'}\prod_{k\ne i,j}x_k$ has coefficient
$-2s_es_{e'}$. In both cases $p_{U,s}\ne p_{U,s'}$.

Distinct monic irreducible polynomials over a field of characteristic zero
have nonzero discriminants and pairwise nonzero resultants, which gives the
statements on eigenvalues for generic coefficients. If an eigenvector $u$
of $(D+A^s)[U]$ for the eigenvalue $\lambda$ vanishes at $x\in U$, its
restriction to some connected component $U'$ of $G[U\setminus x]$ is a
nonzero eigenvector of $(D+A^s)[U']$ for $\lambda$, which is excluded. For
a simple eigenvalue of a complex symmetric matrix the left eigenvector is
$u^{\T}$, so $u^{\T}u\ne0$. Finally $D_x$ is the eigenvalue
of~$(D+A^s)[\{x\}]$.
\end{proof}

\subsection{Degeneration of one edge}
\label{irr:subsec:degeneration}

Let $e_0=ij\in E$ and $H=G-e_0$ (only the edge $e_0$ is removed, not its
parallel edges). We write $p=(p',\varepsilon)$ with $p'\in P_H$ and
$w_{e_0}=\varepsilon$, and $f_0=f_{e_0}$. Splitting the sets
$S\subseteq E$ according to whether $e_0\in S$ gives
\begin{equation}\label{irr:eq:split}
L_\mu(G)=L_\mu(H)+L_{\mu-\ones_{e_0}}(H).
\end{equation}
Put
\[
A=X_\mu(H;p'),\qquad B=X_{\mu-\ones_{e_0}}(H;p'),
\]
and let $\Omega^A_{p'}$, $\Omega^B_{p'}$ be the cofinite sets of
Lemma~\ref{irr:lem:generic} for these two systems on~$H$.

\begin{lemma}[Degeneration]\label{irr:lem:degeneration}
Let $\gamma\subset\Omega^A_{p'}\cap\Omega^B_{p'}$ be compact. There is
$\varepsilon_0>0$ such that every admissible solution $a$ of the system of
$A$ or of $B$ over a point $z\in\gamma$ extends to a family
$a(z',\varepsilon)$ of solutions of the system
\eqref{irr:eq:system-mu} of $G$ at the coefficient $(p',\varepsilon)$, which
is analytic for $z'$ near $z$ and $\abs\varepsilon<\varepsilon_0$, with
$f_0=0$, respectively $f_0=-1$, at $\varepsilon=0$. For
$0<\abs\varepsilon<\varepsilon_0$ these solutions are nondegenerate and
admissible, they are points of $X_\mu(G;p',\varepsilon)$, and they are the
only solutions of the system of $G$ in a neighborhood of $(a,0)$,
respectively $(a,-1)$, which is independent of $\varepsilon$.
\end{lemma}

\begin{proof}
In \eqref{irr:eq:f-system} for $G$ we have
$c_{e_0}=\varepsilon^2/((D_i-z)(D_j-z))$, uniformly $O(\varepsilon^2)$ near
$\gamma$. At $\varepsilon=0$ the equation of $e_0$ is $f_0(1+f_0)=0$ and
the other equations form the system of $H$ with constants
$\mu+f_0\ones_{e_0}$; thus $f_0=0$ gives the solutions of the system of $A$
and $f_0=-1$ those of~$B$. The Jacobian matrix at $\varepsilon=0$ is block
triangular: the row of $e_0$ has the single entry $1+2f_0=\pm1$, and the
complementary block is the Jacobian matrix of the system of $A$ or $B$,
which is invertible over $\Omega^A_{p'}\cap\Omega^B_{p'}$. The implicit
function theorem and the compactness of $\gamma$ give the family and its
uniqueness. The $g_i$ are close to their nonzero values at $\varepsilon=0$,
so the solutions are admissible; being nondegenerate, they lie on branches
of admissible solutions over the $z$-line, hence on
$X_\mu(G;p',\varepsilon)$.
\end{proof}

\subsection{Coordinates near a diverging eigenvector}
\label{irr:subsec:eta}

The following coordinates describe the solutions whose variances diverge on
a whole connected component in the direction of a signed eigenvector. They
are used again in Sections~\ref{sec:classification} and~\ref{sec:genus}.

Let $K$ be a finite loopless multigraph with coefficients, $U$ the vertex
set of one of its connected components, with edge set $E(U)$ and degrees
$d_x$, and $Y=V(K)\setminus U$. Let $s$ be a signing of $E(U)$ and
$\lambda$ a simple eigenvalue of $M_s=(D+A^s)[U]$ with an eigenvector $u$
such that $u^{\T}u=1$ and $u_x\ne0$ for all $x\in U$. Put
\[
y^0_e=-2s_ew_eu_xu_{x'}\qquad(e=xx'\in E(U)).
\]
Let $\mu\in\C^{V(K)}$. We assume that $\lambda\notin\{D_x:x\in Y\}$ and
that $\psi(z)$ is an analytic family, defined for $z$ near $\lambda$, of
nondegenerate admissible solutions of the system \eqref{irr:eq:system-mu} of
$K[Y]$ with constants $\mu_Y$ (no condition if $Y=\emptyset$). The
constants $\mu_U$ will be treated as parameters.

\begin{lemma}[$\eta$-coordinates]\label{irr:lem:eta}
Consider the system
\begin{equation}\label{irr:eq:eta-system}
\begin{aligned}
&(D_x-\lambda-\zeta)\tilde V_x-\eta\Bigl(\mu_x-\frac{d_x}2\Bigr)-\frac12\sum_{e\ni x}y_e=0\quad(x\in U),\\
&y_e^2-\eta^2-4w_e^2\tilde V_x\tilde V_{x'}=0\quad(e=xx'\in E(U)),\qquad
\sum_{x\in U}\tilde V_x=1
\end{aligned}
\end{equation}
in the unknowns $\tilde V\in\C^U$, $y\in\C^{E(U)}$, $\zeta\in\C$, with
parameters $\eta\in\C$ and $\mu_U\in\C^U$.
\begin{enumerate}
\item[\textup{(a)}] For every compact $\Theta\subset\C^U$ there are
$\delta_0>0$, a neighborhood $\Theta'$ of $\Theta$ and an analytic
solution $(\tilde V,y,\zeta)(\eta,\mu_U)$ on
$\{\abs\eta<\delta_0\}\times\Theta'$ with value $(u^2,y^0,0)$ at $\eta=0$,
where $u^2=(u_x^2)_{x\in U}$. It is the only solution in a neighborhood of
$(u^2,y^0,0)$ that does not depend on $(\eta,\mu_U)$.
\item[\textup{(b)}] For $\eta\ne0$, the formulas
\[
z=\lambda+\zeta,\qquad V_x=\frac{\tilde V_x}\eta\ (x\in U),\qquad
f_e=\frac12\Bigl(\frac{y_e}\eta-1\Bigr)\ (e\in E(U)),
\]
together with $\psi(z)$ on $K[Y]$, define a solution of the system
\eqref{irr:eq:system-mu} of $K$, which is admissible if
$z\notin\{D_x:x\in U\}$. Conversely, every solution of the system of $K$
for which $\eta:=(\sum_{x\in U}V_x)^{-1}$ is small,
$(\eta V_U,\eta(1+2f_{E(U)}),z-\lambda)$ is close to $(u^2,y^0,0)$ and the
part on $K[Y]$ is close to $\psi(\lambda)$, arises in this way from the
solution of~\textup{(a)}.
\item[\textup{(c)}] $\zeta=\kappa\,\eta+O(\eta^2)$ with
$\kappa=\kappa(\mu_U)=\abs{E(U)}-\sum_{x\in U}\mu_x$.
\item[\textup{(d)}] If $\kappa(\mu_U)=0$, then $\zeta=C\eta^2+O(\eta^3)$
with
\begin{equation}\label{irr:eq:C}
C=\sum_{e=xx'\in E(U)}s_e\Bigl[\frac1{4w_eu_xu_{x'}}-\frac{w_eu_xu_{x'}}4(\omega_x-\omega_{x'})^2\Bigr],
\end{equation}
where $\omega\in\C^U$ is any solution of
$(M_s-\lambda)(u\circ\omega)=(2\mu_U-d)/u$, with entrywise product and
quotient.
\end{enumerate}
\end{lemma}

For $\mu=\ones$ and $\kappa\ne0$ the family of~(a) is the pole of
Theorem~\ref{thm:local-pole}, $V_U\approx\kappa u^2/(z-\lambda)$,
parametrized by~$\eta$. For $\kappa=0$ and $C\ne0$ it is a band edge:
$z-\lambda\approx C\eta^2$ and $V_U\approx\pm u^2\sqrt{C/(z-\lambda)}$, as at
the edges of the bands of a cycle.

\begin{proof}
(b) The edge equation $f(1+f)=w^2V_xV_{x'}$ is equivalent to
$(1+2f)^2=1+4w^2V_xV_{x'}$, which becomes the second equation of
\eqref{irr:eq:eta-system} after multiplication by $\eta^2$. Multiplying the
vertex equation by $\eta$ and inserting
$\sum_{e\ni x}f_e=\frac12\sum_{e\ni x}(y_e/\eta-1)$ gives the first, and the
normalization $\sum\tilde V_x=1$ is the definition of~$\eta$. The component
$U$ shares no edge with $Y$, so the two parts of the system of $K$ are
coupled only through~$z$. The converse follows from the uniqueness in~(a)
and the nondegeneracy of~$\psi$.

(a) At $\eta=0$ the point $(u^2,y^0,0)$ is a solution: the first equation
is $u_x\bigl((M_s-\lambda)u\bigr)_x=0$ and the second is an identity. Since
$\partial_{y_e}$ of the second equation is $2y^0_e\ne0$, we can eliminate
$y$ to first order,
$dy_e=-s_ew_e\bigl(\tfrac{u_{x'}}{u_x}d\tilde V_x+\tfrac{u_x}{u_{x'}}d\tilde V_{x'}\bigr)$,
and, using $\sum_{e=xx'}s_ew_eu_{x'}=-(D_x-\lambda)u_x$, the Jacobian matrix
of the first equations with respect to $\tilde V$ becomes
\[
\mathcal J=\tfrac12\diag(u)\,(M_s-\lambda)\,\diag(u)^{-1},
\]
with kernel spanned by $u^2$ and left kernel spanned by $\ones$, because
$\lambda$ is simple and $M_s$ is symmetric. The Jacobian matrix with
respect to $(\tilde V,\zeta)$ of the first and third equations is
$\bigl(\begin{smallmatrix}\mathcal J&-u^2\\ \ones^{\T}&0\end{smallmatrix}\bigr)$.
If $\mathcal Jv=u^2\tau$ and $\ones^{\T}v=0$, then applying $\ones^{\T}$
gives $\tau\,u^{\T}u=0$, so $\tau=0$, $v\in\C u^2$ and $v=0$. The matrix is
invertible, and since $\mu_U$ enters only through $\eta\mu_U$ the implicit
function theorem applies uniformly on~$\Theta'$.

(c) Differentiating at $\eta=0$ gives
$\mathcal J\tilde V'-u^2\zeta'=\mu_U-d/2$; applying $\ones^{\T}$ yields
$\zeta'=\abs{E(U)}-\sum_{x\in U}\mu_x$.

(d) Summing the first equations of \eqref{irr:eq:eta-system} and using
$\sum\tilde V_x=1$ and
$y_e=-2s_ew_e\bigl(\tilde V_x\tilde V_{x'}+\eta^2/(4w_e^2)\bigr)^{1/2}$ (the
branch close to $y^0$) we get the exact identity
\[
\zeta=Q(\tilde V)+\sum_{e\in E(U)}2s_ew_e\Bigl[\bigl(\tilde V_x\tilde V_{x'}+\tfrac{\eta^2}{4w_e^2}\bigr)^{1/2}-(\tilde V_x\tilde V_{x'})^{1/2}\Bigr]+\kappa\,\eta,
\]
where
$Q(\tilde V)=\sum_x(D_x-\lambda)\tilde V_x+\sum_e2s_ew_e(\tilde V_x\tilde V_{x'})^{1/2}$
is homogeneous of degree one with $Q(u^2)=0$ and $\nabla Q(u^2)=0$. Write
$\tilde V=u^2+\eta\,u^2\circ\omega+O(\eta^2)$; the first-order equation
of~(c) with $\zeta'=\kappa=0$ is $(M_s-\lambda)(u\circ\omega)=(2\mu_U-d)/u$,
which is solvable because the right-hand side is orthogonal to~$u$, and
determines $\omega$ up to an additive constant. Then
$Q(\tilde V)=\tfrac{\eta^2}2(u^2\circ\omega)^{\T}\nabla^2Q(u^2)(u^2\circ\omega)+O(\eta^3)$,
and the second derivatives of $2s_ew_e(\tilde V_x\tilde V_{x'})^{1/2}$ at
$u^2$, namely $-\tfrac{s_ew_e}2u_{x'}u_x^{-3}$ and $\tfrac{s_ew_e}{2u_xu_{x'}}$,
give the contribution $-\tfrac{s_ew_e}4u_xu_{x'}(\omega_x-\omega_{x'})^2\eta^2$
of the edge~$e$. The bracket contributes $\eta^2s_e/(4w_eu_xu_{x'})+O(\eta^3)$.
\end{proof}

\begin{lemma}\label{irr:lem:band-edge}
Let $U$ be the vertex set of a connected multigraph that contains a cycle
\textup{(}possibly of length two\textup{)}, $\mu\in\Z^U$ with
$\abs{E(U)}=\sum_x\mu_x$, and $s$ a signing. For generic coefficients on
$U$ the number $C$ of \eqref{irr:eq:C} is nonzero at every eigenvalue of
$(D+A^s)[U]$.
\end{lemma}

\begin{proof}
Let $k=\abs U\ge2$. The variety
$\Lambda=\{(D,w,\lambda,u):(M_s-\lambda)u=0,\ u^{\T}u=1,\ u_x\ne0\}$ is
parametrized by $(w,\lambda,u)$ through $D_x=\lambda-(A^su)_x/u_x$; it is
isomorphic to an open subset of the product of $(\C^*)^{E(U)}\times\C$ with
the irreducible quadric $\{u^{\T}u=1\}$, hence irreducible, of the same
dimension as the space of coefficients, onto which it projects with finite
fibers over the generic set of Lemma~\ref{irr:lem:separation}. On the open
subset of $\Lambda$ where $\lambda$ is simple, $C$ is a regular function.
It suffices to find one point where $C\ne0$: its zero set is then a proper
closed subset, whose image in the space of coefficients is not dense.

Take $u=k^{-1/2}\ones$, $w_e=s_e\sigma_e$ with $\sigma_e>0$, and
$D_x=\lambda-\sum_{e\ni x}\sigma_e$. Then $M_s-\lambda=-L_\sigma$, the
Laplacian of the connected weighted graph, and $\lambda$ is simple. With
$b=d-2\mu$, which satisfies $\ones^{\T}b=0$, we get
$\omega=kL_\sigma^+b$ and
\[
C=\frac k4\Bigl[\sum_{e\in E(U)}\frac1{\sigma_e}-b^{\T}L_\sigma^+b\Bigr].
\]
Let $e_1$ be an edge on a cycle and let $\sigma_{e_1}\to0$ with the other
weights fixed. The graph without $e_1$ is connected, so $b^{\T}L_\sigma^+b$
stays bounded, while $1/\sigma_{e_1}\to\infty$. Hence $C\ne0$ for small
$\sigma_{e_1}$.
\end{proof}

\subsection{Joining the two families}
\label{irr:subsec:joining}

We return to the setting of Section~\ref{irr:subsec:degeneration}. In the
next two lemmas $U$ is the vertex set of a connected component of a
subgraph $K$ of $H$, the data $(s,\lambda,u,\psi)$ are as in
Section~\ref{irr:subsec:eta} for $K$ and the coefficient~$p'$, and
$(\tilde V,y,\zeta)(\eta,\mu_U)$ is the solution of
Lemma~\ref{irr:lem:eta}\,(a). For constants $\nu$ on $U$ with
$\zeta(\cdot,\nu)\not\equiv0$ and small $\eta\ne0$ we denote by $P_\nu(\eta)$
the solution of the system of $K$ described in
Lemma~\ref{irr:lem:eta}\,(b); it is admissible for all but finitely many
$\eta$ in a small disc, and lies over $z=\lambda+\zeta(\eta,\nu)$.

\begin{lemma}[Non-bridge edge]\label{irr:lem:stage-N}
Let $e_0=ij$ be an edge of $G$ that is not a bridge, $K=H=G-e_0$, and $U$
the vertex set of the component of $H$ containing $i$ and $j$. Suppose that
\[
\zeta_A(\eta):=\zeta(\eta,\mu_U)\qquad\text{and}\qquad
\zeta_B(\eta):=\zeta(\eta,\mu_U-\ones_{e_0})
\]
do not vanish identically. Let $\eta_1\ne0$ be small and such that
$P_A:=P_{\mu_U}(\eta_1)$ and $P_B:=P_{\mu_U-\ones_{e_0}}(\eta_1)$ are
admissible and lie over points of $\Omega^A_{p'}\cap\Omega^B_{p'}$. Then
for all sufficiently small $\varepsilon\ne0$ the extensions of $P_A$ and of
$P_B$ given by Lemma~\ref{irr:lem:degeneration} lie on the same irreducible
component of $X_\mu(G;p',\varepsilon)$.
\end{lemma}

Here $P_A$ is a point of $A$ and $P_B$ a point of~$B$. By
Lemma~\ref{irr:lem:eta}\,(c) the hypothesis holds if
$\kappa_A:=\abs{E_H(U)}-\sum_{x\in U}\mu_x$ and $\kappa_B=\kappa_A+2$ are
both nonzero; if one of them vanishes it holds when the corresponding
number $C$ of \eqref{irr:eq:C} is nonzero.

\begin{proof}
\emph{Step 1: reduction to a plane curve.} The system of $G$ consists of
the system of $H$ with the constants $\mu+f_0\ones_{e_0}$, in which $f_0$
is regarded as a parameter, together with the equation
$f_0(1+f_0)=\varepsilon^2V_iV_j$ of~$e_0$. Put $a=u_i^2u_j^2\ne0$ and apply
Lemma~\ref{irr:lem:eta} with
$\Theta=\{\mu_U+\tau\ones_{e_0}:\abs\tau\le3\}$. Since
$V_iV_j=\tilde V_i\tilde V_j/\eta^2$, part~(b) of that lemma shows that the
solutions of the system of $G$ in the range of the $\eta$-coordinates
correspond bijectively to the solutions $(\eta,f_0)$, $\eta\ne0$,
$\abs{f_0}<3$, of
\begin{equation}\label{irr:eq:plane-curve}
f_0(1+f_0)\,\eta^2=\varepsilon^2\mathcal A(\eta,f_0),\qquad
\mathcal A(\eta,f_0)=(\tilde V_i\tilde V_j)(\eta,\mu_U+f_0\ones_{e_0})=a+O(\eta).
\end{equation}
Choose $\delta_1<\delta_0$ such that $\abs{\mathcal A-a}\le\abs a$ for
$\abs\eta\le\delta_1$, $\abs{f_0}\le3$, and let
\[
R_\varepsilon=\bigl\{\sqrt{\abs a}\,\abs\varepsilon<\abs\eta<\delta_1\bigr\},\qquad
\Sigma_\varepsilon=\bigl\{(\eta,f_0):\eta\in R_\varepsilon,\ \abs{f_0}<3,\ \eqref{irr:eq:plane-curve}\text{ holds}\bigr\}.
\]

\emph{Step 2: $\Sigma_\varepsilon$ is a connected Riemann surface.} For
$\eta\in R_\varepsilon$ the right-hand side of
$f_0(1+f_0)=(\varepsilon/\eta)^2\mathcal A$ is bounded by $2$ in absolute
value, while $\abs{f_0(1+f_0)}\ge6$ on $\abs{f_0}=3$. By Rouch\'e's theorem
there are exactly two roots $f_0$ in $\abs{f_0}<3$, and their sum
$\sigma_1$ and product $\sigma_2$ are holomorphic functions of
$\eta\in R_\varepsilon$, given by contour integrals over $\abs{f_0}=3$.
Comparing these integrals with those for the constant $a$ in place of
$\mathcal A$ and writing $\eta=\varepsilon\eta'$ we find
$\sigma_1=-1+O(\abs\varepsilon/\abs{\eta'})$,
$\sigma_2=-a/\eta'^2+O(\abs\varepsilon/\abs{\eta'})$, uniformly on
$R_\varepsilon$, so that the discriminant satisfies
\[
\Delta:=\sigma_1^2-4\sigma_2=1+\frac{4a}{\eta'^2}+O\Bigl(\frac{\abs\varepsilon}{\abs{\eta'}}\Bigr).
\]
The function $1+4a/\eta'^2$ has the two simple zeros $\pm2i\sqrt a$, of
modulus $2\sqrt{\abs a}$, is at least $3$ in modulus on
$\abs{\eta'}=\sqrt{\abs a}$, and differs from $1$ by at most $\frac14$ for
$\abs{\eta'}\ge4\sqrt{\abs a}$. By Rouch\'e's theorem, for small
$\varepsilon$ the function $\Delta$ has exactly two zeros in
$R_\varepsilon$, both simple. Now $\Sigma_\varepsilon$ is isomorphic to the
curve $\{(\eta,v):v^2=\Delta(\eta)\}$ over $R_\varepsilon$ through
$v=2f_0-\sigma_1$. This curve is smooth, because the zeros of $\Delta$ are
simple, and connected, because its two sheets are exchanged along a small
loop around a simple zero of~$\Delta$.

\emph{Step 3: the image lies on one component.} By
Lemma~\ref{irr:lem:eta}\,(b), $\Sigma_\varepsilon$ is mapped
biholomorphically onto an open subset $W$ of the solution variety of the
system of $G$; thus this variety is a locally irreducible curve at each
point of $W$, and since $W$ is connected, $W$ is contained in a single
irreducible component~$\Gamma$. On $\Sigma_\varepsilon$ the function
$z=\lambda+\zeta(\eta,\mu_U+f_0\ones_{e_0})$ is not constant for small
$\varepsilon$: on the sheet with $f_0=O(\varepsilon^2)$ over
$\delta_1/2<\abs\eta<\delta_1$ it is uniformly close to
$\lambda+\zeta_A(\eta)$, which is not constant. On $W$ we have
$V_x\ne0$ for $x\in U$, the part on $Y$ is admissible, and
$z\in\{D_x:x\in U\}$ only at finitely many points. Hence $\Gamma$ contains
admissible points and $z$ is not constant on it:
$\Gamma\subseteq X_\mu(G;p',\varepsilon)$.

\emph{Step 4: the two ends.} Over $\eta_1$ the two points of
$\Sigma_\varepsilon$ have $f_0=O(\varepsilon^2)$
and $f_0=-1+O(\varepsilon^2)$. By the continuity of
the solution of Lemma~\ref{irr:lem:eta}\,(a) in $\mu_U$, the corresponding
solutions of the system of $G$ converge, as $\varepsilon\to0$, to $(P_A,0)$
and $(P_B,-1)$. By Lemma~\ref{irr:lem:degeneration}, for small
$\varepsilon$ the first of them therefore belongs to the family
$z'\mapsto a(z',\varepsilon)$ extending $P_A$; this family is an arc of
nondegenerate solutions, so it lies on the unique component through that
point, which is $\Gamma$. In particular the extension of $P_A$ itself lies
on $\Gamma$. The same holds for~$P_B$.
\end{proof}

\begin{lemma}[Pendant edge]\label{irr:lem:stage-P}
Let $e_0=ij$ where $j$ is a vertex of degree one of $G$, $K=G-j$, and $U$
the vertex set of the component of $K$ containing~$i$. Suppose that
$\lambda\ne D_j$, that $\mu_j\notin\{0,1\}$, and that
$\kappa_A:=\abs{E(U)}-\sum_{x\in U}\mu_x\notin\{0,-1\}$. For small
$\eta_1\ne0$ let $P_A$ be the solution $P_{\mu_U}(\eta_1)$ of the system
of $K$ completed by $V_j=\mu_j/(D_j-z)$, and $P_B$ the solution
$P_{\mu_U-\ones_i}(\eta_1)$ completed by $V_j=(\mu_j-1)/(D_j-z)$; assume
that they are admissible and lie over points of
$\Omega^A_{p'}\cap\Omega^B_{p'}$. Then the conclusion of
Lemma~\ref{irr:lem:stage-N} holds.
\end{lemma}

\begin{proof}
In $H=G-e_0$ the vertex $j$ is isolated. The system of $G$ consists of the
system of $K$ with the constants $\mu_U+f_0\ones_i$ on $U$, the relation
$V_j=(\mu_j+f_0)/(D_j-z)$ and the equation of $e_0$, which in the
$\eta$-coordinates reads
\[
f_0(1+f_0)=q\,\tilde b(\eta,f_0)\,(\mu_j+f_0),\qquad q=\frac{\varepsilon^2}\eta,\quad
\tilde b=\frac{\tilde V_i}{D_j-\lambda-\zeta}=b+O(\eta),\quad b=\frac{u_i^2}{D_j-\lambda}\ne0.
\]
For $\tilde b$ replaced by $b$ this is a quadratic equation in $f_0$ with
discriminant $\Delta_0(q)=b^2q^2+2b(2\mu_j-1)q+1$, whose roots $q_\pm$ have
product $b^{-2}$ and are distinct precisely when $(2\mu_j-1)^2\ne1$, that
is, when $\mu_j\notin\{0,1\}$. Let $\rho=2\max\abs{q_\pm}$, choose $F$ with
$F(F-1)>2\rho\abs b(\abs{\mu_j}+F)$, apply Lemma~\ref{irr:lem:eta} with
$\Theta=\{\mu_U+\tau\ones_i:\abs\tau\le F\}$, choose $\delta_1$ such that
$\abs{\tilde b-b}\le\abs b$, hence $\abs{\tilde b}\le2\abs b$, for
$\abs\eta\le\delta_1$ and $\abs{f_0}\le F$, and replace $R_\varepsilon$ by
$\{\abs\varepsilon^2/\rho<\abs\eta<\delta_1\}$, on which $\abs q<\rho$. As
in Step~2 of the proof of Lemma~\ref{irr:lem:stage-N}, for each such $\eta$
there are two roots in $\abs{f_0}<F$, and their discriminant is
$\Delta_0(q)+O(\varepsilon^2)$ because $q(\tilde b-b)=O(\varepsilon^2)$; it
has exactly two zeros, both simple, near $q=q_\pm$. Hence the solutions
form a connected Riemann surface $\Sigma_\varepsilon$. Near the outer
boundary $q\to0$, and the two sheets have $f_0\approx0$, which
corresponds to the constants $\mu$ of $A$, and $f_0\approx-1$, which
corresponds to the constants $\mu-\ones_{e_0}$ of~$B$. By
Lemma~\ref{irr:lem:eta}\,(c) the functions $\zeta(\cdot,\mu_U)$ and
$\zeta(\cdot,\mu_U-\ones_i)$ have the nonzero linear coefficients $\kappa_A$
and $\kappa_A+1$. On $\Sigma_\varepsilon$ the function $f_0$ is not
constant, so $V_j=0$, that is $f_0=-\mu_j$, holds at finitely many points
only. Steps~3 and~4 of the proof of Lemma~\ref{irr:lem:stage-N} now apply
without change.
\end{proof}

\subsection{Existence of an effective edge}
\label{irr:subsec:effective}

\begin{lemma}\label{irr:lem:effective}
Let $\mu\in\Z^V$, $\mu\le\ones$, and $L_\mu(G)\ge2$. Then $G$ has an edge
$e_0=ij$ of one of the following two kinds.
\begin{enumerate}
\item[\textup{(N)}] $e_0$ is not a bridge and, with $U$ the vertex set of
the component of $G-e_0$ containing $i$ and $j$ and
$\kappa_A=\abs{E_{G-e_0}(U)}-\sum_{x\in U}\mu_x$, either
$\kappa_A\notin\{0,-2\}$, or $\kappa_A=0$ and $G-e_0$ contains a cycle
in~$U$.
\item[\textup{(P)}] $j$ has degree one and either $\mu_j\in\{0,1\}$, or,
with $U$ the vertex set of the component of $G-j$ containing $i$ and
$\kappa_A=\abs{E(U)}-\sum_{x\in U}\mu_x$, we have $\kappa_A\notin\{0,-1\}$.
\end{enumerate}
\end{lemma}

In case (P) with $\mu_j\in\{0,1\}$ one of the two curves $A$, $B$ is empty:
for $\mu_j=0$ the isolated vertex $j$ of $H$ has $g_j=0$ in the system
of~$A$, and for $\mu_j=1$ its constant in the system of $B$ is~$0$.

\begin{proof}
We write $\beta(K)$ for the cycle rank of a connected graph~$K$.

Suppose first that some component $K$ of $G$ has $\beta(K)\ge2$, and let
$e_0$ be an edge of $K$ on a cycle. Then $U=V(K)$ and $K-e_0$ is connected
of cycle rank at least one, so
$\abs{E_{G-e_0}(U)}\ge\abs U\ge\sum_{x\in U}\mu_x$ and $\kappa_A\ge0$. If
$\kappa_A=0$, the graph $K-e_0$ still contains a cycle. Thus $e_0$ is of
kind~(N).

Suppose now that all components have cycle rank at most one. Since
$L_\mu(G)=\prod_KL_\mu(K)\ge2$, some component $K$ has $L_\mu(K)\ge2$; it is
not an isolated vertex. If $K$ has a vertex $j$ of degree one with
$\mu_j\in\{0,1\}$, its edge is of kind~(P).

Assume that $K$ has vertices of degree one and that $\mu_j\le-1$ for all of
them. If $K$ is a single edge $ij$, then $U=\{i\}$ and
$\kappa_A=-\mu_i\ge1$. Otherwise let $j$ be a vertex of degree one and
$U=V(K)\setminus\{j\}$, which is connected with
$\abs{E(U)}=\abs U-1+\beta(K)$, so that
$\kappa_A=\sum_{x\in U}(1-\mu_x)+\beta(K)-1$. If $U$ contains another
vertex $j'$ of degree one of $K$, then $1-\mu_{j'}\ge2$ and $\kappa_A\ge1$.
If not, $K$ is not a tree, $\beta(K)=1$ and
$\kappa_A=\sum_{x\in U}(1-\mu_x)\ge0$; when $\kappa_A\ge1$ the edge at $j$
is of kind~(P), and when $\kappa_A=0$, that is $\mu\equiv1$ on $U$, we take
instead an edge $e_0'$ on the cycle of $K$, for which $U'=V(K)$ and
$\kappa_A'=(\abs{V(K)}-1)-(\abs{V(K)}-1+\mu_j)=-\mu_j\ge1$, of kind~(N).

Finally assume that $K$ has no vertex of degree one. Then $K$ is a cycle
$C_\ell$, $\ell\ge2$, and for every edge $U=V(K)$ and
$\kappa_A=\ell-1-\sum_x\mu_x=\sum_x(1-\mu_x)-1\ge-1$. The value
$\kappa_A=0$ occurs only if $\mu_x=0$ at one vertex and $\mu_x=1$ at all
others. But then $L_\mu(C_\ell)=1$: the empty set is excluded by the vertex
with $\mu_x=0$, and a set $S$ with $\emptyset\ne S\ne E$ is a union of
paths and has at least two vertices of degree one, at most one of which has
$\mu_x\ne1$. This contradicts $L_\mu(K)\ge2$, so $\kappa_A\notin\{0,-2\}$
and every edge is of kind~(N).
\end{proof}

\subsection{Proof of Theorems~\ref{irr:thm:general} and~\ref{thm:irreducible}}
\label{irr:subsec:proof}

\begin{proof}[Proof of Theorem~\ref{irr:thm:general}]
We argue by induction on $\abs E$, for all pairs $(G,\mu)$ with
$\mu\le\ones$ simultaneously. If $L_\mu(G)\le1$ the curve has degree at
most one over the $z$-line by Lemma~\ref{irr:lem:infinity} and is
irreducible or empty; this covers $E=\emptyset$. Let $L_\mu(G)\ge2$ and let
$e_0=ij$ be an edge as in Lemma~\ref{irr:lem:effective}, $H=G-e_0$. Both
$\mu$ and $\mu-\ones_{e_0}$ are $\le\ones$, so by induction there is a
nonempty Zariski open set $\mathcal W\subseteq P_H$ on which $A$ and $B$
are irreducible or empty. Let $L_A=L_\mu(H)$ and
$L_B=L_{\mu-\ones_{e_0}}(H)$; by \eqref{irr:eq:split}, $L_\mu(G)=L_A+L_B$.
By Lemma~\ref{irr:lem:generic}\,(iv) it suffices to find one coefficient
$p=(p',\varepsilon)$ at which $X_\mu(G;p)$ is irreducible.

\emph{Choice of $p'$.} If $L_A=0$ or $L_B=0$ take any $p'\in\mathcal W$.
Otherwise let $U$ be as in Lemma~\ref{irr:lem:effective}, let $s$ be the
signing $s\equiv1$ of the edges of $H$ inside $U$, and let $Y$ be the set
of the remaining vertices, not counting $j$ in case~(P). The coefficients
$p'$ consist of the coefficients $p_U$ on $U$, the coefficients $p_Y$ on
$H[Y]$, and $D_j$ in case~(P). The projection of $\mathcal W$ to the space
of the $p_U$ contains a nonempty Zariski open set; we choose $p_U$ in it
such that the conclusions of Lemma~\ref{irr:lem:separation} hold for $U$
and, if $\kappa_A=0$ in case~(N), such that $C\ne0$
(Lemma~\ref{irr:lem:band-edge}, which applies because then $U$ contains a
cycle). Fix an eigenvalue $\lambda$ of $(D+A^s)[U]$ with eigenvector $u$,
$u^{\T}u=1$. For fixed $\lambda$ the map
$p_Y\mapsto c(p_Y,\lambda)\in\C^{E(Y)}$ is dominant, so the set of $p_Y$
with $\lambda\in\Omega^Y_{p_Y}$, the set of
Lemma~\ref{irr:lem:generic} for the system of $H[Y]$ with constants
$\mu_Y$, is nonempty and Zariski open; we choose $p_Y$ (and $D_j\ne\lambda$)
in it such that $p'\in\mathcal W$. Since $L_A>0$ and $L_A$ is the product of
the numbers $L_\mu$ of the components of $H$, the system of $H[Y]$ has
admissible solutions, which over $\Omega^Y_{p_Y}$ are nondegenerate; let
$\psi$ be one of the corresponding analytic families near $z=\lambda$.

\emph{All points of $A$ over $z_0$ lie on one component.} Fix
$z_0\in\Omega^A_{p'}\cap\Omega^B_{p'}$. The part $A^\circ$ of $A$ over
$\Omega^A_{p'}\cap\Omega^B_{p'}$ consists of nondegenerate admissible
solutions and is connected, since $A$ is irreducible. Two points of
$A^\circ$ are joined by a path in $A^\circ$, whose projection to the
$z$-line is compact; by Lemma~\ref{irr:lem:degeneration} the path extends,
for small $\varepsilon\ne0$, to a path of nondegenerate admissible
solutions of the system of $G$, which stays on one irreducible component of
$X_\mu(G;p',\varepsilon)$. Hence the extensions of all points of $A^\circ$
lying on finitely many prescribed paths are on one component $\Gamma_A$
when $\varepsilon$ is small enough, and similarly for $B$ with a
component~$\Gamma_B$. We apply this to the $L_A$ points of $A$ and the
$L_B$ points of $B$ over $z_0$ and, if $L_A,L_B\ge1$, to the points $P_A$
and $P_B$ below.

\emph{The two families are joined.} Let $L_A,L_B\ge1$. In case~(N) of
Lemma~\ref{irr:lem:effective} we have $\kappa_B=\kappa_A+2\ne0$, and either
$\kappa_A\ne0$ or $\kappa_A=0$ and $C\ne0$; thus $\zeta_A$ and $\zeta_B$
do not vanish identically, and Lemma~\ref{irr:lem:stage-N} applies with a
suitable small $\eta_1$. In case~(P) we have $\mu_j\notin\{0,1\}$ because
neither $A$ nor $B$ is empty, hence $\kappa_A\notin\{0,-1\}$, and
Lemma~\ref{irr:lem:stage-P} applies. In both cases the extensions of $P_A$
and $P_B$ lie on one component, so $\Gamma_A=\Gamma_B$.

\emph{Conclusion.} In all cases there is, for small $\varepsilon\ne0$, an
irreducible component $\Gamma$ of $X_\mu(G;p',\varepsilon)$ that contains
the extensions of the $L_A+L_B$ solutions of $A$ and $B$ over $z_0$. These
extensions are pairwise distinct, so $\Gamma$ has degree at least
$L_\mu(G)$ over the $z$-line, which is the degree of
$X_\mu(G;p',\varepsilon)$. Hence $X_\mu(G;p',\varepsilon)=\Gamma$.
\end{proof}

\begin{proof}[Proof of Theorem~\ref{thm:irreducible} and
Corollary~\ref{cor:all-sheets}]
Theorem~\ref{thm:irreducible} is Theorem~\ref{irr:thm:general} for
$\mu=\ones$, together with Lemma~\ref{irr:lem:generic}\,(ii) for the
statement on the monodromy. The generic set contains a dense open set of
real coefficients. For such coefficients and connected $G$, the physical
branch is the branch $S=\emptyset$ at infinity, because
$f^{\mathrm{phys}}_e(z)\to0$ as $z\to-\infty$; it lies on the irreducible
curve $\Xc$, so $\Cphys=\Xc$, and the corollary follows from the
transitivity of the monodromy with $B=\C\setminus\Omega_p$.
\end{proof}

\subsection{The degree of a single Green function}
\label{irr:subsec:green}

Irreducibility of $\Xc$ shows that the vector $(\Vphys,f^{\mathrm{phys}})$
has $\LG$ branches. A single coordinate could still have fewer: for a
disconnected graph $G_1\sqcup G_2$ the function $V_i$, $i\in G_1$, has only
$L(G_1)$ branches although $\Xc$ is irreducible of degree $L(G_1)L(G_2)$.
For connected graphs this does not happen.

\begin{theorem}\label{thm:green-degree}
Let $G$ be connected with $\beta\ge2$. For generic real coefficients and
every $i\in V$, the Green function $G_i$ is an algebraic function of degree
exactly $\LG$ over $\C(z)$: its analytic continuation takes $\LG$ distinct
values over a general point~$z$. In particular it has at least
$2^\beta\ge4$ branches.
\end{theorem}

\begin{proof}
We take the coefficients in the generic sets of
Theorem~\ref{thm:irreducible} and Lemma~\ref{irr:lem:separation}. Let
$L=\LG$ and let $\mathbb F=\C(\Xc)$, of degree $L$ over $\C(z)$. If the
minimal polynomial of $V_i\in\mathbb F$ over $\C(z)$ has degree $d$, its
characteristic polynomial is the $(L/d)$-th power of the minimal
polynomial; hence over a general $z$ the $L$ sheets fall into $d$ groups of
$L/d$ sheets on which $V_i$ defines the same analytic function.

By Theorem~\ref{thm:local-pole}, applied to $\lambda=\sigp$,
$s=\sgn w$ and the Perron vector $\bar u^+$ (the eigenvalue $\sigp$ is
simple, $\bar u^+>0$, $\sigp>\max_iD_i$ and $\kappa_0>0$), there is a point
$P_+$ of the normalization of $\Xc$ over $\sigp$ which is unramified and at
which $V_i$ has a simple pole. Suppose $L/d\ge2$. Then for real $z_1>\sigp$
close to $\sigp$ there is a sheet over $z_1$, different from the sheet
through $P_+$, on which $V_i$ coincides with its values on the sheet
through $P_+$. Following this sheet along the segment from $z_1$ to $\sigp$,
which contains no critical value of $z$ if $z_1$ is close to $\sigp$, we
arrive at a point $Q$ of the normalization over $\sigp$ at which $V_i$ is
unbounded, and $Q\ne P_+$ because $P_+$ is unramified.

We use two lemmas from Section~\ref{sec:classification}. Neither uses
this theorem. By Lemma~\ref{cls:lem:support}, the vertices
at which the pole order at $Q$ is maximal form a set $U$ with $G[U]$
connected, $\sigp$ is an eigenvalue of $(D+A^s)[U]$ for some signing~$s$
with eigenvector $u$, and with $\eta=(\sum_{x\in U}V_x)^{-1}$ we have
$(\eta V_U,\eta(1+2f_{E(U)}))\to(u^2,y^0)$ at~$Q$. By
Lemma~\ref{lem:perron-isolated} this forces $U=V$ and $s$ switching
equivalent to $\sgn w$, so $u=\pm\bar u^+$ after switching. Then
Lemma~\ref{irr:lem:eta}\,(b), applied with $K=G$, $U=V$ and $\mu=\ones$,
shows that the points of $\Xc$ near $Q$ belong to the family of
Lemma~\ref{irr:lem:eta}\,(a), on which
$z-\sigp=\kappa_0\eta+O(\eta^2)$ with $\kappa_0\ne0$. This family is a
single unramified branch, which is the branch through $P_+$. Hence $Q=P_+$,
a contradiction, and $d=L$. The inequality $\LG\ge2^\beta$ is
Proposition~\ref{prop:L-bounds}.
\end{proof}

\subsection{Half-tree functions and the density of states}
\label{irr:subsec:m-functions}

The preceding proof uses the Perron pole to distinguish the sheets. A
different argument, at $z=\infty$, applies to all real coefficients and to
other spectral functions. Recall from \eqref{eq:m-from-Vf} that
$m_{i|e}=V_i/(1+f_e)$ is a rational function on $\Xc$ which on the
physical branch is the half-tree $m$-function. The half-tree in question is
the component, containing the lift of $i$, of the universal cover with one
lift of $e$ removed; we call it the half-tree of $(i,e)$. The Stieltjes
transform of the density of states of $\tH$ is $\frac1n\sum_iG_i$.

\begin{theorem}\label{thm:m-function-degree}
Let $G$ be connected, with arbitrary real coefficients, and let $\mathbb F$
be the function field of the physical component $\Cphys$. Each of the
following functions $q$ generates $\mathbb F$ over $\C(z)$, that is,
$\mathbb F=\C(z,q)$:
\begin{enumerate}
\item[\textup{(a)}] $q=m_{i|e}$ for an edge $e=ij$ such that the half-tree
of $(i,e)$ is infinite;
\item[\textup{(b)}] $q=\sum_ic_iV_i$ with real $c_i\ge0$, not all zero; in
particular a single Green function $G_i$ and the Stieltjes transform
$\frac1n\sum_iG_i$ of the density of states.
\end{enumerate}
Consequently each of these functions is algebraic over $\C(z)$ of degree
exactly $\deg(z\colon\Cphys\to\C)$, and for generic real coefficients this
degree is $\LG$. If the half-tree of $(i,e)$ is finite, $m_{i|e}$ is a
rational function of~$z$.
\end{theorem}

In particular the answer to Problem~10.2 of \cite{ABS} is negative for the
$m$-functions as well: for $\beta\ge2$ and generic coefficients each
$m$-function of an infinite half-tree has $\LG\ge4$ branches. For
part~(b) with a single vertex this gives a second proof of
Theorem~\ref{thm:green-degree}, valid also for $\beta\le1$.

The proof rests on a comparison of the branches at infinity. Let
$t=-1/z$. By Lemma~\ref{lem:infinity} the points of $\Xc$ over $z=\infty$
are unramified and indexed by $S\in\cL(G)$, and
\begin{equation}\label{irr:eq:infinity-expansion}
f^S_e=-\ones_S(e)+O(t^2),\qquad V^S_i=(1-d_S(i))\,t+O(t^2).
\end{equation}
The physical branch is $S=\emptyset$. Let $V(S)$ be the set of vertices
with $d_S(i)>0$. A walk is \emph{non-backtracking} if it never traverses an
edge and then immediately the same edge in the opposite direction.

\begin{lemma}\label{irr:lem:germs}
Let the coefficients be real, $S\in\cL(G)$ nonempty, and let $q$ be $V_i$,
or $m_{i|e}$ with $e=ij\notin S$. Let $r$ be the least length of a
non-backtracking walk from $i$ to a vertex of $V(S)$, whose first step is
not $e$ in the case $q=m_{i|e}$. If $r<\infty$, then
\begin{equation}\label{irr:eq:germ-difference}
q^S-q^\emptyset=-t^{2r+1}\sum_{\gamma}d_S(v(\gamma))\prod_{e'\in\gamma}w_{e'}^2+O(t^{2r+2}),
\end{equation}
where the sum is over these walks $\gamma$ of length $r$, $v(\gamma)$ is
the endpoint of $\gamma$, and the product is over the steps of $\gamma$.
In particular $q^S\ne q^\emptyset$ as germs at $t=0$. If $e=ij\in S$, then
$m^S_{i|e}$ has a pole at $t=0$, whereas $m^\emptyset_{i|e}=t+O(t^2)$.
\end{lemma}

\begin{proof}
If $e=ij\in S$, then by \eqref{eq:m-from-Vf} and
\eqref{irr:eq:infinity-expansion},
$m^S_{i|e}=f^S_e/(w_e^2V^S_j)=\bigl(w_e^2(d_S(j)-1)t\bigr)^{-1}+O(1)$, and
$d_S(j)\ge2$.

On every admissible solution the vertex and edge equations give
$f_{e'}/V_i=w_{e'}^2V_{j'}/(1+f_{e'})=w_{e'}^2m_{j'|e'}$ for $e'=ij'$, hence
the recursion \eqref{eq:m-recursion}:
\[
\frac1{V_i}=D_i-z-\sum_{e'\ni i}w_{e'}^2m_{j'|e'},\qquad
\frac1{m_{i|e}}=\frac{1+f_e}{V_i}=D_i-z-\sum_{e'\ni i,\ e'\ne e}w_{e'}^2m_{j'|e'}.
\]
Subtracting these identities for the branches $S$ and $\emptyset$ we get
\begin{equation}\label{irr:eq:difference-recursion}
q^S-q^\emptyset=q^S\,q^\emptyset\sum_{e'}w_{e'}^2\bigl(m^S_{j'|e'}-m^\emptyset_{j'|e'}\bigr),
\end{equation}
where $e'$ runs over the edges at $i$, except $e$ if $q=m_{i|e}$.

We prove by induction on $k\ge0$: if $r\ge k$ then
$q^S-q^\emptyset=O(t^{2k+1})$, and if $r=k$ then
\eqref{irr:eq:germ-difference} holds. For $k=0$ the first statement is
clear, and if $r=0$, that is $i\in V(S)$, then by
\eqref{irr:eq:infinity-expansion} and $1+f^S_e=1+O(t^2)$ both for $q=V_i$
and for $q=m_{i|e}$, $e\notin S$, we have
$q^S-q^\emptyset=-d_S(i)t+O(t^2)$. Let $r\ge k\ge1$. Then $i\notin V(S)$,
so $q^Sq^\emptyset=t^2+O(t^3)$, and every edge $e'=ij'$ at $i$ is outside
$S$, so the induction hypothesis applies to $m_{j'|e'}$, for which the
least length is at least $k-1$. Inserting it into
\eqref{irr:eq:difference-recursion} gives both statements for~$k$.
All terms of the sum in \eqref{irr:eq:germ-difference} are positive, so the
leading coefficient is negative.
\end{proof}

\begin{proof}[Proof of Theorem~\ref{thm:m-function-degree}]
The field $\mathbb F$ is generated over $\C(z)$ by the $V_i$ and $f_e$, and
all functions $q$ in the statement belong to $\mathbb F$. Let $P_\emptyset$
be the point at infinity of $\Cphys$ on the physical branch.

We claim that the germ of $q$ at $P_\emptyset$, as a function of $t$,
differs from its germ at every other point $P_S$, $S\ne\emptyset$, of
$\Cphys$ over $z=\infty$. In case~(b) every vertex is joined to $V(S)$ by a
path, so Lemma~\ref{irr:lem:germs} applies to each $V_i$; in
$\sum_ic_i(V^S_i-V^\emptyset_i)$ the terms of lowest order have negative
coefficients and cannot cancel. In case~(a), if $e\in S$ the claim is the
last statement of Lemma~\ref{irr:lem:germs}. If $e\notin S$ we must show
that $V(S)$ is reached by a non-backtracking walk from $i$ whose first step
is not~$e$. The half-tree of $(i,e)$ is infinite, so it contains an
infinite path from its root, which projects to an infinite
non-backtracking walk in $G$ starting at $i$ with first step different
from~$e$. This walk traverses some oriented edge twice, hence it contains a
closed non-backtracking walk that can be repeated indefinitely; let $W$ be
its vertex set. Every vertex of $W$ is reached. A vertex outside $W$ is
reached by first walking to the initial vertex of a shortest path from $W$
to it and then following this path: its first edge leaves $W$, so it is not
the reverse of the preceding step. Thus all vertices are reached, in
particular those of $V(S)$.

Let $\mathbb F'=\C(z,q)\subseteq\mathbb F$ and consider the corresponding
maps of normalized projective curves $\Cphys\to C'\to\PP^1_z$. Let $Q$ be
the image of $P_\emptyset$. Every preimage of $Q$ lies over $z=\infty$, so
it is one of the unramified points $P_S$, and $t$ is a local coordinate at
$P_S$ and at~$Q$. If $P_S$ maps to $Q$, the germ of $q$ at $P_S$ in the
coordinate $t$ is the germ of $q$ at $Q$, hence equals the germ at
$P_\emptyset$, and $S=\emptyset$ by the claim. So the fiber of
$\Cphys\to C'$ over $Q$ is the single unramified point $P_\emptyset$, and
the map has degree one: $\mathbb F'=\mathbb F$.

The statement on the degree follows, with Theorem~\ref{thm:irreducible} and
Corollary~\ref{cor:all-sheets} for generic real coefficients. If the
half-tree of $(i,e)$ is finite, the recursion \eqref{eq:m-recursion} on it
terminates at its leaves and expresses $m_{i|e}$ as a diagonal entry of the
inverse of a finite matrix $H_T-z$.
\end{proof}

For special real coefficients the theorem still says that no information is
lost in passing from the whole physical branch to one $m$-function of an
infinite half-tree, or to the density of states: each of them determines
all Green functions $G_i$ and all $m$-functions rationally.

\subsection{Uniform coefficients}
\label{irr:subsec:uniform}

The examples of \cite{CSZ} and the examples computed in \cite{ABS} have
two sheets. For the following class this is a consequence of symmetry.

\begin{proposition}\label{prop:uniform}
Let $G$ be connected with all degrees at least two and real coefficients, and
suppose that a group $\Gamma$ of automorphisms of $G$ acts transitively on
the edges \textup{(}parallel edges forming one orbit\textup{)}. Let $w_e=w$
for all $e$ and let $D$ be constant on the $\Gamma$-orbits of vertices. Then every edge
joins a vertex of degree $a$ and potential $D_a$ with a vertex of degree
$b$ and potential $D_b$ \textup{(}possibly $a=b$ and $D_a=D_b$\textup{)},
and the physical component $\Cphys$ is the curve of the symmetric
solutions $f_e\equiv f$, $V_i=(1+d_if)/(D_i-z)$, where
\begin{equation}\label{irr:eq:uniform}
(1-c\,ab)\,f^2+\bigl(1-c\,(a+b)\bigr)f-c=0,\qquad
c=c(z)=\frac{w^2}{(D_a-z)(D_b-z)}.
\end{equation}
It has exactly two sheets over the $z$-line, so every $G_i$ is an algebraic
function of degree two. The points of the normalization of $\Cphys$ over
finite $z$ at which some $V_i$ is unbounded are the following:
\begin{enumerate}
\item[\textup{(i)}] one unramified point over each of the two roots of
$(D_a-z)(D_b-z)=ab\,w^2$ if $\beta\ge2$. These roots are $\sigm$ and
$\sigp$; there all $V_i$ have simple poles on the nonphysical branch.
If $\beta=1$, then $a=b=2$, and these two roots give ramification points
at which all $V_i$ diverge as $(z-\sigma_\pm)^{-1/2}$ and $V_i\,dz$ has
zero residue;
\item[\textup{(ii)}] if $a\ne b$, one unramified point over $D_a$ at which
the $V_i$ with $d_i=a$ have a simple pole with residue $(a-b)/b$ and the
others vanish, and one over $D_b$ with $a$ and $b$ exchanged;
\item[\textup{(iii)}] if $a=b$ and $D_a\ne D_b$, the two ramification points
over $D_a$ and $D_b$, at which $V_i\sim\mathrm{const}\cdot(D_i-z)^{-1/2}$ for
the vertices of the corresponding orbit; there the differential $V_i\,dz$
has no residue.
\end{enumerate}
\end{proposition}

\begin{proof}
An edge-transitive group on a connected graph with edges has at most two
orbits on the vertices, and if
there are two, every edge joins them; applied to $\Gamma$, this gives the
first statement. An element of $\Gamma$ preserves the coefficients, lifts to the universal
cover and commutes with $\tH$, so $G_i$ is constant on $\Gamma$-orbits of vertices and
$f^{\mathrm{phys}}_e$, which is symmetric in the two ends of $e$, is
constant on edge orbits. Thus $f^{\mathrm{phys}}_e\equiv f$, the quantities
$g_i$ of \eqref{irr:eq:f-system} are $1+af$ and $1+bf$, and
\eqref{irr:eq:f-system} becomes $f(1+f)=c(1+af)(1+bf)$, which is
\eqref{irr:eq:uniform}.

The discriminant of \eqref{irr:eq:uniform} is
$\Delta(c)=1+(4-2a-2b)c+(a-b)^2c^2$, with $\Delta(0)=1$. If $a\ne b$ it has
two roots $c_1\ne c_2$, because its discriminant is $16(a-1)(b-1)\ne0$; the
equation $c(z)=c_k$, that is $(D_a-z)(D_b-z)=w^2/c_k$, has two simple roots
for at least one $k$, at which $\Delta(c(z))$ vanishes to first order. If
$a=b$, $\Delta$ has the single root $c_1=1/(4(a-1))>0$, and
$(D_a-z)(D_b-z)=w^2/c_1>0$ has two simple roots, since the critical value
$-(D_a-D_b)^2/4$ of the left-hand side is not positive. In both cases
$\Delta(c(z))$ is not a square in $\C(z)$, so \eqref{irr:eq:uniform}
defines an irreducible curve with two sheets. It consists of solutions of
\eqref{eq:system} and contains the physical branch; since two irreducible
curves with infinitely many common points coincide, it is~$\Cphys$.

A function $V_i=(1+d_if)/(D_i-z)$ can be unbounded only where $f$ is
unbounded or $z=D_i$. Write $Q=(D_a-z)(D_b-z)/w^2$. The quadratic is
$(Q-ab)f^2+(Q-a-b)f-1=0$, so $f$ is unbounded exactly over $Q=ab$.
These roots in $z$ are simple because $ab w^2>0$. If $\beta\ge2$ then
$ab-a-b\ne0$, and precisely one root of the quadratic has a simple pole
over each of these values. The roots of
$(D_a-z)(D_b-z)=ab\,w^2$ are the eigenvalues of the quotient matrix of
$D\pm\abs A$ for the partition into orbits, hence $\sigm$ and $\sigp$,
whose eigenvectors are constant on orbits. When $\beta\ge2$ the values
lie outside $[\Sigm,\Sigp]$ by Corollary~\ref{cor:gap}, so the physical
branch is analytic there. When $\beta=1$, connectedness and the minimum
degree imply $a=b=2$, giving $f(1+f)=1/(Q-4)$. A local coordinate at
either of the two points then satisfies $z-\sigma_\pm=t^2$ and
$V_i=O(t^{-1})$, so $V_i\,dz$ is holomorphic. This proves (i).
At $z=D_a$ the two roots of
\eqref{irr:eq:uniform} tend to $-1/a$ and $-1/b$. If $a\ne b$, on the sheet
with $f\to-1/a$ the equation $f(1+f)=c(1+af)(1+bf)$ gives $1+af=O(Q)$, so
$V_i$ is bounded for $d_i=a$, while on the sheet with $f\to-1/b$ we get
$V_i=(1-a/b)/(D_a-z)+O(1)$ for $d_i=a$ and $1+bf=O(Q)$, hence
$V_j=O(D_a-z)$ on the other orbit even when $D_a=D_b$. This proves (ii);
when $D_a=D_b$, the two poles are at distinct points above the same
energy. If
$a=b$, then $(1+af)^2=f(1+f)/c$; when $D_a=D_b$ the function $1/c$ has a
double zero and $V_i$ is bounded, and when $D_a\ne D_b$ it has a simple
zero, which gives~(iii).
\end{proof}

The three examples of \cite{CSZ} are of this form: the regular tree
(Example~3.1, $a=b$, $D_a=D_b$, with the branch points
$D\pm2\sqrt{a-1}\,\abs w$), the trees covering complete bipartite graphs
(Example~3.2, $a\ne b$) and the triple edge with an alternating potential
(Example~3.3, $a=b=3$, $D_a\ne D_b$). In case~(ii) the pole over the
potential of the orbit of smaller degree has negative residue; when
$D_a=D_b$ it lies on the physical sheet, because on the other sheet the
functions $V_i$ of the orbit of larger degree have a pole with positive
residue at the same point, which is incompatible with the Herglotz
property. This is the eigenvalue of $\tH$ noted in \cite[Example~3.2]{CSZ};
the pole with positive residue is antibound by
Corollary~\ref{cor:antibound-sign}. By Theorem~\ref{thm:all-signed-poles}
below, all the other $n2^\beta-2$ signed eigenvalues become antibound poles
of the Green function after a generic perturbation of the coefficients,
and the number of sheets jumps from $2$ to $\LG$. The algebraic degree of
the Green function is therefore not a continuous function of the
coefficients.

\subsection{Remarks}
\label{irr:subsec:remarks}

The monodromy group is constant up to conjugacy on a nonempty Zariski
open set of coefficients. Indeed, the proof of
Lemma~\ref{irr:lem:generic}\,(iv) transports the generating loops and all
their permutations on that open set. More generally, a finite set of
loops on which all $L$ admissible sheets are nondegenerate realizes the
same permutations after a small generic perturbation. This follows by
the implicit function theorem along the compact lifted loops. Extra
exceptional values away from those sheets can be avoided by deforming
the loops in their neighborhoods.

\begin{proposition}\label{irr:prop:small-monodromy}
For generic coefficients the monodromy groups of the dipoles $D_3,D_4,D_5$
are $S_5,S_{12},S_{27}$, respectively. Those of the bowtie and $K_4$ are
$S_4$ and $S_{15}$. The group for $K_4$ with one edge doubled is $S_{38}$,
and those of the prism and $K_{3,3}$ are $S_{50}$.
\end{proposition}

\begin{proof}
\emph{Dipoles.}
For $D_k$ with both vertex constants $m\le1$, keep one edge of weight
$1$, set the others to zero, and set both diagonals to zero. Each
choice of $j$ frozen edges with $f_e=-1$ leaves the equation
\[
 f(1+f)=c(\mu+f)^2,\qquad c=z^{-2},\quad\mu=m-j .
\]
Its discriminant is $1+4\mu(\mu-1)c$. For $k\ge3$, at the negative zero
of the discriminant for $\mu=m-k+1$ all other quadratics are regular:
on the negative integers the products $\mu(\mu-1)$ are distinct and
positive, and they vanish for $\mu=0,1$. The extremal choice freezes every deleted
edge at $-1$, so there is just one double root, and it is admissible.
For $\mu=0$ or $1$ the excluded root is constant. All remaining
admissible sheets are therefore regular there. Since $z\mapsto c$ is
locally invertible and the frozen-edge derivatives are $\pm1$, the
generic group contains a transposition.

All these generic groups are transitive by
Theorem~\ref{irr:thm:general}. A transitive group of degree $n$ with a
prime divisor of its order larger than every proper divisor of $n$
is primitive: in a nontrivial system of blocks, an element of that
prime order would fix every block and every letter. A primitive group
containing a transposition is symmetric. Indeed, the graph whose edges
are the conjugates of that transposition has its components permuted by
the group, so it is connected, and its edge transpositions generate
the symmetric group.

The degree for $D_3$, $m=1$, is the prime $5$, giving $S_5$.
Deleting just one edge from $D_4$ gives, for $m=1$ and $m=0$
respectively, pairs of degrees $(5,7)$ and $(7,8)$ on $D_3$ with constants
$m$ and $m-1$. In each case the loops of the simultaneous cover project
onto a transitive group on the seven-sheet family, so $7$ divides their
order. The permutations persist at a small nonzero edge weight.
This gives primitivity, hence symmetric groups, in degrees $12$ and
$15$. For $D_5$, $m=1$, the two limiting degrees are $12$ and $15$;
the $S_{15}$ on the second family gives a prime divisor $13$ of the
order of the simultaneous group. This gives $S_{27}$.

\smallskip
\emph{Bowtie.}
Give the four private vertices diagonal zero, give each edge opposite
the hub weight $1$, and take the squared weights of the two spokes of
triangle $j$ to be $\gamma_j$. Write $s$ for the hub variance and $v_j$
for the common private variance in triangle $j$. Put
$d=-z$, $h=d^2-1$, $x=s/h$. If $P_j=1+2f_{\rm spoke}$ and
$Q_j=1+2f_{\rm opposite}$, the private vertex equation is
$P_j+Q_j=2dv_j$. Subtracting the squared edge equations gives
$hv_j=dP_j-\gamma_j s$. Thus with $y_j=P_j-2\gamma_j d x$ the symmetric
equations are
\[
 y_j^2=1+4\gamma_j^2x^2,\qquad
 t=\frac{y_1+y_2-1}{x},\qquad
 t=(D_{\rm hub}-z)h-2(\gamma_1+\gamma_2)d .
\]
These equations recover $s,v_j,P_j,Q_j$ wherever $d,h,x,v_1,v_2$ are
nonzero. At such points the antisymmetric linearization in each pair
of private vertices is invertible: after eliminating the private
variance difference, the factor on the spoke difference is
$P_j-\gamma_j s/d=hv_j/d$. Hence regular sheets of the displayed
symmetric system are nondegenerate in the full system.

Set $2\gamma_1=1$, $4\gamma_2^2=1-\epsilon^3$ with
$2\gamma_2\to1$. With $r=1/x$, $a=(y_1-y_2)/x$ and
$b=(y_1+y_2)/x$, the equations are
\[
 ab=\epsilon^3,\qquad r^2=(a+b)^2/4-1,\qquad t=b-r .
\]
At $\epsilon=0$, the two components meet at $a=b=0$, $r=i$, over
$t_0=-i$. With $a=\epsilon A$, $b=\epsilon^2/A$ and $r$ near $i$,
\[
 \frac{t-t_0}{\epsilon^2}\longrightarrow
 \frac{1}{A}-\frac{A^2}{8i} .
\]
The degree-three rational function on the right has transitive
monodromy. Its three finite critical points satisfy $A^3=-4i$, are
simple, and have distinct values $3/(2A)$. Thus its group is $S_3$.
Finitely many generating loops in the scaled $t$-plane lift to three
simple sheets with finite nonzero $A$. The expansion, holomorphic
near their compact lifts, preserves those permutations for small
nonzero $\epsilon$. The fourth sheet is regular at $\epsilon=0$:
on $a=0$ it is the other root $r=5t_0/3$ of
$3r^2-2t_0r+4-t_0^2=0$, and $b=t_0+r\ne0$ there.

Choose $D_{\rm hub}=(t_0+2)/3$, so at $\epsilon=0$ the map $z\mapsto t$
is locally invertible at $z=2$ with $t(2)=t_0$. The four sheets just
described give admissible solutions on the lifted $z$-loops. Indeed,
at the limiting points $d=-2$, $h=3$ and $r$ is nonzero and imaginary.
A zero $v_j$ would require $P_j=\pm1$ and $r=\pm\gamma_j h/d$,
which is real there. The fourth sheet is fixed on these loops.
The generic group thus contains an $S_3$ fixing one sheet.
Transitivity gives orbit size four for that sheet, so its stabilizer
of order at least six forces the full group to be $S_4$.

\smallskip
\emph{The complete graph $K_4$ and a doubled edge.}
First let the weights of two opposite edges tend to zero. The four
limiting systems are on a square. The constants are all $1$, zero on
one bipartition class and $1$ on the other (two systems), or all $0$.
The corresponding degrees are $2,3,3,7$. In the two middle cases one
chooses a nonempty subset of the two vertices with constant $1$ and
takes their incident edges. In the last case the edge sets are the
two matchings, the four three-edge sets, and the whole square.
By Theorem~\ref{irr:thm:general} the group of the seven-sheet system
is transitive for generic coefficients. Loops avoiding the exceptional
values of all four systems still act transitively on those seven sheets.
The simultaneous group therefore has order divisible by $7$ and
contains a $7$-cycle, since an element of order $7$ fixes the other
eight sheets. The zero-weight equations have derivatives $\pm1$, so
this permutation persists at generic coefficients of $K_4$.
Since $7$ exceeds the proper divisors of $15$, this proves
primitivity as above.

For $K_4$ with edge $12$ doubled, deleting the extra edge gives degrees
$15$ and $23$. To count the shifted system, require degrees at $1,2$
to be positive and those at $3,4$ to differ from $1$. With edge $34$
absent, a nonempty selection of the two stars from $3,4$ to $1,2$
allows either choice of edge $12$, and the empty selection requires
edge $12$, giving seven possibilities. With $34$ present, each of
$3,4$ must meet at least one of $1,2$. With $12$ present all nine
choices are allowed, and with $12$ absent seven are allowed. This
gives $23$ and a prime divisor $23$ of the group order, proving
primitivity in degree $38$.

For a transposition in either graph, set all diagonals to zero and the weight
of edge $12$ to $1$. Initially all other weights are zero. Put $c=z^{-2}$.
Each choice of $0$ or $-1$ for the deleted edges leaves an equation
\[
 f(1+f)=c(\mu_1+f)(\mu_2+f)
\]
for $f=f_{12}$. Write $M=-1$ for $K_4$ and $M=-2$ for the doubled
graph. The endpoint constants range from $M$ to $1$. The discriminant is
\[
 1+\bigl((2\mu_1-1)(2\mu_2-1)-1\bigr)c+(\mu_1-\mu_2)^2c^2.
\]
At $c_*=-1/(4M(M-1))$ it vanishes only for $\mu_1=\mu_2=M$, where
the coefficient of $c$ is maximal. The other admissible sheets are
regular there: $c_*\ne1$, the excluded roots can only be $0,-1$ and
are constant when they occur, and each isolated vertex has a nonzero
constant on an admissible sheet. There are two double roots in the
full fiber, one for each $f_{34}=k\in\{0,-1\}$, with every deleted
edge incident to $1$ or $2$ set to $-1$. Both double roots are
admissible.

Turn on $\tau=w_{13}^2$ and write $f_{13}=-1+q$. Near either double
root the equations are
\[
 (-1+q)q=\tau c(f+M+q)(k-1+q),\qquad
 F=f(1+f)-c(f+M+q)(f+M)=0 .
\]
The first equation determines $q(f,c,\tau)$ analytically, with
$q_\tau=-c(f+M)(k-1)$ at $\tau=0$. The double roots are simple
ramifications since $F_{ff}=2(1-c_*)$ and $F_c=-(f+M)^2$ there.
Solving $F=F_f=0$ gives derivative of the branch value
$c_k'(0)=-F_\tau/F_c=c_*^2(k-1)$, different for the two choices.
Thus a loop around just one split branch value interchanges two
sheets and fixes all others. The change $c=z^{-2}$ is locally invertible
and the loop persists at generic coefficients. Together with
primitivity this gives $S_{15}$ and $S_{38}$.

\smallskip
\emph{The prism and $K_{3,3}$.}
Each is a hexagon with three disjoint chords. In the prism the two triangles with edges $01,12,20$ and
$34,45,53$ and vertical edges $03,14,25$ contain the hexagon
$0,1,4,3,5,2,0$. In $K_{3,3}$ any alternating ordering of its two
vertex classes gives a hexagon.
Write $H_j$ for the hexagon with $j$ of the chords, and give its
vertices constants $\mu_i=d_{H_j}(i)-2$. The degree counts complements
of subsets with $d_S(i)\ne\mu_i$, hence edge sets with no vertex of
degree $2$. If $k$ chords occur in such a set, the two cycle-edge bits
at each of their $2k$ ends must be equal. At other vertices two
consecutive bits cannot both be $1$. Collapsing the equalities gives
cyclic binary strings without adjacent $1$'s on $6-2k$ positions,
with counts $18,7,3$ for $k=0,1,2$ (matchings on a hexagon, square,
and two parallel edges, including the empty matching).
For $k=3$ all bits coincide, giving two choices. Thus for $j=0,1,2,3$
the degrees are $18,25,35,50$.

For each $j\ge1$ there is a transposition. Keep a chord $uv$ of weight
$1$, put all diagonals to zero and freeze the other weights at zero.
In the single-edge discriminant above the endpoint constants range
from $-1$ to $1$, and the only double roots at $c_*=-1/8$ have both
constants $-1$. The frozen edges incident to $u$ or $v$ are then all
set to $-1$. Turn on those incident edges with a common squared
weight $\tau$. For a frozen edge $b$ incident to $u$ or $v$, denote
its other end by $x_b$ and write $g_{x_b}^0$ for its vertex sum at
$\tau=0$. The calculation for the doubled edge gives the derivative
$c_*^2\sum_b g_{x_b}^0$ of the branch value: to first order in $\tau$
the incident edges enter the equations of $u$ and $v$ additively, so
their contributions add. Every $g_{x_b}^0\ge-2$, and all
are $-2$ when every frozen edge is set to $-1$. This gives an
admissible double root whose derivative is strictly smaller (more
negative) than that of every other one: any edge
omitted from that freezing is disjoint from $uv$, and at least one
of its ends is a cycle neighbor of $u$ or $v$, since the two ends
of a hexagon chord and their cycle neighbors leave at most one
vertex uncovered. A small loop about the separated branch value
therefore gives a single transposition.

Deleting the chord of $H_1$ gives limiting families of degrees $7$
and $18$ (complements containing or omitting the chord). A
transitive action on the seven-sheet family gives a prime divisor
$7$ of the group order, so in degree $25$ the group is primitive
and thus $S_{25}$. For $H_2$ and $H_3$, the shifted family after
deleting a chord is respectively $H_1$ and $H_2$ with the prescribed
constants. Their groups supply the prime divisors $23$ and $31$,
larger than any proper divisors of $35$ and $50$, respectively.
With the transposition this proves $S_{35}$ and $S_{50}$.
On $H_3$ the constants are all $1$.
\end{proof}

\begin{remark}[The monodromy group]\label{rem:monodromy-group}
For a disconnected graph the group is contained in the product of the
groups of the components. Conjecture~\ref{open:monodromy} gives the
expected group for connected $G$ at generic coefficients. In addition to
Proposition~\ref{irr:prop:small-monodromy}, the subdivision result
(Proposition~\ref{blc:prop:subdivision} in the scalar case) gives $S_5$
for the diamond and $K_{2,3}$, both subdivisions of $D_3$.
A refinement of the edge-degeneration induction gives a criterion for the
full symmetric group on more general graphs, conditional on separation and
nondegeneracy properties of the fibers over the signed eigenvalues that
occur in the induction; we do not pursue it here.
\end{remark}

\begin{remark}[The commutative analogue]\label{irr:rem:commutative}
For $\Z^d$-periodic operators on graphs, the analogue of
Theorem~\ref{thm:irreducible} is the irreducibility of Bloch and Fermi
varieties, which Faust and Liu establish for generic parameters
\cite{FaustLiu}. In that setting irreducibility is a tool for excluding
embedded eigenvalues \cite{Kuchment}. For the operators $\tH$ on trees
there is no such implication: by Proposition~\ref{prop:uniform}\,(ii) the
irreducible curve $\Cphys$ can carry a pole of the physical branch, that
is, an eigenvalue of $\tH$, and the classification of
Theorem~\ref{thm:classification} includes such eigenvalues among the poles
with negative residue.
\end{remark}

%% file: sections/classification.tex
\section{Classification of the poles}
\label{sec:classification}

In this section we determine, for generic coefficients, all points of the
spectral curve over finite $z$ at which the differential
\[
\omega=\Bigl(\sum_{i\in V}V_i\Bigr)dz
\]
has a nonzero residue, and we show that for $\beta\ge2$ every eigenvalue
of every signed Jacobi matrix $D+A^s$ is an antibound pole of every Green
function. This proves Theorem~\ref{thm:main-poles}.

The classification concerns residues. A coordinate $V_i$ can also be
unbounded at a point where $\omega$ has no residue, as at the band edges
of a cycle, and there are ramification points at which all coordinates are
finite; these are treated in Section~\ref{sec:genus}. Our Green function is
$\langle\delta_x,(\tH-z)^{-1}\delta_x\rangle$, so a pole of the physical
branch at an eigenvalue $\lambda$ of $\tH$ has the residue
$-\mu_i(\{\lambda\})\le0$, where $\mu_i$ is the spectral measure, and poles
with positive residue are antibound
(Corollary~\ref{cor:antibound-sign}). Gutkin \cite{Gutkin}, who used the
sign of residues in this way for two families of free products, works with
$(z-H)^{-1}$; his signs are the opposite of ours. The poles on the physical
sheet, that is, the point spectrum of $\tH$, were described by Aomoto
\cite{Aomoto91} and by Banks, Garza-Vargas and Mukherjee \cite{BGVM}; they
are among the poles with negative residue classified below.

Recall that the signed eigenvalues are the eigenvalues, on the
antisymmetric subspace, of the two-fold covers of $G$
(Section~\ref{subsec:graphs}); in the language of Bilu and Linial
\cite{BiluLinial} they are the new eigenvalues of the $2$-lifts.

\subsection{Perron values are isolated}
\label{cls:subsec:perron}

\begin{lemma}\label{lem:perron-isolated}
Let $G$ be connected with arbitrary real coefficients. Let $U\subseteq V$
be nonempty, $s$ a signing of the edges of $G[U]$ and $\lambda$ an
eigenvalue of $(D+A^s)[U]$. Then $\sigm\le\lambda\le\sigp$. Moreover
$\lambda=\sigp$ only if $U=V$ and $s$ is switching equivalent to $\sgn w$,
and $\lambda=\sigm$ only if $U=V$ and $s$ is switching equivalent to
$-\sgn w$; in these cases the eigenvalue is simple and the squares of the
entries of its unit eigenvector are $(\bar u^\pm_i)^2$.
\end{lemma}

\begin{proof}
Let $v$ be a real eigenvector. Then
\[
\lambda\,v^{\T}v=v^{\T}(D+A^s)[U]\,v\le\abs v^{\T}(D+\abs A)[U]\,\abs v
\le\lambda_{\max}\bigl((D+\abs A)[U]\bigr)\,v^{\T}v .
\]
If $U\ne V$, the Perron root of the irreducible matrix $D+\abs A+c$, $c$
large, is strictly larger than that of any proper principal submatrix, so
$\lambda<\sigp$. If $U=V$ and $\lambda=\sigp$, equality in the second
inequality gives $\abs v=\bar u^+$ up to a factor, in particular $v_i\ne0$,
and equality in the first gives
$s_ew_ev_iv_j=\abs{w_e}\abs{v_i}\abs{v_j}$ for every edge $e=ij$, that is
$s_e\sgn w_e=\sgn v_i\,\sgn v_j$: switching by $\sgn v$ transforms $s$ into
$\sgn w$. The statements for $\sigm$ follow by applying this to $-D$ and
$-s$.
\end{proof}

\subsection{Pole data and the rescaled system}
\label{cls:subsec:data}

We fix notation for the rest of the section. Let $U\subseteq V$ be
nonempty with $G[U]$ connected, and let
\[
N=N(U)=\{x\notin U:\ x\text{ is adjacent to }U\},\qquad Y=V\setminus(U\cup N).
\]
$E(U)$ is the set of edges inside $U$, $d^U_i$ the number of them at
$i\in U$, $\partial U$ the set of edges between $U$ and $N$, and
$\partial_x$ the set of edges between $x\in N$ and $U$, with
$k_x=\abs{\partial_x}$; for $e\in\partial_x$ we write $i(e)$ for its end
in~$U$. The edges that have an end in $N$ and no end in $U$ are called
\emph{frozen}; they join $N$ to $N$ or $N$ to~$Y$.

\begin{definition}\label{cls:def:datum}
A \emph{pole datum} $\mathfrak d$ of $(G,D,w)$ consists of
\begin{enumerate}
\item[\textup{(1)}] a set $U$ as above, a switching class $[s]$ of
signings of $E(U)$, and an eigenvalue $\lambda$ of $(D+A^s)[U]$ with an
eigenvector $u$, $u^{\T}u=1$;
\item[\textup{(2)}] a choice $F$ of a value $f^0_e\in\{0,-1\}$ for every
frozen edge~$e$; we put
\[
\beta_x=-1-\sum_{e\ni x\text{ frozen}}f^0_e\ (x\in N),\qquad
c_y=1+\sum_{e\ni y\text{ frozen}}f^0_e\ (y\in Y),
\]
\[
\kappa(U,F)=\abs{E(U)}-\abs U-\sum_{x\in N}\beta_x;
\]
\item[\textup{(3)}] for every $x\in N$ a solution
$\bigl((f^*_e)_{e\in\partial_x},\nu_x\bigr)$ with $\nu_x\ne0$ of the
\emph{boundary system}
\begin{equation}\label{cls:eq:boundary}
f^*_e(1+f^*_e)=w_e^2\,u_{i(e)}^2\,\nu_x\quad(e\in\partial_x),\qquad
\sum_{e\in\partial_x}f^*_e=\beta_x;
\end{equation}
\item[\textup{(4)}] an admissible solution $\psi$, at $z=\lambda$, of the
\emph{exterior system}, which is the system \eqref{irr:eq:system-mu} of
$G[Y]$ with the constants $c=(c_y)_{y\in Y}$.
\end{enumerate}
We write $\kappa(\mathfrak d)=\kappa(U,F)$.
\end{definition}

The integers $\beta_x$ are at least $-1$, and the constants $c_y$ are at
most~$1$. The second equation of \eqref{cls:eq:boundary} says that the
right-hand side $1+\sum_{e\ni x}f_e$ of the vertex equation at $x$
vanishes. Summing it over $x$ gives
$\kappa(\mathfrak d)=\abs{E(U)}-\abs U-\sum_{e\in\partial U}f^*_e$. For
$U=V$ a pole datum is a pair $([s],\lambda)$ and $\kappa=\abs E-n=\kappa_0$.
The pole data with $\kappa\ne0$ are the data of
Theorem~\ref{thm:partial-pole}, in a different normalization of $\nu$; the
pole data with arbitrary $\kappa$ reappear as the configurations of
Definition~\ref{def:configuration}.

The next lemma extends the coordinates of Lemma~\ref{irr:lem:eta} to the
situation where $U$ is not a union of components. We keep
$y^0_e=-2s_ew_eu_iu_j$ for $e=ij\in E(U)$.

\begin{lemma}[Rescaled system]\label{cls:lem:chart}
Fix $(U,[s],\lambda,u)$ as in Definition~\ref{cls:def:datum}\,\textup{(1)}
with $\lambda$ simple, $u$ without zero entry, and
$\lambda\notin\{D_y:y\in Y\}$. Consider the following system in the
unknowns $\tilde V\in\C^U$, $y\in\C^{E(U)}$, $\zeta\in\C$, $\nu\in\C^N$ and
$f_e$ for all edges $e\notin E(U)$, with a parameter~$\eta$:
\begin{equation}\label{cls:eq:chart}
\begin{aligned}
&(D_i-\lambda-\zeta)\tilde V_i-\eta\Bigl(1-\frac{d^U_i}2+\sum_{e\in\partial U,\,e\ni i}f_e\Bigr)-\frac12\sum_{e\in E(U),\,e\ni i}y_e=0\quad(i\in U),\\
&y_e^2-\eta^2-4w_e^2\tilde V_i\tilde V_j=0\quad(e=ij\in E(U)),\qquad\sum_{i\in U}\tilde V_i=1,\\
&f_e(1+f_e)-w_e^2\tilde V_{i(e)}\nu_x=0\quad(e\in\partial_x),\qquad
(D_x-\lambda-\zeta)\,\eta\nu_x-1-\sum_{e\ni x}f_e=0\quad(x\in N),\\
&f_e(1+f_e)-w_e^2V_aV_b=0\quad(e=ab\notin E(U)\cup\partial U),
\end{aligned}
\end{equation}
where in the last line $V_x=\eta\nu_x$ for $x\in N$ and
$V_y=\bigl(1+\sum_{e\ni y}f_e\bigr)/(D_y-\lambda-\zeta)$ for $y\in Y$.
\begin{enumerate}
\item[\textup{(a)}] For $\eta\ne0$ the substitution
$z=\lambda+\zeta$, $V_U=\tilde V/\eta$, $f_e=\frac12(y_e/\eta-1)$ for
$e\in E(U)$, $V_N=\eta\nu$ identifies the solutions of
\eqref{cls:eq:chart} with the solutions of \eqref{eq:system} such that
$z\notin\{D_y:y\in Y\}$ and $\sum_{i\in U}V_i=1/\eta$.
\item[\textup{(b)}] The solutions of \eqref{cls:eq:chart} at $\eta=0$ with
$(\tilde V,y,\zeta)=(u^2,y^0,0)$ are the points
\[
b=(u^2,y^0,0,\nu,f^*,f^0,\psi).
\]
Here $f^0_e\in\{0,-1\}$ on the frozen
edges, $(f^*,\nu)$ solves the boundary systems, and $\psi$ solves the
exterior system at $z=\lambda$. Zero $\nu_x$ and inadmissible $\psi$ are
allowed. If every boundary solution $(f^*_{\partial_x},\nu_x)$ and the
exterior solution $\psi$ are nondegenerate in their respective systems,
then the Jacobian matrix of
\eqref{cls:eq:chart} with respect to the unknowns is invertible at~$b$.
Hence there is a unique analytic family $\Xi_b(\eta)$ of solutions through
$b$, and every solution of \eqref{cls:eq:chart} close to $b$ with small
$\eta$ belongs to it.
\item[\textup{(c)}] Along $\Xi_b$ we have $\zeta=\kappa\,\eta+O(\eta^2)$
with $\kappa=\abs{E(U)}-\abs U-\sum_{e\in\partial U}f^*_e$.
\end{enumerate}
\end{lemma}

\begin{proof}
(a) is the computation of Lemma~\ref{irr:lem:eta}\,(b), in which the
constant $\mu_i$ of the vertex $i\in U$ is replaced by
$1+\sum_{e\in\partial U,\,e\ni i}f_e$; for the boundary edges
$V_{i(e)}V_x=\tilde V_{i(e)}\nu_x$, and the vertex equation of $x\in N$ is
unchanged.

(b) The first statement is read off from \eqref{cls:eq:chart} at $\eta=0$.
For the Jacobian matrix, order the unknowns and the equations in four
groups: $(\tilde V,y,\zeta)$ with the first two lines; the frozen edges;
the edges of $G[Y]$; and $(f_{\partial U},\nu)$ with the third line. At
$\eta=0$ the first two lines do not involve the boundary variables, the
equation of a frozen edge involves only its own variable, with derivative
$1+2f^0_e=\pm1$, the equations of the edges of $G[Y]$ involve $\zeta$, the
frozen variables and the variables of $G[Y]$, and the third line does not
involve the variables of $G[Y]$ nor~$\zeta$. The matrix is therefore block
triangular, and its diagonal blocks are the invertible matrix of
Lemma~\ref{irr:lem:eta}\,(a), a diagonal matrix with entries $\pm1$, the
Jacobian matrix of the exterior system at $\psi$, and the Jacobian matrices
of the boundary systems.

(c) As in Lemma~\ref{irr:lem:eta}\,(c), differentiate the first line at
$\eta=0$ and sum over $i\in U$: the contribution of $\tilde V'$ and $y'$
vanishes because $\ones$ spans the left kernel of $\mathcal J$, and we are
left with $\zeta'=-\sum_{i\in U}\bigl(1-d^U_i/2+\sum_{e\in\partial U,\,e\ni i}f^*_e\bigr)$.
\end{proof}

\begin{lemma}[Boundary systems]\label{cls:lem:boundary}
Let $k\ge1$, $\beta\in\Z$ and $a\in(\C^*)^k$ be such that
\begin{equation}\label{cls:eq:nonresonance}
\sum_{e=1}^k\varepsilon_e\sqrt{a_e}\ne0\quad\text{and}\quad
\sum_{e=1}^k\varepsilon_ea_e\ne0\qquad\text{for all }\varepsilon\in\{\pm1\}^k .
\end{equation}
Then the system $f_e(1+f_e)=a_e\nu$ $(1\le e\le k)$, $\sum_ef_e=\beta$ in
the unknowns $(f,\nu)$ has exactly $2^{k-1}$ solutions counted with
multiplicity, of which $\binom k{-\beta}$ have $\nu=0$; these are simple.
For generic $a$ all solutions are simple, so that the number of solutions
with $\nu\ne0$ is
\[
n(k,\beta)=2^{k-1}-\binom k{-\beta},
\]
where $\binom{k}{j}=0$ for $j\notin\{0,\dots,k\}$.
\end{lemma}

\begin{proof}
Eliminating $\nu$ gives the equations
$a_1f_e(1+f_e)=a_ef_1(1+f_1)$, $2\le e\le k$, and $\sum_ef_e=\beta$, of
degrees $(2,\dots,2,1)$. A common zero at infinity satisfies
$f_e=\varepsilon_e\sqrt{a_e/a_1}\,f_1$ and
$f_1\sum_e\varepsilon_e\sqrt{a_e}=0$, so there is none by the first
condition, and B\'ezout's theorem gives $2^{k-1}$ solutions. Those with
$\nu=0$ are the points of $\{0,-1\}^k$ with $\sum_ef_e=\beta$. At such a
point, with $\rho_e=1+2f_e=\pm1$, a vector in the kernel of the Jacobian
matrix satisfies $df_e=a_e\rho_1\,df_1/(a_1\rho_e)$ and
$0=\sum_edf_e=(\rho_1/a_1)\bigl(\sum_e\rho_ea_e\bigr)df_1$, hence vanishes by
the second condition. Finally the map $(f,\nu)\mapsto(f_e(1+f_e)/\nu)_e$
from the smooth variety $\{\sum_ef_e=\beta,\ \nu\ne0\}$ of dimension $k$ to
$\C^k$ is unramified over a nonempty Zariski open set by generic
smoothness, which means that for generic $a$ the solutions with $\nu\ne0$
are simple.
\end{proof}

For example $n(1,\beta)=1$ if $\beta\ge1$ and $n(1,0)=n(1,-1)=0$;
$n(2,-1)=0$, $n(2,0)=1$; $n(3,-1)=1$.

\subsection{Completeness}
\label{cls:subsec:completeness}

The generic coefficients of this subsection are those which satisfy the
following conditions, for all $U$ with $G[U]$ connected, all signings $s$
of $E(U)$, all eigenvalues $\lambda$ of $(D+A^s)[U]$ with eigenvector $u$,
and all choices $F$ of frozen values.
\begin{enumerate}
\item[\textup{(S)}] The conclusions of Lemma~\ref{irr:lem:separation}
hold.
\item[\textup{(B)}] For every $x\in N(U)$ the numbers
$a_e=w_e^2u_{i(e)}^2$, $e\in\partial_x$, satisfy
\eqref{cls:eq:nonresonance}, and all solutions of the boundary systems
\eqref{cls:eq:boundary} are simple.
\item[\textup{(E)}] The point
$\bigl(w_e^2/((D_y-\lambda)(D_{y'}-\lambda))\bigr)_{e=yy'\in E(Y)}$ belongs
to the set $\mathcal U_{\mathrm{good}}$ of Lemma~\ref{irr:lem:generic} for
the graph $G[Y]$ with the constants $c(F)$.
\end{enumerate}
Each condition holds on a nonempty Zariski open set. For (S) this is
Lemma~\ref{irr:lem:separation}. Over its open set the normalized eigenpair
variety has finite fibers, since every eigenvalue is simple and each
eigenvector has two normalizations $u^{\T}u=1$. For each fixed eigenpair on
$U$, the bad boundary and exterior coefficients lie in proper algebraic
subsets of the space of coefficients outside $U$, as we now show. Their
union in the eigenpair variety thus has smaller dimension, and its
projection to coefficient space is contained in a proper Zariski closed
set. For fixed $u$
the sums in \eqref{cls:eq:nonresonance}, namely
$\sum_e\varepsilon_ew_eu_{i(e)}$ up to signs and
$\sum_e\varepsilon_ew_e^2u_{i(e)}^2$, are nonzero polynomials in the
independent weights $w_e$, $e\in\partial_x$, and the map
$(w_e)_{e\in\partial_x}\mapsto a$ is dominant, which gives (B) by
Lemma~\ref{cls:lem:boundary}. For fixed $\lambda$ the point in (E) depends
dominantly on the coefficients of $G[Y]$, and
$\mathcal U_{\mathrm{good}}$ is nonempty and open.

Under (E), every solution of the exterior system at $z=\lambda$ is
nondegenerate and exactly
\[
L_{c(F)}(G[Y])=\#\{S\subseteq E(Y):\ d_S(y)\ne c_y\text{ for all }y\in Y\}
\]
of them are admissible. With Lemma~\ref{cls:lem:boundary} we obtain: for
generic coefficients the number of pole data with given $(U,[s],\lambda,F)$
is
\begin{equation}\label{cls:eq:multiplicity}
M(U,F)=L_{c(F)}(G[Y])\prod_{x\in N}n(k_x,\beta_x).
\end{equation}

The first step of the classification requires only condition~(S), and we
state it separately because it is used in
Theorems~\ref{thm:green-degree} and~\ref{thm:all-signed-poles}.

\begin{lemma}[Support of a divergence]\label{cls:lem:support}
Assume \textup{(S)}. Let $P$ be a point of the normalization of $\Xc$ over
a finite value $z^*$ at which some $V_i$ has a pole, let $t$ be a local
coordinate at $P$, $a=\min_i\ord_tV_i<0$ and
$U=\{i:\ord_tV_i=a\}$. Then $G[U]$ is connected, and there are a signing
$s$ of $E(U)$ and a vector $u$ without zero entry, $u^{\T}u=1$, such that
$z^*$ is a simple eigenvalue of $(D+A^s)[U]$ with eigenvector $u$ and, with
$\eta=(\sum_{i\in U}V_i)^{-1}$,
\[
\eta V_i\to u_i^2\ (i\in U),\qquad\eta\,(1+2f_e)\to y^0_e=-2s_ew_eu_iu_j\ (e=ij\in E(U))
\]
at $P$. For an edge with exactly one end in $U$ we have $\ord_tf_e>a$.
\end{lemma}

\begin{proof}
Write $V_i=v_it^a+\cdots$ with $v_i\ne0$ for $i\in U$. For $e=ij\in E(U)$
the edge equation gives $f_e=b_et^a+\cdots$ with $b_e^2=w_e^2v_iv_j$. For an
edge $e=ij$ with $i\in U$, $j\notin U$ we have $\ord_t(V_iV_j)>2a$, hence
$\ord_tf_e>a$. The coefficient of $t^a$ in the vertex equation of $i\in U$
is therefore
\[
(D_i-z^*)v_i=\sum_{e=ij\in E(U)}b_e .
\]
Choose square roots $\tilde u_i$ of $v_i$ and define $s_e\in\{\pm1\}$ by
$b_e=-s_ew_e\tilde u_i\tilde u_j$. The identity becomes
$\bigl((D+A^s)[U]-z^*\bigr)\tilde u=0$, and $\tilde u$ has no zero entry.
Restricting $\tilde u$ to the components of $G[U]$ and using (S), we see
that $G[U]$ is connected, $z^*$ is a simple eigenvalue and
$\tilde u^{\T}\tilde u=\sum_{i\in U}v_i\ne0$. With
$u=\tilde u/(\tilde u^{\T}\tilde u)^{1/2}$ the limits follow, since
$\eta=t^{-a}\bigl((\sum_Uv_i)^{-1}+O(t)\bigr)$.
\end{proof}

\begin{theorem}[Classification of the residue poles]\label{thm:classification}
Let $G$ be a finite loopless multigraph. For generic coefficients the
following holds.
\begin{enumerate}
\item[\textup{(i)}] Every pole datum $\mathfrak d$ with
$\kappa(\mathfrak d)\ne0$ determines a point $P_{\mathfrak d}$ of the
normalization of $\Xc$ over $z=\lambda$. This point is unramified, and
there
\[
\begin{aligned}
V_i&=\frac{\kappa(\mathfrak d)\,u_i^2}{z-\lambda}+O(1) &&(i\in U),\\
V_x&=O(z-\lambda),\quad V_x\not\equiv0 &&(x\in N),\\
V_y&\longrightarrow V_y(\psi)\ne0 &&(y\in Y),
\end{aligned}
\]
so that $\Res_{P_{\mathfrak d}}\omega=\kappa(\mathfrak d)$. Distinct data
give distinct points.
\item[\textup{(ii)}] Let $P$ be a point of the normalization of $\Xc$ over
a finite value $z^*$ such that $(z-z^*)V_i$ does not tend to $0$ at $P$ for
some~$i$. Then $P=P_{\mathfrak d}$ for a pole datum $\mathfrak d$ with
$\kappa(\mathfrak d)\ne0$. This applies in particular to every point over
finite $z$ at which $\omega$ has a nonzero residue.
\item[\textup{(iii)}] The number of points $P_{\mathfrak d}$ with given
$(U,[s],\lambda,F)$ and $\kappa(U,F)\ne0$ is $M(U,F)$, see
\eqref{cls:eq:multiplicity}. For every $U$ with $G[U]$ connected there are
$\abs U\,2^{\beta(G[U])}$ pairs $([s],\lambda)$, where
$\beta(G[U])=\abs{E(U)}-\abs U+1$.
\end{enumerate}
\end{theorem}

Part~(i) is Theorem~\ref{thm:partial-pole}; we include its short derivation
from Lemma~\ref{cls:lem:chart} because the coordinates are needed for~(ii).
For $U=V$ it is Theorem~\ref{thm:local-pole}.

\begin{proof}
(i) By (B) and (E) the base point $b(\mathfrak d)$ satisfies the hypotheses
of Lemma~\ref{cls:lem:chart}\,(b). Along $\Xi_{b(\mathfrak d)}(\eta)$ we
have $z-\lambda=\kappa\eta+O(\eta^2)$ with $\kappa=\kappa(\mathfrak d)\ne0$,
so $z-\lambda$ is a local coordinate, and for small $\eta\ne0$ the
corresponding solution of \eqref{eq:system} is admissible:
$V_U=\tilde V/\eta\ne0$, $V_N=\eta\nu\ne0$, the exterior part is close to
the admissible $\psi$, and $z\notin\{D_i\}$ because $z$ is injective near
$\eta=0$. This arc lies on $\Xc$ and defines an unramified point of its
normalization at $\eta=0$. The expansions follow from
$V_i=\tilde V_i/\eta$ and $\tilde V_i(0)=u_i^2$, and the residue of
$\omega$ is $\kappa\sum_{i\in U}u_i^2=\kappa$. The datum is recovered from
the limits of the rescaled coordinates at the point; note that $(u^2,y^0)$
determines the switching class of~$s$.

(ii) Let $t$ be a local coordinate at $P$, $r=\ord_t(z-z^*)\ge1$ and
$a=\min_i\ord_tV_i$; the hypothesis is $a\le-r$. The differential $dz$
has order $r-1$, while $\sum_iV_i$ has order at least $a$. Thus $\omega$
is holomorphic if $a>-r$. Let $U$, $s$, $u$,
$\lambda=z^*$ and $\eta$ be as in Lemma~\ref{cls:lem:support}, and $N$, $Y$
as above. We write $a_x=\ord_tV_x$.

\emph{Step 1: all $V_x$, $x\notin U$, are bounded at $P$.} Suppose not,
let $a'=\min_{x\notin U}a_x<0$ and $U'=\{x\notin U:a_x=a'\}$; recall
$a<a'$. If some $x\in U'$ is adjacent to $U$, the variables $f_e$,
$e\in\partial_x$, have order $(a+a')/2<a'$, with leading coefficients
$b_e$ satisfying $b_e^2=w_e^2v_xv_{i(e)}$, whereas all other terms of the
vertex equation of $x$ have order at least $a'$: the left-hand side
because $\ord_tV_x=a'$, the other edges $e=xx'$ at $x$ because
$a_x+a_{x'}\ge2a'$. The coefficient of $t^{(a+a')/2}$ gives
$\sum_{e\in\partial_x}\varepsilon_ew_eu_{i(e)}=0$ for suitable signs,
contradicting~(B). Hence $U'$ is not adjacent to $U$, and the argument of
Lemma~\ref{cls:lem:support} applies to $U'$: it shows that $z^*$ is an
eigenvalue of a signed matrix on a component of $G[U']$, which contradicts
(S) because $z^*$ is already an eigenvalue of $(D+A^s)[U]$.

\emph{Step 2: $a_x=-a$ for $x\in N$.} Let $x\in N$; by Step~1, $a_x\ge0$
and the variables of the edges from $x$ to vertices outside $U$ are
bounded. If $a+a_x<0$, the coefficient of $t^{(a+a_x)/2}$ in the vertex
equation of $x$ comes from the edges in $\partial_x$ only, and we get the
same contradiction with~(B) as in Step~1. If $a+a_x>0$, every edge $e$ at
$x$ has $f_e\to f^0_e\in\{0,-1\}$. For $e\in\partial_x$,
$f_e-f^0_e=w_e^2V_{i(e)}V_x/(1+2f^0_e)+\cdots$ has order exactly $a+a_x$,
while for the other edges at $x$ the difference has order at least $a_x$,
as has the left-hand side of the vertex equation. The constant term
$1+\sum_{e\ni x}f^0_e$ must vanish, and the coefficient of $t^{a+a_x}$
gives $\sum_{e\in\partial_x}\varepsilon_ew_e^2u_{i(e)}^2=0$, again
contradicting~(B).

\emph{Step 3: convergence to a base point.} By
Lemma~\ref{cls:lem:support} and Steps~1 and~2, as $t\to0$: $\eta\to0$,
with $\ord_t\eta=-a$; $(\eta V_U,\eta(1+2f_{E(U)}))\to(u^2,y^0)$;
$V_x/\eta$ tends to a finite nonzero limit $\nu_x$ for $x\in N$; on
$\partial U$ the products $V_{i(e)}V_x$ and hence the $f_e$ have finite
limits $f^*_e$; on the frozen edges $V_aV_b\to0$ and $f_e\to f^0_e\in\{0,-1\}$;
on $G[Y]$ all variables have finite limits; and $\zeta=z-\lambda\to0$. By
Lemma~\ref{cls:lem:chart}\,(a) the points of $\Xc$ near $P$ are solutions
of \eqref{cls:eq:chart}, and they converge to a point $b$ as in
Lemma~\ref{cls:lem:chart}\,(b) with all $\nu_x\ne0$. By (B) and (E) the
Jacobian matrix is invertible at $b$, so the branch of $\Xc$ at $P$ is
parametrized by $t\mapsto\Xi_b(\eta(t))$.

\emph{Step 4: the exterior limit $\psi$ is admissible.} Suppose that
$g_y=1+\sum_{e\ni y}f_e$ vanishes at $b$ for some $y\in Y$. Then
$f_e\in\{0,-1\}$ at $b$ for all edges at $y$, which are frozen edges or
edges of $G[Y]$. In the Jacobian matrix of \eqref{cls:eq:chart} at $b$,
the rows of these edges have nonzero entries only in the columns of these
same edges. At the base point $V_y=0$, and $dV_y=(D_y-\lambda)^{-1}dg_y$
in the unknowns, with $dg_y=\sum_{e\ni y}df_e$. This also covers a frozen
edge, whose other end has $V_x=0$ at the base point. Since the
matrix is invertible, so is the complementary block, which is the Jacobian
matrix of the system obtained from \eqref{cls:eq:chart} by fixing the
variables of the edges at $y$ at their values at $b$ and deleting their
equations. This smaller system has an analytic family of solutions through
$b$; along it $g_y\equiv0$, so the deleted equations hold as well. By the
uniqueness in Lemma~\ref{cls:lem:chart}\,(b) this family is $\Xi_b$. Thus
$V_y$ vanishes identically on the branch of $\Xc$ at $P$. This would make
$V_y$ identically zero on its irreducible component, which has a nonempty
admissible open set. Hence $\psi$
is admissible, and $b=b(\mathfrak d)$ for a pole datum $\mathfrak d$.

\emph{Step 5: $\kappa\ne0$.} Along $\Xi_b$ the function $\zeta(\eta)$ is
not identically zero, since $z$ is not constant near $P$; let $k$ be its
order at $\eta=0$. Then $r=k\cdot\ord_t\eta=-ka\ge kr$, so $k=1$: by
Lemma~\ref{cls:lem:chart}\,(c), $\kappa(\mathfrak d)\ne0$, and $P$ is the
point $P_{\mathfrak d}$ of~(i).

(iii) follows from \eqref{cls:eq:multiplicity} and (S).
\end{proof}

Steps~1 to~4 use only that some $V_i$ has a pole at $P$; they show that
every such point is the limit point of the family $\Xi_{b(\mathfrak d)}$
of a pole datum with arbitrary $\kappa$. This is the starting point of the
genus formula (Lemma~\ref{gen:lem:pole-points} and
Theorem~\ref{thm:genus-formula}), where the data with $\kappa=0$ are
analyzed.

Combined with Proposition~\ref{prop:residue-sum}, according to which the
sum of the residues of $\omega$ over all finite points is
$\sum_{S\in\cL(G)}(2\abs S-n)$, the theorem yields an identity between
integers that depend only on the graph.

\begin{corollary}\label{cls:cor:identity}
For every finite loopless multigraph $G$,
\[
\sum_{S\in\cL(G)}\bigl(2\abs S-n\bigr)
=\sum_{\substack{U\subseteq V\\ G[U]\text{ connected}}}\abs U\,2^{\beta(G[U])}\sum_F\kappa(U,F)\,M(U,F),
\]
where $F$ runs over the choices of frozen values for~$U$.
\end{corollary}

\begin{corollary}[Local identity]\label{cor:local-identity}
For every finite loopless multigraph $G$ and every vertex $v$,
\[
\sum_{S\in\cL(G)}\bigl(d_S(v)-1\bigr)
=\sum_{\substack{U\ni v\\ G[U]\text{ connected}}}2^{\beta(G[U])}\sum_F\kappa(U,F)\,M(U,F).
\]
\end{corollary}

Summing over $v$ gives Corollary~\ref{cls:cor:identity}.

\begin{proof}
Choose generic coefficients with pairwise distinct $D_i$ and small
weights, $\abs{w_e}\le\epsilon$, and let $r$ be less than half the minimal
distance between the $D_i$. On the circle $\abs{z-D_v}=r$ all the
quantities $c_e(z)$ of \eqref{irr:eq:f-system} are $O(\epsilon^2)$, so by
the proof of Lemma~\ref{irr:lem:infinity} the solutions over the circle
lie on $2^{\abs E}$ disjoint branches close to the points $-\ones_S$, the
admissible ones being those with $S\in\cL(G)$, and on the branch $S$ we
have $V_i=(1-d_S(i))/(D_i-z)+O(\epsilon^2)$. Therefore
\[
\frac1{2\pi i}\sum_{S\in\cL(G)}\oint_{\abs{z-D_v}=r}\sum_iV^S_i\,dz
=\sum_{S\in\cL(G)}\bigl(d_S(v)-1\bigr)+O(\epsilon^2).
\]
The left-hand side is the sum of the residues of $\omega$ at the points of
the normalization over the disc. By Theorem~\ref{thm:classification} these
are the numbers $\kappa(\mathfrak d)$ of the pole data with $\lambda$ in
the disc. For small weights the eigenvalues of $(D+A^s)[U]$ are close to
the $D_i$, $i\in U$, one near each; so $\lambda$ lies in the disc exactly
when $v\in U$ and $\lambda$ is the eigenvalue near $D_v$, which happens
once for each switching class. The sum of the residues is the right-hand
side of the corollary. Both sides are integers independent of $\epsilon$,
hence equal.
\end{proof}

\begin{example}\label{cls:ex:K4}
For the triple edge $D_3$ the pole data with $\kappa\ne0$ are the eight
pairs $([s],\lambda)$ with $U=V$, $\kappa=1$; for $U=\{v\}$ the other
vertex has $k_x=3$, $\beta_x=-1$ and $\kappa=0$. Indeed
$\sum_{S\in\cL}(2\abs S-2)=-2+3\cdot2+4=8$. For $K_4$ the points
$P_{\mathfrak d}$ are: $32$ with $U=V$ and $\kappa=2$; $24$ with $U$ a
triangle ($6$ signed eigenvalues each, $n(3,-1)=1$) and $\kappa=1$; $12$
with $U$ an edge, where the opposite edge is frozen at $-1$
($n(2,0)^2=1$) and $\kappa=-1$, the frozen value $0$ giving
$n(2,-1)=0$; and $4$ with $U$ a vertex, where the opposite triangle is
frozen at $-1$ ($n(1,1)^3=1$) and $\kappa=-4$. The residues add up to
$64+24-12-16=60=\sum_{S\in\cL(K_4)}(2\abs S-4)$. For the bowtie there are,
besides the $20$ poles with $U=V$ and $\kappa=1$, eight poles with
$\kappa=-2$ supported on the two vertices of a triangle other than the cut
vertex; for the diamond, besides the $16$ poles with $U=V$, two poles with
residue $-3$ over the potentials of the two vertices of degree two
(Proposition~\ref{prop:degree-two-vertex}).
\end{example}

\begin{remark}\label{cls:rem:census}
Corollary~\ref{cls:cor:identity} is a purely combinatorial statement,
whose proof above goes through the spectral curve. As an independent
check, we verified it by exact enumeration of both sides for all $33867$
labeled simple graphs with at most six vertices and for $300$ random
multigraphs; this computation is not used in any proof.
\end{remark}

\subsection{All signed eigenvalues are antibound poles}
\label{cls:subsec:signed}

\begin{theorem}\label{thm:all-signed-poles}
Let $G$ be connected with $\beta\ge2$. For generic real coefficients the
following holds.
\begin{enumerate}
\item[\textup{(a)}] Let $s$ be a signing and $\lambda$ an eigenvalue of
$D+A^s$ with unit eigenvector~$u$. Then the analytic continuation of the
physical branch along a suitable path that avoids a finite subset of $\C$
reaches a branch that is single-valued near $z=\lambda$ and on which every
$G_i$ has a simple pole at $\lambda$ with residue $\kappa_0u_i^2>0$. There
are $n2^\beta$ pairs $([s],\lambda)$, and the corresponding $n2^\beta$
values $\lambda$ are pairwise distinct and distinct from all eigenvalues
of the matrices $(D+A^{s'})[U]$ with $U\ne V$, $G[U]$ connected.
\item[\textup{(b)}] None of these poles lies on the physical sheet: they
are antibound.
\item[\textup{(c)}] Conversely, if some analytic continuation of some
$G_i$ is unbounded near a finite point $z_0$, then $z_0$ is an eigenvalue
of $(D+A^s)[U]$ for some $U$ with $G[U]$ connected and some signing $s$.
\end{enumerate}
\end{theorem}

Thus, for generic coefficients, the finite singularities at which the
continued Green functions are unbounded lie in
$\bigcup_{U,s}\spec\bigl((D+A^s)[U]\bigr)$, and the part
$\bigcup_s\spec(D+A^s)$ of this set consists of simple poles with positive
residues that all occur. This answers the question raised after
Conjecture~3.4 of \cite{CSZ}, whether there are families of Jacobi matrices
for which the only poles on higher sheets are at the Perron values: such a
family must consist of non-generic coefficients, as in
Proposition~\ref{prop:uniform}. By Lemma~\ref{lem:perron-isolated} and
Theorem~\ref{thm:perron-component}, the two Perron poles are the only ones
among the $n2^\beta$ that persist for all real coefficients.

\begin{proof}
We take real coefficients in the generic sets of
Theorem~\ref{thm:irreducible} and of condition~(S).

(a) By (S), $\lambda$ is a simple eigenvalue, $u$ has no zero entry and
$\lambda\notin\{D_i\}$; since $\abs E>n$,
Theorem~\ref{thm:local-pole} gives a branch of $\Xc$ over a punctured disc
$0<\abs{z-\lambda}<\rho$ that is single-valued and on which
$V_i=\kappa_0u_i^2/(z-\lambda)+O(1)$. Fix a point $z_1$ of this punctured
disc and a path from a base point $z_0<\Sigm$ to $z_1$ in the complement of
the finite set $B$ of Corollary~\ref{cor:all-sheets}. The continuation of
the physical branch along this path ends at some admissible solution over
$z_1$; by the transitivity of the monodromy we can precede the path by a
loop based at $z_0$ so that it ends at the point of the pole branch. The
number of switching classes is $2^\beta$, and the spectrum and the $u_i^2$
depend only on the class. The values are distinct by~(S).

(b) is Corollary~\ref{cor:antibound-sign}.

(c) A continuation of $G_i$ is a path on $\Xc$; if $V_i$ is unbounded
along it as $z\to z_0$, the path converges to a point of the normalization
over $z_0$ at which $V_i$ has a pole, and
Lemma~\ref{cls:lem:support} applies.
\end{proof}

\subsection{The product of the Floquet function over the sheets}
\label{cls:subsec:floquet}

Recall from Proposition~\ref{prop:floquet} that the function
$\Phi=\prod_{e}(1+f_e)\big/\prod_iV_i$ on $\Xc$ satisfies
$\omega=-d\log\Phi$, and that on the physical branch it is the Floquet
function of \cite{BBGSS}. Thus $\Res_P\omega=-\ord_P\Phi$ at every point,
and Theorem~\ref{thm:classification} describes the divisor of $\Phi$ over
the finite part of the $z$-line. We now compute the norm
$\mathrm{Nm}_z\Phi$, the product of the values of $\Phi$ over the $\LG$
sheets of $\Xc$, which is a rational function of~$z$. For $U$ with $G[U]$
connected let
\[
a_U=\sum_F\kappa(U,F)\,M(U,F),
\]
so that $a_V=\kappa_0$, and recall $p_{U,s}(z)=\det(z-(D+A^s)[U])$.

\begin{proposition}\label{prop:floquet-product}
Let $G$ be a finite connected loopless multigraph. For all coefficients
$D\in\C^V$, $w\in(\C^*)^E$,
\begin{equation}\label{cls:eq:norm}
\mathrm{Nm}_z\Phi=K_G(w)\prod_{\substack{U\subseteq V\\ G[U]\text{ connected}}}\ \prod_{[s]}p_{U,s}(z)^{-a_U},
\end{equation}
\[
K_G(w)=\prod_{S\in\cL(G)}\Bigl[(-1)^{n+\abs S}\prod_{e\in S}w_e^2\prod_{i\in V}(1-d_S(i))^{d_S(i)-1}\Bigr],
\]
where $[s]$ runs over the $2^{\beta(G[U])}$ switching classes on $G[U]$.
Moreover $\prod_{[s]}p_{U,s}(z)=\det(z-\widehat H_U)$, where $\widehat H_U$
is the lift of $(D+A)[U]$ to the Galois cover of $G[U]$ with group
$H_1(G[U];\Z/2)$. Consequently
\[
\sum_{\text{sheets}}\ \sum_{i\in V}V_i
=\sum_{U}a_U\,\tr\,(z-\widehat H_U)^{-1}.
\]
\end{proposition}

For special coefficients $\Xc$ may be reducible; the norm is then taken
over all its components. Note the sign: the right-hand side of the last
formula involves the resolvents $(z-\widehat H_U)^{-1}$, whereas the
physical branch is a diagonal entry of $(\tH-z)^{-1}$.

\begin{proof}
First let the coefficients be generic in the sense of
Section~\ref{cls:subsec:completeness}. If $z-\lambda=t^r$ at a point $P$
over $\lambda$ and $\Phi=t^m\cdot(\text{unit})$, the product of $\Phi$ over
the $r$ local sheets has order $m/r$ in $t^r$, that is, order $m$ in
$z-\lambda$ after multiplication of the $r$ values; hence
$\ord_\lambda\mathrm{Nm}_z\Phi=\sum_{P\mapsto\lambda}\ord_P\Phi$. By
Theorem~\ref{thm:classification} the points with $\ord_P\Phi\ne0$ are the
$P_{\mathfrak d}$, with $\ord\Phi=-\kappa(\mathfrak d)$, and by (S) each
$\lambda$ is an eigenvalue for exactly one pair $(U,[s])$ and is a simple
root of $p_{U,s}$. Therefore both sides of \eqref{cls:eq:norm} have the
same zeros and poles in $\C$, and their quotient is constant. To find the
constant, consider the branch $S\in\cL(G)$ at infinity and write
$d_i=d_S(i)$: there $V_i\sim-(1-d_i)/z$, $1+f_e\to1$ for $e\notin S$ and
$1+f_e\sim-w_e^2(1-d_i)(1-d_j)z^{-2}$ for $e=ij\in S$, so the leading term of
$\Phi$ is the bracket in \eqref{cls:eq:norm} times $z^{n-2\abs S}$. The
$p_{U,s}$ are monic.

For arbitrary coefficients $p_0$, the branches $f^S(p,1/z)$ at infinity
are analytic in $(p,1/z)$ jointly for $p$ near $p_0$ and large $\abs z$
(Lemma~\ref{lem:infinity}), and the sheets of $\Xc$ over large $z$ are
those with $S\in\cL(G)$ for every $p$. Both sides of \eqref{cls:eq:norm}
are therefore continuous in $p$ for large $\abs z$, and the identity
extends from the dense set of generic coefficients to $p_0$, first for
large $\abs z$ and then as an identity of rational functions.

For the second statement choose a spanning tree of $G[U]$ and assign to
the edges outside it the elements $h_e$ of a basis of
$\Gamma=(\Z/2)^{\beta(G[U])}$, and $h_e=0$ to the tree edges. The cover has
vertices $(i,g)$, $g\in\Gamma$, and edges from $(i,g)$ to $(j,g+h_e)$ for
$e=ij$. In the basis of the characters $\chi$ of $\Gamma$ the lifted
matrix is the direct sum of the matrices $(D+A^s)[U]$ with
$s_e=\chi(h_e)$, and these signings represent each switching class once.
The last formula is the logarithmic derivative of \eqref{cls:eq:norm},
since the sum of $\omega$ over the sheets is $-d\log\mathrm{Nm}_z\Phi$.
\end{proof}

\begin{example}\label{cls:ex:dipole-norm}
For the triple edge, $a_V=1$ and $a_{\{v\}}=0$, so
$\mathrm{Nm}_z\Phi=-16(w_1w_2w_3)^6\big/\prod_{[s]}p_{V,s}(z)$, a product
over the four switching classes. For $D=0$, $w\equiv1$ this is
$-16/((z^2-9)(z^2-1)^3)$. In this uniform case the curve is reducible: by
Proposition~\ref{prop:uniform} the physical component has two sheets, on
which $f_e\equiv f$ and $\Phi=(1+f)^2/f$ with $z^2f(1+f)=(1+3f)^2$; the
product of $\Phi$ over these two sheets is $-16/(z^2-9)$. On each of the
three remaining sheets one edge has $f_e=1/(z^2-1)$, the two others have
$f_e=-z^2/(z^2-1)$, and $\Phi=1/(z^2-1)$. The poles at $\pm1$, where the three
non-Perron switching classes coincide, are not seen by the physical
component.
\end{example}

\begin{remark}\label{cls:rem:cancellation}
The exponent $a_U$ is a sum of residues of different signs and can
vanish although residue poles are present. For instance, attach two
diamonds to $K_5$, joining the two vertices of degree three of each
diamond to two fixed vertices of $K_5$ by one edge each. Evaluating
\eqref{cls:eq:multiplicity} for $U=V(K_5)$ gives, over each of the
$5\cdot2^6$ signed eigenvalues of $K_5$, sixteen points $P_{\mathfrak d}$
with residues $+1$ (nine times), $-1$ (six times) and $-3$ (once), so that
$a_U=0$: the norm of $\Phi$ and the sum of $\omega$ over the sheets are
regular there. The poles are detected by the second power sum: by the
local computation in the proof above,
$\lim_{z\to\lambda}(z-\lambda)^2\sum_{\text{sheets}}(\sum_iV_i)^2=\sum_{P\mapsto\lambda}(\ord_P\Phi)^2/r_P$,
where $r_P$ is the ramification index, which here equals~$24$.
\end{remark}

\subsection{Poles outside the convex hull of the tree spectrum}
\label{cls:subsec:hull}

By Corollary~\ref{cor:gap} the Perron poles lie outside $[\Sigm,\Sigp]$
when $\beta\ge2$. We now ask which other residue poles can lie there. The
basic tool is a comparison with subgraphs of cycle rank at most one.

\begin{lemma}\label{cls:lem:compare}
Let $G$ be connected with real coefficients, and let $J$ be a connected
subgraph of $G$ \textup{(}a set of vertices with a subset of the edges
joining them, not necessarily all\textup{)} of cycle rank at most one. Then
\[
\Sigm\le\lambda_{\min}(D_J-\abs{A_J})\qquad\text{and}\qquad
\lambda_{\max}(D_J+\abs{A_J})\le\Sigp,
\]
where $A_J$ is formed with the edges of $J$. Consequently, for every
signing $s$ of the edges of $J$, the spectrum of $D_J+A^s_J$ is contained
in $[\Sigm,\Sigp]$.
\end{lemma}

\begin{proof}
We may assume $w_e>0$ (Section~\ref{subsec:cover}). The universal cover of
$J$ embeds into $\tree$ as a subtree, and a subtree of a tree is an induced
subgraph. If $J$ is a tree, the Perron vector of $D_J+A_J$, placed on a
lift of $J$, is a test function whose Rayleigh quotient for $\tH$ is
$\lambda_{\max}(D_J+A_J)$. If $J$ has cycle rank one, its universal cover is
a periodic chain of copies of a fundamental domain; lift the Perron vector
periodically and truncate it to $k$ consecutive copies. The squared norm
grows like $k$, while the quadratic form differs from $k$ times that of
$J$ by a bounded boundary term, so the Rayleigh quotients converge to
$\lambda_{\max}(D_J+A_J)$. This gives the upper bound; the lower one
follows by replacing $D$ with $-D$. The last statement follows from
$v^{\T}(D_J+A^s_J)v\le\abs v^{\T}(D_J+\abs{A_J})\abs v$ and the
corresponding lower bound.
\end{proof}

\begin{proposition}\label{cls:prop:bicyclic-spectrum}
Let $G$ be connected with $\beta=2$ and arbitrary real coefficients. Every
eigenvalue of every signed Jacobi matrix $D+A^s$ lies in $[\Sigm,\Sigp]$,
with exactly two exceptions: $\sigp$, for $s$ switching equivalent to
$\sgn w$, and $\sigm$, for $s$ switching equivalent to $-\sgn w$.
\end{proposition}

\begin{proof}
We may assume $w_e>0$. Let $\lambda$ be an eigenvalue of $D+A^s$ with a
real eigenvector~$v$. Suppose first that $s$ is not switching equivalent
to the all-positive signing. After switching we may assume $v\ge0$. The
signing still has a cycle with negative sign product, hence a negative
edge that is not a bridge; removing it leaves a graph of cycle rank one,
so all components of the subgraph $J$ of the positive edges have cycle
rank at most one. Since $v\ge0$, the negative edges contribute
nonpositively to $v^{\T}(D+A^s)v$, and Lemma~\ref{cls:lem:compare}, applied
to the components of $J$, gives $\lambda\le\Sigp$. Suppose now that $s$ is
all-positive and $\lambda\ne\sigp$. Then $v$ is orthogonal to the positive
Perron vector, so the sets $V^+$ and $V^-$ where $v$ is positive and
negative are both nonempty. The cycle ranks of the components of $G[V^+]$
and $G[V^-]$ add up to at most $2$; exchanging $v$ and $-v$ we may assume
that all components of $G[V^+]$ have cycle rank at most one. Let $v^\pm$
be the positive and negative parts of $v$. Taking the inner product of the
eigenvalue equation with $v^+$,
\[
\lambda\norm{v^+}^2=\langle v^+,(D+A)v^+\rangle-\langle v^+,Av^-\rangle
\le\langle v^+,(D+A)v^+\rangle\le\Sigp\norm{v^+}^2
\]
by Lemma~\ref{cls:lem:compare}. The lower bounds follow by applying these
arguments to $-D$ and $-s$. That $\sigp$ and $\sigm$ lie outside
$[\Sigm,\Sigp]$ is Corollary~\ref{cor:gap}.
\end{proof}

\begin{corollary}\label{cor:bicyclic-hull}
Let $G$ be connected with $\beta=2$. For generic real coefficients, the
points of the normalization of $\Xc$ over $z\in\C\setminus[\Sigm,\Sigp]$,
$z\ne\infty$, at which $\omega$ has a nonzero residue are exactly two: the
Perron poles over $\sigm$ and $\sigp$.
\end{corollary}

\begin{proof}
By Theorem~\ref{thm:classification} such a point is $P_{\mathfrak d}$ for a
pole datum with support $U$ and real $\lambda\in\spec(D+A^s)[U]$. If $U=V$,
Proposition~\ref{cls:prop:bicyclic-spectrum} leaves the two Perron values,
and there is one point over each (Lemma~\ref{lem:perron-isolated},
Theorem~\ref{thm:local-pole}). If $U\ne V$ and $G[U]$ has cycle rank at most
one, then $\lambda\in[\Sigm,\Sigp]$ by Lemma~\ref{cls:lem:compare}.

There remains the case $U\ne V$ with $\beta(G[U])=2$. Counting edges and
vertices, $\beta(G)=\beta(G[U])+\sum_T\bigl(\beta(T)+e(U,T)-1\bigr)$, where
$T$ runs over the components of $G-U$ and $e(U,T)\ge1$ is the number of
edges between $U$ and $T$. Hence every $T$ is a tree joined to $U$ by a
single edge $e=ix$, $i\in U$, $x\in T$. By the recursion
\eqref{eq:m-recursion}, which holds on every admissible solution (see the
proof of Lemma~\ref{irr:lem:germs}), the function $m_{x|e}=V_x/(1+f_e)$ is
the same rational function $h_x(z)$ on all of $\Xc$, namely the diagonal
entry at $x$ of $(D_T+A_T-z)^{-1}$. With \eqref{eq:m-from-Vf},
$f_e=w_e^2h_xV_i$ and
\[
V_x=h_x\,(1+f_e)=h_x+w_e^2h_x^2\,V_i .
\]
The zeros and poles of $h_x$ are eigenvalues of signed matrices on
connected subsets of $T$, so by (S) $h_x(\lambda)$ is finite and nonzero.
Thus $V_x$ has a pole wherever $V_i$ has one over $\lambda$, which
contradicts $V_x\to0$ at $P_{\mathfrak d}$ for $x\in N$. So this case does
not occur.
\end{proof}

For $\beta=3$ the statement fails.

\begin{proposition}\label{cls:prop:rank-three}
Let $G$ consist of two vertices joined by four parallel edges, with $D=0$
and $w=(1,1,1,\epsilon)$, $0<\epsilon<(3-2\sqrt2)/2$. Then
$\lambda_{\max}(D+A^s)=3-\epsilon>\Sigp$ for $s=(+,+,+,-)$, while
$\sigp=3+\epsilon$. Consequently there is a nonempty open set of real
coefficients on this graph for which the continued Green functions have an
antibound pole outside $[\Sigm,\Sigp]$ other than the Perron poles.
\end{proposition}

\begin{proof}
The universal cover is the $4$-regular tree with edges of four colors.
Removing the edges of the fourth color leaves a forest of $3$-regular
trees, whose adjacency operator has norm $2\sqrt2$, and the removed edges
form a perfect matching, of norm $\epsilon$. Hence
$\Sigp\le2\sqrt2+\epsilon<3-\epsilon$. The inequality
$\lambda_{\max}(D+A^s)>\Sigp$ and $\lambda_{\max}(D+A^s)<\sigp$ persist on
an open set of real coefficients, since both sides are continuous; on its
generic points Theorem~\ref{thm:all-signed-poles} applies.
\end{proof}

In general the number of exceptional eigenvalues is bounded as follows.
Let $\taups(G)$ be the least number of vertices whose removal leaves a
graph all of whose components have cycle rank at most one.

\begin{proposition}\label{cls:prop:pseudoforest}
Let $G$ be connected with real coefficients. For every signing $s$, the
matrix $D+A^s$ has at most $\taups(G)$ eigenvalues, counted with
multiplicity, in $(\Sigp,\infty)$, and at most $\taups(G)$ in
$(-\infty,\Sigm)$. Moreover $\taups(G)\le\beta-1$ for $\beta\ge2$.
\end{proposition}

\begin{proof}
Let $T$ be a set of $\taups(G)$ vertices as in the definition. By
Lemma~\ref{cls:lem:compare} the compression of $D+A^s$ to the vectors
supported on $V\setminus T$ has its spectrum in $[\Sigm,\Sigp]$. A spectral
subspace of $D+A^s$ for $(\Sigp,\infty)$ of dimension larger than $\abs T$
would contain a nonzero vector supported on $V\setminus T$, a
contradiction; similarly below. For the last statement remove one end of
each of $\beta-1$ of the $\beta$ edges outside a spanning tree.
\end{proof}

\begin{remark}\label{cls:rem:matching}
Bounds of this kind for one well-chosen signing are classical. The
expectation of $\det(z-(D+A^s))$ over independent uniform signs is the
weighted matching polynomial of $G$, whose roots are real and lie in
$[\Sigm,\Sigp]$ by Godsil's path-tree argument \cite{Godsil}; the method of
interlacing families of Marcus, Spielman and Srivastava \cite{MSS} then
gives a signing whose largest eigenvalue is at most the largest root, see
also \cite{HPS}, and \cite{Spier} for weights and multiple edges. What is
specific to cycle rank two in Proposition~\ref{cls:prop:bicyclic-spectrum}
is that all signings other than the two Perron classes have their whole
spectrum in $[\Sigm,\Sigp]$.
\end{remark}

%% file: sections/real.tex
\section{Real points: counting positive fixed points}\label{sec:real}

In this section the spectral parameter is fixed and absorbed into the
potential: we replace $D-z$ by $D$ and study the positive solutions of
\eqref{eq:system} with $z=0$. By Lemma~\ref{lem:bp-correspondence} these
are the fixed points of Gaussian belief propagation all of whose vertex and
edge beliefs are positive definite. The coefficients are real, and no
assumption is made on the definiteness of the finite matrix $D+A$. We prove
Theorem~\ref{thm:main-real}.

The comparison matrix is $M=D-\abs A$, with $\abs A$ the edgewise absolute
value of Section~\ref{subsec:graphs}. We call the model \emph{walk-summable}
if $M$ is positive definite. For a simple graph this is the notion of
Malioutov, Johnson and Willsky~\cite{MJW}, and by
Proposition~\ref{prop:two-lift-threshold} it holds if and only if all
signed matrices $D+A^s$, equivalently the Jacobi matrices of all two-fold
covers~\cite{RuozziTatikonda}, are positive definite. For connected $G$ and
the diagonal $D-z$ it means $z<\sigm$.

Several facts about the positive fixed points are known. For walk-summable
models the iteration converges~\cite{MJW}. Su and Wu~\cite{SuWu14}
characterized the existence of a positive fixed point by the nonemptiness
of an explicit set $S_1$ and proved that the precision iteration from the
standard initialization converges monotonically to a distinguished fixed
point; this is the maximal fixed point of
Proposition~\ref{prop:max-root}(i). Uniqueness statements in the literature
concern smaller sets of admissible fixed points: a box below the
coefficients of a convex decomposition~\cite{MoallemiVanRoy}, positive
diagonal precisions on the factor graphs of turbo decoding~\cite{RVR}, and
a different factorization, for which the positive definite fixed point is
always unique~\cite{DuEtAl}. To our knowledge the following are new: the
count of \emph{all} positive solutions for the pairwise factorization
(Theorem~\ref{thm:unique-ws}); the converse, that for generic coefficients
a model that is not walk-summable has no positive solution or at least two
(Corollary~\ref{real:cor:generic}); the degree identity and the parity
(Theorem~\ref{thm:degree-identity}, Corollary~\ref{cor:parity}); the bound
on the Morse indices and the fact that all positive solutions but one are
saddles (Proposition~\ref{prop:morse-bound}); the unboundedness of the
number of positive solutions (Proposition~\ref{real:prop:unbounded}); the
compact components of the real curve (Proposition~\ref{prop:circles}); and
models in which every damping of the iteration diverges
(Remark~\ref{real:rem:damping}).

The statements below specify the graph class they use. The degree
identity and parity hold for loopless multigraphs with sums taken over
individual edges. For a simple graph we write $w_{ij}$, $f_{ij}$ for the
quantities of the edge $ij$ and $j\sim i$ for adjacency.

\subsection{Vertex variances and the reduced Bethe function}
\label{real:subsec:reduction}

Let $G$ be a loopless multigraph. Since $f_e=f_+(w_e^2V_iV_j)$ at a
positive solution, the positive solutions with $z=0$ correspond
bijectively, by $V\mapsto(V,f)$, to the vectors $V\in(0,\infty)^V$ with
\begin{equation}\label{real:eq:V}
D_iV_i=1+\sum_{e\ni i,\ e=ij}f_+\bigl(w_e^2V_iV_j\bigr)\qquad(i\in V).
\end{equation}
In particular $D_i>0$ for all $i$ whenever a positive solution exists. We
write $F(V)_i$ for the difference of the two sides of \eqref{real:eq:V} and
$J_V=\partial F/\partial V$. The \emph{reduced Bethe function} is
\begin{equation}\label{real:eq:phi}
\phi(V)=\frac12\sum_{i\in V}D_iV_i-\frac12\sum_{i\in V}\log V_i
-\frac12\sum_{e=ij\in E}g\bigl(w_e^2V_iV_j\bigr),\qquad
g(x)=\sqrt{1+4x}-\log\bigl(1+\sqrt{1+4x}\bigr),
\end{equation}
on $(0,\infty)^V$. It is obtained from the Bethe free energy of the
Gaussian model, written in the moments $V_i$ and the edge covariances, by
eliminating each edge covariance through its stationarity condition and
discarding a constant; we only use the following properties.

\begin{lemma}\label{real:lem:phi}
Let $G$ be a loopless multigraph with real coefficients.
\begin{enumerate}
\item[\textup{(a)}] $2V_i\,\partial\phi/\partial V_i=F(V)_i$. Hence the
positive solutions are the critical points of $\phi$, and
$J_V=2\diag(V)\Hess\phi(V)$ at a critical point.
\item[\textup{(b)}] At a positive solution put $u_e=f_e/(1+f_e)\in(0,1)$.
Then $2\diag(V)\Hess\phi(V)\diag(V)$ is the matrix $H_u$ with entries
\begin{equation}\label{real:eq:Hu}
(H_u)_{ii}=1+\sum_{e\ni i}\frac{u_e^2}{1-u_e^2},\qquad
(H_u)_{ij}=-\sum_{e=ij}\frac{u_e}{1-u_e^2}\quad(i\ne j).
\end{equation}
\item[\textup{(c)}] The Jacobian matrix of the $n+\abs E$ equations
\eqref{eq:system} with respect to $(V,f)$ at a positive solution is
invertible if and only if $\Hess\phi(V)$ is invertible.
\end{enumerate}
\end{lemma}

\begin{proof}
(a) Since $g'(x)=2/(1+\sqrt{1+4x})=f_+(x)/x$, we get
$\partial_{V_i}g(w_e^2V_iV_j)=f_+(w_e^2V_iV_j)/V_i$ and
$2V_i\partial_{V_i}\phi=D_iV_i-1-\sum_{e\ni i}f_+(w_e^2V_iV_j)$. Differentiating
$\partial_{V_i}\phi=F_i/(2V_i)$ at a zero of $F$ gives the last claim.

(b) Write $x_e=w_e^2V_iV_j=f_e(1+f_e)$. Since $f_+'(x)=1/\sqrt{1+4x}$ we
have $x_ef_+'(x_e)=f_e(1+f_e)/(1+2f_e)$. By (a) the matrix
$2\diag(V)\Hess\phi\,\diag(V)=J_V\diag(V)$ has diagonal entries
$D_iV_i-\sum_{e\ni i}x_ef_+'(x_e)=1+\sum_{e\ni i}f_e^2/(1+2f_e)$ and
off-diagonal entries $-\sum_{e=ij}f_e(1+f_e)/(1+2f_e)$. With $f_e=u_e/(1-u_e)$
these are the entries in \eqref{real:eq:Hu}.

(c) The derivative of the edge equations with respect to $f$ is
$\diag(1+2f_e)$, which is invertible at a positive solution. Its Schur
complement in the Jacobian matrix is the Jacobian matrix of the system
obtained by substituting $f_e=f_+(w_e^2V_iV_j)$ into the vertex equations,
that is $J_V$, and $J_V=2\diag(V)\Hess\phi$.
\end{proof}

A positive solution is \emph{nondegenerate} if $\Hess\phi(V)$ is invertible.
Its \emph{Morse index} $n_-$ and \emph{nullity} $n_0$ are the numbers of
negative and of zero eigenvalues of $\Hess\phi(V)$, equivalently of $H_u$.
When the diagonal is $D-z$, a positive solution is nondegenerate if and
only if $z$ is a local coordinate of $\cZ$ at this point (compare
Section~\ref{per:subsec:turning}). By
Proposition~\ref{prop:nondegenerate-bound} and Lemma~\ref{real:lem:phi}(c),
the number of nondegenerate positive solutions is at most $\LG$, for
arbitrary real coefficients.

We also use the cavity precisions $q_{i\to j}=(1+f_e)/V_i$ of
Lemma~\ref{lem:bp-correspondence}, indexed by the set $\vec E$ of oriented
edges (arcs). The \emph{precision update} is the map
\begin{equation}\label{real:eq:update}
\mathcal U(q)_{i\to j}=D_i-\sum_{\substack{e'\ni i,\ e'\ne e\\ e'=ik}}\frac{w_{e'}^2}{q_{k\to i}}
\qquad(e=ij)
\end{equation}
on $(0,\infty)^{\vec E}$, whose fixed points are the solutions of
\eqref{eq:cavity}. Its derivative at $q$ is the nonnegative matrix $C_q$
with entries $w_{e'}^2/q_{k\to i}^2$ in position $(i\to j,\,k\to i)$ for
$e'\ne e$, and zero otherwise. For weights $u\in[0,1)^E$ the weighted
nonbacktracking matrix is
\[
(B_uy)_{i\to j}=\sum_{\substack{e'\ni i,\ e'\ne e\\ e'=ik}}u_{e'}\,y_{k\to i}\qquad(e=ij,\ y\in\R^{\vec E}).
\]
At a positive solution, \eqref{eq:m-from-Vf} gives
$V_i\,(w_{e'}^2/q_{k\to i}^2)/V_k=w_{e'}^2V_iV_k/(1+f_{e'})^2=u_{e'}$, so that
\begin{equation}\label{real:eq:C-B}
C_q=\Psi^{-1}B_u\Psi,\qquad\Psi=\diag\bigl(V_i\bigr)_{(i\to j)\in\vec E} .
\end{equation}

\begin{lemma}\label{real:lem:bass}
Let $G$ be a finite simple graph and $u\in[0,1)^E$.
\begin{enumerate}
\item[\textup{(a)}] $\det(I-B_u)=\det H_u\prod_{e\in E}(1-u_e^2)$; this
holds, with the same definitions of $H_u$ and $B_u$, for all
$u\in(-1,1)^E$. If
$x\in\R^V$ and $y_{i\to j}=(x_i-u_{ij}x_j)/(1-u_{ij}^2)$, then
$\bigl((I-B_u)y\bigr)_{i\to j}=(H_ux)_i$.
\item[\textup{(b)}] $H_u$ is positive definite if and only if $\rho(B_u)<1$.
\item[\textup{(c)}] If every connected component of $G$ has cycle rank at
most one, then $\rho(B_u)<1$.
\end{enumerate}
\end{lemma}

\begin{proof}
(a) Let $S\colon\R^V\to\R^{\vec E}$, $(Sx)_{i\to j}=x_i$, and let
$(Jy)_{i\to j}=u_{ij}y_{j\to i}$. Then $(S^{\T}Jy)_i=\sum_{k\sim i}u_{ik}y_{k\to i}$,
so $B_u=SS^{\T}J-J$, and $P=I+J$ is invertible with
$P^{-1}=(I-J)\diag(1-u_e^2)^{-1}$ and $\det P=\prod_e(1-u_e^2)$. A direct
computation gives $I-S^{\T}JP^{-1}S=H_u$. Hence
$\det(I-B_u)=\det(P-SS^{\T}J)=\det P\,\det(I-S^{\T}JP^{-1}S)=\det P\,\det H_u$.
The vector $y$ in the statement is $P^{-1}Sx$, and
$(I-B_u)P^{-1}S=(P-SS^{\T}J)P^{-1}S=S(I-S^{\T}JP^{-1}S)=SH_u$. This is a
weighted form of Bass's determinant formula; compare
\cite{WatanabeFukumizuZeta}.

(b) If $\rho(B_u)<1$, then $\rho(B_{su})=s\rho(B_u)<1$ for $s\in[0,1]$, so
$\det H_{su}\ne0$ by (a); since $H_0=I$, the matrix $H_u$ is positive
definite. Conversely let $H_u$ be positive definite. It has nonpositive
off-diagonal entries, so $H_u^{-1}$ is entrywise nonnegative with positive
diagonal~\cite{BermanPlemmons}, and $x=H_u^{-1}\ones>0$. The vector $y$ of (a)
satisfies $y=\ones+B_uy$. We claim $y>0$. If $y_{i\to j}\le0$, then some
arc $k\to i$ with $k\ne j$ has $y_{k\to i}<0$, that is $x_k<u_{ik}x_i<x_i$.
Repeating the argument from the arc $k\to i$ produces an infinite strictly
decreasing sequence among the finitely many numbers $x_i$, which is
impossible. Hence $y>0$ and $B_uy=y-\ones\le cy$ with
$c=\max_a(1-1/y_a)<1$, so $\rho(B_u)\le c<1$.

(c) In a graph whose components have cycle rank at most one, a closed
walk without backtracking runs around the cycle of a component in one of
its two directions. Hence, in a suitable order of the arcs, $B_u$ is block
triangular, and each diagonal block is zero or a weighted cyclic
permutation matrix with spectral radius $(\prod_eu_e)^{1/\ell}<1$, the
product being over the $\ell$ edges of a cycle.
\end{proof}

\subsection{The maximal fixed point}\label{real:subsec:max}

Let $G$ be connected and simple, let $\tree$ be its universal cover, and
let $\tH$ be the operator of Section~\ref{subsec:cover} (with $z=0$). For a
finite set $\Lambda$ of vertices of $\tree$ let $\tH_\Lambda$ be the
principal submatrix of $\tH$ on $\Lambda$; we write $K[\Lambda]$ for the
restriction of a kernel $K$ on $\tree$ to $\Lambda\times\Lambda$. Quantities
indexed by vertices and arcs of $G$ are lifted to $\tree$ without change
of notation.

Let $(V,f)$ be a positive solution of \eqref{real:eq:V}, with cavity
precisions $q$. We define a symmetric kernel $K$ on the vertices of
$\tree$ by $K(x,x)=V_x$ and, if $x=x_0,x_1,\dots,x_k=y$ is the path from
$x$ to $y$,
\begin{equation}\label{real:eq:K}
K(x,y)=V_x\prod_{l=1}^{k}\Bigl(-\frac{w_{x_{l-1}x_l}}{q_{x_l\to x_{l-1}}}\Bigr).
\end{equation}
It is symmetric because $V_i/q_{j\to i}=V_iV_j/(1+f_{ij})$ is symmetric in
$i$ and $j$. It is the covariance kernel of the Gaussian Markov field on
$\tree$ obtained by lifting the beliefs of the fixed point; we use only
the following algebraic facts.

\begin{lemma}\label{real:lem:kernel}
Let $G$ be connected and simple and let $(V,f)$ be a positive solution of
\eqref{real:eq:V}.
\begin{enumerate}
\item[\textup{(a)}] For every finite set $\Lambda$ of vertices of $\tree$
the matrix $\tH_\Lambda$ is positive definite and
$K[\Lambda]-\tH_\Lambda^{-1}$ is positive semidefinite.
\item[\textup{(b)}] The limit $G^0(x,y)=\lim_{\Lambda\uparrow\tree}(\tH_\Lambda^{-1})(x,y)$
exists for all $x,y$. It does not depend on the exhaustion or on the
positive solution, and it is invariant under the deck transformations of
$\tree\to G$.
\item[\textup{(c)}] The kernel $H=K-G^0$ is positive semidefinite on every
finite set, and $\sum_y\tH(x,y)H(y,y')=0$ for all $x,y'$.
\end{enumerate}
\end{lemma}

\begin{proof}
By \eqref{eq:m-from-Vf} and the vertex equation, for adjacent $x,y$ in
$\tree$ and for every $x$,
\begin{equation}\label{real:eq:two-identities}
q_{x\to y}-\frac{w_{xy}^2}{q_{y\to x}}=\frac1{V_x},\qquad
\sum_{y\sim x}q_{x\to y}-\frac{d_x-1}{V_x}=\frac{1+\sum_{e\ni x}f_e}{V_x}=D_x .
\end{equation}
Let $F$ be a finite subtree of $\tree$ and let $P_F$ be the symmetric
matrix on $F$ with entries $w_{xy}$ for adjacent $x,y\in F$, with diagonal
entries $\sum_{y\in F,\,y\sim x}q_{x\to y}-(d_F(x)-1)/V_x$, and with zeros
elsewhere. We claim $P_FK[F]=I$. Fix $x,y\in F$. For a neighbor
$x'\in F$ of $x$ that is not on the path from $x$ to $y$,
\eqref{real:eq:K} gives $w_{xx'}K(x',y)=-(w_{xx'}^2/q_{x'\to x})K(x,y)$;
if $x\ne y$ and $x'$ is the neighbor of $x$ on the path, then
$w_{xx'}K(x',y)=-q_{x\to x'}K(x,y)$. By the first identity of
\eqref{real:eq:two-identities}, $(P_FK)(x,y)$ is therefore $K(x,y)$ times
$-(d_F(x)-1)/V_x+(d_F(x)-1)/V_x=0$ if $x\ne y$, and $K(x,x)$ times
$-(d_F(x)-1)/V_x+d_F(x)/V_x=1/V_x$ if $x=y$.

Next, $K[F]$ is positive definite. This is clear if $F$ is a point. If
$F'=F\cup\{y\}$ with $y$ a leaf of $F'$ attached to $x\in F$, then
$K(y,\cdot)=-(w_{xy}/q_{y\to x})K(x,\cdot)$ on $F$ and
$K(y,y)=V_y=(w_{xy}/q_{y\to x})^2V_x+1/q_{y\to x}$, because
$w_{xy}^2V_x/q_{y\to x}=f_{xy}$ and $(1+f_{xy})/q_{y\to x}=V_y$. So the Schur
complement of $K[F]$ in $K[F']$ is $1/q_{y\to x}>0$, and the claim follows
by induction.

(a) Choose a finite subtree $F$ that contains $\Lambda$ and all neighbors
of the vertices of $\Lambda$. By the second identity of
\eqref{real:eq:two-identities}, $P_F[\Lambda]=\tH_\Lambda$. Hence
$\tH_\Lambda=(K[F]^{-1})[\Lambda]$ is positive definite, and
$K[\Lambda]\succeq\bigl((K[F]^{-1})[\Lambda]\bigr)^{-1}=\tH_\Lambda^{-1}$ by
the inequality $(P^{-1})[\Lambda]\succeq(P[\Lambda])^{-1}$ for positive
definite $P$.

(b) By the same inequality, $(\tH_{\Lambda'}^{-1})[\Lambda]\succeq\tH_\Lambda^{-1}$
for $\Lambda\subseteq\Lambda'$. So for a finitely supported vector $a$ the
numbers $a^{\T}\tH_\Lambda^{-1}a$ are nondecreasing in
$\Lambda\supseteq\operatorname{supp}a$ and bounded by $a^{\T}Ka$; they
converge, and polarization gives the limit kernel. It is defined in terms
of $\tH$ alone, which gives the independence and the invariance.

(c) Positivity follows from (a). If $\Lambda$ contains $x$ and its
neighbors, then $\sum_y\tH(x,y)(\tH_\Lambda^{-1})(y,y')=\delta_{xy'}$; in
the limit, $\sum_y\tH(x,y)G^0(y,y')=\delta_{xy'}$. The identity $P_FK[F]=I$
at a vertex $x$ all of whose neighbors are in $F$ gives
$\sum_y\tH(x,y)K(y,y')=\delta_{xy'}$ in the same way.
\end{proof}

\begin{proposition}[Maximal fixed point]\label{prop:max-root}
Let $G$ be connected and simple with real coefficients, and suppose that
\eqref{real:eq:V} has a positive solution.
\begin{enumerate}
\item[\textup{(i)}] The iterates $\mathcal U^n(D)$ of
\eqref{real:eq:update}, started from $q_{i\to j}=D_i$, decrease to a fixed
point $q^{\max}$ of $\mathcal U$ in $(0,\infty)^{\vec E}$, and
$q\le q^{\max}$ entrywise for every fixed point $q$ of $\mathcal U$ in
$(0,\infty)^{\vec E}$. The solution $(V^{\max},f^{\max})$ of
\eqref{eq:system} corresponding to $q^{\max}$ is positive
\textup{\cite{SuWu14}}.
\item[\textup{(ii)}] $V^{\max}_i=G^0(x,x)=\lim_{\varepsilon\downarrow0}G_i(-\varepsilon)$
for $x$ over $i$. Every positive solution satisfies
$V_i-V_i^{\max}=H(x,x)\ge0$ for all $i$, and it equals
$(V^{\max},f^{\max})$ if equality holds for all $i$.
\item[\textup{(iii)}] $\rho(C_{q^{\max}})\le1$, and $\rho(C_q)>1$ for every
other fixed point $q$ of $\mathcal U$ in $(0,\infty)^{\vec E}$.
\end{enumerate}
\end{proposition}

We call $(V^{\max},f^{\max})$ the \emph{maximal fixed point}. It has the
largest cavity precisions and the smallest vertex variances among the
positive solutions. If the diagonal is $D-z$ with $z<\Sigm$, then by (ii)
it is the physical branch $\Vphys(z)$, and at $z=\Sigm$ it is the turning
point of Lemma~\ref{lem:turning}.

\begin{proof}
(i) The map $\mathcal U$ is monotone on $(0,\infty)^{\vec E}$, and a fixed
point $q$ in this set satisfies $q_{i\to j}\le D_i$. Hence
$\mathcal U(D)\le D$ and, by induction, $\mathcal U^{n+1}(D)\le\mathcal U^n(D)$
and $\mathcal U^n(D)\ge q$ for every such $q$. The limit $q^{\max}$ is a
fixed point with $q^{\max}\ge q>0$. If $q$ comes from a positive solution
$(V,f)$, then $1/V^{\max}_i=D_i-\sum_{k\sim i}w_{ik}^2/q^{\max}_{k\to i}\ge1/V_i>0$,
and the edge determinants
$q^{\max}_{i\to j}q^{\max}_{j\to i}-w_{ij}^2=q^{\max}_{i\to j}/V^{\max}_j$ are
positive, so $(V^{\max},f^{\max})$ is positive by
Lemma~\ref{lem:bp-correspondence}.

(ii) Let $x$ lie over $i$ and let $B_n(x)$ be the ball of radius $n$ in
$\tree$. The matrix $\tH_{B_n(x)}$ is positive definite by
Lemma~\ref{real:lem:kernel}(a). Eliminating its vertices from the leaves
towards $x$, one sees by induction on the depth that the pivot at a vertex
$y$ at distance $n-l$ from $x$, with parent $y'$, is
$\mathcal U^{l}(D)_{y\to y'}$. Hence
$1/(\tH_{B_n(x)}^{-1})(x,x)=D_i-\sum_{k\sim i}w_{ik}^2/\mathcal U^{n-1}(D)_{k\to i}$,
which converges to $1/V^{\max}_i$; so $V^{\max}_i=G^0(x,x)$. For
$\varepsilon>0$ the operator $\tH+\varepsilon$ is bounded below by
$\varepsilon$, since $\tH\ge0$ by Lemma~\ref{real:lem:kernel}(a), and
$G_i(-\varepsilon)=\sup_\Lambda\bigl((\tH_\Lambda+\varepsilon)^{-1}\bigr)(x,x)$ by the
variational formula
$\langle a,B^{-1}a\rangle=\sup_\psi\bigl(2\langle a,\psi\rangle-\langle\psi,B\psi\rangle\bigr)$
for a positive operator $B$, restricted to finitely supported $\psi$.
Both suprema are monotone, and exchanging them gives
$\lim_{\varepsilon\downarrow0}G_i(-\varepsilon)=G^0(x,x)$. The identity
$V_i-V_i^{\max}=H(x,x)$ is the definition of $H$, and $H(x,x)\ge0$ by
Lemma~\ref{real:lem:kernel}(c). A positive solution is determined by $V$.

(iii) Put $Y=q^{\max}$ and suppose $\rho=\rho(C_Y)>1$. By the
Perron--Frobenius theorem there is $X\ge0$, $X\ne0$, with $C_YX=\rho X$.
Let $\mu=\max_aX_a/Y_a$ and $0<\varepsilon<(\rho-1)/\mu$. From
$1/Y_a-1/(Y_a+\varepsilon X_a)\ge\varepsilon X_a/\bigl((1+\varepsilon\mu)Y_a^2\bigr)$
we get
$\mathcal U(Y+\varepsilon X)-Y\ge\varepsilon\rho X/(1+\varepsilon\mu)\ge\varepsilon X$.
Moreover $(C_YY)_{i\to j}=D_i-Y_{i\to j}$ and $X\le\mu Y$ give
$\rho X\le\mu(D-Y)$, hence $Z=Y+\varepsilon X\le D$. So $\mathcal U(Z)\ge Z$
and, by monotonicity, $Z\le\mathcal U^n(Z)\le\mathcal U^n(D)\to Y$, which
contradicts $Z\ge Y$, $Z\ne Y$. Now let $q\ne Y$ be another fixed point in
$(0,\infty)^{\vec E}$ and $\Delta=Y-q\ge0$. Then
\[
\Delta_{i\to j}=\sum_{k\sim i,\ k\ne j}\frac{w_{ik}^2\Delta_{k\to i}}{q_{k\to i}Y_{k\to i}},
\qquad(C_q\Delta)_{i\to j}=\sum_{k\sim i,\ k\ne j}\frac{w_{ik}^2\Delta_{k\to i}}{q_{k\to i}^2}
\ge(1+\eta)\Delta_{i\to j},
\]
where $\eta>0$ is the minimum of $\Delta_a/q_a$ over the arcs with
$\Delta_a>0$. Since $\Delta\ge0$ is nonzero, $\rho(C_q)\ge1+\eta$.
\end{proof}

Part (iii) implies that only the maximal fixed point can attract the
precision iteration or any of its dampings: the derivative
$(1-\alpha)I+\alpha C_q$ of the damped update
$(1-\alpha)\mathrm{id}+\alpha\,\mathcal U$, $0<\alpha\le1$, has the
eigenvalue $1-\alpha+\alpha\rho(C_q)>1$ at every other fixed point.

\begin{remark}\label{real:rem:gibbs}
The kernel $H$ has a probabilistic meaning, which we state without proof
and do not use. Let $X$ be the centered Gaussian field on $\tree$ with
covariance $K$; it is a Gibbs measure for the formal Hamiltonian
$\tfrac12\langle\psi,\tH\psi\rangle$. Its conditional expectation given
the tail $\sigma$-field is a Gaussian field $\Psi$ with covariance $H$ and
$\tH\Psi=0$ almost surely, and $X-\Psi$ is independent of $\Psi$ with
covariance $G^0$; this follows from backward martingale convergence. Hence
a positive solution is the maximal fixed point if and only if the tail
$\sigma$-field of $X$ is trivial, if and only if the Gibbs measure of $X$
is extremal. For a non-maximal solution the fraction of the variance at a
vertex over $i$ that is determined by the boundary at infinity is
$1-V_i^{\max}/V_i$.
\end{remark}

\subsection{Walk-summable models}\label{real:subsec:ws}

\begin{theorem}[Uniqueness]\label{thm:unique-ws}
Let $G$ be a finite loopless multigraph with real coefficients such that
$D-\abs A$ is positive definite. Then \eqref{eq:system} with $z=0$ has
exactly one positive solution. On each connected component it is the
physical branch at $z=0$. Equivalently, for connected $G$ and arbitrary
real coefficients, the physical branch is the only positive solution over
each $z<\sigm$.
\end{theorem}

\begin{proof}
The system decouples over the connected components, so let $G$ be
connected.

\emph{Existence.} Let $\bar u=\bar u^->0$ be the eigenvector of $D-\abs A$
for $\sigm>0$, lifted to $\tree$. Since $\tH$ is unitarily equivalent to
the operator with weights $-\abs{w_e}$, for finitely supported $\psi$ on
$\tree$ we have, after this equivalence,
\[
\langle\psi,(\tH-\sigm)\psi\rangle
=\sum_{\tilde e=xy}\abs{w_{\tilde e}}\,\bar u_x\bar u_y
\Bigl(\frac{\psi_x}{\bar u_x}-\frac{\psi_y}{\bar u_y}\Bigr)^2\ge0,
\]
the sum being over the edges of $\tree$; this is verified by expanding the
square and using $\sum_{y\sim x}\abs{w_{xy}}\bar u_y=(D_x-\sigm)\bar u_x$.
Hence $\Sigm\ge\sigm>0$, and the physical branch at $z=0$ is a positive
solution by Proposition~\ref{prop:physical-branch}.

\emph{Uniqueness.} Let $(V,f)$ be a positive solution and
$H$ the kernel of Lemma~\ref{real:lem:kernel}. A positive semidefinite
kernel on a countable set is a Gram kernel: there are vectors $\Psi_x$ in
a real Hilbert space with $\langle\Psi_x,\Psi_y\rangle=H(x,y)$. By
Lemma~\ref{real:lem:kernel}(c),
\[
\Bigl\lVert\sum_y\tH(x,y)\Psi_y\Bigr\rVert^2
=\sum_{y,y'}\tH(x,y)H(y,y')\tH(y',x)=0,
\qquad\text{that is}\qquad
D_x\Psi_x=-\sum_{y\sim x}w_{xy}\Psi_y .
\]
By Lemma~\ref{real:lem:kernel}(b) and \eqref{real:eq:K}, $H(x,x)$ depends
only on the vertex $i$ of $G$ under $x$; put $\sigma_i=\norm{\Psi_x}=H(x,x)^{1/2}$.
Each edge $e=ij$ incident to $i$ lifts to exactly one edge at $x$.
The triangle inequality, applied to these lifts individually, gives
$D_i\sigma_i\le\sum_{e=ij,\,e\ni i}\abs{w_e}\sigma_j$, that is
$(D-\abs A)\sigma\le0$ entrywise. Since $\sigma\ge0$, we get
$\sigma^{\T}(D-\abs A)\sigma\le0$ and therefore $\sigma=0$. By
Proposition~\ref{prop:max-root}(ii), $V=V^{\max}$, and the solution is the
maximal fixed point.
\end{proof}

For vector-valued variables the comparison matrix is replaced by a cone
condition, and the corresponding uniqueness theorem
(Theorem~\ref{thm:block-cone-unique}) contains
Theorem~\ref{thm:unique-ws}, with a different proof of existence.

\begin{remark}\label{real:rem:literature}
The convergence of the variances for walk-summable models~\cite{MJW} and
the convergence to the maximal fixed point~\cite{SuWu14} show that the
iteration from the standard initialization has one possible limit. They
do not exclude other positive fixed points, and Theorem~\ref{thm:unique-ws}
does. The theorem concerns the factorization of the Gaussian density into
vertex and edge factors used here. For the factorization
of~\cite{DuEtAl}, coming from linear observation models, a positive
definite fixed point is unique for all coefficients, and no dichotomy as
in Corollary~\ref{real:cor:generic} below occurs.
\end{remark}

\subsection{Morse indices}\label{real:subsec:morse}

A graph is a \emph{pseudoforest} if each of its connected components has
cycle rank at most one. The \emph{pseudoforest deletion number}
$\taups(G)$ is the least number of vertices whose removal leaves a
pseudoforest. For example $\taups(K_4)=1$, $\taups(K_{2,m})=1$ for
$m\ge3$, every wheel has $\taups=1$, and $\taups(K_5)=2$; the cycle rank is
not bounded on the graphs with $\taups=1$.

\begin{proposition}[Morse indices]\label{prop:morse-bound}
Let $G$ be a finite simple graph with real coefficients.
\begin{enumerate}
\item[\textup{(i)}] For every $u\in[0,1)^E$ the matrix $H_u$ of
\eqref{real:eq:Hu} has at most $\taups(G)$ nonpositive eigenvalues, counted
with multiplicity. Hence $n_-+n_0\le\taups(G)$ at every positive solution.
\item[\textup{(ii)}] At a positive solution, $\Hess\phi$ is positive
definite if and only if $\rho(C_q)<1$.
\item[\textup{(iii)}] Let $G$ be connected. If $\Hess\phi$ is positive
semidefinite and singular at a positive solution, then the corresponding
$q$ is the only fixed point of $\mathcal U$ in $(0,\infty)^{\vec E}$; in
particular there is no other positive solution.
\item[\textup{(iv)}] Let $G$ be connected. At the maximal fixed point
$\Hess\phi$ is positive semidefinite, and it is positive definite if the
maximal fixed point is nondegenerate. At every other positive solution
$\Hess\phi$ has a negative eigenvalue.
\end{enumerate}
\end{proposition}

So the reduced Bethe function has at most one local minimum, the maximal
fixed point, and all other positive solutions are saddles. The implication
from $\rho(C_q)<1$ to positive definiteness in (ii) is, in a different
parametrization, a theorem of Cseke and Heskes~\cite{CsekeHeskes}. Part (i)
does not use the fixed point equations.

\begin{proof}
(i) Let $S\subseteq V$ with $\abs S=\taups(G)$ be such that $G-S$ is a
pseudoforest, and $W=V\setminus S$. The principal submatrix $H_u[W]$ is
the matrix \eqref{real:eq:Hu} of the graph $G-S$, which is positive
definite by Lemma~\ref{real:lem:bass}(b),(c), plus the nonnegative
diagonal matrix with entries $\sum_{j\in S,\,j\sim i}u_{ij}^2/(1-u_{ij}^2)$.
So $H_u$ is positive definite on $\R^W$, and a subspace on which $H_u$ is
nonpositive has dimension at most $\abs S$. The second statement follows
from Lemma~\ref{real:lem:phi}(b).

(ii) By Lemma~\ref{real:lem:phi}(b), \eqref{real:eq:C-B} and
Lemma~\ref{real:lem:bass}(b).

(iii) Let $H_u\succeq0$ be singular. Then $G$ is not a pseudoforest, by
Lemma~\ref{real:lem:bass}(b),(c). Let $G'$ be the subgraph obtained by
removing vertices of degree one repeatedly; it is connected, has minimum
degree at least two and is not a cycle, so the nonbacktracking matrix of
$G'$ is irreducible. If $\ell$ is a leaf attached to $i$ by an edge with
weight $u$, then $(H_u)_{\ell\ell}=1/(1-u^2)$ and
$(H_u)_{\ell i}=-u/(1-u^2)$, and the Schur complement of the entry
$(\ell,\ell)$ in $H_u$ is the matrix \eqref{real:eq:Hu} of $G-\ell$. By
the additivity of inertia~\cite{Haynsworth}, the matrix $H'$ of $G'$ (for
the restriction of $u$) is again positive semidefinite and singular. It
is irreducible with nonpositive off-diagonal entries, so its kernel is
spanned by a positive vector $x$. Let $y_{i\to j}=(x_i-u_{ij}x_j)/(1-u_{ij}^2)$
on the arcs of $G'$ and let $B'$ be the weighted nonbacktracking matrix
of $G'$. By Lemma~\ref{real:lem:bass}(a), $B'y=y$. As in the proof of
Lemma~\ref{real:lem:bass}(b), a negative entry of $y$ would give an
infinite strictly decreasing sequence among the $x_i$, so $y\ge0$; and
$y\ne0$, because $y=0$ gives $x_i=u_{ij}^2x_i$. An irreducible nonnegative
matrix with a nonnegative eigenvector for the eigenvalue $1$ has spectral
radius $1$ and a positive left eigenvector. Let $C'$ be the restriction
of $C_q$ to the arcs of $G'$; by \eqref{real:eq:C-B} it is similar to $B'$
by a positive diagonal matrix. So there is $l>0$ with $l^{\T}C'=l^{\T}$.

Let $\tilde q=q+\delta$ be a fixed point of $\mathcal U$ in
$(0,\infty)^{\vec E}$. On the arcs outside $G'$ that point towards $G'$,
the fixed point equation determines the precisions recursively from the
leaves, so $\delta=0$ there. For an arc $a$ of $G'$ the equation
$\mathcal U(\tilde q)_a-\mathcal U(q)_a=\delta_a$ and
$1/q_b-1/\tilde q_b=\delta_b/q_b^2-\delta_b^2/(q_b^2\tilde q_b)$ give
\[
\delta_a=\sum_{b}(C')_{ab}\,\delta_b-\sum_b(C')_{ab}\,\frac{\delta_b^2}{\tilde q_b},
\]
the sums being over the arcs $b$ of $G'$. Pairing with $l$ gives
$\sum_bl_b\delta_b^2/\tilde q_b=0$, hence $\delta=0$ on $G'$. The
remaining arcs point away from $G'$, and their precisions are determined
by the others. So $\tilde q=q$.

(iv) At the maximal fixed point $\rho(B_u)=\rho(C_q)\le1$ by
Proposition~\ref{prop:max-root}(iii), so $\rho(B_{su})<1$ and
$H_{su}\succ0$ for $s<1$, and $H_u\succeq0$ in the limit. Let $(V,f)$ be
another positive solution. Then $\rho(C_q)>1$, so $\Hess\phi$ is not
positive definite by (ii). If it were positive semidefinite, it would be
singular, and by (iii) there would be only one positive solution.
\end{proof}

\begin{example}[The bound is attained]\label{real:ex:diamond-chain}
Let the diamond be $K_4$ minus the edge $23$, with $D=(73/15,73/15,1,1)$
and all weights $1$. It has $\taups=1$ and exactly two positive solutions
(Theorem~\ref{real:thm:packing} below), of Morse indices $0$ and $1$; in the
cavity precisions on the arcs between the vertices $0$ and $1$ they are
given by $t=5/3$ and $t=6/5$. Join $m$ copies in a chain by edges of a
small weight $\varepsilon$ between vertex $3$ of one copy and vertex $2$
of the next. The resulting graph has $4m$ vertices, maximal degree $3$
and $\taups=m$. The equations \eqref{real:eq:V} and the function $\phi$
depend analytically on $\varepsilon^2$, also at $\varepsilon=0$, where
they are those of $m$ disjoint diamonds. There the $2^m$ choices of a
solution in every copy are positive solutions, and they are nondegenerate
because the Hessian of $\phi$ is the direct sum of the Hessians of the
copies. By the implicit function theorem, for every $m$ there is
$\varepsilon_m>0$ such that for $0<\varepsilon<\varepsilon_m$ each of these
choices continues to a nondegenerate positive solution, whose Morse index
is the number of copies in which the second solution was chosen. So for
such $\varepsilon$ there are at least $\binom mk$ positive solutions of
index $k$ for $0\le k\le m$, and the bound of
Proposition~\ref{prop:morse-bound}(i) is attained. In fact every
$0<\varepsilon\le1/1000$ works for all $m$: each equation involves at most four
terms, so the perturbation of the solutions and of the Hessian is bounded
in the maximum norm by a constant times $\varepsilon^2$ that does not
depend on $m$ or on the choice of solutions
(Theorem~\ref{app:many:thm-chain}).
\end{example}

\subsection{The degree identity and the parity of the number of solutions}
\label{real:subsec:degree}

We say that the real coefficients satisfy \emph{condition \textup{(G)}} if
$\lambda_{\min}(M[S])\ne0$ for every nonempty $S\subseteq V$ that induces a
connected subgraph of $G$, where $M=D-\abs A$ and $M[S]$ is its principal
submatrix. For $S=\{i\}$ this says $D_i\ne0$. Condition (G) excludes
finitely many real hypersurfaces in the space of the potentials $D$; for
the diagonal $D-z$ it fails for at most $2^n-1$ values of $z$.

\begin{theorem}[Degree identity]\label{thm:degree-identity}
Let $G$ be a finite loopless multigraph with real coefficients satisfying
condition \textup{(G)}. Then the set $Z\subset(0,\infty)^V$ of solutions of
\eqref{real:eq:V} is compact, and for every bounded open set $\Omega$ with
$Z\subset\Omega$ and $\overline\Omega\subset(0,\infty)^V$,
\[
\deg(\nabla\phi,\Omega,0)=\begin{cases}1&\text{if $D-\abs A$ is positive definite,}\\
0&\text{otherwise.}\end{cases}
\]
If $Z$ is finite, the left-hand side is the sum of the local indices of
$\nabla\phi$ at the positive solutions, and the index of a nondegenerate
positive solution is $(-1)^{n_-}$.
\end{theorem}

Here $\deg$ is the Brouwer degree; we use its homotopy invariance and
additivity~\cite{Milnor}.

\begin{proof}
Put $v_i=\sqrt{V_i}$ and $\vartheta(x)=\sqrt x-f_+(x)$. Then
$\vartheta(0)=0$ and $0\le\vartheta(x)<\tfrac12$ for $x\ge0$, and
\eqref{real:eq:V} reads
$v_i(Mv)_i=1-\sum_{e\ni i,\,e=ij}\vartheta(w_e^2v_i^2v_j^2)$
(compare \eqref{per:eq:sqrt-form}). For $s\in[0,1]$ let
$\Gamma_s\colon(0,\infty)^V\to\R^V$,
\[
\Gamma_s(v)_i=v_i(Mv)_i-1+s\sum_{e\ni i,\,e=ij}\vartheta\bigl(w_e^2v_i^2v_j^2\bigr).
\]

\emph{Compactness.} We show that the zeros of all $\Gamma_s$ lie in a
compact subset of $(0,\infty)^V$. Let $\Gamma_{s_k}(v^{(k)})=0$. If
$\norm{v^{(k)}}\to\infty$, let $\nu\ge0$, $\norm\nu=1$, be a limit point of
$v^{(k)}/\norm{v^{(k)}}$. The terms other than $v_i(Mv)_i$ are bounded, so
$\nu_i(M\nu)_i=0$ for all $i$. With $S=\operatorname{supp}\nu$ this gives
$M[S]\nu_S=0$ and $\nu_S>0$. On a connected component $S'$ of the subgraph
induced by $S$, the matrix $M[S']$ is irreducible with nonpositive
off-diagonal entries and has a positive vector in its kernel, so
$\lambda_{\min}(M[S'])=0$, which contradicts (G). If $v^{(k)}$ is bounded
and $v^{(k)}_i\to0$ for some $i$, then $v_i(Mv)_i\to0$ and
$\vartheta(w_e^2v_i^2v_j^2)\to0$ on each incident edge, so
$\Gamma_{s_k}(v^{(k)})_i\to-1$, a
contradiction. In particular $Z$ is compact.

\emph{The case $s=0$.} A zero of $\Gamma_0$ is a vector $v>0$ with
$Mv=1/v>0$. A symmetric matrix with nonpositive off-diagonal entries that
maps a positive vector to a positive vector is positive
definite~\cite{BermanPlemmons}. So $\Gamma_0$ has no zero unless $M$ is
positive definite. If it is, $\Gamma_0=\diag(v)\nabla\Xi$ with
$\Xi(v)=\tfrac12v^{\T}Mv-\sum_i\log v_i$, which is strictly convex on
$(0,\infty)^V$ and tends to $+\infty$ at the boundary and at infinity. So
$\Gamma_0$ has exactly one zero, and its derivative there,
$\diag(v)\Hess\Xi$, has positive determinant.

\emph{The case $s=1$.} By Lemma~\ref{real:lem:phi}(a),
$\Gamma_1(v)=2\diag(V)\nabla\phi(V)$ with $V_i=v_i^2$. The change of
variables is an orientation preserving diffeomorphism of $(0,\infty)^V$,
and multiplication of a vector field by a positive diagonal matrix does
not change its degree. So the degree of $\Gamma_1$ on the preimage of
$\Omega$ is $\deg(\nabla\phi,\Omega,0)$, and by additivity it does not
depend on $\Omega$.

By the compactness statement there is a box $[\epsilon,R]^V$ whose
interior contains all zeros of all $\Gamma_s$, and homotopy invariance on
this box gives $\deg\Gamma_1=\deg\Gamma_0$. The local index of
$\nabla\phi$ at a nondegenerate zero is the sign of $\det\Hess\phi$.
\end{proof}

\begin{corollary}[Parity]\label{cor:parity}
Let $G$ be a finite loopless multigraph with real coefficients satisfying
\textup{(G)}, and suppose that all positive solutions are nondegenerate.
Then their number $N$ is finite, and if $c_k$ denotes the number of
positive solutions of Morse index $k$,
\[
\sum_k(-1)^kc_k=[\,D-\abs A\succ0\,],\qquad N\equiv[\,D-\abs A\succ0\,]\pmod 2 .
\]
In particular there is exactly one positive solution if and only if the
model is walk-summable. If $G$ is simple and connected, the model is not
walk-summable and $N>0$, then $c_0=1$, the solution of index $0$ is the
maximal fixed point, and there is a positive solution of odd Morse index.
\end{corollary}

\begin{proof}
A compact set of isolated points is finite, so
Theorem~\ref{thm:degree-identity} gives the first identity, and
$N=\sum_kc_k\equiv\sum_k(-1)^kc_k$. If $N=1$ then $N$ is odd and the model
is walk-summable; the converse is Theorem~\ref{thm:unique-ws}, which needs
neither (G) nor nondegeneracy. In the last statement, the maximal fixed
point exists and has index $0$, and it is the only solution of index $0$,
by Proposition~\ref{prop:morse-bound}(iv); then
$\sum_{k\ \mathrm{odd}}c_k=\sum_{k\ \mathrm{even}}c_k\ge1$.
\end{proof}

Without (G) the last statement can fail for special potentials: as $D$
approaches a hypersurface excluded by (G), a saddle may move to infinity.

\begin{corollary}[Generic coefficients]\label{real:cor:generic}
Let $G$ be a finite loopless multigraph. For generic real coefficients, all
positive solutions are nondegenerate, their number $N$ satisfies
$N\le\LG$, and $N$ is odd if the model is walk-summable and even
otherwise. In particular, for generic real coefficients the positive
solution is unique if and only if the model is walk-summable, and a model
that is not walk-summable has no positive solution or at least two.
\end{corollary}

\begin{proof}
For generic coefficients all admissible solutions of \eqref{eq:system}
with $z=0$ are nondegenerate, and there are $\LG$ of them
(Theorem~\ref{thm:generic-degree}); by Lemma~\ref{real:lem:phi}(c) the
positive solutions are then nondegenerate in the sense used here.
Condition (G) holds as soon as $\det\bigl((D+A^s)[S]\bigr)\ne0$ for all
signings $s$ and all nonempty $S$, which is a nonempty Zariski open
condition. Now apply Corollary~\ref{cor:parity}.
\end{proof}

When the diagonal is $D-z$ and $G$ is connected with $\beta\ge2$, we
obtain the following picture for generic real coefficients. Let $N(z)$
be the number of positive solutions over $z$. By
Corollary~\ref{cor:arc-circles}, $\Xp$ is a smooth curve, so $z$ has
finitely many critical points on it, and outside the finitely many
critical values and the values excluded by (G) all positive solutions
over $z$ are nondegenerate. Then $N(z)=1$ for $z<\sigm$; $N(z)$ is even
and at least $2$ for $\sigm<z<\Sigm$ (Corollary~\ref{cor:two-roots}); and
$N(z)=0$ for $z>\Sigm$ (Lemma~\ref{lem:domain}). For simple graphs we can
characterize when the gap can contain more than two roots and when
compact components can occur.

\begin{theorem}[Packing and positive fixed points]\label{real:thm:packing}
Let $G$ be a connected finite simple graph. Write $p_U(G)$ for the
largest number of pairwise vertex-disjoint connected subgraphs of cycle
rank at least two.
\begin{enumerate}
\item[\textup{(a)}] If $p_U(G)\le1$, there are at most two positive
solutions at each energy for any real coefficients, and $\Xp$ has no compact connected
component. If $G$ is a pseudoforest there is at most one solution at each
$z$. If there are two solutions at a fixed $z$, they are
nondegenerate, of Morse indices $0$ and $1$. If \textup{(G)} holds and
the model is not walk-summable, there is no positive solution, or one
degenerate positive solution, or two.
\item[\textup{(b)}] Suppose $p_U(G)\le1$ and $\beta\ge2$. As $z$ varies,
$N(z)=1$ for $z\le\sigm$, $N(z)=2$ for $\sigm<z<\Sigm$, $N(\Sigm)=1$ and
$N(z)=0$ for $z>\Sigm$. The set $\Xp$ is connected: it is the union of the
physical branch over $(-\infty,\Sigm)$, the turning point over $\Sigm$,
and the graph of a real analytic map $U$ on $(\sigm,\Sigm)$ that is
unbounded as $z\downarrow\sigm$.
\item[\textup{(c)}] If $p_U(G)\ge2$, there is a nonempty open set of
real coefficients on $G$ for which $\Xp$ has a circle of regular points
and there are at least four nondegenerate positive solutions at one
energy.
\end{enumerate}
\end{theorem}

\begin{proof}
(a) For $0<u_e<1$, $H_u$ has nonpositive off-diagonal entries and its
smallest eigenvalue has a positive eigenvector by Perron--Frobenius.
Suppose the second eigenvalue $\lambda_2$ is nonpositive. An eigenvector
$x=x_+-x_-$ orthogonal to the first eigenvector has both signs, with
$x_+,x_-\ge0$ of disjoint support. The eigenvalue equation gives
\[
 x_\pm^{\T}H_ux_\pm
 =\lambda_2\norm{x_\pm}^2+x_+^{\T}H_ux_-\le0.
\]
For a set $S$ of vertices, $H_u[S]$ is the matrix for $G[S]$ plus
nonnegative diagonal entries from the edges leaving $S$. It is positive
definite if $G[S]$ is a pseudoforest, by Lemma~\ref{real:lem:bass}(b),(c).
Thus each sign support contains a connected subgraph of cycle rank at
least two, contradicting $p_U(G)\le1$. So $H_u$ has at most one
nonpositive eigenvalue.

If $G$ is a pseudoforest, then $\rho(C_q)=\rho(B_u)<1$ at every positive
solution by Lemma~\ref{real:lem:bass}(c), so every positive solution is
the maximal fixed point by Proposition~\ref{prop:max-root}(iii). In this
case the spectral bottom is $\sigm$, by switching signs on a tree or by
Proposition~\ref{prop:unicyclic}(a). The unique branch below $\sigm$ is
all of $\Xp$. A positive solution at $\sigm$ would have positive definite
Hessian and continue beyond the spectral bottom, contrary to
Lemma~\ref{lem:domain}. For a single vertex the solution is directly
$V=(D-z)^{-1}$ for $z<D$. In the remaining case $\beta\ge2$ the assertion
about compact components follows from (b).

By the eigenvalue bound a degenerate positive solution has $n_-=0$ and
$n_0=1$, so it is the only positive solution by
Proposition~\ref{prop:morse-bound}(iii). Suppose that there are at
least two positive solutions. Then all are nondegenerate, hence isolated,
and they are finite in number because the solution set is semialgebraic.
Their Morse indices are $0$ or $1$, and $c_0=1$ by
Proposition~\ref{prop:morse-bound}(iv). If there were at least three, we
could perturb $D$ slightly so that (G) holds; three nondegenerate
solutions persist with their indices, so the perturbed model again has no
degenerate solution, $c_0=1$ and $c_1\ge2$, and $c_0-c_1<0$ contradicts
Corollary~\ref{cor:parity}. Under (G) and without walk-summability,
Corollary~\ref{cor:parity} gives $c_0=c_1$ when all solutions are
nondegenerate.

(b) For $z<\sigm$ this is Theorem~\ref{thm:unique-ws}, and for $z>\Sigm$
it is Lemma~\ref{lem:domain}. For $\sigm<z<\Sigm$ there are at least two
positive solutions by Corollary~\ref{cor:two-roots} and at most two by
(a); both are nondegenerate, so by the implicit function theorem they
form two real analytic branches, the physical one and a second one, $U$.
Over $\sigm<\Sigm$ (Corollary~\ref{cor:gap}) the physical branch is a
positive solution. A second positive solution over $\sigm$ would be
nondegenerate and would continue to $z<\sigm$, which is impossible; so
$N(\sigm)=1$. The map $U$ is unbounded as $z\downarrow\sigm$: since
$V_i\ge1/(D_i-z)$ on $\Xp$, a bounded limit point would be a positive
solution over $\sigm$, hence the physical one, and the physical branch is
nondegenerate there (Lemma~\ref{lem:turning}(i)), so $U$ would coincide
with it nearby. Two positive
solutions over $\Sigm$ would be nondegenerate and would continue to
$z>\Sigm$; so $N(\Sigm)\le1$. By Lemma~\ref{lem:blowup}(b) positive solutions
stay bounded as $z\uparrow\Sigm$, so both branches converge to the unique
positive solution over $\Sigm$, which is the turning point of
Lemma~\ref{lem:turning}. Hence $\Xp$ is connected and not compact.

(c) The construction is proved in Appendix~\ref{app:bridge:subsec-packing}.
\end{proof}

For a disconnected simple graph with $p_U(G)\le1$, at most one connected
component has cycle rank at least two. All other components supply a
unique branch on an open half-line of energies bounded above. The common
range of these branches only truncates the response of the remaining
component in $z$, so the bound of two roots and the absence of compact
components also hold. The criterion for circles on a connected simple
graph is the Gaussian counterpart of the packing classification for
binary belief propagation in \cite{WatanabeResponse}.

\begin{example}\label{real:ex:K4}
Let $G=K_4$ with $D=0$ and all weights $1$, so that $\sigm=-3$ and
$\Sigm=-2\sqrt2$. Since $p_U(K_4)=1$, Theorem~\ref{real:thm:packing}
applies. With $a=-z$, the two positive solutions for $2\sqrt2<a<3$ are
uniform, with cavity precisions $t_\pm=\tfrac12(a\pm\sqrt{a^2-8})$ solving
$t=a-2/t$, and vertex variances $1/(t-1/t)$; the solution $t_-$ is
positive exactly for $a<3$, and its vertex precision $t_--1/t_-$ vanishes
at $a=3$, where it leaves $\Xp$ through the Perron pole. For $a=43/15$
the variances are $15/16$ and $30/11$, and
$\rho(C_q)=2/t^2$ equals $18/25$ and $25/18$.
\end{example}

For $\taups\ge2$ the degree identity gives no bound on the number of
positive solutions, and indeed there is none.

\begin{proposition}[The number of positive solutions is unbounded]
\label{real:prop:unbounded}
For every $m\ge1$ there are a connected finite simple graph $G$ with
$\taups(G)=2$ and a nonempty open set of real coefficients on $G$ such
that every model in this set has at least $2m+2$ positive solutions.
Outside a closed subset with empty interior, all positive solutions of
these models are nondegenerate, condition \textup{(G)} holds, and
$(c_0,c_1,c_2)=(1,c_1,c_1-1)$ with $c_1\ge m+1$.
\end{proposition}

\begin{proof}[Outline of the proof]
The graphs consist of two vertices $a$, $b$ joined by an edge, with
triangles attached at $a$ and at $b$. Removing $a$ and $b$ leaves a
forest, and $\taups=2$. Eliminating the triangles, which can be done in
closed form, reduces \eqref{real:eq:V} to two equations in $(V_a,V_b)$,
and then to the zeros of one real function $h$ of $s=V_a$ on an interval,
whose simple zeros correspond to the nondegenerate positive solutions;
the sign of $h'$ determines whether the Morse index is odd. If the
triangles at $a$ have weights on $m+3$ geometrically separated scales,
elementary inequalities show that $h$ changes sign at least $2m+2$ times.
The form of $(c_0,c_1,c_2)$ follows from
Proposition~\ref{prop:morse-bound} and Corollary~\ref{cor:parity}. The
complete proof, with the coefficients and the sign estimates, is given in
Appendix~\ref{app:many} (Theorem~\ref{app:many:thm-unbounded}): the
reduced function alternates signs at
$s=100^{2k}$ and $s=100^{2k+1}$, with $K=m+2$ and indices $k=1,\dots,K$
for the first sequence and $k=1,\dots,K-1$ for the second. The number of
vertices grows exponentially with $m$.
\end{proof}

Small instances can be certified by computer. By a computer-assisted
proof (exact rational interval arithmetic for the function $h$, and an
interval Krawczyk test~\cite{Krawczyk} with interval inertia for each
solution of the full system), a model with rational coefficients on a
graph of this type with $27$ vertices and $34$ edges, with cycles of
lengths $3,3,4,5$ at $a$ and $4,5,4,5$ at $b$, has exactly $6$ positive
solutions, with $(c_0,c_1,c_2)=(1,3,2)$, and a model with $68$ vertices
and $100$ edges ($17$ triangles at $a$, $16$ at $b$) has exactly $10$,
with $(c_0,c_1,c_2)=(1,5,4)$
(Proposition~\ref{app:many:prop-certified}). In particular the bound
$c_k\le\binom{\taups}{k}$, which holds for $\taups\le1$ and is attained in
Example~\ref{real:ex:diamond-chain}, fails for $\taups=2$. The bound
$N\le\LG$ of Proposition~\ref{prop:nondegenerate-bound} is not
contradicted, since $\LG$ grows exponentially with the cycle rank
(Proposition~\ref{prop:L-bounds}); for the graph with $27$ vertices,
$\LG=2^8+15^2=481$.

\subsection{Compact components of the real curve}\label{real:subsec:circles}

We return to the family with diagonal $D-z$. Let $G$ be connected with
$\beta\ge2$. For generic real coefficients, $\Xp$ is a smooth curve with
one non-compact component, the arc $\mathcal A$ from the physical branch
through the turning point to the Perron pole, and finitely many
\emph{circles} (Corollary~\ref{cor:arc-circles}). The arc contains the
physical branch and, by Corollary~\ref{cor:two-roots}, at least one more
point over each $z\in(\sigm,\Sigm)$.

\begin{proposition}[Circles]\label{prop:circles}
Let $G$ be a connected loopless multigraph with $\beta\ge2$ and generic
real coefficients, and let $\mathcal C$ be a circle of $\Xp$, with
$z(\mathcal C)=[a,b]$.
\begin{enumerate}
\item[\textup{(a)}] The number of circles of $\Xp$ is at most the bound
of Proposition~\ref{prop:genus-bound} for the genus of the spectral
curve.
\item[\textup{(b)}] $\sigm\le a<b\le\Sigm$, and $N(z)\ge4$ for all
$z\in(a,b)$.
\item[\textup{(c)}] If $G$ is simple, then $p_U(G)\ge2$, $\taups(G)\ge2$,
and $G$ has at least eight vertices. If $G$ is simple and $\taups(G)=2$,
then for all but
finitely many $z\in(a,b)$ the positive solutions over $z$ are
nondegenerate with $(c_0,c_1,c_2)=(1,c_1,c_1-1)$, $c_1\ge2$, and the
points of $\mathcal C$ over $z$ have Morse index $1$ or $2$.
\item[\textup{(d)}] If $G$ consists of two vertices joined by $k\ge3$
parallel edges, then $N(z)\le2$ for all $z$ and all real coefficients,
and $\Xp$ has no circle.
\end{enumerate}
\end{proposition}

\begin{proof}
(a) is proved in Section~\ref{sec:genus}: a circle is a connected
component of the real points of the smooth projective model of the
spectral curve other than the one containing the arc, and Harnack's
inequality applies.

(b) We have $b\le\Sigm$ by Lemma~\ref{lem:domain}. For $z<\sigm$ the only
positive solution is the physical branch (Theorem~\ref{thm:unique-ws}),
which lies on $\mathcal A$, so $a\ge\sigm$; and $a<b$ because $z$ is not
constant on a component of $\Xp$. A continuous function on a circle takes
every value strictly between its minimum and its maximum at least twice,
and $\mathcal A$, which is disjoint from $\mathcal C$, has at least two
points over each $z\in(a,b)\subseteq(\sigm,\Sigm)$.

(c) Theorem~\ref{real:thm:packing}(a) gives $p_U(G)\ge2$. A set of
vertices whose deletion leaves a pseudoforest must meet each of the two
subgraphs, and each has at least four vertices. Let $\taups(G)=2$ and let
$z\in(a,b)$ be neither a critical
value of $z$ on $\Xp$ nor a value excluded by (G). Then all positive
solutions over $z$ are nondegenerate, their Morse indices are at most
$2$ by Proposition~\ref{prop:morse-bound}(i), $c_0=1$ by
Proposition~\ref{prop:morse-bound}(iv), and $c_0-c_1+c_2=0$ by
Corollary~\ref{cor:parity}. With $N(z)\ge4$ this gives $c_2=c_1-1\ge1$.
The solution of index $0$ is the physical one, which is not on
$\mathcal C$.

(d) Put $t_i=D_i-z>0$. The two vertex equations give
$t_0V_0=t_1V_1=:\gamma$, so the positive solutions over $z$ correspond to
the zeros $\gamma>0$ of
$\varphi(\gamma)=1+\sum_ef_+(w_e^2\gamma^2/(t_0t_1))-\gamma$. This function
is strictly convex with $\varphi(0)=1$, so it has at most two zeros. A
circle would give $N(z)\ge4$ by (b).
\end{proof}

Circles do occur.

\begin{proposition}[Existence of circles]\label{real:prop:circles-exist}
On each of the following graphs there is a nonempty open set of real
coefficients for which $\Xp$ has a smooth compact component consisting of
regular points: the four-vertex multigraph formed by joining two $D_3$'s by
a bridge, and the eight-vertex simple graph formed by joining two diamonds
by a bridge between degree-two vertices.
\end{proposition}

\begin{proof}
Take either $D_3$ or the diamond, with zero potential and unit edge
weights, and choose the vertex for the bridge as in the statement.
The resulting rooted graph is a part in the sense of
Definition~\ref{app:bridge:def-part}. Write $\sigma<\Sigma$ for its lower
Perron value and cover spectral edge; the strict inequality is
Corollary~\ref{cor:gap}. In the second copy add a constant $d$ to all
potentials, where $0<d<\Sigma-\sigma$. The thresholds now satisfy
\[
 \sigma_1=\sigma<\sigma_2=\sigma+d<\Sigma_1=\Sigma<\Sigma_2=\Sigma+d.
\]
By Theorem~\ref{app:bridge:thm-gluing}, joining the copies by a sufficiently
weak nonzero bridge gives exactly one circle, all of whose points are
regular. Lemma~\ref{app:bridge:lem-stable} gives the open coefficient
neighborhood.
\end{proof}

For models with specified bridge strengths and energy bounds, see
Proposition~\ref{app:bridge:prop-M}. The weak-bridge theorem also describes
the limiting shape: in the interlaced case $\sigma_2<\Sigma_1<\Sigma_2$ the
circle follows the two solutions of the first part paired with the second
solution of the other part. Its ends are joined at the turning point
$\Sigma_1$ and near the pole $\sigma_2$, where the weak bridge regularizes
the divergence. The full argument is in Appendix~\ref{app:bridge}.

With the same reduction one obtains,
for every $r\ge1$, a connected simple graph with $\taups=2$ (two families
of triangles attached to the ends of an edge) and a nonempty open set of
real coefficients for which $\Xp$ has at least $r$ circles consisting of
regular points (Theorem~\ref{app:bridge:thm-circles}).
So the number of circles is unbounded already for
$\taups=2$, while by Proposition~\ref{prop:circles}(a) it is bounded in
terms of the graph. Whether circles occur on simple graphs with five to
seven vertices is open. For symmetric coefficients $\Xp$ need not be
smooth, and closed curves of positive solutions may meet the arc (see the
end of Section~\ref{per:subsec:csz}).

\subsection{The Bethe free energy and damping}\label{real:subsec:damping}

Cseke and Heskes~\cite{CsekeHeskes} studied the Bethe free energy of
Gaussian models with unit diagonal and its fractional versions. A
computation shows that, for their parameter $\alpha=1$, the part of their
constrained free energy that depends on the variances agrees with $\phi$
up to an additive constant, after the normalization $D_i=1$; a general
positive $D$ is reduced to this case by the substitution $V_i\mapsto D_iV_i$.
For models that are not walk-summable (not pairwise normalizable, in
their terminology) they put forward the hypothesis that either the free
energy has a unique finite local minimum, to which belief propagation
converges for suitable initial values and sufficient damping, or it has
none. They also give an example of a regular graph with a finite local
minimum without walk-summability.

\begin{proposition}\label{real:prop:CH}
Let $G$ be a connected finite simple graph with real coefficients.
\begin{enumerate}
\item[\textup{(a)}] The function $\phi$ has at most one local minimum on
$(0,\infty)^V$. If it has one, it is the maximal fixed point, and it is a
nondegenerate minimum if $\rho(C_{q^{\max}})<1$.
\item[\textup{(b)}] If the model is walk-summable, then $\phi$ tends to
$+\infty$ at the boundary of $(0,\infty)^V$ and at infinity, and its only
critical point is its global minimum.
\end{enumerate}
\end{proposition}

\begin{proof}
(a) At a local minimum $\Hess\phi$ is positive semidefinite; apply
Proposition~\ref{prop:morse-bound}(ii) and (iv). (b) Since
$\sqrt{1+4x}\le1+2\sqrt x$ and $\log(1+\sqrt{1+4x})\ge0$, we have
$g(x)\le1+2\sqrt x$ and, with $v_i=\sqrt{V_i}$ and
$\mu=\lambda_{\min}(D-\abs A)>0$,
\[
\phi(V)\ge\tfrac12v^{\T}(D-\abs A)v-\tfrac12\sum_i\log V_i-\tfrac12\abs E
\ge\sum_i\bigl(\tfrac\mu2V_i-\tfrac12\log V_i\bigr)-\tfrac12\abs E .
\]
So $\phi$ attains its minimum, and by Theorem~\ref{thm:unique-ws} it has
only one critical point.
\end{proof}

Thus the part of the hypothesis of~\cite{CsekeHeskes} that concerns the
uniqueness of the local minimum holds for the variances, for $\alpha=1$
and all coefficients. We have not shown that a degenerate maximal fixed
point is a local minimum, and we make no statement for fractional
$\alpha$. The part that concerns damping fails.

\begin{remark}[Damping]\label{real:rem:damping}
Belief propagation also iterates linear messages, which determine the
means. With the precisions fixed at a fixed point $q$, their update is
affine, with linear part $T_q$ having the entries $-w_{ij}/q_{i\to j}$ in
position $(i\to j,\,k\to i)$, $k\ne j$. The derivative of the joint update
is block triangular with spectrum $\sigma(C_q)\cup\sigma(T_q)$, and the
derivative of the damped update with parameter $\varepsilon\in(0,1]$ has
the eigenvalues $1+\varepsilon(\lambda-1)$. So if $T_q$ has an eigenvalue
with real part greater than $1$, the fixed point is linearly unstable for
every damping; this is the criterion of Su and Wu~\cite{SuWu15}, whose
examples of failure are models without a positive fixed point. It applies
whether the damping acts on all messages, as in~\cite{CsekeHeskes}, or
only on the linear ones, as in~\cite{SuWu15}.

Let $F_{\ell,k}$ be the simple graph formed by $k\ge2$ cycles of odd
length $\ell\ge3$ with one common vertex, let $0<\varrho<1$, and take
$D_i=1+(d_i-1)\varrho^2$ and all weights equal to $\varrho$. Then $q\equiv1$
is a fixed point of \eqref{real:eq:update}; it is a positive solution with
$V_i=1/(1-\varrho^2)$ and $u_e=\varrho^2$, and $C_q=\varrho^2B$, $T_q=-\varrho B$
for the unweighted nonbacktracking matrix $B$ of $F_{\ell,k}$. Following
an eigenvector of $B$ around the cycles shows that the nonzero eigenvalues
$\lambda$ of $B$ satisfy $\lambda^\ell\in\{2k-1,1,-1\}$ and that all
$\ell$-th roots of $2k-1$ occur. Suppose
\[
\frac1{(2k-1)^{1/\ell}\cos(\pi/\ell)}<\varrho<(2k-1)^{-1/(2\ell)},
\]
which is possible if and only if $(2k-1)^{1/(2\ell)}\cos(\pi/\ell)>1$, for
instance for $(\ell,k)=(3,63)$ with $\varrho=21/50$, or for
$(\ell,k)=(7,3)$. Then $\rho(C_q)=\varrho^2(2k-1)^{1/\ell}<1$, so the fixed
point is the maximal one, it is nondegenerate, and it is the unique local
minimum of $\phi$ (Propositions~\ref{prop:max-root},
\ref{prop:morse-bound}, \ref{real:prop:CH}). Since $\ell$ is odd, the real
eigenvalues of $T_q$ lie in $\{0,\varrho,-\varrho,-\varrho(2k-1)^{1/\ell}\}$
and are less than $1$, so $\det(I-sT_q)\ne0$ for $0\le s\le1$. By
Lemma~\ref{real:lem:bass}(a) with the constant weights $-s\varrho$, for
which $B_u=sT_q$, the matrices \eqref{real:eq:Hu} with these weights are
invertible for $0\le s\le1$, hence positive definite; at $s=1$ this matrix
is $\diag(V)^{1/2}(D+A)\diag(V)^{1/2}$, so $D+A$ is positive definite. By Lemma~\ref{real:lem:bass}(b) with the weights
$\varrho$, the model is walk-summable if and only if
$\varrho(2k-1)^{1/\ell}<1$, which fails here. Finally $T_q$ has the
eigenvalues $\varrho(2k-1)^{1/\ell}e^{\pm i\pi/\ell}$, with real part
greater than $1$. So this is a normalizable Gaussian model whose
precision iteration converges to a nondegenerate fixed point, the unique
local minimum of the free energy, and for which belief propagation is
linearly unstable for every damping $\varepsilon\in(0,1]$; with the
precisions at the fixed point, the means diverge for all initial values
outside a proper affine subspace. The statement concerns synchronous
updates with a uniform convex damping; it does not apply to double-loop
or diagonal-loading schemes such as~\cite{JBD}. For discrete variables,
models with a unique fixed point to which belief propagation does not
converge were reported in~\cite{KnollEtAl}.
\end{remark}

%% file: sections/genus.tex
\providecommand{\Xbar}{\overline{X}}
\providecommand{\cD}{\mathcal D}
\providecommand{\tV}{\widetilde V}
\providecommand{\Efr}{E_{\mathrm{fr}}}

\section{The genus of the spectral curve}\label{sec:genus}

In this section $G$ is a connected finite loopless multigraph and, unless
stated otherwise, the coefficients are generic and complex. We write
$\Xbar$ for the normalization of the projective closure of the spectral
curve $\Xc$. By Theorem~\ref{thm:irreducible} it is a connected compact
Riemann surface, and by Theorem~\ref{thm:generic-degree} the function
$z\colon\Xbar\to\PP^1$ has degree $\LG$. Let $g$ be the genus of $\Xbar$,
let $e_P$ be the ramification index of $z$ at $P\in\Xbar$, and let
$R=\sum_{P}(e_P-1)$ be the total ramification, so that
\begin{equation}\label{gen:eq:RH}
R=2g-2+2\LG
\end{equation}
by the Riemann--Hurwitz formula. We compute $2g-2$ as the degree of the
divisor of an explicit meromorphic differential. The differential has
order $2n-2$ at each of the $\LG$ points over $z=\infty$ and order zero at
every other point at which all $V_i$ are finite, including the points
where some $V_i$ vanishes. The remaining points, where some $V_i$ has a
pole, lie over eigenvalues of signed principal submatrices of the Jacobi
matrix; they are classified by combinatorial data which we call
configurations, and the order of the differential at such a point depends
only on the support of the pole. For generic coefficients no such branch
stays in a single fiber of $z$, so summing the local orders gives the exact
genus (Theorem~\ref{thm:genus-formula}). The genus bounds the number of
compact components of $\Xp$ (Proposition~\ref{prop:genus-bound}).

For a single cycle the spectral curve is the classical hyperelliptic
curve of a periodic Jacobi matrix \cite{vMM}; we recover it in
Proposition~\ref{gen:prop:unicyclic}. Gutkin \cite{Gutkin} found that the
Green function of the free product of two operators with two eigenvalues
each lives on an elliptic curve. Apart from these two cases the genus
formula appears to be new.

\subsection{Configurations and the statement}\label{gen:subsec:statement}

We need the covariance system with integer constants on induced subgraphs.
For a loopless multigraph $K$ with vertex set $Y$, coefficients
$(D_y)_{y\in Y}$, $(w_e)_{e\in E(K)}$ and $c\in\Z^Y$, consider
\begin{equation}\label{gen:eq:system-c}
(D_y-z)V_y=c_y+\sum_{e\ni y}f_e\quad(y\in Y),\qquad
f_e(1+f_e)=w_e^2V_yV_{y'}\quad(e=yy'\in E(K)),
\end{equation}
and put $L_c(K)=\#\{S\subseteq E(K): d_S(y)\ne c_y\text{ for all }y\in Y\}$.
For $c=\ones$ this is \eqref{eq:system} and $L_{\ones}(K)=L(K)$; for
$Y=\emptyset$ the system has one (empty) solution and $L_c(K)=1$. A
solution of a system of equations with as many equations as unknowns is
\emph{nondegenerate} if the Jacobian matrix with respect to the unknowns
is invertible at it; for \eqref{gen:eq:system-c} the unknowns are $(V,f)$
and $z$ is fixed.

\begin{lemma}\label{gen:lem:constants}
Let $K$ and $c$ be as above. There is a nonempty Zariski open set
$\mathcal G_c(K)\subseteq\C^{E(K)}$, containing a neighborhood of the
origin, with the following property. If $z\notin\{D_y\}$ and
$\bigl(w_e^2/((D_y-z)(D_{y'}-z))\bigr)_{e=yy'}\in\mathcal G_c(K)$, then
\eqref{gen:eq:system-c} has exactly $2^{\abs{E(K)}}$ solutions for this
$z$, all of them nondegenerate, and exactly $L_c(K)$ of them have
$V_y\ne0$ for all $y$.
\end{lemma}

\begin{proof}
For $c=\ones$ this is the content of the proof of
Theorem~\ref{thm:generic-degree}: after the elimination
$V_y=(c_y+\sum_{e\ni y}f_e)/(D_y-z)$ the system becomes
$f_e(1+f_e)=\varepsilon_e\,(c_y+\sum f)(c_{y'}+\sum f)$ with couplings
$\varepsilon_e$ as in the statement; at $\varepsilon=0$ its solutions are
the $2^{\abs{E(K)}}$ points $f=-\ones_S$, at which the Jacobian is
$\diag(\pm1)$ and $c_y+\sum_{e\ni y}f_e=c_y-d_S(y)$. The constant $1$
enters that proof only through this last value and through lower order
terms of the equations, so the proof applies verbatim to integer
constants; see also Section~\ref{sec:irreducibility}, where the systems
with $c\le\ones$ are used.
\end{proof}

For $U\subseteq V$ let $E(U)$ be the set of edges with both ends in $U$,
$\beta(U)=\abs{E(U)}-\abs U+1$ the cycle rank of $G[U]$ when $G[U]$ is
connected, $N(U)$ the set of vertices outside $U$ adjacent to $U$, and
$Y_U=V\setminus(U\cup N(U))$. For $x\in N(U)$ let $\partial_x$ be the set
of edges joining $x$ to $U$, $k_x=\abs{\partial_x}$, and
$\partial U=\bigcup_x\partial_x$; for $e\in\partial U$ we write $i(e)$ for
its end in $U$. The \emph{frozen edges} are the edges with no end in $U$
and at least one end in $N(U)$; their set is $\Efr(U)$. Thus
$E=E(U)\sqcup\partial U\sqcup\Efr(U)\sqcup E(Y_U)$.

\begin{definition}[Configurations]\label{def:configuration}
A \emph{configuration} of $(G,D,w)$ is a tuple
$\gamma=(U,s,\lambda,F,\psi)$ of the following data.
\begin{enumerate}
\item $U\subseteq V$ is nonempty and $G[U]$ is connected.
\item $s$ is a switching class of signings of $E(U)$, and $\lambda$ is a
simple eigenvalue of $M_s=(D+A^s)[U]$ whose eigenvector $u$ has no zero
entry and which satisfies $\lambda\ne D_x$ for $x\notin U$. We normalize
$u^{\T}u=1$; the numbers $u_i^2$ and $s_eu_iu_j$ ($e=ij\in E(U)$) depend
only on the class of $s$. For $U=\{v\}$ we have $\lambda=D_v$ and $u_v^2=1$.
\item $F=(f^0,f^*,\nu)$, where $f^0\in\{0,-1\}^{\Efr(U)}$ is a choice of
\emph{frozen values} and, for each $x\in N(U)$, the pair
$((f^*_e)_{e\in\partial_x},\nu_x)$ is a nondegenerate solution with
$\nu_x\ne0$ of the \emph{boundary system}
\begin{equation}\label{gen:eq:boundary}
f_e(1+f_e)=w_e^2u_{i(e)}^2\,\nu_x\quad(e\in\partial_x),\qquad
\sum_{e\in\partial_x}f_e=b_x,\qquad
b_x:=-1-\sum_{e\in\Efr(U),\,e\ni x}f^0_e .
\end{equation}
\item $\psi$ is a nondegenerate solution with all $V_y\ne0$ of the system
\eqref{gen:eq:system-c} on $K=G[Y_U]$ at $z=\lambda$ with the constants
$c_y=1+\sum_{e\in\Efr(U),\,e\ni y}f^0_e$.
\end{enumerate}
We set
\[
\kappa(\gamma)=\abs{E(U)}-\abs U-\sum_{x\in N(U)}b_x,\qquad
\delta(\gamma)=\beta(U)+2\abs{N(U)}-2 .
\]
Lemma~\ref{gen:lem:eta} attaches to $\gamma$ an analytic family of
solutions of \eqref{eq:system}, parametrized by $\eta=1/\sum_{i\in U}V_i$
near $\eta=0$, along which $z=\lambda+\zeta_\gamma(\eta)$ with
$\zeta_\gamma(\eta)=\kappa(\gamma)\eta+O(\eta^2)$. The configuration is
\emph{vertical} if $\zeta_\gamma\equiv0$ and \emph{non-vertical}
otherwise. In particular every configuration with $\kappa(\gamma)\ne0$ is
non-vertical.
\end{definition}

The number $\delta(\gamma)$ depends only on $U$, and we also write
$\delta(U)$. The integers $b_x$ and $c_y$ are determined by $f^0$:
$b_x+1$ and $1-c_y$ are the numbers of frozen edges with value $-1$ at
$x$ and at $y$. For generic coefficients the configurations with a given
$(U,f^0)$ are counted as follows (Lemma~\ref{gen:lem:count}): there are
$2^{\beta(U)}$ switching classes, $\abs U$ eigenvalues for each, and for
each of these
\begin{equation}\label{gen:eq:pi}
\prod_{x\in N(U)}a(k_x,b_x)\cdot L_{c}(G[Y_U]),\qquad
a(k,b)=2^{k-1}-\binom{k}{-b},
\end{equation}
choices of $(f^*,\nu,\psi)$, where $\binom k{j}=0$ for $j\notin[0,k]$. We
write $\pi(U)$ for the sum of \eqref{gen:eq:pi} over
$f^0\in\{0,-1\}^{\Efr(U)}$, and
\begin{equation}\label{gen:eq:Sigma}
\Sigma(G)=\sum_{U}\abs U\,2^{\beta(U)}\,\pi(U)\,\bigl(\beta(U)+2\abs{N(U)}-2\bigr),
\end{equation}
the sum over all nonempty $U\subseteq V$ with $G[U]$ connected. Thus
$\Sigma(G)=\sum_\gamma\delta(\gamma)$ for generic coefficients, and
$\Sigma(G)$ is an integer that depends only on the graph.

\begin{theorem}[Genus formula]\label{thm:genus-formula}
Let $G$ be a connected finite loopless multigraph with $n$ vertices. For
generic coefficients the genus $g$ of the spectral curve satisfies
\begin{equation}\label{gen:eq:genus}
2g-2=2(n-1)\LG+\Sigma(G),
\end{equation}
equivalently $R=2n\LG+\Sigma(G)$. All configurations are non-vertical.
\end{theorem}

\begin{proposition}\label{prop:genus-bound}
Let $G$ be as in Theorem~\ref{thm:genus-formula}, with $\beta\ge2$.
For generic real coefficients the number of compact connected
components of $\Xp$ is at most $g=1+(n-1)\LG+\tfrac12\Sigma(G)$.
\end{proposition}

The proofs occupy Sections~\ref{gen:subsec:divisor}
and~\ref{gen:subsec:poles}. Examples, the table of values and the
behavior under subdivision of edges are in
Sections~\ref{gen:subsec:examples} and~\ref{gen:subsec:subdivision}.

\begin{remark}\label{gen:rem:kappa}
For $\kappa(\gamma)\ne0$ the branch attached to $\gamma$ is the pole branch
of Theorems~\ref{thm:local-pole} and~\ref{thm:partial-pole}: since
$\sum_{i\in U}V_i=1/\eta$ and $z-\lambda=\kappa(\gamma)\eta+O(\eta^2)$, the
differential $(\sum_iV_i)\,dz$ has residue $\kappa(\gamma)$ there, and $z$
is unramified. By Theorem~\ref{thm:classification} these are all points
with nonzero residue. The configurations with $\kappa(\gamma)=0$ carry no
residue; the non-vertical ones are ramification points of $z$. For a cycle
they are the band edges, where $V_i$ behaves like $(z-\lambda)^{-1/2}$.
\end{remark}

\subsection{The canonical divisor away from the poles}\label{gen:subsec:divisor}

Write the system \eqref{eq:system} as $\Phi_i=0$ ($i\in V$), $\Psi_e=0$
($e\in E$) with
\[
\Phi_i=(D_i-z)V_i-1-\sum_{e\ni i}f_e,\qquad
\Psi_e=f_e(1+f_e)-w_e^2V_iV_j\quad(e=ij),
\]
and let $J=\det\partial(\Phi,\Psi)/\partial(V,f)$. We use the following
standard facts about Poincar\'e residues. Let $Z\subset\C^{1+m}$ be the
zero set of $m$ functions $H_1,\dots,H_m$ of $(z,x)$, and suppose that the
Jacobian matrix $\partial H/\partial(z,x)$ has rank $m$ at a point of $Z$.
Then $Z$ is a smooth curve near this point and
\[
\omega_Z=\frac{dz}{\det(\partial H/\partial x)}
=\pm\frac{dx_k}{\det(\partial H/\partial(z,x_{\hat k}))}
\]
(any expression with nonvanishing denominator) is a holomorphic
differential without zeros on $Z$ near this point. If
$\tilde x=\tilde x(z,x)$ is a change of the variables $x$ depending
holomorphically on $z$ and $\tilde H_k=a_kH_k$ with holomorphic functions
$a_k$, then on $Z$
\begin{equation}\label{gen:eq:residue-change}
\det\frac{\partial\tilde H}{\partial\tilde x}
=\Bigl(\prod_ka_k\Bigr)\,\det\frac{\partial H}{\partial x}\,
\Bigl(\det\frac{\partial\tilde x}{\partial x}\Bigr)^{-1},
\end{equation}
because $d(a_kH_k)=a_k\,dH_k$ on $Z$.

On $\Xbar$ we consider the meromorphic differentials
\[
\omega_0=\frac{dz}{J},\qquad\theta=\omega_0\cdot\prod_{i\in V}V_i .
\]
They are not identically zero: at a general point of $\Xc$ the solution is
nondegenerate (Theorem~\ref{thm:generic-degree}), so $J\not\equiv0$. Let
$\cD\subset\Xbar$ be the set of points over finite values of $z$ at which
some $V_i$ has a pole. At the other points over finite $z$ all $V_i$, and
therefore all $f_e$, are finite.

\begin{lemma}\label{gen:lem:no-vertical}
For generic coefficients there is no $z_0\in\C$ for which the solutions of
\eqref{eq:system} with $z=z_0$ and $V_i\ne0$ for all $i$ form a set of
positive dimension. Consequently every solution $(z,V,f)$ of
\eqref{eq:system} with all $V_i\ne0$ is a point of $\Xc$.
\end{lemma}

\begin{proof}
Let $\mathcal Y$ be the set of $(D',w,V,f)$ with $w_e\ne0$, $V_i\ne0$,
$D'_iV_i=1+\sum_{e\ni i}f_e$ and $f_e(1+f_e)=w_e^2V_iV_j$. The projection
to $(V,w)$ exhibits $\mathcal Y$ as the fiber product over
$(\C^*)^V\times(\C^*)^E$ of the double covers $f_e^2+f_e=w_e^2V_iV_j$,
$e\in E$. The branch divisors $1+4w_e^2V_iV_j=0$ in the torus are
irreducible (solve for $V_i$) and distinct (each uses its own $w_e$).
The valuation along the divisor for $e$ is one for its discriminant and
zero for all the others. Thus the discriminants have independent square
classes in the rational function field of the torus. The fiber product
is finite free over the torus and its generic algebra is the resulting
multiquadratic field. Its coordinate ring injects into this field, so
$\mathcal Y$ is irreducible of dimension $n+\abs E$.
The projection $p\colon\mathcal Y\to\C^V\times(\C^*)^E$, $(D',w,V,f)\mapsto(D',w)$,
is dominant with finite general fiber by Theorem~\ref{thm:generic-degree}.
Hence the set $\mathcal Y_1$ of points at which the fiber of $p$ has
positive dimension is a proper closed subset, $\dim\mathcal Y_1\le n+\abs E-1$,
and $B=p(\mathcal Y_1)$ has dimension at most $n+\abs E-2$. The solutions
of \eqref{eq:system} for the coefficients $(D,w)$ and $z=z_0$ with
$V\in(\C^*)^V$ are the fiber of $p$ over $(D-z_0\ones,w)$. If $(D,w)$ lies
outside the closure of $B+\C\cdot(\ones,0)$, which has dimension at most
$n+\abs E-1$, no such fiber has positive dimension. For the last
statement, the solutions with all $V_i\ne0$ form an algebraic set all of
whose components have dimension at least one; on each component $z$ is
not constant, so a general point of it is admissible and the component is
contained in $\Xc$.
\end{proof}

\begin{lemma}[Points with finite nonzero $V$]\label{gen:lem:regular}
For generic coefficients, $\ord_P\theta=0$ at every $P\in\Xbar$ over a
finite value of $z$ at which all $V_i$ are finite and nonzero.
\end{lemma}

\begin{proof}
Let $\mathcal T$ be the set of $(D,w,z,V,f)$ satisfying \eqref{eq:system}
with $w_e\ne0$ and $V_i\ne0$. Since $\partial\Phi_i/\partial D_i=V_i\ne0$
and $\partial\Psi_e/\partial w_e=-2w_eV_iV_j\ne0$, and each equation has
its own coefficient, the Jacobian matrix with respect to $(D,w)$ has full
rank, so $\mathcal T$ is smooth of dimension $n+\abs E+1$. By generic
smoothness in characteristic zero, for generic $(D,w)$ the fiber
$\mathcal T_{(D,w)}$ is a smooth curve, that is, $\partial(\Phi,\Psi)/\partial(z,V,f)$
has full rank at each of its points. By Lemma~\ref{gen:lem:no-vertical}
this curve is an open subset of $\Xc$. Hence $P$ is a smooth point of
$\Xc$, $\omega_0$ is the Poincar\'e residue of the system there, which is
holomorphic without zeros, and the $V_i$ are units at $P$.
\end{proof}

Note that this lemma includes the points over $z=D_i$ with $V_i\ne0$.

\begin{lemma}[Points with a vanishing $V_i$]\label{gen:lem:zero}
For generic coefficients, $\ord_P\theta=0$ at every $P\in\Xbar$ over a
finite value of $z$ at which all $V_i$ are finite and
$I=\{i:V_i(P)=0\}$ is nonempty.
\end{lemma}

\begin{proof}
Let $E_I$ be the set of edges with at least one end in $I$. For
$e=ij\in E_I$ we have $f_e(1+f_e)=w_e^2V_iV_j\to0$ at $P$, so
$f_e(P)=f^0_e\in\{0,-1\}$. Put $\mu_e=V_i$ if $i\in I$, $j\notin I$, and
$\mu_e=V_iV_j$ if $i,j\in I$. Since
$(f_e-f^0_e)(f_e+f^0_e+1)=f_e(1+f_e)$ and $f_e+f^0_e+1\to2f^0_e+1=\pm1$,
the function
\[
h_e=\frac{f_e-f^0_e}{\mu_e}=\frac{w_e^2V_iV_j}{\mu_e\,(f_e+f^0_e+1)}
\]
is holomorphic at $P$. For $i\in I$ the vertex equation at $P$ gives
$1+\sum_{e\ni i}f^0_e=0$. Substituting $f_e=f^0_e+\mu_eh_e$ for
$e\in E_I$ we therefore have identically
\begin{align*}
\Phi_i&=V_i\,\hat\Phi_i,&
\hat\Phi_i&=(D_i-z)-\sum_{e\ni i}\frac{\mu_e}{V_i}h_e&&(i\in I),\\
\Psi_e&=\mu_e\hat\Psi_e,&
\hat\Psi_e&=h_e\bigl(2f^0_e+1+\mu_eh_e\bigr)-\frac{w_e^2V_iV_j}{\mu_e}&&(e\in E_I),
\end{align*}
where $\mu_e/V_i$ is $1$ or $V_j$, and $w_e^2V_iV_j/\mu_e$ is $w_e^2V_j$
or $w_e^2$; these are polynomials. Consider the system $\hat H$ in the
unknowns $(z,V,(h_e)_{e\in E_I},(f_e)_{e\notin E_I})$ consisting of
$\hat\Phi_i$ ($i\in I$), $\hat\Psi_e$ ($e\in E_I$) and of the remaining
$\Phi_j$, $\Psi_e$ after the substitution, and let $\hat J$ be its
Jacobian determinant with respect to all unknowns except $z$. On the open
set where $V_i\ne0$ for $i\in I$, the change of variables
$(V,f)\mapsto(V,h,f)$ has Jacobian determinant $\prod_{e\in E_I}\mu_e^{-1}$,
and \eqref{gen:eq:residue-change} gives
\[
J=\Bigl(\prod_{i\in I}V_i\Bigr)\hat J,\qquad
\theta=\frac{dz}{\hat J}\cdot\prod_{j\notin I}V_j
\]
on a punctured neighborhood of $P$ in $\Xbar$, hence as meromorphic
differentials. The total space of the system $\hat H$ over the
coefficient space, on the open set $V_j\ne0$ ($j\notin I$), is smooth:
$\partial\hat\Phi_i/\partial D_i=1$, $\partial\Phi_j/\partial D_j=V_j\ne0$,
$\partial\hat\Psi_e/\partial w_e$ is $-2w_eV_j$ or $-2w_e$ for
$e\in E_I$, and $\partial\Psi_e/\partial w_e=-2w_eV_iV_j\ne0$ for the other
edges. By generic smoothness, for generic coefficients the zero set
$\hat Z$ of $\hat H$ is a smooth curve; there are finitely many pairs
$(I,f^0)$. The germ of $\Xbar$ at $P$ maps to $\hat Z$ by the holomorphic
functions $(z,V,h,f)$, and this map is injective away from $P$, so it is
an isomorphism onto the germ of $\hat Z$ at the image of $P$. Therefore
$dz/\hat J$, the Poincar\'e residue of $\hat H$, has order zero at $P$,
and the $V_j$ with $j\notin I$ are units.
\end{proof}

\begin{lemma}[Infinity]\label{gen:lem:theta-infinity}
For all coefficients, $\ord_P\theta=2n-2$ at each of the $\LG$ points $P$
of $\Xbar$ over $z=\infty$.
\end{lemma}

\begin{proof}
By Lemma~\ref{lem:infinity} there are exactly $\LG$ such points, they are
unramified, and they are the branches $f=f^S(z)$, $S\in\cL(G)$, analytic
in $\tau=1/z$ with $f^S(\infty)=-\ones_S$. On such a branch
$V_i=(1+\sum_{e\ni i}f_e)/(D_i-z)$ with
$1+\sum_{e\ni i}f_e\to1-d_S(i)\ne0$, so $\ord_\tau V_i=1$. Eliminating the
$V$-columns,
\[
J=\prod_i(D_i-z)\cdot\det\Bigl(\diag(1+2f_e)-C\Bigr),
\]
where the entries of $C$ are $w_e^2V_j/(D_i-z)$ and $w_e^2V_i/(D_j-z)$,
of order $O(\tau^2)$. Since $1+2f_e\to\pm1$, $\ord_\tau J=-n$. With
$dz=-\tau^{-2}d\tau$ this gives $\ord\omega_0=n-2$ and
$\ord\theta=2n-2$.
\end{proof}

\begin{proposition}[Localization]\label{gen:prop:localization}
For generic coefficients
\begin{equation}\label{gen:eq:localization}
2g-2=2(n-1)\LG+\sum_{P\in\cD}\ord_P\theta .
\end{equation}
\end{proposition}

\begin{proof}
The degree of the divisor of the nonzero meromorphic differential
$\theta$ on $\Xbar$ is $2g-2$. Its order is $2n-2$ at the $\LG$ points
over $\infty$ (Lemma~\ref{gen:lem:theta-infinity}) and zero at the points
over finite $z$ outside $\cD$ (Lemmas~\ref{gen:lem:regular}
and~\ref{gen:lem:zero}).
\end{proof}

\subsection{The points where \texorpdfstring{$V$}{V} has a pole}\label{gen:subsec:poles}

Let $\gamma=(U,s,\lambda,F,\psi)$ be a configuration, with eigenvector $u$,
and put $y^0_e=-2s_ew_eu_iu_j$ for $e=ij\in E(U)$. Consider the following
system with parameter $\eta$ in the unknowns $\tV\in\C^U$, $y\in\C^{E(U)}$,
$\zeta\in\C$, $(f_e)_{e\in\partial U\cup\Efr(U)\cup E(Y_U)}$,
$\nu\in\C^{N(U)}$ and $(V_y)_{y\in Y_U}$:
\begin{equation}\label{gen:eq:eta-system}
\begin{aligned}
&(D_i-\lambda-\zeta)\tV_i-\eta\Bigl(1-\tfrac12d_{E(U)}(i)+\!\!\sum_{e\in\partial U,\,e\ni i}\!\!f_e\Bigr)
-\tfrac12\!\!\sum_{e\in E(U),\,e\ni i}\!\!y_e=0&&(i\in U),\\
&\sum_{i\in U}\tV_i=1,\qquad y_e^2-\eta^2-4w_e^2\tV_i\tV_j=0&&(e=ij\in E(U)),\\
&f_e(1+f_e)-w_e^2\,\tV_{i(e)}\,\nu_x=0&&(e\in\partial_x,\ x\in N(U)),\\
&(D_x-\lambda-\zeta)\,\eta\nu_x-1-\sum_{e\ni x}f_e=0&&(x\in N(U)),\\
&f_e(1+f_e)-w_e^2\,\widehat V_x\widehat V_{x'}=0&&(e=xx'\in\Efr(U)),\\
&(D_y-\lambda-\zeta)V_y-1-\sum_{e\ni y}f_e=0,\qquad
f_e(1+f_e)-w_e^2V_yV_{y'}=0&&(y\in Y_U,\ e=yy'\in E(Y_U)),
\end{aligned}
\end{equation}
where $\widehat V_x=\eta\nu_x$ for $x\in N(U)$ and $\widehat V_x=V_x$ for
$x\in Y_U$. The \emph{base point} of $\gamma$ is
$\tV_i=u_i^2$, $y=y^0$, $\zeta=0$, $f_e=f^*_e$ on $\partial U$, $\nu$ as in
$F$, $f_e=f^0_e$ on $\Efr(U)$, and $\psi$ on $Y_U$.

\begin{lemma}[The $\eta$-coordinate]\label{gen:lem:eta}
Let $\gamma$ be a configuration.
\begin{enumerate}
\item The base point solves \eqref{gen:eq:eta-system} at $\eta=0$, and the
Jacobian matrix of \eqref{gen:eq:eta-system} with respect to the unknowns
is invertible there. Hence there is a unique analytic solution
$\gamma(\eta)$, $\abs\eta<\eta_0$, through the base point.
\item For $\eta\ne0$ the substitution
\begin{equation}\label{gen:eq:eta-subst}
z=\lambda+\zeta,\qquad V_i=\tV_i/\eta\ (i\in U),\qquad
f_e=\tfrac12(y_e/\eta-1)\ (e\in E(U)),\qquad V_x=\eta\nu_x\ (x\in N(U))
\end{equation}
identifies the solutions of \eqref{gen:eq:eta-system} with the solutions
of \eqref{eq:system} for which $\sum_{i\in U}V_i=1/\eta$. In particular
every solution of \eqref{eq:system} for which $\eta=1/\sum_{i\in U}V_i$ is
small and the point defined by \eqref{gen:eq:eta-subst} is close to the
base point lies on $\gamma(\eta)$.
\item $\zeta(\eta)=\kappa(\gamma)\,\eta+O(\eta^2)$.
\item If $\gamma$ is non-vertical, then $\gamma(\eta)$ is admissible for
small $\eta\ne0$, and $\eta\mapsto\gamma(\eta)$ extends to a
parametrization of a neighborhood of a point $P_\gamma\in\cD$ of $\Xbar$
over $z=\lambda$ with local coordinate $\eta$. At $P_\gamma$ the function
$V_i$ has a simple pole for $i\in U$, a simple zero for $i\in N(U)$, and
is finite and nonzero for $i\in Y_U$; the ramification index of $z$ is
$\ord_{\eta=0}\zeta_\gamma$. Distinct non-vertical configurations give
distinct points.
\end{enumerate}
Parts (1)--(3) do not use that $V_y\ne0$ in $\psi$.
\end{lemma}

\begin{proof}
(2) Multiply the vertex equation of $i\in U$ by $\eta$ and use
$\eta f_e=\tfrac12(y_e-\eta)$ for $e\in E(U)$; multiply the identity
$(1+2f_e)^2=1+4w_e^2V_iV_j$, equivalent to the edge equation, by $\eta^2$;
and note $V_{i(e)}V_x=\tV_{i(e)}\nu_x$ for $e\in\partial_x$. The other
equations are those of \eqref{eq:system} after the substitution, and
$\sum_U\tV_i=1$ is the definition of $\eta$. The last assertion follows
from the uniqueness in (1).

(1) At $\eta=0$ the first line of \eqref{gen:eq:eta-system} reads
$u_i\bigl((D_i-\lambda)u_i+\sum_{e=ij}s_ew_eu_j\bigr)=0$, the second line
holds because $u^{\T}u=1$ and $(y^0_e)^2=4w_e^2u_i^2u_j^2$, the third and
fourth are the boundary system \eqref{gen:eq:boundary}, the fifth is
$f^0_e(1+f^0_e)=0$, and the last is the system of $\psi$, since
$1+\sum_{e\in\Efr(U),e\ni y}f^0_e=c_y$. Order the unknowns as
$(\tV,y,\zeta)$, the frozen $f_e$, the unknowns on $Y_U$, and
$(f_{\partial U},\nu)$, and the equations correspondingly (lines 1--2,
line 5, line 6, lines 3--4). At $\eta=0$ the equations of lines 1--2
contain only $(\tV,y,\zeta)$; the equation of a frozen edge has
differential $(1+2f^0_e)\,df_e=\pm df_e$; the equations on $Y_U$ contain
$\zeta$, the frozen $f_e$ and the unknowns on $Y_U$; and lines 3--4 do not
contain the unknowns on $Y_U$. So the Jacobian matrix is block lower
triangular. Its diagonal blocks are: the block of $(\tV,y,\zeta)$; a
diagonal matrix with entries $\pm1$; the Jacobian of the system of $\psi$,
invertible by assumption; and the product over $x\in N(U)$ of the
Jacobians of the boundary systems, invertible by assumption. For the
first block, $\partial/\partial y_e$ of the edge equation is
$2y^0_e\ne0$, and eliminating $dy_e=-s_ew_e\bigl(\tfrac{u_j}{u_i}d\tV_i+\tfrac{u_i}{u_j}d\tV_j\bigr)$
and using $\sum_{e=ij}s_ew_eu_j=-(D_i-\lambda)u_i$ turns the differentials
of lines 1--2 into
\[
\mathcal J\,d\tV-u^2\,d\zeta,\quad\ones^{\T}d\tV,\qquad
\mathcal J=\tfrac12\diag(u)\,(M_s-\lambda)\,\diag(u)^{-1},
\]
where $u^2=(u_i^2)_{i\in U}$. Since $\lambda$ is simple,
$\ker\mathcal J=\C u^2$, and $\ones^{\T}\mathcal J=0$. If
$\mathcal Jx=tu^2$ and $\ones^{\T}x=0$, then applying $\ones^{\T}$ gives
$t=0$, so $x\in\C u^2$ and $x=0$. The block is invertible.

(3) Differentiating line 1 in $\eta$ at $\eta=0$ and applying $\ones^{\T}$
gives
$-\zeta'(0)=\sum_{i\in U}\bigl(1-\tfrac12d_{E(U)}(i)\bigr)+\sum_{e\in\partial U}f^*_e
=\abs U-\abs{E(U)}+\sum_xb_x$.

(4) For small $\eta\ne0$ we have $\tV_i\approx u_i^2\ne0$,
$\nu_x\ne0$ and $V_y\ne0$ on $Y_U$, and $z=\lambda+\zeta(\eta)$ avoids
$\{D_i\}$ because $\zeta$ is analytic and not constant. So
$\eta\mapsto\gamma(\eta)$ is an analytic arc of admissible solutions on
which $z$ is not constant; it lies in $\Xc$, and it is injective because
$\eta=1/\sum_UV_i$ is a function on $\Xc$. It extends over $\eta=0$ to a
point $P_\gamma$ of $\Xbar$, with $\eta$ as local coordinate. The orders of
the $V_i$ are read off from \eqref{gen:eq:eta-subst}, and
$z-\lambda=\zeta_\gamma(\eta)$. The data of $\gamma$ are recovered from the
germ at $P_\gamma$: $U$ is the set of poles of $V$, $\lambda=z(P_\gamma)$,
the limits of $\eta(1+2f_e)$ are the $y^0_e$, which determine the
switching class, and $F$ and $\psi$ are limits of coordinates.
\end{proof}

\begin{lemma}[Local term]\label{gen:lem:local-term}
For a non-vertical configuration $\gamma$,
$\ord_{P_\gamma}\theta=\delta(\gamma)=\beta(U)+2\abs{N(U)}-2$.
\end{lemma}

\begin{proof}
Add to \eqref{eq:system} the unknown $\eta$ and the equation
$\eta\sum_{i\in U}V_i-1=0$. The Jacobian determinant of the enlarged
system with respect to $(V,f,\eta)$ is $J\cdot\sum_UV_i=J/\eta$, so its
Poincar\'e residue is $\eta\,\omega_0$. The system
\eqref{gen:eq:eta-system}, in which we now regard $\eta$ as an unknown and
$z=\lambda+\zeta$ as the distinguished variable, is obtained from the
enlarged system by multiplying $\abs U$ equations by $\eta$ and
$\abs{E(U)}$ equations by $4\eta^2$, and by the change of variables
$\tV=\eta V_U$, $y=\eta(1+2f_{E(U)})$, $\nu=V_{N(U)}/\eta$, whose Jacobian
determinant is $\eta^{\abs U}(2\eta)^{\abs{E(U)}}\eta^{-\abs{N(U)}}$. By
\eqref{gen:eq:residue-change} the Poincar\'e residue $\omega'$ of
\eqref{gen:eq:eta-system} satisfies
\[
\omega'=\frac{\eta\,\omega_0}{2^{\abs{E(U)}}\,\eta^{\abs{E(U)}+\abs{N(U)}}}
\]
for $\eta\ne0$. On the other hand
$\omega'=\pm d\eta/\det\bigl(\partial/\partial(\tV,y,\zeta,\dots)\bigr)$,
and this determinant is nonzero at $\eta=0$ by
Lemma~\ref{gen:lem:eta}(1), so $\ord_{P_\gamma}\omega'=0$. Hence
$\ord_{P_\gamma}\omega_0=\abs{E(U)}+\abs{N(U)}-1$, and by
Lemma~\ref{gen:lem:eta}(4)
$\ord_{P_\gamma}\theta=\abs{E(U)}+\abs{N(U)}-1-\abs U+\abs{N(U)}$.
\end{proof}

We now fix the genericity conditions. Besides those of
Theorems~\ref{thm:generic-degree} and~\ref{thm:irreducible},
Lemmas~\ref{gen:lem:no-vertical}--\ref{gen:lem:zero} and
Lemma~\ref{irr:lem:band-edge}, we require the following, for every
$U\subseteq V$ with $G[U]$ connected and every signing $s$ of $E(U)$.
\begin{itemize}
\item[(C1)] The eigenvalues of $M_s=(D+A^s)[U]$ are simple and their
eigenvectors have no zero entry.
\item[(C2)] No eigenvalue of $M_s$ equals $D_x$ with $x\notin U$, and
$M_s$ has no common eigenvalue with $(D+A^{s'})[U']$ if $U'\cap U=\emptyset$.
\item[(C3)] For every eigenvector $u$ of $M_s$, every $x\in N(U)$ and
every $\varepsilon\in\{\pm1\}^{\partial_x}$,
$\sum_{e\in\partial_x}\varepsilon_ew_eu_{i(e)}\ne0$ and
$\sum_{e\in\partial_x}\varepsilon_ew_e^2u_{i(e)}^2\ne0$.
\item[(C4)] For every eigenvector $u$ of $M_s$, every $x\in N(U)$ and
every integer $b$, all solutions with $\nu_x\ne0$ of the boundary system
\eqref{gen:eq:boundary} are nondegenerate.
\item[(C5)] For every eigenvalue $\lambda$ of $M_s$ and every
$f^0\in\{0,-1\}^{\Efr(U)}$, the couplings
$w_e^2/((D_y-\lambda)(D_{y'}-\lambda))$, $e=yy'\in E(Y_U)$, lie in the set
$\mathcal G_c(G[Y_U])$ of Lemma~\ref{gen:lem:constants} for the constants
$c=c(f^0)$.
\end{itemize}
Each condition holds on a nonempty Zariski open set of coefficients. For
(C1) and (C2) the relevant discriminants, resultants and minors are
polynomials in $(D,w)$ which are nonzero near $w=0$ when the $D_i$ are
distinct: the eigenvalues are then close to the distinct numbers $D_i$,
$i\in U$, the leading term of $u_j/u_i$ is a sum of distinct monomials
over the shortest paths from $i$ to $j$ in $G[U]$, the numbers $D_x$,
$x\notin U$, do not enter $M_s$, and matrices on disjoint vertex sets have
no coefficient in common. For (C3)--(C5), the triples
$(\text{coefficients},\lambda,u)$ form a variety that is finite over the
coefficient space and irreducible up to $u\mapsto-u$ (it is parametrized
by $(w,\lambda,u)$ and the $D_x$ with $x\notin U$, through
$D_i=\lambda-(A^su)_i/u_i$ for $i\in U$), the conditions are open on it
and invariant under $u\mapsto-u$, and they hold somewhere because the
weights on $\partial_x$ and on $E(Y_U)$ do not enter $M_s$: (C3) is a
nontrivial polynomial condition on $(w_e)_{e\in\partial_x}$; in (C4) the
map $(w_e)_{e\in\partial_x}\mapsto(w_e^2u_{i(e)}^2)_e$ is dominant, and
for general $a\in(\C^*)^{\partial_x}$ all solutions with $\nu\ne0$ of
$f_e(1+f_e)=a_e\nu$, $\sum f_e=b$ are nondegenerate, by generic smoothness
of the map $(f,\nu)\mapsto(f_e(1+f_e)/\nu)_e$ on the smooth variety
$\{\sum_ef_e=b,\ \nu\ne0\}$; and in (C5) the couplings run through a dense
subset of $\C^{E(Y_U)}$ when the weights on $E(Y_U)$ vary. The image in the
coefficient space of the proper closed subset where a condition fails is
contained in a proper closed subset.

\begin{lemma}\label{gen:lem:pole-points}
For generic coefficients, $\gamma\mapsto P_\gamma$ is a bijection from the
set of non-vertical configurations onto the set $\cD$ of points of $\Xbar$
over finite $z$ at which some $V_i$ has a pole. At $P_\gamma$ the function
$\eta=1/\sum_{i\in U}V_i$ is a local coordinate, $z(P_\gamma)=\lambda$, the
orders of the $V_i$ are $-1$ on $U$, $+1$ on $N(U)$ and $0$ on $Y_U$, and
the ramification index of $z$ is $1$ if $\kappa(\gamma)\ne0$ and at least
$2$ if $\kappa(\gamma)=0$.
\end{lemma}

\begin{proof}
By Lemma~\ref{gen:lem:eta} it remains to show that every $P\in\cD$ is some
$P_\gamma$. Let $t$ be a local coordinate at $P$, $z^*=z(P)$,
$o_i=\ord_tV_i\in\Z$, $o=\min_io_i<0$, $U=\{i:o_i=o\}$, and
$V_i=v_it^{o}+\cdots$ with $v_i\ne0$ for $i\in U$.

\emph{Step 1: $U$ supports a signed eigenvector.} For $e=ij\in E(U)$ the
edge equation gives $f_e=b_et^{o}+\cdots$ with $b_e^2=w_e^2v_iv_j$. For an
edge from $i\in U$ to $j\notin U$ we have $\ord(V_iV_j)>2o$, hence
$\ord f_e>o$. The coefficient of $t^{o}$ in the vertex equation of
$i\in U$ is $(D_i-z^*)v_i=\sum_{e\in E(U),\,e\ni i}b_e$. Choose square
roots $u_i$ of $v_i$ and put $s_e=-b_e/(w_eu_iu_j)\in\{\pm1\}$; then
$(M_s-z^*)u=0$ on each connected component of $G[U]$. By (C2) $G[U]$ is
connected, and by (C1) $z^*=\lambda$ is a simple eigenvalue, so
$\sum_Uv_i=u^{\T}u\ne0$ and $\eta=1/\sum_UV_i$ has order exactly $-o$.

\emph{Step 2: $o_x=-o$ for $x\in N(U)$.} Suppose $o+o_x<0$. For
$e\in\partial_x$ the function $f_e$ has order $(o+o_x)/2$ and leading
coefficient $\varepsilon_ew_eu_{i(e)}\sqrt{v_x}$, where $v_x$ is the
leading coefficient of $V_x$. In the vertex equation of $x$ the term
$(D_x-z)V_x$ has order at least $o_x>(o+o_x)/2$, and an edge from $x$ to
$y\notin U$ has $\ord(V_xV_y)>o+o_x$, so the coefficient of $t^{(o+o_x)/2}$
is $\sqrt{v_x}\sum_{e\in\partial_x}\varepsilon_ew_eu_{i(e)}$, which is
nonzero by (C3): a contradiction. Suppose $o+o_x>0$. Every edge at $x$
then has $\ord(V_xV_y)>0$, so $f_e\to f^0_e\in\{0,-1\}$ and
$f_e-f^0_e=w_e^2V_xV_y/(f_e+f^0_e+1)$. For $y\notin U$ we have
$o_y\ge o+1$, so these differences have order greater than $o+o_x$ unless
$e\in\partial_x$, and $(D_x-z)V_x$ has order at least $o_x>o+o_x$. The
constant term of the vertex equation of $x$ vanishes, and the coefficient
of $t^{o+o_x}$ is
$v_x\sum_{e\in\partial_x}(2f^0_e+1)^{-1}w_e^2v_{i(e)}$, nonzero by (C3): a
contradiction.

\emph{Step 3: no other poles.} Let $o'$ be the minimum of $o_y$ over
$y\in Y_U$ and suppose $o'<0$; let $U'\subseteq Y_U$ be the set where it is
attained. Edges from $U'$ to $N(U)$ have $\ord(V_yV_x)=o'-o>0$, and edges
from $U'$ to $Y_U\setminus U'$ have $\ord f_e>o'$, so the argument of
Step~1 applies to $U'$ and shows that $z^*$ is an eigenvalue of a signed
matrix on a subset of $U'$, contradicting (C2). Hence $o_y\ge0$ on $Y_U$.

\emph{Step 4: the limit point.} As $t\to0$ we have $\eta V_i\to u_i^2$
(with $u$ normalized), $\eta(1+2f_e)\to-2s_ew_eu_iu_j$ on $E(U)$,
$z-\lambda\to0$, $V_x/\eta\to\nu_x\ne0$ on $N(U)$ by Step~2,
$f_e\to f^*_e$ finite on $\partial U$ because $V_{i(e)}V_x$ is finite,
$f_e\to f^0_e\in\{0,-1\}$ on $\Efr(U)$ because $V_xV_y\to0$ there, and the
coordinates on $Y_U$ have finite limits $\psi_0$. The limit point solves
\eqref{gen:eq:eta-system} at $\eta=0$. Thus $(f^*,\nu)$ solves the
boundary systems with $\nu_x\ne0$, nondegenerately by (C4), and $\psi_0$
solves the system on $G[Y_U]$ at $z=\lambda$ with the constants $c(f^0)$,
nondegenerately by (C5) and Lemma~\ref{gen:lem:constants}; note
$\lambda\notin\{D_y\}$ by (C2). By Lemma~\ref{gen:lem:eta}(1)--(2), applied
with $\psi_0$ in place of $\psi$, the germ of $\Xbar$ at $P$ lies on the
unique analytic solution $\gamma(\eta)$ through this limit point.

\emph{Step 5: $\psi_0$ has no vanishing $V_y$.} Let $I$ be the set of
$y\in Y_U$ with $V_y(\psi_0)=0$, and let $E_I$ be the set of edges at $I$;
these are edges of $E(Y_U)$ or frozen edges, and their values at the limit
point lie in $\{0,-1\}$. At the limit point, the differentials of the
equations of \eqref{gen:eq:eta-system} belonging to the vertices of $I$
and to the edges of $E_I$ involve only $dV_y$ ($y\in I$) and $df_e$
($e\in E_I$): the vertex equation of $y$ contributes
$(D_y-\lambda)\,dV_y-\sum_{e\ni y}df_e$, and the equation of an edge
$e=yy'$ with $y\in I$ contributes $(1+2f_e)\,df_e-w_e^2\widehat V_{y'}\,dV_y$,
the coefficient of $d\widehat V_{y'}$ being $-w_e^2V_y=0$. Hence the
Jacobian matrix is block triangular with respect to this splitting, and
the system obtained from \eqref{gen:eq:eta-system} by fixing $V_y=0$
($y\in I$) and $f_e$ ($e\in E_I$) at their limit values and deleting the
corresponding equations has an invertible Jacobian matrix at the limit
point. Its analytic solution, completed by these constant values, solves
all of \eqref{gen:eq:eta-system}: the deleted vertex equations reduce to
$1+\sum_{e\ni y}f_e=0$, which holds at the limit point, and the deleted
edge equations to $f_e(1+f_e)=0$. By uniqueness it equals $\gamma(\eta)$.
If $I\ne\emptyset$, then $V_y\equiv0$ near $P$ for $y\in I$, which is
impossible on $\Xbar$. So $I=\emptyset$, and the limit point is the base
point of a configuration $\gamma$.

Since $z$ is not constant near $P$, $\gamma$ is non-vertical and
$P=P_\gamma$. The remaining assertions are
Lemma~\ref{gen:lem:eta}(3)--(4).
\end{proof}

\begin{lemma}[Counting]\label{gen:lem:count}
For generic coefficients the number of configurations with given $U$ is
$\abs U\,2^{\beta(U)}\,\pi(U)$.
\end{lemma}

\begin{proof}
By (C1) and (C2) each of the $2^{\beta(U)}$ switching classes has $\abs U$
admissible eigenvalues. Fix $f^0$. For $x\in N(U)$ put
$a_e=w_e^2u_{i(e)}^2$ and $k=k_x$. Eliminating $\nu_x$ from
\eqref{gen:eq:boundary} by means of one edge $e_1$ gives the equations
$a_{e_1}f_e(1+f_e)=a_ef_{e_1}(1+f_{e_1})$ ($e\ne e_1$) and $\sum_ef_e=b_x$,
of degrees $(2,\dots,2,1)$. Their leading forms have a common nontrivial
zero only if $\sum_e\varepsilon_e\sqrt{a_e}=0$ for some signs, which (C3)
excludes; so there are $2^{k-1}$ solutions counted with multiplicity. The
solutions with $\nu_x=0$ are the points of $\{0,-1\}^{\partial_x}$ with
sum $b_x$; there are $\binom k{-b_x}$ of them, and they are nondegenerate
because $\sum_e\pm a_e\ne0$ by (C3). The other solutions are nondegenerate
by (C4), so there are $a(k_x,b_x)$ of them. By (C5) and
Lemma~\ref{gen:lem:constants} there are $L_{c}(G[Y_U])$ choices of $\psi$.
\end{proof}

\begin{proof}[Proof of Theorem~\ref{thm:genus-formula}]
Shrink the generic coefficient set to that of
Lemma~\ref{gen:lem:no-vertical}. For any configuration $\gamma$, the
family in Lemma~\ref{gen:lem:eta} has all $V_i$ finite and nonzero for
sufficiently small $\eta\ne0$, by the nonzero limits of the rescaled
coordinates and the exterior variances. The points are distinct because
$\eta=1/\sum_{i\in U}V_i$. If $z$ were constant, these points would belong
to a single positive-dimensional fiber of the system with all $V_i\ne0$,
contrary to Lemma~\ref{gen:lem:no-vertical}. This argument also applies
when $\lambda=D_i$, since that lemma includes these energies.
Proposition~\ref{gen:prop:localization},
Lemma~\ref{gen:lem:pole-points} and Lemma~\ref{gen:lem:local-term} now give
$2g-2=2(n-1)\LG+\sum_\gamma\delta(\gamma)$. The sum is $\Sigma(G)$ by
Lemma~\ref{gen:lem:count}, and the formula for $R$ follows from
\eqref{gen:eq:RH}.
\end{proof}

\begin{proof}[Proof of Proposition~\ref{prop:genus-bound}]
Take real coefficients that are generic in the sense above and for
which the conclusions of Proposition~\ref{prop:generic-smooth} and
Corollary~\ref{cor:arc-circles} hold: $\Xp$ is a smooth real curve with
exactly one non-compact component and finitely many compact ones. Complex
conjugation acts on $\Xbar$ as an antiholomorphic involution; its fixed
point set $\Xbar(\R)$ is a disjoint union of circles, and their number is
at most $g+1$ by Harnack's inequality \cite{Harnack}. A positive solution
has all $V_i\ne0$, so by Lemmas~\ref{gen:lem:no-vertical}
and~\ref{gen:lem:regular} it is a smooth point of $\Xc$; near it the real
points of $\Xc$ form a one-dimensional real manifold, an open subset of
$\Xbar(\R)$, and the positive solutions are an open subset of it. Hence
$\Xp\to\Xbar(\R)$ is an injective open map. The image of a compact
component of $\Xp$ is open, compact and connected, so it is a whole
component of $\Xbar(\R)$. The image of the non-compact component is
disjoint from these and lies in a further component. So the number of
circles of $\Xp$ is at most $g$.
\end{proof}

\begin{remark}\label{gen:rem:scope}
The
points $P_\gamma$ with $\kappa(\gamma)=0$ are ramification points of $z$,
but Theorem~\ref{thm:genus-formula} gives neither the positions of the
other ramification points nor the local monodromy.
\end{remark}

\subsection{Examples and values}\label{gen:subsec:examples}

In the following examples we list, for each $U$ with $\pi(U)\ne0$, the
number $\abs U2^{\beta(U)}\pi(U)$ of configurations, $\kappa$ and $\delta$.
We use repeatedly that $a(1,b)\ne0$ only for $b\ge1$, that $a(2,-1)=0$,
$a(2,0)=1$, $a(2,1)=2$, and $a(3,-1)=1$.

\begin{example}[The cycle]\label{gen:ex:cycle}
Let $G=C_n$, $n\ge2$ ($C_2$ is the double edge), so $\LG=2$. For $U=V$
there are two switching classes, distinguished by the product of the
signs around the cycle, and $n$ eigenvalues for each: $2n$ configurations
with $N(U)=\emptyset$, $\kappa=0$, $\delta=\beta-2=-1$. There are no
others: if $U=V\setminus\{x\}$ then $k_x=2$, $b_x=-1$ and $a(2,-1)=0$, and
otherwise $N(U)$ contains a vertex $x$ with $k_x=1$ and $b_x\in\{-1,0\}$.
Since the configurations with $U=V$ are non-vertical,
\eqref{gen:eq:genus} gives $2g-2=4(n-1)-2n$, that is $g=n-1$ and $R=2n$.
The map $z$ has degree two, so each of the $2n$ points $P_\gamma$ is a
simple ramification point and there are no others. Thus $\Xbar$ is the
double cover of $\PP^1$ branched exactly over the $2n$ eigenvalues of the
two signed matrices $D+A^{s_+}$ and $D+A^{s_-}$,
\[
y^2=\det(z-D-A^{s_+})\,\det(z-D-A^{s_-}),
\]
which are in particular distinct. This is the classical curve
$y^2=\Delta(z)^2-4$ of a periodic Jacobi matrix with discriminant
$\Delta$ \cite{vMM}: the two factors are the characteristic polynomials of
the periodic and antiperiodic problems.
\end{example}

\begin{proposition}[Unicyclic graphs]\label{gen:prop:unicyclic}
Let $G$ be connected with $\beta=1$. For generic coefficients $\Xbar$ is
the double cover of $\PP^1$ branched exactly over the $2n$ eigenvalues of
$H_+=D+A^{s_+}$ and $H_-=D+A^{s_-}$, where $s_\pm$ represent the two
switching classes. In particular these eigenvalues are distinct and
$g=n-1$.
\end{proposition}

\begin{proof}
Here $\LG=2$. The configurations with $U=V$ are the $2n$ pairs of a
switching class and an eigenvalue; they have $\kappa=0$ and
$N(U)=\emptyset$, so they are non-vertical and give $2n$ distinct points
of $\Xbar$ at which $z$ is ramified (Lemma~\ref{gen:lem:pole-points}).
Hence $R\ge2n$. For the opposite inequality let the cycle have vertices
$k\in\Z/\ell$ and edges $e_k=\{k,k+1\}$ with weights $w_k$, and let $F$ be
the set of the other vertices; $G[F]$ is a forest. By
Lemma~\ref{lem:bp-correspondence} the admissible solutions correspond to
the fixed points of \eqref{eq:cavity} for the diagonal $D-z$ with nonzero
edge determinants. The messages on the edges of the forest directed
towards the cycle are determined recursively from the leaves and are
rational functions of $z$; let
$a_k(z)=D_k-z-\sum_{j\in F,\,j\sim k}w_{jk}^2/q_{j\to k}(z)$. The messages
$x_k=q_{k\to k+1}$ and $y_k=q_{k\to k-1}$ satisfy
$x_k=a_k-w_{k-1}^2/x_{k-1}$ and $y_k=a_k-w_k^2/y_{k+1}$. Write
$x_k=w_k\varphi_{k+1}/\varphi_k$ and $y_k=w_{k-1}\chi_{k-1}/\chi_k$ with
functions $\varphi,\chi$ on $\Z$; both solve the recurrence
$w_k\varphi_{k+1}+w_{k-1}\varphi_{k-1}=a_k\varphi_k$ with $\ell$-periodic
coefficients, and $\varphi_{k+\ell}=\mu\varphi_k$,
$\chi_{k-\ell}=\mu'\chi_k$ for some $\mu,\mu'\in\C^*$. The edge
determinant is
$x_ky_{k+1}-w_k^2=w_k^2(\varphi_{k+1}\chi_k-\varphi_k\chi_{k+1})/(\varphi_k\chi_{k+1})$,
so $\varphi$ and $\chi$ are linearly independent. The monodromy $M(z)$ of
the recurrence over one period has determinant one, because the Wronskian
$w_k(\varphi_k\chi_{k+1}-\varphi_{k+1}\chi_k)$ does not depend on $k$;
its eigenvalues are therefore $\mu$ and $1/\mu'=\mu^{-1}$, and
$\mu+\mu^{-1}=\Delta(z):=\tr M(z)\in\C(z)$. All messages, and hence all
$V_i$ and $f_e$, are rational functions of $z$ and $\mu$. Since $z$ has
degree two on $\Xbar$, the function field of $\Xbar$ is $\C(z,\mu)$, and
$z$ is ramified only over points where $\Delta^2-4$ has odd order. At a
pole of $\Delta$ the order of $\Delta^2-4$ is even, so these are zeros of
$\Delta-2$ or of $\Delta+2$. Expanding the determinant of a periodic
Jacobi matrix \cite{vMM} and taking the Schur complement of the forest
block, one finds
\[
\bigl(\Delta(z)\mp2\bigr)\prod_kw_k
=\det\bigl(\diag(a_k(z))-W^{\pm}\bigr)
=\frac{\det(H_\pm-z)}{\det\bigl((D+A)[F]-z\bigr)},
\]
where $W^\pm$ is the weighted adjacency of the cycle with one edge
multiplied by $\pm1$. We choose $H_\pm$ to have off-diagonal cycle part
$-W^\pm$; the signs on the attached trees do not affect the Schur
complement. This convention also works for a cycle of length two:
the two-edge contribution to the determinant is
$-(w_1\pm w_2)^2$. So the finite branch values are among the at most
$2n$ eigenvalues of $H_+$ and $H_-$, each carries one ramification point,
and $z=\infty$ is unramified by Lemma~\ref{lem:infinity}. Hence $R\le2n$.
Therefore $R=2n$, the branch values are exactly the eigenvalues of
$H_\pm$, and $g=n-1$ by \eqref{gen:eq:RH}.
\end{proof}

Comparing with Theorem~\ref{thm:genus-formula}: for a unicyclic graph the
$2n$ configurations with $U=V$ contribute $-2n$, so the non-vertical
configurations with $U\ne V$ contribute zero in total.

\begin{proposition}[Dipoles]\label{gen:prop:dipole}
Let $D_k$ be the graph with two vertices joined by $k\ge1$ parallel edges.
For generic coefficients
\[
g=1-k+(k-1)\,2^{k-1},\qquad R=(k-3)2^{k}+4\,L(D_k),\qquad L(D_k)=2^k-k .
\]
For $k=2,3,4,5$ the genus is $1,6,21,60$.
\end{proposition}

\begin{proof}
Let the vertices be $0,1$, $\sigma=\sum_ef_e$ and
$c(z)=1/((D_0-z)(D_1-z))$. An admissible solution has
$V_i=(1+\sigma)/(D_i-z)$, and \eqref{eq:system} becomes
$f_e(1+f_e)=w_e^2\,t$ with $t=V_0V_1=c(z)(1+\sigma)^2$. With
$y_e=1+2f_e$ this is $y_e^2=1+4w_e^2t$. Let $W$ be the smooth projective
model of the curve $\{(t,y):y_e^2=1+4w_e^2t\ (e\in E)\}$. For distinct
$w_e^2$ it is a connected $(\Z/2)^k$-cover of the $t$-line, branched over
the $k$ points $t=-1/(4w_e^2)$, where the inertia group is generated by
the $e$-th generator, and over $t=\infty$, where it is generated by the
product of all generators. Each of these $k+1$ fibers consists of
$2^{k-1}$ points of ramification index two, so
\[
2g_W-2=-2\cdot2^k+(k+1)2^{k-1},\qquad g_W=1+2^{k-2}(k-3).
\]
On $W$ consider the function $c_W=t/(1+\sigma)^2$,
$\sigma=\tfrac12(\sum_ey_e-k)$. Its zeros are the points over $t=0$ with
$1+\sigma\ne0$: these are the points $y\in\{\pm1\}^E$, where
$\sigma=-\abs S$ for $S=\{e:y_e=-1\}$, with $\abs S\ne1$; $t$ is a local
coordinate there, so the zeros are simple. Indeed at a point over $t=0$
with $\abs S=1$ the function $1+\sigma$ vanishes and $c_W$ has a pole, and
at a point over $t=\infty$, with local coordinate $r$, $t=r^{-2}$ and
$y_e=2\varepsilon_ew_e/r+O(r)$, we get
$c_W\to(\sum_e\varepsilon_ew_e)^{-2}\ne0$ for generic $w$. So
$c_W\colon W\to\PP^1$ has degree $2^k-k=L(D_k)$. The map
$(z,V,f)\mapsto(z,(t,y))$ identifies the admissible solutions with the
points of the fiber product of $c\colon\PP^1_z\to\PP^1$ and
$c_W\colon W\to\PP^1$ at which $t\ne0$, $1+\sigma\ne0$ and
$z\notin\{D_0,D_1,\infty\}$; hence $\Xbar$ is the normalization of this
fiber product. The map $c$ has degree two with critical values $0$ (at
$z=\infty$) and $c_{\mathrm{mid}}=-4/(D_0-D_1)^2$. Over $0$ the map $c_W$
is unramified, as we have seen, and for generic $D$ the value
$c_{\mathrm{mid}}$ is not a branch value of $c_W$, which does not depend
on $D$. Therefore $\Xbar\to W$ is a double cover branched exactly at the
$2L(D_k)$ points of $W$ over $0$ and $c_{\mathrm{mid}}$; it is connected,
and
\[
2g-2=2(2g_W-2)+2L(D_k),\qquad g=2g_W-1+L(D_k)=1-k+(k-1)2^{k-1}.
\]
The value of $R$ follows from \eqref{gen:eq:RH}.
\end{proof}

\begin{example}[The dipole $D_3$]\label{gen:ex:D3}
For $k=3$ the curve $W$ is elliptic, $c_W$ has degree $L(D_3)=5$, and
$\Xbar$ is a double cover of $W$ branched at $10$ points, so $g=6$ and
$R=20$. In terms of Theorem~\ref{thm:genus-formula}: for $U=V$ there are
$4$ switching classes and $2$ eigenvalues each, $8$ configurations with
$N(U)=\emptyset$, $\kappa=3-2=1$ and $\delta=\beta-2=0$; these are the
simple poles of Theorem~\ref{thm:local-pole}. For $U=\{v\}$ the other
vertex $x$ has $k_x=3$ and $b_x=-1$, and $a(3,-1)=1$: one configuration
over $z=D_v$ for each $v$, with $\kappa=0-1+1=0$ and
$\delta=0+2-2=0$. Hence $\Sigma(D_3)=0$, there is no configuration with
$\kappa=0$ and $\delta>0$, and \eqref{gen:eq:genus} gives
$2g-2=2\cdot1\cdot5=10$, in agreement with
Proposition~\ref{gen:prop:dipole}.
\end{example}

\begin{example}[Diamond, bowtie, $K_4$]\label{gen:ex:small}
(a) The diamond ($K_4$ minus an edge; $n=4$, $\LG=5$). For $U=V$: $16$
configurations, $\kappa=1$, $\delta=0$. For $U=\{v\}$ with $v$ of degree
two, $N(U)$ consists of the two vertices of degree three; $a(1,b_x)\ne0$
forces the three frozen edges to have value $-1$, then $b_x=1$ and the
remaining vertex has constant $-1$: one configuration for each such $v$,
with $\kappa=-3$ and $\delta=2$. For $U$ a path through one vertex of
degree three and both vertices of degree two, $N(U)$ is the other vertex
of degree three with $k_x=3$, $b_x=-1$: $2\cdot3$ configurations with
$\kappa=0$ and $\delta=0$. The other $U$ have $\pi(U)=0$. So
$\Sigma=4$, $2g-2=30+4$ and $g=18$.

(b) The bowtie (two triangles sharing a vertex; $n=5$, $\LG=4$). For
$U=V$: $20$ configurations, $\kappa=1$, $\delta=0$. For $U$ the two
vertices of one triangle other than the center $x$: $k_x=2$, and the two
frozen edges must have value $-1$ (otherwise the edge of the other
triangle admits no leafless choice or $a(2,-1)=0$), so $b_x=1$,
$a(2,1)=2$: $2\cdot2$ configurations for each triangle, with $\kappa=-2$
and $\delta=0$. So $\Sigma=0$, $2g-2=32$ and $g=17$.

(c) $K_4$ ($n=4$, $\LG=15$). $U=V$: $32$ configurations, $\kappa=2$,
$\delta=1$. $U$ a triangle: $4\cdot6$ configurations, $\kappa=1$,
$\delta=1$. $U$ an edge (the frozen edge has value $-1$): $6\cdot2$
configurations, $\kappa=-1$, $\delta=2$. $U$ a vertex (the frozen triangle
has values $-1$): $4$ configurations, $\kappa=-4$, $\delta=4$. So
$\Sigma=32+24+24+16=96$, $2g-2=90+96$ and $g=94$, $R=216$.

In (b) and (c) no configuration has $\kappa=0$; in (a) those with
$\kappa=0$ have $\delta=0$. So the values are unconditional.
\end{example}

\begin{table}[t]
\centering
\small\setlength{\tabcolsep}{4pt}\begin{tabular}{@{}lrrrrrr@{}}
\toprule
$G$ & $n$ & $\abs E$ & $\LG$ & $\Sigma(G)$ & $R$ & $g$ \\
\midrule
cycle $C_n$; unicyclic & $n$ & $n$ & $2$ & $-2n$ & $2n$ & $n-1$ \\
dipole $D_3=\Theta(1,1,1)$ & $2$ & $3$ & $5$ & $0$ & $20$ & $6$ \\
dipole $D_k$ & $2$ & $k$ & $2^k-k$ & $(k-3)2^k$ & $(k+1)2^k-4k$ & $1-k+(k-1)2^{k-1}$ \\
$\Theta(1,1,2)$ & $3$ & $4$ & $5$ & $2$ & $32$ & $12$ \\
diamond $=\Theta(1,2,2)$ & $4$ & $5$ & $5$ & $4$ & $44$ & $18$ \\
$K_{2,3}=\Theta(2,2,2)$ & $5$ & $6$ & $5$ & $6$ & $56$ & $24$ \\
theta graph $\Theta(a,b,c)$ & $s-1$ & $s$ & $5$ & $2(s-3)$ & $12s-16$ & $6(s-2)$ \\
bowtie & $5$ & $6$ & $4$ & $0$ & $40$ & $17$ \\
$C_a,C_b$ sharing a vertex & $s-1$ & $s$ & $4$ & $0$ & $8(s-1)$ & $4s-7$ \\
three triangles sharing a vertex & $7$ & $9$ & $8$ & $56$ & $168$ & $77$ \\
dumbbell & $6$ & $7$ & $5$ & $12$ & $72$ & $32$ \\
$K_4$ & $4$ & $6$ & $15$ & $96$ & $216$ & $94$ \\
$K_4$ with one edge subdivided & $5$ & $7$ & $15$ & $134$ & $284$ & $128$ \\
\midrule
prism & $6$ & $9$ & $50$ & $1032$ & $1632$ & $767$ \\
$M_1$: two $D_3$ and a bridge & $4$ & $7$ & $41$ & $324$ & $652$ & $286$ \\
$M_2$: two diamonds and a bridge & $8$ & $11$ & $50$ & $1232$ & $2032$ & $967$ \\
\bottomrule
\end{tabular}
\caption{Genus $g$ and total ramification $R$ of the spectral curve for
generic coefficients. In the rows with $a,b,c$ we put $s=a+b+c$,
respectively $s=a+b$ with $a,b\ge3$; $\Theta(a,b,c)$ consists of three
paths of lengths $a,b,c$ with common ends; the dumbbell is two triangles joined by a bridge. In $M_2$ the bridge joins
vertices of degree two.}
\label{gen:tab:values}
\end{table}

The values of $\LG$ and $\Sigma(G)$ in Table~\ref{gen:tab:values} come from
exact integer enumeration of \eqref{gen:eq:Sigma}. For the diamond,
the bowtie and $K_4$ the enumeration is reproduced by hand in
Example~\ref{gen:ex:small}. The graphs $M_1$ and $M_2$
are graphs on which $\Xp$ has circles for suitable coefficients;
Proposition~\ref{prop:genus-bound} bounds their number by $286$ and
$967$. These bounds are far from the single circle that is known in each
case; we have not investigated how sharp the bound is.

As numerical evidence (floating point, not certified), we located
ramification points of $z$ by computing the monodromy of the sheets around
the cells of a grid in the $z$-plane. Such a survey gives a lower bound
for $R$. It returned exactly the value of $R$ in the table for $D_3$,
$\Theta(1,1,2)$, the diamond, $K_{2,3}$, the bowtie, the dumbbell and
$K_4$, and in no case a larger value; for three triangles sharing a vertex
it found $139$ of the $168$ ramification points.

\subsection{Subdivision of edges}\label{gen:subsec:subdivision}

Let $\Gamma(G)=1+(n-1)\LG+\tfrac12\Sigma(G)$ be the generic genus from
Theorem~\ref{thm:genus-formula}. For an edge $e=uv$ let $G^{(e)}$ be
the graph obtained by subdividing $e$ by a new vertex $x$, and put
\begin{equation}\label{gen:eq:slope}
\Delta_e(G)=\LG+\pi_{G^{(e)}}(\{x\})
+\tfrac12\sum_{U\ni u\ \mathrm{or}\ U\ni v}2^{\beta(U)}\,\pi_G(U)\,\delta(U),
\end{equation}
the sum over the $U\subseteq V$ with $G[U]$ connected that contain $u$ or
$v$.

\begin{proposition}\label{gen:prop:subdivision}
\begin{enumerate}
\item $L(G^{(e)})=\LG$ and $\Gamma(G^{(e)})=\Gamma(G)+\Delta_e(G)$.
\item $G^{(e)}$ has a pair $(U,f^0)$ with $\kappa=0$, $\delta(U)>0$ and
\eqref{gen:eq:pi} nonzero if and only if $G$ has one.
\item $\Delta_{e'}(G^{(e)})=\Delta_e(G)$ for the two edges $e'$ into which
$e$ is subdivided, and $\Delta_{e_1}(G^{(e_2)})=\Delta_{e_1}(G)$ for
distinct edges $e_1,e_2$ of $G$.
\end{enumerate}
\end{proposition}

\begin{proof}
Write $G'=G^{(e)}$. A leafless subset of $E(G')$ contains both or none of
$ux$, $xv$, which gives $L(G')=\LG$. We compare, for each $U'\subseteq V(G')$
with $G'[U']$ connected, the number $t(U')=2^{\beta(U')}\pi_{G'}(U')$ with
the corresponding number for $G$.

(a) $x\notin U'$ and $u,v\in U'$: then $x\in N(U')$ with $k_x=2$,
$b_x=-1$, and $a(2,-1)=0$, so $t(U')=0$.

(b) $x\notin U'$ and exactly one of $u,v$ lies in $U'$: then
$x\in N(U')$ with $k_x=1$ and one frozen edge, so $b_x\in\{-1,0\}$ and
$t(U')=0$.

(c) $x,u,v\notin U'$: then $U'=U\subseteq V$, $G[U]=G'[U]$ and
$N(U)$ is the same in both graphs. If $u,v\in N(U)$, the frozen edge $e$
of $G$ is replaced by two frozen edges and an isolated vertex $x$ of
$G'[Y_{U}]$ with constant $c_x=1+f^0_{ux}+f^0_{xv}$, which must be nonzero;
the pairs $(f^0_{ux},f^0_{xv})=(0,0),(-1,-1)$ correspond to
$f^0_e=0,-1$ with the same $b_u,b_v$. If $u\in N(U)$ and $v\in Y_U$, then
$c_x=1+f^0_{ux}$, and the subsets $S'$ counted by $L_c$ satisfy
$xv\notin S'$ if $f^0_{ux}=0$ and $xv\in S'$ if $f^0_{ux}=-1$; they
correspond to the subsets counted in $G$ for $f^0_e=0$ and $f^0_e=-1$. If
$u,v\in Y_U$, then $G'[Y_U]$ is $G[Y_U]$ with $e$ subdivided and $L_c$ is
unchanged. In all cases $t(U')=t(U)$, with the same $\kappa$ and $\delta$.

(d1) $U'=\{x\}$: $N(U')=\{u,v\}$, $k_u=k_v=1$, $\delta=2$, and
$a(1,b)\ne0$ only for $b\ge1$, so $\kappa=-1-b_u-b_v\le-3$.

(d2) $x\in U'$ and exactly one of $u,v$, say $u$, lies in $U'$: then
$U=U'\setminus\{x\}$ is connected in $G$, contains $u$ and not $v$, and
$G'[U']$ is $G[U]$ with a pendant vertex. The sets $N$, $\Efr$, $Y$ and
the numbers $\beta$, $k_v$, $b$, $c$, $\kappa$, $\delta$ agree, so
$t(U')=t(U)$, while $\abs{U'}=\abs U+1$.

(d3) $x,u,v\in U'$: then $U=U'\setminus\{x\}$ contains $u,v$, $G'[U']$ is
$G[U]$ with $e$ subdivided, and again $t(U')=t(U)$, with the same
$\kappa$ and $\delta$, and $\abs{U'}=\abs U+1$.

Every connected $U\subseteq V$ arises exactly once in (c), (d2) or (d3).
Summing $\abs{U'}\,t(U')\,\delta(U')$ gives
\[
\Sigma(G')=\Sigma(G)+2\pi_{G'}(\{x\})+\sum_{U\ni u\ \mathrm{or}\ U\ni v}2^{\beta(U)}\pi_G(U)\,\delta(U),
\]
and (1) follows, since $n$ increases by one. The correspondence preserves
$\kappa$ and $\delta$, and the new configurations (d1) have $\kappa\ne0$;
this is (2).

(3) Let $e'=ux$ and let $y$ be the vertex subdividing it. Then
$\pi(\{y\})=0$ in the new graph, because $x\in N(\{y\})$ has $k_x=1$ and
one frozen edge. By the classification above, the sets $U'$ of $G'$
containing $u$ or $x$ with $t(U')\ne0$ are $\{x\}$, which contributes
$\tfrac12\cdot\pi_{G'}(\{x\})\cdot2$ to \eqref{gen:eq:slope}, and the
sets $U\cup\{x\}$ of (d2) and (d3), with $U$ containing $u$ or $v$ and
$t(U')\delta(U')=t(U)\delta(U)$. Hence $\Delta_{e'}(G')=\Delta_e(G)$. For
the second statement let $e_1=u_1v_1$ and apply the classification to the
subdivision of $e_2$ by $x_2$: it preserves $t\,\delta$ and the property
of containing $u_1$ or $v_1$, and the new set $\{x_2\}$ contains neither.
Moreover $\pi(\{x_1\})$ in the subdivision of $e_1$ is unchanged when
$e_2$ is subdivided, by case (c) applied to $U=\{x_1\}$. Hence
$\Delta_{e_1}(G^{(e_2)})=\Delta_{e_1}(G)$.
\end{proof}

\begin{corollary}[Linearity in the edge lengths]\label{gen:cor:linear}
For $\ell\in\Z_{\ge1}^E$ let $G_\ell$ be the graph obtained from $G$ by
replacing each edge $e$ by a path of length $\ell_e$. Then
\[
\Gamma(G_\ell)=\Gamma(G)+\sum_{e\in E}(\ell_e-1)\,\Delta_e(G).
\]
\end{corollary}

\begin{proof}
Subdivide one edge at a time and use
Proposition~\ref{gen:prop:subdivision}: by (3) each edge of the path
replacing $e$ has slope $\Delta_e(G)$ at every stage, by (1) each
subdivision adds this slope to $\Gamma$.
\end{proof}

The number of sheets is unchanged under subdivision, while the genus grows
linearly. For the cycle the slope is $\Delta_e=2+0-1=1$
(Example~\ref{gen:ex:cycle}). For $D_k$, the singleton supports have $\delta=0$, while the full
support contributes $2^{k-2}(k-3)$ to \eqref{gen:eq:slope}. At the
new vertex $x$ both neighbors have $k_v=1$, so the frozen edges must
include at least two of the remaining $k-1$ parallel edges.
Thus $\pi(\{x\})=2^{k-1}-k$ and $\Delta_e=2^{k-2}(k+3)-2k$,
giving $g(\Theta(a,b,c))=6(a+b+c-2)$ from $D_3$.

For the bowtie, every edge has an endpoint of degree two, so at the
new vertex $\pi(\{x\})=0$. The supports in
Example~\ref{gen:ex:small}(b) all have $\delta=0$, giving $\Delta_e=4$.
For $K_4$ the four frozen edges from $\{u,v\}$ to the other two
vertices must all take the value $-1$, and the edge between those
vertices is free, giving $\pi(\{x\})=2$. By Example~\ref{gen:ex:small}(c),
the sum in \eqref{gen:eq:slope} is
$8+4\cdot2+5\cdot2+2\cdot4=34$ (supports of sizes $4,3,2,1$).
Thus $\Delta_e=15+2+17=34$. The genera after subdivision are
$4(a+b)-7$ for two cycles of lengths $a,b\ge3$ sharing a vertex,
and $94+34\sum_e(\ell_e-1)$ for $K_4$ with edge lengths $\ell_e$.

%% file: sections/free-products.tex
\section{Free products and free convolution}\label{sec:free}

When the base graph has a single vertex, or two vertices joined by parallel
edges, the universal cover is the Cayley graph of a free product and the
Green function is the Cauchy transform of a free additive convolution. In
this situation the algebraic curve has a long history, and parts of our
picture are known. This section makes the connection precise. We work out
the running example $D_3$ (Example~\ref{ex:D3}); we describe the curve of
the free convolution of finitely atomic measures completely, including
its non-physical sheets (Theorem~\ref{thm:free-convolution-curve}); we
settle the irreducibility of the curve of Aomoto and Kato
(Theorem~\ref{thm:aomoto-kato}); and we show that the poles supported on
proper subgraphs of a general graph can be read through free convolution
with amalgamation over the diagonal algebra~$\C^V$
(Proposition~\ref{prop:operator-valued-poles}).

\paragraph{Conventions.}
Free probability uses the Cauchy transform
$\mathcal C_\mu(\omega)=\int(\omega-a)^{-1}\,d\mu(a)$ and its reciprocal
$F_\mu=1/\mathcal C_\mu$. If $\mu$ is the spectral measure of $\tH$ at a
vertex over $i$, then $G_i(z)=-\mathcal C_\mu(z)$. Consequently an atom
of mass $c$ of $\mu$ is a pole of $\mathcal C_\mu$ with residue $c$ and a
pole of $V_i=G_i$ with residue $-c$, and the poles of $V_i$ with positive
residue of Sections~\ref{sec:local-poles}--\ref{sec:classification}
appear below as poles of $\mathcal C_\mu$ with negative residue, that is,
as ``atoms of negative mass''. Gutkin~\cite{Gutkin} also uses
$(z-\tH)^{-1}$. In Subsections~\ref{free:subsec:atomic}
and~\ref{free:subsec:ak} we use the Cauchy transform; in
Subsection~\ref{free:subsec:operator} we write $g_i=-V_i$.

\paragraph{What is known and what is new.}
For a bouquet of $g$ loops Aomoto~\cite{Aomoto84} wrote the curve
explicitly and found $2^g$ real branches for large~$\abs z$. For free
products of cyclic groups with nearest-neighbor weights, Aomoto and
Kato~\cite[\S3--\S4]{AomotoKato} formed the plane curve of all branches,
computed its total degree, and gave a lower bound for the number of real
solutions at large~$z$; they did not prove irreducibility. Reading the
Green function on its Riemann surface and separating physical from
non-physical poles by the sign of the residue goes back to
Gutkin~\cite{Gutkin}: his Lemma~1 states the sign criterion (as a well
known fact), his Propositions~1 and~2 treat the free product of two
operators with two eigenvalues each (an elliptic curve, with one simple
pole over each of the four sums of eigenvalues), and his Theorem~3 treats
free powers of an operator with two eigenvalues. The atoms on the
physical sheet and their masses are given by the theorem of Bercovici and
Voiculescu~\cite{BercoviciVoiculescu}, which we use in the form restated
in \cite[Theorem~4.1]{Belinschi} and \cite[Theorem~2.1]{BercoviciWang}.
Algebraicity of free convolutions of algebraic measures, with resultant
computations and the observation that other branches may carry apparent
atoms with weight outside $[0,1]$, is in Rao and
Edelman~\cite{RaoEdelman}. The point spectrum of free products of
complete graphs was determined by Cartwright and
Soardi~\cite{CartwrightSoardi} and Kuhn~\cite{Kuhn}. For self-adjoint
convolution operators on discrete groups, Gutkin~\cite[p.~93]{Gutkin}
credits McLaughlin~\cite{McLaughlin} with the inclusion
$\sigma_p(A*B)\subseteq\sigma_p(A)+\sigma_p(B)$.
For universal covers of general graphs,
freeness with amalgamation of the edge operators and the resulting
equations for the physical branch are due to Garza-Vargas and
Kulkarni~\cite{GVK}; the atoms on the physical sheet were classified by
Banks, Garza-Vargas and Mukherjee~\cite{BGVM}; and a necessary condition
for atoms of operator-valued free convolutions is due to Belinschi,
Bercovici and Liu~\cite{BBL}.

To the best of our knowledge the following points are new: for arbitrary
finitely atomic measures, the exact number of sheets, irreducibility for
generic data, the genus, the residue $\sum_jp_j-(k-1)$ at every tuple of
atoms on every sheet together with the ramification when this number
vanishes, and the number of real circles
(Theorem~\ref{thm:free-convolution-curve}); the generic irreducibility
and the exact degree for the model of Aomoto and Kato
(Theorem~\ref{thm:aomoto-kato}); and, in the operator-valued setting, the
poles on non-physical sheets with their directional residues and
rank-one residue matrices
(Proposition~\ref{prop:operator-valued-poles}). The results of
Sections~\ref{sec:local-poles}--\ref{sec:classification} extend this
picture from free products to arbitrary base graphs.

\subsection{The running example}\label{free:subsec:D3}

\begin{example}[The dipole $D_3$]\label{ex:D3}
Let $D_3$ be the graph with two vertices $1,2$ and three parallel edges
with positive weights $w_1,w_2,w_3$, and let $D_1,D_2$ be the potentials. The
leafless subsets are $\emptyset$, the three pairs of edges and the full
edge set, so $L(D_3)=5$; moreover $\beta=2$ and $\kappa_0=1$.

\emph{Closed form.} Put $P=V_1V_2$, $y_e=1+2f_e$ and
$\gamma=1+f_1+f_2+f_3$. The system \eqref{eq:system} is equivalent to
\begin{equation}\label{free:eq:D3}
y_e^2=1+4w_e^2P\quad(e=1,2,3),\qquad
\gamma=\tfrac12(y_1+y_2+y_3-1),\qquad (D_1-z)(D_2-z)=\frac{\gamma^2}{P},
\end{equation}
together with $V_i=\gamma/(D_i-z)$. Let $Y_0$ be the curve defined by the
first three equations; for pairwise distinct $w_e^2$ it is an irreducible
$(\Z/2)^3$-cover of the $P$-line of degree~$8$, branched over
$P=-1/(4w_e^2)$ and $P=\infty$, and of genus~$1$. The spectral curve is
the double cover of $Y_0$ on which
$\Delta=(D_1-D_2)^2+4\gamma^2/P$ has a square root. On the physical
branch all $y_e$ are positive for real $z<\Sigm$.

\emph{Sheets.} As $z\to\infty$ the last equation of \eqref{free:eq:D3}
forces $P\to0$ with $\gamma\not\to0$. Over $P=0$ the curve $Y_0$ has the
eight points $y=(\pm1,\pm1,\pm1)$, at which $P$ is a local coordinate.
The value of $\gamma$ is $1$ at $(+,+,+)$, $-1$ at the three points with
two entries equal to $-1$, and $-2$ at $(-,-,-)$; at each of these five
points $\gamma^2/P$ has a simple pole, the double cover is ramified, and
$z$ has a simple pole in the local coordinate $P^{1/2}$. At the three
points with exactly one entry $-1$ we have $\gamma=0$, so that
$\gamma^2/P\to0$ and $z\to D_1$ or $D_2$; these points are not
admissible. Thus $z$ has degree $5$, and the five sheets at infinity are
indexed by the leafless subsets $S$ through $f_e\to-1$ for $e\in S$ and
$f_e\to0$ for $e\notin S$, as in Lemma~\ref{lem:infinity}. For generic
coefficients the curve is irreducible
(Theorem~\ref{thm:irreducible}). The function $\Delta$ on $Y_0$ has
exactly the five simple poles just described, hence five zeros, and the
double cover is branched at the poles and at the zeros of odd order;
therefore the genus is at most~$6$, with equality when the five zeros
are simple. Theorem~\ref{thm:genus-formula} gives genus $6$ for generic
coefficients. For $D_1=D_2$ the function $4\gamma^2/P$ has a double zero
and the genus drops to~$5$ (Corollary~\ref{free:cor:dipole}).

\emph{Poles.} At finite $z$ the differentials $V_i\,dz$ have poles only
over $P=\infty$. (At the zero of $\gamma$ on $Y_0$ with
$P\ne0$, which lies over $z=D_1$ and $z=D_2$ when $D_1\ne D_2$, the
function $z$ is ramified, one of the $V_i$ grows like
$(z-D_i)^{-1/2}$, and $V_i\,dz$ is holomorphic.) The curve $Y_0$ has four points over $P=\infty$, indexed by
the sign patterns $\varepsilon\in\{\pm1\}^3$ of
$y_e\sim2\varepsilon_e\abs{w_e}P^{1/2}$ modulo a simultaneous change of
sign, that is, by the four switching classes of signings. With
$t=P^{-1/2}$ and $W_\varepsilon=\sum_e\varepsilon_e\abs{w_e}$ we have
$\gamma=W_\varepsilon/t-\tfrac12+O(t)$, hence
$(D_1-z)(D_2-z)=W_\varepsilon^2-W_\varepsilon t+O(t^2)$. If
$W_\varepsilon\ne0$, the two points of the spectral curve over this point
of $Y_0$ lie over
\begin{equation}\label{free:eq:D3-poles}
z=\lambda=\frac{D_1+D_2}2\pm
\sqrt{\Bigl(\frac{D_1-D_2}2\Bigr)^2+(w_1\pm w_2\pm w_3)^2},
\end{equation}
the eigenvalues of the signed Jacobi matrices
$\bigl(\begin{smallmatrix}D_1&W_\varepsilon\\ W_\varepsilon&D_2\end{smallmatrix}\bigr)$.
At each of them $dz/dt=W_\varepsilon/(D_1+D_2-2\lambda)\ne0$, so $z$ is
a local coordinate, and
\[
\lim_{z\to\lambda}(z-\lambda)V_1
=\frac{W_\varepsilon}{D_1-\lambda}\cdot\frac{dz}{dt}
=\frac{D_2-\lambda}{D_1+D_2-2\lambda}=u_1^2,
\]
where $u$ is the unit eigenvector for $\lambda$; similarly for $V_2$.
For generic coefficients this gives the eight simple poles with
residue $\kappa_0u_i^2=u_i^2$ predicted by
Theorem~\ref{thm:local-pole}. The residues are positive, so all eight
poles are antibound (Corollary~\ref{cor:antibound-sign}). The choice
$W_\varepsilon=w_1+w_2+w_3$ gives the two Perron poles of
Theorem~\ref{thm:perron-component}.

\emph{Uniform weights.} If $w_1=w_2=w_3=w$, the physical component lies
in $y_1=y_2=y_3$. It has two sheets, with $y\to\pm1$ at infinity. The
only finite points on its normalization with nonzero residue of
$(V_1+V_2)\,dz$ are the two Perron poles at
$\tfrac12(D_1+D_2)\pm(\tfrac14(D_1-D_2)^2+9w^2)^{1/2}$.
The other three switching classes give the same pair of eigenvalues,
corresponding to $W_\varepsilon^2=w^2$. If $D_1\ne D_2$, the physical
component also has a branch point over each $D_i$. There $V_i$ diverges
as a nonzero multiple of $(D_i-z)^{-1/2}$ and $V_i\,dz$ has zero residue,
as in Proposition~\ref{prop:uniform}.

\emph{Attribution.} For $D_1=D_2=0$ equations \eqref{free:eq:D3} reduce,
with $W=-1/V_1$, to
$2z=\sum_e\pm\sqrt{W^2+4w_e^2}-W$, which is the curve of the weighted
free product $\Z/2*\Z/2*\Z/2$ in \cite[\S4, Example~1]{AomotoKato}.
Example~3.3 of \cite{CSZ} is the specialization to equal weights and
potential $\pm b$, with two sheets. The curve of the running example at
zero potential is therefore a known object; its irreducibility and exact
number of sheets follow from Theorem~\ref{thm:aomoto-kato}, and the
genus, the residues and the correspondence with switching classes are,
as far as we know, not in the literature.
\end{example}

\subsection{Free convolution of finitely atomic measures}
\label{free:subsec:atomic}

Let $\mu_1,\dots,\mu_k$, $k\ge2$, be probability measures on $\R$ with
finite support,
\[
\mu_j=\sum_{\ell=1}^{m_j}p_{j\ell}\,\delta_{a_{j\ell}},\qquad
a_{j1}<\dots<a_{jm_j},\quad p_{j\ell}>0,\quad\sum_\ell p_{j\ell}=1,
\]
and put $F_j=F_{\mu_j}$, $M=\prod_jm_j$. Let
$\mu=\mu_1\boxplus\dots\boxplus\mu_k$ be the free additive convolution.

\begin{lemma}\label{free:lem:F}
The function $F_j$ is a rational map $\PP^1\to\PP^1$ of degree $m_j$
with $F_j(\omega)=\omega+O(1)$ at infinity. The fiber of $F_j$ over
every point of $\R\cup\{\infty\}$ consists of $m_j$ distinct points of
$\R\cup\{\infty\}$. In particular all branch values of $F_j$ are finite
and non-real. Moreover $F_j^{-1}(0)=\{a_{j\ell}\}$ and
$F_j'(a_{j\ell})=1/p_{j\ell}$.
\end{lemma}

\begin{proof}
For non-real $\omega$,
$\im\mathcal C_{\mu_j}(\omega)=-\im\omega\sum_\ell p_{j\ell}\abs{\omega-a_{j\ell}}^{-2}$,
so $F_j$ maps non-real points to non-real points. On the real line
$\mathcal C_{\mu_j}'<0$ away from the atoms, and $\mathcal C_{\mu_j}$
decreases from $+\infty$ to $-\infty$ between consecutive atoms; this
gives one simple zero in each of the $m_j-1$ gaps, and there are no
other zeros because the numerator has degree $m_j-1$. Hence the poles
of $F_j$ are these zeros and $\infty$, all simple, and
$F_j'=-\mathcal C_{\mu_j}'/\mathcal C_{\mu_j}^2>0$ at every real point
where $F_j$ is finite, including the atoms, where
$F_j(\omega)=(\omega-a_{j\ell})/p_{j\ell}+O((\omega-a_{j\ell})^2)$.
Thus the restriction of $F_j$ to the real projective line is an
unramified, orientation preserving covering of degree $m_j$ of the real
projective line, and its fibers exhaust the fibers of~$F_j$.
\end{proof}

Let $Y=Y(\mu_1,\dots,\mu_k)$ be the normalization of the reduced fiber
product of the maps $F_j\colon\PP^1_{\omega_j}\to\PP^1_F$. It is a
smooth projective curve, possibly disconnected, with a map
$F\colon Y\to\PP^1$ of degree $M$, and every irreducible component of
$Y$ maps onto $\PP^1_F$. On $Y$ we consider the rational functions
\begin{equation}\label{free:eq:z}
z=\sum_{j=1}^k\omega_j-(k-1)F\qquad\text{and}\qquad 1/F .
\end{equation}
For a tuple of atoms
$\alpha=(a_{1\ell_1},\dots,a_{k\ell_k})$ let
\[
\lambda_\alpha=\sum_ja_{j\ell_j},\qquad
c_\alpha=\sum_jp_{j\ell_j}-(k-1).
\]

\begin{theorem}\label{thm:free-convolution-curve}
With the notation above, the following hold.
\begin{enumerate}
\item[\textup{(a)}] \emph{Sheets.} The poles of $z$ lie over $F=\infty$.
Of the $M$ points of $Y$ over $F=\infty$, exactly $\sum_j(m_j-1)$ are
not poles of $z$, and the others are simple poles. Hence the degrees of
$z$ on the components of $Y$ add up to
\[
N=\prod_jm_j-\sum_j(m_j-1).
\]
\item[\textup{(b)}] \emph{Physical component.} There is exactly one
point $o\in Y$ with $F(o)=\infty$ and $z/F\to1$. Near $o$ the function
$z$ is a local coordinate and $1/F=\mathcal C_\mu(z)$. If
$Y_{\mathrm{phys}}$ denotes the component of $Y$ containing $o$, then
$\C(Y_{\mathrm{phys}})=\C(z,1/F)$; thus the algebraic degree of
$\mathcal C_\mu$ over $\C(z)$ equals the degree of $z$ on
$Y_{\mathrm{phys}}$, which is at most~$N$.
\item[\textup{(c)}] \emph{Irreducibility and genus.} If the sets of
branch values of $F_1,\dots,F_k$ are pairwise disjoint, then $Y$ is
irreducible, the algebraic degree of $\mathcal C_\mu$ is exactly $N$,
and the genus of $Y$ is
\[
1+M\Bigl(k-1-\sum_{j=1}^k\frac1{m_j}\Bigr).
\]
The hypothesis holds for all positions and masses of the atoms in a
nonempty Zariski open subset of the parameter space; for given measures
it holds after replacing each $\mu_j$ by its dilation by a factor
$t_j>0$, for all $(t_j)$ outside finitely many hyperplanes
$t_ib=t_jb'$.
\item[\textup{(d)}] \emph{Poles on all sheets.} The points of $Y$ over
$F=0$ are in bijection with the tuples of atoms $\alpha$. At the point
$\alpha$ the function $F$ is a local coordinate and
\[
z=\lambda_\alpha+c_\alpha F+O(F^2).
\]
Hence, if $c_\alpha\ne0$, the function $1/F$ has at this point a simple
pole at $z=\lambda_\alpha$ with residue $c_\alpha$ with respect to~$z$.
If $c_\alpha=0$ and $z$ is not constant on the component through
$\alpha$, then $z$ is ramified at $\alpha$ and the differential
$F^{-1}dz$ is holomorphic there. The function $1/F$ has no other poles.
The measure $\mu$ has an atom at $\lambda$ if and only if
$\lambda=\lambda_\alpha$ for a tuple with $c_\alpha>0$; this tuple is
then unique, the mass of the atom is $c_\alpha$, and the physical branch
tends to the point $\alpha$ as $z\to\lambda$ from the upper half-plane.
A point $\alpha$ with $c_\alpha<0$ is not a limit of the physical
branch.
\item[\textup{(e)}] \emph{Real circles.} The preimage in $Y$ of the real
projective line $F\in\R\cup\{\infty\}$ is the set of real points of~$Y$.
It consists of $M/\lcm(m_1,\dots,m_k)$ disjoint circles, each of which
covers the real $F$-line $\lcm(m_1,\dots,m_k)$ times and contains that
many tuples of atoms.
\end{enumerate}
\end{theorem}

\begin{proof}
(a) Over $F=\infty$ all $F_j$ are unramified by Lemma~\ref{free:lem:F},
so the fiber product is smooth there, with $M$ points at which $1/F$ is
a local coordinate. At such a point each $\omega_j$ is either $\infty$,
in which case $\omega_j=F+O(1)$, or one of the $m_j-1$ finite poles
of~$F_j$. If $r$ is the number of indices with $\omega_j=\infty$, then
$z=(r-k+1)F+O(1)$. The coefficient vanishes exactly when $r=k-1$, that
is, when a single factor is at a finite pole; there are
$\sum_j(m_j-1)$ such points. At finite $F$ all $\omega_j$ are finite.
A component on which $z$ is constant contains no pole of $z$, so the
count of poles gives the sum of the degrees.

(b) The point $o$ is the one with $r=k$, and
$z=F+O(1)$ there. By the subordination theorem (see
\cite{BercoviciWang,Belinschi} for two measures; the case of $k$
measures follows by induction, composing the subordination functions of
$\mu_1\boxplus(\mu_2\boxplus\dots\boxplus\mu_k)$) there are analytic
functions $\omega_j$ on the upper half-plane with
\[
F_j(\omega_j(z))=F_\mu(z),\qquad
\sum_j\omega_j(z)=z+(k-1)F_\mu(z),\qquad \omega_j(iy)/iy\to1 .
\]
Hence $z\mapsto(\omega_1(z),\dots,\omega_k(z))$ is a map to the fiber
product which tends to $o$ as $z=iy\to\infty$, and on which the function
\eqref{free:eq:z} takes the value $z$. Since $z$ is a local coordinate
at $o$, this map parametrizes a neighborhood of $o$, and
$1/F=\mathcal C_\mu(z)$ there. Let $Y'$ be the smooth projective curve
with function field $\C(F,z)\subseteq\C(Y_{\mathrm{phys}})$ and
$\phi\colon Y_{\mathrm{phys}}\to Y'$ the corresponding map. Every point
of $\phi^{-1}(\phi(o))$ satisfies $F=\infty$ and $z/F\to1$, so it is
$o$; and $\phi$ is unramified at $o$, because $F$ is. Hence
$\deg\phi=1$.

(c) Let $B_j$ be the set of branch values of $F_j$ and $B=\bigcup_jB_j$.
The monodromy of $F_j$ over $\PP^1\setminus B_j$ is generated by loops
around the points of $B_j$ and is transitive, because
$\PP^1_{\omega_j}$ is connected. A loop in $\PP^1\setminus B$ around a
point of $B_j$ acts trivially on the fibers of $F_i$ for $i\ne j$, since
$B_i\cap B_j=\emptyset$. Therefore the monodromy group of the fiber
product contains the product of the monodromy groups of the factors and
is transitive on the product of the fibers. So $Y$ is connected,
$Y=Y_{\mathrm{phys}}$, and (a) and (b) give the degree. Over a point of
$B_j$ only the factor $F_j$ is ramified, so each ramification point of
$F_j$ of index $e$ gives $M/m_j$ points of $Y$ of index $e$. As the
total ramification of $F_j$ is $2m_j-2$, the Riemann--Hurwitz formula
for $F\colon Y\to\PP^1$ reads
$2g-2=-2M+\sum_j(M/m_j)(2m_j-2)$. For the last assertion, the dilation
of $\mu_j$ by $t_j$ has $F$-transform $t_jF_j(\omega/t_j)$, whose
branch values are those of $F_j$ multiplied by $t_j$; they are nonzero
by Lemma~\ref{free:lem:F}, so a collision between two factors is a
condition $t_ib=t_jb'$ with $b\in B_i$, $b'\in B_j$. The condition is
the nonvanishing of the resultants of the discriminants of the
polynomials obtained from $F_j(\omega)=F$ by clearing denominators, so
it is Zariski open in the positions and masses.

(d) By Lemma~\ref{free:lem:F} the fiber of $F_j$ over $0$ is the set of
atoms and is unramified, with local inverse
$\omega_j=a_{j\ell}+p_{j\ell}F+O(F^2)$. Substituting into
\eqref{free:eq:z} gives the expansion, and the statements on the pole
and on the differential follow from it. For the last statement, let
$x_1,\dots,x_k$ be free self-adjoint elements of a tracial
$W^*$-probability space with distributions $\mu_j$, and let $e_j$ be the
spectral projection of $x_j$ for the atom $a_{j\ell_j}$. The trace of
$e_1\wedge\dots\wedge e_k$ is at least $\sum_jp_{j\ell_j}-(k-1)$, and
the range of this projection lies in the kernel of
$\sum_jx_j-\lambda_\alpha$. Hence $c_\alpha>0$ implies
$\mu(\{\lambda_\alpha\})\ge c_\alpha$. Conversely, suppose $\mu$ has an
atom of mass $c>0$ at $\lambda$. Then
$i\eta\,\mathcal C_\mu(\lambda+i\eta)\to c$, so
$F_\mu(\lambda+i\eta)\to0$. The inverse branches of $F_j$ over a small
disc around $0$ are holomorphic with disjoint images, and the
continuous function $\omega_j(\lambda+i\eta)$ takes values in their
union for small $\eta$, hence converges to an atom of~$\mu_j$. So the
physical branch converges to a point $\alpha'$ of $Y$ over $F=0$, and
the expansion at $\alpha'$ gives $\lambda=\lambda_{\alpha'}$ and
$c=c_{\alpha'}$. Finally, two distinct tuples $\alpha,\alpha'$ with
$\lambda_\alpha=\lambda_{\alpha'}$ differ in $d\ge2$ factors, and in
each of these the two masses add up to at most~$1$; hence
$c_\alpha+c_{\alpha'}\le d+2(k-d)-2(k-1)=2-d\le0$, and at most one of
them is positive. This identifies the tuple reached by the physical
branch. If the physical branch tended to a point $\alpha$ with
$c_\alpha<0$, the same expansion would give
$i\eta\,\mathcal C_\mu(\lambda_\alpha+i\eta)\to c_\alpha<0$, which is
impossible for a positive measure. For two measures the statement about
the atoms of $\mu$ is the theorem of Bercovici and
Voiculescu~\cite{BercoviciVoiculescu,Belinschi,BercoviciWang}.

(e) By Lemma~\ref{free:lem:F} every point of $Y$ over a real value of
$F$ has real coordinates and is unramified, so the preimage of the real
$F$-line is the fiber product of the coverings
$F_j\colon\R\cup\{\infty\}\to\R\cup\{\infty\}$. Each of these is a
connected covering of the circle of degree $m_j$, and one turn around
the base permutes its fiber cyclically. On the product of the fibers one
turn acts by the simultaneous cyclic shift, whose orbits have length
$\lcm(m_j)$.
\end{proof}

\begin{remark}\label{free:rem:circles}
The circles of (e) are the real points of the curve for its natural real
structure. They are not the circles of positive solutions of
Proposition~\ref{prop:circles}: no positivity is imposed, and a circle
of (e) is in general not mapped injectively to the $z$-line.
\end{remark}

\begin{example}\label{free:ex:known-cases}
(i) For $k=2$ and $m_1=m_2=2$ the theorem gives $N=2$ and genus~$1$,
with poles over the four sums of atoms. This is the elliptic curve of
\cite[Propositions~1 and~2]{Gutkin}.
(ii) Let all $\mu_j$ be equal to $p\,\delta_a+q\,\delta_b$. The branch
values coincide, $Y$ is reducible, and $z$ is constant on some
components. The physical component is the diagonal
$\omega_1=\dots=\omega_k$; it meets $F=0$ in the two constant tuples,
where (d) gives the residues $kp-(k-1)=1-kq$ and $1-kp$, in agreement
with \cite[Theorem~3]{Gutkin}.
(iii) For the Cayley graph of $\Z/2*\Z/3$ with generators $a,b,b^{-1}$
the adjacency operator has distribution $\mu_1\boxplus\mu_2$ with
$\mu_1=\tfrac12(\delta_1+\delta_{-1})$ and
$\mu_2=\tfrac13\delta_2+\tfrac23\delta_{-1}$. The branch values are
$\pm2i$ and $1\pm2\sqrt2\,i$, so the curve is irreducible with $N=2$
sheets and genus~$1$. Over $F=0$ we find $c=\tfrac16$ at
$\lambda=0,-2$ and $c=-\tfrac16$ at $\lambda=3,1$. The first two are
the eigenvalues found in \cite[Theorem~4 and the first example following it]{Gutkin};
the other two are poles on the second sheet.
\end{example}

\begin{corollary}[Dipoles with zero potential]\label{free:cor:dipole}
Let $D_k$, $k\ge2$, be the graph with two vertices and $k$ parallel
edges with real weights of pairwise distinct absolute values, and let
the potential be zero. Then $G_1=G_2=-1/F$ on the curve
$Y(\mu_1,\dots,\mu_k)$ for the symmetric Bernoulli measures
$\mu_j=\tfrac12(\delta_{w_j}+\delta_{-w_j})$. This curve is irreducible,
has $2^k-k=L(D_k)$ sheets over the $z$-line and genus
$1+2^{k-1}(k-2)$. The points at which $G_i$ is infinite are the $2^k$
points over $F=0$; they lie over the values
$\sum_j\varepsilon_jw_j$, $\varepsilon\in\{\pm1\}^k$, which are the
eigenvalues of all signed matrices $A^s$. For $k\ge3$ each of them is a
simple pole of $G_i$ with residue $\tfrac k2-1=\kappa_0u_i^2$ and is
antibound. For $k=2$ they are the four ramification points over the
band edges $\pm w_1\pm w_2$. The real points of the curve form
$2^{k-1}$ circles, in bijection with the switching classes of signings:
the circle of the class of $s$ contains the two points over
$\pm\sum_js_jw_j$.
\end{corollary}

\begin{proof}
The universal cover of $D_k$ is the Cayley graph of the free product
$\Gamma$ of $k$ copies of $\Z/2$ with generators $a_1,\dots,a_k$, where
$a_j$ corresponds to the lifts of the $j$th edge, and $\tH$ is the
operator $\sum_jw_ja_j$ in the regular representation. The $a_j$ are
free with respect to the trace of $\Gamma$ and $w_ja_j$ has distribution
$\mu_j$, so the spectral measure of $\tH$ at any vertex is
$\mu_1\boxplus\dots\boxplus\mu_k$. Here $F_j(\omega)=\omega-w_j^2/\omega$
has the branch values $\pm2iw_j$, which are pairwise distinct for
different $j$. Theorem~\ref{thm:free-convolution-curve} applies with
$m_j=2$, $M=2^k$, $\lcm(m_j)=2$ and $c_\alpha=1-\tfrac k2$; the sign of
the residue changes with $G_i=-\mathcal C_\mu$. A circle of (e) carries
two tuples of atoms, $\varepsilon$ and $-\varepsilon$, because one turn
around the real $F$-line exchanges the two real preimages in every
factor. The matrix $A^s$ has the eigenvalues $\pm\sum_js_jw_j$ with
$u_i^2=\tfrac12$, and the switching class of $s$ is $\{s,-s\}$.
\end{proof}

The dictionary with the covariance system is $\omega_j=(1+f_j)F$ and
$V_1=V_2=-1/F$: the equation $F_j(\omega_j)=F$ becomes
$f_j(1+f_j)=w_j^2V_1V_2$, and \eqref{free:eq:z} becomes the vertex
equation $-zV_i=1+\sum_jf_j$. The equation of the curve is that of
\cite[\S4, Example~1]{AomotoKato}, and for a bouquet that of
\cite{Aomoto84}. Note that the number of sheets is $L(D_k)$ although a
constant potential is not generic in the sense of
Theorem~\ref{thm:generic-degree}, while the genus is smaller than for
generic potentials (Example~\ref{ex:D3}). For $k=3$ the four real
circles correspond to the four switching classes of
Example~\ref{ex:D3}. Weights with repeated absolute values are treated
in Proposition~\ref{prop:dipole-inverse}.

\subsection{The model of Aomoto and Kato}\label{free:subsec:ak}

Aomoto and Kato~\cite{AomotoKato} consider the free product
$\Gamma=\Z/n_1*\dots*\Z/n_k$, $n_j\ge2$, with generators $a_j$ of the
factors, and the nearest-neighbor operator
$\sum_jp_j(a_j+a_j^{-1})$ with $p_j>0$ and $2\sum_jp_j=1$. The $j$th
summand has the distribution $\mu_{j,p_j}$ of
$2p_j\cos(2\pi\ell/n_j)$ for $\ell$ uniform in $\Z/n_j$; it has
$m_j=\lfloor n_j/2\rfloor+1$ atoms. In their notation $W=1/\mathcal C$
is our $F$, and their equations (1.16)--(1.18) and (3.1) are
$F_{\mu_{j,p_j}}(\omega_j)=W$ and \eqref{free:eq:z} with
$X_j=\omega_j-W$. In \S3 they form the plane curve in the coordinates
$(W,z)$ from all tuples of branches, find the total degree
$\prod_jm_j$ of the defining product, and show that for large real $z$
there are at least $\prod_jm_j-\sum_j(m_j-1)$ real solutions. They do
not prove that this curve is irreducible.

\begin{theorem}\label{thm:aomoto-kato}
Fix $n_1,\dots,n_k\ge2$. There is a finite union $\Lambda$ of
hyperplanes of the form $p_ib=p_jb'$ in the space of weights such that
for all $p$ with $p_j>0$, $2\sum_jp_j=1$ and $p\notin\Lambda$, the
following hold: the reduced plane curve in $(W,z)$ formed from all
branches is irreducible; it is birational to
$Y(\mu_{1,p_1},\dots,\mu_{k,p_k})$; and the algebraic degree over
$\C(z)$ of the Green function of
$\sum_jp_j(a_j+a_j^{-1})$ at the identity is exactly
\[
\prod_jm_j-\sum_j(m_j-1),\qquad m_j=\lfloor n_j/2\rfloor+1 .
\]
The set of such $p$ is nonempty. For $k=3$, $n_1=n_2=n_3=2$ and
$p_1=p_2=p_3=\tfrac16$, the reduced curve of all branches is the union
of two distinct conics.
\end{theorem}

\begin{proof}
The measure $\mu_{j,p}$ is the dilation of $\mu_{j,1}$ by $p$, so the
branch values of $F_{\mu_{j,p}}$ are those of $F_{\mu_{j,1}}$ multiplied
by~$p$. Let $\Lambda$ be the union of the hyperplanes $p_ib=p_jb'$,
$i\ne j$, with $b,b'$ branch values of $F_{\mu_{i,1}}$ and
$F_{\mu_{j,1}}$. These values are non-real by Lemma~\ref{free:lem:F}, so
each hyperplane meets the set of positive ratios $p_i/p_j$ in at most
one point, and the normalization $2\sum_jp_j=1$ does not restrict the
ratios; thus the complement of $\Lambda$ in the open simplex is
nonempty. For $p\notin\Lambda$,
Theorem~\ref{thm:free-convolution-curve}(c) shows that $Y$ is
irreducible and that the degree is as stated, and (b) shows that
$\C(Y)=\C(W,z)$. The plane curve of all branches is the image of $Y$
under $(W,z)$, so it is irreducible and birational to~$Y$.

In the example each $\mu_j$ is $\tfrac12(\delta_a+\delta_{-a})$ with
$a=\tfrac13$, and $\omega_j=\tfrac12(W\pm\sqrt{W^2+4a^2})$. Summing over
the choices of signs gives
$z+\tfrac12W=\pm\tfrac32\sqrt{W^2+4a^2}$ or
$\pm\tfrac12\sqrt{W^2+4a^2}$, that is, the two conics
$(2z+W)^2=9(W^2+4a^2)$ and $(2z+W)^2=W^2+4a^2$. The physical branch
lies on the first.
\end{proof}

Thus the lower bound of \cite[(3.3)]{AomotoKato} is the exact number of
sheets over the $z$-line for generic weights, and it must be
distinguished from the total degree $\prod_jm_j$ of the plane curve. For
$n_j\ge3$ the Cayley graph is not a tree, so this model lies outside
the setting of Section~\ref{sec:setup}; for $n_j=2$ it is the dipole of
Corollary~\ref{free:cor:dipole}.

\subsection{Poles supported on subgraphs and amalgamation over the
diagonal}\label{free:subsec:operator}

Let $G$ be connected and loopless with real nonzero weights, and write
$\tH=D+A_{\tree}$, where $A_{\tree}$ is the lifted weighted adjacency
operator. Let $B=\C^V$ act on $\ell^2(\tree)$ by multiplication with
$b_{\pi(x)}$. For a bounded operator $T$ commuting with the deck
transformations let $E(T)\in B$ be given by
$E(T)_i=\langle\delta_x,T\delta_x\rangle$, $x\in\pi^{-1}(i)$. For
$e\in E$ let $X_e$ be the part of $A_{\tree}$ supported on the lifts of
$e$, so that $A_{\tree}=\sum_eX_e$; parallel edges give different
operators. The spectral argument is now an element $b\in B$, and the
spectral curve is the restriction to the line $b_i=z-D_i$. For
$\im b_i>0$ we set $g=E((b-A_{\tree})^{-1})$, so that $g_i=-G_i(z)$ on
that line.

\begin{lemma}\label{free:lem:amalgam}
\textup{(i)} The algebras generated by $B$ and $X_e$, $e\in E$, are free
with amalgamation over $B$ with respect to~$E$.

\textup{(ii)} For $\im b_i>0$ the vector $g$, together with suitable
$f\in\C^E$, satisfies
\begin{equation}\label{free:eq:op-system}
f_e(1+f_e)=w_e^2\,g_ig_j\quad(e=ij),\qquad
b_ig_i=1+\sum_{e\ni i}f_e\quad(i\in V),
\end{equation}
which for $b_i=z-D_i$ and $g_i=-V_i$ is the system \eqref{eq:system}.
The edge equation says that $g$ is the $B$-valued Cauchy transform of
$X_e$ at the argument $\omega_e\in B$ with
$\omega_{e,i}=(1+f_e)/g_i$ at the two ends of $e$ and
$\omega_{e,l}=1/g_l$ elsewhere, and the vertex equation is the
additivity of the $B$-valued $R$-transform,
$b=g^{-1}+\sum_e(\omega_e-g^{-1})$.
\end{lemma}

\begin{proof}
(i) The lifts of $e$ form a matching of $\tree$ and
$X_e^2\in B$, so every element of the algebra generated by $B$ and $X_e$
is the sum of an element of $B$ and an operator whose matrix is
supported on the lifts of~$e$; the latter is its part of expectation
zero. In a product of such centered elements for edges
$e_1\ne e_2\ne\dots\ne e_r$, a diagonal matrix entry is a sum over
closed walks in $\tree$ whose consecutive steps are lifts of different
edges. Such a walk does not backtrack and therefore cannot be closed.
(ii) The derivation of \eqref{eq:m-recursion} and of
Proposition~\ref{prop:physical-branch} uses only that the diagonal
operator has positive imaginary part; it applies to $b-A_{\tree}$. The
Cauchy transform of $X_e$ at $\omega\in B$ is the diagonal of the inverse
of a direct sum of $2\times2$ blocks, with value
$\omega_j/(\omega_i\omega_j-w_e^2)$ at $i$ and $1/\omega_l$ at $l\ne i,j$;
equating it with $g$ gives $1+f_e=(1+f_e)^2-w_e^2g_ig_j$.
\end{proof}

Freeness of the edge operators with amalgamation is due to
\cite{GVK}, where the amalgamation is over the matrix algebra
$M_n(\C)$ and the equations for the physical branch are derived from it.
We use \eqref{free:eq:op-system} as an algebraic system on the set where
all $g_i\ne0$, with independent complex arguments $b_i$.

\begin{lemma}\label{free:lem:multi}
\textup{(i)} On every solution of \eqref{free:eq:op-system} with fixed
weights,
$\sum_ig_i\,db_i=d\log\Phi$ with
$\Phi=\prod_e(1+f_e)/\prod_ig_i$. If the weights vary as well, the
right-hand side is $d\log\Phi+2\sum_e(f_e/w_e)\,dw_e$.

\textup{(ii)} If $w_e^2\ne w_{e'}^2$ for parallel edges $e\ne e'$, the
set of solutions $(b,g,f)$ of \eqref{free:eq:op-system} with all
$g_i\ne0$ is an irreducible variety of dimension $n$. Its projection to
$g$ has degree $2^{\abs E}$, and its projection to $b$ is dominant and
generically finite. Hence the physical germ can be continued to every
solution over a general $b$ along paths in $\C^V$.
\end{lemma}

\begin{proof}
(i) Differentiating the vertex equation gives
$g_i\,db_i=\sum_{e\ni i}df_e-(1+\sum_{e\ni i}f_e)\,d\log g_i$. Summing
over $i$ and using
$d\log(g_ig_j)=d\log(f_e(1+f_e))-2\,d\log w_e$ on each edge, the
contribution of $e$ is $df_e/(1+f_e)+2f_e\,d\log w_e$.
(ii) The solution set is isomorphic to the set of $(g,f)$ satisfying the
edge equations with $g_i\ne0$, since $b_i=(1+\sum_{e\ni i}f_e)/g_i$. Its
coordinate ring is obtained from $\C[g_i^{\pm1}]$ by adjoining the roots
of the monic quadratics $f_e^2+f_e-w_e^2g_ig_j$, so it is a free module
and embeds in its localization
$\C(g)\bigl(\sqrt{1+4w_e^2g_ig_j}:e\in E\bigr)$. The radicands are
irreducible and pairwise non-associate, so their square classes are
independent and the localization is a field of degree $2^{\abs E}$. On
the sheet $f_e=f_+(w_e^2g_ig_j)$ near $g=0$ the map
$g\mapsto(1/b_i)_i=(g_i/(1+\sum_{e\ni i}f_e))_i$ is analytic and its
Jacobian at $g=0$ is the identity, so the projection to $b$ is
dominant; the two dimensions agree.
\end{proof}

Part (i) extends the identity of Proposition~\ref{prop:floquet} to
independent vertex arguments. Part (ii) does not imply irreducibility on
a fixed line $b=z-D$; that is Theorem~\ref{thm:irreducible}.

Let $U\subseteq V$ be nonempty with $G[U]$ connected, let $E(U)$ be the
set of edges with both ends in $U$, $\partial U$ the set of edges with
exactly one end in $U$, and $N$ the set of vertices outside $U$ adjacent
to $U$. For a signing $s$ of $E(U)$ let $A^s_U$ be the corresponding
signed adjacency matrix of $G[U]$.

\begin{proposition}\label{prop:operator-valued-poles}
Let $b(t)=b^0+tq+O(t^2)$ be a holomorphic curve in $\C^V$ through a real
point $b^0$, and let $(g(t),f(t))$ be a meromorphic germ of solutions of
\eqref{free:eq:op-system} along it, with all $g_i(t)\ne0$ for small
$t\ne0$. Suppose that $g_i=a_i/t+c_i+O(t)$ with $a_i\ne0$ for $i\in U$
and that $g_i$ is finite at $t=0$ for $i\notin U$. Then:
\begin{enumerate}
\item[\textup{(a)}] \emph{Position.} There are a signing $s$ of $E(U)$
and a vector $\xi$ with $\xi_i^2=a_i$ such that
$(\diag(b^0_U)-A^s_U)\,\xi=0$. On the line $b_i=z-D_i$ the pole
therefore sits at an eigenvalue of the signed principal submatrix
$D_U+A^s_U$; the coefficients outside $U$ do not enter.
\item[\textup{(b)}] \emph{Strength.} The $f_e$, $e\in\partial U$, are
finite at $t=0$, with limits $f_e^*$, and
\begin{equation}\label{free:eq:index}
\sum_{i\in U}q_ia_i=I:=\abs U-\abs{E(U)}+\sum_{e\in\partial U}f_e^*
=\abs U-\abs{E(U)}+\sum_{x\in N}\bigl(d_-(x)-1\bigr),
\end{equation}
where $d_-(x)$ is the number of edges at $x$ with no end in $U$ along
which $f_e\to-1$. In particular $I$ is an integer. If the kernel in
\textup{(a)} is spanned by a real unit vector $u$ and
$\sum_{j\in U}q_ju_j^2\ne0$, then
\begin{equation}\label{free:eq:directional}
\Res_{t=0}g_i\,dt=\frac{I\,u_i^2}{\sum_{j\in U}q_ju_j^2}\qquad(i\in U).
\end{equation}
If the weights also move, $w_e(t)=w_e+tw_e'+O(t^2)$, the denominator is
replaced by
$\sum_{j\in U}q_ju_j^2-2\sum_{e=ij\in E(U)}s_ew_e'u_iu_j$.
\item[\textup{(c)}] \emph{Residue matrix.} Let $b_i=z-D_i$, $t=z-\lambda$,
and let the kernel be spanned by the real unit vector $u$. Let $\tree_U$
be a connected component of $\pi^{-1}(G[U])$, let $\epsilon_x=\pm1$ be
the product of the signs $s_e$ along the path in $\tree_U$ from a base
point to $x$, and $v_x=\epsilon_xu_{\pi(x)}$. Define, for
$x,y\in\tree$ joined by the path $x=x_0,x_1,\dots,x_\ell=y$ with edges
of types $e_1,\dots,e_\ell$,
\[
\mathcal G_{xy}=g_{\pi(x)}\prod_{r=1}^\ell
w_{e_r}\frac{g_{\pi(x_r)}}{1+f_{e_r}} .
\]
On the physical branch $\mathcal G_{xy}$ is the kernel of
$(z-\tH)^{-1}$. For every finite $S\subset\tree_U$,
\[
\Res_{z=\lambda}\mathcal G_{S,S}\,dz=I\,v_Sv_S^{\T}.
\]
\item[\textup{(d)}] \emph{Non-physical sheets.} In the situation of
\textup{(c)}, if $U\ne V$ or $G[U]$ contains a cycle, then the germ is
not the physical branch near $\lambda$. In particular, for $\beta\ge1$
no germ of this kind lies on the physical sheet.
\end{enumerate}
\end{proposition}

In the convention of the other sections the residue of $V_i\,dz$ is
$-Iu_i^2$; for $U=V$ we have $I=n-\abs E=-\kappa_0$, in agreement with
Theorem~\ref{thm:local-pole}. The position depends on the matrix inside
$U$ and the strength on the boundary condition outside; for fixed
$(U,s,\lambda)$ different branches outside $U$ can give different~$I$.

\begin{proof}
(a) For $e=ij\in E(U)$ the edge equation shows that $f_e$ has a simple
pole, $f_e=H_e/t+C_e+O(t)$ with $H_e^2=w_e^2a_ia_j$. Choose $\xi$ with
$\xi_i^2=a_i$ and define $s_e$ by $H_e=s_ew_e\xi_i\xi_j$. For
$e\in\partial U$ the function $f_e(1+f_e)$ has at most a simple pole,
so the meromorphic function $f_e$ is finite. The coefficient of $1/t$
in the vertex equation at $i\in U$ is
$b_i^0a_i=\sum_{e\ni i,\,e\in E(U)}H_e$; dividing by $\xi_i$ gives the
kernel equation.

(b) The next coefficients of the edge equation for $e=ij\in E(U)$ and of
the vertex equation at $i\in U$ are
\[
2C_e+1=H_e\Bigl(\frac{c_i}{a_i}+\frac{c_j}{a_j}\Bigr),\qquad
b_i^0c_i+q_ia_i=1+\sum_{e\ni i,\,e\in E(U)}C_e+
\sum_{e\ni i,\,e\in\partial U}f_e^* .
\]
Summing the second over $i\in U$ and inserting the first, the terms
with $c_i$ cancel because $\sum_{e\ni i,\,e\in E(U)}H_e=b_i^0a_i$, and
the first equality of \eqref{free:eq:index} follows. If the kernel is
spanned by $u$, then $a_i=\rho u_i^2$ with a constant $\rho$, and
\eqref{free:eq:index} gives $\rho$. When the weights move, the right
side of the first displayed identity acquires the term
$2H_ew_e'/w_e=2s_ew_e'\rho u_iu_j$. For $x\in N$ and a boundary edge
$e=ix$ the edge equation gives $g_x=O(t)$. Hence along every edge
$e'=xy$ with $y\notin U$ we have $f_{e'}(1+f_{e'})\to0$, so
$f_{e'}\to0$ or $-1$, and the limit of the vertex equation at $x$ is
$0=1+\sum_{e\in\partial U,\,e\ni x}f_e^*-d_-(x)$. Summing over $x\in N$
gives the second equality.

(c) On the physical branch $g_{\pi(y)}/(1+f_e)$ is the Green function,
for the convention $(z-\tH)^{-1}$, of the half-tree beyond a lift of
$e=\pi(x)\pi(y)$ rooted at $y$ (see \eqref{eq:m-from-Vf}), and the
product formula is the iterated Schur complement along the path. For
$e=ij\in E(U)$ the factor $w_eg_j/(1+f_e)$ tends to
$w_ea_j/H_e=s_eu_j/u_i$. The product telescopes, and with
$\rho=I$ we get $(z-\lambda)\mathcal G_{xy}\to Iv_xv_y$.

(d) Suppose the germ is the physical branch for $\im z>0$ near
$\lambda$. By the spectral theorem $iy(\lambda+iy-\tH)^{-1}$ converges
strongly to the orthogonal projection $P_\lambda$ onto
$\ker(\tH-\lambda)$ as $y\downarrow0$, so by (c) the compression of
$P_\lambda$ to $S$ is $Iv_Sv_S^{\T}$, which must be a positive
contraction. Here $I\ne0$ and all $v_x\ne0$ by \eqref{free:eq:index}.
If $I<0$ this fails for a single point. If $I>0$ and $G[U]$ contains a
cycle, then $\tree_U$ is infinite and $S$ can be chosen with
$I\norm{v_S}^2>1$. If $G[U]$ is a tree, then $\tree_U$ is a copy of it;
taking $S=\tree_U$, the nonzero eigenvalue of the compression is the
integer $I$, so $I=1$ and $\langle v,P_\lambda v\rangle=\norm v^2$ for
the extension $v$ of $v_S$ by zero, whence $\tH v=\lambda v$. If
$U\ne V$ there is a vertex $y\notin\tree_U$ adjacent to $\tree_U$,
through exactly one edge because $\tree$ is a tree, and
$(\tH v)(y)\ne0=\lambda v(y)$, a contradiction. If $\beta\ge1$ and
$U=V$, then $G[U]=G$ contains a cycle.
\end{proof}

The proposition gives necessary conditions for a germ. Existence along a
line $b=b^0+tq$ with all $q_i>0$ reduces to the scalar case: with
$g_i'=q_ig_i$, $w_e'=w_e/\sqrt{q_iq_j}$ and $D_i'=-b_i^0/q_i$, the
system \eqref{free:eq:op-system} becomes the covariance system with
spectral parameter $t$, to which Theorem~\ref{thm:partial-pole} applies
under its nondegeneracy hypotheses. The number $\sum_jq_ju_j^2$ measures
the transversality of the line to the hypersurface
$\det(\diag(b_U)-A^s_U)=0$; tangent directions and kernels of higher
dimension are not covered. Statement (d) concerns germs with a single
connected support; physical atoms supported on induced forests with
several components exist and are classified in~\cite{BGVM}, whose index
theorem also excludes a single proper connected support. The criterion
of \cite[Theorem~3.2]{BBL} for operator-valued atoms is not used.

\begin{example}[A triangle with a pendant vertex]\label{free:ex:triangle-leaf}
Let all weights be $1$ and all potentials $0$ on the triangle with one
pendant vertex $v$ attached at the vertex~$j$. Eliminating the lifts of
$v$ gives a line with period three, and the Floquet determinant yields
\[
g_j=\frac{z^2-1}{\sqrt{\Delta^2-4}},\qquad
g_v=\frac1z+\frac{g_j}{z^2},\qquad \Delta=z^3-4z+z^{-1},
\]
while the other two vertices have $g=(z^2-2)/\sqrt{\Delta^2-4}$. The
point $z=0$ is not in $\sigma(\tH)$, since the Floquet determinant of
the four-vertex matrix, $z(\Delta-2\cos\theta)$, equals $1$ at $z=0$;
by symmetry of the spectrum it lies in a gap. On the physical branch
$\sqrt{\Delta^2-4}\sim\Delta\sim1/z$, so $g_j=-z+O(z^3)$ and $g_v=O(z)$.
On the other branch $g_j=z+O(z^3)$, and $g_v$ has a simple pole at
$z=0$ with residue $2$, while the three other vertices are regular. The
two branches are exchanged by a loop around $z=1$, where $\Delta^2-4$
has a simple zero. Thus the Green function at $v$ has an antibound pole
inside a spectral gap, supported on $U=\{v\}$ at $\lambda=D_v$, with
$d_-(j)=2$ and $I=1+(d_-(j)-1)=2$ in \eqref{free:eq:index}; in agreement with (d), a
tree support with $I\ne1$ is not physical.
\end{example}

\begin{example}[Residues of both signs at one eigenvalue]
\label{free:ex:bridge}
Let the vertices be $0,1,2,3$ with zero potential, four edges of weights
$1,2,3,5$ between $0$ and $1$, an edge of weight $\tfrac14$ between $0$
and $2$, and four edges of weights
$\tfrac1{10},\tfrac2{10},\tfrac3{10},\tfrac4{10}$ between $2$ and $3$.
Take $U=\{0,1\}$ with all inner signs positive, so that $\lambda=11$ and
$u=(1,1)/\sqrt2$. If $r$ of the four outer edges are frozen at $-1$,
\eqref{free:eq:index} gives $f^*=r-1$ on the bridge and $I=r-3$. Such
germs exist for $r=2$ and $r=4$. Indeed, in the chart $x=1/g_0$,
$g_1=h/x$, $g_2=xY$, $g_3=\eta$, with $Y=16f(1+f)$ for the bridge
variable $f$, the inner edge variables are
$f_w=A_w/x$, $A_w=-\tfrac x2+\sqrt{x^2/4+w^2h}$, the outer ones
$B_d$ solve $B_d(1+B_d)=d^2xY\eta$ near $0$ or $-1$, and the vertex
equations on the line $b_i=z$ read
\[
z=x(1+f)+\sum_wA_w,\quad zh=x+\sum_wA_w,\quad
zxY=1+f+\sum_dB_d,\quad z\eta=1+\sum_dB_d .
\]
At $x=0$ the solution is $(z,h,f,\eta)=(11,1,r-1,(1-r)/11)$, and the
Jacobian with respect to $(z,h,f,\eta)$ is block triangular with
determinant $-121$. The implicit function theorem gives a unique
holomorphic solution for each freezing, with $z'(0)=(r-3)/2$. For
$r=2,4$ the function $z$ is a local coordinate, $h,Y,\eta\ne0$, and
the germ is admissible for small $x\ne0$, with residues
$(-\tfrac12,-\tfrac12)$ and $(\tfrac12,\tfrac12)$ of $(g_0,g_1)$. The
value $11$ lies above the spectrum: removing the lifts of the bridge,
the components over $\{0,1\}$ are $4$-regular trees whose operator is
the adjacency operator, of norm at most $2\sqrt3$, plus matchings of
weights $0,1,2,4$, the components over $\{2,3\}$ have norm at most~$1$,
and the bridge contributes $\tfrac14$; so
$\norm{A_{\tree}}\le7+2\sqrt3+\tfrac14<11$. For the germ with $r=4$
the residue matrix on three consecutive vertices of $\tree_U$ has all
entries $\tfrac12$ and the eigenvalue $\tfrac32>1$, which illustrates
how (d) excludes a physical atom even when the diagonal residues have
the sign of a mass. For $r=3$ one has $I=0$, $Y(0)=96$ and
$\eta(0)=-2/11$. Differentiating the outer edge equations and the
equation at $2$ gives
\[
 f'(0)=1056+\frac{192}{11}\sum_d(1+2B_d(0))d^2 .
\]
The first two vertex equations give $h=1-xf/z$ and
$z=x(f-1)+\sum_w\sqrt{w^2h+x^2/4}$. Since $z'(0)=0$, expansion to
second order gives
\[
 z=11+Cx^2+O(x^3),\qquad
 C=\frac{f'(0)}2-\frac1{22}+\frac18\sum_w\frac1w>0 ,
\]
using $\sum_d d^2=3/10$. Thus for all four freezings $g_0,g_1$ grow
like $(z-11)^{-1/2}$ at ramification points and $g_i\,dz$ is holomorphic.

Both germs with $r=2,4$ are points of the irreducible variety of
Lemma~\ref{free:lem:multi}(ii), hence are reached from the physical germ
when the $b_i$ move independently. For the fixed rational coefficients
above we have not decided whether they lie on the physical component of
the line $b_i=z$. They persist under small perturbations of the
coefficients, and for generic real coefficients near the given ones they
lie on the physical component by Theorem~\ref{thm:irreducible}.
\end{example}

%% file: sections/degenerations.tex
\section{Degenerate coefficients, decorations and an inverse problem}
\label{sec:degenerations}

The theorems on sheets and poles hold for generic coefficients, and the
exceptional set is not empty. The basic example is
Proposition~\ref{prop:uniform}: for uniform coefficients on an
edge-transitive graph the physical component has two sheets and its only
poles are the two Perron values. In this section we study three further
kinds of special coefficients, always for the physical component
$\Cphys$ and for real coefficients. For unweighted graphs the natural
candidate for the number of sheets, the count of leafless subsets of the
quotient by the degree refinement, is only an upper bound
(Proposition~\ref{prop:unweighted}). For the dipole with constant
potential and repeated weights we determine sheets, genus and poles, and
find what the poles determine (Proposition~\ref{prop:dipole-inverse}).
Attaching the same finite rooted tree to every vertex pulls the physical
curve back along a rational map and preserves the number of sheets for
all real coefficients (Theorem~\ref{thm:base-change}). Combining the last
two results we obtain pairs of models on the same simple graph with the
same finite poles and residues at every vertex but with different gaps
and densities of states (Theorem~\ref{thm:identical-poles}).

Throughout this section $G$ is connected and loopless, the coefficients
are real, and $\widehat X_G$ denotes the smooth projective model of the
physical component: the normalization of the projective closure of
$\Cphys(G;D,w)$, with the function $z\colon\widehat X_G\to\PP^1$. The
number of sheets of $\Cphys$ is $\deg z$. By
Theorem~\ref{thm:m-function-degree}(b), for arbitrary real coefficients
and every vertex $i$,
\begin{equation}\label{dgn:eq:primitive}
\C(\widehat X_G)=\C(z,G_i).
\end{equation}

\subsection{Unweighted graphs}\label{dgn:subsec:unweighted}

Let all $w_e=1$ and all $D_i=0$. The \emph{degree refinement} of $G$ is
the coarsest partition of $V$ into classes such that for any two classes
$a,b$ the number $b_{ab}$ of edges from a vertex of class $a$ to
vertices of class $b$ does not depend on the vertex. A \emph{type} is an
unordered pair $\{a,b\}$ of classes, possibly equal, with $b_{ab}>0$.
Let $L(B)$ be the number of sets $\mathcal S$ of types such that
$\sum_{b:\{a,b\}\in\mathcal S}b_{ab}\ne1$ for every class $a$; these are
the leafless edge subsets of $G$ that are unions of types, so
$L(B)\le\LG$. The equations for the physical branch on the quotient by
an equitable partition are those of
\cite[Theorem~3.8]{CMMM}, which generalizes \cite[Theorem~6.6]{ABS};
that paper does not count solutions or sheets.

\begin{proposition}\label{prop:unweighted}
Let $G$ be a connected loopless multigraph with all weights $1$ and zero
potential.
\begin{enumerate}
\item[\textup{(a)}] The number of sheets of $\Cphys$ is at most $L(B)$.
\item[\textup{(b)}] Equality can fail: there is a simple graph with $22$
vertices for which $L(B)=3$ and $\Cphys$ has two sheets.
\item[\textup{(c)}] The defect is unbounded: for every integer $c\ge2$
there is a simple graph for which $L(B)=2^{c(c+1)/2}$ and $\Cphys$ has
two sheets.
\end{enumerate}
\end{proposition}

\begin{proof}
(a) Two vertices of the same class are roots of isomorphic universal
covers: the root sees $b_{ab}$ neighbors in each class, and a
non-backtracking walk that arrives at class $a$ from class $b$ continues
to $b_{ac}$ vertices of each class $c$, with one fewer for $c=b$; so the
rooted trees agree in every finite radius. The same holds for the
half-trees obtained by removing one edge from a vertex of class $a$ to
class $b$. Hence on the physical branch $V_i$ depends only on the class
of $i$, $m_{i|e}$ only on the ordered pair of classes of the ends of
$e$, and $f_e=V_i/m_{i|e}-1$ (see \eqref{eq:m-from-Vf}), which is
symmetric in the two ends, only on the type of~$e$. These are polynomial
identities on the physical branch, so they hold on $\Cphys$. By
Lemma~\ref{lem:infinity} the points of $\widehat X_G$ over $z=\infty$ are
unramified and correspond to distinct leafless subsets $S$ through
$f_e\to-\ones_S(e)$. On $\Cphys$ such an $S$ is a union of types, and
there are $L(B)$ leafless sets of this kind. Since $z$ has a simple pole
at each point over $\infty$, its degree is at most $L(B)$.

(b) Let $P_3$ be the path with three vertices. For the convention
$(z-A)^{-1}$ its Green functions at an end and at the middle vertex are
\[
c_{\mathrm{end}}(z)=\frac{z^2-1}{z(z^2-2)},\qquad
c_{\mathrm{mid}}(z)=\frac z{z^2-2},\qquad
2c_{\mathrm{end}}(z)=\frac1z+c_{\mathrm{mid}}(z).
\]
Consider two rooted trees. In $T_A$ the root is joined to an end vertex
of each of two copies of $P_3$. In $T_B$ the root is joined to a single
vertex and to the middle vertex of one copy of $P_3$. By the displayed
identity both trees present the same self-energy to the root, so
\begin{equation}\label{dgn:eq:RA}
R_{T_A}(z)=R_{T_B}(z)=z-2c_{\mathrm{end}}(z)=z-\frac1z-\frac z{z^2-2}
\end{equation}
in the notation of \eqref{dgn:eq:RT} below. The children of the root
have degrees $\{2,2\}$ in $T_A$ and $\{1,3\}$ in $T_B$.

Let $G$ be obtained from $K_4$ by attaching $T_A$ at one vertex and
$T_B$ at each of the other three, identifying roots with vertices. It is
a simple graph with $4+6+3\cdot4=22$ vertices. The four vertices of
$K_4$ have degree $5$ and all others degree at most $3$; among the
former, the degrees of the neighbors separate the vertex carrying $T_A$
(class $a$) from the three carrying $T_B$ (class $b$), which are
permuted by automorphisms. No type inside the decorations belongs to an
admissible set $\mathcal S$: the class of a decoration vertex determines
its distance to $K_4$, each such vertex has exactly one neighbor at
smaller distance, and by descending induction on the distance, once all
types towards larger distance are excluded, the type towards the parent
would contribute exactly $1$. The remaining types are $\{a,b\}$, with
$b_{ab}=3$, $b_{ba}=1$, and $\{b,b\}$, with $b_{bb}=2$. The set
$\{\{a,b\}\}$ is excluded by class $b$, and $\emptyset$, $\{\{b,b\}\}$
and $\{\{a,b\},\{b,b\}\}$ are admissible. Hence $L(B)=3$.

The universal cover of $G$ is the $3$-regular tree with a copy of $T_A$
or $T_B$ at each vertex. Taking the Schur complement onto the vertices
of the $3$-regular tree and using \eqref{dgn:eq:RA}, the Green function
of $G$ at a vertex $i$ of $K_4$ is $G^{K_4}_i(R_{T_A}(z))$, the same
function as for the graph $K_4\triangleright T_A$ in which $T_A$ is
attached at all four vertices. By \eqref{dgn:eq:primitive} the two
physical components have the same function field over $\C(z)$, and by
Theorem~\ref{thm:base-change} the number of sheets for
$K_4\triangleright T_A$ equals that for the unweighted $K_4$, which is
$2$ (Proposition~\ref{prop:uniform}; directly, the half-tree function of
the $3$-regular tree satisfies $2m^2+zm+1=0$ and $V=-1/(z+3m)$).

(c) Divide the vertices of $K_{3c}$ into $c$ groups of three, and attach
at every vertex of the $j$th group, $0\le j<c$, $j$ copies of $T_A$ and
$c-1-j$ copies of $T_B$ with a common root. The vertices of $K_{3c}$
have degree $3c-1+2(c-1)>3$, and the number $2j$ of their neighbors of
degree $2$ separates the groups, so they form exactly $c$ classes, with
$b_{aa}=2$ and $b_{ab}=3$ for $a\ne b$. The types inside the decorations
are excluded as in (b), and every set of the $c(c+1)/2$ remaining types
is admissible. Every vertex of $K_{3c}$ receives the self-energy
$2(c-1)c_{\mathrm{end}}$, the same as when $c-1$ copies of $T_A$ are
attached everywhere, and the argument of (b) gives two sheets.
\end{proof}

In these examples, with $d$ the degree of the complete graph and $R$
the common effective energy after removing the decorations, eliminating
the half-tree function from $(d-1)m^2+Rm+1=0$ and $v=-1/(R+dm)$ gives
$(d^2-R^2)v^2+(d-2)Rv+(d-1)=0$. It is irreducible over $\C(z)$ by the
two-sheet assertion. We do not know a combinatorial formula for the
number of sheets of unweighted graphs, nor whether equality holds in
(a) for leafless graphs.

\subsection{The dipole with constant potential}\label{dgn:subsec:dipole}

Let $D_d$, $d\ge3$, be the graph with two vertices joined by $d$
parallel edges, with potential $D_1=D_2=b$ and positive weights. Let
$a_1,\dots,a_k$ be the distinct values of the weights, $\nu_j$ the
multiplicity of $a_j$, and $s$ the number of $j$ with $\nu_j=1$. Every
vertex of the universal cover sees the same weighted tree, so
$G_1=G_2=:G$. By the \emph{pole data} of a function on $\widehat X$ we
mean the multiset of pairs (value of $z$, residue with respect to $z$)
over all points of $\widehat X$ at which the function has a pole and $z$
is finite; in all cases below these points are unramified for~$z$.

\begin{proposition}\label{prop:dipole-inverse}
For the dipole $D_d$ with constant potential $b$ the following hold.
\begin{enumerate}
\item[\textup{(a)}] $\Cphys$ has $2^k-s$ sheets, and $\widehat X$ has
genus $1+(k-2)2^{k-1}$.
\item[\textup{(b)}] The function $G$ has exactly $2^k$ poles on
$\widehat X$, one for each $\varepsilon\in\{\pm1\}^k$, over
$z=b+\sum_j\varepsilon_j\nu_ja_j$. Each is unramified for $z$ and has
residue $(d-2)/2$; in particular all of them are antibound. The poles
lie only over the Perron values $b\pm\sum_j\nu_ja_j$ if and only if all
weights are equal.
\item[\textup{(c)}] The pole data of $G$ determine $b$, $d$ and the
multiset $\{\nu_ja_j\}$, and are determined by them. If all weights are
distinct, the weights can be recovered from the pole data.
\item[\textup{(d)}] The pole data do not determine the spectrum of
$\tH$. For $d=3$, $b=0$, the weights $(2,2,6)$ and $(3,3,4)$ give curves
with three sheets and genus $1$ and the same pole data: simple poles
over $-10,-2,2,10$ with residue $\tfrac12$. For $(2,2,6)$ the point $0$
lies in a gap of $\sigma(\tH)$; for $(3,3,4)$ it lies in the spectrum
and the density of the spectral measure at $0$ is $1/(10\pi)$.
\end{enumerate}
\end{proposition}

\begin{proof}
(a), (b). Put $W=1/G$ and $\zeta=z-b$, and let $h_j$ be the half-tree
function $m_{i|e}$ for an edge of weight $a_j$. By \eqref{eq:m-recursion},
$W=-\zeta-\sum_j\nu_ja_j^2h_j$ and $h_j^{-1}=W+a_j^2h_j$. With
$y_j=W+2a_j^2h_j$ this becomes
\begin{equation}\label{dgn:eq:dipole}
y_j^2=W^2+4a_j^2\quad(1\le j\le k),\qquad
2\zeta=(d-2)W-\sum_j\nu_jy_j .
\end{equation}
This square-root form of the Green function is known
\cite[Appendix~B]{ABS}, \cite{SBS}. The polynomials $W^2+4a_j^2$ have
simple roots, distinct for different $j$, so their square classes in
$\C(W)$ are independent and $\mathbb M=\C(W,y_1,\dots,y_k)$ has degree
$2^k$ over $\C(W)$, with Galois group the group of sign changes. The
$y_j$ belong to distinct characters, hence are linearly independent
over $\C(W)$, and only the identity fixes $\sum_j\nu_jy_j$. Therefore
$\mathbb M=\C(W,z)=\C(G,z)$. The physical branch has $W\sim-\zeta$ and
$y_j/W\to1$ at infinity, so by \eqref{dgn:eq:primitive} the curve with
function field $\mathbb M$ is $\widehat X$; all $V_i$, $m_{i|e}$ and
$f_e$ lie in $\mathbb M$.

Over $W=\infty$ there are $2^k$ unramified points, indexed by
$\epsilon_j=\lim y_j/W$. There
$2\zeta=(d-2-\sum_j\nu_j\epsilon_j)W+O(1/W)$, and the coefficient
$2(\sum_{\epsilon_j=-1}\nu_j-1)$ vanishes exactly at the $s$ points
where a single weight of multiplicity one has $\epsilon_j=-1$. At the
other $2^k-s$ points $z$ has a simple pole, and $z$ is finite at finite
$W$. This gives the number of sheets. The map $W$ is ramified exactly
over $W=\pm2ia_j$, at $2^{k-1}$ points of index $2$ for each value, and
the Riemann--Hurwitz formula gives
$2g-2=-2^{k+1}+k2^k$.

The poles of $G$ are the $2^k$ points over $W=0$, where
$y_j=2\varepsilon_ja_j$. There $W$ is a local coordinate,
$dy_j/dW=W/y_j=0$, and $dz/dW=(d-2)/2$ by \eqref{dgn:eq:dipole}; so $z$
is a local coordinate, $z=b-\sum_j\varepsilon_j\nu_ja_j$ at the point,
and the residue of $G=1/W$ is $(d-2)/2>0$. The poles are antibound by
Corollary~\ref{cor:antibound-sign}. If $k\ge2$, changing one sign gives
a value strictly between the two Perron values.

(c) The midpoint of the extreme poles is $b$, and $d$ is twice the
residue plus $2$. The map
$\lambda\mapsto(\lambda_{\max}-\lambda)/2$ sends the multiset of pole
positions to the multiset of all subset sums of the positive numbers
$A_j=\nu_ja_j$. A multiset of positive numbers is determined by the
multiset $\Sigma$ of its subset sums: the smallest positive element $A$
of $\Sigma$ is one of the numbers, and $\Sigma=\Sigma'\uplus(A+\Sigma')$,
where $\Sigma'$ is the multiset of subset sums of the others; $\Sigma'$
is obtained by repeatedly moving the smallest remaining element $x$ to
$\Sigma'$ and deleting one copy of $x$ and of $x+A$. Induction on $k$
finishes the argument. The converse is (b).

(d) Let $0<A<B<2A$ with $B>\sqrt2A$, and compare the weights
$(A/2,A/2,B)$ and $(B/2,B/2,A)$. Both have $k=2$, $s=1$ and
$\{\nu_ja_j\}=\{A,B\}$, so by (a)--(c) they have three sheets, genus $1$
and the same pole data. For zero potential, $0\notin\sigma(\tH)$ if and
only if the largest weight $w_*$ satisfies
$w_*^2>\sum_{e\ne*}w_e^2$; this is the criterion of
Fig\`a-Talamanca and Steger, see \cite[Theorem~B.1]{ABS}. For the first
triple $w_*=B$ and $B^2>A^2/2$; for the second $w_*=A$ and
$A^2<B^2/2$. The integer example is $A=4$, $B=6$.

To compute the density, write $a,a,c$ for the weights and
$L=W-2z=2y_1+y_2$ in \eqref{dgn:eq:dipole}. Squaring $y_2=L-2y_1$
gives $Ly_1=W^2-zW+z^2+4a^2-c^2$. Squaring again and substituting
$y_1^2=W^2+4a^2$ gives, with $v=1/W$,
\[
 (z^2-(2a+c)^2)(z^2-(2a-c)^2)v^3
 +2z(4a^2+c^2-z^2)v^2+(4a^2-2c^2-z^2)v+2z=0 .
\]
For the two integer triples this reads $J-56v=0$ and $J+4v=0$,
respectively, where
\[
J(z,v)=(z^4-104z^2+400)v^3-(2z^3-104z)v^2-z^2v+2z.
\]
By symmetry of the spectral measure, $G(i\eta)\in i(0,\infty)$ for $\eta>0$. At $z=0$ both
cubics have leading coefficient $400$ and three simple roots, so the
three branches are holomorphic near $0$. For $(3,3,4)$ the roots are
$0,\pm i/10$; the branch through $0$ is $-z/2+O(z^2)$, which has
negative imaginary part on $i(0,\varepsilon)$, so the physical branch
tends to $i/10$ and the density is
$\pi^{-1}\im G(0+i0)=1/(10\pi)$, continuous near~$0$. For $(2,2,6)$ the
only root on the imaginary axis is $0$, and the physical branch is
$z/28+O(z^2)$, real on the real axis near $0$.
\end{proof}

\begin{remark}\label{dgn:rem:K33}
The same two models live on the simple graph $K_{3,3}$: label both
sides by $\Z/3$ and give the edge from $i$ on the left to $i+j$ on the
right the weight $w_j$. Each $j$ is a perfect matching, and mapping the
three matchings to the three edges of $D_3$ identifies the weighted
universal covers, hence the Green functions and their continuations.
For a second-order invariant that does distinguish the two models note
$\langle\delta,\tH^2\delta\rangle=\sum_ew_e^2=44$ and $34$.
\end{remark}

\begin{remark}[Branch points that are not of square-root type]
\label{dgn:rem:cuberoot}
For $d=3$, $b=0$ and weights $(1,1,\sqrt2)$ the curve has three sheets.
At the point over $W=\infty$ where only the weight $\sqrt2$ has
$\epsilon_j=-1$, equation \eqref{dgn:eq:dipole} gives, with $u=G=1/W$,
\[
z=\frac{1-2\sqrt{1+4u^2}+\sqrt{1+8u^2}}{2u}=-2u^3+O(u^5).
\]
This single point is the whole fiber of the map $z$ of degree $3$ over
$z=0$. So $\widehat X$ is smooth there, $z$ has ramification index~$3$,
and the Green function has a cube-root branch point at $z=0$ on the
physical sheet. This answers in the negative the part of
\cite[Problem~10.1]{ABS} asking whether all branch points are of
square-root type. We do not claim the example: the cube root on this
critical surface $w_*^2=\sum_{e\ne*}w_e^2$ is computed by Singh, Bera
and Shenoy~\cite{SBS}, and the general theory of cubic cusps of the
density is due to Huang, Su and Xu~\cite{HuangSuXu}. What
\eqref{dgn:eq:dipole} adds is that the singularity is a ramification
point of a smooth curve.
\end{remark}

Part (b) answers, within this family, the question raised after
Conjecture~3.4 of \cite{CSZ} whether there are families whose only poles
on higher sheets are the Perron values: here this happens exactly for
uniform weights. Compared with Corollary~\ref{free:cor:dipole}, repeated
weights lower the number of sheets from $2^k-k$ to $2^k-s$ with $k$ now
the number of distinct weights.

\subsection{Attaching a tree at every vertex}\label{dgn:subsec:basechange}

Let $T$ be a finite rooted tree with real potential, equal to $0$ at
the root, and real nonzero weights. The graph $G\triangleright T$ is
obtained from $G$ by attaching a copy of $T$ at every vertex,
identifying its root with the vertex; the vertices of $G$ keep their
potential. Let $B_T$ be the Jacobi matrix of $T$ with the root removed
and $v$ the vector of weights of the edges at the root. Put
\begin{equation}\label{dgn:eq:RT}
R_T(z)=z+v^{\T}(B_T-z)^{-1}v=z-\sum_{j=1}^r\frac{c_j}{z-b_j},
\qquad b_j\in\R\ \text{distinct},\ c_j>0,
\end{equation}
where the $b_j$ are the eigenvalues of $B_T$ on whose eigenspaces $v$
has a nonzero projection, of squared norm $c_j$. Every function of the
form \eqref{dgn:eq:RT} arises, from the star with leaves of potential
$b_j$ and weights $\sqrt{c_j}$. Such a function $R$ has degree $r+1$,
satisfies $\im R(z)=\im z\,(1+\sum_jc_j\abs{z-b_j}^{-2})$, so that
$R^{-1}(\mathbb H)=\mathbb H$ for the upper half-plane $\mathbb H$, has
$R'>0$ on the real line, and all its preimages of $\infty$ and of real
values are real and simple. The formula for the Green function of a
decorated graph and the resulting description of the spectrum are due to
Schenker and Aizenman~\cite[Proposition~3.1]{SchenkerAizenman}; the
composition rule $R_{T\triangleright T'}=R_T\circ R_{T'}$ is the formula
for the reciprocal Cauchy transform of a comb product, that is, monotone
convolution \cite[Proposition~2.5]{AHL}. The new point below is the
statement about all sheets.

\begin{lemma}\label{dgn:lem:fibre}
Let $p\colon X\to\PP^1$ and $q\colon Y\to\PP^1$ be nonconstant maps of
smooth connected projective curves. Suppose that $p$ has a holomorphic
section $\sigma\colon\mathbb H\to X$ and that $q^{-1}(\mathbb H)$ is
connected. Then the normalization $Z$ of the fiber product
$X\times_{\PP^1}Y$ is connected, and $\deg(Z\to Y)=\deg p$,
$\deg(Z\to X)=\deg q$.
\end{lemma}

\begin{proof}
Let $\mathcal B$ be the union of the branch values of $p$ and $q$,
$s_0\in\mathbb H\setminus\mathcal B$ a base point, and
$a=\sigma(s_0)$. Loops in $\mathbb H\setminus\mathcal B$ fix $a$,
because they lift to $\sigma$. They act transitively on $q^{-1}(s_0)$,
because $q^{-1}(\mathbb H\setminus\mathcal B)$ is a connected surface
minus finitely many points and covers
$\mathbb H\setminus\mathcal B$. Loops in $\PP^1\setminus\mathcal B$ act
transitively on $p^{-1}(s_0)$. Hence any pair in
$p^{-1}(s_0)\times q^{-1}(s_0)$ can be moved first to a pair with first
entry $a$ and then, keeping $a$, to any second entry. So the unramified
part of the fiber product is connected. Near the remaining points the
fiber product is given by $x^e=y^h$, all of whose local branches meet
the unramified part, so $Z$ is connected. The degrees are read off a
general fiber.
\end{proof}

The connectedness hypothesis matters: for $p(x)=x^2$ and $q(y)=y^2$ the
map $p$ has a section over $\mathbb H$, but the fiber product splits
into $x=\pm y$.

The physical branch gives a section of $z\colon\widehat X_G\to\PP^1$
over $\mathbb H$: by Proposition~\ref{prop:physical-branch} it is a
holomorphic map from $\mathbb H$ to $\Cphys$ on which $z$ is the
identity, and such a map lifts to the normalization. For $R$ as in
\eqref{dgn:eq:RT} we write $R^*\widehat X_G$ for the normalization of
$\widehat X_G\times_{\PP^1_s}\PP^1_z$, formed with $s\colon\widehat
X_G\to\PP^1$ and $R\colon\PP^1_z\to\PP^1_s$. It carries the section
$z\mapsto(\sigma(R(z)),z)$.

\begin{theorem}\label{thm:base-change}
Let $G$ be a connected loopless multigraph with real potential and real
nonzero weights, and let $T$ be a finite rooted tree as above. Then:
\begin{enumerate}
\item[\textup{(a)}] The curve $R_T^*\widehat X_G$ is connected, and
there is a unique isomorphism over the $z$-line
\[
\widehat X_{G\triangleright T}\;\cong\;R_T^*\widehat X_G
\]
respecting the physical sections. On it, the Green function of
$G\triangleright T$ at a vertex $i$ of $G$ is $G_i\circ R_T$.
\item[\textup{(b)}] The physical components of $G\triangleright T$ and
of $G$ have the same number of sheets, for all real coefficients.
\item[\textup{(c)}] Let $N$ and $g$ be the number of sheets and the
genus of $\widehat X_G$ and $m=\deg R_T$. The genus $g'$ of
$\widehat X_{G\triangleright T}$ satisfies
\[
2g'-2=-2N+\sum_{s\in\PP^1}\ \sum_{P\mapsto s}\ \sum_{z_0\mapsto s}
\bigl(e_P-\gcd(e_P,h_{z_0})\bigr),
\]
where $e_P$ and $h_{z_0}$ are the ramification indices of
$P\in\widehat X_G$ and of $z_0\in\PP^1_z$ over $s$. If the branch values
of $\widehat X_G$ and of $R_T$ are disjoint, then
$g'=1-N+m(g+N-1)$. For given $R_0$ of the form \eqref{dgn:eq:RT}, the
maps $R_\delta(z)=R_0(z-\delta)+\delta$ have this property for all but
finitely many real~$\delta$.
\end{enumerate}
\end{theorem}

\begin{proof}
(a) Since $R_T^{-1}(\mathbb H)=\mathbb H$ is connected,
Lemma~\ref{dgn:lem:fibre} shows that $Z=R_T^*\widehat X_G$ is connected,
of degree $N$ over the $z$-line. Its function field is
$\C(\widehat X_G)\otimes_{\C(s)}\C(z)$, where $s$ is identified with
$R_T(z)$; by \eqref{dgn:eq:primitive} it is generated over $\C(z)$ by
the pull-back of $G_i$, whose germ along the section is
$G_i(R_T(z))$. On the other side, the universal cover of
$G\triangleright T$ is the universal cover of $G$ with a copy of $T$
attached at every vertex. For $z\in\mathbb H$, eliminating the vertices
of the decorations leaves on the vertices of $\tree$ the operator
$\tH-R_T(z)$, with $R_T(z)\in\mathbb H$. Hence the Green function of
$G\triangleright T$ at the vertex $i$ of $G$ is $G_i(R_T(z))$, and by
\eqref{dgn:eq:primitive} for $G\triangleright T$ it generates
$\C(\widehat X_{G\triangleright T})$ over $\C(z)$. Evaluation along the
physical sections embeds both function fields into the field of germs
of meromorphic functions at a point of $\mathbb H$, with the same image
$\C(z,G_i(R_T(z)))$. This gives the isomorphism of the smooth projective
curves, and it is unique because the germ generates. The continuation to
the other sheets is an identity of function fields; no components are
added by clearing denominators at the eigenvalues of $B_T$.

(b) follows from (a) and the lemma.

(c) Over $s$, near $(P,z_0)$, the fiber product is $x^{e_P}=y^{h_{z_0}}$
with $y$ a local coordinate of the $z$-line. Its normalization has
$\gcd(e_P,h_{z_0})$ points, at each of which $y$ vanishes to order
$e_P/\gcd(e_P,h_{z_0})$. The Riemann--Hurwitz formula for the map $z$ of
degree $N$ gives the first formula. If the branch values are disjoint,
each ramification point of $\widehat X_G$ contributes $m$ times and
nothing else contributes, so $2g'-2=-2N+m(2g-2+2N)$. The finite critical
values of $R_\delta$ are those of $R_0$ shifted by $\delta$, and
$\infty$ is not a branch value of $R_\delta$.
\end{proof}

By the universal property of the normalized fiber product,
$(R_1\circ R_2)^*\widehat X\cong R_2^*(R_1^*\widehat X)$, which
corresponds to attaching trees twice. For generic coefficients on $G$
statement (b) is consistent with Theorem~\ref{thm:generic-degree},
because a leafless subset of $G\triangleright T$ contains no edge of a
decoration and $L(G\triangleright T)=\LG$; the content of (b) is that
nothing is lost for special coefficients either.

\subsection{Identical poles, different gaps}\label{dgn:subsec:identical}

Two models on the same graph have \emph{the same pole data} if for
every vertex the pole data of the Green function, in the sense of
Subsection~\ref{dgn:subsec:dipole}, agree.

\begin{theorem}\label{thm:identical-poles}
Let $m\ge1$, let $\xi_1<\dots<\xi_m$ be real numbers and
$d_1,\dots,d_m>0$. There are a finite simple connected graph $\Gamma_m$
with $6m$ vertices and cycle rank~$4$, and two sets of real
coefficients on it, such that:
\begin{enumerate}
\item[\textup{(a)}] the two models have the same pole data at every
vertex, and the Stieltjes transforms of their densities of states also
have the same pole data; none of these poles lies on the physical
sheet;
\item[\textup{(b)}] both physical components have three sheets;
\item[\textup{(c)}] the points $\xi_1,\dots,\xi_m$ lie in $m$ different
bounded gaps of the spectrum of the first model, while the second model
has a continuous density of states near each $\xi_\ell$, equal to
$d_\ell$ at $\xi_\ell$;
\item[\textup{(d)}] if $(d_\ell)$ is replaced by $(td_\ell)$, then for
all but finitely many $t>0$ both curves have genus $3m-2$.
\end{enumerate}
\end{theorem}

\begin{proof}
\emph{Construction.} Put $S=\sum_\ell d_\ell^{-1}$,
$p_\ell=d_\ell^{-1}/S$, $u=\sum_\ell p_\ell\xi_\ell$ and
$a=S/(10m\pi)$. On $K_{3,3}$ with the three matchings of
Remark~\ref{dgn:rem:K33} take the potential $u$ and the weights
$a\cdot(2,2,6)$ for the first model and $a\cdot(3,3,4)$ for the second.
Both have the same Green function $G^{(\nu)}(s)$ at all vertices,
$\nu=1,2$, namely $a^{-1}$ times the Green function of
Proposition~\ref{prop:dipole-inverse}(d) at $(s-u)/a$. Hence: both
curves have three sheets and genus $1$; $G^{(\nu)}$ has exactly four
poles, over $u+a\lambda$, $\lambda\in\{\pm2,\pm10\}$, unramified with
residue $\tfrac12$; over $s=\infty$ there are three unramified points
with $G^{(\nu)}=\gamma/s+O(s^{-2})$, $\gamma\in\{-1,1,2\}$, the physical
one being $\gamma=-1$ (insert $v=\gamma/s$ into the cubics of
Proposition~\ref{prop:dipole-inverse}: the leading term is
$s(\gamma+1)(\gamma-1)(\gamma-2)$); the point $u$ lies in a bounded gap
of the first model, the spectrum being symmetric about $u$; and the
second model has continuous density $1/(10\pi a)$ at $u$.

Let $\mathcal L(z)=\sum_\ell p_\ell/(z-\xi_\ell)$ and
$R(z)=u+1/\mathcal L(z)$. As in Lemma~\ref{free:lem:F}, $\mathcal L$ has
exactly $m-1$ zeros $b_1<\dots<b_{m-1}$, simple and interlacing the
$\xi_\ell$, and $1/\mathcal L(z)=z-u+O(1/z)$ has negative residues
at the $b_j$. So $R$ is of the form \eqref{dgn:eq:RT} with $r=m-1$,
and $R(\xi_\ell)=u$, $R'(\xi_\ell)=1/p_\ell$. Let $\Gamma_m$ be
$K_{3,3}$ with $m-1$ pendant vertices at every vertex, with potentials
$b_j$ and weights $\sqrt{c_j}$, the same in both models. It is simple
and connected with $6m$ vertices and $9+6(m-1)$ edges.

\emph{Sheets.} Theorem~\ref{thm:base-change} gives (b) and identifies
both curves with $R^*\widehat X^{(\nu)}$. On this curve the Green
functions at a vertex $i$ of $K_{3,3}$ and at its pendant vertex of
kind $j$ are, by the block inverse,
\begin{equation}\label{dgn:eq:leaf}
V_i(z)=G^{(\nu)}(R(z)),\qquad
V_{ij}(z)=-\frac1{z-b_j}+\frac{c_j}{(z-b_j)^2}\,G^{(\nu)}(R(z)),
\end{equation}
and the Stieltjes transform of the density of states is
\begin{equation}\label{dgn:eq:dos}
M(z)=\frac1m\Bigl(R'(z)\,G^{(\nu)}(R(z))-\sum_j\frac1{z-b_j}\Bigr).
\end{equation}

\emph{Poles.} Each value $u+a\lambda$ has exactly $m$ preimages $x$
under $R$, all real, simple and different from the $b_j$. Over each
$x$ there is one point at which $G^{(\nu)}\circ R$ has a pole; it is
unramified, and the residues of $V_i$, $V_{ij}$ and $M$ are
$1/(2R'(x))$, $c_j/(2(x-b_j)^2R'(x))$ and $1/(2m)$. They are positive,
so these poles are antibound
(Corollary~\ref{cor:antibound-sign}). Over $z=b_j$ there are three
points, unramified, corresponding to the three points over
$s=\infty$. There $R(z)=-c_j/(z-b_j)+O(1)$ and
$G^{(\nu)}(R(z))=\gamma/R(z)+O((z-b_j)^2)$, so $V_i$ and the $V_{ij'}$
with $j'\ne j$ are regular, while
\[
V_{ij}(z)=\frac{-1-\gamma}{z-b_j}+O(1),\qquad
M(z)=\frac{-1-\gamma}{m(z-b_j)}+O(1).
\]
This is regular at the physical point $\gamma=-1$ and has residue $-2$
or $-3$ (respectively $-2/m$, $-3/m$) at the other two. There are no
further poles at finite $z$, because $G^{(\nu)}$ has no further poles
and \eqref{dgn:eq:leaf}, \eqref{dgn:eq:dos} have no other denominators.
All the data listed are the same for $\nu=1,2$. This proves (a); note
that the poles over $b_j$ have negative residue and are nevertheless
not on the physical sheet.

\emph{Gaps and densities.} For real $x\ne b_j$ the Schur complement
shows that $x$ is in the spectrum of the decorated operator if and only
if $R(x)\in\sigma(\tH)$ \cite[Proposition~3.1]{SchenkerAizenman}. The
points $b_j$ are not in the spectrum: in the equation
$(\tH_{\Gamma_m}-b_j)\psi=\phi$ the rows of the pendant vertices of kind
$j$ determine $\psi$ on the lifts of $K_{3,3}$, the rows of the other
pendant vertices then determine $\psi$ there because
$b_{j'}\ne b_j$, and the rows of the lifts of $K_{3,3}$ determine
$\psi$ on the pendant vertices of kind $j$; each step is bounded. The
function $R$ maps each of the $m$ intervals of $\R\setminus\{b_j\}$
monotonically onto $\R$. The preimage of the gap of the first model
around $u$ in each interval is a bounded open interval free of spectrum
whose end points belong to the spectrum; these are $m$ different
bounded gaps, and by interlacing each contains exactly one $\xi_\ell$.
For the second model, \eqref{dgn:eq:dos} on the physical branch gives
near $\xi_\ell$ a continuous density, with value
\[
\frac1\pi\im M(\xi_\ell+i0)=\frac{R'(\xi_\ell)}{m}\cdot\frac1{10\pi a}
=\frac1{10\pi a\,m\,p_\ell}=d_\ell .
\]

\emph{Genus.} Replacing $d_\ell$ by $td_\ell$ leaves $p_\ell$, $u$ and
$R$ unchanged and replaces $a$ by $a/t$, so the finite branch values of
$\widehat X^{(\nu)}$ become $u+(a/t)\beta$ with $\beta$ in a fixed
finite set. Neither $u$ nor $\infty$ is a branch value of $R$. Hence
each pair of a branch value $\beta\ne0$ and a branch value of $R$
excludes at most one $t$, and for the remaining $t$
Theorem~\ref{thm:base-change}(c) with $N=3$, $g=1$ gives
$g'=3m-2$.
\end{proof}

\begin{example}\label{dgn:ex:leaf}
The smallest decoration with $m=2$ attaches one pendant vertex with
potential $0$ and weight $1$ at every vertex of $K_{3,3}$, with the
weights $(2,2,6)$ and $(3,3,4)$ and zero potential. Then $R(z)=z-1/z$,
$R(\pm1)=0$ and $R'(\pm1)=2$: the points $-1$ and $1$ lie in two
different gaps of the first model, and the second model has density of
states $1/(10\pi)$ at both. The vertices of $K_{3,3}$ carry the eight
poles $\tfrac12(\lambda\pm\sqrt{\lambda^2+4})$,
$\lambda\in\{\pm2,\pm10\}$, and the pendant vertices carry in addition
the two poles over $z=0$ with residues $-2$ and $-3$. The branch values
$\pm2i$ of $R$ are regular values of both cubics: at either value,
substitute $v=zt$ in $J-56v$ and $J+4v$ and divide by $-2z$.
The cubics in $t$ are $1664t^3+224t^2+26t-1$ and
$1664t^3+224t^2-4t-1$. Both reduce modulo $3$ to the squarefree
cubic $2(t+1)(t^2+1)$. Thus the degree-two base change gives
genus~$4$ by Theorem~\ref{thm:base-change}(c).
\end{example}

The case $m=1$ is the pair of Proposition~\ref{prop:dipole-inverse}(d)
on $K_{3,3}$, which is leafless and therefore lies within the framework
of~\cite{CSZ}; for $m>1$ the graphs $\Gamma_m$ have pendant vertices.
The theorem compares diagonal Green functions only. Whether models with
the same poles, residues and branch points, or with the same
off-diagonal data, can still differ in their gaps is not addressed. The
result belongs to the circle of questions around
\cite[Example~9.2, Fact~9.3 and Problem~10.7]{ABS}, where isospectral
periodic Jacobi matrices with different densities of states are
exhibited: here the invariant that fails to determine the density of
states is the set of all poles and residues, on all sheets, of all
diagonal Green functions. The transformation rules for the Green
function under decoration are those of \cite{SchenkerAizenman,AHL}; the
new ingredient is Theorem~\ref{thm:base-change}, which controls all
sheets at once.

%% file: sections/abs.tex
\section{Remarks on problems of Avni, Breuer and Simon}
\label{sec:abs}

Section~10 of \cite{ABS} lists open problems on periodic Jacobi matrices on
trees. Problem~10.2 (is the Green function two-sheeted?) is answered by
Theorems~\ref{thm:green-degree} and~\ref{thm:m-function-degree}, and
Problem~10.1 (are all branch points of square-root type?) by the cube-root
branch point of Proposition~\ref{prop:dipole-inverse}. In this section we
comment on two further items of that list. Problem~10.3 concerns the poles
of the half-tree $m$-functions on non-physical sheets
(Section~\ref{abs:subsec-poles}); Conjectures~10.4 and~10.5, together with
the Borg-type Conjectures~B.4 and~B.5 of the appendix of the arXiv version
of~\cite{ABS}, concern the gaps of $\sigma(\tH)$
(Section~\ref{abs:subsec-gaps}). All references to \cite{ABS} are to the
third arXiv version, arXiv:1911.02612v3; see
Remark~\ref{abs:rem-literature}.

\subsection{Poles of half-tree \texorpdfstring{$m$}{m}-functions}
\label{abs:subsec-poles}

When $G$ is a cycle, $\tree$ is a line with periodic coefficients, the
half-trees are half-lines, and the $m$-function of each half-line is
two-sheeted. It is classical, and recalled in \cite[\S10.1]{ABS}, that the
poles of the $m$-function of one half-line on its second sheet are located
at the poles, on the physical sheet, of the $m$-function of the
complementary half-line. Problem~10.3 of \cite{ABS} asks, in paraphrase,
what relation there is on a general periodic tree between the poles of the
$m$-functions on non-physical sheets and the physical poles of the other
$m$-functions that enter the recursion \eqref{eq:m-recursion}.

Our answer has three parts. The reflection that produces the
correspondence on the line does not preserve the recursion at a vertex of
degree different from two (Proposition~\ref{prop:reflection-residual}).
There is a local criterion, valid on all sheets, which reads off the order
of every $m_{i|e}$ from the orders of the vertex functions
(Proposition~\ref{prop:half-tree-identity}). Finally, explicit families
show that non-physical poles of one $m$-function may or may not sit at the
energy of a physical pole of another one
(Theorem~\ref{abs:thm-parallel}, Proposition~\ref{abs:prop-semiregular}),
and that a non-physical pole can lie at an energy at which no half-tree has
any spectrum at all (Theorem~\ref{abs:thm-bridge}). We do not obtain a
global description of which poles of the various $m$-functions lie over
the same energy; this remains open.

Throughout this subsection $G$ is a connected loopless multigraph with at
least one edge. Let $C$ be the smooth projective model of an irreducible
component of the spectral curve~$\Xc$. The coordinates $z$, $V_i$, $f_e$
are meromorphic functions on $C$; no $V_i$ vanishes identically and no
$f_e$ is identically $0$ or $-1$, because admissible solutions are dense
in $C$ and satisfy $f_e(1+f_e)=w_e^2V_iV_j\ne0$. Hence
$m_{i|e}=V_i/(1+f_e)=f_e/(w_e^2V_j)$ of \eqref{eq:m-from-Vf} is a
meromorphic function on $C$, not identically zero. On the physical
component it is the analytic continuation of the half-tree $m$-function,
by \eqref{eq:G-from-m} and Proposition~\ref{prop:physical-branch}. For a
point $Q\in C$ we write $\ord_Q$ for the order at $Q$, positive at zeros
and negative at poles. A \emph{residue pole} of $m_{i|e}$ is a point $Q$
with $z(Q)$ finite at which the differential $m_{i|e}\,dz$ has a nonzero
residue.

\begin{proposition}[Local criterion]\label{prop:half-tree-identity}
Let $e=ij$ be an edge, $Q\in C$, and $p_i=\ord_QV_i$, $p_j=\ord_QV_j$. Then
\begin{equation}\label{abs:eq-order}
\ord_Q m_{i|e}=
\begin{cases}
(p_i-p_j)/2, & p_i+p_j<0,\\
p_i, & p_i+p_j=0,\\
p_i, & p_i+p_j>0\ \text{and}\ f_e(Q)=0,\\
-p_j, & p_i+p_j>0\ \text{and}\ f_e(Q)=-1,
\end{cases}
\end{equation}
and these cases are exhaustive. Suppose moreover that $z(Q)$ is finite and
that all $V_k$, $k\in V$, are finite at $Q$. Then for every vertex $j$ with
$V_j(Q)=0$ there is exactly one edge $e\ni j$ with $f_e(Q)=-1$; writing
$e=ij$, the function $m_{i|e}$ has a pole of order $p_j$ at $Q$, and
$m_{i'|e'}$ is finite at $Q$ for every other edge $e'=i'j$ at~$j$. If in
addition $z$ is a local coordinate at $Q$ and $V_j$ has a simple zero
there, then
\begin{equation}\label{abs:eq-residue}
\Res_Q m_{i|e}\,dz=-\frac{1}{w_e^2\,(dV_j/dz)(Q)}\ne0 .
\end{equation}
\end{proposition}

\begin{proof}
The edge equation $f_e(1+f_e)=w_e^2V_iV_j$ holds on $C$. If $p_i+p_j<0$,
then $f_e$ has a pole at $Q$, so $\ord_Qf_e=\ord_Q(1+f_e)=(p_i+p_j)/2$ and
$\ord_Qm_{i|e}=p_i-(p_i+p_j)/2$. If $p_i+p_j=0$, then $f_e$ and $1+f_e$ are
finite and nonzero at $Q$ and $\ord_Qm_{i|e}=p_i$. If $p_i+p_j>0$, then
$f_e$ is finite at $Q$ with $f_e(Q)\in\{0,-1\}$. When $f_e(Q)=0$ we use
$m_{i|e}=V_i/(1+f_e)$, and when $f_e(Q)=-1$ we use $m_{i|e}=f_e/(w_e^2V_j)$.

Now let $z(Q)$ and all $V_k(Q)$ be finite and $V_j(Q)=0$. For every edge
$e'=i'j$ the product $V_{i'}V_j$ vanishes at $Q$, so
$f_{e'}(Q)\in\{0,-1\}$, and the vertex equation at $j$ gives
$\sum_{e'\ni j}f_{e'}(Q)=-1$. Hence exactly one edge $e=ij$ at $j$ has
$f_e(Q)=-1$. For it $p_i\ge0$ and $p_j>0$, and the last case of
\eqref{abs:eq-order} gives $\ord_Qm_{i|e}=-p_j$; for the other edges the
third case gives $\ord_Qm_{i'|e'}=p_{i'}\ge0$. The residue formula is the
leading term of $f_e/(w_e^2V_j)$.
\end{proof}

Two consequences illustrate the criterion. At a point where all $V_i$ have
simple poles, as at the Perron poles of Theorem~\ref{thm:local-pole}, the
first case gives $\ord_Qm_{i|e}=0$: every half-tree function is finite and
nonzero there. In the opposite direction, \eqref{abs:eq-residue} produces
residue poles of half-tree functions at points where no vertex Green
function diverges. On the physical branch the identity
$1+f_e=(1-w_e^2m_{i|e}m_{j|e})^{-1}$ and the resulting relations between
zeros of $G_j$ and poles of the $m$-functions are those of
\cite[\S6]{ABS} and \cite[\S2]{BBGSS}; the point of
Proposition~\ref{prop:half-tree-identity} is that they hold on every
sheet.

\begin{proposition}[Reflection]\label{prop:reflection-residual}
Let $(z,V,f)$ be an admissible solution of \eqref{eq:system} and
$m_{i|e}$ as in \eqref{eq:m-from-Vf}. Then all $m_{i|e}$ are finite and
nonzero, and the reflected quantities
\[
m^{\#}_{i|e}:=\frac1{w_e^2\,m_{j|e}}\qquad(e=ij)
\]
satisfy, for every vertex $i$ and every edge $e\ni i$,
\begin{equation}\label{abs:eq-residual}
\frac1{m^{\#}_{i|e}}-\Bigl(D_i-z-\sum_{\substack{e'\ni i,\ e'\ne e\\ e'=ik}}
w_{e'}^2\,m^{\#}_{k|e'}\Bigr)=\frac{d_i-2}{V_i}.
\end{equation}
In particular $m^{\#}$ satisfies the half-tree recursion in
\eqref{eq:m-recursion} at every admissible solution if $G$ is a cycle, and
at no admissible solution otherwise. When $G$ is a cycle, $m^{\#}$ is the
family of half-tree functions of the
admissible solution $(z,-V,-1-f)$.
\end{proposition}

\begin{proof}
An admissible solution has $f_e\notin\{0,-1\}$ by the edge equation, so
the $m_{i|e}$ are finite and nonzero. By \eqref{eq:m-from-Vf},
$1/m^{\#}_{i|e}=w_e^2m_{j|e}$ and $w_{e'}^2m^{\#}_{k|e'}=1/m_{i|e'}
=(1+f_{e'})/V_i=1/V_i+w_{e'}^2m_{k|e'}$ for $e'=ik$. The vertex equation
reads $D_i-z=1/V_i+\sum_{e'\ni i}w_{e'}^2m_{k|e'}$. Therefore the left-hand
side of \eqref{abs:eq-residual} equals
\[
w_e^2m_{j|e}-\frac1{V_i}-\sum_{e'\ni i}w_{e'}^2m_{k|e'}
+\frac{d_i-1}{V_i}+\sum_{e'\ni i,\,e'\ne e}w_{e'}^2m_{k|e'}=\frac{d_i-2}{V_i}.
\]
Since $1/V_i\ne0$, the residual vanishes at every vertex exactly when all
degrees equal two, that is, when the connected graph $G$ is a cycle
(including two vertices joined by two parallel edges). If $d_i=2$ for all
$i$, the map $(V,f)\mapsto(-V,-1-f)$ preserves the edge equations and
transforms the vertex equation into $(D_i-z)V_i=d_i-1+\sum_{e\ni i}f_e$,
which is the vertex equation again. The half-tree functions of the image
are $-V_i/(-f_e)=V_i/f_e=(1+f_e)/(w_e^2V_j)=m^{\#}_{i|e}$.
\end{proof}

For a cycle the spectral curve has $\LG=2$ sheets, and the involution
$(V,f)\mapsto(-V,-1-f)$ exchanges them: near $z=\infty$ the physical
functions satisfy $m_{i|e}\sim-1/z$, whereas $m^{\#}_{i|e}$ grows
linearly. If $e=ik$ and $e'=ij$ are the two edges at a vertex $i$ of the
cycle, the recursion gives
\[
m^{\#}_{j|e'}=\frac1{w_{e'}^2\,m_{i|e'}}
=\frac{D_i-z-w_e^2\,m_{k|e}}{w_{e'}^2}.
\]
Thus the non-physical branch of the $m$-function of the half-line that
starts at $j$ and avoids $i$ is an affine function of the physical
$m$-function of the half-line that starts at $k$ and avoids $i$; their
finite poles are at the same energies, and the residues are in the ratio
$-w_e^2/w_{e'}^2$. This is the correspondence on the line.
Equation~\eqref{abs:eq-residual} shows that at a branching vertex the same
reflection leaves the residual $(d_i-2)/V_i$, so that the reflected family
is not a solution: the correspondence in this form is special to
degree two. The following results show that its conclusion also fails.

\begin{theorem}[Parallel edges]\label{abs:thm-parallel}
Let $G$ have two vertices $1,2$ joined by $k\ge2$ parallel edges, with
real potentials $D_1=a\ne b=D_2$ and pairwise distinct weights $w_e>0$.
Put $c_e=\sum_{e'\ne e}w_{e'}^2-w_e^2$, and assume that $c_e\ne0$ for all
$e$ and $\sum_e\epsilon_ew_e\ne0$ for all $\epsilon\in\{\pm1\}^E$. Then:
\begin{enumerate}
\item $\Xc$ is irreducible; hence every sheet is an analytic continuation
of the physical branch.
\item For each edge $e$ the function $m_{1|e}$ has exactly one residue
pole $Q_e$ on $\Xc$. It lies over $z=a$, the function $z$ is a local
coordinate at $Q_e$, and
\[
\Res_{Q_e}m_{1|e}\,dz=\frac{c_e}{w_e^2}.
\]
At $Q_e$ all $m_{1|e'}$ with $e'\ne e$ and all $m_{2|e'}$ are finite. At
every other point with finite $z$ at which $m_{1|e}$ has a pole, $z=a$,
the map $z$ is ramified, and $m_{1|e}\,dz$ is holomorphic.
\item The point $Q_e$ lies on the physical sheet, that is, in the closure
of the physical branch over the upper half-plane, if and only if
$w_e^2>\sum_{e'\ne e}w_{e'}^2$. In that case the spectral measure of the
root of the half-tree has an atom of mass $-c_e/w_e^2$ at~$a$.
\end{enumerate}
The same statements hold for $m_{2|e}$ with $a$ and $b$ exchanged.
\end{theorem}

\begin{proof}
\emph{The curve.} For an admissible solution put $\Pi=V_1V_2$,
$r_e=1+2f_e$ and $g=1+\sum_ef_e$. The system \eqref{eq:system} gives
\begin{equation}\label{abs:eq-parallel}
r_e^2=1+4w_e^2\,\Pi,\qquad g=(a-z)V_1=(b-z)V_2,\qquad
(a-z)(b-z)=\frac{g^2}{\Pi}.
\end{equation}
Let $Y$ be the smooth projective model of the curve
$\{r_e^2=1+4w_e^2\Pi\ (e\in E)\}$. It is a cover of the $\Pi$-line of
degree $2^k$. The branch values $-1/(4w_e^2)$ are pairwise distinct, and a
small loop around one of them changes the sign of the corresponding $r_e$
only; so the monodromy is transitive and $Y$ is irreducible. On $Y$ the
function $g=1+\tfrac12\sum_e(r_e-1)$ is meromorphic. Let $\widetilde X$ be
the smooth projective model of the double cover of $Y$ defined by
$(2z-a-b)^2=(a-b)^2+4g^2/\Pi$. At the point of $Y$ over $\Pi=0$ with all
$r_e=1$ we have $g=1$ and $\Pi$ is a local coordinate, so the right-hand
side has a simple pole and is not a square in the function field of $Y$.
Hence $\widetilde X$ is irreducible. By \eqref{abs:eq-parallel} the
admissible solutions map injectively to the points with
$\Pi\notin\{0,\infty\}$ and $z\notin\{a,b,\infty\}$ of the curve defined over
$Y$ by this equation, which may be singular; conversely each such
point defines the admissible solution $V_1=g/(a-z)$, $V_2=g/(b-z)$,
$f_e=(r_e-1)/2$. This proves~(1), shows that $\Xc$ is birational to
$\widetilde X$, and so identifies the smooth model of $\Xc$
with $\widetilde X$. On it $m_{1|e}=g/\bigl((a-z)(1+f_e)\bigr)$.

\emph{Points with finite $z$.} We go through the points of $\widetilde X$
according to the value of $\Pi$.

If $\Pi=\infty$, then $t=\Pi^{-1/2}$ is a local coordinate on $Y$,
$r_e=2\epsilon_ew_e/t+O(t)$ for some signs $\epsilon_e$, and
$g=\bigl(\sum_e\epsilon_ew_e\bigr)/t+O(1)$. Thus
$(a-z)(b-z)\to(\sum_e\epsilon_ew_e)^2\ne0$, so $z$ is finite and different
from $a$ and $b$, and $m_{1|e}\to(\sum_{e'}\epsilon_{e'}w_{e'})/
\bigl((a-z)\epsilon_ew_e\bigr)$ is finite.

If $\Pi$ is finite and nonzero, then all $f_e$ are finite and
$1+f_e\ne0$. So $m_{1|e}$ has a pole only if $z=a$, and then $g$ vanishes,
say to order $\ell>0$ in a local coordinate of $\widetilde X$. Since $b-z$
and $\Pi$ are units, \eqref{abs:eq-parallel} gives $\ord(z-a)=2\ell$,
hence $\ord m_{1|e}=-\ell$ and $\ord(m_{1|e}\,dz)=\ell-1\ge0$. The point is
ramified over the $z$-line and carries no residue.

If $\Pi=0$, then $\Pi$ is a local coordinate on $Y$ and $r_e=\pm1$. Let
$N$ be the set of edges with $r_e=-1$, so that $g=1-\abs N$ at this point.
If $\abs N\ne1$, then $g^2/\Pi$ has a pole and $z=\infty$. If $N=\{e\}$,
then
\[
1+f_e=-w_e^2\,\Pi+O(\Pi^2),\qquad f_{e'}=w_{e'}^2\,\Pi+O(\Pi^2)\ (e'\ne e),
\qquad g=c_e\,\Pi+O(\Pi^2),
\]
and $(a-z)(b-z)=c_e^2\,\Pi+O(\Pi^2)$. There are two points of
$\widetilde X$ over this point of $Y$, one with $z=a$ and one with $z=b$,
and $z$ is a local coordinate at both. At the point $Q_e$ with $z=a$,
\[
m_{1|e}=\frac{g}{(a-z)(1+f_e)}=\frac{c_e}{w_e^2\,(z-a)}\bigl(1+O(z-a)\bigr),
\]
while $g/(a-z)\to(b-a)/c_e$ is finite, so $m_{1|e'}$ is finite for
$e'\ne e$; and $m_{2|e'}=g/\bigl((b-z)(1+f_{e'})\bigr)$ tends to
$-c_e/\bigl(w_e^2(b-a)\bigr)$ for $e'=e$ and to $0$ otherwise. At the point
with $z=b$ the function $m_{1|e}$ is finite by the same computation. For
$N=\{e'\}$ with $e'\ne e$ we have just seen that $m_{1|e}$ is finite at
both points. This exhausts the points with finite $z$ and proves~(2).

\emph{The physical sheet.} Let $\mu$ be the spectral measure of the root
$o$ of the half-tree of $m_{1|e}$, so that the physical function is
$m_{1|e}(z)=\int d\mu(x)/(x-z)$ and $-iy\,m_{1|e}(a+iy)\to\mu(\{a\})$ as
$y\downarrow0$. The preimage in $\widetilde X$ of a small punctured disc
around $z=a$ is a disjoint union of punctured discs, one for each point
over $a$, and the physical branch over the upper half of the disc lies in
one of them, around a point $Q$. If $\mu(\{a\})>0$, then $m_{1|e}$ has a
pole at $Q$ and $\abs{m_{1|e}}$ grows like $\abs{z-a}^{-1}$; by~(2) this
excludes the ramified points, where it grows like $\abs{z-a}^{-1/2}$, so
$Q=Q_e$ and $\mu(\{a\})=-\Res_{Q_e}m_{1|e}\,dz=-c_e/w_e^2$; in particular
$c_e<0$. If $\mu(\{a\})=0$, then $y\,m_{1|e}(a+iy)\to0$ and $Q\ne Q_e$. It
remains to show that $c_e<0$ implies $\mu(\{a\})>0$. The root $o$ lies
over the vertex $1$ and has one neighbor for each edge $e'\ne e$; every
other vertex of the half-tree has one neighbor for each edge of $G$. Let
$S$ be the set of vertices $x$ over $1$ such that the labels of the path
from $o$ to $x$ are $(e_1,e,e_2,e,\dots,e_\nu,e)$ with $\nu\ge0$ and
$e_1,\dots,e_\nu\ne e$, and put
$\psi(x)=\prod_{l=1}^{\nu}(-w_{e_l}/w_e)$ for $x\in S$ and $\psi=0$
elsewhere. Then $\psi$ vanishes on the vertices over~$2$. At a vertex $y$
over $2$ whose parent $x_0$ lies in $S$, reached by the edge $e'\ne e$,
the only child in $S$ is the one reached by $e$, and
$w_{e'}\psi(x_0)+w_e\cdot(-w_{e'}/w_e)\psi(x_0)=0$; a vertex over $2$ whose
parent is not in $S$ has no neighbor in $S$. Hence the restricted
operator maps $\psi$ to $a\psi$. Moreover
$\norm\psi^2=\sum_{\nu\ge0}\rho^\nu$ with
$\rho=\sum_{e'\ne e}w_{e'}^2/w_e^2$, which is finite exactly when $c_e<0$.
In that case $\psi$ is an eigenvector with $\psi(o)=1$, so
$\mu(\{a\})\ge1/\norm\psi^2=1-\rho>0$. This proves~(3).
\end{proof}

The theorem assumes $a\ne b$, so the following comparison concerns only
the form of the criterion. For zero potentials the criterion
$w_e^2>\sum_{e'\ne e}w_{e'}^2$ is the condition of Fig\`a-Talamanca and Steger \cite{FTS} for a gap of
$\sigma(\tH)$ around zero, as recalled in \cite[Theorem~B.1]{ABS}, and
\cite[Proposition~B.2]{ABS} states that in this gap exactly one of the
half-tree functions has a physical pole. What
Theorem~\ref{abs:thm-parallel} adds is the computation on all sheets:
every $m_{1|e}$ has a residue pole over $z=a$, whether or not any of them
is physical.

\begin{example}\label{abs:ex-parallel}
Let $k=3$, $a=0$, $b=1$. For $w=(4,5,6)$ we have $c=(45,27,5)$, and the
three functions $m_{1|e}$ have residue poles over $z=0$, at three distinct
non-physical points, with residues $45/16$, $27/25$ and $5/36$. No
half-tree has a physical pole. For $w=(1,2,4)$ we have $c=(19,13,-11)$:
over $z=0$ the physical pole of the half-tree function of the heaviest
edge, with residue $-11/16$, coexists with non-physical residue poles of
the other two, with residues $19$ and $13/4$. In both cases the functions
$m_{2|e}$, which are the ones that appear on the right-hand side of the
recursion \eqref{eq:m-recursion} for $m_{1|e}$, have their residue poles
over $z=1$ only. The same tree operators arise from a simple graph: color
the edges of $K_{3,3}$ by three perfect matchings, give the three weights
to the three colors and the potentials $a,b$ to the two sides; this graph
covers the three-edge dipole compatibly with the coefficients.
\end{example}

\begin{proposition}[Semiregular trees]\label{abs:prop-semiregular}
Let $G$ be a connected bipartite graph in which the vertices of one side
have degree $d_1$ and those of the other side degree $d_2$, with
$2\le d_1<d_2$ (for instance $K_{3,2}$, with $d_1=2$, $d_2=3$). Give all
edges weight $1$ and the two sides the real potentials $a\ne b$. Then the
analytic continuation of the physical half-tree functions is two-sheeted.
Let $m_1$ and $m_2$ denote the functions of the half-trees rooted over the
first and the second side, and $\alpha=d_1-1$, $\beta=d_2-1$. On the
physical sheet the only pole at finite $z$ is a simple pole of $m_1$ at
$z=a$, with residue $-(1-\alpha/\beta)$. On the other sheet the only pole
at finite $z$ is a simple pole of $m_2$ at $z=b$, with residue
$\beta/\alpha-1$. In particular the sets of energies of the physical and
of the non-physical poles are disjoint.
\end{proposition}

\begin{proof}
All half-trees rooted over a given side are isomorphic as rooted trees
with coefficients, so the physical functions depend only on the side of
the root, and
\eqref{eq:m-recursion} reads $1/m_1=a-z-\alpha m_2$,
$1/m_2=b-z-\beta m_1$; for zero potentials this is
\cite[Example~7.2]{ABS}. With $T=m_1m_2$ and $R=(a-z)(b-z)$,
\[
(a-z)m_1=1+\alpha T,\qquad(b-z)m_2=1+\beta T,\qquad
\alpha\beta T^2+(\alpha+\beta-R)T+1=0 .
\]
The discriminant of the quadratic vanishes for
$R=R_\pm=(\sqrt\alpha\pm\sqrt\beta)^2$, two distinct positive numbers, and
each equation $R(z)=R_\pm$ has two distinct real roots. So the discriminant
has four simple zeros in $z$, the quadratic is irreducible over $\C(z)$,
and $T$, $m_1$, $m_2$ are two-valued. For finite $z$ both roots $T$ are
finite, so $m_1$ can have a pole only at $z=a$ and $m_2$ only at $z=b$. At
$z\in\{a,b\}$ we have $R=0$ and the roots are $T=-1/\alpha$ and
$T=-1/\beta$; they are distinct, so both sheets are unramified there. On
the sheet with $T=-1/\beta$ the function $1+\beta T$ vanishes at $z=b$, so
$m_2$ is finite there, and $m_1=(1-\alpha/\beta+O(z-a))/(a-z)$ near $z=a$.
On the sheet with $T=-1/\alpha$ the function $m_1$ is finite at $z=a$ and
$m_2=(1-\beta/\alpha+O(z-b))/(b-z)$ near $z=b$.

It remains to see that the physical sheet takes the value $T=-1/\beta$ at
both $a$ and $b$. The set $\Omega=\{z:R(z)\notin[R_-,\infty)\}$ is the
plane with two real half-lines removed; it is simply connected and
contains the upper half-plane and the interval $[a,b]$. On $\Omega$ the
two roots are distinct holomorphic functions, and the physical $T$ is the
one that tends to $0$ at infinity. On the line $\re z=(a+b)/2$ the value
$R$ is real and runs over $(-\infty,-(a-b)^2/4]$, and on the real segment
between $a$ and $b$ it runs over $[-(a-b)^2/4,0]$. For real $R\le0$ the
two roots are real, negative and distinct, and the one of smaller modulus
depends continuously on $R$, tends to $0$ as $R\to-\infty$, and equals
$-1/\beta$ at $R=0$ because $\beta>\alpha$. Following this path gives the
claim.
\end{proof}

\begin{theorem}[A non-physical pole outside all half-tree spectra]
\label{abs:thm-bridge}
Let $G$ be the graph on the vertices $1,2,3,4$ with three parallel edges
between $1$ and $2$, three parallel edges between $3$ and $4$, and one
edge $13$. Let $\theta_c$ be the coefficients with $D=0$, weight $1$ on
the edges $12$, weight $1/10$ on the edges $34$ and weight $1/20$ on the
edge $13$. There is a neighborhood $W$ of $\theta_c$ in the space of real
coefficients such that for every $\theta\in W$ the following hold.
\begin{enumerate}
\item There are a real number $z_0=z_0(\theta)$, depending continuously on
$\theta$ with $z_0(\theta_c)=3$, and a holomorphic family of solutions of
\eqref{eq:system} over a punctured neighborhood of $z_0$, admissible and
parametrized by $z$, on which $m_{1|13}$ has a simple pole at $z_0$ with
negative residue; the residue is $-1/6$ at $\theta_c$.
\item $z_0(\theta)>\norm{\tH}$. Consequently $z_0$ lies outside the convex
hull of the spectra of $\tH$ and of all its half-tree restrictions, and
every physical half-tree function is holomorphic near $z_0$.
\item If $\theta$ is generic, the family in \textup{(1)} is a branch of the
analytic continuation of the physical branch.
\end{enumerate}
\end{theorem}

\begin{proof}
Write $s=w_{13}$, let $p$ run over the edges $12$ and $q$ over the edges
$34$, and introduce the unknowns
\[
x=-\frac1{V_1},\qquad h=\frac{V_2}{V_1},\qquad F=f_{13},\qquad k=-V_4,
\]
so that $V_2=-h/x$, and the edge equation of the bridge gives $V_3=-xY$
with $Y=F(1+F)/s^2$. The edge equations of the left triple are solved by
$f_p=A_p/x$, where $A_p=-x/2+\sqrt{x^2/4+w_p^2h}$ is the root of
$A^2+xA=w_p^2h$ near $w_p\sqrt h$; those of the right triple by the root
$B_q=f_q$ of $B(1+B)=w_q^2\,xYk$ near $-1$. Both are holomorphic in
$(x,h,F,k)$ and the coefficients near $x=0$, $h=1$. The four vertex
equations become
\begin{equation}\label{abs:eq-bridge}
\begin{aligned}
z-D_1&=x(1+F)+\textstyle\sum_pA_p, &\qquad (z-D_2)\,h&=x+\textstyle\sum_pA_p,\\
(z-D_3)\,xY&=1+F+\textstyle\sum_qB_q, &\qquad (z-D_4)\,k&=1+\textstyle\sum_qB_q .
\end{aligned}
\end{equation}
At $\theta_c$ and $x=0$ we have $A_p=\sqrt h$ and $B_q=-1$, and
$(z,h,F,k)=(3,1,2,-2/3)$ is a solution. The Jacobian matrix of the
differences of the two sides of \eqref{abs:eq-bridge} with respect to
$(z,h,F,k)$ at this point is
\[
\begin{pmatrix}1&-3/2&0&0\\ 1&3/2&0&0\\ 0&0&-1&0\\ -2/3&0&0&3\end{pmatrix},
\]
with determinant $-9$. By the implicit function theorem there is a unique
solution $(z,h,F,k)$ near $(3,1,2,-2/3)$, holomorphic in $x$ and in the
coefficients near $(0,\theta_c)$ and real for real arguments. Put
$z_0(\theta)=z(0;\theta)$. Differentiating the first two equations of
\eqref{abs:eq-bridge} at $\theta_c$ and $x=0$, with
$dA_p/dx=-\tfrac12+\tfrac12h'$, gives $z'=\tfrac32+\tfrac32h'$ and
$z'+3h'=-\tfrac12+\tfrac32h'$, hence $z'(0)=\tfrac12$. So $z$ is a local
coordinate on the family, $V_1=-1/x$ has a simple pole at $z_0$ with
residue $-z'(0)$, and $m_{1|13}=V_1/(1+F)$ has a simple pole with residue
$-z'(0)/(1+F(0))$, equal to $-1/6$ at $\theta_c$. For small $x\ne0$ the
values $V_1=-1/x$, $V_2=-h/x$, $V_3=-xY$ with $Y$ near $2400$, and
$V_4=-k$ are nonzero and $z\ne D_i$, so the solutions are admissible. By
continuity, $z'(0)>0$ and $1+F(0)>0$ on a neighborhood $W$ of
$\theta_c$, which gives~(1).

For~(2), delete all lifts of the edge $13$ from $\tree$. At $\theta_c$ the
remaining components are $3$-regular trees with constant weight $1$ or
$1/10$ and zero potential, of norm $2\sqrt2$ and $2\sqrt2/10$, and the
deleted edges form a matching, which contributes an operator of norm
$1/20$. Hence $\norm{\tH}\le2\sqrt2+1/20<3$ at $\theta_c$. The norm of
$\tH$ depends continuously on $\theta$, so after shrinking $W$ we have
$z_0(\theta)>\norm{\tH}$ on $W$. Every half-tree operator is a compression
of $\tH$ to a subset of the vertices, so its spectrum lies in
$[-\norm{\tH},\norm{\tH}]$.

For~(3), the family consists of admissible solutions and $z$ is not
constant on it, so it lies on $\Xc$; for generic coefficients $\Xc$ is
irreducible by Theorem~\ref{thm:irreducible}, and the generic real
coefficients are dense in $W$.
\end{proof}

At $\theta_c$ the energy $z_0=3$ is the Perron eigenvalue of the Jacobi
matrix of the induced subgraph on $\{1,2\}$, and on the family $V_1,V_2$
have simple poles, $V_3$ has a simple zero and $V_4$ is finite and
nonzero; compare the partially supported poles of
Theorem~\ref{thm:partial-pole}. The orders of all half-tree functions
follow from \eqref{abs:eq-order}; for instance $m_{1|13}$ falls under the
second case. At the
symmetric point $\theta_c$ itself we have not determined whether the
family belongs to $\Cphys$; statement~(3) avoids this question.

A negative residue is the sign of every physical pole, since an atom of
mass $M$ of a spectral measure contributes the residue $-M$. So neither
the position nor the sign of the residue of a non-physical pole of a
half-tree function has to match a physical pole of another half-tree
function: in Theorem~\ref{abs:thm-bridge} no half-tree has spectrum at
$z_0$, in Proposition~\ref{abs:prop-semiregular} the physical and the
non-physical poles are at different energies, and in
Example~\ref{abs:ex-parallel} both coincidence and absence of a physical
partner occur in one family. In this sense the straightforward analogue
of the correspondence on the line fails. What remains is the local
criterion of Proposition~\ref{prop:half-tree-identity}, which together
with the classification of the poles of the $V_i$
(Theorem~\ref{thm:classification}) locates the poles of the half-tree
functions on each sheet separately.

\begin{problem}\label{abs:prob-coincidence}
For a general graph and generic coefficients, describe which poles of the
functions $m_{i|e}$, on which sheets, lie over the same energy.
\end{problem}

\subsection{Open sets of coefficients without gaps}
\label{abs:subsec-gaps}

For a finite graph $G$ with $n$ vertices, a \emph{bounded gap} of
$\sigma(\tH)$ is a bounded connected component of $\R\setminus\sigma(\tH)$.
By Sunada's gap labeling argument \cite{Sunada}, in the form proved for
periodic Jacobi trees in \cite[Theorem~5.1]{ABS}, the integrated density
of states takes a value
$j/n$ with $j\in\{1,\dots,n-1\}$ in each bounded gap; so there are at most
$n-1$ of them, and \cite{ABS} say that \emph{all gaps are open} when there
are exactly $n-1$. The parameter space is $\R^V\times(0,\infty)^E$: all
potentials and all edge weights vary independently. In paraphrase,
Conjecture~10.4 of \cite{ABS} states that for every $G$ the set of
coefficients for which all gaps are open is dense and open in the
parameter space, and Conjecture~10.5 that its complement is a variety of
codimension two. Appendix~B of the arXiv version of \cite{ABS}, following
an observation of Garza-Vargas and Kulkarni \cite{GVK}, already notes that
Conjecture~10.5 fails for graphs with two vertices, where the complement
can have codimension one (see Remark~\ref{abs:rem-two-vertex}). The same
appendix formulates two conjectures of Borg type: if a periodic Jacobi
matrix on a homogeneous tree has all weights equal to $1$ and no gap, then
its potential is constant (Conjecture~B.4); and if a periodic Jacobi
matrix on a tree has no gap, then its potential is constant
(Conjecture~B.5).

\begin{theorem}\label{thm:gapless}
Let $G$ be a connected loopless $d$-regular multigraph with $n$ vertices,
$d\ge3$, and let $\theta_0$ be the coefficients $D=0$, $w\equiv1$.
\begin{enumerate}
\item If $G$ is not bipartite, there is a neighborhood $U$ of $\theta_0$
in $\R^V\times(0,\infty)^E$ such that $\sigma(\tH)$ is an interval for all
coefficients in $U$.
\item If $G$ is bipartite, there is a neighborhood $U$ of $\theta_0$ in
$\R^V\times(0,\infty)^E$ such that for all coefficients in $U$ the
spectrum $\sigma(\tH)$ has at most one bounded gap, and the integrated
density of states equals $1/2$ in it. The coefficients in $U$ with exactly
one bounded gap contain a nonempty open set.
\end{enumerate}
\end{theorem}

The proof occupies the rest of this subsection. We write $\varrho=d-1$,
$\theta=(D,w)$, and $H_\theta$ for the operator $\tH$ with coefficients
$\theta$. The Schur test gives
\begin{equation}\label{abs:eq-norm}
\norm{H_\theta-H_{\theta_0}}\le\max_i\abs{D_i}+d\,\max_e\abs{w_e-1},
\end{equation}
and the same bound holds for the half-tree restrictions.

\emph{The uniform point.} At $\theta_0$ the operator is the adjacency
operator of the $d$-regular tree. All half-tree functions coincide with
the solution $m_0$ of $\varrho m_0^2+zm_0+1=0$ with $m_0\sim-1/z$ at
infinity, and $G_i=G_0=(-z-dm_0)^{-1}$ by \eqref{eq:m-recursion}. The
function $m_0$ extends continuously to the closed upper half-plane, with
\begin{equation}\label{abs:eq-m0}
m_0(E)=\frac{-E+i\sqrt{4\varrho-E^2}}{2\varrho},\qquad
\abs{m_0(E)}^2=\frac1\varrho\qquad(\abs E<2\sqrt\varrho),
\end{equation}
and with real values for $\abs E\ge2\sqrt\varrho$. The denominator
$-E-dm_0(E)$ has negative imaginary part for $\abs E<2\sqrt\varrho$, and
for $\abs E\ge2\sqrt\varrho$ it equals
$-\bigl((d-2)E+d\,\sgn(E)\sqrt{E^2-4\varrho}\bigr)/(2\varrho)\ne0$. So $G_0$
is continuous on the closed upper half-plane, and the density of states
measure $\mu_0$ of $\theta_0$ (the Kesten--McKay law) is absolutely
continuous, with a continuous density that is positive on
$(-2\sqrt\varrho,2\sqrt\varrho)$ and vanishes outside. It is symmetric
about $0$ because the tree is bipartite.

\emph{The linearized recursion.} Let $\vec E$ be the set of the $2\abs E$
oriented edges of $G$. For $\gamma\in\vec E$ let $\bar\gamma$ be the
reversed oriented edge, and let $B$ be the non-backtracking matrix on
$\vec E$: $B_{\gamma\gamma'}=1$ if the head of $\gamma$ is the tail of
$\gamma'$ and $\gamma'\ne\bar\gamma$, and $B_{\gamma\gamma'}=0$ otherwise.
For $\gamma$ oriented from $i$ to $j$ along the edge $e$ we write
$m_\gamma=m_{j|e}$, $w_\gamma=w_e$ and $D_\gamma=D_j$. The second equation
of \eqref{eq:m-recursion} reads
\begin{equation}\label{abs:eq-Phi}
\Phi_\gamma(m;\theta,z):=m_\gamma\Bigl(D_\gamma-z-\sum_{\gamma'}
B_{\gamma\gamma'}\,w_{\gamma'}^2\,m_{\gamma'}\Bigr)-1=0\qquad(\gamma\in\vec E).
\end{equation}

\begin{lemma}\label{abs:lem-nonbacktracking}
The eigenvalues $\lambda$ of $B$ satisfy $\abs\lambda\le\varrho$. If $G$ is
not bipartite, the only eigenvalue of modulus $\varrho$ is $\varrho$; if
$G$ is bipartite, they are $\varrho$ and $-\varrho$. Consequently, for
$\abs E<2\sqrt\varrho$ the matrix $I-m_0(E)^2B$ is invertible, except when
$G$ is bipartite and $E=0$.
\end{lemma}

\begin{proof}
All row sums and column sums of $B$ are $\varrho$, so the spectral radius
is $\varrho$. We show that $B$ is irreducible. Fix a vertex $j$ and
consider the bipartite graph whose two sides are the $d$ oriented edges
with head $j$ and the $d$ oriented edges with tail $j$, with $\gamma$
adjacent to $\gamma'$ when $B_{\gamma\gamma'}=1$. It is $K_{d,d}$ minus a
perfect matching, which is connected for $d\ge3$. An oriented edge from
$i$ to $j$ occurs in the graph of $i$ and in the graph of $j$, and $G$ is
connected; so the directed graph of $B$ is weakly connected. Every node of
it has in-degree and out-degree $\varrho$. A set of nodes with no entering
arc then has no leaving arc either, since the in-degrees and the
out-degrees over the set have the same sum; applied to a strongly
connected component without entering arcs, this shows that the directed
graph is strongly connected.

By the Perron--Frobenius theorem the eigenvalues of modulus $\varrho$ are
$\varrho\,e^{2\pi il/h}$, $l=0,\dots,h-1$, where $h$ is the period of $B$,
and there is a map $c\colon\vec E\to\Z/h$ with $c(\gamma')=c(\gamma)+1$
whenever $B_{\gamma\gamma'}=1$. Two oriented edges with the same head $j$
have a common successor because $d\ge3$, so $c$ takes one value $c_j$ on
the oriented edges with head $j$, and the value $c_j+1$ on those with tail
$j$. For an edge $ij$, its two orientations give $c_j=c_i+1$ and
$c_i=c_j+1$, so $h$ divides~$2$. If $h=2$, then $j\mapsto c_j$ is a proper
$2$-coloring and $G$ is bipartite. Conversely, if $G$ is bipartite, every
closed walk in the directed graph of $B$ projects to a closed walk of the
same length in $G$ and so has even length, whence $h=2$.

Finally let $\abs E<2\sqrt\varrho$ and write
$m_0(E)=e^{i\varphi}/\sqrt\varrho$ with $0<\varphi<\pi$, using
\eqref{abs:eq-m0}. The eigenvalues of $I-m_0(E)^2B$ are $1-m_0(E)^2\lambda$,
and $\abs{m_0(E)^2\lambda}<1$ unless $\abs\lambda=\varrho$. For
$\lambda=\varrho$ we get $1-e^{2i\varphi}\ne0$. For $\lambda=-\varrho$ we
get $1+e^{2i\varphi}$, which vanishes only for $\varphi=\pi/2$, that is,
for $E=0$.
\end{proof}

\begin{lemma}[Stability inside the band]\label{abs:lem-interior}
Let $K\subset(-2\sqrt\varrho,2\sqrt\varrho)$ be compact, with $0\notin K$ if
$G$ is bipartite. Then there is a neighborhood $U_K$ of $\theta_0$ in
$\R^V\times(0,\infty)^E$ such that $K\subset\sigma(H_\theta)$ for all
$\theta\in U_K$.
\end{lemma}

\begin{proof}
Fix $E\in K$ and write $\mathbf m_0(z)$ for the vector with all entries
$m_0(z)$. At $(\mathbf m_0(E);\theta_0,E)$ the map \eqref{abs:eq-Phi}
vanishes, and
\[
\frac{\partial\Phi_\gamma}{\partial m_{\gamma'}}
=\delta_{\gamma\gamma'}\bigl(-E-\varrho m_0(E)\bigr)-m_0(E)B_{\gamma\gamma'}
=\frac1{m_0(E)}\bigl(I-m_0(E)^2B\bigr)_{\gamma\gamma'},
\]
which is invertible by Lemma~\ref{abs:lem-nonbacktracking}. By the
holomorphic implicit function theorem there are $r>0$, $\delta>0$ and a
neighborhood $U$ of $\theta_0$ such that for all $\theta\in U$ and all $z$
in the disc $\Delta=\{\abs{z-E}<\delta\}$ the equation
$\Phi(m;\theta,z)=0$ has exactly one solution $m=M(\theta,z)$ with
$\max_\gamma\abs{m_\gamma-m_0(E)}<r$; moreover $M$ is continuous in
$(\theta,z)$ and holomorphic in $z$. We choose $\delta$ so small that
$\abs{m_0(z)-m_0(E)}<r/2$ for $z\in\Delta$ with $\im z\ge0$, where on the
real axis $m_0$ is the boundary value \eqref{abs:eq-m0}; by uniqueness
$M(\theta_0,z)=\mathbf m_0(z)$ for these $z$.

We claim that, after shrinking $U$, for $\theta\in U$ the physical
half-tree functions $m^\theta(z)$ coincide with $M(\theta,z)$ on
$\Delta^+=\Delta\cap\{\im z>0\}$. Fix $z_*\in\Delta^+$. By the resolvent
identity and \eqref{abs:eq-norm},
$\abs{m^\theta_\gamma(z_*)-m_0(z_*)}\le\norm{H_\theta-H_{\theta_0}}/(\im z_*)^2$,
so for $\theta$ in a smaller neighborhood
$\abs{m^\theta_\gamma(z_*)-m_0(E)}<r$ for all $\gamma$. The physical
functions are continuous in $z$ and solve \eqref{abs:eq-Phi}, so for $z$
near $z_*$ they lie in the ball of radius $r$ and, by uniqueness, equal
$M(\theta,z)$. Both sides are holomorphic on the connected set $\Delta^+$,
so they agree there.

Consequently, for $\theta\in U$ and $x$ over $i$, the Green function
$G^\theta_i(z)=\bigl(D_i-z-\sum_\gamma w_\gamma^2M_\gamma(\theta,z)\bigr)^{-1}$,
where the sum is over the oriented edges with tail $i$, extends
continuously from $\Delta^+$ to the real points of $\Delta$ wherever the
denominator does not vanish. At $\theta_0$ the imaginary part of the
denominator is $-\im z-d\,\im m_0(z)\le-d\,\im m_0(z)$, which is bounded
above by a negative constant on $\{z\in\Delta:\abs{z-E}\le\delta/2,\ \im
z\ge0\}$ if $\delta$ is small. By continuity of $M$ and compactness, after
shrinking $U$ again the imaginary part of the denominator is negative on
this set for all $\theta\in U$. Hence $\im G^\theta_i(E'+i\eta)$ converges,
uniformly in $E'\in J_E=(E-\delta/2,E+\delta/2)$ as $\eta\downarrow0$, to a
continuous positive function. By Stieltjes inversion the spectral measure
of $\delta_x$ has a positive density on $J_E$, so
$J_E\subset\sigma(H_\theta)$. Cover $K$ by finitely many intervals $J_E$
and intersect the corresponding neighborhoods.
\end{proof}

The stability of the physical solution of the recursion in the interior
of a band is a known phenomenon; compare the stability results for vector
Dyson equations in \cite{HuangSuXu}. What we need in addition is a way to
exclude gaps near the band edges, where
Lemma~\ref{abs:lem-interior} says nothing. This is supplied by the
discreteness of the gap labels.

\emph{The integrated density of states.} Let $\mu^\theta_i$ be the spectral
measure of $\delta_x$ for $x$ over $i$,
$\mu_\theta=\frac1n\sum_{i\in V}\mu^\theta_i$ the density of states measure,
and $N_\theta(E)=\mu_\theta((-\infty,E])$. We use three facts.
\begin{enumerate}
\item[(a)] The support of $\mu_\theta$ is $\sigma(H_\theta)$. Indeed, if
$\mu_\theta(I)=0$ for an open interval $I$, the spectral projection $P_I$
of $H_\theta$ satisfies $\norm{P_I\delta_x}^2=0$ for all vertices $x$ of
$\tree$, so $P_I=0$.
\item[(b)] On each bounded gap $N_\theta$ is constant with a value in
$\{1/n,\dots,(n-1)/n\}$ \cite[Theorem~5.1]{ABS}, and distinct bounded gaps
have distinct values, since by (a) the measure $\mu_\theta$ charges every
neighborhood of a point of the spectrum between them. The hypotheses of
\cite[Theorem~5.1]{ABS}, a finite connected graph without leaves, positive
weights and real potentials, are satisfied here.
\item[(c)] For every real $E$, $N_\theta(E)\to N_{\theta_0}(E)$ as
$\theta\to\theta_0$. Indeed, by \eqref{abs:eq-norm} the moments
$\langle\delta_x,H_\theta^l\delta_x\rangle$ are continuous in $\theta$ and
the supports are uniformly bounded, so $\mu_\theta\to\mu_{\theta_0}$
weakly; and $\mu_{\theta_0}=\mu_0$ has no atoms.
\end{enumerate}

\begin{proof}[Proof of Theorem~\ref{thm:gapless}]
(1) Since $N_{\theta_0}$ is continuous with
$N_{\theta_0}(-2\sqrt\varrho)=0$ and $N_{\theta_0}(2\sqrt\varrho)=1$, we
can choose $-2\sqrt\varrho<L<R<2\sqrt\varrho$ with $N_{\theta_0}(L)<1/n$
and $N_{\theta_0}(R)>1-1/n$. By Lemma~\ref{abs:lem-interior} with
$K=[L,R]$ and by (c), there is a neighborhood $U$ of $\theta_0$ such
that for $\theta\in U$
\[
[L,R]\subset\sigma(H_\theta),\qquad N_\theta(L)<\frac1n,\qquad
N_\theta(R)>1-\frac1n .
\]
Let $\theta\in U$ and suppose that $I$ is a bounded gap of
$\sigma(H_\theta)$. Then $I$ is disjoint from $[L,R]$. If
$I\subset(-\infty,L)$, the value of $N_\theta$ on $I$ is at most
$N_\theta(L)<1/n$; if $I\subset(R,\infty)$, it is at least
$N_\theta(R)>(n-1)/n$. Both contradict (b). So $\sigma(H_\theta)$ is a
compact set whose complement has no bounded component, that is, an
interval.

(2) Now $G$ is bipartite and $d$-regular, so the two sides have the same
number of vertices and $n$ is even. By the symmetry of $\mu_0$,
$N_{\theta_0}(0)=1/2$. Choose $L<L_0<0<R_0<R$ in the open band with
\[
N_{\theta_0}(L)<\frac1n,\qquad
\frac12-\frac1n<N_{\theta_0}(L_0),\qquad
N_{\theta_0}(R_0)<\frac12+\frac1n,\qquad N_{\theta_0}(R)>1-\frac1n .
\]
By Lemma~\ref{abs:lem-interior} with $K=[L,L_0]\cup[R_0,R]$ and by (c),
on a neighborhood $U$ of $\theta_0$ these four inequalities hold for
$N_\theta$ and $K\subset\sigma(H_\theta)$. For $\theta\in U$, a bounded gap
cannot lie in $(-\infty,L)$ or in $(R,\infty)$, as in~(1). A bounded gap
in $(L_0,R_0)$ carries a value of $N_\theta$ in
$[N_\theta(L_0),N_\theta(R_0)]\subset(\tfrac12-\tfrac1n,\tfrac12+\tfrac1n)$,
which by (b) is $1/2$; and by (b) there is at most one such gap.

To see that a gap does open, give the two sides of $G$ the potentials
$\varepsilon$ and $-\varepsilon$, with $\varepsilon>0$, and keep
$w\equiv1$. With respect to the bipartition of $\tree$,
\[
H_\theta=\begin{pmatrix}\varepsilon&T\\ T^*&-\varepsilon\end{pmatrix},\qquad
H_\theta^2=\begin{pmatrix}\varepsilon^2+TT^*&0\\0&\varepsilon^2+T^*T\end{pmatrix}
\ge\varepsilon^2,
\]
so $(-\varepsilon,\varepsilon)$ does not meet $\sigma(H_\theta)$, while the
Rayleigh quotients $\pm\varepsilon$ of the basis vectors show that
$\sigma(H_\theta)$ meets both $[\varepsilon,\infty)$ and
$(-\infty,-\varepsilon]$. For every $\theta'$ with
$\norm{H_{\theta'}-H_\theta}<\varepsilon/2$ the interval
$(-\varepsilon/2,\varepsilon/2)$ is free of spectrum and the Rayleigh
quotients of the basis vectors have modulus larger than $\varepsilon/2$
and both signs, so $\sigma(H_{\theta'})$ has a bounded gap. For small
$\varepsilon$ we have $\theta\in U$, and by the analogue of
\eqref{abs:eq-norm} for the pair $\theta',\theta$ these $\theta'$ contain
a neighborhood of $\theta$; its intersection with $U$ is a nonempty open
set of coefficients with exactly one bounded gap.
\end{proof}

For $G=K_4$ the choices in part~(1) can be made explicit: $L=-5/2$ and
$R=5/2$ lie in the band $(-2\sqrt2,2\sqrt2)$; the measure $\mu_0$ is
symmetric with second moment $\norm{H_{\theta_0}\delta_x}^2=3$, so by
Chebyshev's inequality each of the tails $(-\infty,-5/2]$ and
$[5/2,\infty)$ has mass at most $\tfrac12\cdot3\cdot\tfrac4{25}=6/25<1/4$.

\begin{corollary}\label{cor:abs-conjectures}
Let $G$ be a connected loopless non-bipartite $d$-regular multigraph with
$d\ge3$ and $n$ vertices, for instance the simple graph $K_4$.
\begin{enumerate}
\item The set of coefficients in $\R^V\times(0,\infty)^E$ for which
$\sigma(\tH)$ has $n-1$ bounded gaps is not dense, and its complement
contains a nonempty open set. Hence the statements of Conjectures~10.4
and~10.5 of \cite{ABS} do not hold for~$G$; the exceptional set of
Conjecture~10.5 has codimension zero.
\item There are non-constant potentials $D$ such that the periodic Jacobi
operator with these potentials and all weights equal to $1$, which acts on
the $d$-regular tree, has spectrum equal to an interval. Hence the
statements of Conjectures~B.4 and~B.5 of the arXiv version of \cite{ABS}
do not hold.
\end{enumerate}
\end{corollary}

\begin{proof}
Since $d\ge3$ we have $n\ge2$, so $n-1\ge1$ bounded gaps are required,
while the neighborhood $U$ of Theorem~\ref{thm:gapless}(1) consists of
coefficients without any bounded gap. This gives~(1). For~(2) take any
$(D,\ones)\in U$ with $D$ non-constant.
\end{proof}

To our knowledge the failure of Conjecture~10.4, and of
Conjectures~B.4 and~B.5, has not been observed before; the case of two
vertices of Conjecture~10.5 was refuted in the appendix of \cite{ABS} and
by \cite{GVK}. Problem~B.6 of that appendix, which asks for the
codimension of the exceptional set for $n\ge3$, is formulated under the
assumption of Conjecture~10.4. Corollary~\ref{cor:abs-conjectures} does
not contradict the result of \cite[\S10.2]{ABS} that all gaps are open
when the potentials are distinct and the sum of the weights is smaller
than the smallest difference between two potentials. Gap labeling
then gives all labels $1/n,\ldots,(n-1)/n$ simultaneously on a nonempty
open set for every base graph. For $K_4$, the parameter space also has
a nonempty open set where none of the gaps is open. We have not determined
the regions with other numbers of gaps, nor the structure of their
boundaries. For bipartite regular graphs, Theorem~\ref{thm:gapless}(2)
shows that near the uniform point only the gap with label $1/2$ can open,
so the statement of Conjecture~10.4 fails for them as well as soon as
$n\ge4$.

\begin{remark}[Two vertices]\label{abs:rem-two-vertex}
Let $G$ have two vertices joined by $k\ge2$ parallel edges, with
potentials $D_1,D_2$ and weights $w_e>0$. The set of coefficients without
a bounded gap is exactly
\[
\Bigl\{D_1=D_2,\quad 2\max_ew_e^2\le\sum_ew_e^2\Bigr\}.
\]
Indeed, if $D_1\ne D_2$, the block argument in the proof of
Theorem~\ref{thm:gapless}(2), applied to $\tH-\tfrac12(D_1+D_2)$ with
$\varepsilon=\tfrac12\abs{D_1-D_2}$, gives a gap. If $D_1=D_2$, the
spectrum is symmetric about $D_1$ and there is at most one bounded gap,
with label $1/2$; so a gap, if present, contains $D_1$, and by the theorem
of Fig\`a-Talamanca and Steger \cite{FTS}, in the form of
\cite[Theorem~B.1]{ABS}, the point $D_1$ belongs to the spectrum if and
only if $\max_ew_e^2\le\sum_ew_e^2-\max_ew_e^2$. For $k\ge3$ this set has
codimension one, which is the counterexample to Conjecture~10.5 of
\cite[Appendix~B]{ABS}; for $k=2$ it is the set $D_1=D_2$, $w_1=w_2$ of
codimension two. In particular, for loopless graphs with two vertices the
coefficients without a gap never contain an open set; this answers, in
this case, a question raised in a remark of the same appendix. The
codimension-zero example of Corollary~\ref{cor:abs-conjectures} needs at
least three vertices.
\end{remark}

\begin{remark}[What we could not check]\label{abs:rem-literature}
We compared our statements with the text of the third arXiv version of
\cite{ABS}, including its Appendix~B; the numbers B.4 and B.5 refer to
that version, and we have not verified whether the published version
contains this appendix. We searched the papers citing \cite{ABS} that were
accessible to us for a discussion of Conjectures~10.4, B.4 and~B.5 and
found none, but this search was bibliographic and not exhaustive. We have
not read the literature on the stability of absolutely continuous spectrum
on trees, in particular the work of Klein \cite{Klein}, of Aizenman, Sims
and Warzel \cite{AizenmanSimsWarzel} and of Keller, Lenz and Warzel
\cite{KellerLenzWarzel}. As far as we know these works concern mainly
random perturbations, but the argument for Theorem~\ref{thm:gapless} is
short and elementary, and we cannot exclude that the existence of open
sets of periodic coefficients without gaps is noted there or elsewhere.
\end{remark}

%% file: sections/block-real.tex
\providecommand{\Sym}{\mathrm{Sym}}
\providecommand{\nuc}[1]{\lVert#1\rVert_{*}}
\providecommand{\op}[1]{\lVert#1\rVert_{\mathrm{op}}}
\providecommand{\cC}{\mathcal C}
\providecommand{\lamC}{\lambda_{\mathrm C}}
\providecommand{\ran}{\operatorname{ran}}
\providecommand{\cP}{\mathcal P}

\section{Block models: the cone condition and uniqueness}
\label{sec:block-real}

\subsection{Comparison with the scalar theory}\label{blr:subsec:scalar-comparison}

A variable in $\R^{p_i}$ has a symmetric $p_i\times p_i$ potential $D_i$.
An edge $e=ij$ has a $p_i\times p_j$ weight $A_e$, which may be
rectangular or singular. We compare Part~II with the scalar theory through
the threshold for uniqueness, the sheet count and the residue formula.

\emph{Real fixed points.} The cone condition at energy $z$ is
$z<\lamC(D,A)$. The threshold $\lamC$ is obtained by minimizing the
lowest eigenvalue of the finite twisted Jacobi operator over edge
matrices in $O(m)$, with $m=\max_i p_i$ (Lemma~\ref{blr:lem-cone}).
For scalar variables $O(1)=\{\pm1\}$, so this minimization is over
signings and gives $\lambda_{\min}(D-\abs A)$, the lower Perron value
$\sigm$ when $G$ is connected. Below the threshold the positive
solution is unique for arbitrary real coefficients. For generic real
coefficients the converse and the parity rule also hold
(Theorems~\ref{thm:block-cone-unique} and~\ref{thm:block-converse}).
At the threshold, the cover Green function decides whether the unique
possible solution exists (Theorem~\ref{thm:block-boundary}). On a
connected graph with $\beta\ge2$, invertible edge blocks give
$\lamC<\Sigm$. Singular blocks can close this interval on $K_4$ even with
a positive definite finite precision matrix
(Corollary~\ref{blr:cor-cycle-rank} and Example~\ref{blr:ex-K4-singular}).
The scalar real arc to the Perron pole gives at least two positive
solutions throughout the gap for every choice of coefficients. In the
block setting the signed count gives this bound for generic coefficients
and away from finitely many energies
(Corollaries~\ref{cor:two-roots} and~\ref{blr:cor-two}).

\emph{Complex sheets.} The full covariance curve remains generically
irreducible for arbitrary vertex dimensions
(Theorem~\ref{thm:block-irreducible}). At uniform dimension $p$, write
$N_p(G)$ for its degree over the $z$-line. Direct sums of scalar channels
give the lower bound $N_p(G)\ge\LG^p$. This bound is attained on a cycle.
On the bowtie, formed by two triangles sharing a vertex, dimension two
gives $36$ sheets compared with the $16$ supplied by the scalar direct
sum (Theorem~\ref{thm:block-degree}).

\emph{Poles and residues.} The Floquet identity on every sheet uses the
trace differential $(\sum_i\tr V_i)\,dz$ and still gives integer residues
(Proposition~\ref{prop:block-floquet}). On a connected graph with uniform
dimension and invertible edge weights, consider an unramified point at
which all cavity precisions remain finite and invertible and the vertex
covariances have at most simple poles. If some $V_i$ has a pole, the
vertex residues have a common rank $r\ge1$. The energy is an eigenvalue
of a finite Jacobi matrix twisted by an $O(r,\C)$ local system, and the
trace residue is $r(\abs E-n)$.
Rank one gives signed block Jacobi matrices. For higher rank the twist
records orthogonal transport between residue spaces. These are
necessary conditions on such a pole
(Theorem~\ref{thm:block-twisted-poles}).

\subsection{The block covariance system}
\label{blr:subsec-system}

Let $G=(V,E)$ be a finite multigraph without self-loops, with at least one
vertex, and let $p=(p_i)_{i\in V}$ be positive integers. We fix an
orientation of every edge and write $e=ij$ for an edge oriented from $i$
to $j$. The coefficients are
\[
D_i\in\Sym_{p_i}(\R)\quad(i\in V),\qquad A_e\in\R^{p_i\times p_j}\quad(e=ij\in E),
\]
where $\Sym_r$ denotes symmetric $r\times r$ matrices. No rank condition
is imposed on $A_e$. For the two ends of $e=ij$ we write $A_e^{ij}=A_e$ and
$A_e^{ji}=A_e^{\T}$; in a sum $\sum_{e\ni i}$ the other end of $e$ is
denoted by $k$, and parallel edges are counted separately. The space of
real coefficients is $\cP_\R=\prod_i\Sym_{p_i}(\R)\times\prod_{e=ij}\R^{p_i\times p_j}$,
and $\cP$ is its complexification. A statement holds for \emph{generic
real coefficients} if it holds outside a proper real algebraic subset of
$\cP_\R$.

Gaussian belief propagation for the quadratic form
$\tfrac12\sum_ix_i^{\T}D_ix_i+\sum_{e=ij}x_i^{\T}A_ex_j$ iterates one
symmetric matrix for each end of each edge. With the diagonal shifted by
a spectral parameter, a fixed point consists of \emph{cavity precisions}
$Q_{i|e}\in\Sym_{p_i}$, one for each edge $e$ and each end $i$ of $e$,
satisfying
\begin{equation}\label{blr:eq-cavity}
Q_{i|e}=D_i-z-\sum_{\substack{e'\ni i\\ e'\ne e}}
A_{e'}^{ik}\,Q_{k|e'}^{-1}\,A_{e'}^{ki}\qquad(e\in E,\ i\in e).
\end{equation}
Only the edge $e$ itself is omitted from the sum; edges parallel to $e$
are kept. The \emph{vertex precision} and the \emph{edge precision} of a
fixed point are
\begin{equation}\label{blr:eq-RE}
R_i=D_i-z-\sum_{e\ni i}A_e^{ik}Q_{k|e}^{-1}A_e^{ki},\qquad
E_e=\begin{pmatrix}Q_{i|e}&A_e\\ A_e^{\T}&Q_{j|e}\end{pmatrix}\quad(e=ij).
\end{equation}
For $p_i=1$, $Q_{i|e}=q_{i\to j}$ and \eqref{blr:eq-cavity} is
\eqref{eq:cavity} for the diagonal $D-z$. As in Part~I we consider only
the precision part of the algorithm.

The analogue of the unknowns $(V,f)$ of \eqref{eq:system} are the local
covariances. For $V_i\in\Sym_{p_i}(\C)$ and $C_e\in\C^{p_i\times p_j}$ put
\[
S_e=\begin{pmatrix}V_i&C_e\\ C_e^{\T}&V_j\end{pmatrix}\qquad(e=ij),
\]
and write $(M)_{ii}$, $(M)_{ij}$ for the blocks of a matrix $M$ partitioned
as $S_e$.

\begin{definition}\label{def:block-system}
Let $Y_p$ be the set of $y=(V,C)\in\prod_i\Sym_{p_i}(\C)\times
\prod_{e=ij}\C^{p_i\times p_j}$ with $\det V_i\ne0$ for all $i$ and
$\det S_e\ne0$ for all $e$. The \emph{block covariance system} of
$(G,p,D,A)$ is the system
\begin{equation}\label{eq:block-system}
\bigl(S_e^{-1}\bigr)_{ij}=A_e\quad(e=ij\in E),\qquad
\sum_{e\ni i}\bigl(S_e^{-1}\bigr)_{ii}-(d_i-1)\,V_i^{-1}=D_i-z\quad(i\in V)
\end{equation}
in the unknowns $(z,V,C)\in\C\times Y_p$, where $D_i-z$ stands for
$D_i-zI_{p_i}$. We write $\Psi\colon Y_p\to\cP$ for the map that sends
$(V,C)$ to the left-hand sides of \eqref{eq:block-system}, so that
\eqref{eq:block-system} reads $\Psi(V,C)=(D-z,A)$.
\end{definition}

\begin{definition}\label{def:block-admissible}
An \emph{admissible solution} is a point $(z,V,C)\in\C\times Y_p$ that
satisfies \eqref{eq:block-system}; the conditions $\det V_i\ne0$ and
$\det S_e\ne0$ are part of the definition. For real coefficients, a real
admissible solution is \emph{positive} if $V_i\succ0$ for all $i$ and
$S_e\succ0$ for all $e$. A \emph{positive solution of $(D,A)$} is a point
$(V,C)$ such that $(0,V,C)$ is a positive solution; the positive solutions
with a given real $z$ are the positive solutions of $(D-z,A)$.
\end{definition}

\begin{lemma}\label{blr:lem-bp}
Fix $z\in\C$. The map $Q_{i|e}=(S_e^{-1})_{ii}$ is a bijection from the
admissible solutions of \eqref{eq:block-system} with this $z$ onto the
solutions of \eqref{blr:eq-cavity} for which all $Q_{i|e}$, $R_i$ and $E_e$
are nonsingular. Its inverse is $V_i=R_i^{-1}$, $S_e=E_e^{-1}$. For real
coefficients and real $z$ the positive solutions correspond to the fixed
points with $Q_{i|e}\succ0$, $R_i\succ0$ and $E_e\succ0$, that is, to the
fixed points all of whose vertex and edge beliefs are positive definite.
\end{lemma}

\begin{proof}
Let $(z,V,C)$ be admissible and $Q_{i|e}=(S_e^{-1})_{ii}$. Since $S_e$ and
its block $V_j$ are nonsingular, $Q_{i|e}=(V_i-C_eV_j^{-1}C_e^{\T})^{-1}$
is nonsingular, and $E_e:=S_e^{-1}$ has the form \eqref{blr:eq-RE} by the
first equation of \eqref{eq:block-system}. The Schur complement of the
block $Q_{j|e}$ in $E_e$ is the inverse of $(E_e^{-1})_{ii}=V_i$:
\begin{equation}\label{blr:eq-schur}
Q_{i|e}-A_e^{ij}Q_{j|e}^{-1}A_e^{ji}=V_i^{-1}.
\end{equation}
Substituting \eqref{blr:eq-schur} for all edges $e'\ne e$ at $i$ into the
second equation of \eqref{eq:block-system} gives \eqref{blr:eq-cavity},
and substituting it for all edges gives $R_i=V_i^{-1}$. Conversely, let
$Q$ solve \eqref{blr:eq-cavity} with nonsingular $Q_{i|e}$, $R_i$, $E_e$,
and put $V_i=R_i^{-1}$, $S_e=E_e^{-1}$. By \eqref{blr:eq-cavity} and
\eqref{blr:eq-RE}, $Q_{i|e}-A_e^{ij}Q_{j|e}^{-1}A_e^{ji}=R_i$, so the
diagonal blocks of $S_e$ are $R_i^{-1}=V_i$ and $R_j^{-1}=V_j$; thus
$(V,C)$ with $C_e=(S_e)_{ij}$ lies in $Y_p$, $(S_e^{-1})_{ij}=A_e$, and
$\sum_{e\ni i}Q_{i|e}-(d_i-1)R_i=R_i+\sum_{e\ni i}A_e^{ik}Q_{k|e}^{-1}A_e^{ki}=D_i-z$.
The two maps are mutually inverse. In the real case, $S_e\succ0$ gives
$E_e\succ0$ and $Q_{i|e}\succ0$, and $V_i\succ0$ gives $R_i\succ0$; the
converse is the same.
\end{proof}

\begin{remark}\label{blr:rem-scalar}
Let $p_i=1$ for all $i$ and $A_e=w_e\ne0$. If $(z,V,f)$ is an admissible
solution of \eqref{eq:system}, then $C_e=-w_eV_iV_j/(1+f_e)$ defines an
admissible solution of \eqref{eq:block-system}, with
$\det S_e=V_iV_j/(1+f_e)$, and every admissible solution of
\eqref{eq:block-system} with $z\notin\{D_i\}$ arises in this way from
$f_e=C_e^2/\det S_e$. Positive solutions correspond to positive solutions.
The block system records the sign of the covariance $C_e$, which is
determined by the sign of $w_e$.
\end{remark}

The positive solutions are the critical points of a Bethe function. Let
$\cC$ be the open convex cone of real $y=(V,C)$ with $V_i\succ0$ and
$S_e\succ0$ for all $i$ and $e$; its closure $\overline\cC$ is defined by
$V_i\succeq0$, $S_e\succeq0$. For real coefficients put
\begin{equation}\label{blr:eq-bethe}
L(y)=\tfrac12\sum_i\tr(D_iV_i)+\sum_{e}\tr(A_e^{\T}C_e),\qquad
F(y)=L(y)-\tfrac12\sum_e\log\det S_e+\tfrac12\sum_i(d_i-1)\log\det V_i
\end{equation}
for $y\in\cC$. As $(D,A)$ runs over $\cP_\R$, $L$ runs over all linear
forms on the space of real $y$.

\begin{lemma}\label{blr:lem-critical}
The positive solutions of $(D,A)$ are the critical points of $F$ in $\cC$.
More precisely $\nabla F(y)=\iota\bigl((D,A)-\Psi(y)\bigr)$ for a fixed
linear isomorphism $\iota$, so that a positive solution $y$ is a
nondegenerate critical point of $F$ if and only if $d\Psi(y)$ is
invertible.
\end{lemma}

\begin{proof}
For a variation $\delta C_e$ we have
$d(-\tfrac12\log\det S_e)=-\tr\bigl(((S_e^{-1})_{ij})^{\T}\delta C_e\bigr)$,
and for a variation $\delta V_i$,
$dF=\tfrac12\tr\bigl((D_i-\sum_{e\ni i}(S_e^{-1})_{ii}+(d_i-1)V_i^{-1})\delta V_i\bigr)$.
Hence $dF(y)[\delta y]=\langle (D,A)-\Psi(y),\delta y\rangle$ for the
pairing $\langle(D,A),(V,C)\rangle=\tfrac12\sum_i\tr(D_iV_i)+\sum_e\tr(A_e^{\T}C_e)$,
which is nondegenerate.
\end{proof}

We call a positive solution \emph{nondegenerate} if it is a nondegenerate
critical point of $F$, and its \emph{index} is the number of negative
eigenvalues of the Hessian of $F$. For fixed $V$ the function $F$ is
strictly convex in $C$, so the index is also the number of negative
eigenvalues of the Hessian of the reduced function $\min_CF$ of the vertex
covariances; for $p_i=1$ the latter is, up to a constant, the reduced
Bethe function of the vertex variances whose critical points are counted
in Section~\ref{sec:real}.

\subsection{The cover operator and the physical branch}
\label{blr:subsec-physical}

From now on in this section the coefficients are real. Let $G$ be
connected and $\pi\colon\tree\to G$ its universal cover. For a vertex $x$
of $\tree$ we write $p_x$, $D_x$ for $p_{\pi(x)}$, $D_{\pi(x)}$, and for an
edge $\tilde e=xy$ of $\tree$ over $e$ we write $A_{xy}$ for the block
$A_e^{\pi(x)\pi(y)}$. The \emph{cover operator} is the bounded self-adjoint
operator on $\ell^2(\tree)=\bigoplus_x\C^{p_x}$,
\[
(\tH\psi)(x)=D_x\psi(x)+\sum_{\tilde e=xy}A_{xy}\psi(y),
\]
and $\Sigm=\inf\sigma(\tH)$, $\Sigp=\sup\sigma(\tH)$. For a set $\Lambda$ of
vertices, $\tH_\Lambda$ is the compression of $\tH$ to
$\bigoplus_{x\in\Lambda}\C^{p_x}$, and for a kernel $K$ on $\tree$,
$K_\Lambda$ is its restriction to $\Lambda\times\Lambda$. Removing an edge
$\tilde e=xy$ from $\tree$ leaves two half-trees. For $z$ outside the
spectrum of the compression of $\tH$ to the half-tree containing $x$, let
$M_{i|e}(z)\in\C^{p_i\times p_i}$ be the block at the root $x$ of the
resolvent of this compression, where $i=\pi(x)$, $e=\pi(\tilde e)$; it
does not depend on the lift. The block Green function is the block
$G_{xy}(z)$ of $(\tH-z)^{-1}$.

\begin{proposition}[Physical branch]\label{prop:block-physical}
Let $G$ be connected with arbitrary real coefficients, and let $z\in\C\setminus\R$
or $z\in\R\setminus[\Sigm,\Sigp]$. For lifts $x$ of $i$ and $\tilde e=xy$ of
$e=ij$ set
\[
\Vphys_i(z)=G_{xx}(z),\qquad C^{\mathrm{phys}}_e(z)=G_{xy}(z).
\]
\begin{enumerate}
\item $(z,\Vphys(z),C^{\mathrm{phys}}(z))$ is an admissible solution of
\eqref{eq:block-system}. Its cavity precisions are
$Q_{i|e}=M_{i|e}(z)^{-1}$.
\item For real $z<\Sigm$ it is a positive solution. The iteration of the
right-hand side of \eqref{blr:eq-cavity} started from $Q_{i|e}=D_i-z$ is
decreasing in the Loewner order and converges to $M_{i|e}(z)^{-1}$, and
every solution $Q'$ of \eqref{blr:eq-cavity} with all $Q'_{i|e}\succ0$
satisfies $Q'_{i|e}\preceq M_{i|e}(z)^{-1}$.
\end{enumerate}
\end{proposition}

\begin{proof}
We use the Schur complement in the following form. Let $T$ be a bounded
operator on $\mathcal H_1\oplus\mathcal H_2$ with blocks $T_{ab}$, and let
$T$ and $T_{22}$ be invertible. Then $T_{11}-T_{12}T_{22}^{-1}T_{21}$ is
invertible and its inverse is the compression of $T^{-1}$ to $\mathcal H_1$.
For the values of $z$ in the statement, $\tH-z$ and all its compressions to
subsets of vertices are invertible: for $z\notin\R$ because compressions of
a self-adjoint operator are self-adjoint, and for real $z<\Sigm$
(respectively $z>\Sigp$) because the compressions of $\tH-z$ are bounded
below by $\Sigm-z>0$ (above by $\Sigp-z<0$). Moreover every diagonal block
of such an inverse is nonsingular: for $z\notin\R$ its imaginary part is
$\im z$ times the corresponding block of $R^*R$, where $R$ is the inverse,
and this block is positive definite because $R$ is injective; for real $z$
the block is positive or negative definite. In particular the $M_{i|e}(z)$
are nonsingular.

Apply the Schur complement to the compression of $\tH-z$ to the half-tree
with root $x$ obtained by removing $\tilde e$, with $\mathcal H_1=\C^{p_x}$.
The complement of the root in this half-tree is the disjoint union of the
half-trees behind the other edges $\tilde e'\ne\tilde e$ at $x$, so
$T_{22}^{-1}$ is the direct sum of their resolvents and
\[
M_{i|e}^{-1}=D_i-z-\sum_{e'\ni i,\ e'\ne e}A_{e'}^{ik}M_{k|e'}A_{e'}^{ki}.
\]
Thus $Q_{i|e}=M_{i|e}^{-1}$ solves \eqref{blr:eq-cavity}. Applying it to
$\tH-z$ with $\mathcal H_1=\C^{p_x}$ gives $G_{xx}^{-1}=R_i$, and with
$\mathcal H_1=\C^{p_x}\oplus\C^{p_y}$ for an edge $\tilde e=xy$ it shows that
$E_e$ is nonsingular and that $E_e^{-1}$ is the compression of
$(\tH-z)^{-1}$ to $\{x,y\}$. By Lemma~\ref{blr:lem-bp} this is the
admissible solution $(\Vphys,C^{\mathrm{phys}})$. For $z<\Sigm$ the
compression of $(\tH-z)^{-1}$ to $\{x,y\}$ is positive definite, which
proves positivity.

Let $z<\Sigm$, $\delta=\Sigm-z$, and let $\mathcal U$ be the map defined
by the right-hand side of \eqref{blr:eq-cavity}. It is monotone on tuples
of positive definite matrices. Let $Q^{(n)}_{i|e}$ be the Schur complement,
onto the root, of the compression of $\tH-z$ to the vertices of the
half-tree at distance at most $n$ from the root. Then $Q^{(0)}_{i|e}=D_i-z$,
$Q^{(n+1)}=\mathcal U(Q^{(n)})$ by eliminating the subtrees behind the
neighbors of the root, and $Q^{(n)}\succeq\delta$ because a Schur
complement of an operator bounded below by $\delta$ is bounded below by
$\delta$. Since $Q^{(1)}\preceq Q^{(0)}$, the sequence is decreasing and
converges to a solution of \eqref{blr:eq-cavity}. On the other hand
$(Q^{(n)}_{i|e})^{-1}$ is the root block of the inverse of the truncated
compression, and for a positive definite bounded operator $T$ and a
finitely supported $u$,
$\langle u,T_\Lambda^{-1}u\rangle=\sup\{2\langle u,f\rangle-\langle f,Tf\rangle:
\operatorname{supp}f\subseteq\Lambda\}$ increases to
$\langle u,T^{-1}u\rangle$ as $\Lambda$ exhausts the vertex set. Hence
$(Q^{(n)}_{i|e})^{-1}\to M_{i|e}(z)$. Finally, a solution $Q'$ with all
$Q'_{i|e}\succ0$ satisfies $Q'=\mathcal U(Q')\preceq D-z=Q^{(0)}$, hence
$Q'=\mathcal U^n(Q')\preceq Q^{(n)}$ for all $n$.
\end{proof}

We call $(z,\Vphys(z),C^{\mathrm{phys}}(z))$ the \emph{physical branch}.
For $\im z>0$ its cavity precisions are characterized among the solutions
of \eqref{blr:eq-cavity} by $\im Q_{i|e}^{-1}\succ0$: equations of this form
with a positivity constraint have exactly one solution in the upper
half-plane \cite{HRFS,BMS}. The remark in \cite{HRFS} that for fixed $z$
the equation has many other solutions is the starting point of
Section~\ref{sec:block-complex}. For real $z$ below the spectrum the
equation has the opposite sign and is outside the scope of these
results; the rest of this section counts its positive solutions.

\subsection{The cone condition}
\label{blr:subsec-cone}

For a family $X=(X_i)_{i\in V}$ of real symmetric positive semidefinite
matrices $X_i\in\Sym_{p_i}(\R)$ put
\begin{equation}\label{blr:eq-q}
q(X)=\sum_{i\in V}\tr(X_iD_iX_i)-2\sum_{e=ij\in E}\nuc{X_iA_eX_j},
\end{equation}
where $\nuc{\cdot}$ is the nuclear norm, the sum of the singular values.

\begin{definition}\label{def:cone}
The coefficients $(D,A)$ satisfy the \emph{cone condition} if $q(X)>0$ for
every nonzero family $X$ of positive semidefinite matrices. We set
\[
\lamC=\lamC(D,A)=\min\Bigl\{q(X):X_i\succeq0,\ \sum_i\tr X_i^2=1\Bigr\},
\]
so that the cone condition is $\lamC>0$. The coefficients lie on the
\emph{boundary of the cone} if $\lamC=0$.
\end{definition}

Since $\lamC(D-z,A)=\lamC(D,A)-z$, the cone condition for the diagonal
$D-z$ is $z<\lamC(D,A)$. If $p_i=1$ for all $i$, then
$q(x)=x^{\T}(D-\abs A)x$ for $x\in\R_{\ge0}^V$, and since an eigenvector of
the smallest eigenvalue of $D-\abs A$ can be taken nonnegative,
$\lamC=\lambda_{\min}(D-\abs A)$; for connected $G$ this is the Perron
value $\sigm$, and the cone condition is walk-summability. In general
$\lamC$ plays the role of $\sigm$, and orthogonal local systems play the
role of signings:

\begin{lemma}\label{blr:lem-cone}
Let $\ell(y)=\sum_i\tr V_i$.
\begin{enumerate}
\item For $y=(V,C)\in\overline\cC$ we have $2L(y)\ge q(V^{1/2})$, where
$V^{1/2}=(V_i^{1/2})_i$, and for every $V$ equality holds for some $C$.
Hence $\lamC=\min\{2L(y):y\in\overline\cC,\ \ell(y)=1\}$.
\item Let $G$ be connected. Then
$\langle f,\tH f\rangle\ge\lamC\norm f^2$ for every finitely supported $f$
on $\tree$, so $\Sigm\ge\lamC$. The same bound holds on every finite cover
of $G$; in particular the precision matrix of the finite model, with
blocks $D_i$ and $\sum_{e=ij}A_e$, is bounded below by $\lamC$.
\item Let $r=\max_ip_i$. For $O=(O_e)\in O(r)^E$ let $H_O$ be the symmetric
operator on $\bigoplus_i\R^{p_i\times r}$ with quadratic form
$\langle Y,H_OY\rangle=\sum_i\tr(Y_i^{\T}D_iY_i)+2\sum_{e=ij}\tr(Y_i^{\T}A_eY_jO_e)$.
Then $\lamC=\min_{O}\lambda_{\min}(H_O)$.
\item If $G$ is connected with cycle rank $\beta\le1$, then $\Sigm=\lamC$.
\end{enumerate}
\end{lemma}

\begin{proof}
(1) Let $X_i=V_i^{1/2}$. A block matrix
$\bigl(\begin{smallmatrix}X_i^2&C\\ C^{\T}&X_j^2\end{smallmatrix}\bigr)$ is
positive semidefinite if and only if $C=X_iKX_j$ for some $K$ with
$\op K\le1$. Then $\tr(A_e^{\T}C_e)=\tr(X_jA_e^{\T}X_iK_e)\ge-\nuc{X_iA_eX_j}$
by the duality between the operator norm and the nuclear norm, with
equality for a suitable $K_e$.

(2) Since $\tH$ is real, $\langle f,\tH f\rangle$ is the sum of the
quadratic forms of the real and imaginary parts of $f$, so it suffices to
take $f$ real. Let $X_i^2=\sum_{x\in\pi^{-1}(i)}f_xf_x^{\T}$. The lifts of an edge
$e=ij$ form a perfect matching between $\pi^{-1}(i)$ and $\pi^{-1}(j)$; let
$C_e=\sum f_xf_y^{\T}$, the sum over the lifts $\tilde e=xy$ of $e$. Then
$S_e=\sum\binom{f_x}{f_y}\binom{f_x}{f_y}^{\T}\succeq0$, so
$y=(X^2,C)\in\overline\cC$, and
$\langle f,\tH f\rangle=2L(y)\ge q(X)\ge\lamC\sum_i\tr X_i^2=\lamC\norm f^2$
by (1). The same computation applies to a finite cover and to $G$ itself.

(3) Let $\sum_i\tr(Y_i^{\T}Y_i)=1$, $V_i=Y_iY_i^{\T}$ and
$C_e=Y_iO_e^{\T}Y_j^{\T}$. Then $S_e=NN^{\T}$ with
$N=\binom{Y_i}{Y_jO_e}$, so $y=(V,C)\in\overline\cC$, $\ell(y)=1$ and
$\langle Y,H_OY\rangle=2L(y)\ge\lamC$. Conversely let $X$ attain $\lamC$,
let $Y_i=(X_i\ 0)\in\R^{p_i\times r}$, and choose $O_e\in O(r)$ from a
singular value decomposition of $Y_i^{\T}A_eY_j$ so that
$\tr(Y_i^{\T}A_eY_jO_e)=-\nuc{Y_i^{\T}A_eY_j}=-\nuc{X_iA_eX_j}$. Then
$\langle Y,H_OY\rangle=q(X)=\lamC$.

(4) By (2) it remains to show $\Sigm\le\lamC$. Take $Y$ and $O$ as at the
end of the proof of (3). On the tree $\tree$ there are $\Gamma_x\in O(r)$
with $\Gamma_y=O_e\Gamma_x$ for every edge $\tilde e=xy$ over $e=ij$ with
$\pi(x)=i$. Put $f_x=Y_{\pi(x)}\Gamma_x\in\R^{p_x\times r}$; its columns are
vector fields on $\tree$. Then $\tr(f_x^{\T}D_xf_x)=\tr(X_iD_iX_i)$ and
$2\tr(f_x^{\T}A_{xy}f_y)=2\tr(Y_i^{\T}A_eY_jO_e)=-2\nuc{X_iA_eX_j}$. If $G$
is a tree then $\tree=G$ and
$\sum_c\langle f^c,\tH f^c\rangle=q(X)=\lamC=\lamC\sum_c\norm{f^c}^2$,
the sums being over the columns, so some column has Rayleigh quotient at
most $\lamC$. If $\beta=1$, then $\tree$ is a chain of copies of a
fundamental domain, a lift of a spanning tree of $G$, consecutive copies
being joined by one edge. Restricting $f$ to $N$ consecutive copies gives
a field $f^{(N)}$ with $\sum_c\norm{f^{(N),c}}^2=N$ and
$\sum_c\langle f^{(N),c},\tH f^{(N),c}\rangle=N\lamC+2\nuc{X_iA_{e_0}X_j}$,
where $e_0$ is the edge outside the spanning tree. Letting $N\to\infty$
gives $\Sigm\le\lamC$.
\end{proof}

Part (4) fails for $\beta\ge2$, where the boundary of $N$ copies of a
fundamental domain is proportional to its volume; see
Corollary~\ref{blr:cor-cycle-rank}.

\subsection{The cover kernel of a positive solution and uniqueness in the cone}
\label{blr:subsec-kernel}

A positive solution determines a Gaussian Markov field on $\tree$ whose
vertex and edge marginals are the lifted beliefs. We need its covariance
kernel, which we construct directly.

\begin{lemma}\label{blr:lem-kernel}
Let $G$ be connected and let $(V,C)$ be a positive solution of $(D,A)$.
Write $V_x$, $C_{xy}$, $S_{xy}$ for the lifts of $V_i$, $C_e^{ij}$, $S_e$ to
$\tree$. There is a kernel $K(x,y)\in\R^{p_x\times p_y}$ on $\tree$,
symmetric and invariant under deck transformations, with the following
properties.
\begin{enumerate}
\item $K(x,x)=V_x$, $K(x,y)=C_{xy}$ for adjacent $x,y$, and if
$x=x_0,x_1,\dots,x_d=y$ is the path from $x$ to $y$, then
$K(x,y)=V_{x_0}^{1/2}\widehat C_{x_0x_1}\cdots\widehat C_{x_{d-1}x_d}V_{x_d}^{1/2}$
with $\widehat C_{xy}=V_x^{-1/2}C_{xy}V_y^{-1/2}$, $\op{\widehat C_{xy}}<1$.
\item For every finite subtree $F$ of $\tree$, $K_F$ is positive definite,
and $K_F^{-1}=J_F$, where $J_F$ is the sum of the matrices $S_{xy}^{-1}$
over the edges of $F$ minus the sum of $(d_F(x)-1)V_x^{-1}$ over the
vertices of $F$, each placed in the corresponding blocks.
\item $\sum_y\tH(x,y)K(y,w)=\delta_{xw}I$ for all $x,w$; for every finite
set $\Lambda$ of vertices, $\tH_\Lambda\succ0$ and
$\tH_\Lambda^{-1}\preceq K_\Lambda$. In particular $\tH\succeq0$.
\end{enumerate}
\end{lemma}

\begin{proof}
We define $K_F=J_F^{-1}$ and prove (2) and the consistency
$K_{F'}|_{F\times F}=K_F$ for $F\subseteq F'$ by induction on the number of
vertices. For a single vertex, $J_F=V_x^{-1}$. Let $F'=F\cup\{y\}$ with $y$
a leaf attached to $x\in F$ by an edge over $e$, and let
$Q_x=Q_{\pi(x)|e}$, $Q_y=Q_{\pi(y)|e}$, $A=A_{xy}$. Since
$S_{xy}^{-1}=\bigl(\begin{smallmatrix}Q_x&A\\ A^{\T}&Q_y\end{smallmatrix}\bigr)$
and $Q_x-V_x^{-1}=AQ_y^{-1}A^{\T}$ by \eqref{blr:eq-schur},
\[
J_{F'}=\begin{pmatrix}J_F+\epsilon_xAQ_y^{-1}A^{\T}\epsilon_x^{\T}&\epsilon_xA\\
A^{\T}\epsilon_x^{\T}&Q_y\end{pmatrix},
\]
where $\epsilon_x$ is the inclusion of the block of $x$. The Schur
complement of $Q_y\succ0$ is $J_F$, so $J_{F'}\succ0$, the block of
$J_{F'}^{-1}$ on $F$ is $J_F^{-1}=K_F$, and
$K_{F'}(w,y)=-K_F(w,x)AQ_y^{-1}=K(w,x)V_x^{-1}C_{xy}$ for $w\in F$, because
$S_{xy}S_{xy}^{-1}=I$ gives $V_xA+C_{xy}Q_y=0$. Likewise
$K_{F'}(y,y)=Q_y^{-1}+C_{xy}^{\T}V_x^{-1}C_{xy}=V_y$. Every finite subtree is
obtained by adding leaves, and $J_F$ does not depend on the order, so $K$
is well defined; it is invariant because the construction is. The
recursion gives (1), and $\op{\widehat C_{xy}}<1$ because $S_{xy}\succ0$.

(3) The $(x,y)$ block of $J_F$ is $A_{xy}$ for adjacent $x,y\in F$ and zero
for other $x\ne y$, and its $(x,x)$ block is
$\sum Q_{\pi(x)|e}-(d_F(x)-1)V_x^{-1}$, the sum over the edges of $F$ at
$x$. If all neighbors of $x$ in $\tree$ belong to $F$, this is $D_x$ by
\eqref{eq:block-system}. Given a finite $\Lambda$ and a vertex $w$, let $F$ be a finite
subtree containing $w$, $\Lambda$ and all neighbors of the vertices of
$\Lambda$. Then the
rows of $J_F$ indexed by $\Lambda$ are the rows of $\tH$, so $J_FK_F=I$
gives the first claim, and $\tH_\Lambda=(J_F)_\Lambda\succ0$ with
$(J_F)_\Lambda^{-1}\preceq(J_F^{-1})_\Lambda=K_\Lambda$, since the inverse of
a principal submatrix of a positive definite matrix is dominated by the
corresponding submatrix of the inverse.
\end{proof}

\begin{lemma}\label{blr:lem-maximal}
In the situation of Lemma~\ref{blr:lem-kernel}:
\begin{enumerate}
\item The limit $K_0=\lim_\Lambda\tH_\Lambda^{-1}$ along any exhaustion of
$\tree$ by finite sets exists entrywise, is invariant under deck
transformations, and satisfies $0\preceq K_0\preceq K$ on every finite set
and $\sum_y\tH(x,y)K_0(y,w)=\delta_{xw}I$.
\item $V^0_i=K_0(x,x)$ and $C^0_e=K_0(x,y)$, for lifts $x$ of $i$ and $xy$
of $e=ij$, form a positive solution of $(D,A)$. Its cavity precisions are
the limit of the iteration of the right-hand side of \eqref{blr:eq-cavity}
(with $z=0$) started from $Q_{i|e}=D_i$, and they dominate the cavity
precisions of every positive solution. We call it the \emph{maximal
positive solution}. If $\Sigm>0$ it is the physical branch at $z=0$.
\item The kernel $H=K-K_0$ is positive semidefinite on every finite set
and $\sum_y\tH(x,y)H(y,w)=0$. With $H_i=H(x,x)$ and $H_e=H(x,y)$ for
lifts as above,
\begin{equation}\label{blr:eq-energy}
\begin{pmatrix}H_i&H_e\\ H_e^{\T}&H_j\end{pmatrix}\succeq0\quad(e=ij),\qquad
\sum_i\tr(D_iH_i)+2\sum_e\tr(A_e^{\T}H_e)=0.
\end{equation}
If $H_i=0$ for all $i$, then $(V,C)=(V^0,C^0)$.
\end{enumerate}
\end{lemma}

\begin{proof}
(1) For finitely supported $u$, the quantity
$\langle u,\tH_\Lambda^{-1}u\rangle=\sup\{2\langle u,f\rangle-\langle f,\tH f\rangle:
\operatorname{supp}f\subseteq\Lambda\}$ is increasing in $\Lambda$ and
bounded by $\langle u,Ku\rangle$ by Lemma~\ref{blr:lem-kernel}(3). By
polarisation the limit kernel exists and does not depend on the
exhaustion. The identity $\tH_\Lambda\tH_\Lambda^{-1}=I$, read at a row $x$
all of whose neighbors lie in $\Lambda$, passes to the limit because the
rows of $\tH$ are finite.

(2) Let $Q$ be the cavity precisions of $(V,C)$ and $\mathcal U$ the
monotone map of the proof of Proposition~\ref{prop:block-physical}, with
$z=0$. Then $0\prec Q=\mathcal U(Q)\preceq D$, so the iterates
$Q^{(n)}=\mathcal U^n(D)$ decrease and satisfy $Q^{(n)}\succeq Q$; their
limit $Q^0\succeq Q$ solves \eqref{blr:eq-cavity}. Its vertex and edge
precisions dominate those of $Q$, so they are positive definite and $Q^0$
defines a positive solution by Lemma~\ref{blr:lem-bp}; the same argument
with any other positive solution in place of $(V,C)$ shows that $Q^0$
dominates all of them. As in the proof of
Proposition~\ref{prop:block-physical}, $Q^{(n)}_{i|e}$ is the Schur
complement onto the root of the compression of $\tH$ to the vertices of a
half-tree within distance $n$ of the root; these compressions are positive
definite by Lemma~\ref{blr:lem-kernel}(3). Hence, for the ball $B_n(x)$,
the $(x,x)$ block of $\tH_{B_n(x)}^{-1}$ is the inverse of
$D_i-\sum_{e\ni i}A_e^{ik}(Q^{(n-1)}_{k|e})^{-1}A_e^{ki}$ and converges to
$R_i(Q^0)^{-1}$; and for an edge $xy$ over $e$ and the union $\Lambda_n$ of
the two truncated half-trees with roots $x$ and $y$, the compression of
$\tH_{\Lambda_n}^{-1}$ to $\{x,y\}$ is the inverse of the matrix
\eqref{blr:eq-RE} formed with $Q^{(n)}$ and converges to $E_e(Q^0)^{-1}$. By
(1) these limits are $K_0(x,x)$ and the compression of $K_0$ to $\{x,y\}$.
If $\Sigm>0$, then $K_0$ is the kernel of $\tH^{-1}$ by the variational
formula.

(3) The first two statements follow from (1) and
Lemma~\ref{blr:lem-kernel}(3), and the matrix in \eqref{blr:eq-energy} is
a compression of $H$. Taking the trace of the $(x,x)$ block of
$\sum_y\tH(x,y)H(y,x)=0$ and summing over one lift $x$ of each $i\in V$
gives the identity in \eqref{blr:eq-energy}; each edge contributes
$\tr(A_e^{\T}H_e)$ from each of its ends. If all $H_i$ vanish, then $H=0$,
because a positive semidefinite kernel with zero diagonal vanishes, and
the two solutions have the same vertex and edge covariances.
\end{proof}

By Lemma~\ref{blr:lem-kernel}(3) a positive solution with spectral
parameter $z$ can exist only if $\tH-z\succeq0$, that is, for
$z\le\Sigm$, and by Proposition~\ref{prop:block-physical} one exists for
every $z<\Sigm$.

\begin{theorem}[Uniqueness in the cone]\label{thm:block-cone-unique}
Let $G$ be a finite loopless multigraph with vertex dimensions $p_i\ge1$,
and let $D_i$ be real symmetric and $A_e$ arbitrary real matrices. If
$(D,A)$ satisfies the cone condition, then there is exactly one positive
solution. On each connected component it is the physical branch at $z=0$,
and the precision iteration started from $Q_{i|e}=D_i$ converges to its
cavity precisions. Equivalently, for arbitrary real coefficients and every
real $z<\lamC(D,A)$ the block covariance system has exactly one positive
solution with this~$z$.
\end{theorem}

\begin{proof}
The quadratic form $q$ and the positive solutions decompose over the
connected components, so let $G$ be connected. By
Lemma~\ref{blr:lem-cone}(2), $\Sigm\ge\lamC>0$, so the physical branch at
$z=0$ is a positive solution with the stated properties by
Proposition~\ref{prop:block-physical}. Let $(V,C)$ be any positive
solution and $H$ the kernel of Lemma~\ref{blr:lem-maximal}(3). By
\eqref{blr:eq-energy} the family $(H_i,H_e)$ lies in $\overline\cC$ and
$L$ vanishes on it. The minimum of $2L$ among edge blocks with diagonal
blocks $H_i,H_j$ is $q((H_i^{1/2})_{i\in V})$ by
Lemma~\ref{blr:lem-cone}(1). Thus
$0\ge q((H_i^{1/2})_{i\in V})\ge\lamC\sum_i\tr H_i$.
Hence all $H_i$ vanish and $(V,C)$
is the maximal positive solution, which is the physical branch.
\end{proof}

For $p_i=1$ this is the implication from walk-summability to uniqueness
in Theorem~\ref{thm:unique-ws}. For individual coefficients the cone
condition is not necessary for uniqueness: the model of
Example~\ref{blr:ex-K4-boundary} lies on the boundary of the cone and has
exactly one positive solution. For generic coefficients it is necessary,
as we show next.

\subsection{The signed count and the converse}
\label{blr:subsec-degree}

\begin{theorem}[Signed count]\label{thm:block-converse}
Let $G$ be a finite loopless multigraph with vertex dimensions $p_i\ge1$.
\begin{enumerate}
\item For all real coefficients satisfying the cone condition,
the Brouwer degree of $\nabla F$ at $0$ on any bounded open set
$\Omega$ with $\overline\Omega\subset\cC$ that contains the positive
solution is $1$.
\item There is a proper real algebraic subset $Z\subset\cP_\R$ such that
for all $(D,A)\in\cP_\R\setminus Z$ the positive solutions are finite in
number and nondegenerate, $\lamC(D,A)\ne0$, and
\[
\sum_{y}(-1)^{\ind(y)}=\begin{cases}1&\text{if }\lamC(D,A)>0,\\
0&\text{if }\lamC(D,A)<0,\end{cases}
\]
the sum being over the positive solutions $y$ of $(D,A)$ and $\ind(y)$ the
index.
\end{enumerate}
Consequently, for generic real coefficients there is exactly one positive
solution if and only if the cone condition holds; if it fails, the number
of positive solutions is even, and if there is one there is also one of
odd index.
\end{theorem}

For $p_i=1$ this gives, for generic coefficients only,
Theorem~\ref{thm:degree-identity}, Corollary~\ref{cor:parity} and the
converse direction of Theorem~\ref{thm:unique-ws}. The proof deforms
$F$ to a strictly convex barrier while keeping the edge covariances as
variables. For $s\in[0,1]$ let
\[
F_s=L+B_s,\qquad B_s=-\tfrac12\sum_e\log\det S_e+\tfrac12\sum_i(sd_i-1)\log\det V_i
\]
on $\cC$, so that $F_1=F$. Let $K=\{v\in\overline\cC:\ell(v)=1\}$. It is
compact, because for $y\in\overline\cC$ every entry of $C_e$ is bounded by
the geometric mean of two diagonal entries of $V_i$ and $V_j$. Every
$y\in\cC$ is $y=tv$ with $t=\ell(y)>0$ and $v\in K\cap\cC$, and
\begin{equation}\label{blr:eq-scaling}
F_s(tv)=tL(v)-b_s\log t+B_s(v),\qquad b_s=\tfrac12\sum_i\bigl(1+(1-s)d_i\bigr)p_i>0,
\end{equation}
because $\sum_{e=ij}(p_i+p_j)=\sum_id_ip_i$. Hence every critical point
$tv$ of $F_s$ satisfies the radial identity $tL(v)=b_s$. Moreover, for
$v\in K\cap\cC$,
\begin{equation}\label{blr:eq-barrier}
B_s(v)=-\tfrac12\sum_{e=ij}\log\frac{\det S_e}{\det V_i\det V_j}
-\tfrac12\sum_i\bigl(1+(1-s)d_i\bigr)\log\det V_i
\ \ge\ B_1(v)\ \ge0 ,
\end{equation}
since each ratio lies in $(0,1]$ and $\det V_i\le1$ on $K$. When $v$ tends
to a point of $K\setminus\cC$, either some $\det V_i$ tends to $0$, or all
$V_i$ stay positive definite and some $\det S_e$ tends to $0$; in both
cases $B_1(v)\to+\infty$, hence $B_s(v)\to+\infty$ uniformly in $s$.

\begin{lemma}\label{blr:lem-compact}
The set of $(s,y)\in[0,1]\times\cC$ with $\nabla F_s(y)=0$ is compact in
each of the following cases: (a) $(D,A)$ satisfies the cone condition;
(b) $(D,A)$ lies outside a certain proper real algebraic subset
$Z_0\subset\cP_\R$.
\end{lemma}

\begin{proof}
Let $m$ be the dimension of the space of $y$ and regard the directions
$v$ as points of the complex projective space $\PP^{m-1}$. Let $\Delta$ be
the union of the hypersurfaces $\{\ell=0\}$, $\{\det V_i=0\}$ and
$\{\det S_e=0\}$, and let $\rho\colon\widetilde P\to\PP^{m-1}$ be a log
resolution of $\Delta$: a proper morphism from a smooth variety that is an
isomorphism over the complement of $\Delta$ and such that
$\rho^{-1}(\Delta)$ is a divisor with simple normal crossings
\cite{Wlodarczyk}. It depends only on $(G,p)$. The strata of
$\rho^{-1}(\Delta)$, including the open stratum
$\widetilde P\setminus\rho^{-1}(\Delta)$, are finitely many smooth locally
closed subvarieties. By Bertini's theorem for the base-point-free linear
system $\rho^*\abs{\mathcal O(1)}$ restricted to each stratum
\cite[Corollary~III.10.9 and Remark~III.10.9.2]{Hartshorne}, there is a nonempty Zariski open set
$\mathcal O$ of linear forms $L$ such that, on every stratum $S$, a local
equation of $\rho^{-1}\{L=0\}$ restricted to $S$ has nonzero differential
at every point of $S\cap\rho^{-1}\{L=0\}$ (for a stratum of dimension
zero: $\rho^{-1}\{L=0\}$ does not contain it). We call such $L$
transverse, and let $Z_0$ be the set of real $(D,A)$ whose form $L$ does
not belong to $\mathcal O$; it is a proper real algebraic subset of
$\cP_\R$.

On the chart $\ell\ne0$ let $\lambda=L/\ell$,
$g_i=\det V_i/\ell^{p_i}$ and $g_e=\det S_e/\ell^{p_i+p_j}$; on the real
slice $\ell=1$ these are $L$, $\det V_i$ and $\det S_e$. Let
$x_*\in\widetilde P$ with $\rho(x_*)$ in this chart, and let
$u_1,\dots,u_{m-1}$ be local holomorphic coordinates at $x_*$ in which
$\rho^{-1}(\Delta)$ is $\{u_1\cdots u_k=0\}$, $k\ge0$. Each $\rho^*g_i$,
$\rho^*g_e$ is holomorphic near $x_*$ with zero set contained in
$\rho^{-1}(\Delta)$. The local ring at $x_*$ is a unique factorisation
domain in which the $u_j$ are prime, so each of these functions is a unit
times a monomial in $u_1,\dots,u_k$. Therefore, with
$\beta_s=-\tfrac12\sum_ed\log g_e+\tfrac12\sum_i(sd_i-1)d\log g_i$,
\begin{equation}\label{blr:eq-logform}
\rho^*\beta_s=\sum_{j=1}^kc_j(s)\frac{du_j}{u_j}+\eta_s,\qquad
B_s\circ\rho=\sum_{j=1}^kc_j(s)\log\abs{u_j}+h_s ,
\end{equation}
where the $c_j$ are real affine functions of $s$, $\eta_s$ is a
holomorphic $1$-form and $h_s$ a function, both bounded near $x_*$
uniformly in $s\in[0,1]$; the second identity holds at the points over
$K\cap\cC$, where $B_s$ is the corresponding combination of the
$\log\abs{g}$.

Let $(s_n,t_nv_n)$ be a sequence of critical points. By
\eqref{blr:eq-scaling}, on the slice $\ell=1$ the tangential part of the
critical equation is $t_n\,d\lambda+\beta_{s_n}=0$ at $v_n$; a holomorphic
covector that vanishes on the real tangent space vanishes, so this holds
for the holomorphic forms on the chart. After passing to a subsequence,
$s_n\to s_*$, $v_n\to v_*\in K$, and $x_n=\rho^{-1}(v_n)\to x_*$ with
$\rho(x_*)=v_*$, by properness. The radial identity gives $L(v_n)>0$ and
$t_n\ge b_1/\max_KL>0$, since $b_s$ is decreasing in $s$. In the
coordinates at $x_*$,
\begin{equation}\label{blr:eq-crit-coord}
t_n\,\partial_{u_j}(\lambda\circ\rho)(x_n)+\frac{c_j(s_n)}{u_j(x_n)}\,[j\le k]
+\eta_{s_n,j}(x_n)=0\qquad(1\le j\le m-1).
\end{equation}

Suppose first that $t_n$ is bounded and $v_*\notin\cC$, so that $k\ge1$.
Multiplying \eqref{blr:eq-crit-coord} for $j\le k$ by $u_j(x_n)\to0$
gives $c_j(s_*)=0$ for all $j\le k$. By \eqref{blr:eq-logform},
$B_{s_*}(v_n)$ is then bounded, which contradicts the divergence of $B_s$
at $K\setminus\cC$. Suppose next that $t_n\to\infty$. Then
$L(v_n)=b_{s_n}/t_n\to0$, so $\lambda(\rho(x_*))=0$, and dividing
\eqref{blr:eq-crit-coord} for $j>k$ by $t_n$ shows that
$d(\lambda\circ\rho)$ vanishes at $x_*$ on the tangent space of the
stratum $\{u_1=\dots=u_k=0\}$; this is excluded when $L$ is transverse.
Under the cone condition the second case does not occur at all: by
Lemma~\ref{blr:lem-cone}(1), $L\ge\tfrac12\lamC$ on $K$, so
$t_n\le2b_0/\lamC$. In the remaining case $t_n$ is bounded above and below
and $v_*\in\cC$, so the limit is a critical point in $\cC$.
\end{proof}

\begin{proof}[Proof of Theorem~\ref{thm:block-converse}]
In both cases of Lemma~\ref{blr:lem-compact} choose a bounded open
$\Omega$ with $\overline\Omega\subset\cC$ containing all critical points
of all $F_s$. By homotopy invariance,
$\deg(\nabla F_1,\Omega,0)=\deg(\nabla F_0,\Omega,0)$.

The function
$F_0=L-\tfrac12\sum_e\log\det S_e-\tfrac12\sum_i\log\det V_i$ has a positive
definite Hessian everywhere: the second variation of $-\log\det M$ in the
direction $\delta M$ is $\norm{M^{-1/2}\delta M\,M^{-1/2}}_F^2$, and if
these vanish for all $S_e$ and $V_i$ then $\delta y=0$. If $\lamC>0$,
then by \eqref{blr:eq-scaling}, \eqref{blr:eq-barrier} and
Lemma~\ref{blr:lem-cone}(1),
$F_0(tv)\ge\tfrac12\lamC t-b_0\log t+B_0(v)$, so $F_0$ tends to $+\infty$
at the boundary of $\cC$ and at infinity. It has a unique critical point,
a nondegenerate minimum, and the degree is $1$. This proves (1), the
degree being independent of $\Omega$ by excision and
Theorem~\ref{thm:block-cone-unique}. If $\lamC<0$, there is
$w\in\overline\cC$ with $L(w)<0$ by Lemma~\ref{blr:lem-cone}(1), and after a
small perturbation $w\in\cC$. For every $y\in\cC$,
\[
dF_0(y)[w]=L(w)-\tfrac12\sum_e\tr\bigl(S_e(y)^{-1}S_e(w)\bigr)
-\tfrac12\sum_i\tr\bigl(V_i(y)^{-1}V_i(w)\bigr)<0,
\]
so $F_0$ has no critical point and the degree is $0$.

It remains to define $Z$. The map $\Psi\colon Y_p\to\cP$ is a morphism
between smooth irreducible varieties of the same dimension, and its
differential at $V_i=I$, $C_e=0$ is $-\mathrm{id}$. Hence the set where
$d\Psi$ is singular is a proper closed subset of $Y_p$, and by generic
smoothness in characteristic zero its image is contained in a proper
closed subset $Z_1$ of $\cP$. For real $(D,A)\notin Z_1$ all positive
solutions are nondegenerate by Lemma~\ref{blr:lem-critical}. The set
$\{\lamC=0\}$ is semialgebraic and, since $\lamC(D+c,A)=\lamC(D,A)+c$, it
is the graph of a function over a hyperplane of $\cP_\R$; it is therefore
contained in a proper real algebraic subset $Z_2$. Let
$Z=Z_0\cup(Z_1\cap\cP_\R)\cup Z_2$. For $(D,A)\notin Z$ the positive
solutions form a compact set of nondegenerate critical points, hence a
finite set, the local degree of $\nabla F$ at $y$ is $(-1)^{\ind(y)}$, and
their sum is $\deg(\nabla F_1,\Omega,0)$. The consequences follow, since
the summands are $\pm1$.
\end{proof}

\begin{example}\label{blr:ex-K4}
Let $G=K_4$ with all $p_i=p$, $A_e=I_p$ and $D_i=\diag(a_1,\dots,a_p)$,
where $a_k=r_k+2/r_k$ for distinct rational numbers
$r_k\in(\sqrt2,2)$. For each choice of $s_k\in\{r_k,2/r_k\}$,
$Q_{i|e}=\diag(s_1,\dots,s_p)$ solves \eqref{blr:eq-cavity}, since
$K_4$ is $3$-regular and $a_k-2/s_k=s_k$. The vertex precisions
$\diag(s_k-1/s_k)$ are positive definite because $s_k>1$, and the edge
precisions have the eigenvalues $s_k\pm1>0$. These are $2^p$ distinct
positive solutions. They are nondegenerate: by
Lemma~\ref{blr:lem-critical} and the identity \eqref{blr:eq-cavity}, a
kernel vector of $d\Psi$ gives a variation $\delta Q$ fixed by the
linearisation of the right-hand side of \eqref{blr:eq-cavity}, and
$\delta Q=0$ forces the kernel vector to vanish; on the $(k,l)$ entries
this linearisation is $(s_ks_l)^{-1}$ times the nonbacktracking matrix of
$K_4$, a nonnegative integer matrix with row sums $2$, whose rational
eigenvalues are integers of absolute value at most $2$, while
$s_ks_l>1$ is rational and different from $2$ (it is $r_kr_l>2$,
$4/(r_kr_l)<2$, or $2r_k/r_l$ with $k\ne l$). The cone condition fails,
because $X_i=e_ke_k^{\T}$ for all $i$ gives $q(X)=4a_k-12<0$. By the
implicit function theorem the $2^p$ positive solutions persist, and
$\lamC<0$, on an open set of real coefficients, in which the blocks need
not be simultaneously diagonal; at the generic points of this open set the signed count of
Theorem~\ref{thm:block-converse} is zero. For $p=2$ and
$r=(\tfrac32,\tfrac53)$ we get $D_i=\diag(\tfrac{17}6,\tfrac{43}{15})$ and
the four solutions $V_i=\diag(v_1,v_2)$ with $v_1\in\{\tfrac65,\tfrac{12}7\}$
and $v_2\in\{\tfrac{15}{16},\tfrac{30}{11}\}$. The finite precision matrix
of this model is positive definite.
\end{example}

\subsection{The boundary of the cone}
\label{blr:subsec-boundary}

We now assume only $\lamC\ge0$, that is, $q(X)\ge0$ for all $X$. By
Lemma~\ref{blr:lem-cone}(2) the cover operator is then nonnegative, and
$(\tH+\varepsilon)^{-1}$ is defined for $\varepsilon>0$.

\begin{theorem}[Boundary of the cone]\label{thm:block-boundary}
Let $G$ be a connected finite loopless multigraph with vertex dimensions
$p_i\ge1$ and arbitrary real coefficients with $\lamC(D,A)\ge0$.
\begin{enumerate}
\item There is at most one positive solution.
\item There is a positive solution if and only if
\[
\lim_{\varepsilon\downarrow0}\tr\bigl[(\tH+\varepsilon)^{-1}\bigr]_{xx}<\infty
\qquad\text{for every vertex }x\text{ of }\tree ,
\]
that is, if and only if $\int\lambda^{-1}d\mu_i(\lambda)<\infty$ for all
$i\in V$, where $\mu_i$ is the trace of the spectral measure of $\tH$ at a
lift of $i$. In that case the positive solution is the limit of the
physical branch as $z\uparrow0$.
\end{enumerate}
\end{theorem}

The criterion involves only the physical branch. If $\lamC=0<\Sigm$, the
physical branch at $z=0$ is a positive solution, so there is exactly one
(Example~\ref{blr:ex-K4-boundary}). If $\lamC=\Sigm=0$, the criterion asks
whether the traces of the Green function stay bounded at the bottom of the
spectrum. They do not in Corollary~\ref{blr:cor-cycle-rank} and in
Example~\ref{blr:ex-K4-singular}, where there is no positive solution; we
do not know an example with $\lamC=\Sigm=0$ that has one.
The proof of (1) uses a dual description of the closed cone.

\begin{lemma}\label{blr:lem-certificate}
Let every vertex of $G$ have an edge. Then $\lamC\ge0$ if and only if
there are $U_{i,e}\in\Sym_{p_i}(\R)$, one for each edge $e$ and each end $i$
of $e$, such that
\begin{equation}\label{blr:eq-certificate}
Z_e=\begin{pmatrix}U_{i,e}&A_e\\ A_e^{\T}&U_{j,e}\end{pmatrix}\succeq0\quad(e=ij),
\qquad\sum_{e\ni i}U_{i,e}=D_i\quad(i\in V).
\end{equation}
\end{lemma}

\begin{proof}
If \eqref{blr:eq-certificate} holds, then $2L(y)=\sum_e\tr(Z_eS_e)\ge0$ for
$y\in\overline\cC$, and $\lamC\ge0$ by Lemma~\ref{blr:lem-cone}(1). For the
converse we first show that for positive semidefinite $P_i,P_j$ with
square roots $X_i,X_j$,
\begin{equation}\label{blr:eq-mu}
\inf\Bigl\{\tr(P_iU_i)+\tr(P_jU_j):
\begin{pmatrix}U_i&A_e\\ A_e^{\T}&U_j\end{pmatrix}\succeq0\Bigr\}=2\nuc{X_iA_eX_j}.
\end{equation}
The inequality $\ge$ follows by pairing the matrix with
$\bigl(\begin{smallmatrix}P_i&-X_iKX_j\\ -X_jK^{\T}X_i&P_j\end{smallmatrix}\bigr)\succeq0$
for the contraction $K$ with $\tr(X_jA_e^{\T}X_iK)=\nuc{X_iA_eX_j}$. If
$P_i,P_j$ are nonsingular and $X_iA_eX_j=W\Sigma\Theta^{\T}$ is a singular
value decomposition, then $U_i=X_i^{-1}W\Sigma W^{\T}X_i^{-1}$,
$U_j=X_j^{-1}\Theta\Sigma\Theta^{\T}X_j^{-1}$ is feasible with value
$2\tr\Sigma$. In general apply this to $P+\delta I$; the resulting
$U^\delta\succeq0$ satisfy
$\tr(P_iU_i^\delta)\le\tr((P_i+\delta)U_i^\delta)$, and the values converge
as $\delta\to0$.

Let $\mathcal K$ be the set of families $(\sum_{e\ni i}U_{i,e})_{i\in V}$
with $Z_e\succeq0$ for all $e$. It is convex, and closed because
$Z_e\succeq0$ implies $0\preceq U_{i,e}\preceq\sum_{e'\ni i}U_{i,e'}$; it is
stable under adding positive semidefinite families. If $D\notin\mathcal K$,
there is a family $P$ of symmetric matrices with
$\sum_i\tr(P_iD_i)<\inf_{\kappa\in\mathcal K}\sum_i\tr(P_i\kappa_i)$. The
infimum is finite, so $P_i\succeq0$, and since the constraints are
separate for each edge, the infimum is $2\sum_e\nuc{X_iA_eX_j}$ by
\eqref{blr:eq-mu}. Thus $q(P^{1/2})<0$.
\end{proof}

\begin{proof}[Proof of Theorem~\ref{thm:block-boundary}(1)]
If $G$ is a single vertex, the only possible positive solution is
$V=D^{-1}$. Otherwise let $U$ be as in Lemma~\ref{blr:lem-certificate},
let $(V,C)$ be a positive solution, and let $K$, $K_0$, $H=K-K_0$ be the
kernels of Lemmas~\ref{blr:lem-kernel} and~\ref{blr:lem-maximal}. We show
that all $H_i$ vanish. Write $\widehat H_e$ for the matrix in
\eqref{blr:eq-energy}. By \eqref{blr:eq-certificate} and
\eqref{blr:eq-energy},
$\sum_e\tr(Z_e\widehat H_e)=\sum_i\tr(D_iH_i)+2\sum_e\tr(A_e^{\T}H_e)=0$.
All terms are nonnegative, so $Z_e\widehat H_e=0$ for every edge.

The congruence $D_i\mapsto V_i^{1/2}D_iV_i^{1/2}$,
$A_e\mapsto V_i^{1/2}A_eV_j^{1/2}$, with the inverse congruence on $V$,
$C$, the kernels and the corresponding one on $U$, preserves all
hypotheses and conclusions, so we may assume $V_i=I$ for all $i$. Fix an
edge $e=ij$ and let
\[
T=S_e=\begin{pmatrix}I&C\\ C^{\T}&I\end{pmatrix},\qquad J=S^0_e,\qquad
0\prec J\preceq T,\qquad \op C<1,
\]
where $S^0_e$ is the edge covariance of the maximal positive solution. Both
$T^{-1}$ and $J^{-1}$ are edge precisions \eqref{blr:eq-RE} with the same
off-diagonal block $A_e$. Hence
\[
\Phi:=J^{-1}-T^{-1}=\diag(\Delta_i,\Delta_j)\succeq0,\qquad
\widehat H_e=T-J=T\Phi J .
\]
Let $L_i=\ran\Delta_i$ and $L_j=\ran\Delta_j$; these depend on $e$. From
$Z_eT\Phi J=0$ we get $Z_e\binom{x}{C^{\T}x}=0$ for $x\in L_i$ and
$Z_e\binom{Cy}{y}=0$ for $y\in L_j$. Since
$A_e=(T^{-1})_{ij}=-C(I-C^{\T}C)^{-1}$, the first rows of these identities
say that the symmetric matrix
\[
B_{i,e}=U_{i,e}-C(I-C^{\T}C)^{-1}C^{\T}
\]
satisfies
\begin{equation}\label{blr:eq-B}
B_{i,e}x=0\quad(x\in L_i),\qquad B_{i,e}Cy=Cy\quad(y\in L_j).
\end{equation}
The same holds at $j$ with $C$ replaced by $C^{\T}$. In particular
$L_i\perp CL_j$, these being eigenspaces of $B_{i,e}$ for the eigenvalues
$0$ and $1$. Since $\ran\widehat H_e=T(L_i\oplus L_j)$ and
the projection to the first block of the range of a positive semidefinite
block matrix is the range of its first diagonal block,
\begin{equation}\label{blr:eq-W}
W_i:=\ran H_i=L_i\oplus CL_j ,
\end{equation}
an orthogonal decomposition, valid for every edge $e$ at $i$, of a space
that depends only on $i$. By \eqref{blr:eq-B}, $B_{i,e}$ preserves $W_i$,
and its restriction $P_{i,e}$ to $W_i$ is the orthogonal projection of
$W_i$ onto $CL_j$. The second equation of \eqref{eq:block-system} with
$V_i=I$ reads $D_i=I+\sum_{e\ni i}\bigl((I-C_eC_e^{\T})^{-1}-I\bigr)$, and
$(I-CC^{\T})^{-1}-I=C(I-C^{\T}C)^{-1}C^{\T}$, so
$\sum_{e\ni i}B_{i,e}=I$ by \eqref{blr:eq-certificate}. Therefore
$\sum_{e\ni i}P_{i,e}=I_{W_i}$, and orthogonal projections that sum to the
identity have mutually orthogonal ranges: $W_i$ is the orthogonal direct
sum of the spaces $C_eL^{(e)}_j$ over the edges $e=ij$ at $i$.

Since $L_i\perp CL_j$, we have $\Delta_i^{1/2}C\Delta_j^{1/2}=0$, so
$\Phi^{1/2}T\Phi^{1/2}=\Phi$ and
\[
T-J=T\Phi^{1/2}\bigl(I+\Phi^{1/2}T\Phi^{1/2}\bigr)^{-1}\Phi^{1/2}T
=T\diag(\Xi_i,\Xi_j)\,T,\qquad \Xi_i=\Delta_i(I+\Delta_i)^{-1}.
\]
Taking diagonal blocks, $H_i=\Xi_i+C\Xi_jC^{\T}$ and
$H_j=\Xi_j+C^{\T}\Xi_iC$. The two terms of $H_i$ have their ranges in
$L_i$ and $CL_j$, so $H_i$ commutes with $P_{i,e}$, and on the range of
$P_{i,e}$ it equals $C\Xi_jC^{\T}$, of norm at most
$\op{C_e}^2\op{H_j}$ because $\Xi_j\preceq H_j$. As $W_i$ is the
orthogonal sum of these ranges, $\op{H_i}\le c^2\max_j\op{H_j}$ with
$c=\max_e\op{C_e}<1$. Hence all $H_i$ vanish, and by
Lemma~\ref{blr:lem-maximal}(3) the solution $(V,C)$ is the maximal one.
\end{proof}

\begin{proof}[Proof of Theorem~\ref{thm:block-boundary}(2)]
Let $(V,C)$ be a positive solution with kernel $K$. For finitely
supported $u$ and $\varepsilon>0$,
\[
\langle u,(\tH+\varepsilon)^{-1}u\rangle
=\sup_f\bigl\{2\langle u,f\rangle-\langle f,(\tH+\varepsilon)f\rangle\bigr\}
\le\sup_\Lambda\langle u,\tH_\Lambda^{-1}u\rangle\le\langle u,Ku\rangle
\]
by Lemma~\ref{blr:lem-kernel}(3), the first supremum being over finitely
supported $f$. So the limits are finite. Conversely, assume that they are
finite. For $0<\varepsilon\le1$ the physical branch at $z=-\varepsilon$
is a positive solution $(V^\varepsilon,C^\varepsilon)$ of
$(D+\varepsilon,A)$ by Proposition~\ref{prop:block-physical}, and
$S^\varepsilon_e$ is a compression of $(\tH+\varepsilon)^{-1}$. These
compressions increase as $\varepsilon$ decreases, their diagonal blocks
have bounded traces, and they are bounded below by $(\norm{\tH}+1)^{-1}$.
Hence $(V^\varepsilon,C^\varepsilon)$ converges to a point $(V,C)\in\cC$,
and $\Psi(V,C)=\lim\Psi(V^\varepsilon,C^\varepsilon)=(D,A)$ by continuity
of $\Psi$ on $Y_p$. By (1) this limit is the positive solution.
\end{proof}

\begin{corollary}\label{blr:cor-cycle-rank}
Let $G$ be connected with cycle rank $\beta$ and arbitrary real
coefficients.
\begin{enumerate}
\item If $\beta\le1$, then $\Sigm=\lamC$, and a positive solution exists
if and only if the cone condition holds. In particular there is none on
the boundary of the cone, for all vertex dimensions and all real edge
blocks.
\item If $\beta\ge2$ and all edge blocks $A_e$ are invertible, then
$\lamC<\Sigm$. In particular on the boundary of the cone there is exactly
one positive solution.
\end{enumerate}
\end{corollary}

\begin{proof}
(1) The equality is Lemma~\ref{blr:lem-cone}(4). Let $(V,C)$ be a positive
solution with kernel $K$. By Lemma~\ref{blr:lem-kernel}(1),
$\op{K(x,y)}\le c_0c^{d(x,y)}$ with $c=\max_e\op{\widehat C_e}<1$. The tree
$\tree$ is finite if $\beta=0$, and if $\beta=1$ a ball of radius $r$ is
contained in $2r+1$ consecutive copies of a fundamental domain. In both
cases $\sup_x\sum_y\op{K(x,y)}<\infty$, so $K$ is a bounded operator on
$\ell^2(\tree)$ by the Schur test, and $\tH K=I$ by
Lemma~\ref{blr:lem-kernel}(3); taking adjoints, $K\tH=I$. Since
$\tH\succeq0$, it follows that $\tH\succeq\norm K^{-1}$, so
$\lamC=\Sigm>0$. The converse is Theorem~\ref{thm:block-cone-unique}.

(2) The invertibility of the edge blocks and the connectedness force a
common dimension $p_i=p$. Replacing $D$ by $D-\lamC$ we may assume
$\lamC=0$ and have to show $\Sigm>0$; the last statement then follows
from Proposition~\ref{prop:block-physical} and
Theorem~\ref{thm:block-boundary}(1). Let $U$ be as in
Lemma~\ref{blr:lem-certificate}. Each $U_{j,e}$ is positive definite: if
$U_{j,e}v=0$, then $Z_e\succeq0$ gives $A_ev=0$ or $A_e^{\T}v=0$. Summing
\eqref{blr:eq-certificate} over $\tree$ and completing the square,
\begin{equation}\label{blr:eq-transport}
\langle f,\tH f\rangle=\sum_{\tilde e=xy}\binom{f_x}{f_y}^{\!\T}Z_e\binom{f_x}{f_y}
\ \ge\ u_0\sum_{\tilde e=xy}\abs{f_y-T_ef_x}^2,\qquad T_e=-U_{j,e}^{-1}A_e^{\T},
\end{equation}
for real finitely supported $f$, where the edges $\tilde e=xy$ of $\tree$
are oriented as their images $e=ij$, $u_0>0$ is a lower bound for the
$U_{j,e}$, and the $T_e$ are invertible. Suppose $\Sigm=0$. Then there are
real finitely supported $f_n$ with $\norm{f_n}=1$ and
$\langle f_n,\tH f_n\rangle\to0$.

Fix a spanning tree $\tau$ of $G$, a lift $\tilde\tau$ and a root $o$. The
deck group $\Gamma$ is free on generators $g_a$ indexed by the $\beta$
edges $a$ outside $\tau$: every vertex of $\tree$ is $\gamma\tilde\imath$
for unique $\gamma\in\Gamma$ and $i\in V$, the lifts of an edge $e=ij$ of
$\tau$ join $\gamma\tilde\imath$ to $\gamma\tilde\jmath$, and the lifts of
$a=ij$ join $\gamma\tilde\imath$ to $\gamma g_a\tilde\jmath$. Let
$f_n^i(\gamma)=f_n(\gamma\tilde\imath)$, an element of
$\ell^2(\Gamma;\R^p)$, and $b_n=f_n^o$. By \eqref{blr:eq-transport} along
the edges of $\tau$, $\norm{f_n^i-M_ib_n}_2\to0$ for fixed invertible
matrices $M_i$, so $\norm{b_n}_2$ is bounded away from zero, and along
the remaining edges,
\[
\norm{R_{g}b_n-M_{g}b_n}_2\to0\qquad(g=g_a),
\]
where $(R_g\varphi)(\gamma)=\varphi(\gamma g)$ and $M_g=M_j^{-1}T_aM_i$ is
invertible. After passing to a subsequence,
$\sum_\gamma b_n(\gamma)b_n(\gamma)^{\T}\to B\succeq0$, $B\ne0$. Since
$R_g$ permutes the summands, $M_gBM_g^{\T}=B$. Hence $W=\ran B$ is
invariant under $M_g$, and $M_g|_W$ is an isometry for the norm
$\abs{w}_B^2=w^{\T}(B|_W)^{-1}w$. Let $P_W$ be the orthogonal projection onto
$W$. Then $\norm{(I-P_W)b_n}_2^2\to\tr((I-P_W)B)=0$, and because
$P_WM_g-M_gP_W=P_WM_g(I-P_W)$,
\[
\norm{R_gP_Wb_n-M_gP_Wb_n}_2\le\norm{R_gb_n-M_gb_n}_2+\op{M_g}\norm{(I-P_W)b_n}_2\to0 .
\]
Therefore the scalar functions $\alpha_n(\gamma)=\abs{P_Wb_n(\gamma)}_B\ge0$
satisfy $\norm{\alpha_n}_2^2\to\dim W\ge1$ and
$\norm{R_g\alpha_n-\alpha_n}_2\to0$ for each of the $\beta$ free
generators. This is impossible for $\beta\ge2$ by Kesten's bound
\cite{Kesten}: the Cayley graph of $\Gamma$ is the regular tree of degree
$d=2\beta$; orienting it towards an end, each vertex has one parent and
$d-1$ children, and summing
$2\abs{\alpha_x\alpha_y}\le\alpha_x^2/\sqrt{d-1}+\sqrt{d-1}\,\alpha_y^2$
over the edges ($x$ the parent of $y$) bounds the quadratic form of the
adjacency operator by $2\sqrt{d-1}\,\norm\alpha_2^2$, whereas approximate
invariance forces the ratio to tend to $d>2\sqrt{d-1}$.
\end{proof}

For $p_i=1$, part (2) is the inequality $\sigm<\Sigm$ of
Corollary~\ref{cor:gap}, with a different proof. Combined with
Theorem~\ref{thm:block-converse} it gives the block version of
Corollary~\ref{cor:two-roots} for generic coefficients.

\begin{corollary}\label{blr:cor-two}
Let $G$ be connected and $(D,A)\in\cP_\R\setminus Z$, with $Z$ as in
Theorem~\ref{thm:block-converse}. If $\lamC(D,A)<0<\Sigm$, then the
number of positive solutions is even and at least two; one of them is
the physical branch at $z=0$, and at least one of them has odd index. For
every $(D,A)\in\cP_\R\setminus Z$ this applies to the coefficients
$(D-z,A)$ for all but finitely many $z$ in the interval
$(\lamC(D,A),\Sigm)$, which is nonempty when $\beta\ge2$ and all edge
blocks are invertible.
\end{corollary}

\begin{proof}
The first statement follows from Proposition~\ref{prop:block-physical}
and Theorem~\ref{thm:block-converse}. For the second, the line
$\{(D-z,A)\}$ is not contained in the algebraic set $Z$, because it
contains $(D,A)$, so it meets $Z$ in finitely many points.
\end{proof}

\begin{example}\label{blr:ex-K4-boundary}
Let $G=K_4$, $A_e=I_p$ and $D_i=3I_p$. Then $U_{i,e}=I_p$ satisfies
\eqref{blr:eq-certificate} and $q(X)=0$ for $X_i=I_p$, so the model lies on
the boundary of the cone. The unique positive solution has
$Q_{i|e}=2I_p$ and $V_i=\tfrac23I_p$; the other solution $Q_{i|e}=I_p$ of
\eqref{blr:eq-cavity} has vertex precision zero.
\end{example}

\begin{example}\label{blr:ex-K4-singular}
The invertibility in Corollary~\ref{blr:cor-cycle-rank}(2) cannot be
dropped, even when the finite precision matrix is positive definite. Let
$G=K_4$ with vertices $1,2,3,4$ and $p_i=2$. Let $A_e=I_2$ for the three
edges of the triangle $123$ and $A_e=\diag(1,0)$ for the three edges at
the vertex $4$; let $D_i=\diag(4,2)$ for $i\le3$ and $D_4=\diag(4,1)$. The
two coordinates decouple. The finite precision matrix is the direct sum of
$4I+\mathrm{Adj}(K_4)$, $2I+\mathrm{Adj}(K_3)$ and $(1)$, hence positive
definite. A certificate \eqref{blr:eq-certificate} is given by
$U_{i,e}=\diag(\tfrac43,1)$ at the ends of the triangle edges,
$U_{i,e}=\diag(\tfrac43,0)$ at the ends $i\le3$ of the edges to $4$, and
$U_{4,e}=\diag(\tfrac43,\tfrac13)$; and $q(X)=0$ for
$X_1=X_2=X_3=\diag(0,1)$, $X_4=0$. So $\lamC=0$. On the universal cover the
second coordinate lives on the lifts of the triangle, which are
two-sided infinite paths, and for a lift $x$ of a vertex $i\le3$ the
corresponding diagonal entry of $(\tH+\varepsilon)^{-1}$ is
\[
\frac1{2\pi}\int_0^{2\pi}\frac{dt}{2+\varepsilon+2\cos t}
=\frac1{\sqrt{\varepsilon(4+\varepsilon)}}\ \longrightarrow\ \infty .
\]
By Theorem~\ref{thm:block-boundary}(2) there is no positive solution,
although $\beta=3$. Comparing with Example~\ref{blr:ex-K4-boundary}, the
number of positive solutions on the boundary of the cone is not
determined by the graph.
\end{example}

\subsection{Comparison with walk-summability conditions}
\label{blr:subsec-preconditioned}

Several sufficient conditions for the convergence of block Gaussian belief
propagation generalize scalar walk-summability. We compare two of them
with the cone condition. For a simple graph let
$\mathsf J$ be the precision matrix of the finite model, of size
$\sum_ip_i$, with diagonal blocks $D_i$ and off-diagonal blocks $A_e$. For
a real symmetric matrix $S$ the \emph{comparison matrix} $M(S)$ has the
diagonal entries $S_{aa}$ and the off-diagonal entries $-\abs{S_{ab}}$; a
matrix $S$ with positive diagonal is walk-summable in the scalar sense of
\cite{MJW} if and only if $M(S)\succ0$. Following \cite{KSdP}, the model is
\emph{preconditioned walk-summable} if there is a block diagonal matrix
$\Lambda=\diag(\Lambda_i)$ with positive definite blocks
$\Lambda_i\in\Sym_{p_i}(\R)$ such that $\Lambda\mathsf J\Lambda$ is
walk-summable.

\begin{proposition}\label{prop:block-preconditioned}
\begin{enumerate}
\item Let $G$ be a loopless multigraph, $D_i\succ0$ for all $i$, and let
$N$ be the symmetric matrix with entries
$N_{ij}=\sum_{e=ij}\op{D_i^{-1/2}A_eD_j^{-1/2}}$. If the spectral radius
of $N$ is less than $1$, then the cone condition holds.
\item Let $G$ be simple. If the model is preconditioned walk-summable,
then the cone condition holds.
\item Neither implication can be reversed. For (1): on $K_4$ with
$p_i=2$ and $D_i=I_2$ there are rank one edge blocks with $\rho(N)\ge1$
that satisfy the cone condition. For (2): let $G$ be the star with a
center of dimension $2$ and three leaves of dimension $1$, let all $D_i$
be identity matrices, and let the edge blocks be $A_k=tv_k$ with
$t=\sqrt{3/5}$ and $v_k=(\cos\tfrac{k\pi}3,\sin\tfrac{k\pi}3)^{\T}$,
$k=0,1,2$. This model satisfies the cone condition, but
$\Lambda^{\T}\mathsf J\Lambda$ is not walk-summable for any invertible block
diagonal $\Lambda$; since block diagonal matrices with positive definite
blocks are invertible, this is stronger than the failure of
preconditioned walk-summability.
\end{enumerate}
\end{proposition}

\begin{proof}
(1) Let $\varsigma_i=\norm{D_i^{1/2}X_i}_F$. From
$\nuc{PRQ}\le\norm P_F\op R\norm Q_F$ applied to
$X_iA_eX_j=(X_iD_i^{1/2})(D_i^{-1/2}A_eD_j^{-1/2})(D_j^{1/2}X_j)$ we get
$q(X)\ge\varsigma^{\T}(I-N)\varsigma\ge(1-\rho(N))\abs\varsigma^2$.

(2) By Lemma~\ref{blr:lem-cone}(1) the cone condition says that $L>0$ on
$\overline\cC\setminus\{0\}$. The congruence
$D_i\mapsto\Lambda_iD_i\Lambda_i$, $A_e\mapsto\Lambda_iA_e\Lambda_j$
transforms $L$ by the linear automorphism
$V_i\mapsto\Lambda_i^{-1}V_i\Lambda_i^{-1}$,
$C_e\mapsto\Lambda_i^{-1}C_e\Lambda_j^{-1}$ of $\overline\cC$, so the cone
condition is invariant, and we may assume that $S=\mathsf J$ itself
satisfies $M(S)\succ0$. For a family $X$ let $r_{i\alpha}=\abs{X_ie_\alpha}$,
where $e_\alpha$ runs over the coordinate vectors of $\R^{p_i}$. Then
$\tr(X_iD_iX_i)=\sum_{\alpha,\gamma}(D_i)_{\alpha\gamma}\langle X_ie_\alpha,X_ie_\gamma\rangle
\ge\sum_\alpha(D_i)_{\alpha\alpha}r_{i\alpha}^2-\sum_{\alpha\ne\gamma}\abs{(D_i)_{\alpha\gamma}}r_{i\alpha}r_{i\gamma}$,
and expanding
$X_iA_eX_j=\sum_{\alpha,\gamma}(A_e)_{\alpha\gamma}(X_ie_\alpha)(X_je_\gamma)^{\T}$
into rank one matrices gives
$\nuc{X_iA_eX_j}\le\sum_{\alpha,\gamma}\abs{(A_e)_{\alpha\gamma}}r_{i\alpha}r_{j\gamma}$.
Adding up, $q(X)\ge r^{\T}M(S)r>0$ for $X\ne0$.

(3) For (1), choose at each vertex $i$ of $K_4$ unit vectors $u_{ij}\in\R^2$,
one for each neighbor $j$, with pairwise angles $2\pi/3$, so that
$\sum_ju_{ij}u_{ij}^{\T}=\tfrac32I$, and let $A_{ij}=c\,u_{ij}u_{ji}^{\T}$.
Then
$2\nuc{X_iA_{ij}X_j}=2c\abs{X_iu_{ij}}\abs{X_ju_{ji}}\le c(\abs{X_iu_{ij}}^2+\abs{X_ju_{ji}}^2)$,
so $q(X)\ge(1-\tfrac32c)\sum_i\tr X_i^2$, while $N=c\,\mathrm{Adj}(K_4)$ has
$\rho(N)=3c$. For $\tfrac13\le c<\tfrac23$ the cone condition holds and
$\rho(N)\ge1$.

For (2), let $X\succeq0$ be the matrix at the center and $x_k\ge0$ the
numbers at the leaves. Since $\sum_kv_kv_k^{\T}=\tfrac32I$,
\[
q(X,x)=\tr X^2+\sum_kx_k^2-2t\sum_kx_k\abs{Xv_k}
=\sum_k\bigl(x_k-t\abs{Xv_k}\bigr)^2+\bigl(1-\tfrac32t^2\bigr)\tr X^2 ,
\]
which is positive for $(X,x)\ne0$ because $1-\tfrac32t^2=\tfrac1{10}$. Let
$\Lambda$ be invertible and block diagonal. Multiplying $\Lambda$ by a
positive diagonal matrix does not change the definiteness of the
comparison matrix, so we may assume that the leaf blocks are $\pm1$ and that
the two columns $b_1,b_2$ of the center block are unit vectors; let
$c=\langle b_1,b_2\rangle$. The comparison matrix of
$\Lambda^{\T}\mathsf J\Lambda$ has unit diagonal, the entry $-\abs c$ between
the two center coordinates and the entries $-t\abs{\langle b_\alpha,v_k\rangle}$
between the center and the leaves. Eliminating the leaves leaves the
$2\times2$ matrix with diagonal entries
$1-t^2\sum_k\langle b_\alpha,v_k\rangle^2=\tfrac1{10}$ and off-diagonal entry
$-\abs c-t^2\sum_k\abs{\langle b_1,v_k\rangle\langle b_2,v_k\rangle}$. If
$\abs c\ge\tfrac1{10}$ it is not positive definite. If
$\abs c<\tfrac1{10}$, write $b_2=cb_1+\sqrt{1-c^2}\,b_\perp$ with a unit
vector $b_\perp\perp b_1$. By the double angle formula
$\sum_k\langle b_1,v_k\rangle^2\langle b_\perp,v_k\rangle^2=\tfrac38$, so
some $k$ has
$\abs{\langle b_1,v_k\rangle\langle b_\perp,v_k\rangle}\ge1/(2\sqrt2)$, and
for this $k$
\[
\abs{\langle b_1,v_k\rangle\langle b_2,v_k\rangle}
\ge\sqrt{1-c^2}\,\abs{\langle b_1,v_k\rangle\langle b_\perp,v_k\rangle}-\abs c>\tfrac15 .
\]
The off-diagonal entry then has absolute value greater than
$\tfrac35\cdot\tfrac15>\tfrac1{10}$, and the matrix is not positive
definite.
\end{proof}

Under either hypothesis of
Proposition~\ref{prop:block-preconditioned},
Theorem~\ref{thm:block-cone-unique} therefore gives a unique positive
solution. The paper \cite{KSdP} proves convergence of the means and
covariances of block Gaussian belief propagation for preconditioned
walk-summable models and does not discuss the set of fixed points. For
linear Gaussian observation models, \cite{DuEtAl} prove that the
information matrices of the factor graph algorithm have a unique positive
definite fixed point; the set of fixed points considered there is
different from the set of all fixed points with positive definite
pairwise beliefs counted here. We have not compared the cone condition
with the vector walk-summability of \cite{MalioutovThesis}. By
Lemma~\ref{blr:lem-cone}(3) the cone condition is the positivity of all
orthogonally twisted finite operators $H_O$, in the same way as scalar
walk-summability is the positivity of $D+A^s$ for all signings~$s$.

%% file: sections/block-complex.tex
\section{Block models: sheets, permutation groups and twisted poles}
\label{sec:block-complex}

\providecommand{\blcSym}{\operatorname{Sym}}
\providecommand{\blcMat}{\operatorname{Mat}}
\providecommand{\blcGL}{\mathrm{GL}}
\providecommand{\blcSp}{\mathrm{Sp}}
\providecommand{\blcO}{\mathrm{O}}
\providecommand{\blcp}{\mathbf p}
\providecommand{\blcq}{\mathbf q}
\providecommand{\blcMon}{\operatorname{Mon}}
\providecommand{\blcgraph}{\operatorname{graph}}

We use the block covariance system of Section~\ref{sec:block-real} with
complex coefficients. After proving generic irreducibility, we give
degree formulas and compute degrees and permutation groups for selected
graphs. In particular, the generic group for a cycle of length at least
three and uniform dimension is hyperoctahedral
(Theorem~\ref{thm:block-cycle}). We then prove the trace version of the
Floquet identity and derive the orthogonal twist at a pole under the
hypotheses of Theorem~\ref{thm:block-twisted-poles}.

For a fixed spectral parameter in the upper half-plane, the equations
studied here are of the type of operator-valued quadratic equations with a
positivity constraint, for which the solution with positive imaginary part
is unique \cite{HRFS,BMS}; Helton, Rashidi Far and Speicher remark that
such equations have many other solutions at a fixed $z$, without counting
them \cite{HRFS}. The count and the connectedness of all solutions are the
subject of this section.

\subsection{Notation and the coefficient map}
\label{blc:subsec:notation}

Throughout, $G=(V,E)$ is a finite loopless multigraph with at least one
vertex, and $\blcp=(p_i)_{i\in V}$ are positive integers, the \emph{vertex
dimensions}. We fix an orientation of each edge and write $e=ij$ for an
edge oriented from $i$ to~$j$. The coefficient space is
$\mathcal D\times\mathcal A$ with
\[
\mathcal D=\prod_{i\in V}\blcSym_{p_i}(\C),\qquad
\mathcal A=\prod_{e=ij\in E}\blcMat_{p_i\times p_j}(\C);
\]
a point is $(D,A)=((D_i)_i,(A_e)_e)$, and for $e=ij$ we write $A_{ij}=A_e$,
$A_{ji}=A_e^{\T}$, the edge being understood when there are parallel edges.
All symmetric matrices are complex symmetric (not Hermitian). As in
Definition~\ref{def:block-system}, the unknowns of the block system are
one symmetric \emph{cavity precision} $Q^{e}_{i\to j}\in\blcSym_{p_i}(\C)$
for each edge $e$ and each of its two directions, subject to
\begin{equation}\label{blc:eq:cavity}
Q^{e}_{i\to j}=D_i-zI-\sum_{\substack{e'\ni i,\ e'\ne e\\ e'=ik}}
A_{ik}\bigl(Q^{e'}_{k\to i}\bigr)^{-1}A_{ki}
\end{equation}
(this is the cavity form \eqref{blr:eq-cavity} of \eqref{eq:block-system},
equivalent to it by Lemma~\ref{blr:lem-bp}; for $p_i=1$ it is
\eqref{eq:cavity} with
$D$ replaced by $D-z$). The \emph{vertex precision} is
$R_i=D_i-zI-\sum_{e'\ni i}A_{ik}(Q^{e'}_{k\to i})^{-1}A_{ki}$, so that
$Q^e_{i\to j}=R_i+A_{ij}(Q^e_{j\to i})^{-1}A_{ji}$, and the \emph{edge
precision} is
$E_e=\bigl(\begin{smallmatrix}Q^e_{i\to j}&A_{ij}\\ A_{ji}&Q^e_{j\to i}\end{smallmatrix}\bigr)$.
A solution is \emph{admissible} (Definition~\ref{def:block-admissible},
transported to precisions by Lemma~\ref{blr:lem-bp}) if
all $Q^e_{i\to j}$, $R_i$ and $E_e$ are invertible; when the cavity and
vertex precisions are invertible, so are the edge precisions, because
$R_i$ is the Schur complement of $Q^e_{j\to i}$ in $E_e$. The
\emph{vertex covariances} are $V_i=R_i^{-1}$ and the \emph{edge
covariances} are $S_e=E_e^{-1}$.

It is convenient to parametrize solutions by their covariances. Let
\[
Y_{\blcp}=\Bigl\{(V,C):\ V_i\in\blcSym_{p_i}(\C),\ C_e\in\blcMat_{p_i\times p_j}(\C),\
\det V_i\ne0,\ \det S_e\ne0\Bigr\},\qquad
S_e=\begin{pmatrix}V_i&C_e\\ C_e^{\T}&V_j\end{pmatrix},
\]
a nonempty open subset of an affine space of dimension
\[
m=\sum_{i}\tfrac12p_i(p_i+1)+\sum_{e=ij}p_ip_j=\dim(\mathcal D\times\mathcal A),
\]
and define the \emph{coefficient map}
$\Psi_{\blcp}\colon Y_{\blcp}\to\mathcal D\times\mathcal A$,
$\Psi_{\blcp}(V,C)=(D,A)$, by
\begin{equation}\label{blc:eq:Psi}
A_e=(S_e^{-1})_{ij},\qquad
D_i=\sum_{e\ni i}(S_e^{-1})_{ii}-(d_i-1)\,V_i^{-1},
\end{equation}
where $(S_e^{-1})_{ii}$ and $(S_e^{-1})_{ij}$ are the diagonal block at $i$
and the off-diagonal block of $S_e^{-1}$.

\begin{lemma}\label{blc:lem:recovery}
\textup{(i)} For every $(D,A)$, the assignment
$Q^e_{i\to j}=(S_e^{-1})_{ii}$ is a bijection from the fiber
$\Psi_{\blcp}^{-1}(D,A)$ onto the set of admissible solutions of
\eqref{blc:eq:cavity} with $z=0$; its inverse is $V_i=R_i^{-1}$,
$S_e=E_e^{-1}$.

\textup{(ii)} $\Psi_{\blcp}$ is dominant, and there is a nonempty Zariski
open set $W\subseteq\mathcal D\times\mathcal A$ such that
$\Psi_{\blcp}^{-1}(W)\to W$ is finite \'etale of some degree
$N_{\blcp}(G)\ge1$.
\end{lemma}

We call $N_{\blcp}(G)$ the \emph{generic degree}: it is the number of
admissible complex solutions for generic coefficients, each of them
nondegenerate. For uniform dimension $p_i=p$ we write $N_p(G)$. By
Remark~\ref{blr:rem-scalar}, applied with $z=0$, which is not among the
generic $D_i$, and Theorem~\ref{thm:generic-degree}, $N_1(G)=\LG$.

\begin{proof}
(i) Let $(V,C)\in Y_{\blcp}$ and write
$S_e^{-1}=\bigl(\begin{smallmatrix}Q_{i\to j}&A_{ij}\\ A_{ji}&Q_{j\to i}\end{smallmatrix}\bigr)$.
By Jacobi's identity for complementary minors,
$\det Q_{i\to j}=\det V_j/\det S_e\ne0$, and the block inversion formula
gives $Q_{i\to j}-A_{ij}Q_{j\to i}^{-1}A_{ji}=V_i^{-1}$. Put $R_i=V_i^{-1}$.
Summing $Q^e_{i\to j}=R_i+A_{ij}(Q^e_{j\to i})^{-1}A_{ji}$ over the $d_i$
edges at $i$ and comparing with \eqref{blc:eq:Psi} gives
$D_i=R_i+\sum_{e\ni i}A_{ij}(Q^e_{j\to i})^{-1}A_{ji}$, which is the
definition of the vertex precision and, after subtracting the term of one
edge, the cavity equation \eqref{blc:eq:cavity}. All precisions are
invertible. Conversely, at an admissible solution the upper left block of
$E_e^{-1}$ is the inverse of the Schur complement $R_i$, so
$V_i=R_i^{-1}$, $S_e=E_e^{-1}$ lies in $Y_{\blcp}$ and is mapped to $(D,A)$.
For an isolated vertex the formula reads $D_i=V_i^{-1}$.

(ii) At the point $V_i=I$, $C_e=0$ we have
$S_e^{-1}=I-\bigl(\begin{smallmatrix}\delta V_i&\delta C_e\\ \delta C_e^{\T}&\delta V_j\end{smallmatrix}\bigr)+\cdots$,
so the differential of $\Psi_{\blcp}$ is
$(\delta V,\delta C)\mapsto(-\delta V,-\delta C)$. Hence $\Psi_{\blcp}$ is a
dominant map between irreducible varieties of the same dimension, and it
is finite and \'etale over a nonempty open subset of the target.
\end{proof}

Replacing $D_i$ by $D_i-zI$ gives, for fixed $(D,A)$, the set
\[
\cZ_{\blcp}(D,A)=\bigl\{(z,V,C)\in\C\times Y_{\blcp}:\ \Psi_{\blcp}(V,C)=(D-zI,A)\bigr\}
\]
of admissible solutions over the $z$-line. In analogy with
Definition~\ref{def:curve}, the \emph{block spectral curve}
$X_{\blcp}(G;D,A)$ is the union of the irreducible components of the
Zariski closure of $\cZ_{\blcp}(D,A)$, taken in the affine space with
coordinates $z$ and all entries of $(V,C)$, on which $z$ is not constant.
Points of the closure at which one of the determinants $\det V_i$,
$\det S_e$ vanishes or which lie at infinity are not counted as solutions.
By Lemma~\ref{blc:lem:recovery}\textup{(i)} the same curve is obtained, up
to a birational map over $\C(z)$, from the cavity precisions.

\subsection{Irreducibility}
\label{blc:subsec:irreducible}

\begin{theorem}\label{thm:block-irreducible}
Let $G$ be a finite loopless multigraph with at least one vertex and let
$\blcp$ be arbitrary vertex dimensions. For generic
$(D,A)\in\mathcal D\times\mathcal A$ the block spectral curve
$X_{\blcp}(G;D,A)$ is an irreducible curve, the map $z\colon
X_{\blcp}\to\C$ has degree $N_{\blcp}(G)$, and over all but finitely many
$z$ its fiber is the set of admissible solutions, all nondegenerate. The
set of such $(D,A)$ contains a dense open subset of the real coefficients,
and for these the analytic continuation in $z$ of the physical branch of
Proposition~\ref{prop:block-physical} reaches every point of a general
fiber.
\end{theorem}

The statement concerns the curve in all covariance coordinates. The degree
over the $z$-line of a single matrix entry of some $V_i$ may be smaller,
and we do not determine it. The proof does not use signings or transfer
matrices: it moves all vertex diagonals independently and then restricts to
a line. It uses two general facts. The first is the Bertini
irreducibility theorem for a morphism in the form of
\cite[Theorem~1.5]{PoonenSlavov}: if $\mathcal X$ is a geometrically
irreducible variety and $\phi\colon\mathcal X\to\PP^d$ is a morphism whose
nonempty fibers all have the same dimension, then the hyperplanes $H$ for
which $\phi^{-1}(H)$ is not geometrically irreducible form a constructible
set of dimension at most $\operatorname{codim}\phi(\mathcal X)+1$. The
second is that for a morphism of finite type to an irreducible base whose
generic fiber is geometrically irreducible, the fibers over a nonempty open
subset of the base are geometrically irreducible
\cite[Lemma 37.27.5]{StacksProject}.

\begin{proof}
\emph{Step 1: fixing the edge coefficients.} Irreducibility of $Y_{\blcp}$
does not imply irreducibility of the fibers of the projection to
$\mathcal A$. To control them we extend this projection across
$\{\det V_i=0\}$. Write $C_e=V_iT_eV_j$ with
$T_e\in\blcMat_{p_i\times p_j}$. Then
$\det S_e=\det V_i\det V_j\det(I-V_jT_e^{\T}V_iT_e)$ and the block
inversion formula gives
\begin{equation}\label{blc:eq:extended}
A_e=(S_e^{-1})_{ij}=-T_e\,(I-V_jT_e^{\T}V_iT_e)^{-1}.
\end{equation}
Let $\bar Y$ be the open subset of the affine space of all $(V,T)$,
$V_i\in\blcSym_{p_i}$, where the determinants on the right of
\eqref{blc:eq:extended} do not vanish; the open subset
$\{\det V_i\ne0\ \forall i\}$ of $\bar Y$ is isomorphic to $Y_{\blcp}$, and
\eqref{blc:eq:extended} defines a morphism $a\colon\bar Y\to\mathcal A$
extending the edge part of $\Psi_{\blcp}$. It has the section $V=0$,
$T=-A$, which consists of points that are not solutions. Let
$K=\C(\mathcal A)$. Since $a$ is dominant and $\bar Y$ is integral, the
generic fiber $\bar Y_K$ is integral. The $K$-rational point $V=0$, $T=-A$
is a smooth point of $\bar Y_K$, because at $V=0$ the derivative of
\eqref{blc:eq:extended} with respect to $T$ is $-I$. An integral variety
over a field of characteristic zero with a smooth rational point is
geometrically integral: an element of its function field algebraic over
$K$ is integral over the (normal) local ring of that point, hence lies in
it, and the residue map embeds the field of such elements into $K$. By
Lemma~\ref{blc:lem:recovery}\textup{(ii)} the open subsets $Y_{\blcp}$ and
$\Psi_{\blcp}^{-1}(W)$ of $\bar Y$ are nonempty and therefore meet the
generic fiber, so $\Psi_{\blcp}^{-1}(W)_K$ is geometrically irreducible. It
is a finite \'etale cover of degree $N_{\blcp}(G)$ of the open subset $W_K$
of $\mathcal D_K$.

\emph{Step 2: restriction to a general line.} Let $d=\dim\mathcal D$. If
$d=1$ there is nothing to do. If $d\ge2$, apply the Bertini theorem, over
an algebraic closure of $K$, to the composition
$\Psi_{\blcp}^{-1}(W)_K\to W_K\subseteq\mathcal D_K\cong\mathbb A^d\subset\PP^d$.
Its image is dense and its nonempty fibers are finite, so the exceptional
hyperplanes form a constructible set of dimension at most one; since
$d\ge2$, this set does not contain the generic point of the dual
projective space, so the preimage of the generic hyperplane is
geometrically irreducible; it is again a finite
\'etale cover of degree $N_{\blcp}(G)$ of a nonempty open subset of an
affine space, now of dimension $d-1$. Cutting successively by the generic
members of a flag of hyperplanes, we arrive at the generic line. Since
geometric irreducibility is insensitive to extension of the base field, we
may parametrize lines as $z\mapsto D^0-zM$ with
$(D^0,M)\in\mathcal D\times\mathcal D$: the map from pairs $(D^0,M)$ with
$M\ne0$ to the variety of affine lines in $\mathcal D$ is dominant, so
$K(D^0,M)$ is an extension of the field of definition of the generic line,
and the preimage of the line stays geometrically irreducible over it. We
conclude: for the family
\[
\mathcal Y=\bigl\{(A,D^0,M,z,y):\ y\in\Psi_{\blcp}^{-1}(W),\ \Psi_{\blcp}(y)=(D^0-zM,A)\bigr\}
\longrightarrow\mathcal A\times\mathcal D\times\mathcal D
\]
the generic fiber is geometrically irreducible. By the second fact quoted
above, the fibers over a nonempty Zariski open set
$\mathcal U\subseteq\mathcal A\times\mathcal D\times\mathcal D$ are
irreducible.

\emph{Step 3: the direction $M=I$.} The identity direction is a special
point of the space of directions, and irreducibility could be lost in a
specialization $M\to I$. Instead we use the action of
$P=(P_i)\in\prod_i\blcGL_{p_i}(\C)$ by congruence. The morphism
\[
(P,D',A')\longmapsto\bigl(A_e=P_i^{\T}A'_eP_j,\ D^0_i=P_i^{\T}D'_iP_i,\ M_i=P_i^{\T}P_i\bigr)
\]
is onto the open set $\{\det M_i\ne0\ \forall i\}$, because a nonsingular
complex symmetric matrix is congruent to the identity. Hence the preimage
of $\mathcal U$ is nonempty and open, and for a suitable fixed $P$ it
contains $\{P\}\times\mathcal U'$ with $\mathcal U'$ a nonempty Zariski open
set of pairs $(D',A')$. For such a pair and every $z$, congruence by $P$
transforms the precisions of the model $(D^0-zM,A)$ into those of
$(D'-zI,A')$ and the covariances accordingly; it preserves
\eqref{blc:eq:Psi} and all admissibility conditions and commutes with~$z$.
Therefore the part of $\cZ_{\blcp}(D',A')$ over
$\Omega=\{z:(D^0-zM,A)\in W\}$ is isomorphic to the fiber of $\mathcal Y$
over $(A,D^0,M)$: it is irreducible and finite \'etale of degree
$N_{\blcp}(G)$ over $\Omega$, and its fiber over $z\in\Omega$ consists of
all admissible solutions of the model $(D'-zI,A')$. After shrinking
$\mathcal U$ we may assume that $\Omega$ is nonempty; it is then the
complement of a finite set.

\emph{Step 4: conclusion.} Every component of the closure of
$\cZ_{\blcp}(D',A')$ on which $z$ is not constant meets the part over
$\Omega$, which is a smooth irreducible curve; so $X_{\blcp}$ is the closure
of that part. A connected finite \'etale cover has transitive monodromy,
so the continuation of any germ of solution reaches the whole fiber. The
real points of a nonempty Zariski open subset of $\mathcal D\times\mathcal A$
are dense and open in the real coefficients, and for real coefficients and
$z\ll0$ the physical branch is an admissible germ by
Proposition~\ref{prop:block-physical}; when $G$ is not connected, the
physical branch is understood as the product of the physical branches of
its connected components.
\end{proof}

\subsection{The generic degree: supermultiplicativity and first values}
\label{blc:subsec:degree}

In the next theorem the \emph{bowtie} is the graph of two triangles
sharing one vertex ($\LG=4$) and the \emph{diamond} is $K_4$ minus an edge
($\LG=5$).

\begin{theorem}\label{thm:block-degree}
Let $G$ be a finite loopless multigraph with at least one vertex.
\begin{enumerate}
\item[\textup{(i)}] For all vertex dimensions $\blcp,\blcq$,
\begin{equation}\label{blc:eq:supermult}
N_{\blcp+\blcq}(G)\ \ge\ N_{\blcp}(G)\,N_{\blcq}(G).
\end{equation}
In particular $N_p(G)\ge\LG^{\,p}$ for uniform dimension $p$, and
$N_{\blcp}(G)\ge\LG$ for all $\blcp$.
\item[\textup{(ii)}] If $G$ is a tree, $N_{\blcp}(G)=1$. If $G$ is connected
of cycle rank one, $N_{\blcp}(G)=2^{r}$, where $r$ is the smallest
dimension of a vertex on the cycle.
\item[\textup{(iii)}] For the bowtie, $N_1=4$ and, computer-assisted,
$N_2=36$.
\item[\textup{(iv)}] For the diamond, $N_1=5$. Computer-assisted lower
bounds: for the bowtie, $N_3\ge1079$; for the diamond, $N_2\ge66$.
\end{enumerate}
Consequently, for the bowtie $N_{a+2b+3c}\ge4^{a}36^{b}1079^{c}$ and for the
diamond $N_{2b+r}\ge66^{b}5^{r}$ $(r\in\{0,1\})$; for both graphs
$N_p(G)>\LG^{\,p}$ for every $p\ge2$.
\end{theorem}

Appendix~\ref{app:bowtie} proves $N_2=36$ by counting branches at
infinity of a reduced system. Its genericity checks use exact algebra
at explicit coefficients. The lower bounds in \textup{(iv)} use interval
certificates at explicit coefficients: each box contains a distinct,
nondegenerate admissible solution, which persists to generic
coefficients. For $p=3$ on the bowtie the certificate is for a reduced
system and the solutions are lifted to the full graph. The coefficient
data and certificates are in Appendix~\ref{app:cert}.
Part \textup{(ii)} is proved with
Theorem~\ref{thm:block-cycle} below. By
Theorem~\ref{thm:block-irreducible}, every lower bound for $N_{\blcp}(G)$ is
also a lower bound for the number of solutions reached from the physical
branch.

\begin{proof}[Proof of \textup{(i)}]
Taking direct sums of covariances, $V_i\oplus V_i'$ and $C_e\oplus C_e'$,
defines a morphism
$j\colon Y_{\blcp}\times Y_{\blcq}\to Y_{\blcp+\blcq}$, and by
\eqref{blc:eq:Psi}
$\Psi_{\blcp+\blcq}\circ j$ is the direct sum of the two coefficient
points. The function $\det d\Psi_{\blcp+\blcq}\circ j$, the Jacobian being
taken with respect to all coordinates of the large space, is regular on the
irreducible variety $Y_{\blcp}\times Y_{\blcq}$ and does not vanish at the
point where all $V_i=I$ and $C_e=0$ (proof of
Lemma~\ref{blc:lem:recovery}). Its zero set $Z$ is a proper closed subset.
Over the product of the two open sets of Lemma~\ref{blc:lem:recovery} the
map $\Psi_{\blcp}\times\Psi_{\blcq}$ is finite, so the image of $Z$ is a
proper closed subset. For a generic pair of coefficient points
$(c,c')$ all $N_{\blcp}(G)N_{\blcq}(G)$ points $j(y,y')$ of the fiber are
therefore nondegenerate solutions for the block diagonal coefficients
$c\oplus c'$, and they are distinct. By the implicit function theorem they
persist, in disjoint neighborhoods, when $c\oplus c'$ is perturbed in the
large coefficient space, including the couplings between the two blocks.
Choosing the perturbation in the open set of
Lemma~\ref{blc:lem:recovery}\textup{(ii)} gives \eqref{blc:eq:supermult}.
The bound $N_p\ge\LG^{\,p}$ follows by induction from $N_1(G)=\LG$
(Theorem~\ref{thm:generic-degree}). The same argument applies when some
$q_i$ are zero, with the conventions that an empty determinant is $1$ and
that an edge with an end of dimension zero carries no unknown, its
contribution to \eqref{blc:eq:Psi} canceling against one term
$V_i^{-1}$; the model with dimensions $\blcq$ is then the model on the
subgraph induced by $\{q_i>0\}$ and has degree at least one. With
$\blcp=\ones$ this gives $N_{\ones+\blcq}(G)\ge\LG$.
\end{proof}

The consequences stated after \textup{(iv)} follow from
\eqref{blc:eq:supermult} and $36>4^2$, $66>5^2$.

\begin{example}[Many positive solutions on $K_4$]\label{blc:ex:k4}
Supermultiplicativity has a real counterpart that shows that the
uniqueness statement of Theorem~\ref{thm:block-cone-unique} fails outside
the cone. On $K_4$ take $A_e=I_p$ for all edges and
$D_i=\diag(a_1,\dots,a_p)$ for all vertices, where $a_k=r_k+2/r_k$ for
distinct rational numbers $r_k\in(\sqrt2,2)$. \emph{This model has at
least $2^p$ distinct nondegenerate positive solutions at $z=0$, and so
does every real model in a neighborhood of it; for all of them the cone
condition of Definition~\ref{def:cone} fails.} Indeed, choose
$s_k\in\{r_k,2/r_k\}$ for each $k$ and put $Q_{i\to j}=\diag(s_k)$ on all
twelve arcs. Since $K_4$ is $3$-regular, \eqref{blc:eq:cavity} reads
$D-2Q^{-1}=Q$, which holds because $s_k+2/s_k=a_k$. The vertex precisions
$\diag(s_k-1/s_k)$ and the edge precisions, with eigenvalues $s_k\pm1$, are
positive definite, since $s_k>1$. In the entry $(k,l)$ of the symmetric
matrices the derivative of the cavity update is $B/(s_ks_l)$, where $B$ is
the non-backtracking matrix of $K_4$ on its twelve arcs. This integer
matrix has row sums $2$, so its rational eigenvalues are integers of
absolute value at most $2$, whereas $s_ks_l$ is a rational number greater
than $1$ and different from $2$: it equals $r_kr_l>2$, or $4/(r_kr_l)<2$, or
$2r_k/r_l\ne2$ with $k\ne l$, and $s_k^2=2$ is impossible. Hence
$I-B/(s_ks_l)$ is invertible for all $k\le l$, each of the $2^p$ solutions
is nondegenerate with respect to all block variations, and they persist on
an open set of real coefficients, on which the channels are no longer
decoupled. Finally the quadratic form
$q(X)=\sum_i\tr(X_iD_iX_i)-2\sum_{e=ij}\norm{X_iA_eX_j}_*$ of the cone
condition ($\norm{\cdot}_*$ the nuclear norm) takes the value
$4a_k-12<0$ at $X_i=e_ke_k^{\T}$, because $r+2/r<3$ on $(\sqrt2,2)$. The
finite precision matrix of this model is positive definite.
\end{example}

\subsection{Two finite formulas for the degree}
\label{blc:subsec:formulas}

We record two expressions of $N_{\blcp}(G)$ in terms of the graph and the
vertex dimensions alone. Neither is practical; their interest is that the
degree is determined by a finite computation that involves no choice of
generic coefficients.

Let $x_1,\dots,x_m$ be the independent entries of $(V,C)$, let
$H=\prod_i\det V_i\prod_e\det S_e$, of degree
$h=\sum_ip_i+\sum_{e=ij}(p_i+p_j)$, and let $P_a=H\,\Psi_a\in\Z[x]$
$(1\le a\le m)$ be the numerators of the $m$ components of
$\Psi_{\blcp}$; the admissible region is $\{H\ne0\}$. Put
$b_0=(m+1)(h+1)^m$, introduce indeterminates $v_1,\dots,v_m$ and
$t=\sum_av_ax_a$. For $0\le k\le b_0$ let $M_k$ be the matrix over
$\mathbb Q(v)$ whose columns are indexed by the pairs
$(\alpha,j)\in\N^m\times\{0,\dots,k\}$ with $\abs\alpha+j\le b_0$, the
column $(\alpha,j)$ being the coefficient vector, in the monomial basis of
$\mathbb Q(v)[x]$, of $H^{\,b_0-\abs\alpha}\prod_aP_a^{\alpha_a}\,t^{\,j}$.

\begin{proposition}[Rank formula]\label{blc:prop:rank-formula}
$N_{\blcp}(G)$ is the number of $k\in\{0,\dots,b_0-1\}$ for which $M_k$ has
full column rank.
\end{proposition}

\begin{proof}
Let $F=\mathbb Q(v)$ and consider the closure of the image of
$x\mapsto(y,T)=(\Psi_{\blcp}(x),t(x))$ in $\mathbb A^{m+1}_F$. Two distinct
points of a geometric generic fiber of $\Psi_{\blcp}$ differ by a nonzero
vector, whose pairing with the independent $v_a$ is nonzero; so $t$
separates the $N=N_{\blcp}(G)$ points of the fiber and is a primitive
element. The closure is an irreducible hypersurface with equation
$R(y,T)$, and $\deg_TR=N$ by Gauss's lemma, the $y_a$ being algebraically
independent. Substituting the monomials $y^\alpha T^j$ of total degree at
most $b_0$ and multiplying by $H^{b_0}$ gives polynomials in $x$ of degree
at most $(h+1)b_0$, which span a space of dimension at most
$\binom{(h+1)b_0+m}{m}\le(h+1)^m\binom{b_0+m}{m}<\binom{b_0+m+1}{m+1}$,
the number of monomials substituted. Hence there is a nontrivial relation
of total degree at most $b_0$; $R$ divides it, so $R$ has total degree at
most $b_0$ and $N\le b_0$. The kernel of $M_k$ is the space of relations
of total degree at most $b_0$ and $T$-degree at most~$k$; since all columns
were cleared by the same power of $H$, no spurious solutions with $H=0$
enter. A nonzero relation is a multiple of $R$, so the kernel vanishes for
$k<N$, and it contains $R$ for $k\ge N$.
\end{proof}

The second formula removes ranks and elimination altogether, at the price
of a range of summation that is larger still. To each vertex attach
symmetric matrices of letters $V_i,W_i$ of size $p_i$, to each edge $e=ij$
a matrix of letters $C_e$ and a symmetric matrix of letters $T_e$ of size
$p_i+p_j$; let $X_1,\dots,X_\nu$ be all the independent letters. Number
the $m$ output positions (the upper triangular entries of the $D_i$ and
all entries of the $A_e$) by $a=1,\dots,m$, and let $D_i(t)$, $A_e(t)$ be
the matrices with $t^{w_a}$ at position $a$, for weights $w_a$ fixed below.
Let $(f_r)_r$ be the list of the following polynomials in $\Z[X,t]$: all
entries of $V_iW_i-I$ and of $S_eT_e-I$, the upper triangular entries of
$\sum_{e\ni i}(T_e)_{ii}-(d_i-1)W_i-D_i(t)$, and all entries of
$(T_e)_{ij}-A_e(t)$. Define
\[
\begin{gathered}
L_0=\prod_ip_i!\prod_{e=ij}(p_i+p_j)!,\quad u_1=2+m+(2\abs E+2)L_0,\\
b_1=(2m+1)(h+2)^{2m},\quad c_1=\binom{b_1+2m+1}{2m+1},\\
k_1=(4+2h)(b_1+1)^2,\quad w_a=(k_1+1)^{a-1},\quad j_1=(k_1+1)^{2m},
\end{gathered}
\]
let $\ell$ be any prime larger than $c_1!\,u_1^{2b_1c_1}+b_1+2$, and put
$s=2^{\,j_1\lceil\log_2\ell\rceil}$ and $q=\ell^{\,s\cdot b_1!}$.

\begin{proposition}[Multinomial formula]\label{blc:prop:multinomial}
With these choices, $N_{\blcp}(G)$ is the residue in $\{0,\dots,\ell-1\}$
modulo $\ell$ of the integer
\[
(-1)^{\nu}\sum_{\gamma\in((q-1)\Z_{>0})^{\nu}}\ \sum_{j\ge0}(-1)^j\,
\bigl[X^\gamma t^{sj}\bigr]\prod_r\bigl(1-f_r(X,t)^{\,q-1}\bigr),
\]
where $[X^\gamma t^{a}]$ denotes the coefficient of a monomial. The sum is
finite.
\end{proposition}

Expanding the powers by the multinomial theorem, the right-hand side is a
signed sum of products of multinomial coefficients over the ways of
choosing, for each row $r$ that is used, $q-1$ monomials of $f_r$ such that
every letter occurs a positive multiple of $q-1$ times in total. For
$\blcp=\ones$ this should be compared with the sum over edge subsets in
Theorem~\ref{thm:generic-degree}, where each edge is either chosen or not.

\begin{proof}
First, the generic degree is the same over $\mathbb F_\ell$. In every
characteristic $\{H\ne0\}$ is irreducible and $d\Psi_{\blcp}=-I$ at
$V_i=I$, $C_e=0$, so $\Psi_{\blcp}$ is generically finite and separable. As
in the proof of Proposition~\ref{blc:prop:rank-formula}, now with $v$ as
$m$ further variables, the closure of the image of
$(x,v)\mapsto(\Psi_{\blcp}(x),v,\sum_av_ax_a)$ is an irreducible hypersurface
$R(y,v,Z)=0$ with $\deg_ZR$ equal to the generic degree, and the same
monomial count (with $h+2$ and $2m$ in place of $h+1$ and $m$) gives
$\deg R\le b_1$, in characteristic $0$ and $\ell$. The coefficient vectors
of the $c_1$ substituted monomials, cleared by $H^{b_1}$, are integral with
absolute sum at most $u_1^{2b_1}$ (the coefficients of $H$ have absolute
sum at most $L_0$, those of $P_a$ at most $(2\abs E+1)L_0$), so every minor
of the matrix they form is smaller than $\ell$ in absolute value, and a
set of columns is independent over $\mathbb Q$ if and only if it is
independent over $\mathbb F_\ell$. The smallest $Z$-degree of a relation
of total degree at most $b_1$, which is $\deg_ZR$, is therefore the same in
both characteristics: the degree over $\mathbb F_\ell$ is
$N:=N_{\blcp}(G)\le b_1<\ell$.

Second, work over $\mathbb F_\ell$. Differentiating
$R(\Psi_{\blcp}(x),v,\sum_a v_ax_a)=0$ with respect to $v_a$ gives
$x_a=-R_{v_a}/R_Z$ in the function field, and $H(-R_v)=R_Z^{\,h}H(x)\ne0$.
Hence
$b=\operatorname{lc}_Z(R)\operatorname{Res}_Z(R,R_Z)\operatorname{Res}_Z(R,H(-R_v))$
is a nonzero polynomial in $(y,v)$, of degree less than $k_1$, and where
$b\ne0$ the $N$ distinct roots $Z$ of $R$ are in bijection, by
$x_a=-R_{v_a}/R_Z$, with the admissible points of the fiber of
$\Psi_{\blcp}$ over $y$. Substituting $t^{(k_1+1)^{a-1}}$ for the $a$th of
the $2m$ coordinates of $(y,v)$ sends distinct monomials of degree less
than $k_1$ to distinct powers of $t$, so $b$ becomes a nonzero polynomial
of degree less than $j_1$ in $t$. If $\zeta$ is a primitive $2s$th root of
unity, then $\zeta\in\mathbb F_{\ell^s}$, and $b(\zeta)\ne0$ because
$2s>\ell^{j_1}$ forces the degree of $\zeta$ over $\mathbb F_\ell$ to exceed
$j_1$. The roots of $R$ lie in extensions of $\mathbb F_{\ell^s}$ of degree
at most $b_1$, hence in $\mathbb F_q$, and since the inverses $W_i,T_e$ are
determined by the covariances, the system $f_r(X,\zeta)=0$ has exactly $N$
solutions in $\mathbb F_q^{\nu}$.

Third, $1-f^{q-1}$ is the indicator function of $\{f=0\}$ on
$\mathbb F_q^\nu$, so
$\sum_{X}\prod_r(1-f_r(X,\zeta)^{q-1})=N$ in $\mathbb F_\ell$. As
$\sum_{x\in\mathbb F_q}x^{u}$ is $-1$ if $u$ is a positive multiple of
$q-1$ and $0$ otherwise, only the monomials $X^\gamma$ with
$\gamma\in((q-1)\Z_{>0})^\nu$ contribute, with the factor $(-1)^\nu$.
Averaging over the $s$ roots $\zeta$ of $t^s=-1$, which are all primitive
$2s$th roots of unity, replaces $t^a$ by $0$ unless $a=sj$, and then by
$(-1)^j$. This gives the stated integer modulo $\ell$, and $0<N<\ell$.
\end{proof}

The constants are upper bounds chosen to make the argument uniform, not
recommended values; already for a triangle with $p=2$ they are far beyond
any computation, and we have not evaluated the formula in any nontrivial
case.

\subsection{Permutation groups}
\label{blc:subsec:groups}

For generic $(D,A)$ let $\blcMon_z(G,\blcp)\subseteq S_N$,
$N=N_{\blcp}(G)$, be the monodromy group of the finite \'etale cover of
Theorem~\ref{thm:block-irreducible} over the punctured $z$-line. It is
defined up to conjugacy, it is transitive, and it is the same group for
all coefficients in a nonempty Zariski open set, by the spreading-out
argument at the end of the proof of
Proposition~\ref{blc:prop:line-group}. For a finite \'etale cover
$Z\to B$ of degree $N$ of a smooth connected variety we use the cover
$\mathcal O(Z)\to B$ of \emph{orderings}, whose fiber is the set of
bijections from $\{1,\dots,N\}$ to the fiber of $Z$; it is the open and
closed subset of the $N$-fold fiber product of $Z$ where the points are
pairwise distinct. The group $S_N$ acts simply transitively on its fibers
and transitively on its connected components, and the stabilizer of a
component is the monodromy group of $Z\to B$. Consequently, if
$B_1\subseteq B$ is a smooth connected locally closed subvariety, the
monodromy group of $Z|_{B_1}$ is a subgroup of that of $Z$, with equality if
and only if every component of $\mathcal O(Z)$ remains connected over
$B_1$. Over another algebraically closed field of characteristic zero, as
in the proof below, the same statements hold, with the same proof, when the
monodromy group is read as the Galois group of the Galois closure of the
extension of function fields, acting on the generic fiber.

\begin{proposition}\label{blc:prop:line-group}
For generic $(D,A)$, the group $\blcMon_z(G,\blcp)$ is equal, up to
conjugacy in $S_N$, to the monodromy group of the cover
$\Psi_{\blcp}^{-1}(W)\cap\{A\text{ fixed}\}$ of an open subset of
$\mathcal D$, in which the edge coefficients are fixed and all vertex
diagonals vary independently. Thus no permutation is lost by restricting
the diagonals to the line $D-zI$.
\end{proposition}

\begin{proof}
Fix the edge coefficients at the generic point of $\mathcal A$ as in
Step~1 of the proof of Theorem~\ref{thm:block-irreducible}, and let
$\mathcal O_1$ be a connected component of the cover of orderings of
$\Psi_{\blcp}^{-1}(W)_K\to W_K$ over an algebraic closure of~$K$. It is
smooth and connected, hence irreducible, and finite \'etale over $W_K$. The
Bertini argument of Step~2 applies to $\mathcal O_1$ in place of
$\Psi_{\blcp}^{-1}(W)_K$ and shows that its restriction to the generic line
$D^0-zM$ is irreducible. By the criterion above, the monodromy group over
the generic line equals the monodromy group over $W_K$. Congruence by
$P\in\prod_i\blcGL_{p_i}$ maps the family of all diagonals for the edge
coefficients $A$ to the family of all diagonals for $A'$, and the line
$D^0-zM$ to the line $D'-zI$, compatibly with the covers; as in Step~3
this transports the equality to the direction $M=I$. To pass from generic
points to a nonempty open set of complex coefficients, one spreads out the
finitely many components of the two covers of orderings over the
parameter space, after a finite base change over which they are defined,
and applies \cite[Lemma 37.27.5]{StacksProject} to each of them: on a
nonempty open set of parameters each component stays irreducible, so the
number of components, and with it the order of each of the two groups, is
constant. Since the group of the line is always contained in the group of
the family of all diagonals, the two groups coincide.
\end{proof}

Equivalently, $\blcMon_z(G,\blcp)$ is the Galois group over
$\overline{\mathbb Q(v,D,A)}(z)$ of the minimal polynomial of the primitive
element $t$ of Proposition~\ref{blc:prop:rank-formula} along the line
$D-zI$; like the rank formula, this is a finite characterization that we
have not been able to evaluate beyond the examples below.

\begin{proposition}[Subdivision]\label{blc:prop:subdivision}
Let $G'$ be obtained from $G$ by replacing an edge $e_0=ij$ by a path
$i\,w\,j$ through a new vertex $w$ of dimension $p_w\ge\min(p_i,p_j)$.
Then the generic degrees of $G'$ and of $G$ are equal, and the generic
group $\blcMon_z$ of $G'$ contains that of $G$, up to conjugacy. In
particular the two groups are equal if the group of $G$ is the full
symmetric group.
\end{proposition}

\begin{proof}
Let $J$ be the diagonal coefficient at $w$ and
$B_1\in\blcMat_{p_i\times p_w}$, $B_2\in\blcMat_{p_w\times p_j}$ the
coefficients of the edges $iw$, $wj$; generically $\det J\ne0$. Eliminating
$w$ by a Schur complement gives effective coefficients on $G$,
\[
A'_{e_0}=-B_1J^{-1}B_2,\qquad D'_i=D_i-K_i,\quad D'_j=D_j-K_j,\qquad
K_i=B_1J^{-1}B_1^{\T},\quad K_j=B_2^{\T}J^{-1}B_2,
\]
all other coefficients being unchanged. This map of coefficient spaces is
dominant: for $J=I$ every $p_i\times p_j$ matrix is of the form $-B_1B_2$,
since its rank is at most $p_w$. Let $X_i,X_j$ be the cavity precisions on
$e_0$ at a solution for $G$ with the effective coefficients, and
$E'=\bigl(\begin{smallmatrix}X_i&A'_{e_0}\\ A_{e_0}'^{\T}&X_j\end{smallmatrix}\bigr)$.
Put
\[
Q_{i\to w}=X_i+K_i,\quad Q_{j\to w}=X_j+K_j,\quad
Q_{w\to j}=J-B_1^{\T}Q_{i\to w}^{-1}B_1,\quad
Q_{w\to i}=J-B_2Q_{j\to w}^{-1}B_2^{\T},
\]
and $R_w=J-B_1^{\T}Q_{i\to w}^{-1}B_1-B_2Q_{j\to w}^{-1}B_2^{\T}$, keeping all
other cavity precisions. By the Woodbury identity
$B_1Q_{w\to i}^{-1}B_1^{\T}=K_i+A'_{e_0}X_j^{-1}A_{e_0}'^{\T}$, so the message
entering $i$ along $iw$ differs from the message entering $i$ along $e_0$
exactly by the shift $K_i$ of the diagonal; hence the cavity equations of
$G'$ and of $G$ correspond, and the vertex precisions at $i$ and $j$ are
the same. Computing the determinant of the matrix with diagonal blocks
$(Q_{i\to w},J,Q_{j\to w})$ and off-diagonal blocks $B_1,B_2$, and of its
leading two-by-two block, by eliminating either the middle or the outer
blocks first, gives
\[
\det J\,\det X_i=\det Q_{i\to w}\,\det Q_{w\to j},\qquad
\det J\,\det E'=\det Q_{i\to w}\,\det Q_{j\to w}\,\det R_w,
\]
and the first identity with $i$ and $j$ exchanged. Hence an admissible
solution for $G'$ gives an admissible solution for $G$, injectively.
Conversely an admissible solution for $G$ lifts to an admissible solution
for $G'$ as soon as $\det(X_i+K_i)\det(X_j+K_j)\ne0$: the identities then
show that $Q_{w\to j}$, $Q_{w\to i}$ and $R_w$ are invertible, and
$\det E_{iw}=\det Q_{i\to w}\det R_w$. Replacing $(B_1,B_2)$ by
$(\lambda B_1,\lambda^{-1}B_2)$ and correcting $D_i,D_j$ so that the
effective coefficients are unchanged, the first condition has a nonzero
limit as $\lambda\to0$ and the second as $\lambda\to\infty$, for each of
the finitely many solutions of a generic fiber for $G$. Since $K_i$ is
proportional to $\lambda^2$ and $K_j$ to $\lambda^{-2}$, both conditions are
Laurent polynomials in $\lambda^2$ that are not identically zero, so they
hold simultaneously for all but finitely many values of $\lambda^2$. So for generic
coefficients of $G'$ all solutions of the effective model lift, and the
degrees are equal.

For the groups, fix generic coefficients $(D,A)$ of $G$, a base point and
finitely many loops in the $z$-plane generating $\blcMon_z(G,\blcp)$.
Factor $A_{e_0}=-U_0V_0$ through $\C^{p_w}$ and, for $b>0$, give $w$ the
diagonal $bI$ and the edges the coefficients $\sqrt b\,U_0$, $\sqrt b\,V_0$,
adding $U_0U_0^{\T}$ and $V_0^{\T}V_0$ to $D_i$ and $D_j$. With the spectral
parameter, the effective edge coefficient is $A_{e_0}\,b/(b-z)$ and the
effective diagonals are $D_i-zI-zU_0U_0^{\T}/(b-z)$ and similarly at $j$.
As $b\to\infty$ these converge to the original path $(D-zI,A)$, uniformly
on the chosen loops. Hence for large $b$, and for generic coefficients of
$G'$ near this point, the effective paths are homotopic to the original
loops inside the open set over which the cover for $G$ is finite \'etale.
By the first part all solutions lift outside finitely many values of $z$,
which the loops avoid after a small deformation. Thus every generating
permutation occurs for $G'$.
\end{proof}

\begin{example}\label{blc:ex:symmetric-groups}
Take the bowtie, the graph $D_3$ of three parallel edges, or $K_4$, with
scalar vertices, and subdivide every edge by one vertex of arbitrary
dimension $s_e\ge1$, with generic coefficients. By
Proposition~\ref{blc:prop:subdivision} the resulting block models have
$4$, $5$ and $15$ sheets, and their groups are $S_4$, $S_5$ and $S_{15}$ by
Proposition~\ref{irr:prop:small-monodromy}. Starting from the bowtie with
uniform dimension $2$ and subdividing with $s_e\ge2$ preserves the
$36$ sheets of Theorem~\ref{thm:block-degree}.
\end{example}

\begin{proposition}\label{blc:prop:product-group}
Let $\blcp=\sum_a\blcp_a$ be a finite sum of vertex dimensions on $G$, with
generic degrees $N_a$ and generic groups $M_a\subseteq S_{N_a}$. Then
$\blcMon_z(G,\blcp)$ has a subgroup that preserves a set of $\prod_aN_a$
sheets and acts on it as the product $\prod_aM_a$ in its product action.
In particular $\prod_a\abs{M_a}$ divides $\abs{\blcMon_z(G,\blcp)}$.
\end{proposition}

\begin{proof}
Take independent generic coefficients for the summands and form the block
diagonal model. The argument for Theorem~\ref{thm:block-degree}\textup{(i)}
applies to any finite number of summands: for all but finitely many $z$,
all direct sums of solutions are nondegenerate solutions of the large
model. Indeed, nondegeneracy is an open condition on the fiber product over
the $z$-line of the finite \'etale covers of the summands; each component of
this fiber product has a point over $z=0$, where the coefficients form a
generic tuple and the argument for \textup{(i)} gives nondegeneracy, so the
condition fails only over finitely many $z$. Shifting the diagonals of the summands by independent generic
constants, we may assume that their finite branch values are pairwise
disjoint. The group $M_a$ is generated by lassos around the branch values
of the $a$th summand; along such a lasso the sheets of the other summands
return to themselves, so on the direct sum sheets the lassos generate
$\prod_aM_a$. The lifts of these finitely many loops to the direct sum
sheets form a compact set of nondegenerate solutions, so after a small
generic perturbation of all coefficients of the large model, the
continuation of these sheets along the loops, and the correspondence of
their endpoints, is unchanged. The remaining sheets have finitely many
exceptional values, which the loops avoid after a small deformation. The
permutations so obtained generate a subgroup of $\blcMon_z(G,\blcp)$
preserving the set of direct sum sheets, with image $\prod_aM_a$ in the
symmetric group of that set.
\end{proof}

The action on the complementary sheets is not determined by this
construction, so the product is a quotient of a subgroup. By
Proposition~\ref{irr:prop:small-monodromy}, for the bowtie with uniform
dimension $2$ a subgroup of the group of the $36$ sheets acts on $16$
of them as $S_4\times S_4$, and the order of the group is a multiple
of $576$.

\subsection{Scalar separators}
\label{blc:subsec:separators}

There is a large class of block models whose degree is still given by
leafless edge subsets.

\begin{theorem}\label{thm:block-separator}
Let $G$ be a connected loopless multigraph with vertex dimensions $\blcp$.
Suppose that there is a nonempty set $U\subseteq V$ of vertices of
dimension one such that every connected component $C$ of $G-U$ is a tree
and is joined to $U$ by at most three edges. Then
\[
N_{\blcp}(G)=\LG .
\]
The dimensions of the vertices outside $U$ are arbitrary, several of the
edges leaving a component may end at the same vertex of $U$, and there is
no restriction on the cycle rank.
\end{theorem}

For instance, raising the dimensions on an independent set of vertices of
a cubic graph does not change the degree: $K_4$ with one vertex of
arbitrary dimension has $15$ sheets, and $D_3$ with dimensions $(1,p)$ has
$5$. The case $U=V$ is Theorem~\ref{thm:generic-degree}.

\begin{proof}
The inequality $N_{\blcp}(G)\ge\LG$ is part of
Theorem~\ref{thm:block-degree}\textup{(i)}. For the reverse inequality we
replace each component by a scalar vertex.

\emph{The response of a tree.} Let $C$ be a component of $G-U$, with
internal precision matrix $J_C$ (the block matrix of the $D_i$, $i\in C$,
and the $A_e$ of the edges inside $C$), generically invertible. Call the
$k\le3$ edges from $C$ to $U$ its \emph{ports}; for a port $a$ let $b_a$ be
its coefficient, a column vector supported on the block of its end in $C$,
and put $K_{ab}=b_a^{\T}J_C^{-1}b_b$. At a solution let $x_a$ be the
(scalar) cavity precision entering $C$ through the port $a$, and $m_a$ the
message returned through $a$, that is, the quantity subtracted from the
diagonal of the outer end of $a$. Because belief propagation on a tree is
Gaussian elimination, attaching to $C$ the other ports as leaves with
diagonals $x_b$ and eliminating gives
\begin{equation}\label{blc:eq:tree-response}
m_a=K_{aa}+\sum_{b,b'\ne a}K_{ab}\,\bigl[(\diag(r)-P)_{-a}^{-1}\bigr]_{bb'}\,K_{b'a},
\qquad r_b=x_b-K_{bb},
\end{equation}
where $P$ is the off-diagonal part of $K$ and the subscript $-a$ denotes
the principal submatrix with $a$ removed. The inputs $x_a$ determine all
cavity precisions inside $C$. Ports ending at the same vertex of $U$ are
treated as separate ports.

\emph{Replacement by a star.} Let $G_*$ be the scalar graph in which each
$C$ is replaced by one vertex with diagonal $h_C\ne0$, joined to $U$ by
the same ports with weights $g_a$; its response matrix is
$K^*_{ab}=g_ag_b/h_C$. For $k\le3$ and nonzero off-diagonal entries of $K$
one can choose $(g,h_C)$ with $K^*_{ab}=K_{ab}$ for all $a\ne b$: for $k=3$
take $h_C=P_{12}P_{13}/P_{23}$ and $g=(h_C,P_{12},P_{13})$; for $k=2$ take
$h_C=1$, $g=(u,P_{12}/u)$; for $k=1$, $h_C=1$, $g=(u)$, with $u\ne0$ free;
in all cases $(g,h_C)$ may further be replaced by $(\tau g,\tau^2h_C)$. Now
add $\delta_a=K^*_{aa}-K_{aa}$ to the diagonal of the outer end of each
port $a$, and replace $x_a$ by $x_a+\delta_a$. In
\eqref{blc:eq:tree-response} the quantities $r_a$ and $P$ are unchanged, so
$m_a$ changes by $\delta_a$; the vertex precisions on $U$ and the cavity
precisions on edges inside $U$ are unchanged, and the cavity precision
entering each port changes by its own $\delta_a$, consistently. Thus an
admissible solution for $G$ is mapped to a solution of the scalar system
for $G_*$ with the coefficients just described. The map is injective: the
$x_a$ are recovered by subtraction, and the interior of each tree is
determined by them.

\emph{Genericity.} We must check that, for generic coefficients of $G$ and
generic auxiliary parameters $u,\tau$, (a) the coefficients of $G_*$ are
generic and (b) every solution of the fiber is mapped to an admissible
solution of $G_*$. For (a): the map from the coefficients of $G$ and the
auxiliary parameters to the coefficients of $G_*$ is dominant. Indeed, the
off-diagonal entries $K_{ab}$ can be made nonzero (take all diagonals of
$C$ equal to $I$, small negative couplings between the first coordinates
along the edges of $C$, and $b_a$ in the first coordinate, and use the
Neumann series); scaling $b_a$ by $\lambda_a$ scales $P_{ab}$ by
$\lambda_a\lambda_b$, which for $k\le3$ reaches all nonzero off-diagonal
data, and together with $u,\tau$ all values of $(g,h_C)$; the outer
diagonals are free. For (b): all inverses used, the invertibility of all
cavity, vertex and edge precisions of the star model, and the
invertibility of the intermediate Schur complements in the trees, are open
conditions on the irreducible variety of triples (coefficients, solution,
auxiliary parameters), which by Lemma~\ref{blc:lem:recovery} is the
product of $Y_{\blcp}$ with a torus. They hold simultaneously at one point:
in the real example of (a) with small couplings on all ports and inside
$U$, positive outer diagonals, and the positive solution close to the
decoupled one, all eliminations are positive definite. There $P_{ab}>0$, so
$h_C>0$, and for real positive $u,\tau$ the star model is weakly coupled
and positive definite as well. Hence the bad
triples form a proper closed subset, and since the projection to the
coefficients and auxiliary parameters is generically finite, for a generic
choice of these no solution of the fiber is bad.

It follows that $N_{\blcp}(G)\le N_{\ones}(G_*)=L(G_*)$. Finally
$L(G_*)=\LG$. Fix which ports of each $C$ belong to an edge subset. Inside
the tree $C$, the edge subsets that create no vertex of degree one are:
only the empty set if no port is chosen; none if one port is chosen; and,
if two or three ports are chosen, only the minimal subtree joining their
ends, because every component of the chosen forest needs at least two
chosen ports and at most three are available. These are exactly the
numbers of possibilities, $1,0,1,1$, at the center of the star, and the
conditions at the vertices of $U$ are the same in $G$ and in~$G_*$.
\end{proof}

The bound of three ports is sharp for this method and for the statement.
A single scalar star with four ports has
$P_{12}P_{34}=P_{13}P_{24}=P_{14}P_{23}$, which fails for a two-dimensional
vertex with $J_C=I$ and couplings $(1,0),(1,1),(1,2),(1,3)$, where the
three products are $7,4,3$. The next proposition shows that the degree
itself increases.

\begin{proposition}\label{prop:block-four-port}
Let $D_k$ be the graph with two vertices joined by $k$ parallel edges, so
that $L(D_k)=2^k-k$, with vertex dimensions $(1,p)$.
\begin{enumerate}
\item[\textup{(i)}] If $p\ge2$ and $k=a+b$ with $a,b\ge1$, then
$N_{(1,p)}(D_k)\ge L(D_k)+(2^a-a-1)(2^b-b-1)$; hence
$N_{(1,p)}(D_k)>L(D_k)$ for $k\ge4$.
\item[\textup{(ii)}] \textup{(Computer-assisted.)} For four parallel edges,
\[
N_{(1,p)}(D_4)=\begin{cases}12,&p=1,\\ 14,&p\ge2.\end{cases}
\]
\end{enumerate}
For three parallel edges $N_{(1,p)}(D_3)=5$ for all $p$, by
Theorem~\ref{thm:block-separator}.
\end{proposition}

\begin{proof}
(i) Let $H_{a,b}$ be the scalar graph with vertices $o,u,v$, an edge $uv$,
$a$ parallel edges $ou$ and $b$ parallel edges $ov$. Generic scalar
coefficients on $H_{a,b}$ define special coefficients on $D_k$ with
dimensions $(1,p)$: the two-by-two precision matrix of $\{u,v\}$, extended
by $I_{p-2}$, is the diagonal coefficient of the second vertex, and each
edge keeps its weight, placed in the coordinate of $u$ or of~$v$. In both
models the unknowns can be reduced to the $k$ scalar cavity precisions
entering the ports, and the reduced equations coincide, because the
response \eqref{blc:eq:tree-response} of the tree $uv$ and of the block
vertex is the same Schur complement. The eliminations are regular
substitutions near each solution, of the triangular form ``unknown equals a
function of the remaining ones'', so they preserve simple solutions; the
invertibility they require is an open condition that holds at the positive
solution of a weakly coupled real model, hence at all solutions for
generic scalar coefficients, because the scalar incidence variety
$Y_{\ones}(H_{a,b})$, an open subset of an affine space, is irreducible and
the condition therefore fails only on a proper closed subset of it. The $L(H_{a,b})$ simple solutions persist
under a perturbation to generic coefficients of $D_k$, so
$N_{(1,p)}(D_k)\ge L(H_{a,b})$. If the edge $uv$ is not chosen, a leafless
subset avoids choosing exactly one edge in either group:
$(2^a-a)(2^b-b)$ subsets; if it is chosen, at least one edge of each group
is needed: $(2^a-1)(2^b-1)$ subsets. The sum is
$2^k-k+(2^a-a-1)(2^b-b-1)$.

(ii) The value for $p=1$ is Theorem~\ref{thm:generic-degree}. Let $p\ge2$.
The symbols introduced now are used in this proof only. Let $d$ be the
diagonal coefficient of the scalar vertex, $\Theta\in\blcSym_p$ that of the
other vertex, generically invertible, $b_1,\dots,b_4\in\C^p$ the edge
coefficients, and $M_{ij}=b_i^{\T}\Theta^{-1}b_j$. Let $y_i$ be the cavity
precision from the scalar vertex into the $i$th edge. The precision of the
block vertex is $P_0=\Theta-\sum_ib_ib_i^{\T}/y_i$, its cavity precisions
are $P_0+b_ib_i^{\T}/y_i$, and the precision $R$ of the scalar vertex
satisfies $3R=\sum_iy_i-d$. Put $x_i=y_i-M_{ii}$, $c=d-\sum_iM_{ii}$,
$\sigma=\sum_ix_i-c$ and
$\mathcal H=\det(\diag(y)-M)$, a polynomial in $x$. By the
Sherman--Morrison formula,
$\partial_i\mathcal H/\mathcal H=1/y_i+b_i^{\T}P_0^{-1}b_i/y_i^2$, and the
cavity equation at the $i$th edge is equivalent to this quantity being
$1/R$. Hence the admissible solutions correspond, with multiplicities, to
the solutions of
\begin{equation}\label{blc:eq:four-port}
\sigma\,\partial_i\mathcal H=3\mathcal H\quad(1\le i\le4),\qquad
\mathcal H\,\sigma\prod_iy_i\ne0 ,
\end{equation}
a system in four unknowns that depends on the coefficients only through
the symmetric matrix $M$ and the number~$c$.

Write $u=x_3$, $v=x_4$, $\mu=M_{34}$, $\Delta=uv-\mu^2$, let $C$ be the
block $M_{\{1,2\},\{3,4\}}$,
$K=C\bigl(\begin{smallmatrix}v&\mu\\ \mu&u\end{smallmatrix}\bigr)C^{\T}$,
$\alpha=K_{11}+K_{22}$ and $\eta=-M_{12}\Delta-K_{12}$. Then
$\Delta\mathcal H=(\Delta x_1-K_{11})(\Delta x_2-K_{22})-\eta^2$. The
equations for $i=1,2$ give $\Delta x_1-K_{11}=\Delta x_2-K_{22}=:t$ and
\[
t^2-\bigl(\alpha+\Delta S\bigr)t-3\eta^2=0,\qquad S=u+v-c .
\]
With $\Lambda_w=2\alpha\,\partial_w\Delta+\Delta(S\,\partial_w\Delta-\Delta-\partial_w\alpha)$
and $\Pi_w=\Delta\,\partial_w\eta-2\eta\,\partial_w\Delta$ for
$w\in\{u,v\}$, the equations for $i=3,4$ become $\Lambda_wt=2\eta\,\Pi_w$.
On the open set $\Delta\eta\Lambda_u\ne0$ we have
$t=2\eta\Pi_u/\Lambda_u$, the solution is recovered from $(u,v)$, and
$(u,v)$ lies on the two plane curves
\[
\mathcal E=(\Lambda_u\Pi_v-\Lambda_v\Pi_u)/\Delta=0,\qquad
\mathcal F=4\eta\Pi_u^2-2(\alpha+\Delta S)\Pi_u\Lambda_u-3\eta\Lambda_u^2=0 ,
\]
where $\mathcal E$ is a polynomial. All solutions of a generic fiber lie
in this open set: the condition is open on the irreducible variety
$Y_{(1,p)}$, and it holds at the solutions of the example below.

By symbolic expansion, $\mathcal E$ and $\mathcal F$ have bidegrees at
most $(3,3)$ and $(4,7)$ in $(u,v)$, with equality and without common
component for generic coefficients. Their intersection number in
$\PP^1\times\PP^1$ is $3\cdot7+3\cdot4=33$. At least $19$ of it lies
outside the admissible set, for coefficients in a neighborhood of the
example below: at least $7$ at the four points $(\infty,0)$,
$(\infty,v_0)$, $(0,\infty)$, $(\infty,\infty)$ at infinity, where $v_0$ is
the root of the leading coefficient in $u$ (with multiplicities at least
$3,1,2,1$, read off from the supports of $\mathcal E$ and $\mathcal F$); at
least $8$ at the four transversal common zeros of $\Lambda_u$ and $\Pi_u$
with $\Delta\eta\ne0$, where $\mathcal E$ vanishes and
$\mathcal F\in(\Lambda_u,\Pi_u)^2$; and at least $4$ at the two points of
$\Delta=0$ where $\alpha=\pm2\eta$, at each of which $\mathcal E$ vanishes
and $\mathcal F$ vanishes to second order, since on $\Delta=0$ one has
$\alpha\pm2\eta=(\xi_1\mp\xi_2)^2/v$ with $\xi_i=C_{i1}v+C_{i2}\mu$ and
$\mathcal F=4\eta v^2(4\eta^2-\alpha^2)-4v\Delta[(4\eta^2+\alpha^2)\partial_u\eta-2\alpha\eta\,\partial_u\alpha]+O(\Delta^2)$.
Since this neighborhood meets the open set where the degree is generic,
$N_{(1,p)}(D_4)\le33-19=14$.

For the lower bound take $\Theta=I_2$, $d=17$ and
$(b_1\,b_2\,b_3\,b_4)=\bigl(\begin{smallmatrix}1&2&1&2\\ 0&1&3&2\end{smallmatrix}\bigr)$.
The resultant of $\mathcal E$ and $\mathcal F$ with respect to $u$ has an
irreducible factor of degree $14$ and multiplicity one; in its residue
field one recovers $u$, $t$, $x_1$, $x_2$ and checks
\eqref{blc:eq:four-port}, $\Delta\eta\Lambda_u\ne0$ and the nonvanishing of
the Jacobian of the four equations. The four zeros of
$(\Lambda_u,\Pi_u)$ and the two points on $\Delta=0$ used above are also
verified at this point. These are exact computations over $\mathbb Q$.
For $p>2$ take $\Theta=I_p$ and extend the $b_i$ by zeros, which gives the
same system \eqref{blc:eq:four-port}. So there are $14$ simple admissible
solutions, which persist to generic coefficients.
\end{proof}

The number $33=19+14$ is the only place in this paper where a block degree
beyond the range of Theorem~\ref{thm:block-separator} is obtained by an
intersection count.

\subsection{Trees and single cycles: selections and the hyperoctahedral group}
\label{blc:subsec:cycle}

For a cycle the block system is a periodic matrix Riccati recursion, and
its solutions are described by a symplectic transfer matrix. Let $C_n$,
$n\ge3$, have vertices $0,\dots,n-1$ (indices modulo $n$), uniform
dimension $p$, and edge coefficients $A_i=A_{i,i+1}\in\blcGL_p(\C)$; the
coefficient space is
$\mathcal P_n=\blcSym_p(\C)^n\times\blcGL_p(\C)^n\ni c$. Since every vertex
has degree two, \eqref{blc:eq:cavity} decouples into a recursion for the
clockwise cavity precisions $X_i=Q_{i\to i+1}$ and one for the
counterclockwise ones $X'_i=Q_{i+1\to i}$:
\begin{equation}\label{blc:eq:riccati}
X_{i+1}=D_{i+1}-zI-A_i^{\T}X_i^{-1}A_i,\qquad
X'_{i}=D_{i+1}-zI-A_{i+1}X_{i+1}'^{-1}A_{i+1}^{\T}.
\end{equation}
Let
\[
T_{i+1}(z)=\begin{pmatrix}(D_{i+1}-zI)A_i^{-1}&-A_i^{\T}\\ A_i^{-1}&0\end{pmatrix},
\qquad M_c(z)=T_0(z)\,T_{n-1}(z)\cdots T_1(z).
\]
Each $T_{i+1}$ preserves the standard symplectic form on $\C^{2p}$, so
$M_c(z)\in\blcSp(2p,\C[z])$ and its eigenvalues come in pairs
$\{\lambda,\lambda^{-1}\}$. Hence
$\det(\lambda I-M_c(z))=\lambda^pP_c(z,\lambda+\lambda^{-1})$ with
$P_c(z,x)$ monic of degree $p$ in $x$. We write $\Sigma_c$ for the
\emph{eigenvalue curve} $\{\det(\lambda I-M_c(z))=0\}$, $\Gamma_c$ for
$\{P_c(z,x)=0\}$, and
\[
B_\pm(z)=P_c(z,\pm2),\qquad\delta_c(z)=\operatorname{disc}_xP_c(z,x).
\]
The zeros of $B_\pm$, where $M_c$ has the eigenvalue $\pm1$, are the
\emph{band edges}, and the zeros of $\delta_c$ the \emph{collisions}. When
the eigenvalues are distinct and different from $\pm1$, a \emph{selection}
is a choice of one eigenvalue from each of the $p$ pairs. The permutations
of the $2p$ eigenvalues commuting with $\lambda\mapsto\lambda^{-1}$ form
the hyperoctahedral group $W(C_p)\cong(\Z/2)^p\rtimes S_p$, of order
$2^pp!$, the Weyl group of $\blcSp_{2p}$; it acts faithfully on the $2^p$
selections. Its normal subgroup $(\Z/2)^p$ is generated by the
\emph{flips} of single pairs, and a \emph{parallel swap} of two pairs
($\lambda_j\leftrightarrow\lambda_k$ and
$\lambda_j^{-1}\leftrightarrow\lambda_k^{-1}$) maps to a transposition in
$S_p$.

\begin{lemma}\label{blc:lem:lagrangian}
Write $\blcgraph(X)=\{(Xw;w):w\in\C^p\}$ and
$\Lambda_i=\{(A_iv;\,A_i^{-\T}X'_iv):v\in\C^p\}$. For invertible symmetric
$X_i,X'_i$ satisfying \eqref{blc:eq:riccati}, both are Lagrangian,
$T_{i+1}\blcgraph(X_i)=\blcgraph(X_{i+1})$ and $T_{i+1}\Lambda_i=\Lambda_{i+1}$;
in particular $\blcgraph(X_0)$ and $\Lambda_0$ are invariant under
$M_c(z)$. The edge precision $E_i$ of the edge $(i,i+1)$ is singular if and
only if $\blcgraph(X_i)\cap\Lambda_i\ne0$, and this holds for one $i$ if and
only if it holds for all.
\end{lemma}

\begin{proof}
With $w'=A_i^{-1}X_iw$ we get
$T_{i+1}(X_iw;w)=((D_{i+1}-zI)w'-A_i^{\T}X_i^{-1}A_iw';\,w')=(X_{i+1}w';w')$,
and
$T_{i+1}(A_iv;A_i^{-\T}X'_iv)=((D_{i+1}-zI-X'_i)v;\,v)=(A_{i+1}v'';A_{i+1}^{-\T}X'_{i+1}v'')$
with $v''=X_{i+1}'^{-1}A_{i+1}^{\T}v$, by \eqref{blc:eq:riccati}. The graph
of a matrix is Lagrangian exactly when the matrix is symmetric, and
$\Lambda_i$ is the image of such a graph under a symplectic map. If
$E_i(u;v)=0$ then $(X_iu;u)=(-A_iv;u)=(A_i(-v);A_i^{-\T}X'_i(-v))$ is a
common vector, and conversely. The last assertion holds because the
$T_{i+1}$ transport both subspaces.
\end{proof}

\begin{theorem}\label{thm:block-cycle}
\begin{enumerate}
\item[\textup{(i)}] \textup{(Degree.)} For a tree and arbitrary dimensions,
$N_{\blcp}(G)=1$. For a connected loopless multigraph of cycle rank one
and arbitrary dimensions, $N_{\blcp}(G)=2^{r}$, where $r$ is the smallest
dimension of a vertex on the cycle; the cycle may have length two.
\item[\textup{(ii)}] \textup{(Selections.)} Let $G=C_n$, $n\ge3$, with
uniform dimension $p$. For generic $c\in\mathcal P_n$ and all but finitely
many $z$, the admissible solutions are in bijection with the selections of
eigenvalues of $M_c(z)$: $\blcgraph(X_0)$ is spanned by the eigenvectors of
the selected eigenvalues and $\Lambda_0$ by those of the complementary
selection. The bijection is compatible with analytic continuation in $z$.
\item[\textup{(iii)}] \textup{(Group.)} In the situation of \textup{(ii)},
$\blcMon_z(C_n,p)$ is the full hyperoctahedral group $W(C_p)$ in its action
on the $2^p$ selections. The same holds for $C_n$ with trees attached to
its vertices, with $M_c$ replaced by the transfer matrix of the cycle with
the tree messages absorbed into the diagonals.
\item[\textup{(iv)}] \textup{(Branching and genus.)} For generic
$c\in\mathcal P_n$: $\deg B_\pm=np$, $\deg\delta_c=np(p-1)$, and
$B_+B_-\delta_c$ is squarefree. At a band edge exactly one pair passes
through $\pm1$ and the local monodromy is the flip of that pair; at a
collision exactly two pairs meet and the local monodromy is their parallel
swap; all three curves are unramified over $z=\infty$. The curves
$\Sigma_c$ and $\Gamma_c$ are irreducible with groups $W(C_p)$ and $S_p$,
and the genera of the normalizations are
\[
g(X_c)=1+2^{p-3}\bigl(np(p+3)-8\bigr),\qquad g(\Sigma_c)=np^2-2p+1,\qquad
g(\Gamma_c)=1-p+\tfrac12np(p-1),
\]
where $X_c$ is the block spectral curve.
\end{enumerate}
\end{theorem}

For $p=1$ the three genera are $n-1$, $n-1$, $0$, and $X_c=\Sigma_c$ is the
hyperelliptic curve $y^2=(\tr M_c)^2-4$. For $p\ge2$ the group is much
smaller than the symmetric group $S_{2^p}$, in contrast with
Remark~\ref{rem:monodromy-group}: the pairing of the eigenvalues restricts
it, and $W(C_p)$ is the largest group compatible with this pairing.

\begin{remark}[Attribution]\label{blc:rem:attribution}
Several ingredients of Theorem~\ref{thm:block-cycle} are classical or
known, and we do not claim them. That the symmetric solutions of a
(periodic) matrix Riccati equation correspond to invariant Lagrangian
subspaces of the associated Hamiltonian or symplectic (monodromy) matrix,
and that there are $2^p$ of them when the eigenvalues are distinct and off
the unit circle or imaginary axis, is standard in the theory of algebraic
Riccati equations; see \cite{LancasterRodman} and the account of the
parametrization of all solutions in \cite{AlpagoFerrante}. We have not
consulted the original sources for the discrete periodic case, and state
Lemma~\ref{blc:lem:lagrangian} with its short proof for this reason. The
curve $\Gamma_c$ is the Riemann surface of the Lyapunov function of the
periodic matrix-valued Jacobi operator of \cite{KorotyaevKutsenko}, who
show that it may be reducible but do not treat generic coefficients. The
irreducibility of the eigenvalue curve $\Sigma_c$ for generic coefficients
is contained in the theorem of \cite{FaustLiu} on Bloch varieties of
$\Z$-periodic graphs, and its genus $np^2-2p+1$ follows from the genus
formula of \cite{vMM} for regular periodic band matrices, applied with
period $np$ and bandwidth $p$ on either side after a triangularization of
the couplings. What is specific to the present setting is that for generic
coefficients every selection gives an admissible fixed point, the value
$2^r$ for mixed dimensions and the extension to attached trees, the
determination of the group $W(C_p)$, and the genus of the curve $X_c$ of
fixed points, which is a $2^p$-sheeted cover of the $z$-line and not the
eigenvalue curve.
\end{remark}

\begin{proof}[Proof of Theorem~\ref{thm:block-cycle}]
\emph{\textup{(ii)}.} Let $(*)$ be the following condition on
$(D-zI,A)$: the eigenvalues of $M=M_c(z)$ are distinct and different from
$\pm1$; the $2^p$ subspaces $L_S$ spanned by the eigenvectors of a
selection $S$, and their transports to the other vertices, are transversal
to $\C^p\times0$ and to $0\times\C^p$; and the cavity precisions so defined
are invertible. This is a Zariski open condition. Under $(*)$ the
invariant subspaces of $M$ are spanned by eigenvectors; since the
symplectic form pairs the eigenvectors of $\lambda$ and $\mu$ trivially
unless $\lambda\mu=1$, the invariant Lagrangian subspaces are exactly the
$L_S$, and $L_S\cap L_{S'}=0$ if and only if $S'$ is the complementary
selection. By Lemma~\ref{blc:lem:lagrangian}, the solutions of
\eqref{blc:eq:riccati} with invertible cavity precisions are the pairs
$(S,S')$, and since $\det E_e=\det Q^e_{j\to i}\det R_i$ ($R_i$ being the
Schur complement of $Q^e_{j\to i}$ in $E_e$), the admissible
ones are those with invertible edge precisions, that is, $S'$
complementary to $S$. On the open set $\Omega\subseteq\mathcal P_n\times\C$
where $(*)$ holds the eigenprojections are analytic, so the admissible
solutions form a covering of degree $2^p$ compatible with the labeling by
selections. $\Omega$ is nonempty: for diagonal $D_i$ and $A_i$, $M$ is the
direct sum of $p$ scalar transfer matrices $m_k(z)$; the entries
$(m_k)_{12}$, $(m_k)_{21}$ of the transfer matrix started at any vertex are
polynomials of degree $n-1$, $\tr m_k$ has degree $n$ and leading
coefficient $(-1)^n/\prod_ia^{(k)}_i$, and an eigenvector $(\alpha;\beta)$ of
a two-by-two matrix with nonzero off-diagonal entries has
$\alpha\beta\ne0$; so $(*)$ holds for all but finitely many $z$ as soon as
the products $\prod_ia^{(k)}_i$ are distinct. For $c$ with
$\Omega_c\ne\emptyset$ the set $\Omega_c$ is cofinite, and as in Step~4 of
the proof of Theorem~\ref{thm:block-irreducible} $X_c$ is the closure of
the cover over $\Omega_c$.

\emph{\textup{(i)}.} Let $j$ be a leaf with neighbor $i$. The cavity
precision from $j$ is $D_j$, so the model on $G-j$ with $D_i$ replaced by
$D_i-K$, $K=A_{ij}D_j^{-1}A_{ji}$, has the same cavity equations, and
$Q_{i\to j}=R_i+K$ with $R_i$ the vertex precision of the reduced model.
Keeping the reduced diagonal fixed and scaling $A_{ij}$ to zero shows that
$R_i+K$ can be made invertible at all solutions of a generic reduced
fiber; then the precision at $j$ and the edge precision are invertible by
$\det E=\det D_j\det R_i=\det(R_i+K)\det R_j$. Hence
$N_{\blcp}(G)=N_{\blcp}(G-j)$, and removing leaves gives $1$ for a tree
(a single vertex has the solution $V=D^{-1}$). For cycle rank one we are
reduced to a cycle. Let $v$ be a vertex of minimal dimension $r$ on it and
insert a vertex of dimension $r$ on each of the two edges at $v$, say $a$
and $b$; by Proposition~\ref{blc:prop:subdivision} the degree is unchanged.
The same graph is obtained from the triangle $v,a,b$ with uniform dimension
$r$ by subdividing the edge $ab$ repeatedly, inserting the remaining
vertices in their order starting from $a$, so that the edge subdivided
always has the end $b$ of dimension $r$ and the new vertex has dimension at
least $r$. Thus $N_{\blcp}(G)=N_r(C_3)$, which is $2^r$ by \textup{(ii)}.

\emph{\textup{(iii)}.} By \textup{(ii)} the group $G_c$ of a coefficient
$c$ with $\Omega_c\ne\emptyset$ is induced by the monodromy of the
eigenvalues, which commutes with inversion; so $G_c\subseteq W(C_p)$. If
loops $\gamma_1,\dots,\gamma_r\subset\Omega_{c_1}$ generate a subgroup $H$
of $G_{c_1}$, then for $c$ near $c_1$ the same loops lie in $\Omega_c$ and
induce the same permutations, since $\Omega$ is open and the cover over it
is locally trivial; so $\abs{G_c}\ge\abs H$. As the group is constant on a
nonempty Zariski open set, it suffices to exhibit one $c$ whose group
contains all flips and, for every $j<k$, an element mapping to the
transposition $(jk)$; these generate $W(C_p)$.

Start from the \emph{diagonal model} $c_0$:
$D_i=\diag_k(d^{(k)}_i+s_k)$, $A_i=\diag_k(a^{(k)}_i)$, with distinct
products $\prod_ia^{(k)}_i$ and shifts $s\in\C^p$. Then
$M_{c_0}(z)=\bigoplus_km_k(z-s_k)$ and the pairs are given by
$x_k(z)=t_k(z-s_k)$, $t_k=\tr m_k$. No $t_k^2-4$ is a square in $\C[z]$
(if $t^2-4=q^2$ then $(t-q)(t+q)=4$ forces $t$ to be constant), so it has a
zero $\beta_k$ of odd order, around which the two eigenvalues of channel
$k$ are exchanged. For $s$ outside finitely many hyperplanes
$\{s_k-s_j=\mathrm{const}\}$ the other channels are unbranched at
$\beta_k+s_k$, and a lasso around this point acts as the flip of the
$k$th pair. In particular $X_{c_0}$ is already irreducible, with group
$(\Z/2)^p$.

To obtain transpositions we split a collision. For $s$ outside finitely
many further hyperplanes, $t_j(z-s_j)-t_k(z-s_k)$ has $n$ simple zeros for
each $j<k$, and at such a zero $z^*$ the common value $x^*$ is different
from $\pm2$, the other $x_l(z^*)$ are different from $x^*$, from $\pm2$
and from each other, and $(m_j)_{21}(m_k)_{21}\ne0$. (Multiple zeros occur
only if $s_k-s_j$ is a difference $u-u'$ of a common solution of
$t_j(u)=t_k(u')$, $t_j'(u)=t_k'(u')$; these are finite in number, since
along a common component one would have $du=du'$ and then
$t_j(u)\equiv t_k(u-\sigma_0)$, contradicting the distinct leading
coefficients. The other conditions exclude finitely many values in the
same way.) Let $C$ be a symmetric matrix with nonzero off-diagonal entries
and let $c_\varepsilon$ be $c_0$ with $D_0$ replaced by $D_0+\varepsilon C$.
Then $M_{c_\varepsilon}(z)=(I+\varepsilon N)M_{c_0}(z)$ with
$N=\bigl(\begin{smallmatrix}0&C\\0&0\end{smallmatrix}\bigr)$. Near
$(z^*,0)$ the two eigenvalues close to $\lambda^*$ (where
$\lambda^*+1/\lambda^*=x^*$) are the eigenvalues of an analytic two-by-two
matrix $E(z,\varepsilon)$, the compression of $M_{c_\varepsilon}$ to the
range of the Riesz projection, with
$E(z,0)=\diag(\lambda_j(z),\lambda_k(z))$ and off-diagonal entries
$\varepsilon e_{12}$, $\varepsilon e_{21}$. First-order perturbation theory
gives $e_{12}e_{21}(z^*,0)=\lambda^{*2}C_{jk}^2\gamma_j\beta_j\gamma_k\beta_k$,
where $(\alpha_j;\beta_j)$ and $(\gamma_j,\delta_j)$ are right and left
eigenvectors of $m_j$; this is nonzero because $(m_j)_{21}\ne0$ forces
$\beta_j\gamma_j\ne0$. The discriminant of $E$ is
$(\lambda_j-\lambda_k+\varepsilon a)^2+4\varepsilon^2e_{12}e_{21}$ with $a$
analytic; in the variable $w=(z-z^*)/\varepsilon$ it converges, after
division by $\varepsilon^2$, to a quadratic polynomial with two distinct
roots, because $\lambda_j-\lambda_k$ has a simple zero at $z^*$. By
Hurwitz's theorem, for small $\varepsilon\ne0$ the discriminant has exactly
two simple zeros in a small disc around $z^*$ and no others there. A loop
around one of them exchanges the two eigenvalues near $\lambda^*$, and
simultaneously their inverses near $1/\lambda^*$, and fixes all other
eigenvalues: a parallel swap of the pairs $j$ and $k$. Conjugating by a
path from the base point, which for small $\varepsilon$ carries the same
labeling as for $c_0$, gives an element of $G_{c_\varepsilon}$ mapping to
$(jk)$. The flips persist for small $\varepsilon$ by the local constancy
above. Hence $G_{c_\varepsilon}=W(C_p)$.

For $C_n$ with trees attached, the sum $S_i(z)$ of the messages entering
the cycle vertex $i$ from its trees is a rational function of $z$
determined by the recursion from the leaves, independent of $D_i$, and the
argument applies to the transfer matrix with diagonals $D_i-zI-S_i(z)$:
the traces $t_k$ of the diagonal model are rational functions with the
same leading behavior at infinity, the scalar double covers are branched
at finite points by Theorem~\ref{thm:irreducible}, and the perturbation
$\varepsilon C$ of $D_0$ still multiplies the transfer matrix by
$I+\varepsilon N$. Here the dimension is uniform, and the shifts $s_k$ of
the diagonal model are applied to the diagonals of the tree vertices as
well, so that the model remains a direct sum of shifted scalar channels.
The differences $t_j(z-s_j)-t_k(z-s_k)$ are now rational functions; the
argument above uses one simple zero with the stated properties for each
pair $j<k$, not their number $n$.

\emph{\textup{(iv)}.} Since
$M_c(z)/(-z)^n\to\bigl(\begin{smallmatrix}\mathsf H&0\\0&0\end{smallmatrix}\bigr)$
with $\mathsf H=A_{n-1}^{-1}\cdots A_0^{-1}$, the large eigenvalues are
$\lambda_k\sim h_k(-z)^n$ with $h_k$ the eigenvalues of $\mathsf H$. Hence
$\deg B_\pm=np$ always, and $\deg\delta_c=np(p-1)$ when the $h_k$ are
distinct; in that case the eigenvalues of $M_c/(-z)^n$ and their inverses
are analytic in $1/z$ near $z=\infty$, so nothing ramifies there.
Squarefreeness of $B_+B_-\delta_c$ is a Zariski open condition; it holds at
$c_\varepsilon$ if the channels of $c_0$ are chosen with squarefree
$t_k^2-4$ (which holds for generic scalar coefficients) and generic
shifts: the band edges of $c_0$ are then simple and stay simple, and by
the splitting above each of the $n\binom p2$ double zeros of $\delta_{c_0}$
becomes two simple zeros, which accounts for the full degree $np(p-1)$.
At a simple zero of $B_\pm$ exactly one $x_k\mp2$ vanishes simply and
$\lambda_k=\frac12(x_k+\sqrt{x_k^2-4})$ has a square root branch point;
at a simple zero of $\delta_c$ exactly two roots $x_j,x_k$ are exchanged,
with $x\ne\pm2$, so the eigenvalues undergo a parallel swap. On the $2^p$
selections a flip is a product of $2^{p-1}$ transpositions and a parallel
swap a product of $2^{p-2}$ transpositions; there are $2np$ band edges and
$np(p-1)$ collisions. The Riemann--Hurwitz formula gives
$2g(X_c)-2=-2^{p+1}+2np\,2^{p-1}+np(p-1)2^{p-2}$,
$2g(\Sigma_c)-2=-4p+2np+2np(p-1)$ and $2g(\Gamma_c)-2=-2p+np(p-1)$. The
curves $\Sigma_c$ and $\Gamma_c$ are irreducible because $W(C_p)$ is
transitive on the $2p$ eigenvalues and $S_p$ on the $p$ pairs.
\end{proof}

The diagonal model is irreducible but has group $(\Z/2)^p$ only, so the
exceptional set in \textup{(iii)} is not empty. For the double edge
$D_2$, which is a cycle of length two, we have proved the degree $2^r$ but
not determined the group.

\subsection{The block Floquet identity and the points at infinity}
\label{blc:subsec:floquet}

Define, on admissible solutions,
\[
\Phi_B=\frac{\prod_{e\in E}\det E_e}{\prod_{i\in V}(\det R_i)^{d_i-1}},
\qquad
\omega=\Bigl(\sum_{i\in V}\tr V_i\Bigr)dz .
\]
For $p_i=1$, $\det E_e=(1+f_e)/(V_iV_j)$ in the variables of
\eqref{eq:system}, so $\Phi_B=\prod_e(1+f_e)/\prod_iV_i$ is the function of
Proposition~\ref{prop:floquet}.

\begin{proposition}[Block Floquet identity]\label{prop:block-floquet}
Let $G$ be a finite loopless multigraph with arbitrary vertex dimensions.
\begin{enumerate}
\item[\textup{(i)}] At an admissible solution,
$S_e=\bigl(\begin{smallmatrix}V_i&C_e\\ C_e^{\T}&V_j\end{smallmatrix}\bigr)$
with $C_e=-V_iA_{ij}(Q^e_{j\to i})^{-1}=-(Q^e_{i\to j})^{-1}A_{ij}V_j$,
and
$\sum_{e\ni i}Q^e_{i\to j}-(d_i-1)R_i=D_i-zI$. On the variety of admissible
solutions with varying coefficients,
\begin{equation}\label{blc:eq:floquet}
\frac{d\Phi_B}{\Phi_B}=\sum_i\tr(V_i\,dD_i)+2\sum_e\tr(C_e^{\T}dA_e)
-\Bigl(\sum_i\tr V_i\Bigr)dz .
\end{equation}
In particular $\omega=-d\log\Phi_B$ on every branch of the block spectral
curve, and for fixed coefficients and fixed $z$ the function $\Phi_B$ is
constant on every connected analytic family of admissible solutions.
\item[\textup{(ii)}] For arbitrary coefficients and every component $C$ of
the block spectral curve, with smooth compactification $\widehat C$:
$\Res_P\omega=-\ord_P\Phi_B\in\Z$ for all $P\in\widehat C$,
$\sum_{P\in\widehat C}\Res_P\omega=0$, and
$\ord_P\Phi_B=\tfrac12\bigl[\sum_{\text{arcs }a}\ord_P\det Q_a+\sum_i(2-d_i)\ord_P\det R_i\bigr]$.
\item[\textup{(iii)}] Let $G$ be connected and simple, with uniform
dimension $p$, real symmetric $D_i$ and real invertible $A_e$, let $\tH$
be the lift of $D+A$ to the universal cover and $\nu$ the sum over
$i\in V$ of the traces of the matrix-valued spectral measures of $\tH$ at
a lift of $i$ (total mass $np$). On the physical branch, for real
$z<\inf\sigma(\tH)$,
\[
\log\Phi_B(z)=\int\log(t-z)\,d\nu(t).
\]
\item[\textup{(iv)}] For uniform dimension $p$, invertible $A_e$ and
arbitrary coefficients, there is exactly one branch with
$Q_a/(-z)\to I$ on all arcs as $z\to\infty$; it is unramified at
$z=\infty$, has $\Res\omega=np$ there, and if $G$ is connected, for real
coefficients it is the physical branch.
\item[\textup{(v)}] For $C_n$ $(n\ge3)$ with uniform dimension $p$ and
generic coefficients,
$\Phi_B=\prod_i\det X_i=\prod_e\det A_e\cdot\prod_{\lambda\in S}\lambda$ on
the sheet of the selection $S$. All residues of $\omega$ at finite points
vanish. Over $z=\infty$ there are exactly $2^p$ points, unramified,
labeled by the set $S_-$ of pairs from which the small eigenvalue is
selected, with $\Res\omega=n(p-2\abs{S_-})$.
\end{enumerate}
\end{proposition}

Identity \eqref{blc:eq:floquet} is the matrix form of the stationarity of
the Gaussian Bethe free energy: the gradient of the Bethe log-determinant
$\log\Phi_B$ with respect to $D_i$ is the belief covariance $V_i$. It holds
at all admissible solutions, physical or not. For $p=1$,
\textup{(iii)} is the identification with the Floquet function of
\cite{BBGSS} in Proposition~\ref{prop:floquet}; we have not compared the
block statement with the literature. The argument for \textup{(iii)} uses only
Proposition~\ref{prop:block-physical} and the fact that the cover is a
tree.

\begin{proof}
(i) The first formulas are the block inversion of $E_e$, whose Schur
complements are $R_i$ and $R_j$, and the sum of
$Q^e_{i\to j}=R_i+A_{ij}(Q^e_{j\to i})^{-1}A_{ji}$ over the edges at $i$.
From $d\log\det E_e=\tr(S_e\,dE_e)$,
\[
d\log\det E_e=\tr(V_i\,dQ^e_{i\to j})+\tr(V_j\,dQ^e_{j\to i})+2\tr(C_e^{\T}dA_e).
\]
Summing over the edges, each arc occurs once, and
$\sum_{e\ni i}dQ^e_{i\to j}=(d_i-1)dR_i+dD_i-I\,dz$; the terms
$(d_i-1)\tr(V_i\,dR_i)=(d_i-1)\,d\log\det R_i$ cancel against the
denominator of $\Phi_B$. The last assertion is the case
$dD=dA=dz=0$.

(ii) On $C$ the function $\Phi_B$ is rational, finite and nonzero at
admissible points, so $\omega=-d\Phi_B/\Phi_B$ is a meromorphic
differential on $\widehat C$ with residues $-\ord_P\Phi_B$, whose sum
vanishes. The formula for the order follows from
$(\det E_e)^2=\det Q_{i\to j}\det Q_{j\to i}\det R_i\det R_j$.

(iii) On the physical branch $V_i$ is the diagonal block of
$(\tH-z)^{-1}$ and the inverses of the cavity precisions are diagonal
blocks of resolvents of restrictions of $\tH-z$ to half-trees
(Proposition~\ref{prop:block-physical}); hence all precisions are positive
definite, $Q_a\succeq(\inf\sigma(\tH)-z)I$, and $Q_a,R_i=-zI+O(1)$ as
$z\to-\infty$. By \textup{(i)},
$\frac d{dz}\log\Phi_B=-\sum_i\tr V_i=-\int(t-z)^{-1}d\nu(t)$, which is the
derivative of the right-hand side. Both sides are
$np\log(-z)+o(1)$ as $z\to-\infty$, because
$2p\abs E-p\sum_i(d_i-1)=np$ and $\nu$ has compact support and mass $np$.

(iv) Put $Q_a=-zP_a$ and $\tau=1/z$. Then \eqref{blc:eq:cavity} reads
$P^e_{i\to j}=I-\tau D_i-\tau^2\sum_{e'\ne e}A_{ik}(P^{e'}_{k\to i})^{-1}A_{ki}$,
which at $\tau=0$ has the solution $P=I$ with identity Jacobian. The
implicit function theorem gives a unique analytic solution $P(\tau)$,
admissible for small $\tau\ne0$, with $\sum_i\tr V_i=-np\,\tau+O(\tau^2)$,
so $\omega=np\,d\tau/\tau+\cdots$. The physical branch has this asymptotic
behavior by the proof of \textup{(iii)}.

(v) The Schur complement of the block $X_i=Q_{i\to i+1}$ in $E_i$ is
$R_{i+1}$, so $\det E_i=\det X_i\det R_{i+1}$, and $d_i=2$ gives
$\Phi_B=\prod_i\det X_i$. In the coordinate $w$ of
$\blcgraph(X_0)$ the restriction of $M_c$ is the ordered product of the
maps $w\mapsto A_i^{-1}X_iw$, whose determinant is
$\prod_i\det X_i/\prod_i\det A_i$ and also the product of the selected
eigenvalues. At finite $z$ the eigenvalues of the polynomial matrix $M_c$
and their inverses are bounded, so $\Phi_B$ is a unit at every point of
$\widehat C$ over a finite $z$. At infinity, by the proof of
Theorem~\ref{thm:block-cycle}\textup{(iv)}, the eigenvalues are
$\lambda_k=h_k(-z)^n(1+O(1/z))$ and their inverses, analytic in $1/z$, so
there are $2^p$ unramified points and
$\Phi_B=\mathrm{const}\cdot(-z)^{n(p-2\abs{S_-})}(1+O(1/z))$.
\end{proof}

For graphs of higher cycle rank we have no description of the points over
$z=\infty$. If at such a point each $Q_a$ has $p-k_a$ eigenvalues of the
order of $z$ and $k_a$ of the order of $z^{-1}$, each $\det R_i$ has the
order of $z^p$, and the ramification index is $e_P$, then
\textup{(ii)} gives $\Res_P\omega=e_P(np-\sum_ak_a)$, which for $p=1$ is
the residue $n-2\abs S$ of Lemma~\ref{lem:infinity}. Numerical continuation
for the bowtie, the diamond and $D_3$ with $p=2$ is consistent with this
splitting and with unramified points at infinity, but we have no proof; a
proof would turn the classification of the points at infinity into a
formula for $N_p(G)$, as in the scalar case.

\subsection{Poles and twisted operators}
\label{blc:subsec:twisted}

In this subsection $G$ is connected, the dimension is uniform, $p_i=p$,
the $D_i$ are complex symmetric and the $A_e$ are invertible. For $r\ge1$
and a family $O=(O_e)_{e\in E}$ of complex orthogonal matrices
$O_e\in\blcO(r,\C)$ define the \emph{twisted operator} $H_O$ on families
$U=(U_i)_{i\in V}$ of $p\times r$ matrices by
\[
(H_OU)_i=D_iU_i+\sum_{e=ij}A_{ij}U_jO_e+\sum_{e=ji}A_{ij}U_jO_e^{\T},
\]
where the first sum is over the edges oriented away from $i$ and the
second over those oriented towards $i$. It is a symmetric matrix of size
$npr$: the block Jacobi matrix twisted by the orthogonal local system with
holonomies $O_e$. A change of basis $U_i\mapsto U_ig_i$,
$g_i\in\blcO(r,\C)$, replaces $O_e$ by $g_j^{-1}O_eg_i$, so its spectrum
depends only on the conjugacy class of the holonomy representation of
$\pi_1(G)$, and its lift to the universal cover is $\tH\otimes I_r$. For
$r=1$, $O_e=s_e=\pm1$ and $H_O=D+A^s$ is a signed block Jacobi matrix.

Let $C$ be a component of the block spectral curve and $\widehat C$ its
smooth compactification. A point $P\in\widehat C$ over a finite
$z=\lambda$ is \emph{of type~I} if all cavity precisions are finite and
invertible at~$P$. Write
$T_{j\leftarrow i}=-(Q^e_{j\to i})^{-1}A_{ji}$ for $e$ joining $i$ and $j$.

\begin{theorem}\label{thm:block-twisted-poles}
\begin{enumerate}
\item[\textup{(i)}] At a point $P$ of type~I, the order
$m_P=\ord_P\det R_i\ge0$ and the dimension $r_P=\dim\ker R_i(P)$ do not
depend on $i$, $\ord_P\det E_e=m_P$ for all edges, $r_P\le m_P$, the
transfers $T_{j\leftarrow i}$ restrict to isomorphisms
$\ker R_i(P)\to\ker R_j(P)$ with inverse $T_{i\leftarrow j}$, and
\[
\Res_P\omega=m_P\,(\abs E-n).
\]
\item[\textup{(ii)}] Let $P$ be an unramified point of type~I at which
every $V_i$ has at most a simple pole and at least one has a pole, and let
$W_i=\Res_PV_i$. Then all $W_i$ have the same rank $r=r_P=m_P\ge1$, with
image $\ker R_i(P)$; writing $W_i=U_iU_i^{\T}$ with
$U_i\in\blcMat_{p\times r}$ of rank $r$, there are unique
$O_e\in\blcO(r,\C)$ with
$T_{j\leftarrow i}U_i=U_jO_e$ and $T_{i\leftarrow j}U_j=U_iO_e^{\T}$ for
$e=ij$, and
\[
H_OU=\lambda U,\qquad\sum_i\tr W_i=\sum_i\tr(U_i^{\T}U_i)=r\,(\abs E-n).
\]
\item[\textup{(iii)}] If $r=1$ in \textup{(ii)}, then $\lambda$ is an
eigenvalue of a signed matrix $D+A^s$, and for $\abs E\ne n$
\[
V_i(z)=\frac{(\abs E-n)\,u_iu_i^{\T}}{z-\lambda}+O(1),\qquad
(D+A^s)u=\lambda u,\quad\sum_iu_i^{\T}u_i=1 .
\]
\end{enumerate}
\end{theorem}

For $p=1$ every pole of type~I has $r=1$, and \textup{(iii)} is the
location and residue of Theorem~\ref{thm:local-pole}. For $p\ge2$ the
poles of type~I are stratified by $r\in\{1,\dots,p\}$: those with $r=1$ lie
at the finitely many signed eigenvalues, those with $r\ge2$ at eigenvalues
of operators twisted by a local system that is not known in advance and
varies in a continuous family. Theorem~\ref{thm:block-twisted-poles} gives
necessary conditions only. We do not prove the existence of such poles,
which for $p=1$ is Theorem~\ref{thm:all-signed-poles}, and we do not
classify the poles at which a cavity precision degenerates or diverges
(compare Theorem~\ref{thm:classification}). Numerical continuation for
$D_3$ with $p=2$ exhibits poles of both ranks, several poles of rank one
over the same signed eigenpair, and poles with negative residue, which by
\textup{(i)} are not of type~I when $\abs E>n$; these are observations, not
theorems. The role of these points for positive solutions, where the
relevant twists are unitary, is discussed with
Theorem~\ref{thm:block-boundary}.

\begin{proof}
(i) Since $\det Q_a$ is a unit at $P$, the identities
$\det E_e=\det Q_{j\to i}\det R_i=\det Q_{i\to j}\det R_j$ show that
$\ord_P\det R_i=\ord_P\det E_e=\ord_P\det R_j$ along every edge, hence
everywhere by connectedness. Then
$\ord_P\Phi_B=m_P\abs E-m_P\sum_i(d_i-1)=m_P(n-\abs E)$, and
Proposition~\ref{prop:block-floquet}\textup{(ii)} gives the residue. If
$E_e(P)(x;y)=0$, the second row gives $y=T_{j\leftarrow i}x$ and the first
then gives $R_ix=0$; symmetrically $x=T_{i\leftarrow j}y$ and $R_jy=0$.
So $x\mapsto(x;T_{j\leftarrow i}x)$ identifies $\ker R_i(P)$ with
$\ker E_e(P)$, and likewise for $j$, which gives the isomorphisms and the
independence of $r_P$. In a basis adapted to $\ker R_i(P)$, $r_P$ columns
of $R_i$ vanish at $P$, so $r_P\le m_P$.

(ii) Let $\zeta=z-\lambda$, $R_i=R_i^0+\zeta R_i^1+\cdots$ and
$V_i=W_i/\zeta+V_i^0+\cdots$. The coefficient of $\zeta^{-1}$ in
$V_iR_i=I$ gives $W_iR_i^0=0$, hence $\im W_i\subseteq K_i:=\ker R_i^0$ by
symmetry; the coefficient of $\zeta^0$, applied to $k\in K_i$, gives
$W_iR_i^1k=k$. So $\im W_i=K_i$, the bilinear form $R_i^1$ is nondegenerate
on $K_i$, and $\det R_i$ vanishes to order exactly $r_P=\dim K_i$; thus
$\rank W_i=r_P=m_P$, positive because some $W_i\ne0$. The symmetry of the
edge covariance in Proposition~\ref{prop:block-floquet}\textup{(i)} reads
$T_{j\leftarrow i}V_i=V_jT_{i\leftarrow j}^{\T}$; multiplying by $\zeta$,
\begin{equation}\label{blc:eq:residue-symmetry}
T_{j\leftarrow i}W_i=W_jT_{i\leftarrow j}^{\T}\qquad\text{at }P.
\end{equation}
A complex symmetric matrix of rank $r$ is $U_iU_i^{\T}$ with $U_i$ of rank
$r$, and the columns of $U_i$ span $K_i$. By \textup{(i)} there are unique
invertible $O_e,O_e'$ with $T_{j\leftarrow i}U_i=U_jO_e$ and
$T_{i\leftarrow j}U_j=U_iO_e'$. Substituting in
\eqref{blc:eq:residue-symmetry} gives $O_e=O_e'^{\T}$, and
$T_{i\leftarrow j}T_{j\leftarrow i}=\mathrm{id}$ on $K_i$ gives
$O_e'O_e=I$; so $O_e$ is orthogonal. Expanding $R_i(P)U_i=0$,
\[
0=(D_i-\lambda)U_i-\sum_{e\ni i}A_{ij}(Q^e_{j\to i})^{-1}A_{ji}U_i
=(D_i-\lambda)U_i+\sum_{e\ni i}A_{ij}\,T_{j\leftarrow i}U_i ,
\]
and $T_{j\leftarrow i}U_i$ is $U_jO_e$ or $U_jO_e^{\T}$ according to the
orientation of $e$; this is $H_OU=\lambda U$. Finally
$\sum_i\tr W_i=\Res_P\omega=m_P(\abs E-n)$ by \textup{(i)}.

(iii) $\blcO(1,\C)=\{\pm1\}$; put $u=U/\sqrt{\abs E-n}$.
\end{proof}

We have not found the twisted operators $H_O$ or the block identity
\eqref{blc:eq:floquet} in the literature on Gaussian belief propagation or
on matrix-valued Jacobi operators on trees that we consulted. Our search
of the prior art is incomplete here: we have not surveyed the formulas of
Ihara--Bass type for graph zeta functions twisted by representations of
the fundamental group, nor the work on non-backtracking operators with
matrix coefficients in the line of Bordenave and Collins, where related
objects may appear.

\begin{remark}[Open questions]\label{blc:rem:open}
\textup{(a)} For uniform dimension $p\ge2$ and a general graph we have no
formula for $N_p(G)$ of the type of $\LG$, that is, a sum over edge
subsets, or over other combinatorial data, of products of local weights.
Theorems~\ref{thm:block-degree}, \ref{thm:block-separator},
\ref{thm:block-cycle} and Proposition~\ref{prop:block-four-port} show that
such a formula must be supermultiplicative in $p$, reduce to $\LG$ when
the block vertices are separated by scalar ones with at most three ports,
give $2^p$ for a cycle, $36$ for the bowtie with $p=2$, and $14$ rather
than $12$ for four parallel edges with dimensions $(1,p)$, $p\ge2$. The
exact values of $N_3$ for the bowtie and of $N_2$ for the diamond are not
known; only the lower bounds $1079$ and $66$ are. The formulas of
Propositions~\ref{blc:prop:rank-formula} and~\ref{blc:prop:multinomial}
determine $N_{\blcp}(G)$ in principle but give no insight.

\textup{(b)} The isomorphism type of the group $\blcMon_z$ of the $36$
sheets of the bowtie with $p=2$ is not known. It is transitive, and the
scalar $S_4$ result of Proposition~\ref{irr:prop:small-monodromy} gives a
subgroup acting as $S_4\times S_4$ on a set of $16$ sheets by
Proposition~\ref{blc:prop:product-group}.

\textup{(c)} The points of the block spectral curve over $z=\infty$, the
existence of the poles of Theorem~\ref{thm:block-twisted-poles}, and the
poles that are not of type~I are not understood beyond the single cycle.
\end{remark}

%% file: sections/formalization.tex
\section{Formal verification}
\label{sec:formal}

Part of the scalar theory has been formalized in Lean~4 \cite{Lean} with
Mathlib \cite{Mathlib}. Table~\ref{tab:lean} gives representative
declarations in the namespace \path{BetheGeometry.Gaussian}.
The Gaussian declarations are audited at build time for axiom
dependencies. The audit allows \texttt{propext}, \texttt{Classical.choice}
and \texttt{Quot.sound}.

The Lean source is available in the
\href{https://github.com/ywatanabex/gaussian-bp-jacobi-lean}{GitHub repository
\texttt{ywatanabex/gaussian-bp-jacobi-lean}}.
This paper refers to commit
\href{https://github.com/ywatanabex/gaussian-bp-jacobi-lean/commit/18e6de182b6d2601c36167a20ba86eaca9ba9d1d}{\texttt{18e6de1}}.

\subsection{Graph and solution conventions}

A multigraph has an edge type with two endpoint maps, so parallel edges
are distinct. The graph convention counts each self-loop twice in the
vertex sum. The universal cover is the type of reduced walks from a
root. Its Jacobi operator acts on the real Hilbert space $\ell^2$ of this
type, and its spectrum is independent of the root.

Lean's \path{IsAdmissible} requires the covariance equations and
$V_i\ne0$ for every vertex. The fixed-coefficient count and the covering
construction exclude a finite set of energies that contains the $D_i$.
On this domain the admissibility condition agrees with
Definition~\ref{def:system}. For local positivity Lean uses $V_i>0$ and
$f_e\ge0$. For nonzero real edge weights the edge equation then gives
$f_e>0$, as used in the manuscript. The comparison matrix uses the
edgewise absolute values of the weights.

\subsection{Verified statements}

\begin{table}[htbp]
\centering
\small
\begin{tabular}{@{}>{\raggedright\arraybackslash}p{0.43\textwidth}>{\raggedright\arraybackslash}p{0.52\textwidth}@{}}
\toprule
Statement & Lean declaration \\
\midrule
Strict spectral gaps at both Perron values for connected multigraphs
with $\abs E>n$ and nonzero real weights
& \path{comparisonEigenvalue_lt_coverSpectrumInf_of_connected},
  \path{coverSpectrumSup_lt_comparisonEigenvalue_of_connected} \\
Continuation of a single physical germ to both Perron poles, with
simple pole orders and the vertex residues
& \path{Paper.perron_analytic_continuation} \\
Analytic normalization of the physical projective component, carrying
both real Perron arcs
& \path{Paper.physical_projective_normalization} \\
Two positive solutions throughout the lower gap
& \path{Paper.two_lower_gap_solutions} \\
The generic complex count and the count at fixed complex coefficients
away from finitely many energies, both with simple roots
& \path{Paper.generic_solution_count},
  \path{Paper.fixed_coefficient_solution_count} \\
Existence of the local signed pole branch with the specified vertex
and edge limits
& \path{Paper.signed_eigenvalue_local_pole} \\
Uniqueness under walk-summability and its generic converse
& \path{Paper.walk_summability_implies_uniqueness},
  \path{Paper.generic_uniqueness_iff_walk_summability} \\
The signed count at nearby regular values, the Morse index sum and
parity on nondegenerate slices
& \path{Paper.local_degree},
  \path{Paper.morse_index_sum}, \path{Paper.root_count_parity} \\
At most two positive solutions per energy and at most one nonpositive
Hessian eigenvalue under the packing hypothesis
& \path{Paper.packing_positive_root_bound},
  \path{Paper.packing_hessian_bound} \\
\bottomrule
\end{tabular}
\caption{Selected formal statements for the scalar theory.}
\label{tab:lean}
\end{table}

The continuation of Theorem~\ref{thm:main-perron}(b) is expressed by a
path of locally agreeing holomorphic germs. It starts at the
cover Green function and ends at a germ in the explicit pole chart.
The normalization is constructed by filling the punctures of the
finite covering, including those over infinity. The projective map
from the physical component has a compact smooth source, is proper
with finite fibers, and admits a unique local lift of each nonconstant
holomorphic germ in its image.

For Theorem~\ref{thm:main-sheets}, the formal count uses the Jacobian of
the reduced edge system. The Schur complement in
Lemma~\ref{deg:lem-reduction} relates its invertibility to that of the
original covariance system. The fixed-coefficient theorem counts all
admissible points outside its finite exceptional set.

For the local pole theorem, the formal existence result assumes
$\abs E>n$, a one-dimensional kernel spanned by $u$, $u_i\ne0$ for every
$i$, and $u^{\T}u\ne0$. The chart is
$(z,V,f)=(\lambda+t,W(t)/t,(Y(t)/t-1)/2)$ with $W$ and $Y$ analytic
at zero. It gives
\[
W_i(0)=\frac{\kappa_0u_i^2}{u^{\T}u},\qquad
Y_e(0)=-\frac{2\kappa_0s_e w_e u_i u_j}{u^{\T}u}\quad(e=ij).
\]
A separate real statement gives local uniqueness in the scaled
coordinates when the kernel vector is positive, as needed for the
Perron arc.

The formal degree identity uses the gradient in logarithmic variance
coordinates. For the analogue of condition~(G), Lean tests the
principal blocks whose off-diagonal support is connected. This agrees
with the graph condition for loopless multigraphs with nonzero weights.
Nearby regular values give the signed count even when the central
fiber is singular. At zero, the Hessian sign sum and root parity are
verified under the nondegeneracy hypothesis of
Theorem~\ref{thm:main-real}(b).

\subsection{Scope}

The analytic normalization concerns the physical irreducible component.
For Theorem~\ref{real:thm:packing}(a), \path{Paper.packing_positive_root_bound}
verifies that $p_U(G)\le1$ gives at most two positive solutions at every
energy for arbitrary real coefficients with nonzero edge weights. The
Lean proof orders the cavity precisions by applying the packing sign
lemma to the difference of two fixed points and rules out a chain of
three. Separately, \path{Paper.packing_hessian_bound} verifies that the
normalized Hessian has at most one nonpositive eigenvalue, counted with
multiplicity. The branch and circle classification uses the written proofs.

Generic irreducibility of the full spectral curve, the function-field
degree of individual Green functions, the classification of all residue
poles, the cycle-rank-one band-edge expansion, the genus formula, and
the vector-valued theory are outside the formalization described here.

%% file: sections/open-problems.tex
\section{Open problems}
\label{sec:open}

\begin{conjecture}[Full symmetric monodromy]\label{open:monodromy}
For every finite connected loopless multigraph $G$, the monodromy group
of the scalar spectral curve over the $z$-line is $S_{\LG}$ on its
$\LG$ sheets for generic coefficients $(D,w)$.
\end{conjecture}

Theorem~\ref{thm:main-sheets}(b) gives transitivity, which proves the
conjecture when $\LG\le2$. Proposition~\ref{irr:prop:small-monodromy}
and Remark~\ref{rem:monodromy-group} prove it in ten further examples.

\begin{problem}[Exceptional coefficients]\label{open:exceptional}
Theorems~\ref{thm:main-sheets}(b),(c) and~\ref{thm:main-poles}(b),(c) hold
for generic coefficients, and we do not describe the exceptional set. For
uniform coefficients on an edge-transitive graph the physical component
has two sheets (Proposition~\ref{prop:uniform}). Determine the number of
sheets of the physical component, and which signed eigenvalues are poles
on it, in terms of the symmetries of $(G,D,w)$.
\end{problem}

\begin{problem}[Half-tree poles]\label{open:half-tree}
Describe, for a general finite graph, how the poles of the half-tree
functions $m_{i|e}$ on all sheets are distributed over the energies, in
the sense of \cite[Problem~10.3]{ABS}. Section~\ref{sec:abs} gives a local
criterion and shows that the correspondence valid on the line fails.
\end{problem}

\begin{problem}[Gaps]\label{open:gaps}
Theorem~\ref{thm:gapless} shows that for regular nonbipartite graphs no
gap opens near the uniform coefficients. Does every connected nonbipartite
graph of minimum degree at least three admit a nonempty open set of
coefficients for which the spectrum of $\tH$ is an interval?
\end{problem}

\begin{problem}[Degree of block models]\label{open:block-degree}
Find a formula for the generic block degree at uniform dimension $p\ge2$
that, like $\LG$, sums products of local weights over combinatorial
choices. Determine the permutation group of the $36$ sheets of the
two-dimensional model on two triangles sharing a vertex
(Remark~\ref{blc:rem:open}).
\end{problem}

%% file: appendices/many-solutions.tex
\section{Families with many positive solutions}\label{app:many}

This appendix contains the proofs of three statements of
Section~\ref{real:subsec:degree} and Example~\ref{real:ex:diamond-chain}:
a bridge strength that works for diamond chains of every length
(Theorem~\ref{app:many:thm-chain}), the unboundedness of the number of
positive solutions for pseudoforest deletion number two
(Theorem~\ref{app:many:thm-unbounded}, which contains
Proposition~\ref{real:prop:unbounded}), and the two certified small
instances (Proposition~\ref{app:many:prop-certified}). As in
Section~\ref{sec:real}, the spectral parameter is absorbed into the
potential, so that we study the positive solutions of \eqref{eq:system}
with $z=0$; all graphs are finite and simple, $D\in\R^V$ and
$w\in(\R\setminus\{0\})^E$. We work with the cavity precisions $q_{i\to j}$
of Lemma~\ref{lem:bp-correspondence}, indexed by the arcs of $G$, and with
the precision update $\mathcal U$ of \eqref{real:eq:update}. For a vector
$q\in(0,\infty)^{\vec E}$ we call
\[
R_i(q)=D_i-\sum_{k\sim i}\frac{w_{ik}^2}{q_{k\to i}}\qquad(i\in V)
\]
the \emph{vertex precisions}. Norms of vectors and matrices are maximum
norms and maximal absolute row sums, $\norm\cdot_\infty$, unless stated
otherwise.

\subsection{Cavity precisions of positive solutions}
\label{app:many:subsec-cavity}

\begin{lemma}\label{app:many:lem-cavity}
Let $G$ be a connected simple graph with real coefficients and let
$q\in(0,\infty)^{\vec E}$ be a fixed point of $\mathcal U$. Then for every
edge $ij$
\begin{equation}\label{app:many:eq-det}
q_{i\to j}q_{j\to i}-w_{ij}^2=R_i(q)\,q_{j\to i}=R_j(q)\,q_{i\to j},
\end{equation}
so all $R_i(q)$ have the same sign. If one of them is positive, then
$V_i=1/R_i(q)$, $f_{ij}=w_{ij}^2V_i/q_{j\to i}$ is a positive solution of
\eqref{eq:system} with $z=0$, and at this solution the numbers
$u_{ij}=f_{ij}/(1+f_{ij})$ of Lemma~\ref{real:lem:phi}\textup{(b)} are
$u_{ij}=w_{ij}^2/(q_{i\to j}q_{j\to i})$. Every positive solution arises in
this way from exactly one such~$q$.
\end{lemma}

\begin{proof}
The fixed point equation for the arc $i\to j$ reads
$q_{i\to j}=R_i(q)+w_{ij}^2/q_{j\to i}$; multiplying by $q_{j\to i}$ gives the
first equality in \eqref{app:many:eq-det}, and exchanging $i$ and $j$ the
second. Since $q>0$, $R_i$ and $R_j$ have the same sign for adjacent $i,j$,
and $G$ is connected. If the $R_i$ are positive, the vertex precisions and,
by \eqref{app:many:eq-det}, the edge determinants are positive, so all
vertex and edge beliefs of $q$ are positive definite, and
Lemma~\ref{lem:bp-correspondence} gives a positive solution with the stated
$V$ and $f$. Its inverse map $q_{i\to j}=(1+f_{ij})/V_i$ gives
$u_{ij}=f_{ij}/(1+f_{ij})=w_{ij}^2V_i/(q_{j\to i}q_{i\to j}V_i)$. Conversely,
by the same lemma a positive solution comes from a unique fixed point
whose edge beliefs are positive definite; their diagonal entries
$q_{i\to j}$ are positive, and $R_i=1/V_i>0$.
\end{proof}

\subsection{The diamond chain}\label{app:many:subsec-chain}

Let $\Delta$ be the diamond of Example~\ref{real:ex:diamond-chain}: vertices
$0,1,2,3$, edges $01,02,03,12,13$, all weights $1$, and
$D=(\tfrac{73}{15},\tfrac{73}{15},1,1)$. We call $0,1$ the \emph{hubs} and
$2,3$ the \emph{terminals}. Put $t_+=\tfrac53$, $t_-=\tfrac65$, the roots of
$t^2-\tfrac{43}{15}t+2=0$, and for $t\in\{t_+,t_-\}$ let $q^t$ be the vector
\[
q^t_{0\to1}=q^t_{1\to0}=t,\qquad q^t_{h\to s}=t+1,\qquad
q^t_{s\to h}=\frac t{t+1}\qquad(h\in\{0,1\},\ s\in\{2,3\}).
\]

\begin{lemma}\label{app:many:lem-diamond}
For $t\in\{t_+,t_-\}$:
\begin{enumerate}
\item[\textup{(a)}] $q^t$ is a fixed point of $\mathcal U$, with vertex
precisions $R_0=R_1=t-1/t$ and $R_2=R_3=(t-1)/(t+1)$, that is $\tfrac{16}{15}$
and $\tfrac14$ for $t_+$, $\tfrac{11}{30}$ and $\tfrac1{11}$ for $t_-$. It
gives a positive solution with $u_{01}=t^{-2}$ and $u_{hs}=t^{-1}$.
\item[\textup{(b)}] The matrix $H_u$ of \eqref{real:eq:Hu} at this solution has the
eigenvalues $(t^2+1)/(t^2-1)$ (eigenvector $e_2-e_3$),
$(t^2+2)/(t^2-1)$ (eigenvector $e_0-e_1$) and the two eigenvalues of its
restriction to the span of $e_0+e_1$, $e_2+e_3$. These are
$\tfrac{17}8,\tfrac{43}{16},\tfrac{271\pm3\sqrt{7261}}{136}$ for $t_+$ and
$\tfrac{61}{11},\tfrac{86}{11},\tfrac{7167\pm5\sqrt{2146321}}{1342}$ for $t_-$.
So the solution of $t_+$ has Morse index $0$, that of $t_-$ has Morse index
$1$, and every eigenvalue has absolute value greater than $\tfrac1{10}$.
\item[\textup{(c)}] Let $C^t=C_{q^t}$ be the derivative of $\mathcal U$ at $q^t$
(Section~\ref{real:subsec:reduction}), a $10\times10$ matrix, and let
$\eta^t_s\in\R^{\vec E(\Delta)}$, for a terminal $s$, be the row vector with
entries $(t+1)^{-2}$ at the two arcs $h\to s$ and $0$ elsewhere. Then
$I-C^t$ is invertible and
\[
\norm{(I-C^t)^{-1}}_\infty=\frac{2025}{119},\ \frac{3355565}{201971},
\qquad 1+\norm{\eta^t_s(I-C^t)^{-1}}_1=\frac{20017}{3808},\ \frac{366717}{51667}
\]
for $t=t_+,t_-$ respectively (for both terminals $s$). In particular all
four numbers are at most $K:=\tfrac{2025}{119}<17.02$.
\end{enumerate}
\end{lemma}

\begin{proof}
(a) For an arc $h\to s$ from a hub to a terminal,
$\mathcal U(q^t)_{h\to s}=\tfrac{73}{15}-\tfrac1t-\tfrac{t+1}t$, and for the arc
$0\to1$, $\mathcal U(q^t)_{0\to1}=\tfrac{73}{15}-\tfrac{2(t+1)}t$; by
$t+2+2/t=\tfrac{73}{15}$ these equal $t+1$ and $t$. For $s\to h$,
$\mathcal U(q^t)_{s\to h}=1-1/(t+1)=t/(t+1)$. The vertex precisions follow
from $R_i=q_{i\to j}-w_{ij}^2/q_{j\to i}$, they are positive, and
Lemma~\ref{app:many:lem-cavity} gives the positive solution and $u$.

(b) With $\varphi_1(u)=u/(1-u^2)$ and $\varphi_2(u)=u^2/(1-u^2)$, the matrix
$H_u$ has diagonal entries $1+\varphi_2(t^{-2})+2\varphi_2(t^{-1})$ at the hubs
and $1+2\varphi_2(t^{-1})$ at the terminals, the entry $-\varphi_1(t^{-2})$
at $01$, $-\varphi_1(t^{-1})$ at the four hub--terminal pairs, and zero at
$23$. It commutes with the exchanges $0\leftrightarrow1$ and
$2\leftrightarrow3$, which gives the eigenvectors $e_2-e_3$, $e_0-e_1$ and
the invariant plane. The trace and determinant on that plane give the
remaining two eigenvalues. Substituting the two values of $t$ gives the
stated list. The bounds follow from $(271-13.6)^2>9\cdot7261$ and
$25\cdot2146321>(7167+134.2)^2$.
The eigenvalues are nonzero, so the solutions are nondegenerate.

(c) Index the two hubs and the two terminals separately by signs
$h,s\in\{\pm1\}$, with bars denoting sign reversal. Write the unknown
perturbations on the arcs $h\to\bar h$, $h\to s$, $s\to h$ as $A_h,B_{hs},C_{hs}$
and their right-hand sides in $(I-C^t)(A,B,C)=(x,y,z)$ as $x_h,y_{hs},z_{hs}$.
Put $j=t^{-2}$, $k=(t+1)^2/t^2$, $\ell=(t+1)^{-2}$, so $k\ell=j$. The equations are
\[
 A_h-k\sum_s C_{hs}=x_h,\qquad
 B_{hs}-j A_{\bar h}-kC_{h\bar s}=y_{hs},\qquad
 C_{hs}-\ell B_{\bar h s}=z_{hs}.
\]
Expand each array in the characters of the two sign exchanges, with
Fourier coefficients given by averaging against the character. A mode
$(\epsilon,\delta)$ has signs $\epsilon,\delta\in\{\pm1\}$ under hub and
terminal exchange; the $A,x$ arrays only have $\delta=+1$. In that mode,
write $A,B,C,x,y,z$ for the coefficients. Eliminating $C=z+\ell\epsilon B$
and, when $\delta=+1$, $A=x+2kC$, gives
\[
 B_{\epsilon,+}=
 \frac{y+j\epsilon x+k(1+2j\epsilon)z}{1-j\epsilon-2j^2},
 \qquad
 B_{\epsilon,-}=\frac{y-kz}{1+j\epsilon}.
\]
Every denominator is nonzero at the two stated values of $t$.

For clarity, the row sums from these two formulas are listed below. A
row is determined by its arc type, and the last row uses the functional
$\ell(B_{+s}+B_{-s})$ for a fixed terminal. To recover its coefficient
on any entry of an input array, sum the modal coefficients with the
signs of that entry's characters and divide by the array size (two for
$x$, four for $y,z$). In the last row only the modes with $\epsilon=+1$
survive and their coefficients are multiplied by $2\ell$.
\[
\begin{array}{c|cc}
 &t=5/3&t=6/5\\ \hline
 A &2025/119&2954934/201971\\
 B &1801/119&3355565/201971\\
 C &23825/7616&104848866/24438491\\
 1+\|\eta^t_s(I-C^t)^{-1}\|_1&20017/3808&366717/51667
\end{array}
\]
For example, the $A$ row at $t=5/3$ has absolute sums
$200/119$, $225/119$ and $1600/119$ on $x,y,z$, respectively.
Taking the maximum over the first three rows gives the stated inverse
norm.
\end{proof}

For $m\ge1$ and $\varepsilon>0$ let $\Delta_m(\varepsilon)$ be the chain of
Example~\ref{real:ex:diamond-chain}: copies $\Delta^{(1)},\dots,\Delta^{(m)}$
of $\Delta$, with vertices $0_j,1_j,2_j,3_j$, and the \emph{bridges}
$3_j2_{j+1}$ of weight $\varepsilon$ for $1\le j<m$. It has $4m$ vertices
and maximal degree~$3$. For $\tau\in\{+,-\}^m$ let $q^\tau$ be the vector on
the arcs of $\Delta_m(\varepsilon)$ equal to $q^{t_{\tau_j}}$ on the arcs of
$\Delta^{(j)}$, and equal on a bridge arc $s\to s'$ to the terminal
precision $(t-1)/(t+1)$ of the copy that contains~$s$, where $t=t_{\tau_j}$.

\begin{theorem}\label{app:many:thm-chain}
Let $m\ge1$ and $0<\varepsilon\le\tfrac1{1000}$. For every
$\tau\in\{+,-\}^m$ the model $\Delta_m(\varepsilon)$ has a nondegenerate
positive solution whose cavity precisions $q$ satisfy
$\norm{q-q^\tau}_\infty\le375\,\varepsilon^2$ and whose Morse index is the
number of $j$ with $\tau_j=-$. The $2^m$ solutions obtained in this way are
pairwise distinct. Consequently $\Delta_m(\varepsilon)$ has at least
$\binom mk$ nondegenerate positive solutions of Morse index $k$ for
$0\le k\le m$, and $\taups(\Delta_m(\varepsilon))=m$.
\end{theorem}

The constants do not depend on $m$, on $\tau$ or on $\varepsilon$ in this
range; only the existence statement for $m$-dependent $\varepsilon_m$ is
proved in Example~\ref{real:ex:diamond-chain}. Floating point computations
of the solutions for $m\le4$ suggest that $\varepsilon=\tfrac1{100}$ also
works, with all eigenvalues of $H_u$ of absolute value at least $0.106$;
this is evidence only.

\begin{proof}
Fix $m$, $\tau$ and $\varepsilon$, and write $t_j=t_{\tau_j}$. Let $N$ be the
set of arcs, split into the internal arcs (those of the copies) and the
bridge arcs, and let $\mathcal R^\varepsilon(q)=q-\mathcal U(q)$ on
$(0,\infty)^N$, where $\mathcal U$ is the update of
$\Delta_m(\varepsilon)$. The fixed points of $\mathcal U$ are the zeros of
$\mathcal R^\varepsilon$.

\emph{Step 1: the linearization.} Let $\mathcal R^0$ be the same map with
$\varepsilon=0$, that is with the bridge terms $\varepsilon^2/q$ deleted
from the update. Then $\mathcal R^0(q^\tau)=0$: on the internal arcs this is
Lemma~\ref{app:many:lem-diamond}(a), and for a bridge arc $3_j\to2_{j+1}$
the update is $1-1/q_{0_j\to3_j}-1/q_{1_j\to3_j}$, the terminal precision
of $\Delta^{(j)}$ (similarly for $2_{j+1}\to3_j$). Let
$L_0=D\mathcal R^0(q^\tau)$. In $\mathcal R^0$ the internal arcs of a copy
depend only on internal arcs of the same copy, and a bridge arc out of a
terminal $s$ depends on itself (coefficient $1$) and on the two internal
arcs $h\to s$, with coefficients $-(t_j+1)^{-2}$. Ordering the internal arcs
first,
\[
L_0=\begin{pmatrix}\bigoplus_j(I-C^{t_j})&0\\ -\mathcal H&I\end{pmatrix},
\qquad
L_0^{-1}=\begin{pmatrix}\bigoplus_j(I-C^{t_j})^{-1}&0\\
\mathcal H\bigoplus_j(I-C^{t_j})^{-1}&I\end{pmatrix},
\]
where each row of $\mathcal H$ is a vector $\eta^{t_j}_s$ of
Lemma~\ref{app:many:lem-diamond}(c) placed on the arcs of one copy. By that
lemma, $\norm{L_0^{-1}}_\infty\le K=\tfrac{2025}{119}$, for every $m$ and
$\tau$.

\emph{Step 2: estimates.} The map $\mathcal R^\varepsilon$ differs from
$\mathcal R^0$ only on the four internal arcs out of a bridged terminal: for
example
$\mathcal R^\varepsilon(q)_{3_j\to0_j}=\mathcal R^0(q)_{3_j\to0_j}+\varepsilon^2/q_{2_{j+1}\to3_j}$.
The bridge arcs have $q^\tau$-values $\tfrac14$ or $\tfrac1{11}$, so
\begin{equation}\label{app:many:eq-res}
\norm{\mathcal R^\varepsilon(q^\tau)}_\infty\le11\varepsilon^2 .
\end{equation}
Let $0<r\le\tfrac1{2500}$ and $B_r=\{q:\norm{q-q^\tau}_\infty\le r\}$. The
internal arcs have $q^\tau$-values at least $\tfrac6{11}$, so on $B_r$ every
internal arc has $q\ge\tfrac6{11}-r\ge0.545$ and every bridge arc has
$q\ge\tfrac1{11}-r\ge0.0905$; in particular $B_r\subset(0,\infty)^N$. The
nonzero off-diagonal entries of $D\mathcal R^\varepsilon(\xi)$ are
$-1/\xi_b^2$ for internal arcs $b$ and $-\varepsilon^2/\xi_b^2$ for bridge
arcs $b$; those of $L_0$ are $-1/(q^\tau_b)^2$ for internal $b$. For
$\abs{\xi_b-c}\le r$ and $c\ge\tfrac6{11}$,
$\abs{\xi_b^{-2}-c^{-2}}\le2r/(c-r)^3\le12.36\,r$. Each row has at most two
entries of the first kind, or one of the first kind and one of the second
(a terminal has at most three neighbors, one of which is the head of the
arc). Hence for $\xi\in B_r$
\begin{equation}\label{app:many:eq-lip}
\norm{D\mathcal R^\varepsilon(\xi)-L_0}_\infty\le\delta(r):=
\max\bigl(24.72\,r,\ 12.36\,r+122.1\,\varepsilon^2\bigr),
\end{equation}
using $(\tfrac1{11}-r)^{-2}\le122.1$.

\emph{Step 3: contraction.} Put $\rho=22K\varepsilon^2$. Since
$22K=\tfrac{44550}{119}<374.4$ and $\varepsilon\le10^{-3}$, we have
$\rho\le3.744\cdot10^{-4}<\tfrac1{2500}$ and
$122.1\,\varepsilon^2\le0.33\,\rho$, so $\delta(\rho)\le24.72\,\rho$ and
$K\delta(\rho)\le17.02\cdot24.72\cdot3.744\cdot10^{-4}<\tfrac12$. Let
$\Gamma(q)=q-L_0^{-1}\mathcal R^\varepsilon(q)$ on $B_\rho$. For
$q,q'\in B_\rho$, by the mean value theorem on the convex set $B_\rho$ and
\eqref{app:many:eq-lip},
\[
\norm{\Gamma(q)-\Gamma(q')}_\infty
=\Bigl\lVert L_0^{-1}\int_0^1\bigl(L_0-D\mathcal R^\varepsilon(q'+s(q-q'))\bigr)ds\,(q-q')\Bigr\rVert_\infty
\le\tfrac12\norm{q-q'}_\infty,
\]
and with \eqref{app:many:eq-res},
$\norm{\Gamma(q)-q^\tau}_\infty\le\tfrac12\rho+K\cdot11\varepsilon^2=\rho$. So
$\Gamma$ is a contraction of $B_\rho$ into itself, and its fixed point
$q^*\in B_\rho$ is a zero of $\mathcal R^\varepsilon$, that is a fixed point
of $\mathcal U$ in $(0,\infty)^N$, with
$\norm{q^*-q^\tau}_\infty\le\rho\le375\,\varepsilon^2$.

\emph{Step 4: positivity and distinctness.} The vertex precision of $0_1$
is $q^*_{0_1\to1_1}-1/q^*_{1_1\to0_1}\ge(t_1-\rho)-(t_1-\rho)^{-1}>0$, since
$t_1-\rho>1$. By Lemma~\ref{app:many:lem-cavity}, $q^*$ gives a positive
solution. If $\tau\ne\tau'$ differ at the copy $j$, then
$\abs{q^\tau_{0_j\to1_j}-q^{\tau'}_{0_j\to1_j}}=t_+-t_-=\tfrac7{15}>2\rho$, so the
balls are disjoint and the solutions are distinct.

\emph{Step 5: the Hessian.} Let $H(0)$ be the matrix \eqref{real:eq:Hu} of
$\Delta_m(\varepsilon)$ for the weights $u$ of $q^\tau$, with $u=0$ on the
bridges, and $H(\varepsilon)$ the one for the weights $u^*$ of $q^*$ given
by Lemma~\ref{app:many:lem-cavity}. Then $H(0)$ is the direct sum of the
matrices of Lemma~\ref{app:many:lem-diamond}(b), so it has
$\#\{j:\tau_j=-\}$ negative eigenvalues and all its eigenvalues have
absolute value greater than $\tfrac1{10}$. Both $\varphi_1$ and $\varphi_2$
are increasing with $\varphi_1'+\varphi_2'=(1-u)^{-2}$, so the row of a
vertex $i$ of $H(\varepsilon)-H(0)$ has absolute sum at most
$\sum_{e\ni i}\abs{u^*_e-u_e}\sup(1-u)^{-2}$, the supremum over the segment
between $u_e$ and $u^*_e$. For an internal edge with arcs of $q^\tau$-values
$\alpha,\beta$ we have $u=1/(\alpha\beta)$, $\abs{\partial u/\partial\alpha}=u/\alpha$, and
on $B_\rho$
\[
u\le\bar u=\frac1{(\alpha-\rho)(\beta-\rho)},\qquad
\abs{u^*_e-u_e}\le\rho\,\bar u\Bigl(\frac1{\alpha-\rho}+\frac1{\beta-\rho}\Bigr).
\]
For the edge $01$ ($\alpha=\beta=t$) and $t=t_-$ this gives $\bar u<0.6949$,
$\abs{u^*_e-u_e}\le1.1585\rho$ and $(1-\bar u)^{-2}<10.75$, so a contribution
less than $12.46\rho$; for a hub--terminal edge ($\alpha=t+1$,
$\beta=t/(t+1)$) and $t=t_-$ it gives $\bar u<0.83405$,
$\abs{u^*_e-u_e}\le1.9094\rho$ and $(1-\bar u)^{-2}<36.32$, so less than
$69.35\rho$. For $t=t_+$ the two contributions are smaller, less than
$1.06\rho$ and $7.44\rho$. For a bridge,
$u^*=\varepsilon^2/(q^*_{s\to s'}q^*_{s'\to s})\le122.1\,\varepsilon^2\le1.23\cdot10^{-4}$,
and its contribution to the rows of its ends is
$\varphi_1(u^*)+\varphi_2(u^*)=u^*/(1-u^*)\le1.24\cdot10^{-4}$. A hub row
contains one edge $01$ and two hub--terminal edges, a terminal row two
hub--terminal edges and at most one bridge, so
\[
\norm{H(\varepsilon)-H(0)}_2\le\norm{H(\varepsilon)-H(0)}_\infty
\le\max\bigl(151.2\,\rho,\ 138.7\,\rho+1.24\cdot10^{-4}\bigr)\le0.057<\tfrac1{10},
\]
where the first inequality holds because the matrix is symmetric. By
Weyl's inequality the eigenvalues of $H(\varepsilon)$, in increasing
order, differ from those of $H(0)$ by less than $\tfrac1{10}$, so they are
nonzero and have the same signs. By Lemma~\ref{real:lem:phi}(b) the
solution is nondegenerate with Morse index $\#\{j:\tau_j=-\}$.

Finally, the choices $\tau$ with exactly $k$ minus signs give $\binom mk$
distinct solutions of index $k$. The solution for $\tau=(-,\dots,-)$ has
index $m$, so $\taups\ge m$ by Proposition~\ref{prop:morse-bound}(i).
Removing the hubs $0_j$ leaves the path
$2_1,1_1,3_1,2_2,1_2,3_2,\dots,3_m$, so $\taups\le m$.
\end{proof}

\subsection{Bridges of two bouquets}\label{app:many:subsec-bouquets}

A \emph{bouquet} with center $z$ is a graph formed by $n\ge1$ cycles
through $z$, of lengths at least $3$, any two of which share only $z$. A
\emph{bridge of two bouquets} is the graph $G$ obtained from two disjoint
bouquets $G_a$ and $G_b$ with centers $a$ and $b$ by adding the edge $ab$.
It is connected and simple, and every vertex other than $a,b$ has degree
$2$. We write $\mathcal C_a$, $\mathcal C_b$ for the sets of cycles, $n_a$,
$n_b$ for their numbers and $\kappa=w_{ab}^2$.

\emph{One cycle.} Let $C=(z,j_1,\dots,j_L,z)$, $L\ge2$, be a cycle of $G$ through
a center $z$, with $c_0=w_{zj_1}^2$, $c_k=w_{j_kj_{k+1}}^2$ for
$1\le k<L$, and $c_L=w_{j_Lz}^2$. For $P\in\R$, the \emph{cycle equations}
$(\mathrm E_C(P))$ are the following equations for the $2(L+1)$ cavity
precisions on the arcs of $C$:
\begin{gather*}
q_{z\to j_1}=P+\frac{c_0}{q_{j_1\to z}},\qquad
q_{j_1\to j_2}=D_{j_1}-\frac{c_0}{q_{z\to j_1}},\qquad
q_{j_{k}\to j_{k+1}}=D_{j_k}-\frac{c_{k-1}}{q_{j_{k-1}\to j_k}}\ (2\le k\le L),
\end{gather*}
where $j_{L+1}=z$, together with the same equations for the reversed cycle
$\bar C=(z,j_L,\dots,j_1,z)$. A fixed point of $\mathcal U$ on $G$
satisfies $(\mathrm E_C(R_z(q)))$, because
$q_{z\to j_1}=D_z-\sum_{k\ne j_1}w_{zk}^2/q_{k\to z}=R_z(q)+c_0/q_{j_1\to z}$ and the
interior vertices of $C$ have degree~$2$.

For $s>0$ let $\kappa_1(s)=D_{j_1}-c_0/s$ and
$\kappa_{k+1}(s)=D_{j_{k+1}}-c_k/\kappa_k(s)$; these are the cavity
precisions along $C$ when $q_{z\to j_1}=s$. Let
$U_C=\{s>0:\kappa_1(s),\dots,\kappa_L(s)>0\}$, and define $U_{\bar C}$ for
$\bar C$ in the same way. With
\[
T_0=\begin{pmatrix}0&c_0\\1&0\end{pmatrix},\qquad
T_k=\begin{pmatrix}0&c_k\\-1&D_{j_k}\end{pmatrix}\ (1\le k\le L),\qquad
M_C=T_L\cdots T_1T_0=\begin{pmatrix}m_{11}&m_{12}\\m_{21}&m_{22}\end{pmatrix},
\]
the M\"obius transformation $M_C(s)=(m_{11}s+m_{12})/(m_{21}s+m_{22})$ is the
composition of $s\mapsto c_0/s$ and of the maps $\mu\mapsto c_k/(D_{j_k}-\mu)$;
so $M_C(s)=c_L/\kappa_L(s)$ for $s\in U_C$, the incoming message at $z$
from $j_L$. Note $\det M_C=-\prod_{k=0}^Lc_k<0$. If $m_{21}\ne0$ we put
\[
A_C=\frac{m_{11}}{m_{21}},\quad C_C=\frac{m_{22}}{m_{21}},\quad
B_C=-\frac{\det M_C}{m_{21}^2}>0,\quad b_C=2\sqrt{B_C},\quad E_C=A_C-C_C,
\]
so that $M_C(s)=A_C+B_C/(s+C_C)$.

\begin{lemma}\label{app:many:lem-cycle}
Let $C$ be a cycle as above.
\begin{enumerate}
\item[\textup{(a)}] $U_C$ is empty or an interval $(s_*,\infty)$ with
$s_*\ge0$; on $U_C$ all $\kappa_k$ are strictly increasing and $M_C$ is
finite, positive and strictly decreasing. The M\"obius transformation of
$\bar C$ is $M_{\bar C}=M_C^{-1}$.
\item[\textup{(b)}] Suppose $m_{21}\ne0$ and let $P>0$. The equations
$(\mathrm E_C(P))$ have at most one solution with all cavity precisions
positive. A solution exists if and only if
\[
s_C(P)=\tfrac12\bigl(P+\sqrt{P^2+b_C^2}\bigr)-C_C\in U_C
\quad\text{and}\quad
s'_C(P)=\tfrac12\bigl(P+\sqrt{P^2+b_C^2}\bigr)+A_C\in U_{\bar C} .
\]
Then $q_{z\to j_1}=s_C(P)$, $q_{z\to j_L}=s'_C(P)$, and the two incoming
messages at $z$ have the sum
\begin{equation}\label{app:many:eq-messages}
\frac{c_0}{q_{j_1\to z}}+\frac{c_L}{q_{j_L\to z}}=\sqrt{P^2+b_C^2}-P+E_C .
\end{equation}
\item[\textup{(c)}] If $m_{21}\ne0$, $s_C(0)=\tfrac12b_C-C_C\in U_C$ and
$s'_C(0)=\tfrac12b_C+A_C\in U_{\bar C}$, then the conditions of \textup{(b)}
hold for every $P>0$.
\end{enumerate}
\end{lemma}

We say that $C$ is \emph{valid} if the hypotheses of (c) hold.

\begin{proof}
(a) If $s\in U_C$ and $s'>s$, then $\kappa_1(s')>\kappa_1(s)>0$, and if
$\kappa_k(s')>\kappa_k(s)>0$ then
$\kappa_{k+1}(s')=D_{j_{k+1}}-c_k/\kappa_k(s')>\kappa_{k+1}(s)$. So $U_C$ is
closed upwards and the $\kappa_k$ increase on it; it is open, hence empty or
of the stated form, and $M_C=c_L/\kappa_L$ is positive and decreasing on
it. For the inverse, let $\mu_0=c_0/s$ and $\mu_k=c_k/(D_{j_k}-\mu_{k-1})$, so that
$y=M_C(s)=\mu_L$. The reversed recursion starts with $c_L/y$, and the
identity $D_{j_k}-c_k/\mu_k=\mu_{k-1}$ shows by induction that its messages
are $c_{k-1}/\mu_{k-1}$, $k=L,\dots,1$, and that its output is
$c_0/\mu_0=s$. This holds for all but finitely many $s$, so
$M_{\bar C}\circ M_C=\mathrm{id}$ as M\"obius transformations.

(b) Let a positive solution be given and put $x=c_0/q_{j_1\to z}$,
$y=c_L/q_{j_L\to z}$. The forward equations say that $s=q_{z\to j_1}=P+x$ lies
in $U_C$ and $y=M_C(s)$; the backward ones that $s'=q_{z\to j_L}=P+y$ lies in
$U_{\bar C}$ and $x=M_{\bar C}(s')$. By (a), $M_C(x)=s'=P+M_C(P+x)$, and $x$ is
not the pole $-C_C$ of $M_C$. With $\sigma=x+C_C$ this reads
$B_C/\sigma-B_C/(\sigma+P)=P$, that is $\sigma(\sigma+P)=B_C$. The pole of $M_C$
is not in the interval $U_C\ni s$ by (a), so $-C_C<s$ and
$\sigma+P=s+C_C>0$. Of the two roots of $\sigma^2+P\sigma-B_C$ only
$\sigma=\tfrac12(\sqrt{P^2+b_C^2}-P)$ has $\sigma+P>0$. Hence $s=s_C(P)$,
$y=A_C+B_C/(\sigma+P)=A_C+\sigma$ and $s'=s'_C(P)$, and all cavity
precisions are determined by $s$ and $s'$. This proves uniqueness and
\eqref{app:many:eq-messages}, since $x+y=2\sigma+A_C-C_C$. Conversely, if
$s=s_C(P)\in U_C$ and $s'=s'_C(P)\in U_{\bar C}$, define the cavity
precisions along $C$ and $\bar C$ from $s$ and $s'$. With $\sigma$ as above,
$\sigma>0$ and $\sigma(\sigma+P)=B_C$, so $M_C(s)=A_C+B_C/(\sigma+P)=s'-P$ and
$M_C(\sigma-C_C)=A_C+B_C/\sigma=s'$, that is $M_{\bar C}(s')=s-P$. These are the
two equations at $z$, so $(\mathrm E_C(P))$ holds with positive cavity
precisions.

(c) $s_C(P)$ and $s'_C(P)$ increase with $P$, and $U_C$, $U_{\bar C}$ are
closed upwards.
\end{proof}

\begin{lemma}\label{app:many:lem-triangle}
For a triangle $(z,u,v)$ with $c_0=w_{zu}^2$, $c_1=w_{uv}^2$, $c_2=w_{vz}^2$,
\[
M_C=\begin{pmatrix}c_2D_u&-c_0c_2\\ D_uD_v-c_1&-c_0D_v\end{pmatrix}.
\]
If $D_u=D_v=1$ and $c_1<1$, then $m_{21}=1-c_1$,
$b_C=2\sqrt{c_0c_1c_2}/(1-c_1)$, $E_C=(c_0+c_2)/(1-c_1)$, and the triangle is
valid.
\end{lemma}

\begin{proof}
The matrix is $T_2T_1T_0$. For $D_u=D_v=1$, $A_C=c_2/(1-c_1)$,
$C_C=-c_0/(1-c_1)$ and $B_C=c_0c_1c_2/(1-c_1)^2$. With
$\varrho=\sqrt{c_0c_1c_2}$ we get $s_C(0)=(\varrho+c_0)/(1-c_1)>0$,
$\kappa_1=1-c_0/s_C(0)=(\varrho+c_0c_1)/(\varrho+c_0)>0$ and
$\kappa_2=1-c_1/\kappa_1=\varrho(1-c_1)/(\varrho+c_0c_1)>0$. The reversed
triangle is obtained by exchanging $c_0$ and $c_2$, which exchanges $A_C$ and
$-C_C$, so $s'_C(0)\in U_{\bar C}$ in the same way.
\end{proof}

For a bridge of two bouquets all of whose cycles have $m_{21}\ne0$ we define
the \emph{side functions}
\begin{equation}\label{app:many:eq-Phi}
\Phi_A(s)=D_as+n_a-\sum_{C\in\mathcal C_a}\bigl(E_Cs+\sqrt{1+b_C^2s^2}\bigr),\qquad
\Phi_B(t)=D_bt+n_b-\sum_{C\in\mathcal C_b}\bigl(E_Ct+\sqrt{1+b_C^2t^2}\bigr),
\end{equation}
and, on $J_A=\{s>0:\Phi_A(s)>1\}$,
\begin{equation}\label{app:many:eq-h}
T(s)=\frac{\Phi_A(s)\bigl(\Phi_A(s)-1\bigr)}{\kappa s},\qquad
h(s)=\Phi_B\bigl(T(s)\bigr)-\Phi_A(s).
\end{equation}
Since $b_C>0$, $\Phi_A$ and $\Phi_B$ are strictly concave on $[0,\infty)$ and
vanish at $0$, so $J_A$ is an open interval or empty, and $h$ is real
analytic on $J_A$. The function $T$ and the number $\kappa=w_{ab}^2$ are not
related to the matrices $T_k$ and the functions $\kappa_k$ above.

\begin{proposition}\label{app:many:prop-reduction}
Let $G$ be a bridge of two bouquets with real coefficients such that every
cycle is valid.
\begin{enumerate}
\item[\textup{(a)}] The map $(V,f)\mapsto V_a$ is a bijection from the set of
positive solutions of \eqref{eq:system} with $z=0$ onto the set of zeros of
$h$ in $J_A$. At a positive solution, $V_b=T(V_a)$ and
$f_{ab}=\Phi_A(V_a)-1=\Phi_B(V_b)-1$.
\item[\textup{(b)}] Let $(V,f)$ be a positive solution, $s=V_a$,
$\nu=\Phi_A(s)$ and $\alpha=\Phi_A'(s)$. It is nondegenerate if and only if
$h'(s)\ne0$. If $h'(s)<0$, its Morse index is $1$. If $h'(s)>0$, its Morse
index is $0$ if $s\alpha>\nu(\nu-1)/(2\nu-1)$ and $2$ if
$s\alpha<\nu(\nu-1)/(2\nu-1)$, and one of the two holds.
\end{enumerate}
\end{proposition}

\begin{proof}
(a) Let $(V,f)$ be a positive solution and $q$ its fixed point
(Lemma~\ref{app:many:lem-cavity}), and put $P=R_a(q)=1/V_a$. Every cycle
$C\in\mathcal C_a$ satisfies $(\mathrm E_C(P))$ with positive cavity
precisions, so by \eqref{app:many:eq-messages}
\[
P=D_a-\sum_{C\in\mathcal C_a}\Bigl(\sqrt{P^2+b_C^2}-P+E_C\Bigr)-\frac\kappa{q_{b\to a}} .
\]
Multiplying by $V_a$ and using $\kappa V_a/q_{b\to a}=f_{ab}$
(Lemma~\ref{app:many:lem-cavity}) gives $1=\Phi_A(V_a)-f_{ab}$; in the same
way $1=\Phi_B(V_b)-f_{ab}$. Since $f_{ab}>0$, $V_a\in J_A$, and the edge
equation $f_{ab}(1+f_{ab})=\kappa V_aV_b$ gives $V_b=T(V_a)$, hence
$h(V_a)=0$. The map is injective: $V_a$ determines $V_b$ and
$f_{ab}$, hence $q_{b\to a}=\kappa V_a/f_{ab}$ and $q_{a\to b}=\kappa V_b/f_{ab}$,
and by Lemma~\ref{app:many:lem-cycle}(b) the cavity precisions on every
cycle are determined by $R_a=1/V_a$ or $R_b=1/V_b$; and $q$ determines the
solution.

Conversely let $s\in J_A$ with $h(s)=0$, and put $\nu=\Phi_A(s)>1$,
$t=T(s)>0$, so that $\Phi_B(t)=\nu$ and $\nu(\nu-1)=\kappa st$. On every
cycle of $\mathcal C_a$ take the positive solution of
$(\mathrm E_C(1/s))$, on every cycle of $\mathcal C_b$ that of
$(\mathrm E_C(1/t))$ (Lemma~\ref{app:many:lem-cycle}(b),(c)), and put
$q_{a\to b}=\nu/s$, $q_{b\to a}=\kappa s/(\nu-1)$. By
\eqref{app:many:eq-messages}, $\Phi_A(s)=\nu$ says that the sum of the
cycle messages at $a$ is $D_a-\nu/s$, so
$R_a(q)=\nu/s-\kappa/q_{b\to a}=1/s$. Similarly the cycle messages at $b$ sum to
$D_b-\nu/t$ and $R_b(q)=\nu/t-\kappa s/\nu=1/t$, using $\kappa s/\nu=(\nu-1)/t$.
Hence the equations of $\mathcal U$ hold on the arcs out of $a$ towards the
cycles (they are part of $(\mathrm E_C(R_a))$), on
$a\to b$ ($R_a+\kappa/q_{b\to a}=\nu/s$) and on $b\to a$
($R_b+\kappa/q_{a\to b}=1/t+\kappa s/\nu=\kappa s/(\nu-1)$), and likewise at
$b$; the other arcs are covered by the cycle equations. So $q>0$ is a fixed
point with $R_a>0$, and Lemma~\ref{app:many:lem-cavity} gives a positive
solution with $V_a=s$.

(b) Let $S=\{a,b\}$ and $W=V\setminus S$. A direct computation from
\eqref{real:eq:phi}, using $g'(x)=f_+(x)/x$ and $f_+'(x)=1/(1+2f_+(x))$, shows
that $2\diag(V)\Hess\phi(V)\diag(V)$ equals the matrix $H_u$ of
\eqref{real:eq:Hu} with $u_e=f_e/(1+f_e)$, $f_e=f_+(w_e^2V_iV_j)$, at
\emph{every} $V\in(0,\infty)^V$; the vertex equations are not needed for the
second derivatives. Since $G-S$ is a forest, the proof of
Proposition~\ref{prop:morse-bound}(i) shows that $H_u[W]$, hence
$\Hess_{WW}\phi$, is positive definite everywhere. By the implicit function
theorem the equations $\nabla_W\phi=0$ define near the given solution a
real analytic map $V_S\mapsto V_W^*(V_S)$, and
$\psi(V_S)=\phi(V_W^*(V_S),V_S)$ satisfies
$\nabla\psi=\nabla_S\phi(V_W^*(V_S),V_S)$ and
$\Hess\psi=\Hess_{SS}\phi-\Hess_{SW}\phi\,(\Hess_{WW}\phi)^{-1}\Hess_{WS}\phi$ at the
solution. By the additivity of inertia~\cite{Haynsworth}, $\Hess\phi$ has
$\abs W$ positive eigenvalues more than $\Hess\psi$, and the same numbers of
negative and zero eigenvalues.

At $(V_W^*(V_S),V_S)$ the vertex equations at $W$ and the edge relations
hold with $V>0$. Put $q_{i\to j}=(1+f_{ij})/V_i>0$. The computation in the
proof of Lemma~\ref{lem:bp-correspondence}, which at a vertex uses only
the vertex equation at that vertex and the edge equations, gives the
cavity equations on the arcs out of the vertices of $W$; at $a$ the
equation of $(\mathrm E_C(1/V_a))$ is
$(1+f_{aj_1})/V_a=1/V_a+c_0V_{j_1}/(1+f_{aj_1})$, which is the edge
equation. So every cycle at $a$ satisfies $(\mathrm E_C(1/V_a))$ with
positive cavity precisions. By \eqref{app:many:eq-messages} and
$f_{aj}=V_a\,c/q_{j\to a}$, the sum of $f_e$ over the cycle edges at $a$ is
$\sum_C(E_CV_a+\sqrt{1+b_C^2V_a^2}-1)$. With Lemma~\ref{real:lem:phi}(a),
\[
2V_a\,\partial_{V_a}\psi=\Phi_A(V_a)-1-f_+(\kappa V_aV_b),\qquad
2V_b\,\partial_{V_b}\psi=\Phi_B(V_b)-1-f_+(\kappa V_aV_b).
\]
Differentiating at the solution, where both right-hand sides vanish, and
writing $\beta=\Phi_B'(V_b)$ and $p=\kappa f_+'(\kappa V_aV_b)=\kappa/(2\nu-1)>0$,
\[
2\Hess\psi=\begin{pmatrix}(\alpha-pV_b)/V_a&-p\\-p&(\beta-pV_a)/V_b\end{pmatrix},
\qquad
4V_aV_b\det\Hess\psi=\alpha\beta-p(\alpha V_a+\beta V_b).
\]
On the other hand, $T(s)=\nu(\nu-1)/(\kappa s)$ gives
$T'(s)=((2\nu-1)\alpha-\kappa T)/(\kappa s)=(\alpha-pV_b)/(pV_a)$, hence
$pV_ah'(s)=\beta(\alpha-pV_b)-p\alpha V_a=4V_aV_b\det\Hess\psi$. So
$\det\Hess\psi$ has the sign of $h'(s)$. If $h'(s)<0$, $\Hess\psi$ has one
negative and one positive eigenvalue; if $h'(s)>0$ it is definite, and its
sign is that of $\alpha-pV_b$, which is nonzero; and
$V_a(\alpha-pV_b)=s\alpha-\kappa V_aV_b/(2\nu-1)=s\alpha-\nu(\nu-1)/(2\nu-1)$.
\end{proof}

\subsection{The number of positive solutions is unbounded}
\label{app:many:subsec-unbounded}

Let $\pi(x)=1+x-\sqrt{1+x^2}=2x/(1+x+\sqrt{1+x^2})$ for $x\ge0$.

\begin{lemma}\label{app:many:lem-pi}
$\pi$ is increasing, $\pi(0)=0$, $\pi(1)=2-\sqrt2$, and for $x>0$
\[
\max\Bigl(x-\frac{x^2}2,\ 1-\frac1{2x}\Bigr)\le\pi(x)<\min(x,1).
\]
\end{lemma}

\begin{proof}
$\pi'(x)=1-x/\sqrt{1+x^2}>0$. The bounds follow from
$\sqrt{1+x^2}\le1+x^2/2$, $\sqrt{1+x^2}\le x+1/(2x)$ and
$\sqrt{1+x^2}>\max(1,x)$.
\end{proof}

\emph{The family.} Let $\theta=100$. A \emph{triangle of type $\beta>0$ at a
center $z$} consists of two new vertices $u,v$ with $D_u=D_v=1$ and the
edges $zu$, $uv$, $vz$ with $w_{zu}=w_{vz}=\tfrac12\sqrt{3\beta}$ and
$w_{uv}=\tfrac12$. By Lemma~\ref{app:many:lem-triangle} it is valid, with
$b_C=\beta$ and $E_C=2\beta$. For integers $K\ge2$ and
$\kappa\in[\tfrac14,\tfrac{101}{400}]$ let $\mathfrak M_K(\kappa)$ be the
bridge of two bouquets with
\begin{itemize}
\item at $a$: for each $k=0,1,\dots,K$, $\theta^k$ triangles of type
$\theta^{-2k}$, and $D_a=3\sum_{k=0}^K\theta^{-k}-\theta^{-3K}$;
\item at $b$: $n_B=40\,\theta^{K-1}$ triangles of type $1$, and
$D_b=\tfrac52n_B$;
\item $w_{ab}=\sqrt\kappa$.
\end{itemize}
Its graph $G_K$ has $2+2\bigl(\sum_{k=0}^K\theta^k+n_B\bigr)$ vertices.
Since $\sum_CE_C=2\sum_C\beta_C$ and
$\beta s+1-\sqrt{1+\beta^2s^2}=\pi(\beta s)$, the side functions
\eqref{app:many:eq-Phi} are
\begin{equation}\label{app:many:eq-family}
\Phi_A(s)=\sum_{k=0}^K\theta^k\pi(\theta^{-2k}s)-\theta^{-3K}s,\qquad
\Phi_B(t)=n_B\,\gamma(t),\quad \gamma(t)=1+\frac t2-\sqrt{1+t^2}.
\end{equation}
In this subsection $K$ is the number of scales, not the constant of
Section~\ref{app:many:subsec-chain}; in part (a) of the following theorem
and its proof, $c_k$ and $j_k$ are points of $J_A$, not the squared weights
of Lemma~\ref{app:many:lem-cycle}, while in part (b), as in
Corollary~\ref{cor:parity}, $c_k$ is the number of positive solutions of
Morse index~$k$.

\begin{theorem}\label{app:many:thm-unbounded}
\begin{enumerate}
\item[\textup{(a)}] Let $K\ge2$ and $\kappa\in[\tfrac14,\tfrac{101}{400}]$. The
graph $G_K$ is connected and simple with $\taups(G_K)=2$. The points
$c_k=\theta^{2k}$ $(1\le k\le K)$ and $j_k=\theta^{2k+1}$ $(1\le k\le K-1)$
lie in $J_A$, and $h(c_k)<0<h(j_k)$. Hence $h$ has at least $2K-2$ zeros in
$J_A$, and $\mathfrak M_K(\kappa)$ has at least $2K-2$ positive solutions.
\item[\textup{(b)}] Let $m\ge1$ and $K=m+2$. There is a nonempty open set
$\mathcal O$ of real coefficients on $G_K$ such that every model in
$\mathcal O$ has at least $2m+2$ positive solutions, and a closed subset
$\mathcal N\subset\mathcal O$ with empty interior such that every model in
$\mathcal O\setminus\mathcal N$ satisfies condition \textup{(G)}, has only
nondegenerate positive solutions, and has
$(c_0,c_1,c_2)=(1,c_1,c_1-1)$ with $c_1\ge m+1$ and $c_k=0$ for $k\ge3$.
\end{enumerate}
In particular Proposition~\ref{real:prop:unbounded} holds.
\end{theorem}

\begin{proof}
(a) \emph{The graph.} $G_K-\{a,b\}$ is a disjoint union of edges, so
$\taups\le2$. If one vertex $x$ is removed, then either $x\ne a$ and two
triangles at $a$ not containing $x$ remain, or $x=a$ and the triangles at
$b$ remain; in both cases a component of cycle rank at least $2$ remains.
So $\taups(G_K)=2$.

\emph{The points $c_k$.} Keeping only the term $k$ in
\eqref{app:many:eq-family} and using $\pi\ge0$,
$\Phi_A(c_k)\ge\theta^k\pi(1)-\theta^{2k-3K}$. Since $2k-3K\le k-4$, we get
$\Phi_A(c_k)\ge\theta^kL_c$ with $L_c=2-\sqrt2-10^{-4}>0.5856$, in particular
$c_k\in J_A$. The function $x\mapsto x(x-1)$ increases for $x>\tfrac12$, so
\[
\frac{\Phi_A(c_k)(\Phi_A(c_k)-1)}{c_k}\ge L_c(L_c-\theta^{-k})
\ge L_c(L_c-0.01)>0.3371,
\qquad T(c_k)>\frac{0.3371}{\kappa}\ge\frac{0.3371}{0.2525}>\frac43 .
\]
Now $\gamma(t)<0$ if and only if $(1+t/2)^2<1+t^2$, that is $t>\tfrac43$.
So $\Phi_B(T(c_k))<0<\Phi_A(c_k)$ and $h(c_k)<0$.

\emph{The points $j_k$.} Let $1\le k\le K-1$. Using $\pi<1$ for the terms
$i\le k$ and $\pi(x)<x$ for $i>k$,
\[
\Phi_A(j_k)<\sum_{i\le k}\theta^i+\sum_{i>k}\theta^{2k+1-i}
<\frac{\theta^{k+1}}{\theta-1}+\frac{\theta^{k+1}}{\theta-1}=U_j\theta^k,
\qquad U_j=\frac{2\theta}{\theta-1}=\frac{200}{99}.
\]
Using $\pi(x)\ge1-1/(2x)$ for $i\le k$ (with $x=\theta^{2(k-i)+1}$),
$\pi(x)\ge x-x^2/2$ for $i=k+1$ (with $x=\theta^{-1}$), $\pi\ge0$ for $i>k+1$,
$\sum_{i\le k}\theta^{3i}<\theta^{3k+3}/(\theta^3-1)$ and
$\theta^{-3K}j_k=\theta^{2k+1-3K}\le\theta^{k-2K}\le\theta^k\cdot10^{-8}$ (as
$k\le K-1$ and $K\ge2$),
\[
\frac{\Phi_A(j_k)}{\theta^k}\ge\frac\theta{\theta-1}+1-\frac1{2\theta}
-\frac{\theta^2}{2(\theta^3-1)}-\frac1{\theta(\theta-1)}-10^{-8}
>1.9999=:L_j ;
\]
the term $1/(\theta(\theta-1))$ bounds $\theta^{-k}/(\theta-1)$, which comes
from $\sum_{i\le k}\theta^i=(\theta^{k+1}-1)/(\theta-1)$, and the middle
expression equals $2-1.5\cdot10^{-8}$ to within $10^{-14}$. Hence
$j_k\in J_A$, and
\[
\frac{\Phi_A(j_k)(\Phi_A(j_k)-1)}{j_k}\in
\Bigl[\frac{L_j(L_j-0.01)}\theta,\ \frac{U_j^2}\theta\Bigr]\subset[0.039796,\ 0.040813],
\]
so $T(j_k)\in[0.039796/0.2525,\ 0.040813/0.25]\subset[0.1576,\ 0.16326]$. The
function $\gamma$ is concave, and $\gamma(0.1576)>0.06645$,
$\gamma(0.16326)>0.06839$; so $\gamma(T(j_k))>0.0664$ and
\[
\Phi_B(T(j_k))>0.0664\cdot40\,\theta^{K-1}=2.656\,\theta^{K-1}\ge2.656\,\theta^k
>U_j\theta^k>\Phi_A(j_k),
\]
that is $h(j_k)>0$.

\emph{Counting.} The points satisfy $c_1<j_1<c_2<\dots<j_{K-1}<c_K$. Since
$\Phi_A$ is concave and exceeds $1$ at $c_1$ and $c_K$, the interval
$[c_1,c_K]$ lies in $J_A$, where $h$ is continuous. The $2K-1$ points
carry alternating signs, so each of the $2K-2$ open intervals between
consecutive points contains a zero of $h$. By
Proposition~\ref{app:many:prop-reduction}(a) these give $2K-2$ distinct
positive solutions.

(b) Let $p_0$ be the coefficient point of $\mathfrak M_K(\tfrac14)$, which has
all weights positive. We use the following finitely many conditions on the
coefficients $p=(D,w)$ of $G_K$: for every triangle,
$m_{21}=D_uD_v-w_{uv}^2\ne0$ and the four cavity precisions (two forward,
two backward) at the inputs $s_C(0)$ and $s'_C(0)$ of
Lemma~\ref{app:many:lem-cycle}(c) are positive; further
$h(c_k)<0$ for $1\le k\le K$, $h(j_k)>0$ for $1\le k\le K-1$, and
$\Phi_A(c_1)>1$, $\Phi_A(c_K)>1$. By (a) and
Lemma~\ref{app:many:lem-triangle} they hold at $p_0$. Where $m_{21}\ne0$,
the numbers $A_C,B_C,C_C,b_C,E_C$ are continuous functions of $p$, hence so
are $s_C(0)$, $s'_C(0)$, the cavity precisions computed from them (whose
denominators are positive at $p_0$), and, for fixed $s$, the values
$\Phi_A(s)$, $T(s)$ and $\Phi_B(T(s))$. So the conditions hold on an open
neighborhood $\mathcal O$ of $p_0$, on which we also keep $w>0$. For
$p\in\mathcal O$ every triangle is valid, Proposition~\ref{app:many:prop-reduction}
applies, $\Phi_A$ is concave, and the argument of (a) gives at least
$2K-2=2m+2$ positive solutions.

By Theorem~\ref{thm:generic-degree} and Lemma~\ref{real:lem:phi}(c) there is
a nonempty Zariski open set of complex coefficients on which all positive
solutions are nondegenerate (compare the proof of
Corollary~\ref{real:cor:generic}). Its complement is contained in the zero
set of a nonzero complex polynomial, whose real or imaginary part is a
nonzero real polynomial, and the zero set of a nonzero real polynomial is
closed with empty interior. Since $w>0$ on $\mathcal O$, condition (G)
fails only where $\det\bigl((D-A)[S]\bigr)=0$ for some nonempty connected
$S$; each of these determinants is a polynomial in $p$ that does not vanish
for large $D$. Let $\mathcal N$ be the intersection of $\mathcal O$ with the
union of the zero set of the real polynomial above and of the finitely many
zero sets of these determinants.

Let $p\in\mathcal O\setminus\mathcal N$. By
Corollary~\ref{cor:parity} the number $N$ of positive solutions is finite,
and since $N\ge2m+2>1$ the model is not walk-summable
(Theorem~\ref{thm:unique-ws}). By Proposition~\ref{prop:morse-bound}(i) the
Morse indices are at most $\taups(G_K)=2$; by the last statement of
Corollary~\ref{cor:parity}, $c_0=1$; and by its identity
$c_0-c_1+c_2=0$. Hence $c_2=c_1-1$ and $2m+2\le N=c_0+c_1+c_2=2c_1$.
\end{proof}

\begin{remark}\label{app:many:rem-mechanism}
Heuristically, the $\theta^k$ triangles of type $\theta^{-2k}$ switch on
near $s=\theta^{2k}$ and raise $\Phi_A$ by about $\theta^k$, so that the
quotient $\Phi_A(\Phi_A-1)/s$ oscillates, as $\log s$ increases, between
values of at least about $0.34$ at the points $c_k$ and about $4/\theta$ at
the points $j_k$.
The function $\Phi_B=n_B\gamma$ is negative for $t>\tfrac43$ and of order
$n_B$ near $t=0.16$, and the many triangles at $b$ make $\Phi_B(T(j_k))$
larger than $\Phi_A(j_k)$. The graph $G_{m+2}$ has about
$2.8\cdot10^{2m+4}$ vertices, so the lower bound $c_2\ge m$ is of the
order of the logarithm of the number of vertices.
\end{remark}

\subsection{Certified small instances}\label{app:many:subsec-certified}

Both models have $D_i=1$ for all $i$ and rational weights; they are bridges
of two bouquets with $a=0$ and $b=1$.

\emph{The model $\mathfrak S_6$} has $27$ vertices and $34$ edges:
$w_{ab}=1868/785871$, and the cycles of Table~\ref{app:many:tab-S6}, each
given by its weights in the order of a walk around it starting at its
center.

\begin{table}[ht]
\centering\small
\begin{tabular}{@{}cl@{}}
\toprule
center & weights along the cycle\\
\midrule
$a$ & $17797/902086$, $424644/920567$, $5266/63711$\\
$a$ & $133669/889500$, $718576/892037$, $94025/777589$\\
$a$ & $79138/955265$, $135205/930357$, $373096/994207$, $65041/430754$\\
$a$ & $80096/821269$, $264419/353155$, $101007/997834$, $584612/607337$, $39225/219422$\\
$b$ & $307983/999116$, $68033/380581$, $323645/416872$, $17358/399337$\\
$b$ & $100331/819923$, $5283/50176$, $777902/983059$, $115211/985213$, $36749/974461$\\
$b$ & $147411/499304$, $466238/748569$, $433457/609660$, $131637/942167$\\
$b$ & $28211/553112$, $537/1633$, $236712/357209$, $260713/834717$, $30577/915607$\\
\bottomrule
\end{tabular}
\caption{The cycles of the model $\mathfrak S_6$.}
\label{app:many:tab-S6}
\end{table}

\emph{The model $\mathfrak S_{10}$} has $68$ vertices and $100$ edges. All its
cycles are triangles $(z,u,v)$ with $w_{zu}=w_{vz}=p$ and $w_{uv}=\tfrac12$;
at $a$ there are $11$ triangles with $p=130923/870130$ and $6$ with
$p=10788/839159$; at $b$ there are $2$, $1$, $4$ and $9$ triangles with
$p=276918/802963$, $35431/637504$, $15763/362293$ and $11425/892976$; and
$w_{ab}=671/876039$.

\begin{proposition}\label{app:many:prop-certified}
The model $\mathfrak S_6$ has exactly $6$ positive solutions and
$\mathfrak S_{10}$ has exactly $10$. All of them are nondegenerate, and
$(c_0,c_1,c_2)=(1,3,2)$ for $\mathfrak S_6$ and $(1,5,4)$ for
$\mathfrak S_{10}$. In increasing order of $V_a$ the Morse indices are
$0,1,2,1,2,1$, respectively $0,1,2,1,2,1,2,1,2,1$.
\end{proposition}

The proof is computer-assisted. It consists of a complete count by
Proposition~\ref{app:many:prop-reduction}, which uses exact rational
arithmetic and the rational interval arithmetic (I3) of
Appendix~\ref{app:cert:subsec-methods}, and of an independent
certification of every solution in the full system with the arithmetic (I1)
of that appendix.

\begin{proof}
\emph{Step 1: the reduction (exact).} For each cycle the matrix $M_C$ was
computed in rational arithmetic; in all cases $m_{21}\ne0$, so $b_C^2=4B_C$
and $E_C$ are rational. The cavity precisions of
Lemma~\ref{app:many:lem-cycle}(c) at the inputs $s_C(0)=\sqrt{B_C}-C_C$ and
$s'_C(0)=\sqrt{B_C}+A_C$ were enclosed with the arithmetic (I3) and are
positive, so every cycle is valid (for the triangles of $\mathfrak S_{10}$
this is also Lemma~\ref{app:many:lem-triangle}). Hence
Proposition~\ref{app:many:prop-reduction} applies, with $\Phi_A$, $\Phi_B$ of
\eqref{app:many:eq-Phi} having explicit rational coefficients $D_a-\sum E_C$
and $b_C^2$. For the evaluation we write
$\Phi_A(s)=\sum_C\pi(b_Cs)-\epsilon_As$ with
$\epsilon_A=\sum_Cb_C-D_a+\sum_CE_C$, and similarly for
$\Phi_B$; the form $\pi(x)=2x/(1+x+\sqrt{1+x^2})$ avoids cancellation.

\emph{Step 2: interval enclosures.} For an interval $I=[s_1,s_2]\subset(0,\infty)$,
enclosures of $\Phi_A(I)$ and $\Phi_A'(I)$ were computed with (I3) from the
monotonicity of $\pi$ (increasing) and of
$\pi'(x)=1/(r(r+x))$, $r=\sqrt{1+x^2}$ (decreasing): $\sum_C\pi(b_Cs_1)$ and
$\sum_C\pi(b_Cs_2)$ bound the first sum, and the term $-\epsilon_As$ is an
interval product. Enclosures of $T(I)$, $h(I)$ and of
$h'=\Phi_B'(T)\,T'-\Phi_A'$ with
$T'=((2\Phi_A-1)\Phi_A's-\Phi_A(\Phi_A-1))/(\kappa s^2)$ follow by interval
arithmetic.

\emph{Step 3: the complete count.} Rational numbers
$s_{\mathrm{lo}}<s_{\mathrm{hi}}$ were chosen with $\Phi_A(s_{\mathrm{lo}})<1$,
$\Phi_A'(s_{\mathrm{lo}})>0$, $\Phi_A(s_{\mathrm{hi}})<1$ and
$\Phi_A'(s_{\mathrm{hi}})<0$, verified with (I3)
($s_{\mathrm{lo}}\approx4.53$, $s_{\mathrm{hi}}\approx3.33\cdot10^3$ for
$\mathfrak S_6$; $s_{\mathrm{lo}}\approx1.92$, $s_{\mathrm{hi}}\approx2.68\cdot10^5$
for $\mathfrak S_{10}$). Since $\Phi_A$ is concave, $\Phi_A'\ge\Phi_A'(s_{\mathrm{lo}})>0$
on $(0,s_{\mathrm{lo}}]$, so $\Phi_A<1$ there, and similarly on
$[s_{\mathrm{hi}},\infty)$; hence $J_A\subset(s_{\mathrm{lo}},s_{\mathrm{hi}})$.
The interval $[s_{\mathrm{lo}},s_{\mathrm{hi}}]$ was bisected until each piece
$I$ satisfied one of the following, checked with the enclosures of Step~2:
\begin{enumerate}
\item[(i)] $\sup\Phi_A(I)\le1$; then $I\cap J_A=\emptyset$;
\item[(ii)] $0\notin h(I)$; then $h$ has no zero in $I$;
\item[(iii)] $\inf\Phi_A(I)>1$, $0\notin h'(I)$, and $h$ has the same strict
sign at both ends of $I$; then $h$ is monotone without zero on $I$;
\item[(iv)] $\inf\Phi_A(I)>1$, $0\notin h'(I)$, and $h$ has strictly
opposite signs at the ends of $I$; then $I\subset J_A$ and $h$ has exactly
one zero in $I$, at which $h'\ne0$.
\end{enumerate}
The pieces cover $[s_{\mathrm{lo}},s_{\mathrm{hi}}]$ and every zero of $h$ in
$J_A$ lies in some piece, which must be of type (iv). So the number of
zeros of $h$ in $J_A$ is the number of pieces of type (iv). The bisection
ended with $193$ pieces for $\mathfrak S_6$, of which $6$ are of type (iv),
and $505$ pieces for $\mathfrak S_{10}$, of which $10$ are of type (iv). By
Proposition~\ref{app:many:prop-reduction}(a) the models have exactly $6$ and
exactly $10$ positive solutions, and by part (b) all of them are
nondegenerate, since $h'\ne0$ on the pieces of type (iv). On each piece of
type (iv) the sign of $h'$ and, when $h'>0$, the sign of
$s\Phi_A'(s)-\nu(\nu-1)/(2\nu-1)$ with $\nu=\Phi_A(s)$ were enclosed with (I3)
and are constant; by Proposition~\ref{app:many:prop-reduction}(b) this gives
the indices stated in the proposition. The completeness of the list, and not
only a lower bound, comes from the covering in this step: no part of
$J_A$ is left unexamined, and each piece contains at most one zero.

\emph{Step 4: independent certification in the full system.} With $D=1$ the
system \eqref{eq:system} at $z=0$ is equivalent, by $V_i=g_i$, to the system
\[
F_e(f)=f_e(1+f_e)-w_e^2g_ig_j=0\quad(e=ij\in E),\qquad g_i=1+\sum_{e\ni i}f_e,
\]
in the unknowns $f\in\R^E$, and its Jacobian matrix is the Schur complement
of the identity block in the Jacobian matrix of \eqref{eq:system} with
respect to $(V,f)$. So a zero $f>0$ gives the positive solution $(g,f)$,
and by Lemma~\ref{real:lem:phi}(c) this solution is nondegenerate if
$F'(f)$ is invertible. For each of the $6$, respectively $10$, solutions, a box $X$
around a floating point approximation, of relative radius between
$10^{-12}$ and $10^{-8}$, satisfied the Krawczyk test of
Lemma~\ref{app:cert:lem-krawczyk} in the arithmetic (I1), with the rational
numbers $w_e^2$ enclosed by floating point intervals and with an
approximate inverse $Y$ of $F'$ at the center. By that lemma, $X$ contains
exactly one zero of $F$, and $F'$ is invertible there. The lower end points
of all coordinates of all boxes are positive (the smallest is
$5.0\cdot10^{-4}$ for $\mathfrak S_6$ and $7.8\cdot10^{-6}$ for
$\mathfrak S_{10}$), so the zeros are nondegenerate positive solutions, and
the boxes are pairwise disjoint, so they are distinct. For each box the
entries of $H_u$ in \eqref{real:eq:Hu}, written as $1+\sum f_e^2/(1+2f_e)$ and
$-f_e(1+f_e)/(1+2f_e)$, are increasing functions of the $f_e$ and were
enclosed by their values at the end points; the test of
Lemma~\ref{app:cert:lem-inertia} succeeded, with margins $\min_i\abs{d_i}-e$
between $0.0086$ and $0.053$ for $\mathfrak S_6$ and between $0.0025$ and
$0.31$ for $\mathfrak S_{10}$, and gave the Morse indices. Since the
total numbers of positive solutions are $6$ and $10$ by Step~3, the boxes
contain all positive solutions, and this step confirms independently
the nondegeneracy and the indices. For every box the value of $V_a=g_a$
lies in the piece of type (iv) of Step~3 with the same index.
\end{proof}

The approximate values of $(V_a,V_b)$ at the solutions are, in the order
of the proposition, $(7.44,11.9)$, $(14.8,6.83\cdot10^3)$,
$(95.4,4.13\cdot10^3)$, $(330,1.75\cdot10^3)$, $(1.13\cdot10^3,311)$,
$(1.78\cdot10^3,36.7)$ for $\mathfrak S_6$, and $(3.15,4.22)$,
$(7.89,6.63\cdot10^5)$, $(1.32\cdot10^3,1.77\cdot10^5)$,
$(8.10\cdot10^3,4.19\cdot10^4)$, $(2.05\cdot10^4,1.60\cdot10^4)$,
$(4.94\cdot10^4,4.50\cdot10^3)$, $(9.67\cdot10^4,853)$,
$(1.29\cdot10^5,200)$, $(1.48\cdot10^5,47.2)$, $(1.58\cdot10^5,10.1)$ for
$\mathfrak S_{10}$. In both models the solution of index $0$ is the maximal
fixed point (Proposition~\ref{prop:morse-bound}(iv)), and the identity
$c_0-c_1+c_2=0$ of Corollary~\ref{cor:parity} holds, although we have not
verified condition (G) for these models. The model $\mathfrak S_6$ shows
that the bound $c_k\le\binom{\taups}k$, which holds for $\taups\le1$ and is
attained by Theorem~\ref{app:many:thm-chain}, fails for $\taups=2$.

%% file: appendices/weak-bridge.tex
\section{Weak bridges and compact components}\label{app:bridge}

We prove the weak-bridge gluing theorem, construct graphs with arbitrarily
many compact components of $\Xp$, and use a bridge construction to prove
the existence direction of the packing theorem.

Throughout, the diagonal is $D-z$, and $g(z,V)$ is the map of
\eqref{per:eq:g}. Since $f_e=f_+(w_e^2V_iV_j)$ at a positive solution, we
identify $\Xp$ with $\{(z,V)\in\R\times(0,\infty)^V:\ g(z,V)=0\}$; a point
with $V>0$ and $g=0$ automatically has $f_e>0$ and $z<\min_iD_i$. A point of
$\Xp$ is \emph{regular} if the $n\times(n+1)$ matrix $(\partial_zg\mid J_V)$
has rank $n$. By the implicit function theorem, near a regular point $\Xp$
is a real analytic curve. A \emph{circle} of $\Xp$ is a compact connected
component of $\Xp$; when all its points are regular it is a compact
connected one-dimensional manifold, hence diffeomorphic to a circle. Some
letters have a second, local meaning fixed by the context: $g$ also denotes
the function of \eqref{real:eq:phi} in the proof of
Lemma~\ref{app:bridge:lem-hess} and the number $-S'(\sigma)$ in
Lemma~\ref{app:bridge:lem-response}; $\varphi$ is a convex function in the
proof of Lemma~\ref{app:bridge:lem-roots} and $\Phi_1(\sigma_2,\cdot)$ in
\eqref{app:bridge:eq-Z}; and $\Phi$ denotes the response function, except
for the map $\Phi_z$ of \eqref{per:eq:Phi} in
Lemma~\ref{app:bridge:lem-left}. We use the following gluing principle.

\begin{lemma}\label{app:bridge:lem-glue}
Let $X$ be a topological space, let $\Gamma_1,\dots,\Gamma_m\subseteq X$ be
connected and closed in $X$, and let $O_1,\dots,O_m\subseteq X$ be open in
$X$ with $O_i\subseteq\Gamma_i$. Suppose that $\Gamma:=\bigcup_i\Gamma_i$ is
contained in $\bigcup_iO_i$ and that the graph on $\{1,\dots,m\}$ in which
$i$ and $j$ are adjacent when $\Gamma_i\cap\Gamma_j\ne\emptyset$ is
connected. Then $\Gamma$ is open and closed in $X$ and is a connected
component of $X$. If every point of $\Gamma$ has a neighborhood in $X$
homeomorphic to an open interval, then $\Gamma$ is homeomorphic to a circle
if it is compact and to $\R$ otherwise.
\end{lemma}

\begin{proof}
$\Gamma$ is closed as a finite union of closed sets, and
$\Gamma=\bigcup_iO_i$ is open. It is connected because the $\Gamma_i$ are
connected and the adjacency graph is connected. A connected set that is
open and closed is a component. The last statement is the classification
of connected one-dimensional manifolds~\cite{Milnor}.
\end{proof}

\subsection{Parts and their response functions}\label{app:bridge:subsec-parts}

\begin{definition}\label{app:bridge:def-part}
A \emph{part} is a pair $(G,a)$, where $G=(V,E)$ is a connected loopless
multigraph with $\beta(G)\ge2$ which is either simple or a \emph{dipole}
(two vertices joined by $k\ge3$ parallel edges), with real coefficients
$D$, $w$, and $a\in V$ is a vertex such that $G-a$ is a pseudoforest. We
write $V^\circ=V\setminus\{a\}$, $b_j=\sum_{e=aj}\abs{w_e}$ for
$j\in V^\circ$,
\[
M(z)=(D-z-\abs A)[V^\circ],\qquad \rho=\lambda_{\min}\bigl((D-\abs A)[V^\circ]\bigr),
\]
and $\sigma=\sigm(G)$, $\Sigma=\Sigm(G)$ for the lower Perron value and the
bottom of the spectrum of $\tH$.
\end{definition}

For a dipole, $V^\circ$ is a single vertex. A simple part has
$\taups(G)\le1$. The diamond is a part for every choice of $a$, and so is a
bouquet of at least two cycles of length at least $3$ with its center. By Corollary~\ref{cor:gap}, $\sigma<\Sigma$.

We use the reduced Bethe function $\phi$ of \eqref{real:eq:phi} with $D$
replaced by $D-z$. The formula of Lemma~\ref{real:lem:phi}(b) holds at every
point, not only at positive solutions.

\begin{lemma}\label{app:bridge:lem-hess}
Let $G$ be a loopless multigraph and $V\in(0,\infty)^V$. Put
$f_e=f_+(w_e^2V_iV_j)$ and $u_e=f_e/(1+f_e)$. Then
$2\diag(V)\Hess\phi(V)\diag(V)=H_u$, with $H_u$ as in \eqref{real:eq:Hu}.
\end{lemma}

\begin{proof}
Let $x_e=w_e^2V_iV_j$ for $e=ij$. As in the proof of
Lemma~\ref{real:lem:phi}, $\partial_{V_i}g(x_e)=f_e/V_i$, so
$2V_i^2\partial^2_{V_iV_i}\phi=1+\sum_{e\ni i}(f_e-x_ef_+'(x_e))$ and
$2V_iV_j\partial^2_{V_iV_j}\phi=-\sum_{e=ij}x_ef_+'(x_e)$ for $i\ne j$.
Since $f_e(1+f_e)=x_e$ we have $x_ef_+'(x_e)=f_e(1+f_e)/(1+2f_e)$ and
$f_e-x_ef_+'(x_e)=f_e^2/(1+2f_e)$, and with $f_e=u_e/(1-u_e)$ these are the
entries of \eqref{real:eq:Hu}.
\end{proof}

\begin{lemma}\label{app:bridge:lem-rho}
Let $(G,a)$ be a part and let $\widehat G\supseteq G$ be a loopless
multigraph with coefficients extending those of $G$, such that every edge
of $\widehat G$ with an end in $V^\circ$ is an edge of $G$. Then every
positive solution $(z,V,f)$ of the covariance system of $\widehat G$
satisfies $z<\rho$. In
particular $\Sigma<\rho$.
\end{lemma}

\begin{proof}
Let $C$ be a connected component of $G-a$; it has cycle rank at most one,
and it is joined to $a$ by at least one edge because $G$ is connected. For
$j\in C$ the vertex equation reads
$(D'_j-z)V_j=1+\sum_{e\ni j,\ e\subseteq C}f_e$ with
$D'_j=D_j-\sum_{e=ja}f_e/V_j$, so $(V_C,f_C)$ is a positive solution of the
covariance system of $C$ with the potential $D'$. If $C$ is a single
vertex, this gives $z<D'_j$. Otherwise Lemma~\ref{lem:domain} gives
$z\le\Sigm(C,D')$, and $\Sigm(C,D')=\lambda_{\min}((D'-\abs A)[C])$: for a
tree the covering tree is $C$ itself and the signs of the weights can be
removed, and for cycle rank one this is
Proposition~\ref{prop:unicyclic}(a). In both cases
$z\le\lambda_{\min}((D'-\abs A)[C])$. Now $D'\le D$ on $C$, with strict
inequality at the neighbors of $a$, and $(D-\abs A)[C]$ is irreducible with
a positive Perron vector $x$; hence
$\lambda_{\min}((D'-\abs A)[C])\le x^{\T}(D'-\abs A)x/\norm x^2<\lambda_{\min}((D-\abs A)[C])$.
Taking the minimum over $C$ gives $z<\rho$. The last statement follows by
applying this with $\widehat G=G$ to the turning point of
Lemma~\ref{lem:turning}(ii), a positive solution over $\Sigma$.
\end{proof}

For $z<\rho$ and $s>0$ the \emph{slice system} of a part consists of the
vertex equations at $j\in V^\circ$, with $f_e=f_+(w_e^2V_iV_j)$ and $V_a=s$
fixed:
\begin{equation}\label{app:bridge:eq-slice}
F_j(V^\circ):=(D_j-z)V_j-1-\sum_{e\ni j,\ e=jj'}f_+\bigl(w_e^2V_jV_{j'}\bigr)=0\qquad(j\in V^\circ).
\end{equation}

\begin{lemma}\label{app:bridge:lem-slice}
Let $(G,a)$ be a part, $z<\rho$ and $s>0$. The slice system has exactly one
solution $V^\circ(z,s)\in(0,\infty)^{V^\circ}$. Its Jacobian matrix
$J^\circ=\partial F/\partial V^\circ$ is invertible with
$\det J^\circ>0$ and $(J^\circ)^{-1}\ge0$ entrywise, and $V^\circ(z,s)$ is a
real analytic function of $(z,s)$, nondecreasing in $z$ and in $s$.
\end{lemma}

\begin{proof}
\emph{Convexity.} By Lemma~\ref{app:bridge:lem-hess}, at every
$V\in(0,\infty)^V$ the matrix $2\diag(V)\Hess\phi\,\diag(V)$ restricted to
$V^\circ$ is $H_u[V^\circ]$, which is the matrix \eqref{real:eq:Hu} of the
pseudoforest $G-a$ plus the nonnegative diagonal matrix with entries
$\sum_{e=ja}u_e^2/(1-u_e^2)$. It is positive definite: for simple $G$ by
Lemma~\ref{real:lem:bass}(b),(c), and for a dipole because it is the
$1\times1$ matrix $1+\sum_eu_e^2/(1-u_e^2)$. Hence $V^\circ\mapsto\phi(V^\circ,s)$
is strictly convex on the convex set $(0,\infty)^{V^\circ}$. Since
$F_j=2V_j\partial_{V_j}\phi$ (Lemma~\ref{real:lem:phi}(a)), the solutions
of the slice system are its critical points, so there is at most one.

\emph{Existence.} The matrix $M=M(z)$ is positive definite with nonpositive
off-diagonal entries, so $M^{-1}\ge0$ with positive
diagonal~\cite{BermanPlemmons}. Let $m=\min_j(M^{-1}\ones)_j>0$,
$c=m^{-1/2}$, $y=M^{-1}(\sqrt s\,b+c\ones)$ and $\tilde V=y^2$ entrywise.
Then $y_j\ge cm$, and $f_+(x)\le\sqrt x$ gives
\[
F_j(\tilde V)\ge y_j\Bigl((D_j-z)y_j-\sum_{e\ni j,\ j'\in V^\circ}\abs{w_e}y_{j'}-\sqrt s\,b_j\Bigr)-1=cy_j-1\ge c^2m-1=0 .
\]
The map $\mathcal M(V)_j=(1+\sum_{e\ni j}f_+(w_e^2V_jV_{j'}))/(D_j-z)$ is
nondecreasing (note $z<\rho\le\min_{j\in V^\circ}D_j$), and
$\mathcal M(\tilde V)\le\tilde V$. Its iterates from $0$ increase, stay
below $\tilde V$, and converge to a fixed point $V^\circ\ge\mathcal M(0)>0$,
which solves \eqref{app:bridge:eq-slice}.

\emph{Derivatives.} Differentiating $F_j=2V_j\partial_{V_j}\phi$ and using
$F_j=0$ for $j\in V^\circ$ (only these rows occur in $J^\circ$), we get at
the solution, by Lemma~\ref{app:bridge:lem-hess},
$J^\circ=2\diag(V^\circ)\Hess_{V^\circ}\phi=H_u[V^\circ]\diag(V^\circ)^{-1}$.
Its determinant is positive, and
$(J^\circ)^{-1}=\diag(V^\circ)H_u[V^\circ]^{-1}\ge0$ because
the inverse of a positive definite matrix with nonpositive off-diagonal
entries is nonnegative. The implicit function theorem gives analyticity,
and $\partial_zV^\circ=(J^\circ)^{-1}V^\circ\ge0$,
$\partial_sV^\circ=-(J^\circ)^{-1}\partial_sF\ge0$ because
$\partial_sF_j=-\sum_{e=ja}f_+'(w_e^2sV_j)w_e^2V_j\le0$.
\end{proof}

\begin{definition}\label{app:bridge:def-response}
The \emph{response function} of the part $(G,a)$ is, for $z<\rho$ and
$s>0$,
\[
\Phi(z,s)=(D_a-z)s-\sum_{e\ni a,\ e=aj}f_+\bigl(w_e^2sV^\circ_j(z,s)\bigr),
\]
and $\mathcal J(z)=\{s>0:\Phi(z,s)>1\}$.
\end{definition}

The vertex equation at $a$ of $G$ reads $\Phi(z,V_a)=1$; in the joined
graph below it will read $\Phi(z,V_a)=1+f_{\mathrm{bridge}}$.

\begin{lemma}[Response at infinity]\label{app:bridge:lem-response}
Let $(G,a)$ be a part.
\begin{enumerate}
\item[\textup{(a)}] $\partial_z\Phi\le-s<0$, and $\Phi(z,s)\to0$ as
$s\to0$, locally uniformly in $z<\rho$.
\item[\textup{(b)}] The function $\Psi(z,u)=u\,\Phi(z,1/u)$, $u>0$,
extends to a real analytic function on an open set containing
$(-\infty,\rho)\times[0,\infty)$, with
\[
\Psi(z,0)=S(z):=D_a-z-b^{\T}M(z)^{-1}b .
\]
The function $S$ is strictly decreasing, its only zero in $(-\infty,\rho)$
is $\sigma$, $g:=-S'(\sigma)=1+\norm{M(\sigma)^{-1}b}^2>0$, and
$\partial_u\Psi(\sigma,0)=\beta(G)$. Thus
$\Phi(z,s)=S(z)s+\beta(G)+O(\abs{z-\sigma})+O(1/s)$ as $s\to\infty$, locally
uniformly in $z$.
\end{enumerate}
\end{lemma}

\begin{proof}
(a) $\partial_z\Phi=-s-\sum_ef_+'\,w_e^2s\,\partial_zV^\circ_j\le-s$ by
Lemma~\ref{app:bridge:lem-slice}. For $s\le1$,
$V^\circ(z,s)\le V^\circ(z,1)$, so
$0\le\sum_ef_+(w_e^2sV^\circ_j)\le\sum_e\abs{w_e}\sqrt{sV^\circ_j(z,1)}\to0$.

(b) Put $u=1/s$, $Y_j=uV_j$ for $j\in V^\circ$, $Y_a=1$, and
$y_e=u(1+2f_e)$. Since $(1+2f_e)^2=1+4w_e^2V_iV_{i'}$, the slice system
with $f_e=f_+$ is, for $u>0$, equivalent to $y_e=\sqrt{u^2+4w_e^2Y_iY_{i'}}$
and
\[
P_j(z,u,Y):=(D_j-z)Y_j-u\Bigl(1-\frac{d_j}2\Bigr)-\frac12\sum_{e\ni j}\sqrt{u^2+4w_e^2Y_jY_{j'}}=0
\qquad(j\in V^\circ),
\]
and $\Psi=D_a-z-\frac12\sum_{e\ni a}(y_e-u)$. At $u=0$ the positive
solutions are $Y=x^2$ with $x_j((D-z-\abs A)x)_j=0$ for $j\in V^\circ$ and
$x_a=1$, that is $x_{V^\circ}=M(z)^{-1}b$, which is positive since $M(z)$ is
block diagonal over the components of $G-a$, each block is an irreducible
nonsingular M-matrix and so has an entrywise positive
inverse~\cite{BermanPlemmons}, and $b\ge0$ is nonzero on every component. Using
$(D_j-z)x_j=\sum_{e\ni j}\abs{w_e}x_{j'}$ one finds
$\partial P/\partial Y=\frac12\diag(x)M(z)\diag(x)^{-1}$ at $u=0$, which is
invertible. By the implicit function theorem, and by uniqueness of the
local solutions, $Y$ extends to a real analytic function of $(z,u)$ on a
neighborhood of $(-\infty,\rho)\times\{0\}$. For small $u>0$ it gives
$V_j=Y_j/u>0$ and $f_e=(y_e/u-1)/2>0$ solving the slice system, so by
Lemma~\ref{app:bridge:lem-slice} it is $V^\circ(z,1/u)$. Hence $\Psi$
extends analytically, and $\Psi(z,0)=D_a-z-\sum_{e=aj}\abs{w_e}x_j=S(z)$.

We have $S'(z)=-1-b^{\T}M^{-2}b<0$ and $\det(D-z-\abs A)=S(z)\det M(z)$
with $\det M(z)>0$ on $(-\infty,\rho)$. By Lemma~\ref{app:bridge:lem-rho}
and Corollary~\ref{cor:gap}, $\sigma<\Sigma<\rho$, so $S(\sigma)=0$, and
$\sigma$ is the only zero since $S$ is strictly decreasing.

For the first-order term write $Y_j=x_j^2+u\,x_j\omega_j+O(u^2)$, with
$\omega_a=0$. Then
$y_e=2\abs{w_e}x_jx_{j'}+u\abs{w_e}(x_j\omega_{j'}+x_{j'}\omega_j)+O(u^2)$,
and the terms of order $u$ in $P_j$, together with
$(D_j-z)x_j=\sum_{e\ni j}\abs{w_e}x_{j'}$, give
$((D-z-\abs A)\omega)_j=(2-d_j)/x_j$ for $j\in V^\circ$; the expression for
$\Psi$ gives $\partial_u\Psi(z,0)=\frac12((D-z-\abs A)\omega)_a+\frac{d_a}2$.
At $z=\sigma$ the matrix $L=D-\sigma-\abs A$ satisfies $Lx=0$, the entry
at $a$ being $S(\sigma)=0$. Since $L$ is symmetric,
$0=x^{\T}L\omega=\sum_{j\in V^\circ}(2-d_j)+(L\omega)_a$, hence
$\partial_u\Psi(\sigma,0)=\frac12\bigl(\sum_{j\in V}d_j-2(n-1)\bigr)=\abs E-n+1=\beta(G)$.
The expansion of $\Phi=s\Psi(z,1/s)$ follows.
\end{proof}

\begin{lemma}[Roots of one part]\label{app:bridge:lem-roots}
Let $(G,a)$ be a part and $z<\rho$.
\begin{enumerate}
\item[\textup{(a)}] The map $s\mapsto(V^\circ(z,s),s)$ is a bijection from
$\{s>0:\Phi(z,s)=1\}$ onto the set of positive solutions of $G$ over $z$, and
the solution is nondegenerate ($J_V$ invertible) if and only if
$\partial_s\Phi(z,s)\ne0$.
\item[\textup{(b)}] For $z\le\sigma$ the equation $\Phi(z,\cdot)=1$ has
exactly one root $L(z)$, and $\partial_s\Phi(z,L(z))>0$. For
$\sigma<z<\Sigma$ it has exactly two roots $L(z)<U(z)$, with
$\partial_s\Phi(z,L)>0>\partial_s\Phi(z,U)$. For $z=\Sigma$ it has exactly
one root $s^*$, at which $\partial_s\Phi=0$, and for $\Sigma<z<\rho$ it has
none. Here $L(z)$ is the $a$-entry of the physical branch for $z<\Sigma$, $L$
and $U$ are real analytic, $U(z)\to\infty$ as $z\downarrow\sigma$, and
$L(z),U(z)\to s^*$ as $z\uparrow\Sigma$.
\item[\textup{(c)}] $\mathcal J(z)=(L(z),\infty)$ for $z\le\sigma$,
$\mathcal J(z)=(L(z),U(z))$ for $\sigma<z<\Sigma$, and $\mathcal J(z)=\emptyset$
for $\Sigma\le z<\rho$.
\end{enumerate}
\end{lemma}

\begin{proof}
(a) A positive solution restricts on $V^\circ$ to a positive solution of the
slice system with $s=V_a$, hence to $V^\circ(z,s)$, and its vertex
equation at $a$ is $\Phi(z,s)=1$; conversely such $s$ gives a positive
solution. Order the coordinates as $(V^\circ,V_a)$. The Schur complement
of $J^\circ$ in $J_V$ is the derivative of $s\mapsto g_a(z,V^\circ(z,s),s)=\Phi(z,s)-1$,
so $\det J_V=\det J^\circ\cdot\partial_s\Phi$ with $\det J^\circ>0$.

(b) We first note that the conclusions of
Theorem~\ref{real:thm:packing}(b) hold for every part. For simple $G$ this
follows because $G-a$ is a pseudoforest and hence $p_U(G)\le1$.
For a dipole with vertices $0,1$
put $t_i=D_i-z$; as in the proof of Proposition~\ref{prop:circles}(d),
positive solutions correspond to the zeros $\gamma>0$ of the strictly convex
function $\varphi(\gamma)=1+\sum_ef_+(w_e^2\gamma^2/(t_0t_1))-\gamma$, with
$t_0V_0=t_1V_1=\gamma$. Writing $c=\sum_ef_+'(w_e^2V_0V_1)w_e^2$, a direct
computation gives $\det J_V=t_0t_1-c(t_0V_0+t_1V_1)=-t_0t_1\varphi'(\gamma)$.
Since $\varphi(0)=1>0$, there are at most two positive solutions, and if
there are two, both are simple zeros of $\varphi$ and hence nondegenerate.
The proof of Theorem~\ref{real:thm:packing}(b) uses part (a) of that
theorem only through these two facts, and its other ingredients
(Theorem~\ref{thm:unique-ws}, Lemmas~\ref{lem:domain}, \ref{lem:blowup}
and~\ref{lem:turning}, Corollary~\ref{cor:two-roots}) hold for loopless
multigraphs; so it applies verbatim.

Thus $N(z)=1$ for $z\le\sigma$, $N(z)=2$ for $\sigma<z<\Sigma$, $N(\Sigma)=1$
and $N(z)=0$ for $z>\Sigma$; all positive solutions over $z<\Sigma$ are
nondegenerate; and the two branches converge to the turning point $p^*$
as $z\uparrow\Sigma$, whose $J_V$ is singular by Lemma~\ref{lem:turning}(ii).
By (a), $\Phi(z,\cdot)=1$ has the stated number of roots, simple for
$z<\Sigma$ and with $\partial_s\Phi=0$ at $z=\Sigma$. By
Lemma~\ref{app:bridge:lem-response}, $\Phi(z,0^+)=0$, while as
$s\to\infty$ the function $\Phi(z,s)$ tends to $+\infty$ for $z<\sigma$, to
$\beta(G)\ge2$ for $z=\sigma$, and to $-\infty$ for $z>\sigma$. Since all
roots are simple for $z<\Sigma$, the signs of $\partial_s\Phi$ follow, and
for $z>\Sigma$ we get $\Phi<1$. Analyticity of $L$, $U$ follows from the
implicit function theorem. The $a$-entry of the physical branch is a
continuous root on $(-\infty,\Sigma)$; it equals $L$ for $z\le\sigma$,
it stays bounded near $\sigma$ while $U\to\infty$ (shown next), and $L<U$
on $(\sigma,\Sigma)$, so it equals $L$ throughout. Let $z_m\downarrow\sigma$ and suppose that
$U(z_m)$ stays bounded. By Lemma~\ref{app:bridge:lem-slice},
$V^\circ(z_m,U(z_m))$ is then bounded, and a subsequence of the solutions
converges to a positive solution over $\sigma$ (note $V_i\ge1/(D_i-z)$),
which is the physical one. This point is nondegenerate, so near it the
positive solutions form the physical branch, and $U(z_m)=L(z_m)$ for
large $m$, a contradiction.

(c) follows from (b) and the sign of $\Phi-1$ near $0$.
\end{proof}

\subsection{Two parts joined by an edge}\label{app:bridge:subsec-reduction}

Let $(G_1,a_1)$ and $(G_2,a_2)$ be two parts with disjoint vertex sets, and
let $G_w$ be the graph obtained from $G_1\sqcup G_2$ by adding an edge
$a_1a_2$ of weight $w\ne0$; we write $\kappa=w^2$. We use the notation of
Section~\ref{app:bridge:subsec-parts} with an index $k\in\{1,2\}$:
$\Phi_k$, $\Psi_k$, $S_k$, $g_k$, $\beta_k=\beta(G_k)$, $\rho_k$,
$\sigma_k$, $\Sigma_k$, $L_k$, $U_k$, $\mathcal J_k$. Put
$z_+=\min(\Sigma_1,\Sigma_2)$.

\begin{lemma}[Reduction]\label{app:bridge:lem-reduction}
For every $\kappa>0$ the map $(z,V)\mapsto(z,s,t)=(z,V_{a_1},V_{a_2})$ is a
homeomorphism, real analytic with real analytic inverse, from $\Xp(G_w)$
onto
\[
\mathcal C_\kappa=\bigl\{(z,s,t):\ z<\min(\rho_1,\rho_2),\ s,t>0,\
\Phi_1(z,s)=\Phi_2(z,t)=:\nu>1,\ \nu(\nu-1)=\kappa st\bigr\}.
\]
Moreover $\mathcal C_\kappa\subseteq\{z<z_+\}$, and $(z,V)$ is regular if
and only if the $2\times3$ Jacobian matrix of
$E_1=\Phi_1(z,s)-\Phi_2(z,t)$ and $E_2=\Phi_1(z,s)(\Phi_1(z,s)-1)-\kappa st$
with respect to $(z,s,t)$ has rank $2$ at the image point.
\end{lemma}

\begin{proof}
By Lemma~\ref{app:bridge:lem-rho}, applied to each part with $\widehat G=G_w$,
every point of $\Xp(G_w)$ has $z<\min\rho_k$. The restriction of a
positive solution to $V_k^\circ$ solves the slice system of $G_k$ with
$V_{a_k}$ given, so it equals $V^\circ_k(z,V_{a_k})$. The vertex equations
at $a_1$ and $a_2$ read $\Phi_1(z,s)=1+f$ and $\Phi_2(z,t)=1+f$ with
$f=f_+(\kappa st)>0$; so $\nu=1+f>1$ and $\nu(\nu-1)=\kappa st$. Conversely,
for $(z,s,t)\in\mathcal C_\kappa$ the number $\nu-1>0$ is the positive root of
$x(1+x)=\kappa st$, and $(V^\circ_1(z,s),s,t,V^\circ_2(z,t))$ is a positive
solution. Both maps are real analytic. Since $\nu>1$ means
$\mathcal J_k(z)\ne\emptyset$, Lemma~\ref{app:bridge:lem-roots}(c) gives
$z<\Sigma_k$.

For the regularity, order the columns of $(\partial_zg\mid J_V)$ as
$(V_1^\circ,V_2^\circ\mid z,s,t)$ and the rows as $(V_1^\circ,V_2^\circ\mid a_1,a_2)$.
The upper left block is $\diag(J_1^\circ,J_2^\circ)$, invertible by
Lemma~\ref{app:bridge:lem-slice}, so the rank is
$\abs{V_1^\circ}+\abs{V_2^\circ}$ plus the rank of the Schur complement,
which by the chain rule is the Jacobian matrix of
$G_{a_1}=\Phi_1(z,s)-1-f_+(\kappa st)$ and $G_{a_2}=\Phi_2(z,t)-1-f_+(\kappa st)$.
At a point of $\mathcal C_\kappa$ we have $f_+'(\kappa st)=1/(2\nu-1)$, so
$(2\nu-1)\,dG_{a_1}=dE_2$ and $dE_1=dG_{a_1}-dG_{a_2}$: the two pairs of
differentials span the same space.
\end{proof}

\begin{corollary}[Level sets]\label{app:bridge:cor-level}
Let
$\mathcal S=\{(z,s,t):z<\min\rho_k,\ s,t>0,\ \Phi_1(z,s)=\Phi_2(z,t)>1\}$
and $\Lambda(z,s,t)=\nu(\nu-1)/(st)$ with $\nu=\Phi_1(z,s)$. For every
$\kappa>0$, $\mathcal C_\kappa=\Lambda^{-1}(\kappa)\cap\mathcal S$. At a
point $p\in\mathcal C_\kappa$ where $d(\Phi_1-\Phi_2)\ne0$, the set
$\mathcal S$ is a smooth surface near $p$, and the corresponding point of
$\Xp(G_w)$ is regular if and only if $p$ is not a critical point of
$\Lambda|_{\mathcal S}$.
\end{corollary}

\begin{proof}
On $\mathcal S$ we have $E_2=st(\Lambda-\kappa)$, hence $dE_2=st\,d\Lambda$
on $\Lambda^{-1}(\kappa)$. The rank of $(dE_1,dE_2)$ is $2$ if and only if
$d\Lambda$ does not vanish on $T_p\mathcal S=\ker dE_1$.
\end{proof}

\subsection{The weak bridge limit}\label{app:bridge:subsec-weak}

\begin{theorem}[Weak bridge]\label{app:bridge:thm-gluing}
Let $(G_1,a_1)$, $(G_2,a_2)$ be parts with arbitrary real coefficients
such that $\sigma_1,\sigma_2,\Sigma_1,\Sigma_2$ are distinct, numbered so
that $\sigma_1<\sigma_2$. There is $\bar\kappa>0$ such that for
$0<\kappa=w^2<\bar\kappa$:
\begin{enumerate}
\item[\textup{(i)}] every point of $\Xp(G_w)$ is regular, so $\Xp(G_w)$ is
a smooth curve, and exactly one connected component of $\Xp(G_w)$ is not
compact;
\item[\textup{(ii)}] $\Xp(G_w)$ has exactly one circle if
$\sigma_2<\Sigma_1<\Sigma_2$, and none if $\Sigma_2<\Sigma_1$ or
$\Sigma_1<\sigma_2$;
\item[\textup{(iii)}] if $\sigma_2<\Sigma_1<\Sigma_2$, then as $\kappa\to0$
the image of the circle under $(z,V)\mapsto(z,V_{a_1},V_{a_2}/(1+V_{a_2}))$
converges in the Hausdorff distance to the closure of
\[
\{(z,X_1(z),\hat U_2(z)):\ \sigma_2<z\le\Sigma_1,\ X\in\{L,U\}\}\ \cup\
\{(\sigma_2,s,1):\ L_1(\sigma_2)\le s\le U_1(\sigma_2)\},
\]
where $\hat U_2=U_2/(1+U_2)$; in particular its $z$-projection converges to
$[\sigma_2,\Sigma_1]$. Moreover, for every compact interval
$I\subset\mathcal J_1(\sigma_2)=(L_1(\sigma_2),U_1(\sigma_2))$ there are
$\delta,\kappa_I>0$ such that for $0<\kappa<\kappa_I$ the points of the
circle with $V_{a_1}\in I$ and $V_{a_2}>1/\delta$ form the graph of a real
analytic map $s\mapsto(z,V_{a_2})$ on $I$ with
\begin{equation}\label{app:bridge:eq-Z}
z=\sigma_2+\kappa Z(s)+O(\kappa^2),\qquad
\frac1{V_{a_2}}=\frac{\kappa s}{\varphi(s)(\varphi(s)-1)}+O(\kappa^2),\qquad
Z(s)=\frac{s\,(\beta_2-\varphi(s))}{g_2\,\varphi(s)(\varphi(s)-1)},
\end{equation}
with $\varphi(s)=\Phi_1(\sigma_2,s)$, the error terms being bounded together
with their $s$-derivatives uniformly on $I$.
\end{enumerate}
\end{theorem}

Thus the circle consists of the fiber products of the two positive
solutions of $G_1$ with the second solution of $G_2$, joined at $\Sigma_1$
by the turning point of $G_1$ and near $\sigma_2$ by a piece on which $G_2$
is close to its Perron pole and $z-\sigma_2=O(\kappa)$. Since
$\varphi\to1$ at both ends of $\mathcal J_1(\sigma_2)$ and $\beta_2\ge2$,
$Z\to+\infty$ there; no further condition is needed for this piece to join
the two sheets. The constant $\bar\kappa$ depends on the coefficients of
the parts and on the distances between the four thresholds; we do not
estimate it.

The rest of this subsection is the proof. We work in the space $\R^3$ with
coordinates $(z,\hat s,\hat t)=(z,s/(1+s),t/(1+t))$, in which points with
$s=\infty$ or $t=\infty$ have $\hat s=1$ or $\hat t=1$; near such points
$v=1/s$ and $u=1/t$ are analytic coordinates, also for $v,u\le0$. We fix
$Z_*<\min(Z_0,\sigma_1-1)$, with $Z_0$ as in Lemma~\ref{app:bridge:lem-left}
below.

\emph{The skeleton.} Let
$F=\{(z,s,t):z\le z_+,\ s,t>0,\ \Phi_1(z,s)=1=\Phi_2(z,t)\}$. By
Lemma~\ref{app:bridge:lem-roots}, $F$ is the union of the four sets
$F_{XY}=\{(z,X_1(z),Y_2(z))\}$, $X,Y\in\{L,U\}$, over the values of $z\le z_+$
where they are defined, together with the turning points of the part with
$\Sigma_k=z_+$. Let
\[
\mathcal D_2=\{(\sigma_2,s,\infty):\Phi_1(\sigma_2,s)\ge1\},\qquad
\mathcal D_1=\{(\sigma_1,\infty,t):\Phi_2(\sigma_1,t)\ge1\}\cup\{P_\infty\},
\]
with $P_\infty=(\sigma_1,\infty,\infty)$, where $\mathcal D_2=\emptyset$ if
$\sigma_2\ge\rho_1$,
and $\mathcal K=F\cup\mathcal D_1\cup\mathcal D_2$,
$\mathcal K'=\mathcal K\cap\{z\ge Z_*\}$. Since $\sigma_1<\sigma_2$, we
have $S_2(\sigma_1)>0$ and $\mathcal D_1$ is the segment
$t\in[L_2(\sigma_1),\infty]$. If $\sigma_2<\Sigma_1$ then $\mathcal D_2$ is
the segment $s\in[L_1(\sigma_2),U_1(\sigma_2)]$, and if $\Sigma_1<\sigma_2$
it is empty. The points $(\sigma_2,L_1(\sigma_2),\infty)$,
$(\sigma_2,U_1(\sigma_2),\infty)$ and $(\sigma_1,\infty,L_2(\sigma_1))$ are
the \emph{nodes}; at them $F_{LU}$, $F_{UU}$ and $F_{UL}$, respectively,
end as $z\downarrow\sigma_2$ or $z\downarrow\sigma_1$, by
Lemma~\ref{app:bridge:lem-roots}(b). From that lemma one reads off the
connected components of the compact set $\mathcal K'$:
\begin{itemize}
\item \emph{interlaced} ($\sigma_2<\Sigma_1<\Sigma_2$, so $z_+=\Sigma_1$):
$\mathcal K'_1=F_{LL}\cup F_{UL}\cup\mathcal D_1$ (cut at $Z_*$), where
$F_{LL}$ and $F_{UL}$ meet at the turning point of $G_1$ over $\Sigma_1$
and $F_{UL}$ meets $\mathcal D_1$ at a node; and
$\mathcal K'_2=F_{LU}\cup F_{UU}\cup\mathcal D_2$, where $F_{LU}$ and
$F_{UU}$ meet at the turning point over $\Sigma_1$ and end at the two nodes
of $\mathcal D_2$. These two compact connected sets are disjoint.
\item \emph{nested} ($\Sigma_2<\Sigma_1$, so $z_+=\Sigma_2$ and
$\sigma_1<\sigma_2<\Sigma_2<\Sigma_1$): $\mathcal K'$ is connected; it runs
$F_{LL}$, turning point of $G_2$, $F_{LU}$, node, $\mathcal D_2$, node,
$F_{UU}$, turning point of $G_2$, $F_{UL}$, node, $\mathcal D_1$, $P_\infty$.
\item \emph{separated} ($\Sigma_1<\sigma_2$): $F_{\cdot U}=\mathcal D_2=\emptyset$,
and $\mathcal K'=F_{LL}\cup F_{UL}\cup\mathcal D_1$ is connected.
\end{itemize}
In each case, the component containing the end $z=Z_*$ of $F_{LL}$ also
contains $P_\infty$.

\begin{lemma}[Left end]\label{app:bridge:lem-left}
There is $Z_0$ such that for all $\kappa\in[0,1]$ and $z<Z_0$ the set
$\Xp(G_w)$ (for $\kappa=0$, of $G_1\sqcup G_2$) has exactly one point over $z$; it is regular, depends
continuously on $(z,\kappa)\in(-\infty,Z_0)\times[0,1]$, and for $\kappa=0$
its image is $(z,L_1(z),L_2(z))$.
\end{lemma}

\begin{proof}
Let $\omega=\max(1,\max_e\abs{w_e})$ over the edges of the parts, $d$ the
maximal degree of $G_w$, $D_0=\min_iD_i$, and $z<D_0-d\omega$. If $V$ is a
positive solution and $V_i=\max V$, then $(D_i-z)V_i\le1+d\omega V_i$ by
$f_+(x)\le\sqrt x$, so $V\in Q_z=[0,1/(D_0-z-d\omega)]^V$. The map
$\Phi_z$ of \eqref{per:eq:Phi} maps $Q_z$ into itself, and its Jacobian
matrix has row sums at most $2d\omega^2/((D_0-z)(D_0-z-d\omega))$ on $Q_z$
since $0<f_+'\le1$. For $z<Z_0$, with $Z_0$ depending only on $D_0,d,\omega$,
this is less than $\frac12$, so $\Phi_z$ is a contraction of the convex set
$Q_z$ in the maximum norm, uniformly in $\kappa\in[0,1]$. Its unique fixed
point is positive and depends continuously on $(z,\kappa)$; the spectral
radius of $D\Phi_z$ is less than $1$ there, so
$J_V=\diag(D-z)(I-D\Phi_z)$ is invertible. For $\kappa=0$ the fixed point is
the pair of the physical solutions of the parts.
\end{proof}

\begin{lemma}[Compactness]\label{app:bridge:lem-compact}
For every open set $N\supseteq\mathcal K'$ of $\R^3$ there is $\kappa_N>0$
such that $\mathcal C_\kappa\cap\{z\ge Z_*\}\subseteq N$ for
$0<\kappa<\kappa_N$.
\end{lemma}

\begin{proof}
Otherwise there are $\kappa_m\to0$ and $p_m=(z_m,s_m,t_m)\in\mathcal C_{\kappa_m}\setminus N$
with $z_m\in[Z_*,z_+]$ (Lemma~\ref{app:bridge:lem-reduction}); pass to a
subsequence with $z_m\to z^*$, $s_m\to s^*\in[0,\infty]$,
$t_m\to t^*\in[0,\infty]$, and write $\nu_m=\Phi_1(z_m,s_m)>1$,
$v_m=1/s_m$, $u_m=1/t_m$. If $s^*=0$ or $t^*=0$, then $\nu_m\to0$ by
Lemma~\ref{app:bridge:lem-response}(a), which is impossible. If
$s^*,t^*\in(0,\infty)$, then $\nu_m(\nu_m-1)=\kappa_ms_mt_m\to0$, so the
limit lies in $F$. If $s^*=\infty$ and $t^*\in(0,\infty)$, then
$\nu_m=\Phi_2(z_m,t_m)\to\nu^*\ge1$ and $\Psi_1(z_m,v_m)=\nu_mv_m\to0$, so
$S_1(z^*)=0$, $z^*=\sigma_1$, and the limit lies in $\mathcal D_1$; the
case $s^*<\infty=t^*$ is symmetric and gives a point of $\mathcal D_2$
(and is impossible if $\Sigma_1<\sigma_2$). If $s^*=t^*=\infty$ and
$\nu_m$ has a bounded subsequence, then $S_1(z^*)=S_2(z^*)=0$, so
$\sigma_1=\sigma_2$, a contradiction; hence $\nu_m\to\infty$ and
$\Psi_1(z_m,v_m)\Psi_2(z_m,u_m)=\nu_m^2u_mv_m=\kappa_m\nu_m/(\nu_m-1)\to0$.
Both factors are positive, so $S_1(z^*),S_2(z^*)\ge0$ and one of them
vanishes. Since $S_1(\sigma_2)<0$, we get $z^*=\sigma_1$, and the limit is
$P_\infty$. In all cases the limit lies in $\mathcal K'\subseteq N$, which
contradicts $p_m\notin N$ for large $m$.
\end{proof}

We now describe $\mathcal C_\kappa$ near each point of $\mathcal K'$. A
\emph{graph box} is a compact set $\bar B=\psi([-1,1]\times[-1,1]^2)$, where
$\psi$ is an analytic diffeomorphism from a neighborhood of the cube onto
an open subset of $\R^3$, together with $\kappa_B>0$ and continuous maps
$\gamma_\kappa\colon[-1,1]\to(-1,1)^2$ for $0\le\kappa<\kappa_B$, such that
$\psi^{-1}(\mathcal K\cap\bar B)$ is the graph of $\gamma_0$,
$\psi^{-1}(\mathcal C_\kappa\cap\bar B)$ is the graph of $\gamma_\kappa$ for
$\kappa>0$, $\gamma_\kappa\to\gamma_0$ uniformly, and all points of
$\mathcal C_\kappa\cap\bar B$ are regular in the sense of
Lemma~\ref{app:bridge:lem-reduction}. We write
$B=\psi((-1,1)^3)$; it is an open neighborhood of the points
$\psi(c,\gamma_0(c))$, $\abs c<1$.

\begin{lemma}[Local structure]\label{app:bridge:lem-local}
\begin{enumerate}
\item[\textup{(a)}] Every point of $\mathcal K'$ other than $P_\infty$ is of
the form $\psi(0,\gamma_0(0))$ for graph boxes of arbitrarily small
diameter.
\item[\textup{(b)}] There are an open neighborhood $B_\infty$ of $P_\infty$
with compact closure, $u_0>0$ and $\kappa_\infty>0$ such that for
$0<\kappa<\kappa_\infty$ the set $\mathcal C_\kappa\cap\bar B_\infty$ is
the image of a continuous injective map $c\mapsto\gamma_\kappa(c)$ on
$(0,u_0]$ which converges, uniformly on $(0,u_0]$, to
$c\mapsto(\sigma_1,\infty,1/c)$ as $\kappa\to0$, and along which
$s,t\to\infty$ as $c\to0$; its points are regular, those with $c<u_0$ lie
in $B_\infty$, and
$\mathcal K\cap\bar B_\infty=\{(\sigma_1,\infty,t):t\ge1/u_0\}\cup\{P_\infty\}$.
\item[\textup{(c)}] The expansion \eqref{app:bridge:eq-Z} holds for the
graph boxes at the points of $\mathcal D_2$ with $s\in I$.
\end{enumerate}
\end{lemma}

\begin{proof}
\emph{Points of $F$.} At $\kappa=0$ and $\nu=1$ we have $dE_2=d\Phi_1$, and
the Jacobian matrix of $(E_1,E_2)$ with respect to $(z,s,t)$ is
\[
\begin{pmatrix}\partial_z\Phi_1-\partial_z\Phi_2&\partial_s\Phi_1&-\partial_t\Phi_2\\
\partial_z\Phi_1&\partial_s\Phi_1&0\end{pmatrix}.
\]
Its minors for the columns $(z,s)$ and $(z,t)$ are
$-\partial_s\Phi_1\,\partial_z\Phi_2$ and $\partial_t\Phi_2\,\partial_z\Phi_1$,
and $\partial_z\Phi_k<0$. Both vanish only if
$\partial_s\Phi_1=\partial_t\Phi_2=0$, that is, by
Lemma~\ref{app:bridge:lem-roots}(b), at a point over $\Sigma_1=\Sigma_2$,
which is excluded. So at every $p\in F$ one of $c=t$, $c=s$ can serve as a
parameter, and $c=z$ can where $\partial_s\Phi_1\,\partial_t\Phi_2\ne0$
(the minor for the columns $(s,t)$), and the implicit function theorem, applied to $(E_1,E_2)$ as
functions of $(z,s,t,\kappa)$, gives a box around $p$, of arbitrarily small
size, in which for $0\le\kappa<\kappa_B$ the zero set of $(E_1,E_2)$ is the
graph of a continuous function of $c$, analytic in $(c,\kappa)$, with
$\nu\in(\frac12,\frac32)$ and with rank-$2$ Jacobian. For $\kappa=0$ the
zero set is $F$ near $p$, because $\nu(\nu-1)=0$ and $\nu$ is close to $1$.
For $\kappa>0$ the zeros have $\nu(\nu-1)=\kappa st>0$, hence $\nu>1$, so
they are exactly the points of $\mathcal C_\kappa$ in the box.

\emph{Interior points of $\mathcal D_2$.} Let $p=(\sigma_2,s_0,\infty)$ with
$\nu_0=\Phi_1(\sigma_2,s_0)>1$. In the coordinates $(z,s,u)$, $u=1/t$, the
functions
$E_1'=\Psi_2(z,u)-u\Phi_1(z,s)=-uE_1$ and
$E_2'=\Phi_1(\Phi_1-1)u-\kappa s=uE_2$ are analytic near $p$, also for
$u\le0$, and for $u>0$ they define $\mathcal C_\kappa$ with the same rank of
the Jacobian matrix. At $(p,\kappa=0)$ their Jacobian matrix with respect
to $(z,u)$ is
$\bigl(\begin{smallmatrix}-g_2&\beta_2-\nu_0\\0&\nu_0(\nu_0-1)\end{smallmatrix}\bigr)$
by Lemma~\ref{app:bridge:lem-response}(b), which is invertible. The
implicit function theorem gives a box in which the zero set, for
$0\le\kappa<\kappa_B$, is the graph of an analytic map
$(s,\kappa)\mapsto(z(s,\kappa),u(s,\kappa))$ with $z(s,0)=\sigma_2$ and
$u(s,0)=0$ (at $\kappa=0$, $E_2'=0$ forces $u=0$ because $\nu>1$, and then
$E_1'=S_2(z)=0$). For $\kappa>0$, $u=\kappa s/(\nu(\nu-1))>0$, so the graph
lies in $\mathcal C_\kappa$; and $\mathcal K$ meets the box only in
$\mathcal D_2$, since $\nu$ is close to $\nu_0>1$. Differentiating
$E_2'=0$ and $E_1'=0$ with respect to $\kappa$ at $\kappa=0$ gives
$\nu_0(\nu_0-1)\partial_\kappa u=s$ and
$-g_2\,\partial_\kappa z+(\beta_2-\nu_0)\,\partial_\kappa u=0$, which is
\eqref{app:bridge:eq-Z}; the remainder is analytic in $(s,\kappa)$, which
gives the uniformity on compact $I$. This also proves (c). Interior points
of $\mathcal D_1$ are treated in the same way with the roles of the parts
exchanged, in the coordinates $(z,v,t)$.

\emph{Nodes.} Let $p=(\sigma_2,s_0,\infty)$ with $\Phi_1(\sigma_2,s_0)=1$;
by Lemma~\ref{app:bridge:lem-roots}(b) applied at $\sigma_2\in(\sigma_1,\Sigma_1)$,
$\partial_s\Phi_1(\sigma_2,s_0)\ne0$. The map $(z,s,u)\mapsto(e,p_1,q)=(E_1',\Phi_1(z,s)-1,u)$
has Jacobian determinant $-g_2\,\partial_s\Phi_1(\sigma_2,s_0)\ne0$ at $p$,
so it is an analytic coordinate system near $p$. In these coordinates
$\mathcal C_\kappa$ is, near $p$, the set of points with $e=0$, $p_1>0$
(that is $\nu>1$), $q>0$ (that is $t>0$) and $(1+p_1)p_1q=\kappa s$, that is
\[
p_1q=\kappa\,m(p_1,q),\qquad m(p_1,q)=\frac{s(0,p_1,q)}{1+p_1},\qquad m(0,0)=s_0>0,
\]
while $\mathcal K$ is $\{e=0,\ q=0,\ p_1\ge0\}\cup\{e=0,\ p_1=0,\ q\ge0\}$:
the first set is $\mathcal D_2$, and the second is $F_{XU}$, since $e=0$,
$\Phi_1=1$ give $\Psi_2(z,u)=u$, that is $\Phi_2(z,t)=1$. Choose $r>0$ such
that on $\abs{p_1},\abs q\le2r$ we have $m_0\le m\le m_1$ with $m_0>0$ and
$\abs{\partial_{p_1}m}+\abs{\partial_qm}\le M_1$. Put $\xi=p_1-q$,
$\eta=p_1+q$ and $H(\xi,\eta)=(\eta^2-\xi^2)/4-\kappa m$, so that
$\mathcal C_\kappa=\{e=0,\ \eta>\abs\xi,\ H=0\}$ near $p$. For
$\abs\xi\le r$ and $\abs\xi<\eta\le3r$: if $\eta\le4\kappa M_1$ then
$H\le4\kappa^2M_1^2-\kappa m_0<0$ for $\kappa<m_0/(4M_1^2)$; if
$\eta\ge4\kappa M_1$ then $\partial_\eta H\ge\eta/2-\kappa M_1\ge\eta/4>0$;
and $H(\xi,3r)\ge2r^2-\kappa m_1>0$ for small $\kappa$. Hence for each
$\xi$ there is exactly one solution $\eta=\gamma_\kappa(\xi)$, it depends
continuously on $\xi$, $\partial_\eta H\ne0$ there (so the points are
regular, since $de\ne0$ and $dE_2'=(1+p_1)\,dH$ on $\{e=0,H=0\}$), and
$\abs\xi<\gamma_\kappa(\xi)\le(\xi^2+4\kappa m_1)^{1/2}\le\abs\xi+2(\kappa m_1)^{1/2}$.
So the cube $\{\abs\xi\le r,\ \abs\eta\le3r,\ \abs e\le r\}$ is a graph box
with parameter $\xi$, $\gamma_0(\xi)=\abs\xi$. The node of $\mathcal D_1$
is treated in the same way, with $\partial_t\Phi_2(\sigma_1,L_2(\sigma_1))>0$.

\emph{The point $P_\infty$.} In the coordinates $(z,v,u)$ the set
$\mathcal C_\kappa$ is defined, for $u,v>0$, by $\Psi_1(z,v)/v=\Psi_2(z,u)/u=\nu$
and $\nu(\nu-1)=\kappa/(uv)$. Put $\mu=\Psi_2(z,u)$, so that $\nu=\mu/u$;
then the equations are equivalent to $v=\kappa u\,\omega$ and
$\Psi_1(z,v)=\kappa\,\omega\mu$ with $\omega=1/(\mu(\mu-u))$. Write
$z=\sigma_1+\kappa\zeta$ and
$\Theta(\zeta,\tilde v;\kappa)=\kappa^{-1}\Psi_1(\sigma_1+\kappa\zeta,\kappa\tilde v)$,
which is analytic up to $\kappa=0$ because $\Psi_1(\sigma_1,0)=0$, with
$\Theta(\zeta,\tilde v;0)=-g_1\zeta+\beta_1\tilde v$. The system becomes
\[
\omega\,\Psi_2(\sigma_1+\kappa\zeta,u)\bigl(\Psi_2(\sigma_1+\kappa\zeta,u)-u\bigr)-1=0,\qquad
\Theta(\zeta,u\omega;\kappa)-\omega\,\Psi_2(\sigma_1+\kappa\zeta,u)=0,
\]
in the unknowns $(\zeta,\omega)$ with parameters $(u,\kappa)$. At
$(u,\kappa)=(0,0)$ it has the unique solution $\omega_0=\mu_0^{-2}$,
$\zeta_0=-1/(g_1\mu_0)$, where $\mu_0=S_2(\sigma_1)>0$, and its Jacobian
matrix with respect to $(\zeta,\omega)$ there is
$\bigl(\begin{smallmatrix}0&\mu_0^2\\-g_1&-\mu_0\end{smallmatrix}\bigr)$, with
determinant $g_1\mu_0^2\ne0$. So $(\zeta,\omega)$ are analytic functions of
$(u,\kappa)$ near $(0,0)$. For $\kappa>0$ the passage from $(z,v,u)$ to
$(\zeta,\omega,u)$ is a diffeomorphism on $u>0$, and the new system is
obtained from $(E_1,E_2)$ by multiplying with a matrix of functions of
determinant $-uv^2/\kappa^2\ne0$ (the second equation is $(v/\kappa)E_1$, the
first is $(uv/\kappa)(E_2-(\Phi_1+\Phi_2-1)E_1)$); so its points are
regular. A priori, let $(z,v,u)\in\mathcal C_\kappa$ with
$\abs{z-\sigma_1}$, $v$ and $u$ small. Then $\mu$ is close to $\mu_0$,
$\omega$ is close to $\omega_0$, $v=\kappa u\omega$, and
$\abs{S_1(z)}\le\abs{\Psi_1(z,v)}+Cv\le C'\kappa$, so $\zeta$ is bounded;
then $\Theta(\zeta,u\omega;\kappa)=-g_1\zeta+O(u+\kappa)$ shows that
$(\zeta,\omega)$ is close to $(\zeta_0,\omega_0)$, inside the domain of
uniqueness of the implicit function theorem. Hence for the box
$\bar B_\infty=\{\abs{z-\sigma_1}\le\delta_1,\ \abs v\le\delta_1,\ \abs u\le u_0\}$
with small $\delta_1$ and $u_0$, its interior $B_\infty$, and small $\kappa$, the
set $\mathcal C_\kappa\cap\bar B_\infty$ is
$\{(\sigma_1+\kappa\zeta(u,\kappa),\ \kappa u\,\omega(u,\kappa),\ u):0<u\le u_0\}$,
which has the stated properties with $c=u$.
\end{proof}

\begin{proof}[Proof of Theorem~\ref{app:bridge:thm-gluing}]
Let $\mathcal K'_1,\dots,\mathcal K'_\ell$ ($\ell\in\{1,2\}$) be the
components of $\mathcal K'$, $\mathcal K'_1$ the one containing $P_\infty$,
and let $\delta>0$ be less than a third of their mutual distances. Choose
graph boxes $\bar B_1,\dots,\bar B_N$ of diameter less than $\delta$ by
Lemma~\ref{app:bridge:lem-local}(a), together with $\bar B_\infty$ of
diameter less than $\delta$, such that $B_\infty\cup B_1\cup\dots\cup B_N$
covers $\mathcal K'$, and such that $\bar B_1$ is centered at
$(Z_*,L_1(Z_*),L_2(Z_*))$ with parameter $z$ (which is possible since
$\partial_s\Phi_1,\partial_t\Phi_2\ne0$ on $F_{LL}$). We may assume that each
$B_i$ meets $\mathcal K'$; then it meets exactly one component. Now take
$\kappa\in(0,1]$, so that Lemma~\ref{app:bridge:lem-left} applies, smaller
than all $\kappa_{B_i}$, than $\kappa_\infty$ and than
$\kappa_N$ of Lemma~\ref{app:bridge:lem-compact} for
$N=B_\infty\cup\bigcup B_i$, and so small that the following finitely many
conditions hold.
\begin{itemize}
\item[(1)] For each $i$, the two end points
$\psi_i(\pm1,\gamma_{\kappa}(\pm1))$ of the arc $\mathcal C_\kappa\cap\bar B_i$
lie in $B_\infty\cup\bigcup_jB_j\cup\{z<Z_*\}$. This holds for small
$\kappa$ because the end points $\psi_i(\pm1,\gamma_0(\pm1))$ of
$\mathcal K\cap\bar B_i$ lie in this open set, by the choice of the cover
and because near $z=Z_*$ the set $\mathcal K$ is $F_{LL}$. The same holds
for the end point at $c=u_0$ of the arc in $\bar B_\infty$, which is close
to the point $(\sigma_1,\infty,1/u_0)$ of $\mathcal D_1$.
\item[(2)] Whenever $\mathcal K\cap B_i\cap B_j\ne\emptyset$ (with $i$ or
$j$ possibly $\infty$), the arcs in $\bar B_i$ and $\bar B_j$ meet. Indeed,
choose $x=\psi_i(c,\gamma_0(c))$ in $\mathcal K\cap B_i\cap B_j$ with
$x\ne P_\infty$; then $\psi_i(c,\gamma_\kappa(c))\to x$, so for small $\kappa$
it lies in $\mathcal C_\kappa\cap B_j$, which is contained in the arc of
$\bar B_j$.
\item[(3)] The point of $\mathcal C_\kappa$ over $z=Z_*$, which by
Lemma~\ref{app:bridge:lem-left} converges to $(Z_*,L_1(Z_*),L_2(Z_*))$, lies
in $B_1$.
\end{itemize}
Apply Lemma~\ref{app:bridge:lem-glue} to $X=\mathcal C_\kappa$, which by
Lemma~\ref{app:bridge:lem-reduction} is homeomorphic to $\Xp(G_w)$, with the
following sets $\Gamma\supseteq O$: the arcs $\mathcal C_\kappa\cap\bar B_i$
and $\mathcal C_\kappa\cap B_i$, the arc $\mathcal C_\kappa\cap\bar B_\infty$
and $\mathcal C_\kappa\cap B_\infty$, and $\mathcal C_\kappa\cap\{z\le Z_*\}$
and $\mathcal C_\kappa\cap\{z<Z_*\}$; the last is connected, being a graph
over $(-\infty,Z_*]$ by Lemma~\ref{app:bridge:lem-left}. Group them into
$\Gamma^{(1)}$, consisting of the last two kinds and of the boxes meeting
$\mathcal K'_1$, and $\Gamma^{(2)}$, consisting of the boxes meeting
$\mathcal K'_2$ (if $\ell=2$). The covering condition of
Lemma~\ref{app:bridge:lem-glue} holds by (1), (3) and the interior points
of the boxes; connectivity of the adjacency graphs follows from (2), (3)
and the connectedness of $\mathcal K'_1$, $\mathcal K'_2$, which are covered
by the open sets $\mathcal K\cap B_i$. By Lemma~\ref{app:bridge:lem-compact}
every point of $\mathcal C_\kappa$ with $z\ge Z_*$ lies in some box, so
$\mathcal C_\kappa=\Gamma^{(1)}\cup\Gamma^{(2)}$, and the two are disjoint
by the choice of $\delta$ (note $\mathcal K'_2\subseteq\{z\ge\sigma_2\}$). All
points are regular, by Lemmas~\ref{app:bridge:lem-left}
and~\ref{app:bridge:lem-local}, and $\Xp(G_w)$ is a smooth curve. So the
components of $\Xp(G_w)$ correspond to $\Gamma^{(1)}$, which is not compact
since it contains points with $z\to-\infty$, and to $\Gamma^{(2)}$, which
is a finite union of compact arcs, hence a circle. By the list of
components of $\mathcal K'$, $\ell=2$ exactly in the interlaced case. This
proves (i) and (ii).

For (iii), the circle lies in the $\delta$-neighborhood of
$\mathcal K'_2$, and every point of $\mathcal K'_2$ is within $\delta$ of
the circle for small $\kappa$, since the arcs converge uniformly to the
pieces of $\mathcal K'_2$ in the boxes; as $\delta$ was arbitrary, this
is the Hausdorff convergence, and the $z$-projection of $\mathcal K'_2$ is
$[\sigma_2,\Sigma_1]$. For a compact $I\subset\mathcal J_1(\sigma_2)$, the
points of $\mathcal K$ with $s\in I$ and $t$ large lie in $\mathcal D_2$;
cover $\mathcal D_2\cap\{s\in I\}$ by finitely many of the boxes of
Lemma~\ref{app:bridge:lem-local}(c) and apply
Lemma~\ref{app:bridge:lem-compact}; the local graphs coincide on overlaps
and give \eqref{app:bridge:eq-Z}.
\end{proof}

The same computation in the nested case shows that the piece of
\eqref{app:bridge:eq-Z} lies on the non-compact component. Near
$P_\infty$, $z\to\sigma_1+\kappa\zeta(0,\kappa)$ along the non-compact
component while $V\to\infty$; by Lemma~\ref{lem:blowup}(b) this limit is
the lower Perron value of $G_w$, so
$\sigm(G_w)=\sigma_1-\kappa/(g_1S_2(\sigma_1))+O(\kappa^2)$.

\subsection{Unboundedly many circles}\label{app:bridge:subsec-many}

We now use parts which are bouquets of triangles. For $\beta>0$, a
\emph{triangle of type $\beta$} at a vertex $c$ consists of two further
vertices $x,y$ with $D_x=D_y=1$ and edges $cx$, $xy$, $yc$ with
$w_{cx}^2=w_{yc}^2=\frac34\beta$ and $w_{xy}^2=\frac14$ (the signs of the
weights are arbitrary). Let $\varpi(x)=1+x-\sqrt{1+x^2}$; it is strictly
increasing and strictly concave on $[0,\infty)$, with $0<\varpi(x)<\min(1,x)$
and $\varpi(x)\ge1-1/(2x)$ for $x>0$.

\begin{lemma}[Bouquets of triangles]\label{app:bridge:lem-triangles}
Let $X$ be the graph formed by at least two triangles, of types
$\beta_C$ ($C\in\mathcal T$), at a center $c$, sharing only $c$, with any $D_c$, and let
$B_X=\sum_C\beta_C$. Then $(X,c)$ is a part with $\rho=\frac12$,
$\beta(X)=\abs{\mathcal T}$, and for $z<\frac12$ and $s>0$
\[
\Phi_X(z,s)=\sum_{C\in\mathcal T}\varpi\bigl(\beta_Cb_0(z)s\bigr)-\epsilon_X(z)s,\qquad
b_0(z)=\frac{3/4}{(1-z)^2-1/4},\qquad
\epsilon_X(z)=\frac{3B_X}{1-2z}-D_c+z .
\]
In particular $\Phi_X(z,\cdot)$ is strictly concave, $S_X=-\epsilon_X$,
$\sigma_X$ is the unique zero of the strictly increasing function
$\epsilon_X$ on $(-\infty,\frac12)$, and $b_0$ is increasing on
$[0,\frac12)$ with $b_0(0)=1$.
\end{lemma}

\begin{proof}
$X-c$ is a disjoint union of edges $xy$, and $(D-\abs A)[\{x,y\}]$ has
smallest eigenvalue $\frac12$. The slice system decouples over the
triangles. For one triangle put $d=1-z>\frac12$, $\gamma=\frac34\beta$,
$\gamma'=\frac14$ and $m=d^2-\gamma'>0$. The slice solution is unique and
the triangle has the symmetry $x\leftrightarrow y$, so $V_x=V_y=:y_0$ and
$f_{cx}=f_{yc}$. Let $p=1+2f_{cx}$ and $q=1+2f_{xy}$. The equations are
$p^2-1=4\gamma sy_0$, $q^2-1=4\gamma'y_0^2$ and $dy_0=(p+q)/2$. Eliminating
$q=2dy_0-p$ gives $y_0=(dp-\gamma s)/m$ and then
$mp^2-4\gamma ds\,p+4\gamma^2s^2-m=0$, whose roots are
$p_\pm=(2\gamma ds\pm(m^2+4\gamma^2\gamma's^2)^{1/2})/m$. They are distinct
for all $s$, $p$ depends continuously on $s$, and $p\to1=p_+(0)$ as
$s\to0$ because $y_0$ stays bounded (Lemma~\ref{app:bridge:lem-slice});
so $p=p_+$. The contribution of the triangle to $\Phi_X$ is
$-2f_{cx}=-(p-1)=-(Es+\sqrt{1+b^2s^2}-1)$ with $E=2\gamma d/m$ and
$b=2\gamma\sqrt{\gamma'}/m=\beta b_0(z)$, and
$E+b=\frac34\beta(1+2d)/((d-\frac12)(d+\frac12))=3\beta/(1-2z)$. Summing,
$\Phi_X=(D_c-z)s-\sum_C(E_C+b_C)s+\sum_C\varpi(b_Cs)$, which is the formula.
The remaining statements follow from it and from
Lemma~\ref{app:bridge:lem-response}; $\beta(X)=3\abs{\mathcal T}-(2\abs{\mathcal T}+1)+1$.
\end{proof}

For two bouquets of triangles $(G_1,a_1)$, $(G_2,a_2)$ as parts, let
$\Phi_k^{\max}(z)=\max_{s>0}\Phi_k(z,s)$ (attained when $z>\sigma_k$, by
concavity and $\Phi_k\to-\infty$), and for $\nu<\Phi_2^{\max}(z)$ let
$U_2(z,\nu)$ be the larger root of $\Phi_2(z,\cdot)=\nu$.

\begin{lemma}[Separation of the descending sheet]\label{app:bridge:lem-sheet}
Let $(G_1,a_1)$, $(G_2,a_2)$ be bouquets of triangles with
$\sigma_1<\sigma_2<\Sigma_1<\Sigma_2$ and
$\Phi_1^{\max}(\sigma_2)<\Phi_2^{\max}(\Sigma_1)$. Let
$W=\{(z,s):z>\sigma_2,\ \Phi_1(z,s)>1\}$ and
$\Lambda_W(z,s)=\nu(\nu-1)/(s\,U_2(z,\nu))$ with $\nu=\Phi_1(z,s)$. Then:
\begin{enumerate}
\item[\textup{(a)}] $\partial_t\Phi_2\ne0$ on $\mathcal S$, so
$\mathcal S^-=\{p\in\mathcal S:\partial_t\Phi_2<0\}$ is open and closed in
$\mathcal S$, and $(z,s,t)\mapsto(z,s)$ is an analytic diffeomorphism from
$\mathcal S^-$ onto $W$ with inverse $t=U_2(z,\Phi_1(z,s))$, carrying
$\Lambda$ to $\Lambda_W$;
\item[\textup{(b)}] $\bar W\subseteq[\sigma_2,\Sigma_1]\times[L_1(\sigma_2),U_1(\sigma_2)]$,
and $\Lambda_W$ extends continuously to $\bar W$ by $0$ on $\partial W$;
\item[\textup{(c)}] for every $\kappa>0$, the subset of $\Xp(G_w)$ that
corresponds to $\Lambda_W^{-1}(\kappa)$ under
Lemma~\ref{app:bridge:lem-reduction} and \textup{(a)} is compact and open
and closed in $\Xp(G_w)$; so each of its connected components is a circle
of $\Xp(G_w)$. If $\kappa$ is a regular value of $\Lambda_W$, these circles
consist of regular points.
\end{enumerate}
\end{lemma}

\begin{proof}
(a) For $z\le\sigma_2$, $\Phi_2(z,\cdot)$ is strictly concave with
asymptotic slope $-\epsilon_2(z)\ge0$, since $\varpi'(x)\to0$; so
$\partial_t\Phi_2>0$. Points of $\mathcal S$ have $z<\Sigma_1$, since
$\mathcal J_1(z)=\emptyset$ for $z\ge\Sigma_1$. For $\sigma_2<z<\Sigma_1$,
$\nu=\Phi_1(z,s)\le\Phi_1^{\max}(z)\le\Phi_1^{\max}(\sigma_2)<\Phi_2^{\max}(\Sigma_1)\le\Phi_2^{\max}(z)$,
because $\partial_z\Phi_k<0$; so $t$ is one of the two roots of
$\Phi_2(z,\cdot)=\nu$, neither of which is the maximum point. The rest
follows from the implicit function theorem and from $\mathcal S^-\subseteq\{z>\sigma_2\}$.

(b) For $z>\sigma_2$ we have $\Phi_1(z,s)<\Phi_1(\sigma_2,s)$, so the
$s$-coordinates of $W$ lie in $\mathcal J_1(\sigma_2)=(L_1(\sigma_2),U_1(\sigma_2))$,
a bounded interval with $L_1(\sigma_2)>0$. A point of $\partial W$ has
$z=\sigma_2$ or $\Phi_1=1$ (a limit point with $z=\Sigma_1$ has
$\Phi_1\le\Phi_1^{\max}(\Sigma_1)=1$). Let $\tau_2(z)$ be the maximum point
of $\Phi_2(z,\cdot)$ for $\sigma_2<z\le\Sigma_1$; it is continuous and
positive, and $\tau_2(z)\to\infty$ as $z\downarrow\sigma_2$ because
$\partial_t\Phi_2(\sigma_2,\cdot)>0$ on every $[0,T]$. Since
$U_2\ge\tau_2$, $\Lambda_W\le\nu(\nu-1)/(L_1(\sigma_2)\inf\tau_2)\to0$ as
$\nu\to1$. Near $z=\sigma_2$: $\nu\le\nu_1:=\Phi_1^{\max}(\sigma_2)<\Phi_2^{\max}(\Sigma_1)<\beta_2$
(because $\varpi<1$ and $\epsilon_2>0$ for $z>\sigma_2$), and
$\Phi_2(\sigma_2,t)=\sum_C\varpi(\beta_Cb_0t)\uparrow\beta_2$; given $T$
choose $T'\ge T$ with $\Phi_2(\sigma_2,T')>\nu_1$; then
$\Phi_2(z,T')>\nu_1\ge\nu$ for $z$ close to $\sigma_2$, so $U_2(z,\nu)>T$,
and $\Lambda_W\le\nu_1^2/(L_1(\sigma_2)T)$.

(c) By Corollary~\ref{app:bridge:cor-level} and (a), the set in question is
open and closed in $\mathcal C_\kappa\cong\Xp(G_w)$, and by (b) it is
homeomorphic to the compact set $\{\Lambda_W=\kappa\}\subset W$. A compact
open and closed subset is a union of components. On $\mathcal S^-$,
$d(\Phi_1-\Phi_2)\ne0$, so by Corollary~\ref{app:bridge:cor-level} its
regular points are those where $d\Lambda_W\ne0$.
\end{proof}

\begin{definition}\label{app:bridge:def-Nr}
Let $\theta=100$ and $r\ge1$. The graph $\mathfrak N_r$ consists of an edge
$ab$ (the bridge, of weight $w$, $\kappa=w^2$), of $\theta^{2k}$ triangles
of type $\theta^{-3k}$ at $a$ for each $k=1,\dots,r$, and of $n$ triangles of
type $1$ at $b$, where $n\ge2.1\,\theta^{2r}$ is an integer; all triangles
share only their center. The potentials are $D=1$ on the vertices of the
triangles, $D_a=3B_a-\theta^{-2r}$ with $B_a=\sum_{k=1}^r\theta^{-k}$, and
$D_b=3n$.
\end{definition}

$\mathfrak N_r$ is a connected simple graph with $2+2\sum_k\theta^{2k}+2n$
vertices. Removing $a$ and $b$ leaves disjoint edges, and removing one
vertex leaves a component of cycle rank at least $2$; so
$\taups(\mathfrak N_r)=2$. Both sides are parts (bouquets of triangles), with
$\epsilon_a(z)=\theta^{-2r}+z+6B_az/(1-2z)$ and $\epsilon_b(z)=z+6nz/(1-2z)$
by Lemma~\ref{app:bridge:lem-triangles}; hence $\sigma_b=0$ and
$\sigma_a<0$. Let $\nu_*=\sum_k\theta^{2k}=\beta_a$ and
$S_\theta=(1-\theta^{-2})^{-1}$, so that $\nu_*<S_\theta\theta^{2r}$. All
numerical constants below are rational bounds for algebraic numbers and
were checked in exact rational arithmetic (and can be checked by hand).

\begin{proposition}\label{app:bridge:prop-Nr}
For every $r\ge1$ and every $\kappa\in[\frac1{100},\frac9{50}]$, the set
$\Xp(\mathfrak N_r)$ has at least $r$ compact connected components. For
almost every such $\kappa$ they are circles consisting of regular points.
\end{proposition}

\begin{proof}
We apply Lemma~\ref{app:bridge:lem-sheet} with $G_1$ the side of $a$ and
$G_2$ the side of $b$. Put $z_1=\frac1{50}$; then
$b_0(z_1)=\frac{625}{592}<1.0558$.

\emph{Step 1: $\Sigma_a<z_1$.} Since $\epsilon_a(z_1)\ge z_1$ and
$\varpi(x)\le\min(1,x)$, with $y=b_0(z_1)s$ and $e_1=z_1/b_0(z_1)=0.018944$,
\[
\Phi_a(z_1,s)\le g(y):=\sum_{k=1}^r\theta^{2k}\min(1,\theta^{-3k}y)-e_1y .
\]
If $y\le\theta^3$, then $g(y)\le y(\sum_k\theta^{-k}-e_1)\le y(\frac1{99}-e_1)<0$.
If $y\ge\theta^{3i}$ with $1\le i\le r$, bound the terms $k\le i$ by
$\theta^{2k}$ and the terms $k>i$ by $\theta^{-k}y$; since
$\sum_{k>i}\theta^{-k}\le\theta^{-i}/99<e_1$,
$g(y)\le S_\theta\theta^{2i}+\theta^{3i}(\theta^{-i}/99-e_1)\le\theta^{2i}(S_\theta+\frac1{99}-100e_1)<-0.88\,\theta^{2i}$.
So $\Phi_a(z_1,\cdot)<1$, hence $\Phi_a(z,\cdot)<1$ for $z\ge z_1$, and
$\Sigma_a<z_1$ by Lemma~\ref{app:bridge:lem-roots}.

\emph{Step 2: the descending root of $b$.} Let $0<z<z_1$ and $1<\nu\le\nu_*$.
Then $6nz\le\epsilon_b(z)\le z(1+\frac{25}4n)$ and $n-\nu\ge c_1n$ with
$c_1=1-S_\theta/2.1>0.5237$. Let $t_0=(n-\nu)/(2\epsilon_b)$. Using
$\varpi(x)\ge1-1/(2x)$ and $b_0\ge1$,
\[
\Phi_b(z,t_0)\ge\nu+\frac{n-\nu}2-\frac{n\epsilon_b}{n-\nu}>\nu
\quad\text{since}\quad(n-\nu)^2\ge c_1^2n^2>0.25n^2+0.04n\ge2n\epsilon_b,
\]
the middle inequality because $(c_1^2-\frac14)n\ge(c_1^2-\frac14)\cdot21000>0.04$.
So $\Phi_b(z,\cdot)=\nu$ has two roots $L_b<t_0<U_b$, and $\varpi<1$ gives
$\epsilon_bU_b<n-\nu$:
\begin{equation}\label{app:bridge:eq-Ub}
\frac{n-\nu}{2\epsilon_b(z)}<U_b(z,\nu)<\frac{n-\nu}{\epsilon_b(z)} .
\end{equation}
In particular $\Phi_b^{\max}(z)>\nu_*$ for $0<z<z_1$.

\emph{Step 3: hypotheses of Lemma~\ref{app:bridge:lem-sheet}.} By Step 2,
$\mathcal J_b(z)\ne\emptyset$ for $z<z_1$, so $\Sigma_b\ge z_1>\Sigma_a$. The
points $p_k$ of Step 5 lie in $W$ and have $z>0$, so $\Sigma_a>0=\sigma_b$.
Thus $\sigma_a<\sigma_b<\Sigma_a<\Sigma_b$. Moreover
$\Phi_a^{\max}(0)<\nu_*<\Phi_b^{\max}(\Sigma_a)$, since $\varpi<1$ and
$\epsilon_a(0)>0$. On $W$ we have $0<z<z_1$ and $1<\nu<\nu_*$, and
\eqref{app:bridge:eq-Ub} gives
\begin{equation}\label{app:bridge:eq-Kbounds}
\frac{6z\,\nu(\nu-1)}{s}<\Lambda_W(z,s)<\frac{2\nu(\nu-1)\epsilon_b}{s(n-\nu)}\le24\,\frac{z\nu^2}{s},
\end{equation}
because $2(\frac{25}4+\frac1n)/c_1<23.9$.

\emph{Step 4: separating lines.} For $1\le k\le r-1$ let
$s_k^\sharp=\theta^{3k+3/2}$. On $W\cap\{s=s_k^\sharp\}$, using
$\epsilon_a(z)\ge z$ and $b_0(z)\le b_0(z_1)$,
$\nu\le A_k-zs$ with
$A_k=\sum_{j\le k}\theta^{2j}+b_0(z_1)s_k^\sharp\sum_{j>k}\theta^{-j}\le(S_\theta+\frac{1000}{99}b_0(z_1))\theta^{2k}<11.6642\,\theta^{2k}$.
Since $\nu>0$, $z<A_k/s$ and $z\nu^2\le\max_{0\le z\le A_k/s}z(A_k-zs)^2=4A_k^3/(27s)$.
By \eqref{app:bridge:eq-Kbounds},
$\Lambda_W<24\cdot\frac4{27}\cdot11.6642^3\theta^{6k}/\theta^{6k+3}<0.0057<\frac1{100}$.

\emph{Step 5: peaks.} For $1\le k\le r$ let $s_k=0.6\,\theta^{3k}$ and let
$z_k$ be the solution of $\epsilon_a(z_k)=0.24\,\theta^{-k}$, which exists
since $\epsilon_a$ is continuous and increasing with
$\epsilon_a(0)=\theta^{-2r}<0.24\,\theta^{-k}$; since $\epsilon_a(z)>z$ for
$z>0$, $0<z_k<0.0024$. Keeping only the triangles of type $\theta^{-3k}$
and using $b_0\ge1$,
$\nu(p_k)\ge\theta^{2k}\varpi(0.6)-0.24\theta^{-k}s_k\ge0.2898\,\theta^{2k}\ge2898$,
since $\varpi(0.6)=1.6-\sqrt{1.36}>0.4338$. So $p_k=(z_k,s_k)\in W$. For
$z\le0.0024$, $\epsilon_a(z)\le\theta^{-2r}+1.0609\,z$, and
$\theta^{-2r}\le\theta^{-k}/100$; hence $z_k\ge0.23\,\theta^{-k}/1.0609\ge0.2167\,\theta^{-k}$,
and \eqref{app:bridge:eq-Kbounds} gives
\[
\Lambda_W(p_k)>\frac{6\cdot0.2167\cdot0.2898^2\,(1-\frac1{2898})}{0.6}>0.1819>\frac9{50}.
\]
Finally $s_{k-1}^\sharp<s_k<s_k^\sharp$, with $s_0^\sharp=0$, $s_r^\sharp=\infty$.

\emph{Step 6.} Let $\kappa\in[\frac1{100},\frac9{50}]$ and let $O_k$ be the
connected component of $\{\Lambda_W>\kappa\}\cap W$ containing $p_k$. By
Step 4, $O_k$ lies in the strip $\Pi_k=\{s_{k-1}^\sharp<s<s_k^\sharp\}$. It
is bounded and nonempty, so its boundary is nonempty; by
Lemma~\ref{app:bridge:lem-sheet}(b) a boundary point is not in $\partial W$
(where $\Lambda_W=0$), so it lies in $W$, and there $\Lambda_W=\kappa$ (if
$\Lambda_W>\kappa$ the point would belong to $O_k$). By Step 4 it is not on
the lines $s=s_j^\sharp$. Hence the compact set
$\Lambda_W^{-1}(\kappa)$ meets each of the $r$ strips $\Pi_k$, and none of
its connected components meets two strips, because it does not meet the
lines. By Lemma~\ref{app:bridge:lem-sheet}(c) this gives at least $r$
compact components of $\Xp(\mathfrak N_r)$. By Sard's
theorem~\cite{Milnor}, almost every $\kappa$ is a regular value of the
analytic function $\Lambda_W$, and then these components are circles of
regular points.
\end{proof}

Since $\sigma_a<\sigma_b<\Sigma_a<\Sigma_b$ and the four numbers are
distinct, Theorem~\ref{app:bridge:thm-gluing} shows that for the same
coefficients and sufficiently small $\kappa>0$ there is exactly one circle:
the number of circles changes with the strength of the bridge.

\begin{lemma}[Stability]\label{app:bridge:lem-stable}
Let $G$ be a loopless multigraph with real coefficients $c_0=(D,w)$, and let
$\mathcal C$ be a compact connected component of $\Xp(c_0)$ all of whose
points are regular. Then there are a bounded open set
$N\supseteq\mathcal C$ with $\bar N\subset\R\times(0,\infty)^V$ and a
neighborhood $\mathcal N$ of $c_0$ such that for every $c\in\mathcal N$ the
set $\Xp(c)\cap N$ is nonempty, open and closed in $\Xp(c)$, and a finite
union of circles of $\Xp(c)$ consisting of regular points. The set $N$ can
be chosen in any given neighborhood of $\mathcal C$. The same
conclusions hold for the zero set of any sufficiently small $C^1$
perturbation of $g_{c_0}$ on $\bar N$.
\end{lemma}

\begin{proof}
Near a regular point $\Xp(c_0)$ is an arc, so $\Xp(c_0)$ is locally
connected near $\mathcal C$, and the component $\mathcal C$ is open in
$\Xp(c_0)$; it is also compact. Since $\Xp(c_0)$ is closed in
$\R\times(0,\infty)^V$, the set $\Xp(c_0)\setminus\mathcal C$ has positive
distance from $\mathcal C$, and a small $\varepsilon$-neighborhood $N$ of
$\mathcal C$ satisfies $\bar N\cap\Xp(c_0)=\mathcal C$. Let
$\mu=\min_{\partial N}\abs{g_{c_0}}>0$. The map $g_c(z,V)$ depends
continuously on $c$, uniformly on $\bar N$, so for $c$ near $c_0$ it has no
zero on $\partial N$, and $\Xp(c)\cap N=\Xp(c)\cap\bar N$ is compact and
open and closed in $\Xp(c)$. If there were $c_m\to c_0$ and non-regular
zeros $x_m\in\bar N$ of $g_{c_m}$, a limit point $x$ would be a zero of
$g_{c_0}$ in $\bar N$, hence a point of $\mathcal C$, hence regular, and
regularity is an open condition in $(x,c)$; so for $c$ near $c_0$ all
points of $\Xp(c)\cap N$ are regular. At a point $x_0\in\mathcal C$,
restricting the residual to a transverse slice gives a nondegenerate
zero, which persists under small perturbations. Thus $\Xp(c)\cap N$ is
a nonempty compact one-dimensional
manifold, a finite union of circles, each of which is a component of
$\Xp(c)$. The boundary estimate and the rank argument apply to any
uniformly $C^1$-close residual, and the same transverse slice gives a
zero of that residual.
\end{proof}

\begin{theorem}[Unboundedly many circles]\label{app:bridge:thm-circles}
For every $r\ge1$ there are a connected simple graph with $\taups=2$,
consisting of an edge $ab$ with triangles attached at $a$ and at $b$, and
a nonempty open set of real coefficients on it, such that for every
coefficient in this set $\Xp$ has at least $r$ circles consisting of
regular points.
\end{theorem}

\begin{proof}
Take $\mathfrak N_r$ with the coefficients of
Definition~\ref{app:bridge:def-Nr} and a bridge weight with
$\kappa\in[\frac1{100},\frac9{50}]$ a regular value of $\Lambda_W$
(Proposition~\ref{app:bridge:prop-Nr}). It has at least $r$ circles of
regular points. Apply Lemma~\ref{app:bridge:lem-stable} to each of them,
with pairwise disjoint sets $N$, and intersect the neighborhoods
$\mathcal N$.
\end{proof}

The open set meets the generic set of Proposition~\ref{prop:generic-smooth},
so for generic coefficients the number of circles is unbounded on the
graphs with $\taups=2$. The graphs $\mathfrak N_r$ are large: the number of
vertices grows like $\theta^{2r}$, so that the number of circles obtained is
of the order of the logarithm of the number of vertices.

\subsection{From packing to a circle}\label{app:bridge:subsec-packing}

\begin{proof}[Proof of Theorem~\ref{real:thm:packing}(c)]
In each of the two vertex-disjoint subgraphs take a spanning tree and
two additional edges, then prune leaves. This gives connected graphs
$H_1,H_2$ of cycle rank two and minimum degree at least two. Let $a_k$ be
the end in $H_k$ of a shortest path joining the two graphs in $G$. Then
$H_k-a_k$ is a pseudoforest. Indeed, if $r_C$ is the number of edges
from $a_k$ to a component $C$ of $H_k-a_k$, then
\[
 2=\sum_C\bigl(\beta(C)+r_C-1\bigr).
\]
If one component had cycle rank two, every $r_C$ would be one, and any
other component would be a tree attached by a single edge, giving a
leaf of $H_k$. If there were no other component, $a_k$ itself would be a
leaf. Thus $(H_k,a_k)$ is a part.

First join $a_1,a_2$ by a single edge. By shifting the potentials of
the second part we can arrange $\sigma_1<\sigma_2<\Sigma_1<\Sigma_2$.
Theorem~\ref{app:bridge:thm-gluing} then gives a regular circle for a
sufficiently weak nonzero bridge. At any fixed energy $z_0$ strictly
between $\sigma_2$ and $\Sigma_1$, the two roots of each part give four
nondegenerate roots at bridge weight zero. They persist for sufficiently
small nonzero bridge weight by the implicit function theorem.

We next replace this edge by a path while preserving a regular circle
and the four nondegenerate roots at $z_0$. Subdivide an edge $ij$ by $x$,
with weights $a,b$ on $ix,xj$, and put $t=D_x-z>0$. Writing $p=f_{ix}$
and $q=f_{xj}$, the equations at $x$ and its two edges are
\[
 tV_x=1+p+q,\qquad
 p(1+p)=a^2V_iV_x,\qquad q(1+q)=b^2V_jV_x.
\]
With $f=pq/(1+p+q)$ these equations give
\[
 p=f+\frac{a^2V_i}{t},\qquad
 q=f+\frac{b^2V_j}{t},\qquad
 f(1+f)=\frac{a^2b^2}{t^2}V_iV_j.
\]
Conversely, these formulas and $V_x=(1+p+q)/t$ recover the original
positive solution. Eliminating $x$ changes the endpoint potentials
to $D_i-a^2/t,D_j-b^2/t$ and gives the edge weight $ab/t$. This
elimination preserves regularity and nondegeneracy at fixed $z$. With
the positive edge roots eliminated, the residual at $x$ is a convex
function of $V_x$ equal to $-1$ at zero, so its derivative at any
positive zero is at least $1/V_x>0$.

If the old edge has absolute weight $w>0$, take $D_x=M$ and
$a=b=\sqrt{Mw}$, and add $w$ to each endpoint potential. On compact
sets in $(z,V)$, the reduced equations on the old vertices converge in
$C^1$ to the old equations as $M\to\infty$. Apply the perturbation form of
Lemma~\ref{app:bridge:lem-stable} to the circle and the implicit function
theorem to the four roots at $z_0$. For sufficiently large $M$ the perturbed
circle lifts to a regular compact component of the subdivided system.
The elimination is a local analytic graph over the reduced system, so
the lifted circle is open in the full solution set and closed by
compactness. The four roots lift with nondegenerate Jacobians as well.
Repeating the subdivision gives the chosen path in $G$.

Add each unused vertex as an isolated vertex with $D_i$ above the compact
range of energies in the construction. Its equation $(D_i-z)V_i=1$
gives a unique positive analytic graph there with nonzero derivative
in $V_i$. Add all remaining edges with initially zero weights, then
turn them on with sufficiently small nonzero weights. The reduced
vertex residuals extend analytically at zero weights, so the circle
and the four roots persist. The same stability argument gives the
claimed open set of coefficients.
\end{proof}

%% file: appendices/bowtie.tex
\section{The two-dimensional model on the bowtie has 36 sheets}\label{app:bowtie}

\providecommand{\appbwSym}{\operatorname{Sym}}
\providecommand{\appbwadj}{\operatorname{adj}}
\providecommand{\appbwdisc}{\operatorname{disc}}
\providecommand{\appbwgr}{\operatorname{gr}}
\providecommand{\appbwT}{\mathcal T}
\providecommand{\appbwL}{\mathcal L}
\providecommand{\appbwJ}{\mathcal J}
\providecommand{\appbwN}{N^{\mathrm{nf}}}

In this appendix $G$ is the bowtie: the vertex $0$ and, for $t=1,2$, the
vertices $v_1^t,v_2^t$, with the six edges $0v_1^t$, $v_1^tv_2^t$, $v_2^t0$
oriented in this order. All vertex dimensions are equal to $p$, and we use
the cavity form \eqref{blc:eq:cavity} of the block system at $z=0$; by
Lemma~\ref{blc:lem:recovery}, $N_p(G)$ is the number of its admissible
solutions for generic $(D,A)$. As there are no parallel edges we write
$Q_{i\to j}$ for $Q^e_{i\to j}$, and we put
$m_{k\to i}=A_{ik}Q_{k\to i}^{-1}A_{ki}$ (the message from $k$ into $i$),
so that $Q_{i\to j}=R_i+m_{j\to i}$ and $R_i=D_i-\sum_{k\sim i}m_{k\to i}$.
We determine $N_2(G)$ by reducing the two triangles to a system at
their shared vertex and counting its branches at infinity.

\begin{theorem}[Computer-assisted]\label{app:bowtie:thm-36}
For the bowtie, $N_2(G)=36$. More precisely: for generic $(\widetilde D,B_1,B_2)$ the normal form
\eqref{app:bowtie:eq-nf} with $p=2$ has exactly $36$ admissible solutions,
all nondegenerate; its curve over the line $z\mapsto(\widetilde D-zI,B)$ (the
covering of the $z$-line, outside a finite set, formed by the admissible
solutions at $(\widetilde D-zI,B)$; see the paragraph ``Branches'' in
Subsection~\ref{app:bowtie:subsec-branches}) has
exactly $36$ points over $z=\infty$, all unramified, and they are of the nine
types listed in Proposition~\ref{app:bowtie:prop-branches} with the numbers
$1,1,1,1,2,2,4,4,20$.
\end{theorem}

Proposition~\ref{app:bowtie:prop-normal-form} identifies the generic count
with that of a system with one symmetric unknown $U_t$ for each triangle
and the covariance $V$ at the shared vertex. For $p=2$, eight types of
branches at infinity contribute $16$ sheets in total; type $(1,1)$
contributes $20$ by a plane intersection count. We then exclude branches
that escape these charts. The genericity checks in
Subsection~\ref{app:bowtie:subsec-comp} use exact rational computations
at explicit coupling points.

\subsection{Reduction to the normal form}\label{app:bowtie:subsec-nf}

For $Y\in\appbwSym_p(\C)$ let $\appbwgr(Y)=\{(Yu,u):u\in\C^p\}\subset\C^{2p}$.
A $p$-dimensional subspace of $\C^{2p}$ is a graph if and only if its
projection to the second factor is bijective. Let
$\appbwJ=\bigl(\begin{smallmatrix}0&I\\-I&0\end{smallmatrix}\bigr)$. For
$A\in\mathrm{GL}_p(\C)$ and $D'\in\appbwSym_p(\C)$ put
\[
T(D',A)=\begin{pmatrix}D'A^{-1}&-A^{\T}\\ A^{-1}&0\end{pmatrix},
\]
a symplectic matrix ($T^{\T}\appbwJ T=\appbwJ$). For a symplectic $K$ with
blocks $K_{ij}$ and symmetric $Y$ with $K_{21}Y+K_{22}$ invertible put
$k_K(Y)=(K_{11}Y+K_{12})(K_{21}Y+K_{22})^{-1}$.

\begin{lemma}\label{app:bowtie:lem-transfer}
Let $A$ be invertible and $D',Y$ symmetric with $Y$ invertible.
\textup{(a)} $T(D',A)\appbwgr(Y)=\appbwgr(D'-A^{\T}Y^{-1}A)$.
\textup{(b)} $T(D',A)\appbwgr(AY^{-1}A^{\T})=\appbwgr(D'-Y)$.
\textup{(c)} Let $K$ be symplectic. If $K\appbwgr(Y)$ is a graph for a
symmetric $Y$, then $K_{21}Y+K_{22}$ is invertible and
$K\appbwgr(Y)=\appbwgr(k_K(Y))$. If $K_{21}$ is invertible, then
$B=K_{21}^{-1}$, $s=K_{11}B$ and $s'=BK_{22}$ satisfy $s=s^{\T}$,
$s'=s'^{\T}$, and for symmetric $Y$ the space $K\appbwgr(Y)$ is a graph if
and only if $Y+s'$ is invertible, in which case
$k_K(Y)=s-B^{\T}(Y+s')^{-1}B$.
\end{lemma}

\begin{proof}
(a) $T(Yu,u)=(D'A^{-1}Yu-A^{\T}u,A^{-1}Yu)$; with $w=A^{-1}Yu$ this is
$((D'-A^{\T}Y^{-1}A)w,w)$. (b) $\appbwgr(AY^{-1}A^{\T})$ consists of the
vectors $(Aw,A^{-\T}Yw)$, which $T$ maps to $((D'-Y)w,w)$.
(c) $K(Yu,u)=((K_{11}Y+K_{12})u,(K_{21}Y+K_{22})u)$, which gives the first
statement. The relations $K^{\T}\appbwJ K=\appbwJ$ and
$K\appbwJ K^{\T}=\appbwJ$ say that $K_{11}^{\T}K_{21}$ and $K_{21}K_{22}^{\T}$
are symmetric and $K_{11}K_{22}^{\T}-K_{12}K_{21}^{\T}=I$. The first two give
the symmetry of $s$ and $s'$, the third gives $K_{12}=sB^{-1}s'-B^{\T}$, and
multiplying out shows
$K=\bigl(\begin{smallmatrix}I&s\\0&I\end{smallmatrix}\bigr)
\bigl(\begin{smallmatrix}0&-B^{\T}\\B^{-1}&0\end{smallmatrix}\bigr)
\bigl(\begin{smallmatrix}I&s'\\0&I\end{smallmatrix}\bigr)$. The outer
factors add $s$, $s'$ to the argument of $\appbwgr$, and the middle factor
maps $(Xu,u)$ to $(-B^{\T}u,B^{-1}Xu)$, which is a graph if and only if $X$
is invertible, and then equals $\appbwgr(-B^{\T}X^{-1}B)$.
\end{proof}

For $t=1,2$ the \emph{loop matrix} of the triangle $t$ is
\[
K_t=T(0,A_{v_2^t0})\,T(D_{v_2^t},A_{v_1^tv_2^t})\,T(D_{v_1^t},A_{0v_1^t}).
\]

\begin{lemma}\label{app:bowtie:lem-loop}
Let $Q$ satisfy the four equations \eqref{blc:eq:cavity} \textup{(}$z=0$\textup{)}
at the arcs leaving $v_1^t$ and $v_2^t$, with all cavity precisions
of the triangle $t$ invertible. Then
\[
K_t\appbwgr(Q_{0\to v_1^t})=\appbwgr(-m_{v_2^t\to0}),\qquad
K_t\appbwgr(m_{v_1^t\to0})=\appbwgr(-Q_{0\to v_2^t}).
\]
\end{lemma}

\begin{proof}
Write $v_1,v_2$ for $v_1^t,v_2^t$. The equations at the internal arcs are
$Q_{v_1\to v_2}=D_{v_1}-A_{v_10}Q_{0\to v_1}^{-1}A_{0v_1}$,
$Q_{v_2\to0}=D_{v_2}-A_{v_2v_1}Q_{v_1\to v_2}^{-1}A_{v_1v_2}$, and
$Q_{v_2\to v_1}=D_{v_2}-A_{v_20}Q_{0\to v_2}^{-1}A_{0v_2}$,
$Q_{v_1\to0}=D_{v_1}-A_{v_1v_2}Q_{v_2\to v_1}^{-1}A_{v_2v_1}$. Three
applications of Lemma~\ref{app:bowtie:lem-transfer}(a) map
$\appbwgr(Q_{0\to v_1})$ to $\appbwgr(Q_{v_1\to v_2})$,
$\appbwgr(Q_{v_2\to0})$ and $\appbwgr(-A_{0v_2}Q_{v_2\to0}^{-1}A_{v_20})$.
Three applications of (b) map
$\appbwgr(m_{v_1\to0})=\appbwgr(A_{0v_1}Q_{v_1\to0}^{-1}A_{0v_1}^{\T})$ to
$\appbwgr(D_{v_1}-Q_{v_1\to0})=\appbwgr(A_{v_1v_2}Q_{v_2\to v_1}^{-1}A_{v_1v_2}^{\T})$,
then to $\appbwgr(A_{v_20}Q_{0\to v_2}^{-1}A_{v_20}^{\T})$ and to
$\appbwgr(-Q_{0\to v_2})$.
\end{proof}

Let $\mathcal P^\circ$ be the set of coefficients $(D,A)$ of the bowtie for
which all $A_e$ and both $(K_t)_{21}$ are invertible, and for
$(D,A)\in\mathcal P^\circ$ let $B_t,s_t,s'_t$ be the matrices of
Lemma~\ref{app:bowtie:lem-transfer}(c) for $K_t$ and
$k_t(Y)=s_t-B_t^{\T}(Y+s'_t)^{-1}B_t$. Put
\[
\varphi(D,A)=(\widetilde D,B_1,B_2),\qquad \widetilde D=D_0+\textstyle\sum_t(s_t+s'_t).
\]
The \emph{normal form} with coefficients
$(\widetilde D,B_1,B_2)\in\appbwSym_p(\C)\times\mathrm{GL}_p(\C)^2$ is the
system
\begin{equation}\label{app:bowtie:eq-nf}
U_tVU_t-U_t-B_tVB_t^{\T}=0\quad(t=1,2),\qquad
\sum_{t=1,2}\bigl(U_t+B_t^{\T}U_t^{-1}B_t\bigr)-V^{-1}=\widetilde D
\end{equation}
in the unknowns $U_1,U_2,V\in\appbwSym_p(\C)$; it is
\eqref{app:cert:eq-rose} with $m=2$. A solution is \emph{admissible} if
$U_1,U_2,V$ are invertible, and \emph{nondegenerate} if the Jacobian matrix
of \eqref{app:bowtie:eq-nf} with respect to $(U,V)$ is invertible there.
Given an admissible solution at $\varphi(D,A)$, define for $t=1,2$
\begin{gather*}
Y^+_0=U_t-s'_t,\quad Y^+_1=D_{v_1^t}-A_{v_1^t0}(Y^+_0)^{-1}A_{0v_1^t},\quad
Y^+_2=D_{v_2^t}-A_{v_2^tv_1^t}(Y^+_1)^{-1}A_{v_1^tv_2^t},\\
Y^-_0=V^{-1}-s_t+B_t^{\T}U_t^{-1}B_t,\quad Y^-_1=D_{v_2^t}-A_{v_2^t0}(Y^-_0)^{-1}A_{0v_2^t},\quad
Y^-_2=D_{v_1^t}-A_{v_1^tv_2^t}(Y^-_1)^{-1}A_{v_2^tv_1^t}.
\end{gather*}
We say that the solution satisfies the \emph{chain condition} if these
twelve matrices are defined and invertible.

\begin{proposition}\label{app:bowtie:prop-normal-form}
Let $(D,A)\in\mathcal P^\circ$.
\begin{enumerate}
\item[\textup{(i)}] The map $Q\mapsto(U_t=Q_{0\to v_1^t}+s'_t,\ V=R_0^{-1})$
is a bijection from the admissible solutions of \eqref{blc:eq:cavity} with
$z=0$ onto the admissible solutions of \eqref{app:bowtie:eq-nf} at
$\varphi(D,A)$ that satisfy the chain condition. The inverse is
$(Q_{0\to v_1^t},Q_{v_1^t\to v_2^t},Q_{v_2^t\to0})=(Y_0^+,Y_1^+,Y_2^+)$ and
$(Q_{0\to v_2^t},Q_{v_2^t\to v_1^t},Q_{v_1^t\to0})=(Y_0^-,Y_1^-,Y_2^-)$.
Both maps are isomorphisms between the corresponding open subschemes of the
two solution schemes; in particular they preserve multiplicities and
nondegeneracy.
\item[\textup{(ii)}] $\mathcal P^\circ$ is a nonempty Zariski open set and
$\varphi(\mathcal P^\circ)=\appbwSym_p(\C)\times\mathrm{GL}_p(\C)^2$; in
particular $\varphi$ is dominant.
\item[\textup{(iii)}] $N_p(G)=\appbwN_p$, where $\appbwN_p$ is the generic
number of admissible solutions of \eqref{app:bowtie:eq-nf}
\textup{(}Lemma~\textup{\ref{app:bowtie:lem-generic})}.
\end{enumerate}
\end{proposition}

\begin{proof}
We use the identity, valid for invertible $U,R,U-R,B$ with $V=R^{-1}$,
\begin{equation}\label{app:bowtie:eq-stein}
UVU-U-BVB^{\T}=BV\,\Phi\,B^{-1}(U-R)VU,\qquad
\Phi=R-B^{\T}U^{-1}R(U-R)^{-1}B,
\end{equation}
which follows by expanding the right-hand side
($BVRB^{-1}(U-R)VU=UVU-U$ and the second term gives $BVB^{\T}$). Note that
$\Phi=R-k(U-s')+k(U-s'-R)$ for $k(Y)=s-B^{\T}(Y+s')^{-1}B$, by
$U^{-1}-(U-R)^{-1}=-U^{-1}R(U-R)^{-1}$.

(i) Let $Q$ be admissible; put $Q_t=Q_{0\to v_1^t}$, $S_t=m_{v_1^t\to0}$,
$R=R_0$. The equations at the arcs leaving $0$ read $Q_t=R+S_t$ and
$Q_{0\to v_2^t}=R+m_{v_2^t\to0}$. By Lemma~\ref{app:bowtie:lem-loop} the
spaces $K_t\appbwgr(Q_t)$ and $K_t\appbwgr(S_t)$ are graphs, so by
Lemma~\ref{app:bowtie:lem-transfer}(c) $U_t=Q_t+s'_t$ and
$S_t+s'_t=U_t-R$ are invertible, $k_t(Q_t)=-m_{v_2^t\to0}$ and
$k_t(S_t)=-Q_{0\to v_2^t}$. Hence
\begin{equation}\label{app:bowtie:eq-hub}
Q_t-S_t=R,\qquad k_t(Q_t)-k_t(S_t)=R\quad(t=1,2),\qquad
R=D_0-\textstyle\sum_t\bigl(S_t-k_t(Q_t)\bigr),
\end{equation}
the last being the definition of $R_0$. By \eqref{app:bowtie:eq-stein} the
second relation is the first equation of \eqref{app:bowtie:eq-nf} with
$V=R^{-1}$, and inserting $S_t=U_t-s'_t-R$ and
$k_t(Q_t)=s_t-B_t^{\T}U_t^{-1}B_t$ into the third gives the second equation.
Thus $(U,V)$ is an admissible solution; $Y^\pm_j$ are the cavity precisions
listed in the statement (for $Y_0^-$: $R-k_t(Q_t)=R+m_{v_2^t\to0}$), so the
chain condition holds, and $Q$ is recovered from $(U,V)$.

Conversely let $(U,V)$ be admissible with the chain condition, and put
$R=V^{-1}$, $Q_t=U_t-s'_t$, $S_t=Q_t-R$. The first equation of
\eqref{app:bowtie:eq-nf} gives $UV(U-R)=BVB^{\T}$, so $U_t-R$ is invertible
and \eqref{app:bowtie:eq-stein} gives $k_t(Q_t)-k_t(S_t)=R$. Define $Q$ by the
inverse formulas. The equations at the arcs leaving $v_1^t,v_2^t$ hold by
construction, so Lemma~\ref{app:bowtie:lem-loop} applies:
$m_{v_2^t\to0}=-k_t(Q_t)$, hence $Q_{0\to v_2^t}=Y^-_0=R+m_{v_2^t\to0}$, and
$K_t\appbwgr(m_{v_1^t\to0})=\appbwgr(-Y_0^-)=\appbwgr(k_t(S_t))=K_t\appbwgr(S_t)$,
hence $m_{v_1^t\to0}=S_t$. The second equation of \eqref{app:bowtie:eq-nf} is
the third relation of \eqref{app:bowtie:eq-hub}, so
$R_0=D_0-\sum_t(m_{v_1^t\to0}+m_{v_2^t\to0})=R$, and then
$Q_{0\to v_1^t}=R_0+m_{v_1^t\to0}$, $Q_{0\to v_2^t}=R_0+m_{v_2^t\to0}$ are the
two remaining equations. All cavity precisions and $R_0$ are invertible. For
an edge $ij$ with invertible cavity precisions,
$\det E_e=\det Q_{j\to i}\det R_i=\det Q_{i\to j}\det R_j$, because
$Q_{i\to j}-A_{ij}Q_{j\to i}^{-1}A_{ji}=R_i$; as $G$ is connected, all $R_i$
and $E_e$ are invertible and $Q$ is admissible.

For the scheme statement we eliminate variables in steps that do not change
the scheme. (a) The equations at the arcs leaving $v_1^t,v_2^t$ express the
four internal cavity precisions as regular functions of $Q_{0\to v_1^t}$
and $Q_{0\to v_2^t}$ (triangular form), so they can be eliminated. (b) The
map $\mu_t:Q_{0\to v_2^t}\mapsto m_{v_1^t\to0}$ obtained from the
counterclockwise chain is a composition of isomorphisms between open
subsets of $\appbwSym_p$ of the form $Y\mapsto D'-A^{\T}Y^{-1}A$, the last one
of the form $Y\mapsto A^{\T}Y^{-1}A$, and by
Lemma~\ref{app:bowtie:lem-loop} its inverse is $S\mapsto-k_t(S)$; we use
$S_t=\mu_t(Q_{0\to v_2^t})$ as variable. The four equations at the arcs
leaving $0$ become $Q_t-S_t-\rho=0$ and $k_t(Q_t)-k_t(S_t)-\rho=0$ with
$\rho=D_0-\sum_s(S_s-k_s(Q_s))$. (c) Replace the second equation by the
difference of the two, introduce $R$ with the equation $R=\rho$, eliminate
$S_t=Q_t-R$, and pass to $U_t=Q_t+s'_t$, $V=R^{-1}$. The equations become
$\Phi_t=0$ and the second equation of \eqref{app:bowtie:eq-nf}; by
\eqref{app:bowtie:eq-stein} the entries of $\Phi_t$ and of
$U_tVU_t-U_t-B_tVB_t^{\T}$ generate the same ideal, the two factors being
invertible matrices of regular functions.

(ii) $\mathcal P^\circ$ is open, and nonempty by the following computation.
Take $D_{v_1^t}=D_{v_2^t}=0$, $A_{0v_1^t}=A_{v_1^tv_2^t}=I$ and
$A_{v_2^t0}=-B_t$. Then $T(0,I)^2=-I$ and $K_t=-T(0,-B_t)
=\bigl(\begin{smallmatrix}0&-B_t^{\T}\\B_t^{-1}&0\end{smallmatrix}\bigr)$, so
$(K_t)_{21}=B_t^{-1}$, $s_t=s'_t=0$ and $\varphi(D,A)=(D_0,B_1,B_2)$ with
$D_0$ arbitrary. (This is the realization of
Lemma~\ref{app:cert:lem-windmill}.)

(iii) Let $W$ be the open set of Lemma~\ref{blc:lem:recovery}(ii) and
$\Omega$ that of Lemma~\ref{app:bowtie:lem-generic}. By (ii),
$\mathcal P^\circ\cap\varphi^{-1}(\Omega)\cap W$ is nonempty and open; for
$(D,A)$ in it the $N_p(G)$ admissible solutions inject by (i) into the
$\appbwN_p$ admissible solutions at $\varphi(D,A)$.

For the reverse inequality choose $(\widetilde D,B)\in\Omega$ and its
realization in (ii). Every admissible normal-form solution satisfies
the chain condition there. Indeed, writing $R=V^{-1}$, the first
normal-form equation gives
\[
 U_tV(U_t-R)=B_tVB_t^{\T}.
\]
Thus $U_t-R$ is invertible, and \eqref{app:bowtie:eq-stein} gives
$R+B_t^{\T}U_t^{-1}B_t=B_t^{\T}(U_t-R)^{-1}B_t$. In the realization from
(ii), $Y_0^+=U_t$ and $Y_0^-=R+B_t^{\T}U_t^{-1}B_t$ are therefore
invertible. The other chain matrices are obtained by inversion and
congruence by invertible edge blocks, since the two internal diagonals
are zero. By (i) all $\appbwN_p$ simple solutions lift to simple
admissible solutions of the full system. They persist at nearby
coefficients by the implicit function theorem, including on $W$.
This proves $N_p(G)\ge\appbwN_p$.
\end{proof}

\begin{lemma}\label{app:bowtie:lem-generic}
\textup{(i)} There are a nonempty Zariski open set
$\Omega\subseteq\appbwSym_p(\C)\times\mathrm{GL}_p(\C)^2$ and an integer
$\appbwN_p\ge0$ such that for $(\widetilde D,B)\in\Omega$ the system
\eqref{app:bowtie:eq-nf} has exactly $\appbwN_p$ admissible solutions, all
nondegenerate. \textup{(ii)} For $p=2$ one may choose $\Omega$ such that in
addition $\tr(VU_t)\ne1$ $(t=1,2)$ at every admissible solution.
\end{lemma}

\begin{proof}
Let $\sigma=p(p+1)/2$ and let $\mathcal I$ be the set of
$(U_1,U_2,V,B_1,B_2)$ with $U_t,V$ invertible symmetric, $B_t$ invertible,
satisfying the first equation of \eqref{app:bowtie:eq-nf}. For fixed
$(U,V)$ the derivative of $B\mapsto BVB^{\T}$ is onto $\appbwSym_p$
($\delta B=\tfrac12HB^{-\T}V^{-1}$ gives $H$), so $\mathcal I$ is smooth of
pure dimension $3\sigma+2p^2-2\sigma=\sigma+2p^2$, the dimension of the
coefficient space. Let $\pi:\mathcal I\to\appbwSym_p\times\mathrm{GL}_p^2$
send a point to $(\widetilde D,B)$ with $\widetilde D$ given by the second
equation; its fibers are the sets of admissible solutions. By generic
smoothness \cite[Corollary~III.10.7]{Hartshorne}, applied to each
component, there is a dense open set over which $\pi$ is smooth, hence
\'etale since the dimensions agree. A tangent vector $(\delta U,\delta V,\delta B)$
of $\mathcal I$ is killed by $d\pi$ only if $\delta B=0$ and it is killed by
the Jacobian of \eqref{app:bowtie:eq-nf} in $(U,V)$; so at points where
$\pi$ is \'etale the solutions are nondegenerate. The components of
$\mathcal I$ that do not dominate have images in a proper closed set, and
over a dense open set the dominating ones are finite
\cite[Exercise~II.3.7]{Hartshorne}. This gives (i). For (ii), the map
$\mathcal I\to\{(U,V)\}$ is surjective onto an open subset of
$\appbwSym_2^3\cong\C^9$ (complex nondegenerate symmetric matrices are
congruent) with fibers of dimension $2\dim\mathrm O_2=2$, so the preimage of
the hypersurface $\{\tr(VU_t)=1\}$ has dimension $10<11$ and its image under
$\pi$ is not dense.
\end{proof}

\subsection{The case \texorpdfstring{$p=2$}{p=2}: a five-variable system and one coupling}\label{app:bowtie:subsec-five}

From now on $p=2$. For a solution of \eqref{app:bowtie:eq-nf} put $W=V^{-1}$,
$C_t=VB_tVB_t^{\T}$, $k_t=\tr C_t$, $L_t(V)=B_tVB_t^{\T}+B_t^{\T}VB_t$,
$\beta_t=\det B_t$, $Z_t=VU_t$, $s_t=\tr Z_t-1$, $\delta_t=\det Z_t$.

\begin{lemma}\label{app:bowtie:lem-ch}
Let $Z\in M_2(\C)$, $C=Z^2-Z$, $s=\tr Z-1$, $\delta=\det Z$, $k=\tr C$.
\textup{(i)} $C=sZ-\delta I$, $k=s^2+s-2\delta$, $\det C=\delta(\delta-s)$,
and $s^4-(2k+1)s^2+k^2-4\det C=0$. \textup{(ii)} If $s\ne0$ then
$Z=(C+\delta I)/s$. \textup{(iii)} Conversely, if $C\in M_2(\C)$ and $s\ne0$
satisfy $s^4-(2k+1)s^2+k^2-4\det C=0$ with $k=\tr C$, then
$\delta=(s^2+s-k)/2$ and $Z=(C+\delta I)/s$ satisfy $Z^2-Z=C$,
$\tr Z=s+1$, $\det Z=\delta$.
\end{lemma}

\begin{proof}
(i) Cayley--Hamilton gives $Z^2=(s+1)Z-\delta I$, i.e. $C=sZ-\delta I$; the
eigenvalues of $C$ are $s\zeta_i-\delta$, whose sum and product are
$s(s+1)-2\delta$ and $\delta(\delta-s)$; substituting
$\delta=(s^2+s-k)/2$ gives $4\det C=(s^2-k)^2-s^2$. (ii) is (i). (iii) $Z$
is a polynomial in $C$ and $s^2(Z^2-Z-C)=C^2-kC+\delta(\delta-s)I$, which
vanishes by Cayley--Hamilton for $C$ since the quartic equation says
$\delta(\delta-s)=\det C$; $\tr Z=(k+2\delta)/s=s+1$ and
$\det Z=(\det C+\delta k+\delta^2)/s^2=\delta$.
\end{proof}

\begin{proposition}\label{app:bowtie:prop-five}
Put $\Delta_t(V)=k_t^2-4(\det V)^2\beta_t^2$ $(=\appbwdisc C_t)$ and
$N_t(V)=VL_t(V)-k_tI$. The map $(U,V)\mapsto(V,s_1,s_2)$ is a bijection from
the admissible solutions of \eqref{app:bowtie:eq-nf} with $s_1s_2\ne0$ onto
the solutions of
\begin{equation}\label{app:bowtie:eq-five}
E_t:\ s_t^4-(1+2k_t)s_t^2+\Delta_t(V)=0\quad(t=1,2),\qquad
F:\ \sum_t\frac{L_t(V)-k_tW}{s_t}+(s_1+s_2-1)W=\widetilde D
\end{equation}
with $\det V\ne0$, $s_t\ne0$ and $\delta_t:=(s_t^2+s_t-k_t)/2\ne0$; the
inverse is $U_t=W(C_t+\delta_tI)/s_t$. Moreover $N_t(V)$ is traceless, and
$VF$ reads $\sum_tN_t(V)/s_t+(s_1+s_2-1)I=V\widetilde D$.
\end{proposition}

\begin{proof}
Multiplying the first equation of \eqref{app:bowtie:eq-nf} by $V$ gives
$Z_t^2-Z_t=C_t$, and from $U_t(VU_t-I)=B_tVB_t^{\T}$ we get
$B_t^{\T}U_t^{-1}B_t=W\,\Xi_t(Z_t-I)\Xi_t^{-1}$ with $\Xi_t=VB_t^{\T}$. By
Lemma~\ref{app:bowtie:lem-ch}, $Z_t=(C_t+\delta_tI)/s_t$, and
$\Xi_tC_t\Xi_t^{-1}=VB_t^{\T}VB_t$, so $W(C_t+\Xi_tC_t\Xi_t^{-1})=L_t(V)$.
With $U_t=WZ_t$ the second equation of \eqref{app:bowtie:eq-nf} becomes
$\sum_tW(C_t+\Xi_tC_t\Xi_t^{-1}+2\delta_tI)/s_t-3W=\widetilde D$, i.e.
$\sum_t(L_t+2\delta_tW)/s_t-3W=\widetilde D$, which is $F$ because
$2\delta_t/s_t=s_t+1-k_t/s_t$; $E_t$ is Lemma~\ref{app:bowtie:lem-ch}(i) with
$\det C_t=(\det V)^2\beta_t^2$, and $\delta_t=\det Z_t\ne0$ is the
invertibility of $U_t$. Conversely, given a solution of
\eqref{app:bowtie:eq-five}, Lemma~\ref{app:bowtie:lem-ch}(iii) gives $Z_t$
and $U_t=WZ_t$ satisfies $U_tVU_t-U_t-B_tVB_t^{\T}=W(Z_t^2-Z_t-C_t)=0$;
$U_t$ is symmetric because $WC_t=B_tVB_t^{\T}$ and
$WC_t^{j+1}=(WC_t)V(WC_t^j)$ are symmetric, and the computation above runs
backwards. Finally $\tr(VL_t(V))=2k_t$.
\end{proof}

\paragraph{One coupling.} For $B\in\mathrm{GL}_2(\C)$ write
$B=\Sigma+aJ$ with $\Sigma$ symmetric and
$J=\bigl(\begin{smallmatrix}0&1\\-1&0\end{smallmatrix}\bigr)$, and put
$L_B(V)=BVB^{\T}+B^{\T}VB$, $\kappa_B(V)=\tr(VBVB^{\T})$,
$\Delta_B(V)=\kappa_B(V)^2-4(\det V)^2(\det B)^2$,
$N_B(V)=VL_B(V)-\kappa_B(V)I$. Let $\appbwT$ be the plane of traceless
symmetric $2\times2$ matrices, $q=q_B=\appbwdisc\Sigma$, $Y=Y_B=2\Sigma-(\tr\Sigma)I$,
$X=X_B=YJ$, $\lambda_X=2(\det\Sigma-a^2)$, $\lambda_Y=\lambda_X+q$ and
$\gamma=(\tr\Sigma)^2-\lambda_X$. Note $\tr B=\tr\Sigma$, $\det B=\det\Sigma+a^2$
and $\appbwdisc B=q-4a^2$.

\begin{lemma}\label{app:bowtie:lem-coupling}
For $V\in\appbwSym_2(\C)$ let $\theta_B(V)=\tr(\Sigma V)$ and
$M_B(V)=2V\Sigma-\theta_B(V)I$. Then
\begin{gather*}
L_B(V)=2\Sigma V\Sigma+2a^2\appbwadj V,\qquad
\kappa_B(V)=\theta_B(V)^2+2(a^2-\det\Sigma)\det V,\\
N_B(V)=\theta_B(V)\,M_B(V),\qquad
\Delta_B(V)=\bigl(\theta_B(V)^2-4\det\Sigma\det V\bigr)\bigl(\theta_B(V)^2+4a^2\det V\bigr).
\end{gather*}
If $q\ne0$, then $X,Y$ are traceless symmetric, $X^2=Y^2=qI$, $XY=-YX$,
$\tr(XY)=0$, $\{I,X,Y\}$ is a basis of $\appbwSym_2$, and
$L_B(X)=\lambda_XX$, $L_B(Y)=(\tr\Sigma)qI+\lambda_YY$,
$L_B(I)=\gamma I+(\tr\Sigma)Y$, $\Delta_B(X)=-16a^2q^2\det\Sigma$,
$\Delta_B(Y)=(\tr\Sigma)^2(q-4a^2)q^2$.
\end{lemma}

\begin{proof}
The cross terms of $(\Sigma+aJ)V(\Sigma-aJ)+(\Sigma-aJ)V(\Sigma+aJ)$ cancel,
and $JVJ^{\T}=\appbwadj V$. Then
$\kappa_B=\tfrac12\tr(VL_B(V))=\tr((V\Sigma)^2)+2a^2\det V$, and with
$P=V\Sigma$, Cayley--Hamilton gives $\tr P^2=(\tr P)^2-2\det P$ and
$N_B=2P^2-(\tr P^2)I=\tr P\,(2P-\tr P\,I)$. The factorization of $\Delta_B$
follows from $\Delta_B=(\kappa_B-2\det B\det V)(\kappa_B+2\det B\det V)$. The
last assertions follow from $\Sigma=\frac12(Y+\tr\Sigma\,I)$, $Y^2=qI$,
$JY=-YJ$ and $\appbwadj V=(\tr V)I-V$.
All identities of this lemma have also been verified by computer as
polynomial identities in the entries of $B$ and $V$.
\end{proof}

\begin{corollary}\label{app:bowtie:cor-coupling}
Assume $a\,q\det\Sigma\tr\Sigma\ne0$. \textup{(a)} $N_B(V)=0$ if and only
if $\theta_B(V)=0$ or $V\in\C\appbwadj\Sigma$. \textup{(b)} On
$\C\appbwadj\Sigma$ the matrix $VBVB^{\T}$ is scalar, $\Delta_B$ vanishes to
second order there, and $\theta_B(\appbwadj\Sigma)=2\det\Sigma$.
\textup{(c)} On $\{\theta_B=0\}$, $L_B(V)=\lambda_X(V-(\tr V)I)$,
$\kappa_B=-\lambda_X\det V$ and $\Delta_B=-16a^2\det\Sigma(\det V)^2$; and
$\{\theta_B=0,\ \det V=0\}=\C w_+w_+^{\T}\cup\C w_-w_-^{\T}$, where
$w_\pm\ne0$ are the two isotropic vectors of $\Sigma$
\textup{(}$w^{\T}\Sigma w=0$, $w\ne0$\textup{)}; on these lines $VBVB^{\T}=0$. In
particular $\Delta_B\in(\theta_B,\det V)^2$ and
$\{N_B=0\}\cap\{\Delta_B=0\}$ is the union of the three lines
$\C\appbwadj\Sigma$, $\C w_\pm w_\pm^{\T}$. \textup{(d)} For
$V\in\appbwT\setminus\{0\}$: $N_B(V)=0$ iff $V\in\C^*X$, and
$N_B(V)\in\C^*V$ iff $V\in\C^*Y$; if also $\appbwdisc B\ne0$, then
$\Delta_B(V)\ne0$ in both cases. \textup{(e)} If $V=\sigma ww^{\T}$ has rank
one, then $N_B(V)w=\kappa_B(V)w$ and $\Delta_B(V)=\kappa_B(V)^2$.
\end{corollary}

\begin{proof}
(a) $M_B(V)=0$ iff $V\Sigma$ is scalar. (b) $\appbwadj\Sigma\,\Sigma=\det\Sigma\,I$;
with $K=\Sigma^{-1}J$ one has $K^2=-I/\det\Sigma$ (since
$J\Sigma^{-1}J=-\appbwadj(\Sigma^{-1})$), so for $V=c\,\appbwadj\Sigma$,
$VBVB^{\T}=c^2(\det\Sigma)^2(I+aK)(I-aK)$ is scalar; $\Delta_B=\appbwdisc(VBVB^{\T})$
is a quadratic form in the traceless part of $VBVB^{\T}$, which vanishes on
the line. (c) On $\theta_B=0$, $\tr((V\Sigma)^2)=-2\det V\det\Sigma$, and the
formulas follow from the lemma and the identity
$\Sigma V\Sigma=\tr(\Sigma V)\Sigma-\det\Sigma\,\appbwadj V$ (both sides are
linear in $V$; check on a basis). If $\det V=0$ then $V=\sigma ww^{\T}$,
$\theta_B(V)=\sigma w^{\T}\Sigma w$ and $VBVB^{\T}=\sigma^2(w^{\T}Bw)\,ww^{\T}B^{\T}$
with $w^{\T}Bw=w^{\T}\Sigma w$. (d) For traceless $V$,
$\theta_B(V)=\frac12\tr(YV)$, so $\theta_B(V)=0$ iff $V\in\C X$, and
$\C\appbwadj\Sigma\cap\appbwT=0$ since $\tr\appbwadj\Sigma=\tr\Sigma\ne0$.
If $N_B(V)=\lambda V$ with $\lambda\ne0$, then $\theta=\theta_B(V)\ne0$ and
$V(\Sigma-cI)=\frac\theta2I$ with $c=\lambda/(2\theta)$, so
$V\in\C\appbwadj(\Sigma-cI)=\C\bigl((\tr\Sigma-c)I-\Sigma\bigr)$, and
$\tr V=0$ forces $c=\frac12\tr\Sigma$, i.e. $V\in\C Y$. Conversely
$\theta_B(Y)=q$ and $N_B(Y)=q\,M_B(Y)=q(\tr\Sigma)Y$. The values of
$\Delta_B(X)$, $\Delta_B(Y)$ are in the lemma. (e) $V\Sigma w=\sigma(w^{\T}\Sigma w)w$,
so $N_B(V)w=\theta(2\theta-\theta)w$ and $\kappa_B=\theta^2$ as $\det V=0$.
\end{proof}

\subsection{Branches over \texorpdfstring{$z=\infty$}{z at infinity}}\label{app:bowtie:subsec-branches}

Fix $B=(B_1,B_2)$, write $\Sigma_t,a_t,q_t,X_t,Y_t,\theta_t,M_t,\kappa_t,\Delta_t,N_t,L_t$
for the quantities of $B_t$, and consider the line
$z\mapsto(\widetilde D-zI,B)$. Put $\tau=1/z$ and $V'=-zV$. For $\tau\ne0$,
a solution $(V,s)$ of \eqref{app:bowtie:eq-five} at $(\widetilde D-zI,B)$
corresponds to a solution $(V',s)$ of
\begin{equation}\label{app:bowtie:eq-scaled}
s_t^4-(1+2\tau^2\kappa_t)s_t^2+\tau^4\Delta_t=0,\qquad
\tau^2\sum_t\frac{L_t-\kappa_tW'}{s_t}+(s_1+s_2-1)W'=I-\tau\widetilde D,
\end{equation}
where $\kappa_t,\Delta_t,L_t$ are evaluated at $V'$ and $W'=V'^{-1}$; indeed
$L_t(V)=-\tau L_t(V')$, $k_t=\tau^2\kappa_t(V')$, $\Delta_t(V)=\tau^4\Delta_t(V')$,
and $-\tau F$ is the second equation. Multiplying by $V'$,
\begin{equation}\label{app:bowtie:eq-scaledV}
\tau^2\sum_tN_t(V')/s_t+(s_1+s_2-1)I=V'(I-\tau\widetilde D).
\end{equation}
The admissibility conditions are $\det V'\ne0$, $s_t\ne0$ and
$s_t^2+s_t-\tau^2\kappa_t\ne0$.

\paragraph{Branches.} Let $(\widetilde D,B)\in\Omega$, with $\Omega$ as in
Lemma~\ref{app:bowtie:lem-generic}. Then $(\widetilde D-zI,B)\in\Omega$ for
all $z$ outside a finite set $F$. Over $\C\setminus F$ the admissible
solutions form a finite unramified covering of degree $\appbwN_2$, and
$s_t\ne0$ on it. Choose $\rho$ with $F\subset\{\abs z<\rho\}$. The connected
components of the part over $\{\abs z>\rho\}$ are the \emph{branches at
infinity}; a branch $\beta$ is a cyclic covering of the punctured disc of
some degree $e_\beta$, parametrized by $u\mapsto z=u^{-e_\beta}$, on which
$V'$ and $s_t$ are holomorphic in $u$ and, being algebraic functions of $z$,
meromorphic at $u=0$ (Puiseux). Hence
\begin{equation}\label{app:bowtie:eq-sum}
\appbwN_2=\textstyle\sum_\beta e_\beta .
\end{equation}
On a branch we use the order $\nu$ in $\tau$ ($\nu(\tau)=1$, rational values),
and write $s_t\to s_t^0$, $V'\to V'_0$ for limits as $\tau\to0$. Two
elementary facts about the roots of $E_t$ are used throughout. If
$\tau^2C^0_t\to0$, where $C^0_t=V'B_tV'B_t^{\T}$, then the roots
$\sigma_\pm$ of $\sigma^2-(1+2\tau^2\kappa_t)\sigma+\tau^4\Delta_t=0$ satisfy
$\sigma_+\to1$ and $\sigma_-=\tau^4\Delta_t/\sigma_+$; so $s_t^0\in\{1,-1,0\}$,
and if $s_t^0=0$ then $s_t^2=\tau^4\Delta_t(V')/\sigma_+$ with
$\Delta_t(V')\not\equiv0$. Second, in the unscaled variable, if
$V=\tau^{-b}(V_*+o(1))$ with $b\ge0$, then $s_t=\tau^{-b}\varsigma$ where
$\varsigma^4-(\tau^{2b}+2\kappa_t(V_*+o(1)))\varsigma^2+\Delta_t(V_*+o(1))=0$;
thus $s_t=O(\tau^{-b})$, and if $\Delta_t(V_*)\ne0$ then $\tau^bs_t$ is
bounded away from $0$.

\paragraph{Charts and limit systems.} A \emph{type} is
$r=(r_1,r_2)\in\{0,1,2\}^2$; loops with $r_t=1$ are of \emph{rank one}. In
the chart $\mathcal R^r$ the unknowns are $(V',y_1,y_2)$ with $s_t=\tau^2y_t$
if $r_t=1$ and $s_t=y_t$ otherwise. Dividing $E_t$ by $\tau^4$ when $r_t=1$,
the chart equations are
\begin{gather*}
\tau^4y_t^4-(1+2\tau^2\kappa_t)y_t^2+\Delta_t=0\ (r_t=1),\qquad
y_t^4-(1+2\tau^2\kappa_t)y_t^2+\tau^4\Delta_t=0\ (r_t\ne1),\\
\sum_{r_t=1}\frac{L_t-\kappa_tW'}{y_t}+\tau^2\sum_{r_t\ne1}\frac{L_t-\kappa_tW'}{y_t}
+\Bigl(\tau^2\sum_{r_t=1}y_t+\sum_{r_t\ne1}y_t-1\Bigr)W'=I-\tau\widetilde D .
\end{gather*}
Put $\lambda_r=1-\sum_{r_t\ne1}(1-r_t)$. At $\tau=0$ we take $y_t=1-r_t$ for
$r_t\ne1$ (simple roots of $y^4-y^2$), and obtain the \emph{limit system}
$\appbwL^r(B)$ in the unknowns $V'$ and $y_t$ ($r_t=1$):
\begin{equation}\label{app:bowtie:eq-limit}
y_t^2=\Delta_t(V')\ (r_t=1),\qquad
\sum_{r_t=1}\frac{L_t(V')-\kappa_t(V')W'}{y_t}-\lambda_rW'=I .
\end{equation}
It does not involve $\widetilde D$. A solution is \emph{admissible} if
$\det V'\ne0$ and $y_t\ne0$ for $r_t=1$, and \emph{nondegenerate} if the
Jacobian of \eqref{app:bowtie:eq-limit} is invertible; since the equations
for $y_t$ ($r_t\ne1$) have derivative $\pm2$ and the other equations do not
involve these $y_t$ at $\tau=0$, this is the invertibility of the Jacobian of
the chart equations at $\tau=0$. The values are
$\lambda_{(0,0)}=-1$, $\lambda_{(0,2)}=\lambda_{(2,0)}=1$, $\lambda_{(2,2)}=3$,
$\lambda_{(0,1)}=\lambda_{(1,0)}=0$, $\lambda_{(1,2)}=\lambda_{(2,1)}=2$,
$\lambda_{(1,1)}=1$.

A branch is \emph{regular} if $V'\to V'_0$ with $V'_0$ invertible and
$\Delta_t(V'_0)\ne0$ for every $t$ with $s_t^0=0$.

\begin{lemma}\label{app:bowtie:lem-regular}
Let $(\widetilde D,B)\in\Omega$.
\textup{(a)} A regular branch with $r_t=1-s_t^0$ for $s_t^0=\pm1$ and
$r_t=1$ for $s_t^0=0$ converges in the chart $\mathcal R^r$ to an admissible
solution $x_0$ of $\appbwL^r(B)$.
\textup{(b)} If $x_0$ is a nondegenerate admissible solution of
$\appbwL^r(B)$, then exactly one branch converges to $x_0$ in $\mathcal R^r$,
and it has $e_\beta=1$.
\textup{(c)} Hence, if every branch is regular and every admissible solution
of every limit system is nondegenerate, $\appbwN_2$ is the total number of
admissible solutions of the nine limit systems.
\end{lemma}

\begin{proof}
(a) As $V'$ is bounded, $\tau^2C^0_t\to0$, so $s_t^0\in\{1,-1,0\}$. If
$s_t^0=0$ then $y_t=s_t/\tau^2$ satisfies $y_t^2=\Delta_t(V')/\sigma_+\to\Delta_t(V'_0)\ne0$.
The chart equations hold for small $\tau\ne0$ and pass to the limit.
(b) By the implicit function theorem the chart equations have a unique
analytic solution $x(\tau)$ near $(x_0,0)$. For small $\tau\ne0$ it is
admissible: $\det V'\ne0$ and $s_t\ne0$ are clear; for $r_t=1$,
$s_t^2+s_t-\tau^2\kappa_t=\tau^2(\tau^2y_t^2+y_t-\kappa_t)$ and
$(y_t-\kappa_t)(y_t+\kappa_t)\to-4(\det V'_0)^2\beta_t^2\ne0$; for $r_t=0$ the
expression tends to $2$; for $r_t=2$, $\delta_t(\delta_t-s_t)=\det C_t\ne0$
by Lemma~\ref{app:bowtie:lem-ch}. Hence for $\abs z$ large $x(1/z)$ is one of
the admissible solutions, and $x(\tau)$ traces a branch with $e_\beta=1$.
Points of any branch converging to $x_0$ lie, for small $\tau$, in the
neighborhood where $x(\tau)$ is the only solution. (c) follows from (a),
(b) and \eqref{app:bowtie:eq-sum}.
\end{proof}

Writing $M=\sum_{r_t=1}L_t(V')/y_t$ and taking the trace of $V'$ times
\eqref{app:bowtie:eq-limit} (using $\tr(V'L_t(V'))=2\kappa_t$) shows that
\eqref{app:bowtie:eq-limit} is equivalent to
\begin{equation}\label{app:bowtie:eq-limit2}
\tr V'=-2\lambda_r\quad\text{and}\quad M-I\in\C\,\appbwadj V'.
\end{equation}

\begin{proposition}\label{app:bowtie:prop-branches}
There is a nonempty Zariski open set $\mathcal U\subseteq\mathrm{GL}_2(\C)^2$
such that for $B\in\mathcal U$ and every $\widetilde D$ with
$(\widetilde D,B)\in\Omega$:
\textup{(i)} every branch at infinity is regular;
\textup{(ii)} the limit systems of the types
$(0,0),(0,2),(2,0),(2,2),(0,1),(1,0),(1,2),(2,1),(1,1)$ have exactly
$1,1,1,1,2,2,4,4,20$ admissible solutions, all nondegenerate;
\textup{(iii)} the map from branches to pairs (type, limit) is a bijection
onto the $36$ admissible limit solutions, and every branch is unramified.
Consequently $\appbwN_2=36$.
\end{proposition}

\begin{proof}
(ii) is Propositions~\ref{app:bowtie:prop-single} and~\ref{app:bowtie:prop-11},
(i) is Proposition~\ref{app:bowtie:prop-noirr}, and (iii) and the last
statement follow from Lemma~\ref{app:bowtie:lem-regular}. The set
$\mathcal U$ is the intersection of the finitely many nonempty open sets
introduced in those propositions.
\end{proof}

\begin{proof}[Proof of Theorem~\ref{app:bowtie:thm-36}]
By Proposition~\ref{app:bowtie:prop-branches}, $\appbwN_2=36$ (take
$(\widetilde D,B)\in\Omega\cap(\appbwSym_2\times\mathcal U)$, which is
nonempty and open), and by
Proposition~\ref{app:bowtie:prop-normal-form}(iii), $N_2(G)=36$.
\end{proof}

\begin{remark}\label{app:bowtie:rem-bij}
Since $N_2(G)=\appbwN_2$, for generic coefficients of the bowtie the
injection of Proposition~\ref{app:bowtie:prop-normal-form}(i) is a
bijection: every admissible solution of the normal form at $\varphi(D,A)$
satisfies the chain condition. The $z$-line used above is that of the normal
form; the $z$-line of the bowtie, which shifts all five vertex diagonals,
is mapped by $\varphi$ to a curve in the coefficient space of
\eqref{app:bowtie:eq-nf} that is in general not a line.
\end{remark}

\subsection{The limit systems}\label{app:bowtie:subsec-limit}

\begin{proposition}\label{app:bowtie:prop-single}
Assume that $a_tq_t\det\Sigma_t\tr\Sigma_t(q_t-4a_t^2)\ne0$ for $t=1,2$.
\begin{enumerate}
\item[\textup{(a)}] For $r\in\{(0,0),(0,2),(2,0),(2,2)\}$ the only solution
of $\appbwL^r(B)$ is $V'=-\lambda_rI$, and it is nondegenerate.
\item[\textup{(b)}] For $r=(0,1)$ or $(1,0)$, with rank-one loop $t$, the
admissible solutions are exactly the two points
$V'=\pm Y_t/\sqrt{q_t-4a_t^2}$, $y_t=\pm(\tr\Sigma_t)q_t/\sqrt{q_t-4a_t^2}$
\textup{(}same sign\textup{)}.
\item[\textup{(c)}] For $r=(1,2)$ or $(2,1)$, with rank-one loop $t$, every
admissible solution has $V'=-2I+\eta Y_t$, where $\eta$ is a root of
\[
Q_{B_t}(\eta)=q^2(q-4a^2)\eta^4-4(\tr\Sigma)q^2\eta^3+\bigl(3q^2+16q\det B_t\bigr)\eta^2
+4(\tr\Sigma)q\eta-4(\tr\Sigma)^2
\]
\textup{(}quantities of $B_t$\textup{)}, and $y_t$ is determined by $\eta$;
so there are at most $4$.
\item[\textup{(d)}] For $B$ in a nonempty Zariski open set, the solutions in
\textup{(b)} are nondegenerate and in \textup{(c)} there are exactly four,
all nondegenerate.
\end{enumerate}
\end{proposition}

\begin{proof}
(a) There are no rank-one loops and \eqref{app:bowtie:eq-limit} reads
$-\lambda_rW'=I$ with $\lambda_r\ne0$; its derivative is invertible.

(b) Here $\lambda_r=0$, so $V'=T\in\appbwT$, $\appbwadj T=-T$, and by
\eqref{app:bowtie:eq-limit2} $L_t(T)/y_t-I\in\C T$. Write $T=\xi X+\eta Y$.
By Lemma~\ref{app:bowtie:lem-coupling},
$L_t(T)=\eta(\tr\Sigma)qI+\xi\lambda_XX+\eta\lambda_YY$, so the condition is
$y_t=\eta(\tr\Sigma)q$ and $\xi\eta(\lambda_Y-\lambda_X)=\xi\eta q=0$. If
$\eta=0$ then $y_t=0$; so $\xi=0$, and $y_t^2=\Delta_t(\eta Y)=\eta^4\Delta_t(Y)$
gives $\eta^2=1/(q-4a^2)$. Both points are admissible:
$\det V'=-q/(q-4a^2)$.

(c) Here $\lambda_r=2$, $V'=-2I+T$ with $T=\xi X+\eta Y\in\appbwT$,
$\appbwadj V'=-(2I+T)$, and the condition is
$L_t(V')=\alpha I+\alpha'T$ with $\alpha-2\alpha'=y_t$. By the lemma,
$L_t(V')=(-2\gamma+\eta(\tr\Sigma)q)I+\xi\lambda_XX+(\eta\lambda_Y-2\tr\Sigma)Y$.
Parallelism of the $\appbwT$-parts gives $\xi(q\eta-2\tr\Sigma)=0$. If
$\xi\ne0$, then $\eta=2\tr\Sigma/q$, $\alpha'=\lambda_X$ and
$y_t=-2\gamma+2(\tr\Sigma)^2-2\lambda_X=0$. So $\xi=0$, $\eta\ne0$ (for
$\xi=\eta=0$ the condition would put $L_t(-2I)=-2\gamma I-2(\tr\Sigma)Y$ in
$\C I$), and, writing $s=\tr\Sigma$ and $\theta=q\eta-2s$,
\[
 \eta y_t=sq\eta^2-2(\gamma+\lambda_Y)\eta+4s=(s\eta-2)\theta .
\]
Since $q=s^2-4\det\Sigma$ and $\det V'=4-q\eta^2$, the first factor of
$\Delta_t(V')$ in Lemma~\ref{app:bowtie:lem-coupling} is
$\theta^2-4\det\Sigma\det V'=q(s\eta-2)^2$. Hence the equation
$\eta^2y_t^2=\eta^2\Delta_t(V')$ factors as
\[
 0=-(s\eta-2)^2
 \bigl[q\eta^2\bigl(\theta^2+4a^2(4-q\eta^2)\bigr)-\theta^2\bigr].
\]
The bracket is $Q_{B_t}(\eta)$. The factor $s\eta-2$ would give $y_t=0$,
so the admissible solutions come from roots of $Q_{B_t}$.

(d) Nondegeneracy of the two points of (b) and of the solutions of (c), and
the existence of four distinct roots of $Q_{B_t}$ with $\eta y_t\ne0$ and
$\det V'\ne0$, are expressed by the nonvanishing of polynomials in the
entries of $B$ (for (b), of the product of the Jacobian determinants at the
two points, which is a rational function of $B$; for (c), of the
discriminant of $Q_{B_t}$ and of the resultants of $Q_{B_t}$ with the
numerators of $\eta y_t$, $\det V'$ and of the Jacobian determinant on the
line). It suffices to check them at one point. At
\[
B_0=(B_1,B_2)=\left(\begin{pmatrix}-2&4\\3&-3\end{pmatrix},\begin{pmatrix}0&4\\2&5\end{pmatrix}\right)
\]
one has $(\tr\Sigma_t,q_t,q_t-4a_t^2,\det B_t)=(-5,50,49,-6)$ and $(5,61,57,-8)$.
In the coordinates $(V'_{11},V'_{12},V'_{22},y_t)$ and the equations
$((\cdot)_{11},(\cdot)_{12},(\cdot)_{22},y_t^2-\Delta_t)$ the Jacobian
determinant at both points of (b) is $98$ for $t=1$ and $114$ for $t=2$.
For (c), $Q_{B_1}=100(1225\eta^4+500\eta^3+27\eta^2-10\eta-1)$ and
$Q_{B_2}=212097\eta^4-74420\eta^3+3355\eta^2+1220\eta-100$ have nonzero
discriminants, and are coprime to $\eta y_t$, to $\det V'$ and to the
numerator and the denominator of the Jacobian determinant restricted to the
line; the residuals of \eqref{app:bowtie:eq-limit} vanish modulo $Q_{B_t}$.
\end{proof}

\paragraph{Type $(1,1)$.} Here $\lambda_r=1$, $V'=-I+T$ with
$T=\bigl(\begin{smallmatrix}a&b\\b&-a\end{smallmatrix}\bigr)$, and
$\appbwadj V'=-(V'+2I)$. Identify $\appbwSym_2$ with $\C^3$ by
$M\mapsto(M_{11},M_{12},M_{22})$ and write $[M,M',M'']$ for the
determinant of three vectors. Put
\begin{gather*}
D_1=[L_2V',I,V'],\quad D_2=[L_1V',I,V'],\quad n=-[L_1V',L_2V',V'+2I],\\
\xi_t=\tr(TX_t),\qquad \ell_3=[V',\appbwadj\Sigma_1,\appbwadj\Sigma_2],
\end{gather*}
where $L_tV'=L_t(V')$, and note $\theta_t(V')=\tr(\Sigma_tV')=\frac12\tr(TY_t)-\tr\Sigma_t$.
The identity
$\Sigma_tV'\Sigma_t=\theta_t\Sigma_t-\det(\Sigma_t)\appbwadj V'$
and Lemma~\ref{app:bowtie:lem-coupling} give
$L_tV'=2\theta_t\Sigma_t+\lambda_{X,t}(V'+2I)$.
The terms proportional to $V'+2I$ disappear from the three determinants.
Also $2[\Sigma_t,I,V']=-\tr(TX_t)=-\xi_t$, by expanding the determinant.
Finally adjugation on $\appbwSym_2$ preserves the determinant of three
vectors, since it acts as $+1$ on scalars and $-1$ on the traceless plane.
It sends $V'+2I$ to $-V'$. Thus
\begin{equation}\label{app:bowtie:eq-plane}
D_2=-\xi_1\theta_1(V'),\qquad D_1=-\xi_2\theta_2(V'),\qquad
n=4\,\theta_1(V')\theta_2(V')\,\ell_3 .
\end{equation}
Put $f_t=D_t^2\Delta_t(V')-n^2$ and
\[
g_1=\xi_2^2\Delta_1(V')-16(\theta_1\ell_3)^2,\qquad
g_2=\xi_1^2\Delta_2(V')-16(\theta_2\ell_3)^2,
\]
so that $f_1=\theta_2^2g_1$, $f_2=\theta_1^2g_2$, and $\deg g_t\le6$.

\begin{lemma}\label{app:bowtie:lem-plane}
Let $B$ satisfy $a_tq_t\det\Sigma_t\tr\Sigma_t\ne0$.
\textup{(i)} On $\{\det V'\ne0,\ y_1y_2\ne0\}$ the scheme $\appbwL^{11}(B)$
is isomorphic to the scheme in $(V',y,\nu)$ given by $\tr V'=-2$,
$H:=\rho_1L_1V'+\rho_2L_2V'-\nu(V'+2I)-I=0$ and $y_t^2=\Delta_t(V')$, where
$\rho_t=1/y_t$.
\textup{(ii)} Near a point with $n\ne0$ this scheme is isomorphic to
$\{f_1=f_2=0\}$ in the $(a,b)$-plane, with $\rho_1=D_1/n$,
$\rho_2=-D_2/n$. Hence the admissible solutions with $n\ne0$ are in
bijection with the points $P$ of $\{g_1=g_2=0\}$ with $n(P)\ne0$ and
$\det V'(P)\ne0$, and the multiplicity \textup{(}length of the local
ring\textup{)} of the solution equals the intersection number
$I_P(g_1,g_2)$; in particular the solution is nondegenerate iff
$I_P(g_1,g_2)=1$.
\textup{(iii)} If $n=0$ at a solution, then $L_1V',L_2V',I,V'$ span a
subspace of dimension at most $2$.
\end{lemma}

\begin{proof}
(i) Let $\Phi$ be the left side minus the right side of
\eqref{app:bowtie:eq-limit}. Then $\tr(V'\Phi)=-(\tr V'+2)$, so the ideal
contains $\tr V'+2$. Modulo $\tr V'+2$, $W'=-(V'+2I)/\det V'$ and
$\Phi=H|_{\nu=\nu_0}$ with $\nu_0=-(\rho_1\kappa_1+\rho_2\kappa_2+1)/\det V'$;
conversely, using $\tr V'^2=(\tr V')^2-2\det V'$, $\tr(V'H)\equiv2\det V'(\nu-\nu_0)$,
so $\nu-\nu_0$ lies in the ideal of $(H,\tr V'+2)$.
(ii) $H=A(\rho_1,\rho_2,\nu)^{\T}-I$ with columns
$L_1V',L_2V',-(V'+2I)$, and $\det A=n$. Where $n\ne0$, $H=0$ is the graph
$(\rho,\nu)=A^{-1}I$, and Cramer's rule gives $\rho_1=D_1/n$, $\rho_2=-D_2/n$.
Then $y_t^2=\Delta_t$ is $\rho_t^2\Delta_t=1$, i.e. $f_t=0$. Where $n\ne0$,
$\theta_1\theta_2\ne0$ by \eqref{app:bowtie:eq-plane}, so $(f_1,f_2)=(g_1,g_2)$
locally; and $f_t=0$, $n\ne0$ force $D_t\ne0$, i.e. $y_t$ finite and
nonzero. (iii) If $n=0$ then $I$ lies in the span of the columns of the
singular matrix $A$.
\end{proof}

\begin{proposition}\label{app:bowtie:prop-11}
For $B$ in a nonempty Zariski open set, $\appbwL^{11}(B)$ has exactly $20$
admissible solutions, all nondegenerate.
\end{proposition}

\begin{proof}
\emph{Six special points.} For $t=1,2$ let $M_t^\circ=-2\appbwadj\Sigma_t/\tr\Sigma_t$
and let $A_t^\pm=\sigma w_\pm w_\pm^{\T}$ with $w_\pm$ the isotropic vectors
of $\Sigma_t$ and $\sigma$ chosen with $\tr A_t^\pm=-2$ (possible when
$w_\pm^{\T}w_\pm\ne0$). They lie in the plane $\tr V'=-2$. At $A_t^\pm$,
$\theta_t=\det V'=0$, so $n=0$, and by
Corollary~\ref{app:bowtie:cor-coupling}(c) $\Delta_t\in(\theta_t,\det V')^2$,
hence $g_t\in(\theta_t,\det V')^2$ has multiplicity at least $2$; moreover
$g_{t'}(A_t^\pm)=0$ for $t'\ne t$. To see this, choose orthonormal
coordinates in which $V'=\operatorname{diag}(-2,0)$. They exist because
$w_\pm^{\T}w_\pm\ne0$. In these coordinates write
$\Sigma_t=\bigl(\begin{smallmatrix}0&b\\b&c\end{smallmatrix}\bigr)$ and
$(\Sigma_{t'})_{11}=d$. Then $\xi_t=4b$, $\ell_3=2bd$ and
$\theta_{t'}=-2d$; the relation
$16\ell_3^2=\xi_t^2\theta_{t'}^2$ is unchanged by orthogonal conjugation.
Since $\det V'=0$, one has $\Delta_{t'}=\theta_{t'}^4$, yielding the
claimed vanishing. At $M_t^\circ$, $\ell_3=0$ because
$M_t^\circ\in\C\appbwadj\Sigma_t$, and $\xi_t=0$ because
$M_t^\circ=-I+Y_t/\tr\Sigma_t$ and $\tr(X_tY_t)=0$; so $n=0$ and
$g_{t'}\in(\xi_t,\ell_3)^2$; and $\Delta_t$ vanishes to second order at
$M_t^\circ\in\C\appbwadj\Sigma_t$ (Corollary~\ref{app:bowtie:cor-coupling}(b)),
so $g_t$ has multiplicity at least $2$ too. By
$I_P(g_1,g_2)\ge m_P(g_1)m_P(g_2)$ \cite[Section~3.3]{Fulton}, if the six
points are distinct and $g_1,g_2$ have no common component, they contribute
at least $2+2+2+2+4+4=16$ to the intersection.

\emph{Upper bound.} Let $\mathcal U_1$ be the set of $B$ such that
$a_tq_t\det\Sigma_t\tr\Sigma_t\ne0$, $w_\pm^{\T}w_\pm\ne0$, the six points
are distinct, $X_1,X_2$ are linearly independent, $n(T=0)\ne0$,
$\tr(Y_1Y_2)\ne0$, $\Phi_1\Phi_2\ne0$ (defined below), the coefficient of
$b^6$ in $g_1$ is nonzero, and $\Res_b(g_1,g_2)$, formed with the formal
degrees $(6,6)$, is not the zero polynomial in $a$. The last two conditions
imply that $g_1,g_2$ have no common component (a common factor would divide
the constant coefficient of $b^6$ in $g_1$ if it did not involve $b$). All
conditions are nonvanishing conditions of polynomials in $B$ (for the
distinctness, of discriminants and resultants of the defining equations),
so $\mathcal U_1$ is open; it contains $B_0$ (see below). Let
$B\in\mathcal U_1$. By B\'ezout's theorem \cite[Section~5.3]{Fulton},
$\sum_{P\in\C^2}I_P(g_1,g_2)\le36$, so by Lemma~\ref{app:bowtie:lem-plane}(ii)
the admissible solutions with $n\ne0$ have total multiplicity at most
$36-16=20$. There are no admissible solutions with $n=0$. Indeed, by
Lemma~\ref{app:bowtie:lem-plane}(iii) all $3\times3$ minors of
$[L_1V',L_2V',I,V']$ vanish there, in particular $D_1=D_2=0$, i.e.
$\xi_1\theta_1=\xi_2\theta_2=0$. If $\xi_1=\xi_2=0$ then $T=0$, where
$n\ne0$. If $\theta_1=\theta_2=0$, Corollary~\ref{app:bowtie:cor-coupling}(c)
gives $L_tV'=\lambda_{X,t}(V'+2I)$, so
$H=(\rho_1\lambda_{X,1}+\rho_2\lambda_{X,2}-\nu)(V'+2I)-I\ne0$ since
$T\ne0$ on $\{\theta_t=0\}$. If $\xi_1=0$ and $\theta_2=0$, then $T=\eta Y_1$
is the point $T_1^*$ of $\C Y_1$ on $\{\theta_2=0\}$, $L_1$ preserves
$\mathrm{span}(I,Y_1)$, so $L_1V'=\alpha_1I+\alpha_1'V'$, and $H=0$ forces
$y_1=\alpha_1-2\alpha_1'$; then $y_1^2=\Delta_1$ requires
$\Phi_1:=(\alpha_1-2\alpha_1')^2-\Delta_1(V')=0$ at $T_1^*$, a rational
function of $B$. The case $\xi_2=0$, $\theta_1=0$ is symmetric
($\Phi_2$).

\emph{Lower bound.} At $B_0$ the twenty points of $\{g_1=g_2=0\}$ found in
the next paragraph are transversal intersections ($I_P=1$) with
$n\det V'\Delta_1\Delta_2D_1D_2\ne0$, hence twenty nondegenerate admissible
solutions of $\appbwL^{11}(B_0)$. By the implicit function theorem they
persist, nondegenerate and admissible, for $B$ in a Euclidean neighborhood
of $B_0$. The set of $B$ with at least twenty nondegenerate admissible
solutions is constructible (Chevalley) and contains a nonempty Euclidean open
set, so it contains a nonempty Zariski open set $\mathcal U_2$. For
$B\in\mathcal U_1\cap\mathcal U_2$ there are twenty solutions of
multiplicity one and total multiplicity at most $20$: exactly twenty, all
nondegenerate.

\emph{The point $B_0$.} The exact computation
gives polynomials $g_1,g_2\in\Z[a,b]$ of degree $6$; the coefficient of
$b^6$ in $g_1$ is $115200$, $\deg_bg_2=5$ (so the formal resultant of
degrees $(6,6)$ is $\pm115200$ times the following one), and
\[
\Res_b(g_1,g_2)=c\,(5a-3)^2(5a+4)^2(61a-11)^2(5a+1)^4(a+1)^6\,p_{20}(a)
\]
with $c\in\mathbb Q^*$ and $p_{20}$ irreducible over $\mathbb Q$ of degree
$20$. The six special points are
$M_1^\circ,A_1^\pm\leftrightarrow(a,b)=(-\tfrac15,-\tfrac75),(-\tfrac45,-\tfrac35),(\tfrac35,-\tfrac45)$
and $M_2^\circ,A_2^\pm\leftrightarrow(-1,\tfrac65),(\tfrac{11}{61},\tfrac{60}{61}),(-1,0)$,
distinct, accounting for the linear factors. The subresultant of degree one
in $b$ is $c_1(a)b+c_0(a)$ with $\gcd(c_1,p_{20})=1$; with
$b(a)=-c_0/c_1\bmod p_{20}$, both $g_1(a,b(a))$ and $g_2(a,b(a))$ vanish
modulo $p_{20}$, so above each of the twenty roots of the squarefree
$p_{20}$ there is exactly one common point. The values at $(a,b(a))$ of $n$,
$\det V'$, $\Delta_1$, $\Delta_2$, $D_1$, $D_2$ and of the Jacobian
$\partial(g_1,g_2)/\partial(a,b)$ have numerators coprime to $p_{20}$. The
remaining conditions of $\mathcal U_1$ at $B_0$: $X_1,X_2$ are independent,
$n(0)=2050$, $\Phi_1=3439461284/1874161$, $\Phi_2=-5612855040/1874161$,
$\tr(Y_1Y_2)=74$, and $w_\pm^{\T}w_\pm\ne0$ (the points $A_t^\pm$ exist).
\end{proof}

\subsection{No irregular branches}\label{app:bowtie:subsec-noirr}

Let $\mathcal G$ be the set of $B$ with $a_tq_t\det\Sigma_t\tr\Sigma_t(q_t-4a_t^2)\ne0$
$(t=1,2)$ for which the quartic forms $\Delta_1|_{\appbwT}$ and
$\Delta_2|_{\appbwT}$ have no common zero in $\appbwT\setminus\{0\}$; it is
open (a resultant of binary forms), and contains $B_0$ and the point $P_0$
below. In this subsection $\widetilde D$ is arbitrary with
$(\widetilde D,B)\in\Omega$, and a branch is a branch at infinity.

\begin{lemma}\label{app:bowtie:lem-A}
For $B\in\mathcal G$, on every branch $V'$ has a finite nonzero limit.
\end{lemma}

\begin{proof}
Let $\alpha=-\nu(V')$ and $V'=\tau^{-\alpha}(V'_*+o(1))$, $V'_*\ne0$. We
exclude $\alpha\ne0$ using \eqref{app:bowtie:eq-scaledV}, whose trace reads
$2(s_1+s_2-1)=\tr(V'(I-\tau\widetilde D))$ because $N_t$ is traceless.

\emph{$\alpha<0$.} Then $\tau^2C^0_t\to0$, the trace gives $s_1^0+s_2^0=1$,
so $\{s_1^0,s_2^0\}=\{1,0\}$, say $s_t\to0$, $s_{t'}\to1$. Then
$s_{t'}-1$ and $s_t$ are $O(\tau^{2+2\abs\alpha})$, and comparing the trace
at order $\tau^{\abs\alpha}$ gives $\tr V'_*=0$. The term
$\tau^2N_{t'}/s_{t'}$ is $O(\tau^{2+2\abs\alpha})$, while
$\tau^2N_t/s_t=\pm N_t(V')\sqrt{\sigma_+}/\sqrt{\Delta_t(V')}$ is homogeneous
of degree $0$ in $V'$ and must tend to $0$; hence $N_t(V'_*)=0$ (if
$\Delta_t(V'_*)=0$ and $N_t(V'_*)\ne0$ it would diverge). By
Corollary~\ref{app:bowtie:cor-coupling}(d), $V'_*=cX_t$ with $c\ne0$ and
$\Delta_t(X_t)\ne0$. Write $V'=\tau^{\abs\alpha}(cX_t+R)$ with $R\to0$. As $\theta_t(X_t)=0$,
$N_t(V')=\tau^{2\abs\alpha}\theta_t(R)M_t(cX_t+R)$, and
$M_t(X_t)=(\tr\Sigma_t)X_t-q_tJ$ has the nonzero antisymmetric part
$-q_tJ$. The antisymmetric part of the right side of
\eqref{app:bowtie:eq-scaledV} is $-\tau\,\mathrm{Anti}(V'\widetilde D)=O(\tau^{1+\abs\alpha})$,
so $\theta_t(R)=O(\tau^{1+\abs\alpha})$; then the symmetric traceless part of
the left side is $O(\tau^{1+\abs\alpha})$, while that of the right side is
$\tau^{\abs\alpha}(cX_t+o(1))$, a contradiction.

\emph{$0<\alpha<1$.} Again $\tau^2C^0_t=O(\tau^{2-2\alpha})\to0$ and $s_t$
is bounded, so the trace gives $\tr V'_*=0$, and the traceless part of the
right side is $\tau^{-\alpha}(V'_*+o(1))$. Terms with $s_t^0=\pm1$ are
$O(\tau^{2-2\alpha})=o(\tau^{-\alpha})$; terms with $s_t^0=0$ and
$\Delta_t(V'_*)\ne0$ are bounded. Hence some $t$ has $s_t^0=0$ and
$\Delta_t(V'_*)=0$; by $\mathcal G$ there is only one, and
$\tau^2N_t/s_t=\tau^{-\alpha}(V'_*+o(1))$. As
$N_t(V')=\tau^{-2\alpha}(N_t(V'_*)+o(1))$, this says
$N_t(V'_*)+o(1)=\varepsilon(\tau)(V'_*+o(1))$ for the scalar
$\varepsilon=\tau^{\alpha-2}s_t$, which therefore has a finite limit
$\varepsilon_0$. If $\varepsilon_0=0$ then $N_t(V'_*)=0$ and
$V'_*\in\C^*X_t$; otherwise $V'_*\in\C^*Y_t$. Both contradict
$\Delta_t(V'_*)=0$ by Corollary~\ref{app:bowtie:cor-coupling}(d).

\emph{$\alpha=1$.} In the unscaled variable $V=-\tau V'\to V_\infty\ne0$,
and \eqref{app:bowtie:eq-five} times $V$ reads
$\sum_tN_t(V)/s_t+(s_1+s_2-1)I=V(\widetilde D-zI)$, with $s_t$ bounded. The
trace gives $\tr V_\infty=0$, and the traceless part of the right side is
$-z(V_\infty+o(1))+O(1)$, unbounded. So some $s_t\to0$, which by $E_t$ forces
$\Delta_t(V_\infty)=0$; by $\mathcal G$ only one $t$, and
$N_t(V)=s_t(-z(V_\infty+o(1))+O(1))$; the scalar $-zs_t$ has a finite limit and the
argument of the previous case gives a contradiction.

\emph{$\alpha>1$.} Then $V=\tau^{-b}(V_*+o(1))$ with $b=\alpha-1>0$,
$s_t=O(\tau^{-b})$, and the trace gives $\tr V_*=0$. The traceless part of
the right side is $-\tau^{-1-b}(V_*+o(1))$. Loops with $\Delta_t(V_*)\ne0$
contribute $O(\tau^{-b})$ since $\tau^bs_t$ is bounded away from $0$. Hence
exactly one loop has $\Delta_t(V_*)=0$, the scalar $-\tau^{b-1}s_t$ has a finite limit,
and we conclude as before.
\end{proof}

Now let $V'\to V'_0\ne0$. Then $s_t^0\in\{1,-1,0\}$; define the type by
$r_t=1-s_t^0$ ($s_t^0=\pm1$), $r_t=1$ ($s_t^0=0$). The trace of
\eqref{app:bowtie:eq-scaledV} gives $\tr V'_0=-2\lambda_r$. For rank-one
loops $y_t=s_t/\tau^2$ satisfies $y_t^2\to\Delta_t(V'_0)$; a rank-one loop
with $\Delta_t(V'_0)=0$ is called \emph{bad}.

\begin{lemma}\label{app:bowtie:lem-B}
Let $B\in\mathcal G$ and $V'\to V'_0\ne0$ on a branch.
\textup{(a)} If there is no bad loop, $V'_0$ is invertible, so the branch is
regular.
\textup{(b)} For $B$ in a nonempty Zariski open set there is at most one bad
loop $t$, and then $N_t(V'_0)=0$, $\Delta_t(V'_0)=0$.
\textup{(c)} In the situation of \textup{(b)}, $V'_0$ is one of twelve
candidate points: for $r=(1,1)$, $M_t^\circ$ and $A_t^\pm$ (as in
Proposition~\ref{app:bowtie:prop-11}); for $r=(1,2)$, $2M_1^\circ$ and
$2A_1^\pm$; for $r=(2,1)$, $2M_2^\circ$ and $2A_2^\pm$. Types $(0,1)$ and
$(1,0)$ have no candidates.
\end{lemma}

\begin{proof}
(a) Without bad loops the limit of \eqref{app:bowtie:eq-scaledV} is
$\sum_{r_t=1}N_t(V'_0)/y_t^0-\lambda_rI=V'_0$ with $y^0_t\ne0$. Suppose
$V'_0=\sigma ww^{\T}$. Then $V'_0w=(\tr V'_0)w=-2\lambda_rw$, and by
Corollary~\ref{app:bowtie:cor-coupling}(e) $N_t(V'_0)w=\kappa_tw$ with
$(y_t^0)^2=\Delta_t(V'_0)=\kappa_t^2\ne0$, so $\kappa_t/y_t^0=\pm1$. Applying the limit equation to $w$ gives
$\sum_{r_t=1}(\pm1)=-\lambda_r$, impossible: with no rank-one loop
$\lambda_r\in\{-1,1,3\}$; with one, $\lambda_r\in\{0,2\}$; with two,
$\lambda_r=1$.

(b) A single bad loop $t$: the other terms of \eqref{app:bowtie:eq-scaledV}
are bounded, so $N_t(V')/y_t$ is bounded with $y_t\to0$, whence
$N_t(V'_0)=0$. Two bad loops occur only for $r=(1,1)$; then
$\Delta_1(V'_0)=\Delta_2(V'_0)=0$, and since $y_2/y_1$ has a limit
$c\in\C\cup\{\infty\}$ (Puiseux series), say $c\ne\infty$,
$y_2(N_1/y_1+N_2/y_2)=N_2+(y_2/y_1)N_1\to0$ gives
$N_2(V'_0)=-cN_1(V'_0)$, so $N_1(V'_0)\wedge N_2(V'_0)=0$. In the
projective closure of the plane $\tr V'=-2$ the conditions
$\Delta_1=\Delta_2=0$, $N_1\wedge N_2=0$ define a closed subset of
$\PP^2\times\mathrm{GL}_2^2$; its projection to $\mathrm{GL}_2^2$ is closed.
At the point $P_0$ below the affine part is empty (the Gr\"obner basis of
the ideal generated by $\Delta_1,\Delta_2$ and the $2\times2$ minors of
$(N_1,N_2)$, as $2\times4$ matrix of entries, is $\{1\}$), and at infinity
the top-degree forms $\Delta_t|_{\appbwT}$ have no common zero. So for $B$ off a proper closed set there are no
two bad loops.

(c) By Corollary~\ref{app:bowtie:cor-coupling}(c), $V'_0$ lies on one of
the lines $\C\appbwadj\Sigma_t$, $\C w_\pm w_\pm^{\T}$ and has trace
$-2\lambda_r$. For $\lambda_r=0$ there is no nonzero such point:
$\tr\appbwadj\Sigma_t=\tr\Sigma_t\ne0$, and a traceless $V$ with
$\theta_t(V)=0$ lies in $\C X_t$ by (d) of the corollary, while
$\det(cX_t)=-c^2q_t$; so $\det V=0$ forces $V=0$. For $\lambda_r=1,2$ the
intersections are the listed points.
\end{proof}

For a candidate $V'_0$ with bad loop $t$ let
$\mathcal R=V'_0+\lambda_rI-\sum N_o(V'_0)/y_o^0$, the sum over good rank-one
loops $o$ (at most one), with $y_o^0=\pm\sqrt{\Delta_o(V'_0)}$. The limit of
\eqref{app:bowtie:eq-scaledV} gives $N_t(V')/y_t\to\mathcal R$.

\begin{lemma}\label{app:bowtie:lem-balance}
\textup{(a)} At $A$-type candidates $V'_0=A_t^\pm$ or $2A_t^\pm$
\textup{(}$\theta_t(V'_0)=0$, $M_t(V'_0)\ne0$\textup{)}, the scalar
$\omega=\theta_t(V')/y_t$ converges to $\omega_0$ with
$\mathcal R=\omega_0M_t(V'_0)$.
\textup{(b)} At $M$-type candidates $V'_0=M_t^\circ$ or $2M_t^\circ$
\textup{(}$M_t(V'_0)=0$, $\theta_t(V'_0)\ne0$\textup{)}, there are a sign
and $d\in\appbwT\setminus\{0\}$ with $\mathcal H(d)\ne0$ and
$\theta_t(V'_0)M_t(d)=\pm\sqrt{\mathcal H(d)}\,\mathcal R$, where
$\mathcal H(d)$ is the second-order term of $\Delta_t(V'_0+d)$. Equivalently:
$\mathcal R\ne0$, $\mathcal R=\theta_t(V'_0)M_t(d_{\mathcal R})$ for some
$d_{\mathcal R}\in\appbwT$, and $\mathcal H(d_{\mathcal R})=1$.
\end{lemma}

\begin{proof}
(a) $N_t/y_t=\omega M_t(V')$ and $M_t(V')\to M_t(V'_0)\ne0$. (b) Write
$V'=(1+c)V'_0+\delta'$ with $\delta'\in\appbwT$; then
$M_t(V')=M_t(\delta')$, and by homogeneity and
Corollary~\ref{app:bowtie:cor-coupling}(b),
$\Delta_t(V')=(1+c)^2\mathcal H(\delta')+O(\abs{\delta'}^3)$ and
$y_t^2=\Delta_t(V')(1+O(\tau^2))$. If $\delta'\equiv0$ then $y_t\equiv0$.
Otherwise $\delta'=\tau^\gamma(d+o(1))$, and if $\mathcal H(d)=0$ then
$\abs{M_t(\delta')/y_t}\to\infty$ since $M_t$ is injective on $\appbwT$
($\ker M_t=\C\appbwadj\Sigma_t$). The equivalence uses this injectivity:
$\theta M_t(d)=\lambda\mathcal R$ with $\lambda^2=\mathcal H(d)\ne0$ means
$d=\lambda d_{\mathcal R}$ and $\mathcal H(d_{\mathcal R})=1$.
\end{proof}

\begin{proposition}\label{app:bowtie:prop-noirr}
For $B$ in a nonempty Zariski open set and every $\widetilde D$ with
$(\widetilde D,B)\in\Omega$, every branch at infinity is regular.
\end{proposition}

\begin{proof}
By Lemmas~\ref{app:bowtie:lem-A} and~\ref{app:bowtie:lem-B} it remains to
exclude branches converging to the twelve candidates. Each candidate and
each sign of $y_o^0$ is defined over a finite covering $\mathcal W$ of the
space of $B$ (adjoining $w_\pm$ and $\sqrt{\Delta_o(V'_0)}$). The necessary
conditions of Lemma~\ref{app:bowtie:lem-balance} are closed conditions on
$\mathcal W$, or on $\mathcal W\times\PP(\appbwT)$ in case (b) (written as
$\theta M_t(d)\wedge\mathcal R=0$ and
$(\theta M_t(d))_{ij}^2=\mathcal H(d)\mathcal R_{ij}^2$ for all $i,j$, which
by injectivity of $M_t$ on $\appbwT$ is equivalent to the condition of the
lemma). Their images in the space of $B$ are closed. They fail at the
rational point
\[
P_0=\left(\begin{pmatrix}1&\frac57\\-\frac57&-\frac49\end{pmatrix},
\begin{pmatrix}1&-\frac{11}{36}\\-\frac{43}{36}&\frac{161}{400}\end{pmatrix}\right)
\]
for all candidates except the four points $A_t^\pm$ of type $(1,1)$: for $2A_t^\pm$, $\mathcal R=V'_0+2I$ and
$M_t(V'_0)$ are not parallel; for $2M_t^\circ$ and for $M_t^\circ$ (both
signs of $y_o^0$), $\mathcal R\in\theta M_t(\appbwT)$ but
$\mathcal H(d_{\mathcal R})\ne1$. At the four $A_t^\pm$ of type $(1,1)$
the condition (a) has, for exactly one sign, the solution
$y_o^0=-\kappa_o(A_t^\pm)$ and $\omega_0\ne0$; at $P_0$ the values of
$\omega_0$ are $\{-4401/13298,-99/1102\}$ for $t=1$ and
$\{-7475/3161,-2275/1159\}$ for $t=2$. These four candidates are excluded by
Lemma~\ref{app:bowtie:lem-blowup}.
\end{proof}

\begin{lemma}\label{app:bowtie:lem-blowup}
For $B$ in a nonempty Zariski open set and every $\widetilde D$, no branch
of type $(1,1)$ converges to $A_t^\pm$ with bad loop $t$.
\end{lemma}

\begin{proof}
Let $o\ne t$, $V_A=A_t^\pm=\sigma ww^{\T}$ ($w^{\T}\Sigma_tw=0$,
$\sigma w^{\T}w=-2$), and $K=\kappa_o(V_A)$; then $\Delta_o(V_A)=K^2$ by
Corollary~\ref{app:bowtie:cor-coupling}(e), and $K\ne0$ because the loop $o$
is not bad. Consider the unknowns $x=(V',\omega,y_o)$, put
$y_t=\theta_t(V')/\omega$, and let
\begin{gather*}
\widetilde G=\omega M_t(V')+\frac{N_o(V')}{y_o}+\bigl(\tau^2(y_t+y_o)-1\bigr)I-V'(I-\tau\widetilde D),\\
\widetilde E_t=\omega^2\theta_t^2(1+2\tau^2\kappa_t)-\tau^4\theta_t^4-\omega^4\Delta_t,\qquad
E_o=y_o^2(1+2\tau^2\kappa_o-\tau^4y_o^2)-\Delta_o .
\end{gather*}
These are \eqref{app:bowtie:eq-scaledV} in the chart $\mathcal R^{11}$ with
$N_t/y_t=\omega M_t$, and $\omega^4$ times the chart equation of $y_t$. They
are analytic near $p_A=(V_A,\omega_0,y_o^0)$, $\tau=0$, and $p_A$ solves
them at $\tau=0$ by Lemma~\ref{app:bowtie:lem-balance}(a); by (c) below
(at $\tau=0$, $s'=1$) this forces $y_o^0=-K$ and determines $\omega_0$.

(a) An admissible branch converging to $V_A$ gives, by
Lemma~\ref{app:bowtie:lem-balance}(a), a solution $x(\tau)\to p_A$ of
$\widetilde G=0$, $\widetilde E_t=E_o=0$ with $\theta_t(V')=\omega y_t\ne0$
(as $\omega_0\ne0$).

(b) Let $\Gamma$ be defined by $\widetilde G_{11}=\widetilde G_{12}=\widetilde G_{22}=0$
and $E_o=0$ near $(p_A,0)$; it may be larger than the solution set of
$\widetilde G=0$, $E_o=0$, since the argument only uses that it contains the
points of the branch and $x_{\mathrm{sp}}$. At $P_0$, for all
four points, the $4\times5$ Jacobian of these equations with respect to $x$
at $(p_A,0)$ has rank $4$, so $\Gamma$ is a smooth surface with coordinates
$(s,\tau)$, $s$ a coordinate along the kernel vector $v$; and
$Q:=\theta_t(v)^2-\omega_0^2\mathcal H(v)\ne0$, where $\mathcal H$ is the
second-order term of $\Delta_t$ at $V_A$ (the first-order term vanishes
since $\Delta_t\in(\theta_t,\det V')^2$). As
$\theta_t(V_A)=0$, $f:=\widetilde E_t|_\Gamma$ satisfies
$f(s,0)=\omega_0^2Qs^2+O(s^3)$, because $\nabla\theta_t^2$ and
$\nabla\Delta_t$ vanish at $V_A$. The conditions $\omega_0\ne0$, rank $4$ and
$Q\ne0$ are open conditions on the finite covering $\mathcal W$ that hold
over $P_0$ (at all four points), so they hold for generic $B$. By the Weierstrass
preparation theorem $f(\cdot,\tau)$ has exactly two zeros near $0$, counted
with multiplicity, for every small $\tau$.

(c) On $S=\{V'=s'V_A\}$ we have $\theta_t=\det V'=0$, hence $y_t=0$,
$\Delta_t=0$ and $\widetilde E_t=0$. Here $M_t(V_A)=2\sigma w(\Sigma_tw)^{\T}=c\,w(Jw)^{\T}$
with $c\ne0$ (as $\Sigma_tw\perp w$ and $\det\Sigma_t\ne0$), so $M_t(V_A)w=0$
and $M_t(V_A)Jw=c(w^{\T}w)w$; and $N_o(V_A)w=Kw$. Since $w,Jw$ is a basis,
$\widetilde G=0$ on $S$ is equivalent to the vanishing of the coefficients of
$w$ in $\widetilde Gw$ and of $w$, $Jw$ in $\widetilde GJw$:
\begin{gather*}
\frac{s'^2K}{y_o}+\tau^2y_o-1+2s'+\tau s'\sigma\,w^{\T}\widetilde Dw=0,\qquad
-\frac{s'^2K}{y_o}+\tau^2y_o-1=0,\\
\omega s'c\,(w^{\T}w)+\frac{s'^2\sigma\,w^{\T}L_o(V_A)Jw}{y_o}+\tau s'\sigma\,w^{\T}\widetilde DJw=0,
\end{gather*}
and $E_o$ follows from the second equation, as
$E_o=y_o^2-(s'^2K-\tau^2y_o^2)^2$ on $S$. At $\tau=0$ the solution is
$s'=1$, $y_o=-K$, $\omega=\omega_0$, and the Jacobian in $(s',y_o,\omega)$
is block triangular with blocks
$\bigl(\begin{smallmatrix}0&-1/K\\2&1/K\end{smallmatrix}\bigr)$ and
$c\,w^{\T}w\ne0$. So there is an analytic solution $x_{\mathrm{sp}}(\tau)\in S$,
$x_{\mathrm{sp}}(0)=p_A$, for every $\widetilde D$; it is not admissible
($y_t\equiv0$).

(d) $x_{\mathrm{sp}}$ lies on $\Gamma$ as $s=s_{\mathrm{sp}}(\tau)$, and
$\theta_t|_\Gamma$, $\det V'|_\Gamma$ vanish there, so they are divisible by
$s-s_{\mathrm{sp}}(\tau)$ in $\C\{s,\tau\}$. Since
$\widetilde E_t\in(\theta_t,\det V')^2$, $f\in(s-s_{\mathrm{sp}}(\tau))^2$;
by (b) the two zeros of $f(\cdot,\tau)$ are both at $s_{\mathrm{sp}}(\tau)$.
The points of the branch of (a) lie on $\Gamma$ and are zeros of $f$, hence
equal to $x_{\mathrm{sp}}(\tau)$, contradicting $\theta_t\ne0$.
\end{proof}

\subsection{The exact computations}\label{app:bowtie:subsec-comp}

All computations are exact (SymPy, rational arithmetic), and each was
carried out twice, by independently written programs. We list what was
computed.

The polynomial identities used in the count, including
\eqref{app:bowtie:eq-plane} and the vanishing at $A_t^\pm$, were derived
above; they were also checked by expansion.

\emph{(C1) At $B_0$.} The data of Proposition~\ref{app:bowtie:prop-single}(d)
and of the last paragraph of the proof of
Proposition~\ref{app:bowtie:prop-11}: the polynomials $g_1,g_2$, the
resultant $\Res_b(g_1,g_2)$ of degree $36$ and its factorization, the
irreducibility and squarefreeness of $p_{20}$, the degree-one subresultant,
the reductions modulo $p_{20}$ and the gcds. The admissibility and
transversality of the twenty common points are obtained from the exact gcd
computation described above (subresultant, reduction modulo $p_{20}$, gcds
with $p_{20}$, including the Jacobian). The single-loop data of
Proposition~\ref{app:bowtie:prop-single}(d), including the Jacobian
determinants, are computed exactly in the same way. The conditions defining
$\mathcal G$ hold at $B_0$: the two quartic forms, as binary forms, have no
common zero; in the affine coordinate used, one of them vanishes at the
point at infinity, where the other does not.

\emph{(C2) At $P_0$.} The Gr\"obner basis $\{1\}$ and the resultant of the
quartic forms in Lemma~\ref{app:bowtie:lem-B}(b); for each of the twelve
candidates of Lemma~\ref{app:bowtie:lem-B}(c) the exact test of
Lemma~\ref{app:bowtie:lem-balance} (linear algebra over $\mathbb Q$ or over
$\mathbb Q(\sqrt{\Delta_o(V'_0)})$); at the four points $A_t^\pm$ the
residuals of the blown-up system, the rank $4$ of the Jacobian (over
$\mathbb Q$), the kernel vector, $Q\ne0$ and $c\ne0$. The same checks also succeed at $B_0$
(for instance $\omega_0\in\{-39/115,\,1/160\}$ for $t=1$ and
$\{17/4416,\,7/24\}$ for $t=2$).

%% file: appendices/certificates.tex
\section{Computer-assisted statements and their certificates}\label{app:cert}

\providecommand{\appcJ}{\mathcal J}
\providecommand{\appcK}{\mathcal K}
\providecommand{\appcSym}{\operatorname{Sym}}

This appendix describes the computations used in the paper. The proofs
using them are given alongside the statements or in
Appendices~\ref{app:many}--\ref{app:bowtie}. The computer-assisted examples
and numerical values are listed in Table~\ref{app:cert:tab}. The exact
coefficients and boxes for the block root-count certificates are
available in the ancillary files linked from the arXiv abstract page.
The accompanying \path{README.md} lists the datasets and describes their format.

\subsection{Methods}\label{app:cert:subsec-methods}

\subsubsection{Exact computations}
A computation is called \emph{exact} if every number in it is an integer, a
rational number, or an element of an explicitly given number field
$\mathbb Q[v]/(P)$ with $P\in\mathbb Q[v]$ irreducible, and every operation
is carried out without rounding. The exact computations of this paper are of
four kinds: enumeration of finite sets (edge subsets, vertex subsets,
configurations, permutation groups by the Schreier--Sims algorithm), with
integer arithmetic in the Python standard library; linear algebra over
$\mathbb Q$ (inverses, determinants, characteristic polynomials, $LDL^{\T}$
decompositions), in rational arithmetic; polynomial arithmetic over $\mathbb Q$ (expansion,
resultants and subresultants, factorization into irreducible factors,
gcds, Gr\"obner bases in lexicographic order, implicit differentiation at a
rational point), with SymPy~1.14.0; and arithmetic in a residue field
$\mathbb Q[v]/(P)$, where an element is represented by its remainder modulo
$P$ and inverses are computed by the extended Euclidean algorithm. A
statement that a polynomial expression is nonzero in such a residue field
means that it is nonzero at every complex root of~$P$. SymPy is used only
through these exact operations; no floating point evaluation enters an
exact certificate.

\subsubsection{Interval certificates}
An \emph{interval} is a compact interval $[a,b]\subset\R$ and a
\emph{complex interval} is a rectangle $[a,b]+i[c,d]\subset\C$; vectors and
matrices of (complex) intervals are products of such sets, and are called
boxes. An interval extension of a function $\phi$ is a procedure that
returns, for every box $X$ in its domain, a box containing $\phi(X)$. All
interval certificates below rest on the following existence and uniqueness
test for Krawczyk's operator~\cite{Krawczyk}. We include the proof because we use it for a general
convex enclosure of the Jacobian and for an arbitrary center.

\begin{lemma}[Krawczyk test]\label{app:cert:lem-krawczyk}
Let $X\subset\R^N$ be a box with nonempty interior, $x_0\in X$, let $F$ be a
$C^1$ map from an open neighborhood of $X$ to $\R^N$, let $\appcJ$ be a
closed convex set of real $N\times N$ matrices with $F'(x)\in\appcJ$ for all $x\in
X$, and let $Y$ be a real $N\times N$ matrix. Suppose that
\begin{equation}\label{app:cert:eq-krawczyk}
x_0-YF(x_0)+(I-YJ)(x-x_0)\in\operatorname{int}X\qquad
\text{for all }J\in\appcJ\text{ and }x\in X .
\end{equation}
Then $Y$ and every $J\in\appcJ$ are invertible, $F$ has exactly one zero
$x^*$ in $X$, $x^*$ lies in the interior of $X$, and $F'(x^*)$ is invertible.
The same holds for a holomorphic map $F$ from an open neighborhood of a
box $X\subset\C^N$ of complex intervals to $\C^N$, with $\appcJ$ a closed
convex set of complex matrices and $Y$ complex.
\end{lemma}

\begin{proof}
Put $c=-YF(x_0)$ and $D=X-x_0$, a compact convex set. For $J\in\appcJ$ let
$M=I-YJ$ and $\phi(d)=c+Md$. By \eqref{app:cert:eq-krawczyk},
$\phi(D)\subseteq\operatorname{int}X-x_0=\operatorname{int}D$. By Brouwer's
theorem $\phi$ has a fixed point $d^*\in D$. If $YJ$ were singular, there
would be $v\ne0$ with $Mv=v$; the set $\{s\in\R:d^*+sv\in D\}$ is a compact
interval $[s_-,s_+]$, the point $d^*+s_+v$ lies on the boundary of $D$, and
$\phi(d^*+s_+v)=d^*+s_+v$, contradicting $\phi(D)\subseteq\operatorname{int}D$.
Hence $YJ$ is invertible, and so are $Y$ and~$J$.

For $x\in X$ the segment from $x_0$ to $x$ lies in $X$, and
$F(x)-F(x_0)=\bar J(x)(x-x_0)$ with
$\bar J(x)=\int_0^1F'(x_0+t(x-x_0))\,dt$. The matrix $\bar J(x)$ lies in
$\appcJ$, since it is a limit of Riemann sums, which are convex combinations
of elements of the closed convex set $\appcJ$. Then
$g(x)=x-YF(x)=x_0+c+(I-Y\bar J(x))(x-x_0)$ lies in $\operatorname{int}X$
by \eqref{app:cert:eq-krawczyk}. So the continuous map $g$ sends $X$ into
$\operatorname{int}X$ and has a fixed point $x^*\in\operatorname{int}X$,
which is a zero of $F$ because $Y$ is invertible. If $y^*\in X$ is another
zero, then $0=F(y^*)-F(x^*)=\bar J'(y^*-x^*)$ with
$\bar J'=\int_0^1F'(x^*+t(y^*-x^*))\,dt\in\appcJ$, which is invertible, so
$y^*=x^*$. Finally $F'(x^*)\in\appcJ$ is invertible.

In the holomorphic case identify $\C^N$ with $\R^{2N}$. A box of complex
intervals is a real box, the real derivative of $F$ is the realification of
the complex derivative, the realifications of the matrices in $\appcJ$ form a
closed convex set, and a complex matrix is invertible if and only if its
realification is. The real case applies.
\end{proof}

In practice $F(x_0)$ cannot be evaluated exactly and is replaced by an
interval enclosure, and the left-hand side of
\eqref{app:cert:eq-krawczyk} is replaced by an interval vector
$\appcK(X)$ computed by interval arithmetic from an enclosure of $F(x_0)$, a
floating point matrix $Y$ (usually an approximate inverse of $F'(x_0)$; any
matrix is allowed), an interval matrix $\appcJ\supseteq F'(X)$, which is
closed and convex, and the box $X-x_0$ or a box containing it. If
$\appcK(X)\subseteq\operatorname{int}X$, then \eqref{app:cert:eq-krawczyk}
holds. In the \emph{parametric} form the coefficients of $F$, or the
spectral parameter $z$, range over a box $C$, the enclosures of $F(x_0)$ and
of $F'$ are taken over $C$, respectively over $X\times C$, and the
conclusion of Lemma~\ref{app:cert:lem-krawczyk} holds for every parameter
in~$C$.

The signs and inertia of symmetric matrices are certified as follows.

\begin{lemma}[Interval inertia]\label{app:cert:lem-inertia}
Let $\mathbf H$ be a box of real symmetric $n\times n$ matrices, let $Q$ be
a real $n\times n$ matrix with $\norm{Q^{\T}Q-I}_\infty<1$, let $\mathbf M$
be a box of matrices containing $Q^{\T}HQ$ for all $H\in\mathbf H$, let
$d_1,\dots,d_n$ be real numbers and $e\ge\max_i\sum_j\abs{M_{ij}-d_i\delta_{ij}}$
for all $M\in\mathbf M$. If $\min_i\abs{d_i}>e$, then every symmetric
$H\in\mathbf H$ is invertible and has exactly $\#\{i:d_i<0\}$ negative
eigenvalues.
\end{lemma}

\begin{proof}
$Q^{\T}Q$ is invertible, so $Q$ is. Let $H\in\mathbf H$ be symmetric and
$M=Q^{\T}HQ$, a symmetric matrix in $\mathbf M$. For $t\in[0,1]$ the symmetric
matrix $M_t=\diag(d)+t(M-\diag(d))$ has, by Gershgorin's theorem, all its
eigenvalues in $\bigcup_i[d_i-e,d_i+e]$, which does not contain $0$. The
eigenvalues depend continuously on $t$, so $M=M_1$ has as many negative
eigenvalues as $M_0=\diag(d)$, and none is zero. By Sylvester's law of
inertia the same holds for~$H$.
\end{proof}

The hypothesis on $Q$ is used only for the invertibility of $Q$, which also
follows from the other hypotheses: by Gershgorin's theorem every matrix in
$\mathbf M$, in particular $Q^{\T}HQ$, is invertible.

\begin{corollary}\label{app:cert:cor-certificate}
Let $F$ be one of the polynomial or rational systems of this paper at a
fixed coefficient point $P^*$ (or for all coefficients in a box), and let
$X_1,\dots,X_K$ be boxes for which the hypothesis of
Lemma~\ref{app:cert:lem-krawczyk} has been verified.
\begin{enumerate}
\item[\textup{(a)}] If the boxes are pairwise disjoint, the zeros $x_k^*\in X_k$
are $K$ distinct nondegenerate solutions.
\item[\textup{(b)}] If, in addition, interval enclosures over $X_k$ of the
determinants (or coordinates) that define admissibility exclude $0$, then
$x_k^*$ is admissible; if enclosures of the coordinates that define
positivity are positive, then $x_k^*$ is positive.
\item[\textup{(c)}] If an interval matrix $\mathbf H$ that contains
$\Hess(x_k^*)$, or a matrix congruent to it such as $H_u$ of
\eqref{real:eq:Hu} (Lemma~\ref{real:lem:phi}(b)), satisfies the hypothesis of Lemma~\ref{app:cert:lem-inertia}, the Morse
index of $x_k^*$ is $\#\{i:d_i<0\}$.
\item[\textup{(d)}] If the coefficient space is irreducible and the number
of admissible solutions equals a constant $N$ on a nonempty Zariski open
subset, then \textup{(a)} and \textup{(b)} at one point $P^*$ give $N\ge K$.
Indeed, by the implicit function theorem the $x_k^*$ continue to
nondegenerate solutions on a Euclidean neighborhood of $P^*$, which stay
distinct and admissible, and this neighborhood meets the Zariski open set.
\end{enumerate}
\end{corollary}

Three implementations of interval arithmetic are used.
\begin{enumerate}
\item[\textup{(I1)}] IEEE binary64 arithmetic with rounding to nearest,
with each elementary operation ($+,-,\times,\div,\sqrt{\ }$) on interval
endpoints followed by a step to the adjacent floating-point number
outwards. Coefficients are binary64 numbers used as point intervals, or
rational numbers enclosed by the two neighboring binary64 numbers.
For the full cavity system, the stored binary64 endpoints themselves
define the box. The test uses an outward enclosure of that box minus
its center, with outward rounding of each operation in the matrix
products and sums. Inclusion and separation are checked by direct
comparisons of endpoints.
\item[\textup{(I2)}] The rectangular complex interval arithmetic of mpmath~1.3.0
with $160$-bit end points and outward rounding.
\item[\textup{(I3)}] Intervals with rational end points, rounded outward
after every operation to $96$ significant bits, with square roots enclosed
by integer square roots.
\end{enumerate}

\subsubsection{Floating point computations}
Floating point computations (path tracking, monodromy, Newton's method,
grid searches) are used to find candidate solutions, coefficient points and
curves, which are then certified; the certificates do not depend on how the
candidates were found. Floating point results are never used as proofs. The
paper reports floating point results as evidence in several places, each
marked as such where it occurs; among them are the survey of ramification
points after Table~\ref{gen:tab:values} and the numbers of solutions found
by monodromy that exceed the certified lower bounds, which are mentioned
below.

\subsection{Catalog}\label{app:cert:subsec-catalogue}

Table~\ref{app:cert:tab} lists the items; each is described in one of the
following subsections. ``Exact'' and ``interval'' refer to the methods
above.

\begin{table}[ht]
\centering\small
\begin{tabular}{@{}>{\raggedright\arraybackslash}p{0.08\textwidth}>{\raggedright\arraybackslash}p{0.56\textwidth}>{\raggedright\arraybackslash}p{0.26\textwidth}@{}}
\toprule
Item & Statement & Method\\
\midrule
\ref{app:cert:item-block-lower} &
Thm.~\ref{thm:block-degree}(iv): bowtie $N_3\ge1079$;
diamond $N_2\ge66$ & interval (I1), (I2)\\
\ref{app:cert:item-bowtie-upper} &
Thm.~\ref{app:bowtie:thm-36}: bowtie $N_2=36$ & exact\\
\ref{app:cert:item-four-port} &
Prop.~\ref{prop:block-four-port}(ii): $N_{(1,p)}(D_4)=14$, $p\ge2$ &
exact\\
\ref{app:cert:item-K4} &
Examples~\ref{blr:ex-K4}, \ref{blc:ex:k4} ($p=2$ instance) & exact (check)\\
\ref{app:cert:item-many-solutions} &
Prop.~\ref{app:many:prop-certified}: exactly $6$, resp.\ $10$, positive
solutions & interval (I3), (I1)\\
\ref{app:cert:item-circles} &
Prop.~\ref{app:bridge:prop-M}: circles for $M_1$, $M_2$ & interval
(I1), parametric\\
\ref{app:cert:item-genus} &
Table~\ref{gen:tab:values}: $\LG$, $\Sigma(G)$, $g$, $R$ & exact\\
\ref{app:cert:item-census} &
Rem.~\ref{cls:rem:census}: Cor.~\ref{cls:cor:identity} for $33867$ graphs &
exact (check)\\
\ref{app:cert:item-L} &
Table~\ref{deg:table}: $\LG$, $Q(G)$ & exact\\
\bottomrule
\end{tabular}
\caption{The computer-assisted statements of the paper and the methods
used.}
\label{app:cert:tab}
\end{table}

\subsubsection{Lower bounds for block degrees}\label{app:cert:item-block-lower}

\emph{Statements.} Theorem~\ref{thm:block-degree}(iv):
$N_3(\text{bowtie})\ge1079$ and $N_2(\text{diamond})\ge66$.
The certificates for $N_2(\text{bowtie})\ge36$ also recorded here are
independent checks of Theorem~\ref{app:bowtie:thm-36}.
The ancillary data give the coefficient points and all accepted
boxes for the two full cavity systems and for the reduced system with $p=3$.
All real and imaginary parts are given as exact dyadic numbers.

\emph{Diamond, $p=2$.} The unknowns are the ten symmetric $2\times2$ cavity
precisions, $30$ complex numbers. The coefficient point $P^*$ consists of
complex symmetric $D_i$ and complex $A_e$, given as binary64 numbers.
For each candidate $\tilde x$, form the real and imaginary endpoints
by adding and subtracting $r_j=\rho(1+\abs{\tilde x_j})$ in binary64
arithmetic, trying $\rho=10^{-9}$ first and then
$10^{-13},\dots,10^{-6}$. The endpoints define $X$ as in (I1).
Take $Y$ to be an approximate inverse of the Jacobian at $\tilde x$.
The Jacobian enclosure over $X$ uses the derivative
$-A_{ik}Q^{-1}EQ^{-1}A_{ki}$ of the messages. On the same box the enclosures of
$\det Q^e_{i\to j}$, $\det R_i$ and $\det E_e$ (the latter as
$\det Q_{i\to j}\det(Q_{j\to i}-A_{ji}Q_{i\to j}^{-1}A_{ij})$) exclude $0$.
The Krawczyk image is strictly inside $X$ for $66$ of the $67$ candidates.
For every pair of accepted boxes, some real or imaginary coordinate
intervals are disjoint by direct comparison of their endpoints.
Hence $N_2(\text{diamond})\ge66$. Floating-point monodromy found
$67$ candidates.

\emph{Bowtie, $p=2$, full system.} The same test, applied to $36$
candidates at one coefficient point ($36$ complex unknowns), certifies all
$36$, with pairwise disjoint boxes. The bound
$N_2(\text{bowtie})\ge36$ is also obtained below with the arithmetic (I2).

\emph{Reduced normal form.} For larger $p$ we use the system for $m$
triangles sharing the vertex $0$ (the bowtie is $m=2$). Its unknowns are
symmetric $p\times p$ matrices $U_1,\dots,U_m,V$, its coefficients are
$\widetilde D\in\appcSym_p(\C)$ and $B_1,\dots,B_m\in\mathrm{GL}_p(\C)$, and it
reads
\begin{equation}\label{app:cert:eq-rose}
U_tVU_t-U_t-B_tVB_t^{\T}=0\quad(1\le t\le m),\qquad
\sum_{t=1}^m\bigl(U_t+B_t^{\T}U_t^{-1}B_t\bigr)-(m-1)V^{-1}=\widetilde D .
\end{equation}
For $p=3$, $m=2$ it has $18$ unknowns instead of the $72$ of the cavity
system.

\begin{lemma}\label{app:cert:lem-windmill}
Let $G$ consist of $m$ triangles $0v_1^tv_2^t$ sharing the vertex $0$, with
uniform dimension $p$. Given $(\widetilde D,B)$, define coefficients on $G$
by $D_0=\widetilde D$, $D_{v_1^t}=D_{v_2^t}=0$, $A_{0v_1^t}=A_{v_1^tv_2^t}=I$
and $A_{v_2^t0}=-B_t$. Let $(U,V)$ solve \eqref{app:cert:eq-rose} with
$U_t$, $V$ and $W_t=V^{-1}+B_t^{\T}U_t^{-1}B_t$ invertible. Then
\begin{gather*}
Q_{0\to v_1^t}=U_t,\qquad Q_{v_1^t\to v_2^t}=-U_t^{-1},\qquad Q_{v_2^t\to0}=U_t,\\
Q_{0\to v_2^t}=W_t,\qquad Q_{v_2^t\to v_1^t}=-B_tW_t^{-1}B_t^{\T},\qquad
Q_{v_1^t\to0}=B_t^{-\T}W_tB_t^{-1}
\end{gather*}
is an admissible solution of \eqref{blc:eq:cavity} at $z=0$, distinct
solutions of \eqref{app:cert:eq-rose} give distinct solutions, and the
solution is nondegenerate if the Jacobian matrix of \eqref{app:cert:eq-rose}
is invertible at $(U,V)$.
\end{lemma}

\begin{proof}
We drop the index $t$ where possible. The equations at the internal arcs
hold by direct substitution: $Q_{v_1\to v_2}=0-A_{v_10}U^{-1}A_{0v_1}$,
$Q_{v_2\to0}=0-A_{v_2v_1}Q_{v_1\to v_2}^{-1}A_{v_1v_2}$, and similarly in the
other direction. The messages into $0$ are
$A_{0v_2}Q_{v_2\to0}^{-1}A_{v_20}=B^{\T}U^{-1}B$ and
$A_{0v_1}Q_{v_1\to0}^{-1}A_{v_10}=BW^{-1}B^{\T}$. From the first equation of
\eqref{app:cert:eq-rose}, $UV(U-V^{-1})=BVB^{\T}$, and one checks
the identity $UVB=B+(UVU-U)U^{-1}B$; by the first equation its right-hand
side is $B+BVB^{\T}U^{-1}B=BVW$, and $UVB=BVW$ is equivalent to
$W=V^{-1}+B^{\T}U^{-1}B=B^{\T}(U-V^{-1})^{-1}B$; hence
$BW^{-1}B^{\T}=U-V^{-1}$. The vertex precision at $0$ is therefore
$R_0=\widetilde D-\sum_t(U_t-V^{-1}+B_t^{\T}U_t^{-1}B_t)=V^{-1}$ by the second
equation, and the equations at the arcs $0\to v_1^t$ and $0\to v_2^t$ read
$U_t=V^{-1}+B_tW_t^{-1}B_t^{\T}$ and $W_t=V^{-1}+B_t^{\T}U_t^{-1}B_t$, which
hold. The vertex precisions are $R_0=V^{-1}$, $R_{v_1}=-U^{-1}+B^{-\T}WB^{-1}
=B^{-\T}V^{-1}B^{-1}$ and $R_{v_2}=U-BW^{-1}B^{\T}=V^{-1}$, all invertible,
and the cavity precisions are invertible; hence so are the edge precisions,
and the solution is admissible. Distinct $(U,V)$ give distinct $U_t=Q_{0\to v_1^t}$
or $R_0^{-1}=V$.

Nondegeneracy. The equation of each internal arc has the form ``unknown $=$
function of one earlier unknown'' along the chains
$0\to v_1\to v_2\to0$ and $0\to v_2\to v_1\to0$; eliminating these unknowns
does not change the Jacobian determinant (the block of the internal arcs is
the identity plus a nilpotent matrix). There remain the $2m$ unknowns
$U_t=Q_{0\to v_1^t}$, $W_t=Q_{0\to v_2^t}$ and the equations
$U_t=\widetilde D-\sum_s\mu_s+B_tW_t^{-1}B_t^{\T}$,
$W_t=\widetilde D-\sum_s\mu_s+B_t^{\T}U_t^{-1}B_t$, with
$\mu_s=B_sW_s^{-1}B_s^{\T}+B_s^{\T}U_s^{-1}B_s$. Adjoin the unknown $R$ and the
equation $E_0\colon R-\widetilde D+\sum_s\mu_s=0$, which has derivative $I$ in
$R$, and subtract $E_0$ from the other equations. The equations
$W_t-R-B_t^{\T}U_t^{-1}B_t=0$ eliminate $W_t=\omega_t(U_t,R)$ in the same way.
What remains is
$\Phi_t=U_t-R-B_t\omega_t^{-1}B_t^{\T}=0$ and $E_0=0$; adding $\sum_t\Phi_t$ to
$E_0$ gives
$\Psi=\sum_t(U_t+B_t^{\T}U_t^{-1}B_t)-(m-1)R-\widetilde D=0$. After the change of
variables $R=V^{-1}$, a direct computation gives
\[
U_tVU_t-U_t-B_tVB_t^{\T}=(U_tV)\,\Phi_t\,(B_t^{-\T}\omega_tVB_t^{\T})
\]
identically near the solution, with invertible outer factors. At a common
zero the differentials of the two left-hand sides differ by these invertible
factors, so the two square systems have Jacobian matrices with the same
kernel. All operations were invertible, so the Jacobian matrix of the cavity
system is invertible if and only if that of \eqref{app:cert:eq-rose} is.
\end{proof}

The identity displayed in the proof follows from multiplying
$\Phi_tB_t^{-\T}\omega_t=(U_t-R)B_t^{-\T}R-RU_t^{-1}B_t$ on the left by
$U_tR^{-1}$ and on the right by $VB_t^{\T}$; we also checked it in exact arithmetic for random rational
$2\times2$ and $3\times3$ matrices.

\emph{Certificates for \eqref{app:cert:eq-rose}.} The entries of
$\widetilde D$ on and above the diagonal and all entries of $B_t$ are
given as binary64 numbers. The lower triangle of $\widetilde D$ is
obtained by transpose. Each candidate is refined by Newton's method
with $160$ bits. An approximate inverse $Y$ of the Jacobian is computed
in $160$-bit floating point. The dyadic center and radii of each box define
$X=\tilde x+r([-1,1]+i[-1,1])$, with radii computed from
$10^{-24}(1+\abs{\tilde x_j})$. The arithmetic (I2) encloses
$\appcK(X)-\tilde x$, which is compared with the exact radius intervals
in every coordinate. On $X$ the enclosures of $\det U_t$, $\det V$ and
$\det W_t$ exclude $0$. The boxes are pairwise disjoint. The results are
$36$ of $36$ candidates for $p=2$ and $1079$ of $1079$ for $p=3$.
By Lemma~\ref{app:cert:lem-windmill} and
Corollary~\ref{app:cert:cor-certificate}(d), $N_3(\text{bowtie})\ge1079$ and,
again, $N_2(\text{bowtie})\ge36$. The candidates for $p=3$ were found by a
floating point monodromy in the system \eqref{app:cert:eq-rose} that stopped
after six loops without new solutions; a second run, stopped after 30
minutes, had found $1069$. Whether $N_3(\text{bowtie})=1079$ is not known.
The same test certifies $310$ solutions for three triangles sharing a
vertex with $p=2$; this bound is not used in the paper.

\subsubsection{Exact count for the bowtie}\label{app:cert:item-bowtie-upper}

Theorem~\ref{app:bowtie:thm-36} uses exact arithmetic at two coupling
points, listed in Subsection~\ref{app:bowtie:subsec-comp}: the resultant at
$B_0$ and the single-loop polynomials and Jacobians; the admissibility and
transversality checks at $B_0$, which use the linear subresultant and gcds
modulo $p_{20}$; and the genericity checks at $P_0$ for excluding irregular
branches. These computations, over $\mathbb Q$ with SymPy, supply all
computational inputs to the equality $N_2=36$.

\subsubsection{Four parallel edges}\label{app:cert:item-four-port}

\emph{Statement.} Proposition~\ref{prop:block-four-port}(ii). The upper bound
$33-19$ is proved in the text from the bidegrees and the supports of
$\mathcal E,\mathcal F$; the lower bound uses one witness.

\emph{Computation} (exact, SymPy). (1) With $\alpha=a_1u+a_2v+a_3$ and
$\eta=k_1uv+k_2u+k_3v+k_4$ for indeterminates $a_i,k_i$, with
$\Delta=uv-\mu^2$ and $S=u+v-c$ for indeterminates $\mu^2,c$, the
polynomials $\mathcal E$, $\mathcal F$ of the proof are expanded identically;
one records, for each power of $u$, the range of powers of $v$, and
checks that the leading coefficients in $u$ are $v(k_1v+k_2)$ and
$2v^3(k_1v+k_2)(a_1+v^2-cv+\mu^2)$. These are the supports and the root
$v_0$ used in the text. (2) At the witness $\Theta=I_2$, $d=17$,
$b_i=(1,0),(2,1),(1,3),(2,2)$ one checks that the contents of
$\mathcal E$ and $\mathcal F$ with respect to $u$ are coprime, computes the subresultant
sequence with respect to $u$, and factors the resultant over $\mathbb Q$: its
irreducible factors have degrees $1,1,1,1,4,14$ with multiplicities
$1,2,2,3,2,1$. In the residue field of the factor of degree $14$ the
subresultant of degree one in $u$ gives $u$; then $t$, $x_1$, $x_2$ are
recovered, and one checks in this field that the four equations
\eqref{blc:eq:four-port} vanish, that $\mathcal H$, $\sigma$, all $y_i$,
$\Delta$, $\eta$ and $\Lambda_u$ are nonzero, and that the Jacobian
determinant of the four equations is nonzero. Since the factor is
irreducible over $\mathbb Q$, these hold at all $14$ complex roots, which
are distinct. (3) The resultant of $\Lambda_u,\Pi_u$ with respect to $u$ is a
nonzero constant times $v(5v^4+111v^3+848v^2+7424v+49152)$; in the residue
field of the quartic $\Lambda_u=\Pi_u=0$, $\Delta\eta\ne0$ and the Jacobian
determinant of $(\Lambda_u,\Pi_u)$ is nonzero; at $(u,v)=(-8,-8)$ and
$(-6,-\tfrac{32}3)$ one has $\Delta=0$, $\eta\ne0$ and
$\mathcal E=\mathcal F=\partial_u\mathcal F=\partial_v\mathcal F=0$.

\emph{Certified.} At the witness there are $14$ simple admissible solutions,
and the base points used in the upper bound are present with the stated
local properties. The extension to $p>2$ by zeros is by hand.

\subsubsection{Positive solutions on \texorpdfstring{$K_4$}{K4}}\label{app:cert:item-K4}

Examples~\ref{blr:ex-K4} and~\ref{blc:ex:k4} are proved by hand. We also
checked the instance
$p=2$, $r=(\tfrac32,\tfrac53)$, $D_i=\diag(\tfrac{17}6,\tfrac{43}{15})$, $A_e=I_2$
in exact rational arithmetic: for the four choices $Q=\diag(s,t)$,
$s\in\{\tfrac32,\tfrac43\}$, $t\in\{\tfrac53,\tfrac65\}$, the identity
$D-2Q^{-1}=Q$ and $s,t>1$; the characteristic polynomial of the
non-backtracking matrix of $K_4$ on its $12$ arcs,
$(\lambda-2)(\lambda-1)^3(\lambda+1)^2(\lambda^2+\lambda+2)^3$; the
nonvanishing of $\det(I-B/(s_ks_l))$ for the three entries $(k,l)$; and the
values $4a_k-12=-\tfrac23,-\tfrac8{15}$ of the quadratic form of the cone
condition. The identity \eqref{blc:eq:extended} and the
factorization of $\det S_e$ were checked for one edge with dimensions $2$
and $3$ at an explicit rational point. These are checks of hand proofs.

\subsubsection{Many positive solutions}\label{app:cert:item-many-solutions}

\emph{Statement.} Proposition~\ref{app:many:prop-certified}, quoted after
Proposition~\ref{real:prop:unbounded}: the model on $27$ vertices and $34$
edges has exactly $6$ positive solutions, with $(c_0,c_1,c_2)=(1,3,2)$, and
the model on $68$ vertices and $100$ edges has exactly $10$, with
$(c_0,c_1,c_2)=(1,5,4)$, at $z=0$.

\emph{Input.} All $D_i=1$ and rational weights. The graph consists of the
vertices $a=0$, $b=1$, the edge $ab$, and cycles through $a$ or through $b$
that are disjoint otherwise. For the $27$-vertex model the cycles have
lengths $3,3,4,5$ at $a$ and $4,5,4,5$ at $b$, $w_{ab}=\tfrac{1868}{785871}$,
and the $33$ cycle weights are listed in Appendix~\ref{app:many}. For the $68$-vertex model all cycles are
triangles $zuvz$ with $w_{zu}=w_{vz}=p$ and $w_{uv}=\tfrac12$; at $a$ there are
$11$ triangles with $p=\tfrac{130923}{870130}$ and $6$ with
$p=\tfrac{10788}{839159}$, at $b$ there are $2$, $1$, $4$, $9$ triangles with
$p=\tfrac{276918}{802963},\tfrac{35431}{637504},\tfrac{15763}{362293},\tfrac{11425}{892976}$,
and $w_{ab}=\tfrac{671}{876039}$.

\emph{Complete count} (exact and (I3)). One computes, in exact
rational arithmetic, the parameters of the reduction of
Appendix~\ref{app:many} for every cycle, and checks the conditions under
which the reduction is valid. The positive solutions then correspond to the
zeros of an explicit real function $h$ of $s=V_a$ on the set where a second
explicit function $\Phi_A$ exceeds $1$, the simple zeros to the
nondegenerate solutions, and the sign of $h'$, together with the sign of a
second explicit quantity when $h'>0$, gives the Morse index. With the
arithmetic (I3) one checks at two rational points
$s_{\mathrm{lo}}<s_{\mathrm{hi}}$ that $\Phi_A<1$, $\Phi_A'(s_{\mathrm{lo}})>0$ and
$\Phi_A'(s_{\mathrm{hi}})<0$; by concavity of $\Phi_A$ there is no zero outside
$(s_{\mathrm{lo}},s_{\mathrm{hi}})$. This interval is bisected until every piece
$I$ satisfies one of: $\sup\Phi_A(I)\le1$; $0\notin h(I)$; $\Phi_A>1$ on $I$,
$0\notin h'(I)$ and $h$ has the same sign at both ends (no zero); or
$\Phi_A>1$ on $I$, $0\notin h'(I)$ and $h$ has opposite signs at the ends
(exactly one zero, simple). No piece remained undecided, and on every piece
with a zero the index test was decided. The output is $6$ zeros with indices
$0,1,1,1,2,2$, resp.\ $10$ zeros with indices $0$ (once), $1$ (five times) and
$2$ (four times).

\emph{Independent check of every root} (arithmetic (I1)). The full system is
\eqref{irr:eq:f-system} at $z=0$, $D=1$, that is $f_e(1+f_e)=w_e^2g_ig_j$ with
$g_i=V_i=1+\sum_{e\ni i}f_e$, in the $\abs E$ unknowns $f$. For each root a
box with relative radius $10^{-12}$ to $10^{-8}$ passes the test of
Lemma~\ref{app:cert:lem-krawczyk} (with an equilibrated approximate inverse $Y$ if the plain one fails), $f>0$ on the box (smallest lower end $5\cdot10^{-4}$ for
the $27$-vertex model), the boxes are pairwise disjoint (compared on their
floating point end points), and the value of $V_a$ on the box lies in the
interval of the one-dimensional count. The Morse index is certified with
Lemma~\ref{app:cert:lem-inertia} for the matrix $H_u$ of
\eqref{real:eq:Hu}: since $f^2/(1+2f)$ and $f(1+f)/(1+2f)$ are increasing for
$f>0$, their values at the ends of the interval of $f_e$ enclose the entries;
$Q$ is the orthogonal matrix of eigenvectors of the midpoint computed in
floating point, $\mathbf M$ is enclosed by (I1), and the margins
$\min\abs{d_i}-e$ lie between $0.0025$ and $0.31$. The certified indices
agree with those of the one-dimensional count. As a further check, the
inertia of $H_u$ at the midpoint of every box was recomputed by an exact
$LDL^{\T}$ decomposition over $\mathbb Q$.

\subsubsection{Circles}\label{app:cert:item-circles}

\begin{proposition}[Circles at specified coefficients]\label{app:bridge:prop-M}
Consider the following two models.
\begin{enumerate}
\item[\textup{(M1)}] The multigraph with vertices $0,1,2,3$, three
parallel edges $01$ with weights $1,\tfrac78,\tfrac98$, three parallel
edges $23$ with the same weights, and an edge $12$ of weight
$\tfrac1{16}$; $D=(0,0,\tfrac5{64},\tfrac5{64})$; $\eta=10^{-9}$.
\item[\textup{(M2)}] The simple graph formed by two diamonds with edges
$01,02,03,12,13$, joined by an edge between their vertices $2$, of weight
$\tfrac1{32}$; the first diamond has
$64D=(75,94,167,139)$ and $64w=(-38,-60,-63,-42,79)$, the second has
$64D=(76,95,168,140)$ and $64w=(-39,-58,-63,-43,77)$, in the order of the
edges above; $\eta=10^{-12}$.
\end{enumerate}
For every perturbation of at most $\eta$ in each coefficient, $\Xp$ has a
smooth circle of regular points. The lower bound on every variance and the
enclosure of the energy on the certified circle are
\[
\begin{array}{c|cc}
 &\text{lower bound on } V_i & \text{enclosure of }z\\\hline
 \mathrm{M1}&0.9485&[-2.91740,-2.83794]\\
 \mathrm{M2}&0.8543&[-0.54680,-0.53480].
\end{array}
\]
\end{proposition}

\begin{proof}
The certificate uses arithmetic (I1) with the coefficients
enclosed in intervals of radius $\eta$. It works with
$G_i=(D_i-z)V_i-1-\sum_{e\ni i}f_+(w_e^2V_iV_j)$ and evaluates
$f_+(y)$ as $2y/(\sqrt{1+4y}+1)$. A numerically traced loop supplies a
cyclic list of boxes in affine coordinates. In each box one coordinate
$p$ of $(V,z)$ varies in an interval $P$ and all others have the form
$q=\ell(p)+u$ with $\ell(p)=\hat q+s(p-\hat p)$, $\hat p\in P$ and
$u\in U=[-r,r]^n$.
The following checks certify the loop for each fixed coefficient vector
in the coefficient box.

\begin{enumerate}
\item Each box has positive variances. A mean value enclosure on the
line $L=\{(p,\ell(p)):p\in P\}$ is
$G(\hat p,\hat q)+(J_q(L)s+J_p(L))(P-\hat p)$.
Using it for $G(p,\ell(p))$ and enclosing $J_q$ over the whole box,
the Krawczyk image for $u$ at center $0$ is strictly inside $U$.
Lemma~\ref{app:cert:lem-krawczyk}, applied for each $p$, gives a unique
zero with invertible $J_q$. By the implicit function theorem these
zeros form one regular arc, with its interior in the interior of the box.
\item At each endpoint $p_*$, evaluate $G(p_*,\ell(p_*))$ directly
by interval arithmetic. Apply four Krawczyk steps with center $u=0$,
starting from $U_0=U$. Each $U_{j+1}$ is the symmetric box about $0$
enclosing the intersection of $U_j$ with its Krawczyk image. The true
endpoint remains in the intersection by the mean value formula. The
enclosure from the fourth intersection lies strictly inside the
neighboring box, checked in that box's affine coordinates.
\end{enumerate}

The checks pass on all $2228$ boxes for (M1) and $21916$ for (M2).
Thus the union of the arcs is connected, compact, and covered by their
open chart interiors, including at the junctions. By
Lemma~\ref{app:bridge:lem-glue} it is a compact component of $\Xp$; full
rank in each chart makes it a smooth circle. The coordinate bounds in
the proposition enclose the boxes.
\end{proof}

The certificate gives a circle in each model. Floating-point continuation
at $z=-2.88$ in (M1) and $z=-0.54$ in (M2) finds four positive roots with
indices $0,1,1,2$; two belong to the traced circle. No count of all circles
or roots for these fixed coefficient vectors follows from the certificate.
With $\eta=10^{-9}$ the subdivision for (M2) did not terminate near the
lower turning point; only the smaller radius was certified.

\subsubsection{Genus table}\label{app:cert:item-genus}

\emph{Statement.} The values of $\LG$ and $\Sigma(G)$ in
Table~\ref{gen:tab:values}; $g$ and $R$ follow from
Theorem~\ref{thm:genus-formula} and the Riemann--Hurwitz formula
$R=2g-2+2\LG$.

\emph{Computation} (exact integer enumeration). One enumerates all
nonempty $U\subseteq V$ with $G[U]$ connected and all frozen values on the
edges at $N(U)$, and evaluates \eqref{gen:eq:Sigma} for the cycles, dipoles
$D_1$ to $D_6$, theta graphs, cycles sharing a vertex, three triangles
sharing a vertex, the dumbbell, $K_4$, $K_4$ with one edge subdivided, the
prism, $K_{2,3}$, $M_1$ and $M_2$. A second, independent implementation
written directly from the definitions in Section~\ref{sec:genus}
reproduces $(\LG,\Sigma,g)$ for $D_3$, $D_4$, the diamond,
the bowtie, $K_4$, the dumbbell, the prism, $M_1$ and $M_2$, and gives
$\Sigma(G)=-2n$ for $84$ random unicyclic graphs with attached trees. The
closed forms in the rows with parameters are proved in Section~\ref{sec:genus};
the enumeration confirms them for small parameters. The survey of
ramification points reported after the table is floating point evidence.

\subsubsection{The residue identity on small graphs}\label{app:cert:item-census}

\emph{Statement.} Remark~\ref{cls:rem:census}. An exact computation evaluates both sides
of Corollary~\ref{cls:cor:identity}, with $M(U,F)$ computed from the numbers
of nondegenerate boundary roots and of leafless subsets with prescribed
forbidden degrees, for all $1+2+8+64+1024+32768=33867$ labeled simple graphs
with at most six vertices (including disconnected graphs, trees and isolated
vertices) and for $300$ random multigraphs with at most six vertices and nine
edges; all agree. Corollary~\ref{cor:local-identity} was checked in the same
way for all $5405$ pairs of a labeled simple graph with at most five
vertices and a vertex. This is a check of a proved statement and is not
used in any proof.

\subsubsection{Leafless subsets and binary fixed points}\label{app:cert:item-L}

\emph{Statement.} The values of $\LG$ and $Q(G)$ in Table~\ref{deg:table} and
of $L(K_n)$, $n\le7$. They are obtained by enumerating all
$2^{\abs E}$ edge subsets and evaluating \eqref{deg:L-sum} and
\eqref{deg:Q-sum} exactly. The rows for trees and for $\beta=1$ are immediate (the
leafless subsets are $\emptyset$, resp.\ $\emptyset$ and the cycle), and so
are those for the theta graphs and the dumbbell: the leafless subsets are
$\emptyset$, the three cycles and $E$, resp.\ $\emptyset$, the two cycles,
their union and $E$, which gives $\LG=5$ and $Q=1+1+1+1+4=8$.

%% file: sections/declarations.tex
\section*{Declarations}
Large language model tools were used in exploring candidate statements,
in drafting, reorganizing, and checking the proofs and the text, in the
literature search, in writing the Lean formalization, and in preparing
the exact and interval computations referred to as computer-assisted. The
author assumes responsibility for all content, including the definitions,
the written proofs, the correspondence between the written and the formal
statements, and the citations.